\documentclass{gtmon_a}
\pdfoutput=1

\usepackage{pinlabel}
\usepackage[matrix, arrow, curve]{xy}
\usepackage{stmaryrd}
\input{diagxy}

\proceedingstitle{The interaction of finite-type and Gromov--Witten
invariants (BIRS 2003)}
\conferencestart{15 November 2003}
\conferenceend{20 November 2003}
\conferencename{The interaction of finite-type and Gromov--Witten
invariants}
\conferencelocation{Banff International Research Station, Banff, Alberta,
Canada}

\editor{David Auckly}
\givenname{David}
\surname{Auckly}

\editor{Jim Bryan}
\givenname{Jim}
\surname{Bryan}

\title{Introduction to the Gopakumar--Vafa Large $N$ Duality}

\author{Dave Auckly}
\givenname{Dave}
\surname{Auckly}
\address{Department of Mathematics\\
Kansas State University\\\newline
Manhattan KS 66506\\
USA}
\email{dav@math.ksu.edu}
  \urladdr{http://www.math.ksu.edu/~dav/}

\author{Sergiy Koshkin}
\givenname{Sergiy}
\surname{Koshkin}
\address{Department of Mathematics\\
Northwestern University\\\newline
2033 Sheridan Road\\
Evanston IL 60208-2730\\
USA}
\email{koshkin@math.northwestern.edu}
\urladdr{}

\volumenumber{8}
\issuenumber{}
\publicationyear{2006}
\papernumber{9}
\startpage{195}
\endpage{456}

\doi{}
\MR{}
\Zbl{}

\arxivreference{math.GT/0701568}


\received{19 May 2006}
\revised{}
\accepted{}
\published{21 September 2007}
\publishedonline{21 September 2007}
\proposed{}
\seconded{}
\corresponding{}
\version{}

\let\xysavmatrix\xymatrix
\def\xymatrix{\disablesubscriptcorrection\xysavmatrix}
\AtBeginDocument{\let\bar\wbar\let\tilde\wtilde\let\hat\what\let\wcheck\check}

\renewcommand{\to}{\rightarrow}
\let\barrsquare\square
\let\square\undefined

\makeop{val}
\makeop{vir}
\makeop{virdim}
\makeop{Aut}
\makeop{inv}
\makeop{Mor}
\makeop{Ob}
\makeop{Hom}
\newcommand{\Ext}{\mathbb{E}\mathrm{xt}}
\newcommand{\II}{\mathit{II}}
\newcommand{\III}{\mathit{III}}
\newcommand{\IV}{\mathit{IV}}
\newcommand\sbinom[2]{\genfrac{[}{]}{0pt}{}{#1}{#2}}
\newcommand\pling{!}
\newcommand\bthin{\!}
\def\underC{{}\mskip1mu\underline{\mskip-1mu C \mskip-2mu}\mskip2mu}
\def\underCC{{}\mskip1mu\underline{\mskip-1mu \C \mskip-1mu}\mskip1mu}
\def\underCP{{}\mskip2mu\underline{\mskip-2mu \CP \mskip-2mu}\mskip2mu}
\def\underV{{}\mskip2mu\underline{\mskip-2mu V \mskip-2mu}\mskip2mu}
\def\underT{{}\mskip2mu\underline{\mskip-2mu T \mskip-2mu}\mskip2mu}
\def\underN{{}\mskip1mu\underline{\mskip-1mu N \mskip-5mu}\mskip5mu}
\def\underM{{}\mskip1mu\underline{\mskip-1mu M \mskip-5mu}\mskip5mu}
\def\underMM{{}\mskip2mu\underline{\mskip-2mu \mathcal{M} \mskip-2mu}\mskip2mu}

\makeatletter
\def\cnewtheorem#1[#2]#3{\newtheorem{#1}{#3}[section]
\expandafter\let\csname c@#1\endcsname\c@thm}

\newtheorem{thm}{Theorem}[section]
\cnewtheorem{lemma}[thm]{Lemma}
\cnewtheorem{prop}[thm]{Proposition}
\cnewtheorem{cor}[thm]{Corollary}
\cnewtheorem{conjecture}[thm]{Conjecture}

\theoremstyle{definition}
\cnewtheorem{aside}[thm]{Aside}
\cnewtheorem{example}[thm]{Example}
\cnewtheorem{remark}[thm]{Remark}
\cnewtheorem{exm}[thm]{Exercise}
\cnewtheorem{defn}[thm]{Definition}
\makeatother  
\makeautorefname{exm}{Exercise}
\makeautorefname{defn}{Definition}

\numberwithin{figure}{section}
\numberwithin{table}{section}

\newdimen\tableauside\tableauside=0.8ex
\newdimen\tableaurule\tableaurule=0.3pt
\newdimen\tableaustep
\def\phantomhrule#1{\hbox{\vbox to0pt{\hrule height\tableaurule width#1\vss}}}
\def\phantomvrule#1{\vbox{\hbox to0pt{\vrule width\tableaurule height#1\hss}}}
\def\sqar{\vbox{%
  \phantomhrule\tableaustep
  \hbox{\phantomvrule\tableaustep\kern\tableaustep\phantomvrule\tableaustep}%
  \hbox{\vbox{\phantomhrule\tableauside}\kern-\tableaurule}}}
\def\squares#1{\hbox{\count0=#1\noindent\loop\sqar
  \advance\count0 by-1 \ifnum\count0>0\repeat}}
\def\tableau#1{\vcenter{\offinterlineskip
  \tableaustep=\tableauside\advance\tableaustep by-\tableaurule
  \kern\normallineskip\hbox
    {\kern\normallineskip\vbox
      {\gettableau#1 0 }%
     \kern\normallineskip\kern\tableaurule}%
  \kern\normallineskip\kern\tableaurule}}
\def\gettableau#1 {\ifnum#1=0\let\next=\null\else
  \squares{#1}\let\next=\gettableau\fi\next}

\def\CP{{\mathbb{C}\mathrm{P}}}

\def\l{\text{\ell}}

\def\baru{{\bar {u}}}

\def\barz{{\bar {z}}}
\def\dbar{{\bar {\partial}}}
\let\mathds\mathbb
\def\bbone{{\mathds{1}}}

\def\varlambda{\lambda}
\def\ve{\varepsilon}
\def\eps{\epsilon}
\def\g{\mathfrak{g}}
\def\h{\mathfrak{h}}
\def\SL{\mathfrak{sl}}
\def\Sl{\mathfrak{sl}}
\def\Sum{\sum\nolimits}
\def\Prod{\prod\nolimits}
\def\Tilt{\overline{\mathcal{T}\!\mathit{ilt}}_\epsilon}
\def\End{\mathop{\rm End}\nolimits}
\def\Tr{\mathop{\rm Tr}\nolimits}
\def\id{\mathop{\rm id}\nolimits}
\def\rotimes{\wwbar{\otimes}}
\def\AW{\wwtilde{W}^l}

\def\CC{{\mathbf C}}
\def\DD{{\mathbf D}}

\def\calI{{\mathcal I}}
\def\calA{{\mathcal A}}
\def\calB{{\mathcal B}}

\def\calD{{\mathcal D}}
\def\calL{{\mathcal L}}
\def\calM{{\mathcal M}}
\def\calO{{\mathcal O}}
\def\calQ{{\mathcal Q}}

\def\calU{{\mathcal U}}
\def\calV{{\mathcal V}}

\def\calW{{\mathcal W}}
\def\A{{\mathbb A}}
\def\Q{{\mathbb Q}}
\def\C{{\mathbb C}}

\def\R{{\mathbb R}}

\def\E{{\mathbb E}}
\def\F{{\mathbb F}}
\def\H{{\mathbb H}}
\def\Z{{\mathbb Z}}

\def\tildeS{\wwtilde{\Sigma}}
\def\l{{\ell}}

\def\lra{\longrightarrow}
\def\disj{\coprod}

\makeindex

\begin{document}

\begin{htmlabstract}
<p class="noindent">
Gopakumar&ndash;Vafa Large N Duality is a correspondence between
Chern&ndash;Simons invariants of a link in a 3&ndash;manifold and relative
Gromov&ndash;Witten invariants of a 6&ndash;dimensional symplectic manifold
relative to a Lagrangian submanifold. We address the correspondence
between the Chern&ndash;Simons free energy of S<sup>3</sup> with no link and the
Gromov&ndash;Witten invariant of the resolved conifold in great detail.
This case avoids mathematical difficulties in formulating a
definition of relative Gromov&ndash;Witten invariants, but includes all of
the important ideas.
</p>
<p class="noindent">
There is a vast amount of background material related to this
duality. We make a point of collecting all of the background
material required to check this duality in the case of the
3&ndash;sphere, and we have tried to present the material in a way
complementary to the existing literature. This paper contains a
large section on Gromov&ndash;Witten theory and a large section on quantum
invariants of 3&ndash;manifolds. It also includes some physical
motivation, but for the most part it avoids physical terminology.
</p>
\end{htmlabstract}

\begin{abstract}   
Gopakumar--Vafa Large $N$ Duality is a correspondence between
Chern--Simons invariants of a link in a 3--manifold and relative
Gromov--Witten invariants of a 6--dimensional symplectic manifold
relative to a Lagrangian submanifold. We address the correspondence
between the Chern--Simons free energy of $S^3$ with no link and the
Gromov--Witten invariant of the resolved conifold in great detail.
This case avoids mathematical difficulties in formulating a
definition of relative Gromov--Witten invariants, but includes all of
the important ideas.

There is a vast amount of background material related to this
duality. We make a point of collecting all of the background
material required to check this duality in the case of the
3--sphere, and we have tried to present the material in a way
complementary to the existing literature. This paper contains a
large section on Gromov--Witten theory and a large section on quantum
invariants of 3--manifolds. It also includes some physical
motivation, but for the most part it avoids physical terminology.
\end{abstract}

\maketitle

\tableofcontents

\newpage

\part*{Introduction}
\setobjecttype{Part}

Large $N$ duality is a conjectural correspondence between two very
different types of mathematical objects: the large $N$ limit of a
gauge theory with structure group $U(N)$ and a string theory. Since
gauge theories and string theories are both meant to describe the
same universe it is natural to expect a correspondence between their
predictions. There are several examples in physics and mathematics literature of this apparent
duality. The original description of this correspondence goes back
to the 1974 theoretical physics paper \cite{tH} by 't Hooft. At the time the main
computational tool in both theories was perturbative expansion and
't Hooft noticed intriguing similarities that occur in those
expansions. One of the first mathematical papers related to a Large $N$ Duality is a 1992 paper
where Kontsevich introduced a matrix model to resolve the
conjectures of Witten about $2$--dimensional gravity \cite{KontAiry}.

In oversimplified terms, gauge theory studies moduli spaces of
connections on principal bundles while string theory studies spaces
of maps from a particular class of domains into different targets.
Both theories lead to invariants (of the base manifold in the gauge
case and the target manifold in the string case).  These invariants
can be conveniently assembled into generating functions (called
'partition functions' by physicists). 't Hooft was considering
expansions for gauge theories with SU($N$) structure groups and
noticed that as $N\to\infty$ they turn into partition functions one
expects from a string theory.

In this generality 't Hooft's principle remains beyond the reach of
mathematical theory for the foreseeable future. Ten years after
Kontsevich's paper people began to understand the Large $N$ Duality
relating the Chern--Simons SU($N$) gauge theory of $3$--manifolds to
the Gromov--Witten string theory of complex $3$--folds.

The aspect of Large $N$ Duality that we address in this survey is
a duality between Chern--Simons theory and Gromov--Witten
theory, the Gopakumar--Vafa Large $N$ Duality duality. This duality states that the
Chern--Simons (Reshetikhin--Turaev) invariants of a link in a
$3$--manifold are related to the relative Gromov--Witten invariants of
a $6$--dimensional symplectic manifold relative to a $3$--dimensional
Lagrangian submanifold. We address the correspondence between the
Chern--Simons free energy of $S^3$ with no link and the Gromov--Witten
invariant of the resolved conifold in great detail.

The key trait of both theories is that they are 'topological' in the
sense that they do not depend on a background metric on the manifold
in question. This greatly simplifies the setting and makes it
possible to explain the Large $N$ Duality in mathematical terms,
something that remains impossible for other examples.

Mathematically the most that can be done is to compute the invariant
on the Gromov--Witten side and compute the invariant on the
Chern--Simons side and compare the two answers. This is what we do.
The first part of this paper (\fullref{GW}) is an exposition
of Gromov--Witten theory up to the point of the computation of the
Gromov--Witten free energy of the resolved conifold (the full
multiple cover formula). The definition of Gromov--Witten invariants
is given in \fullref{coarse} with intuitive descriptions of some of
the more technical elements. Formal definitions are given in the following
subsections as they
are motivated by and required for ever more complicated sample
computations. We have included all of the relevant definitions. For
example, the definition of a stack is given in \fullref{app:a}. Most
expositions on Gromov--Witten theory avoid this definition.

One of the unique things that we do in this paper is a computation of
the genus zero, degree two Gromov--Witten invariant of the resolved
conifold directly from the definition; see \fullref{mc}. This case is
addressed via localization in the book \cite{CK} by Cox and Katz.

The second part of this paper (\fullref{CS}) is an
exposition of Chern--Simons theory up to the point of the computation
of the Chern--Simons free energy of the $3$--sphere. The Chern--Simons
invariants were motivated by a path integral expression. We outline
the progression from this heuristic definition to formal definitions
of perturbative invariants in Sections \ref{pnhdesc} and \ref{pCS}
and the first two parts of \fullref{TQFT}. The motivation for
introducing the free energy is explained in \fullref{pnhdesc}.

The skein theory approach is described in \fullref{app:c}; this is
the easiest way to describe the Chern--Simons invariants. It is
however very difficult to compute from the resulting expressions, so
we rely on the quantum group approach instead.

The main subsection in the second part is \fullref{TQFT}. It is here
that we motivate and give a formal definition of the quantum
group invariants. The second part of the paper ends with the
computation of the Chern--Simons partition function.

\fullref{part3} begins with \fullref{comparison} where we use
special function techniques to derive the formal relation between
the Gromov--Witten free energy and the Chern--Simons free energy. The
remainder of this part is overview and history.

To get a fast introduction to the Gopakumar--Vafa Large $N$ Duality
one may just read the overview or read the first definition of the
Gromov--Witten invariants from \fullref{coarse}, the skein theory
definition of the Chern--Simons invariants from \fullref{app:c} and
the comparison of the two in \fullref{comparison}. Some physical
intuition may be obtained from the description of the perturbative
expansion in \fullref{pCS}.


Chern--Simons theory is defined for real 3--manifolds while the
relevant Gromov--Witten theory is defined for Calabi--Yau complex
3--manifolds (in general one can define Gromov--Witten invariants for
arbitrary symplectic manifolds, see \fullref{GW}). Thus the
mathematically oriented reader can see Large $N$ Duality as an
interesting correspondence between 3--dimensional real and complex
geometries and topologies. Physically the importance of the
Calabi--Yau condition is that in string theory Calabi--Yau 3--folds
(ie 6--real dimensional manifolds) provide complementary
'compactified' dimensions to the 4 observed ones of the classical
space-time.

The existing literature on Large $N$ Duality is vast and is growing
exponentially so it would be impossible to survey it here. However, it
appears to fall mostly into two categories: one written by physicists
with extensive use of physical terminology (see the survey by Mari\~no
\cite{Mar} and references therein), another by or for algebraic
geometers (see Cox and Katz \cite{CK} and the second half of Hori et al
\cite{Hori}). This reflects the fact that the physical insight and the
complex side of the duality are at present the most developed parts.

One of the difficulties that impedes further progress in this
subject is the amount of background material required to comprehend
all the mathematical aspects of the conjectural duality (as one can
judge, for example, from the size of the background chapters in Hori et al
\cite{Hori} and Turaev \cite{T}).

We tried to provide a relatively self-contained introduction to
the existing ideas and methods involved in Large $N$ Duality that
was complementary to the existing literature. We fill in details
where we had trouble finding them and leave well-documented
computations as exercises.

In order to keep the size of the the paper manageable we stuck to
topics that could be formalized mathematically. In particular we
cover the duality between the $3$--sphere and the resolved conifold
without including any knots. We also provide a number of
computational and illustrative examples to make the matters clearer
to a non-specialist. It is hoped that this paper will be accessible
to advanced graduate students and will help to bring new blood into
the field. Exercises are spread all over the text along with
references to their solutions (or ideas for such). Although formally
not necessary to get through the paper they are important for those
who plan to acquire a working understanding of the subject matter.

The idea of writing this paper dates back to the workshop
Interaction of Finite-type and Gromov--Witten Invariants at the Banff
International Research Station in November 2003 co-organized by the
first author and Jim Bryan. M Mari\~no gave a mini-series of
lectures at this workshop on physics and mathematics of the
Gopakumar--Vafa conjecture (as the Large $N$ Duality was dubbed at
the time). The plan was to write up lecture notes accessible to
mathematics graduate students.  However, M Mari\~no wrote an
introduction that followed his lectures fairly closely \cite{Mar}
and it became clear that more details could not be given without
writing a fair amount of background material. The final version of
the paper emerged out of friendly discussions between the two
authors as we learned this material.

We would like to thank the directors and staff of the Banff
International Research Station for providing an amazing place to
share mathematics. We would like to thank Marcos Mari\~no for his
illuminating lectures and Arthur Greenspoon for his careful
editing of an earlier draft of this paper. The first author would
like to thank Jim Bryan for first telling him about this duality.
The second author would like to thank C-C\,M Liu and D Karp for
fruitful discussions and suggestions related to the content of this
work. Dave Auckly and Sergiy Koshkin were partially supported by NSF
grant DMS-0204651.

\section{Mathematical history of Large $N$ Duality}\label{02}


It is instructive to trace the development of elements involved in
the modern picture of Large $N$ Duality. Hopefully this will also
serve as a non-technical step-by-step introduction into the field.
The reader should keep in mind that this is a mathematician's take
on this task and a sporadic one at that. For instance, we almost
completely ignore the physical undercurrent of the process except
for a few landmark papers. A different perspective can be found in
the introductory parts of Grassi--Rossi \cite{GR}.

The history can be divided into four periods separated by physical
breakthroughs into the mathematical realm: prehistory (1974--1989),
formation of concepts and tools (1989--1998), the Gopakumar--Vafa
conjecture (1998--2003), life after the vertex (2003--present). The
dates in the text refer to arxiv submissions while references are
given wherever possible to journal publications.

This historical overview describes the state of the subject in 2003
before the seminal paper `The topological vertex' \cite{vertex}
re-shifted the perspective. Accordingly in this essay we only
address the developments in 1974--2003. New results and directions
that appeared after the November 2003 workshop are surveyed (or
rather sketched) in the last section of this paper.

\subsection{Prehistory (1974--1989)}\label{021}

When 't Hooft first formulated the Large $N$ Duality principle
neither  Chern--Simons theory nor  Gromov--Witten theory existed.
However, this period saw the appearance of many ingredients that
later fit into the picture.

In 1983 H Clemens in a paper called `Double Solids' (ie complex
threefolds) \cite{Cl} studied extremal transitions between
threefolds that include deformations of complex structure into a
singular one with subsequent resolution of the singularity. The
conifold transition that connects cotangent bundles of 3--manifolds
to their large $N$ dual threefolds is an example of such a
transition. The study of threefold singularities led to a conjecture
known as Reid's fantasy (1987) that the moduli space of Calabi--Yau
threefolds forms a single family related through such transitions
(see Grassi and Rossi \cite{GR} and Reid \cite{Rd}).

At about the same time, physicists realized that string theories remain
well-defined on varieties with certain singularities (see Dixon, Harvey,
Vafa and Witten \cite{DHVW}) and can therefore change smoothly as the
underlying manifolds undergo a singular transition. This idea that the
same theory can be described on topologically different bases is at the
heart of Large $N$ Duality.

On the other end, in 1985 V Jones introduced his polynomial
invariant of knots \cite{Jn}. However, his motivation came from
operator algebras and no connection to Chern--Simons theory was known
at the time. Already in 1985 several groups of authors generalized
the Jones polynomial. Six authors published a joint paper
\cite{HOMFLY}  giving the new HOMFLY polynomial invariant (named
with the initials of their last names). Two other authors \cite{PT}
operated behind the iron curtain and their work remained
unrecognized until somewhat later. To give full credit HOMFLY is now
sometimes expanded to HOMFLYPT or THOMFLYP. It was later discovered
that the Jones polynomial corresponds to the SU($2$) and THOMFLYP to
the SU($N$) Chern--Simons theories respectively.

\subsection{Formation of concepts and tools (1989--1998)}\label{022}

This period starts with the appearance of the first \cite{W1} of
three seminal papers by E Witten.  These papers revolutionized both
the Chern--Simons and the Gromov--Witten theories and then tied them
together. In `Quantum field theory and the Jones polynomial' Witten
introduced the (quantum) Chern--Simons theory as a gauge theory on
3--manifolds with Lagrangian density given by the Chern--Simons form
(see Witten \cite{W1} and Baez--Muniain \cite{BzM}). He then
`solved' it, that is, found explicit mathematical expressions for
expectations of observables by reducing the computation to conformal
field theory on Riemann surfaces (see Di
Francesco--Mathieu--S\'en\'echal \cite{DMS} or Appendix
\ref{app:e}). The observables turned out to be the so-called `Wilson
loops', that is, holonomies of connections over knots and links in
the manifold. Witten's heuristic computation showed that expectation
values of Wilson loops are the Jones and THOMFLYP polynomials with
minor renormalizations. One can find a more mathematical account of
Witten's ideas in M Atiyah's book \cite{At}.

Shortly after (1990--1991) Witten's results were put on a firm
mathematical basis via quantum groups and conformal field theory.
The quantum group approach was led by N Reshetikhin and V Turaev
\cite{RT1,RT2,KRT,T}, who redefined Witten's quantum invariants
using the machinery of quantum groups and modular tensor categories
(see \fullref{TQFT}).

One first had to study perturbative expansions of the invariants
rather than their exact expressions before the gauge/string duality
would become apparent. Coefficients of such expansions represent
another kind of knot invariants known as finite-type or Vassiliev
invariants introduced in 1990 by V Vassiliev from completely
different considerations. The connection between Chern--Simons theory
and Vassiliev invariants along this line was established in the 1991
PhD thesis of D Bar-Natan (see also the subsequent papers
\cite{bn0,bn1,bn2}). The idea was to apply the the `Feynman rules'
from quantum field theory to the 'path integral' on the space of
connections describing Chern--Simons theory (see Sections
\ref{pnhdesc} and \ref{pCS} for more details). Mathematically the
universal finite type invariant now known as the Kontsevich integral
was given in 1992 by M Kontsevich \cite{KontU1}.

Similar developments also occurred in the field of $3$--manifolds. In
1991 S Axelrod and I Singer formalized the Feynman integral construction
for perturbative Chern--Simons $3$--manifold invariants and in 1994
M Kontsevich demonstrated that his integral proved to be a universal
finite-type invariant \cite{KontU2}. Finite-type invariants for
$3$--manifolds were introduced by T Ohtsuki in 1996 \cite{Oht} and in
in 1998 T Le, J Murakami and T Ohtsuki gave a detailed construction of a
Kontsevich-type invariant for $3$--manifolds and proved its universality
among the Ohtsuki invariants for integral homology spheres \cite{LMO}. Now
known as the LMO invariant it is believed to capture the trivial
connection contribution to Witten's quantum invariant (see Rozansky
\!\cite{Roz} and Bar-Natan--Garoufalidis--Rozansky--Thurston \cite{bar}).

Rapid progress on the complex side of the conjecture was initiated by
Witten's 1991 paper `2D Gravity and Intersection Theory on Moduli Space'
\cite{2D}. In this paper Witten defines what are now called tautological
classes on the moduli space of stable algebraic curves and gives string,
dilaton and divisor equations that are sufficient to compute all of
the corresponding intersection numbers (see \fullref{GWco}, Hori et al
\cite{Hori} or Vakil \cite{vakil}). These intersection numbers could
be included in the framework of invariants of symplectic manifolds
introduced in M Gromov's 1985 paper \cite{Gromov} on pseudoholomorphic
curves (Gromov was more interested in topological applications).

The next year M Kontsevich provided a proof of Witten's conjectured
equations \cite{KontAiry} using a presentation of the moduli space by
ribbon graphs and reducing the generating function for the intersection
numbers to a matrix integral of Gaussian type \cite{KontU2}. This can be
considered the first mathematical instance of Large $N$ Duality. The
techniques of matrix models and graph combinatorics have become
indispensable in computations related to Large $N$ Duality.

In 1994 M Kontsevich generalized the definition of stable algebraic curves
to stable maps \cite{konloc} and paved the way to general definitions of
(closed) Gromov--Witten invariants introduced by various authors in 1996
(see the discussion in McDuff--Salamon \cite{MdS} or Cox--Katz \cite{CK}).
Another important achievement of this paper is the discovery of the
`localization', technique for computing Gromov--Witten invariants (see
\fullref{locGW}). The basic idea dates back to the Atiyah--Bott paper
\cite{AB} that shows how to compute an integral of functions equivariant
under a torus action by localizing to an integral over the fixed points
of the action. In 1997 this was generalized to `virtual localization'
by Graber and Pandharipande \cite{virloc}.

The third and final Witten paper that we ought to discuss
`Chern--Simons gauge theory as a string theory' was made available
in 1992 although the first printed version appeared only in 1995 in
the Floer Memorial Volume \cite{Wcss}. It contains an outline of the
first step in the Large $N$ Duality between the Chern--Simons gauge
theory and a theory of open strings, so-called `holomorphic
instantons at infinity'. The second and more complicated step in 't
Hooft's large $N$ program that involves transition from open to
closed strings had to wait until later.

As described above, quantum invariants of knots and links correspond
to expectation values in Chern--Simons theory. Witten took the next
step and along with a 3--manifold $M$ considered its cotangent bundle
$T^*M$ with its natural symplectic structure. This allows one to
define (pseudo-)holomorphic curves (stable maps) in $T^*M$ and
`count' them with Gromov--Witten invariants. For this idea to work it
is important  that the dual threefold is a Calabi--Yau so that
holomorphic curves (stable maps) are generically isolated and can be
counted. Holomorphic curves ending on the zero section of $T^*M$ are
supposed to be the holomorphic instantons of Witten's theory. Witten
argues moreover that the duality is exact, that is, it holds not only for
large but for all $N$. The last property is due to the topological
nature of the Chern--Simons theory. The partition function of this
theory corresponds to the generating function of the Gromov--Witten
invariants of holomorphic curves which is dual to the generating
function of Chern--Simons invariants of knots and links.

There is one major problem with this picture. A cotangent bundle is
a Calabi--Yau but a very degenerate one, in particular there are no
non-constant holomorphic curves there either closed or ending on the
zero section. This circumstance was well known to Witten and is
explicitly pointed out in \cite{Wcss}. Recall that in 1992 there was
no notion of Gromov--Witten invariants and even the modern
`translation' of holomorphic instantons as stable maps is
incomplete. The physical notion is broader and includes objects like
framed trivalent graphs, in particular framed knots and links in the
zero section, a.k.a. instantons at infinity. One can think of them
as made of infinitely thin ribbons thus representing degenerate
Riemann surfaces with boundary. Even today we do not have a grand
theory that would incorporate such degeneracies. Some potential
contenders have emerged recently (for example symplectic field
theory, see more in the last \fullref{last}). In a way, one can view
the theory of Chern--Simons link invariants as `the Gromov--Witten
theory' of cotangent bundles. Building up on this idea R Gopakumar
and C Vafa completed the 't Hooft's program by finding in 1998 a
dual theory of closed strings that turned out to be a `true'
Gromov--Witten theory but not on $T^*M$.

\subsection{The Gopakumar--Vafa conjecture (1998--2003)}\label{023}

Recall that string theories may change smoothly as the target space
passes through some mild singular transitions. Gopakumar and Vafa
conjectured that the dual theory of closed instantons lives on a
threefold obtained from $T^*M$ via such a transition. They succeeded
in finding the dual threefold for $M=S^3$ \cite{GV}. The
corresponding transition known only to few algebraic geometers from
Clemens' paper \cite{Cl} was very explicitly described in 1990 by P
Candelas and X De La Ossa \cite{Cnf}. It involves the zero section
for $T^*S^3$ shrinking to a nodal point and then getting resolved
into an exceptional $\CP^1$. Their terminology of deformed and
resolved conifolds and the schematic picture of the conifold
transition have since migrated from one paper to the next. They also
showed that all three actors in the conifold transition admit
Calabi--Yau metrics and computed them explicitly. The manifold on the
resolved side of the transition is just the sum
$\calO(-1)\oplus\calO(-1)$ of two tautological bundles over $\CP^1$.
It is now called the resolved conifold even by mathematicians.

Intuitively, as the zero section in $T^*S^3$ collapses into the
nodal point the `open curves' that end on it close up and stay
closed after the resolution of singularity. Thus, Gopakumar and Vafa
conjectured that Witten's open string theory on $T^*S^3$ turns into
the usual Gromov--Witten theory of closed holomorphic curves on
$\calO(-1)\oplus\calO(-1)$, in short:
$$\bfig
  \morphism<1600,0>[\text{CS theory on } S^3`
  \text{GW theory on } \calO(-1)\oplus\calO(-1);
  \text{large $N$ Duality}]
  \efig
$$
This was the original meaning of the Gopakumar--Vafa conjecture on
the level of partition functions.

Physicists make many claims like this; one reason this particular
conjecture generated so much excitement in the mathematical
community is that mathematical machinery was just mature enough to
define both sides of the duality (if not the connection between
them) in rigorous terms.

Another reason is that Gopakumar and Vafa did not stop at a general
physical claim but made two important and completely mathematical
predictions. First, based on Witten's Chern--Simons computations they
predicted an explicit form of the Gromov--Witten free energy function
(and thus all closed Gromov--Witten invariants) for the resolved
conifold. This can be compared to the prediction made by another
physical duality, Mirror Symmetry, for the Gromov--Witten invariants
of projective quintics (see Cox--Katz \cite{CK} and Hori et al
\cite{Hori}). Second, although
Gromov--Witten invariants themselves are rational numbers (see
\fullref{locstr}) they can be represented as combinations of
different numbers, later named the Gopakumar--Vafa invariants, that
are integers. This integrality prediction later replaced the
original meaning of `the Gopakumar--Vafa conjecture' among algebraic
geometers. The environment was so ripe that the prediction of the
Gromov--Witten free energy was verified the same year by C Faber and
R Pandharipande \cite{FP} by an explicit localization computation.
The integrality conjecture for the resolved conifold was proved by
A Okounkov and R Pandharipande in 2003 \cite{OP}.

On the Chern--Simons side of the duality significant progress was made
in computational techniques. In 1999 S Garoufalidis, D Bar-Natan,
L Rozansky and D Thurston discovered an effective algorithm for
computing the LMO invariant \cite{bar,wheels}. It was called the
{\AA}rhus integral (in honor of the city in Denmark where they started
their work in 1995) and uses the graphical calculus of Bar-Natan and
formal Gaussian integration. The same year R Lawrence and L Rozansky
gave a representation for the SU($2$) LMO invariant in terms of
integrals and residues distinguishing contributions from different
flat connections. This type of representation led to the discovery by
Mari\~no \cite{Smatrix} of a relation between the LMO invariants and
matrix integrals.

The next year D Bar-Natan and R Lawrence used surgery and the
{\AA}rhus integral to give an explicit formula for the LMO invariant
of Seifert fibered homology spheres \cite{BL}. Contributions from
nontrivial flat connections for Seifert fibered spaces were
determined in 2002 by M Mari\~no using physical arguments
\cite{Smatrix}. Finally, in 2002 combinatorial expressions for the
Reshetikhin--Turaev version of Witten's invariants associated to an
arbitrary compact Lie group were derived for Seifert fibered spaces
by S Hansen and T Takata \cite{HT}.

The predictions in Gopakumar--Vafa \cite{GV} did not involve mathematical
expectations of Chern--Simons observables, that is, quantum link
invariants. This was rectified in a 2000 paper by H Ooguri and C
Vafa \cite{OV}, who conjectured that to each framed knot in $S^3$
there corresponds a Lagrangian submanifold in the resolved conifold
and polynomial invariants of the knot give the generating function
of open Gromov--Witten invariants. The latter `count' open
holomorphic curves with boundaries on the Lagrangian submanifold.
For the case of the unknot Ooguri and Vafa gave an explicit
construction of the corresponding Lagrangian submanifold in
$\calO(-1)\oplus\calO(-1)$, and predicted the generating function of
open invariants and their integral structure analogous to the
Gopakumar--Vafa invariants in the closed case. A number of explicit
computations for nontrivial knots followed on the Chern--Simons side
(see Labastida--Mari\~no--Vafa \cite{LMV} and references therein).

Unfortunately, the situation on the Gromov--Witten side did not
develop as successfully. For one thing, unlike in the closed case
the definition of the open Gromov--Witten invariants was lacking.
Also, it was not clear how to generalize the Ooguri--Vafa
construction of the Lagrangian to nontrivial knots. Nevertheless,
in 2001 two groups of researchers succeeded in verifying the unknot
predictions of Ooguri and Vafa. To compute the as-of-yet undefined
invariants S Katz and C-C\,M Liu used the virtual localization
technique  as applied by Faber and Pandharipande in the closed case.
T Li and Y\,S Song on the other hand avoided the use of open
invariants altogether \cite{li-song} by replacing them with relative
Gromov--Witten invariants, the theory of which was developed by J Li
at the same time \cite{litian}. A year later M Liu gave a rigorous
definition of open invariants for the case when the Lagrangian
submanifold is invariant under a torus action \cite{liu}.

On the other front J\,M\,F Labastida, M Mari\~no and C Vafa
generalized the Ooguri--Vafa construction to all algebraic knots in
2000 \cite{LMV} and in 2002 C\,H Taubes extended it even to all
symmetric knots \cite{Tb}. S Koshkin later gave a different construction
for a Lagrangian associated to a knot that is valid for any knot
\cite{koshkin}. By the time of the Banff workshop in 2003 the time
seemed right for a general theory of open invariants to emerge and
predictions of the Large $N$ Duality to be verified. However, this
did not happen.

The difficulties proved to be much more significant than in the
closed case: Lagrangian submanifolds for nontrivial knots do not
admit convenient torus actions so the standard computational
techniques do not apply; moreover, the available definition of
invariants does not work in those cases. At the same time (2003) M
Aganagic, A Klemm, M Mari\~no and C Vafa discovered that all
closed Gromov--Witten invariants of toric Calabi--Yau threefolds could
be computed by `slicing' them along Lagrangians corresponding to
framed unknots. Their algorithm known as `the topological vertex'
\cite{vertex} captured the attention of mathematicians and shifted
the focus away from general knots. Shortly after, the topological
vertex was restated in a rigorous form by J Li, C-C\,M Liu, K Liu
and J Zhou \cite{Mvertex} using relative invariants and thus
eliminating a need for open ones altogether as far as toric
varieties are concerned.

Thus 2003 closes the Gopakumar--Vafa period of research and we
conclude our historical excursion. Progress made after 2003 went in
a different direction and we shall give some indications about this
in the last section of this paper.

\section{Overview}\label{03}

In this section we introduce some minimal notation to explain in a
nutshell what the rest of this fat paper is all about. The
Gopakumar--Vafa Large $N$ Duality \cite{GV} is a correspondence
between two theories, one defined on $S^3$ and the other on
$\calO(-1)\oplus\calO(-1)$. Here $\calO(-1)$ is the tautological
line bundle over $\CP^1$ defined by
\index{$O$@$\calO(-1)$ the tautological line bundle}
$$\calO(-1):=\biggl\{ ([z_1:z_2], w_1,w_2)\in \CP^1\times\C^2
   ~\bigg|~ \biggl|\begin{matrix}w_1&w_2\\
     z_1&z_2 \end{matrix}\biggr| =0 \biggr\}.$$
The projective line $\CP^1$ is the collection of complex lines
through the origin in $\C^2$ and the $\calO(-1)$ bundle is simply
the collection of pairs consisting of a line together with a point
on that line (hence the name `tautological line bundle'). The
projection just maps each point to the corresponding line. The
number $-1$ refers to the fact that the first Chern class of this
bundle evaluates to $-1$ on the $\CP^1$ cycle.

Intuitively, the large $N$ correspondence can be traced to a
geometric relation between the two spaces called the conifold
transition $T^*S^3\leadsto\calO(-1)\oplus\calO(-1)$. To understand
this transition consider $T^*S^3\simeq TS^3$ realized as
$\text{SL}_2\C$ in $\C^4\simeq$End($\C^2$):
\begin{equation}\label{SL2C}
T^*S^3\simeq\biggl\{
  W=\biggl(\begin{matrix}w_1&w_2\\w_3&w_4\end{matrix}\biggr)
  \in\text{End}(\C^2)
   ~\bigg|~ \det(W)=1 \biggr\}=\text{SL}_2\C.
\end{equation}

\begin{figure}[ht!]
\centering
\labellist\small
\pinlabel {Deformation} [b] at 245 148
\pinlabel {Resolution} [b] at 495 148
\endlabellist
\includegraphics[width=4truein]{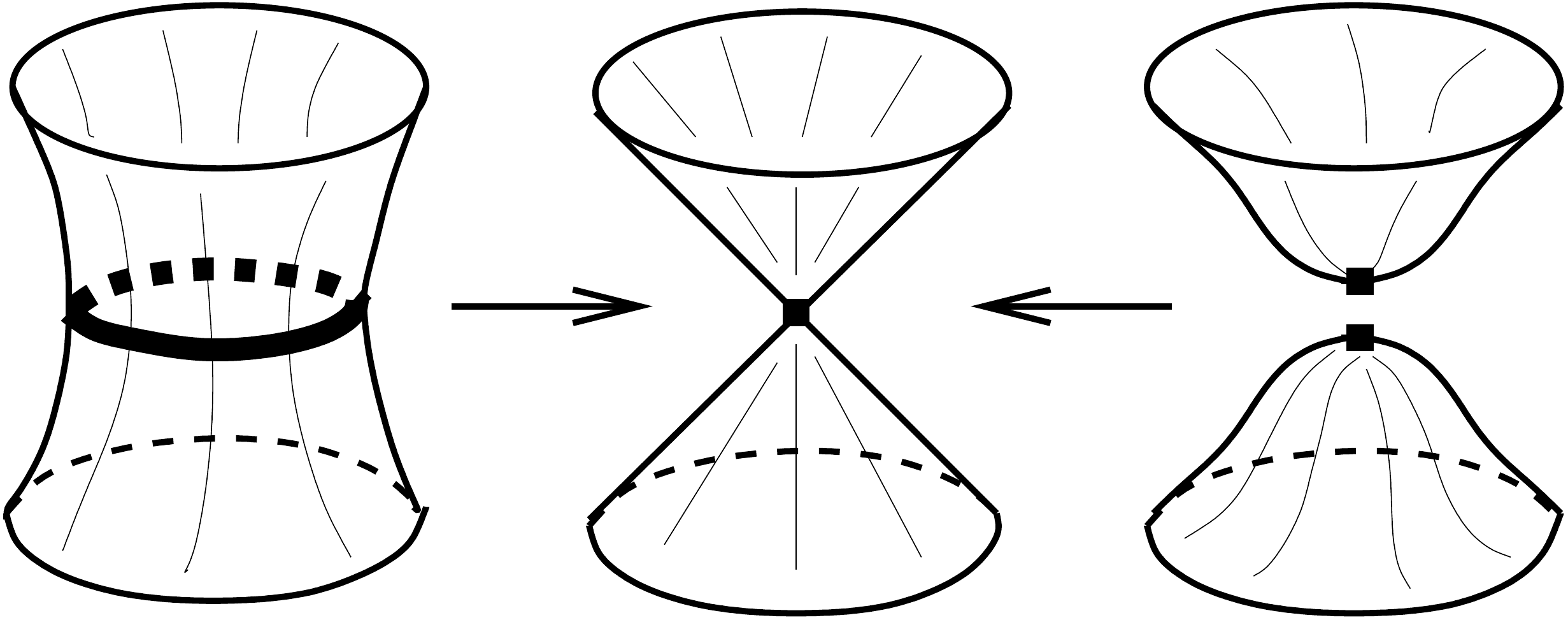}
\caption{The `conifold' transition two dimensions down
$S^1{\times}\R^1{\leadsto}S^0{\times}\R^2$}\label{F111}
\end{figure}

\begin{exm} Consider the standard embedding $S^3\hookrightarrow\R^4$ as the unit sphere
inducing the embedding $TS^3\hookrightarrow\R^4\times\R^4\simeq\C^4$
and find an explicit automorphism $\C^4\to\C^4$ that restricts to a
diffeomorphism $TS^3\to\text{SL}_2\C$ and $S^3\to\text{SU}(2)$ (see
Koshkin \cite{koshkin}).
\end{exm}
Now set $\det(W)=\mu$ in equation \eqref{SL2C}. Taking $\mu\to0$
produces a complex deformation of $T^*S^3$ into a singular variety
that physicists and mathematicians call the conifold (it has an
ordinary double point at the origin and is a higher dimensional
analog of the usual double cone as in \fullref{F111}):
$$\wcheck{X}_{S^3}:=\biggl\{
  \begin{pmatrix}w_1&w_2\\w_3&w_4\end{pmatrix}\in\text{End}(\C^2)
   ~\bigg|~ \det(W)=0 \biggr\}.$$
The conifold admits a small resolution
$X_{S^3}\xrightarrow{\rho}\wcheck{X}_{S^3}$ of the singularity
(see Harris \cite{AG1}) small meaning that the exceptional locus $\rho^{-1}(0)$
is a curve rather than a surface as in the case of a blow-up. This
resolution is the resolved conifold.
\begin{defn} The resolved conifold is
$$X_{S^3}=\biggl\{
   ([z_1:z_2],w_1,w_2,w_3,w_4)\in\CP^1\times\C^4
   ~\bigg|~ \begin{vmatrix}w_1&w_2\\z_1&z_2\end{vmatrix}
    =\begin{vmatrix}w_3&w_4\\z_1&z_2\end{vmatrix} =0 \biggr\}.$$
\end{defn}\index{$X_{S^3}$ the resolved conifold}\index{resolved conifold}
The resolved conifold $X_{S^3}$ is easily seen to be biholomorphic
to $\calO(-1)\oplus\calO(-1)$.

Geometrically the deformation shrinks the zero section (ie $S^3$)
into the double point and the resolution replaces it with an
exceptional $\CP^1$ so topologically we have the surgery
$S^3\times\R^3\leadsto \R^4\times S^2$. A low-dimensional analog of
the conifold transition is shown in \fullref{F111}.

The resolved conifold $X_{S^3}$ is a K\"ahler manifold and one can
show that the expected dimension of the space of stable maps
(holomorphic curves) from a Riemann surface into $X_{S^3}$ is zero.
Thus such maps are (formally) isolated and it makes sense to count
them. The Gromov--Witten invariant $N_{g,d}(X_{S^3})$ is intuitively
the number of maps of genus $g$ surfaces representing the homology
class $d[\CP^1]$. The actual numbers can be fractional because one
must assign fractional weights to curves with automorphisms.
Invariants defined in this way can be conveniently collected into
the full Gromov--Witten free energy (see \fullref{GW} for
details):
$$F^{\text{GW}}_{X_{S^3}}(t,y):=\sum_{g=0}^\infty\sum_{d=0}^\infty
N_{g,d}(X_{S^3})e^{-td}y^{2g-2}.$$
On the other side, the unnormalized Chern--Simons free energy is the
logarithm of its partition function $F^{\text{CS}}_{S^3}:=\ln
Z^{\text{CS}}(S^3)$. Again intuitively, the latter is the value of
the path integral:
$$Z^{\text{CS}}(S^3)=\int_{{\mathcal A}} e^{iCS(A)}{\mathcal D}A,$$
taken over the space of connections on a trivial $\text{SU}(N)$
bundle over $S^3$. In this formula
$$CS(A):=\frac{k}{4\pi}\int_M \text{Tr}(A\wedge dA
  +\tfrac23 A\wedge A\wedge A)$$
is the $SU(N)$ Chern--Simons action, where $k$ is an arbitrary
integer called the level. Thus the unnormalized Chern--Simons free
energy is a function of two parameters $k,N$ encoded as $N$ and
$x:=\frac{2\pi}{k+N}$ just as the full Gromov--Witten free energy is
a function of $t,y$.

Today there is a multitude of rigorous constructions that define
invariants $Z^{\text{CS}}(M)$ for any 3--manifold $M$ that have the
properties that one would conjecture based on heuristic path
integral manipulations. Most of them use the fact due to Lickorish
and Kirby that any $3$--manifold $M$ can be presented as a surgery
on a framed link $L_M$ and two different links present the same
3--manifold if and only if they are related by a sequence of the
so-called Kirby moves (see Prasolov--Sossinsky \cite{prasolov} and
Turaev \cite{T}). Therefore, if one can come up with a framed link
invariant that does not change under the Kirby moves one gets a
3--manifold invariant. N Reshetikhin and V Turaev were the first ones
to come up with a systematic procedure for constructing Kirby-move
invariant link invariants. Their invariants were based on the theory
of quantum groups \cite{RT1,RT2}. We use a version of their
construction in Sections \ref{modinv} and \ref{RTinv}. The
corresponding invariant which is just the THOMFLYP polynomial at
roots of unity can also be constructed in a number of other ways:
skeins, TQFT, etc that all lead to the same quantity identified with
the above partition function. We briefly touch on the skein and TQFT
approaches in Appendices \ref{app:c} and \ref{app:e} respectively.

We now give an intuitive idea how the Gopakumar--Vafa duality arises
from the physics of string theory (we are grateful to M Mari\~no
for explaining this to us). As mentioned in the history section
physicists work with a very broad notion of `holomorphic instantons'
described by a topological version of string theory known as the
`topological A--model'. Holomorphic instantons live in Calabi--Yau
threefolds and can be closed or open, that is, have boundary. In the
latter case their boundary lies on `D-branes', located at (`wrapped
around') special Lagrangian submanifolds of the threefold. Open and
closed holomorphic curves and stable maps are examples of
holomorphic instantons but there are more degenerate ones as well,
for example `instantons at infinity' (see Witten \cite{Wcss}) -- trivalent ribbon graphs
in Lagrangian submanifolds representing infinitely thin `curves with
boundary'. Physical quantities produced by the theory are called
`string amplitudes' and in good cases they can be identified with
Gromov--Witten invariants of the threefolds. The Calabi--Yau condition
is needed to make sure that holomorphic instantons are isolated and
can be `counted' with finite amplitudes.

Now consider two extreme cases of this picture. The first case is
when there are no `honest' holomorphic curves as in cotangent
bundles $T^*M$ to 3--manifolds since the symplectic form on them is
exact. At the same time, the zero section is a special Lagrangian
submanifold and knotted trivalent graphs in $M$ (framed knots and
links in the simplest case) can be seen as degenerate instantons
with boundary. Witten discovered in \cite{Wcss} that the string
amplitudes of a cotangent bundle with `$N$ D-branes wrapped around
the zero section' can be recovered from the quantum $\text{SU}(N)$
Chern--Simons invariants of $M$ computable via the surgery
prescription from his earlier paper \cite{W1}. Moreover, one can
also consider instantons ending on conormal bundles to links in $M$
that are also special Lagrangian in $T^*M$. This time the string
amplitudes coincide with the the quantum $\text{SU}(N)$ link
invariants. In other words, the topological A--model reduces to the
quantum Chern--Simons theory in this case and can be viewed as `the
Gromov--Witten theory' of cotangent bundles. K Fukaya gives this
idea a more precise meaning in terms of Floer homology in
\cite{Fuk}.

In the second case there are no D-branes in the picture and the only
holomorphic instantons that remain are closed stable maps. This is
the case of the resolved conifold and its usual Gromov--Witten
theory. A striking feature of string theory discovered by Dixon, Harvey,
Vafa and Witten \cite{DHVW} is that physically equivalent models can be set in
different `geometric backgrounds', that is, live on different
threefolds. This occurs when the geometric backgrounds are related
by special geometric transitions. Whereas physical quantities do not
change, the underlying threefold may undergo a singular transition
as some of the D-brains and/or holomorphic homology classes collapse
or appear. The conifold transition is the simplest example of such a
geometric transition.

Open instantons that end on the zero section close up as they shrink
to the nodal point and transform into closed holomorphic curves in
the resolved conifold. Unlike the zero section conormal bundles to
links do not collapse and reappear as Lagrangian submanifolds in the
resolved conifold. Instantons at infinity that ended on them
therefore transform into open holomorphic curves. String amplitudes
computed on both sides of the transition should be the same since
the physics does not change. In a nutshell, this is the insight
behind the Gopakumar, Ooguri and Vafa predictions of equality
between the Chern--Simons 3--manifold and link invariants on one side
and closed and open Gromov--Witten invariants on the other
(see Gopakumar--Vafa \cite{GV} and Ooguri--Vafa \cite{OV}).

The name Large $N$ Duality comes from the specific way string
amplitudes for instantons at infinity are recovered from the
Chern--Simons invariants. One needs to consider the latter not for a
specific rank but for {\it all} ranks $N$. The resulting function
turns out to be analytic in $\frac{1}{N}$ around $0$ modulo some
logarithmic terms and can be expanded into a Laurent series. The
coefficients of this series are the string amplitudes in question.
From the perspective of knot theory this means that string
amplitudes are given not by the `exact' invariants (such as Jones or
THOMFLYP polynomials) but by the so-called Vassiliev (or finite-type
or perturbative) invariants \cite{bn2,wheels}. One should not be
misled by the name into believing that the duality holds at large
$N$ only. The Laurent coefficients match the Gromov--Witten
invariants of the resolved conifold at all powers of $\frac1N$ and
the duality is exact.

One of the most attractive traits of the Large $N$ Duality is its
computational power. It is usually much easier to compute gauge
theory quantities (partition functions, correlators, etc) than the
corresponding string amplitudes. This is due to the apparatus of
informal but effective path integral manipulations successfully
applied by physicists for quite some time. String theory techniques
(such as equivariant localization) are more recent and much more
cumbersome. It turns out that contributions from instantons at
infinity can be reduced to path integrals and even those
contributions of honest holomorphic curves can be represented by
path integrals with extra insertions. Thus the computational
machinery of the gauge theory becomes available for string theories
as well. One remarkable achievement of this approach is the
topological vertex algorithm originally derived from the
Chern--Simons path integral.  This algorithm computes Gromov--Witten
invariants of all toric Calabi--Yau threefolds (see
Aganagic--Klemm--Mari\~no--Vafa \cite{vertex} and Mari\~no
\cite{Mar}).

At present we are very far removed from a mathematical definition of
the topological A--model in anywhere near the generality used by
physicists. However, the prediction of the equality of the
Gromov--Witten and Chern--Simons free energies or partition functions
on $X_{S^3}$ and $S^3$ is a well-defined mathematical statement.
This equality should not be taken too literally, one has to
renormalize and change variables to make it work. But it is true
that one function can be recovered from the other as Gopakumar and
Vafa convincingly demonstrated by correctly predicting the values of
the Gromov--Witten invariants of $X_{S^3}$ in \cite{GV}.

The computation that verifies the Gopakumar--Vafa predictions was
originally done by C Faber and R Pandharipande in \cite{FP} but
it does not cover the Chern--Simons side relying on the formulas
obtained by path integral methods. Later papers \cite{GR,Mar} that
compute and compare both free energies skip many of the details to
stay within a limited length. In this paper we provide the
background material on both theories necessary to understand the
structures behind these computations and reproduce the computations
themselves (\fullref{comparison}). The result of comparison can be
packaged in the following form:
\begin{thm}\label{sphereGV}
The full Gromov--Witten free energy and the unnormalized Chern--Simons
free energy are related by
$$\text{Re}\bigl(F^{\text{GW}}_{X_{S^3}}(iNx,x)-\,
  F^{\text{CS}}_{S^3}(N,x)\bigr)=\tfrac{5}{12}\ln
  x+\zeta(3)x^{-2}-\tfrac12\ln(2\pi)-\zeta^\prime(-1).$$
\end{thm}
Some comments are in order about the form of this formula. First of
all, even though the free energies are complex-valued the relevant
coefficients in their expansions, that is, the invariants themselves are
real so it suffices to consider only the real parts. Secondly, the
free energies fail to be holomorphic in $x$ at zero where the
expansions are taken, but they do so in a very minor way. The terms
on the right appear as a result of regularization and do not
indicate any meaningful discrepancy.

To appreciate how powerful this theorem is note that the full
Gromov--Witten free energy encapsulates the Gromov--Witten invariants
of the resolved conifold in {\it all degrees} and {\it all genera}.
In its turn, the Chern--Simons partition function for the Hopf link
contains the $\Sl_n\C$ Vassiliev invariants of {\it all knots and
links} in the three-sphere. It should come as no surprise that the
duality for `just' this one example led to computation of the
Gromov--Witten invariants for all local Calabi--Yau threefolds
\cite{vertex,Mar}. Despite its somewhat unappealing form this
formula is a very strong confirmation of the Gopakumar--Vafa
conjecture. We finish this introduction by stating a far-reaching
generalization of \fullref{sphereGV} suggested by M Mari\~no in
his Banff lectures.
\begin{conjecture}
For every rational homology $3$--sphere $M$ there exists a large
$N$ dual Calabi--Yau threefold $X_M$ such that the Chern--Simons
theory on $M$ is equivalent to the Gromov--Witten theory on $X_M$. In
particular, the corresponding invariants can be recovered from each
other.
\end{conjecture}
For a reader who may think that the task of learning so much
algebraic geometry and quantum algebra is overwhelming we promise
that learning these complex but remarkable structures is well worth
the effort. While navigating the deep waters of abstraction the
reader should always keep in mind that we are merely
computing two complex-valued functions -- the free energies of
$X_{S^3}$ and $S^3$.

\newpage

\part{Gromov--Witten invariants}
\setobjecttype{Part}
\label{GW}

The theory of Gromov--Witten invariants is the mathematical theory
closest to string theory in physics. These invariants arise as
generalizations of enumerative invariants. In this part, we will
outline the definition of Gromov--Witten invariants and give some
sample computations.

The first ingredient in understanding these invariants is the cohomological
interpretation of intersection theory. As a simple example consider counting
the number of zeros of a degree $d$ polynomial in $\C[x]$, say
$p(x)$. The answer is easier when we are using
complex coefficients. In more general counting problems the answers
will be more uniform if we work in complex projective spaces.
Given a degree $d$ polynomial $p$ we can define a function,
$f_p\co\CP^1\to\CP^1$ given by $f_p([z:w])=[p(z/w)w^d:w^d]$. This
induces a map on the second cohomology, $f^*_p\co H^2(\CP^1;\Z)\to
H^2(\CP^1;\Z)$. Since $ H^2(\CP^1;\Z)\cong \Z$, this map is just
multiplication by some integer. This integer is known as the degree
of the map and it coincides with the degree of the original
polynomial. We can write this as
$$\#(p^{-1}(0))=\#(f_p^{-1}([0:1]))=\int_{[\CP^1]}f^*_p\omega_{\CP^1}.$$
Here $\#$ represents a signed count of a set of points in general
position. We will later describe methods to address non-generic
situations. The integral represents the cap product pairing of
homology and cohomology, so $[\CP^1]$ is the fundamental homology
cycle ($\CP^1$ has a natural orientation coming from the complex
structure) and $\omega_{\CP^1}$ is the orientation class. Using the
de Rham model for cohomology, the integral will become an honest
integral.

\fullref{intvco} provides a correspondence between geometric
intersections and cohomological operations. We will describe various
lines in this table as we use them. General topological folk wisdom
suggests thinking via intersections and proving via cohomology.
\begin{table}[ht!]
\centering
\begin{small}
\begin{tabular}{|c|c|} \hline
Intersections & Cohomology\\ \hline
A codimension $k$ homology submanifold, $A$ &
   A cohomology class $\alpha\in H^k(M;\Z)$\\
$\#(A\cap F)$ & $\alpha([F])=\alpha\cap[F]=\int_{[F]}\alpha$ \\
$A\cap B$ & $\alpha\cup\beta=\alpha\wedge\beta$ \\
$f^{-1}(A)$ & $f^*\alpha$ \\
$\sigma_1^{-1}(\sigma_0(X))$ & $c_1(L)$ \, \text{or\ }$e(E)\in H^r(X)$ \\
$\sigma_0(X)$ & \text{Thom class \ }$\Phi\in H^r(E,E-\sigma_0(X))$ \\ \hline
\end{tabular}
\vskip.1in
\end{small}
\label{intvco}
\caption{Geometric intersections vs cohomology}
\end{table}

We now turn to a more serious motivating question: how many lines
pass through two generic points in a plane? More generally, how
many degree $d$ parameterized curves pass through the `right' number
of points in a plane modulo reparameterization of the domain? We may
describe the space of lines with two marked points in the plane as
\begin{multline*}
M_{0,2}(\CP^2,1[\CP^1]) = \\[-1ex]
\{u\co\CP^1\to \CP^2, p_1, p_2 ~|~ u_*[\CP^1]=1[\CP^1], \wbar\partial u=0,
p_1\ne p_2\in\CP^1\}/\sim.
\end{multline*}
Maps in this set are explicitly given by
$u([z:w])=[az+bw:cz+dw:ez+fw]$. Now count the dimension of this
space. There are $6$ complex parameters in the definition of our
degree one parameterized curve. However the points in the projective
plane are only defined up to a scale, so we subtract one parameter.
We wish to count two maps as equivalent if they are related by a
reparameterization of the domain (this is what the $\sim$ represents
in the equation for $ M_{0,2}(\CP^2,1[\CP^1])$). The holomorphic
isomorphisms (reparameterizations) of $\CP^1$ are just the linear fractional
transformations. More explicitly, we
define $(u, p_1, p_2)\sim(v, q_1, q_2)$ to hold if and only if there
is a linear fractional transformation, say $\varphi$, such that
$u=v\circ\varphi$ and $\varphi(p_k)=q_k$. A similar count allows one to conclude that the
space of linear fractional  transformations has complex dimension
$3$. It follows that the space of complex projective lines in the
complex projective plane has complex dimension two. Adding two
points in the domain adds two more complex parameters, so the
complex dimension of $ M_{0,2}(\CP^2,1[\CP^1])$ is four.

Now we have two natural evaluation maps taking $M_{0,2}(\CP^2,1[\C
P^1])$ to $\CP^2$, given by $\text{ev}_k([u, p_1, p_2]):= u(p_k)$.
Since $\CP^2$ is two complex-dimensional the following integral
makes sense and represents the number of lines passing through two
points:
$$\int_{[\wwbar{M}_{0,2}(\CP^2,1[\CP^1])]} \text{ev}_1^*\omega_{\C
P^2} \wedge \text{ev}_2^*\omega_{\CP^2}$$
There is an infinite number of lines passing through one fixed
point, and there are no lines passing through three generic points
in the plane, thus two is the `right' number of points to mark when
counting lines.

\begin{exm}
Count the dimension of the space of degree $d$ genus zero
parameterized curves into $\CP^2$ modulo reparameterization of the
domain. Write an expression similar to the integral above
representing the number of such curves through the `right' number of
points. We will outline two different ways to compute these numbers
later in \fullref{GW}.
\end{exm}

In the next section we will define the Gromov--Witten invariant
$\langle \gamma_1,\ldots, \gamma_n\rangle_{g,\beta}^X$ The intuitive
interpretation of $\langle \gamma_1,\ldots,
\gamma_n\rangle_{g,\beta}^X$ is the number of genus $g$ curves in
the class $\beta$ that intersect the cycles
$\Gamma_1,\ldots,\Gamma_n$ Poincar\'e dual to $\gamma_1,\ldots,
\gamma_n$. In general, this interpretation fails. In fact, $\langle
\gamma_1,\ldots, \gamma_n\rangle_{g,\beta}^X$ are only rational
numbers, not integers. This is because curves must be counted with
fracitonal weights to get a correct definition.

\section{The coarse moduli space}\label{coarse}
Gromov--Witten invariants extend the ideas described in the above
problem. These invariants may be defined for symplectic manifolds or
for projective algebraic varieties. We mostly use the symplectic
definition in this section, but give some idea of the
algebraic definition at the end.

\subsection{The symplectic construction}

Recall that a \index{symplectic manifold}\index{symplectic form}
\index{$\omega$ a symplectic form} symplectic manifold is a (real)
$2n$--dimensional manifold with a $2$--form $\omega$ such that
$d\omega=0$ and $(n!)^{-1}\omega^{\wedge n}$ is a volume form on
$X$. Any symplectic manifold admits a compatible almost complex
structure. An almost complex structure $J\in\Gamma(\text{End}(TX))$ \index{almost complex
structure} \index{$J$ an almost complex structure} is an
endomorphism of the tangent bundle which squares to negative one $J^2=-I$. Such is compatible with a
symplectic form, say $\omega$, if the tensor defined by
$g(X,Y)=\omega(X,JY)$ is a Riemannian metric. Any two of $g$,
$\omega$ or $J$ uniquely determine the third via the compatibility
condition. The standard symplectic structure on $\C$ is given by
$\omega= \frac{i}{2}dz\wedge d\wbar{z}$ and the standard
symplectic structure on $\CP^1$ is given by
$$\omega=\frac{i(dz\wedge d\wbar{z}+ dw\wedge d\wbar{w})}{2(|z|^2+|w|^2)^2}.
$$
Together with the standard complex structure this produces a round
metric of radius of $\frac12$ on the Riemann sphere $\CP^1$. There
are similar symplectic structures on all complex projective spaces.
The standard complex structure on $\C$ considered as a real vector
space is just multiplication by $i$. In the basis $\{1,i\}$ it is
given by the matrix
$$J=\begin{pmatrix}0&-1\\1&0  \end{pmatrix}.$$
In tensor notation it is given by
$$J=\partial_y\otimes dx-\partial_x\otimes dy = i\partial_z\otimes dz
-i\partial_{\wbar{z}}\otimes d{\wbar{z}}.$$
In a local chart the complex structure on any complex manifold takes
this form (generalized in the obvious way to $\C^n$). The product of
two symplectic manifolds, $(X_k,\omega_k)$, $k=1,2$, is the
symplectic manifold $(X_1\times X_2,p_1^*\omega_1+p_2^*\omega_2)$
where $p_k$ are the natural projections. Thus $\CP^1\times\C^4$
inherits a natural symplectic structure as a product. This leads to
a symplectic structure on our main example, the {\it resolved
conifold}
$$X_{S^3}:=\biggl\{[z_1,z_2],(w_1,w_2,w_3,w_4)\in \CP^1\times\C^4 ~\bigg|~
\biggl|\begin{matrix}w_1 &w_2\\z_1 &z_2\end{matrix}\biggr|=
\biggl|\begin{matrix}w_3 &w_4\\z_1 &z_2\end{matrix}\biggr|=0\biggr\}.$$
This is our main example because it will turn out to be the large
$N$ dual of $S^3$. We will describe the corresponding moduli space
in \fullref{vc}, compute the Gromov--Witten invariants in
\fullref{fmc} and see that this is the large $N$ dual of
$S^3$ in \fullref{comparison}.

The notion of a holomorphic curve can be generalized from algebraic
varieties to symplectic manifolds. A symplectic manifold with a
compatible almost complex structure provides exactly the data needed
for the target of a pseudoholomorphic curve. A Riemann surface
$\Sigma$ has an associated almost complex structure which is
automatically a complex structure denoted by $j$. The
Cauchy--Riemann operator \index{Cauchy--Riemann operator}
\index{$\dbar$ the Cauchy--Riemann operator} is defined on a map
$u\co\Sigma\to X$ by
$$\dbar u=\tfrac12\left(du+J_u\circ du\circ j\right).$$
By definition, a map $u$ is pseudoholomorphic
\index{pseudoholomorphic} when $\dbar u \equiv 0$.
\begin{exm}
Write out the Cauchy--Riemann operator on maps $f\co\C\to\C$ using
$x$--$y$ coordinates on the first factor and $u$--$v$ coordinates on
the second assuming the natural symplectic and almost complex
structures.
\end{exm}

\begin{defn}
A smooth, genus $g$, $n$--marked, pseudoholomorphic curve in $X$ is a
tuple $(\Sigma,j,u,p_1,\ldots,p_n)$ where $\Sigma$ is an oriented
genus $g$ surface, $j$ is an almost complex structure on $\Sigma$,
$u$ is a pseudoholomorphic map $u\co\Sigma\to X$ and $p_k\in\Sigma$
are distinct. A morphism between $(\Sigma,j,u,p_1,\ldots,p_n)$ and
$(\Sigma^\prime,j^\prime,u^\prime,p_1^\prime,\ldots,p_n^\prime)$ is
a holomorphic map, $\varphi\co\Sigma\to\Sigma^\prime$, such that
$u^\prime\circ\varphi=u$ and $\varphi(p_k)=p_k^\prime$. A genus $g$,
$n$--marked, pseudoholomorphic curve in $X$ will be called stable if
it has a finite automorphism group.
\end{defn} \index{stable}

We can now define \index{$M$@$M_{g,n}(X,\beta)$ the coarse moduli
space} \index{coarse moduli space} the (coarse) moduli space of
genus $g$, $n$--marked, stable pseudoholomorphic curves in a homology
class $\beta\in H_2(X;\Z)$ to be the set of equivalence classes of
such,
$$
M_{g,n}(X,\beta)\,=\, \{[\Sigma,j,u,p_1,\ldots,p_n] |
u_*[\Sigma]=\beta\}/\sim.
$$
There are natural evaluation maps, $\text{ev}_k\co M_{g,n}(X,\beta)\to
X$ given by \index{$E$@$\text{ev}_k$ evaluation at point $k$}
$$
\text{ev}_k([\Sigma,j,u,p_1,\ldots,p_n])=u(p_k).
$$
Using these, we have our first definition of the Gromov--Witten
invariants. Let $\gamma_1,\ldots, \gamma_n\in H^*(X;\Q)$, $\beta\in
H_2(X;\Z)$  and define the Gromov--Witten invariants by
\index{$\langle  \rangle_{g,\beta}^X$ the Gromov--Witten
invariants}\index{Gromov--Witten invariants}
$$
\langle \gamma_1,\ldots, \gamma_n \rangle_{g,\beta}^X:=\int_{[\wwbar
{\calM}_{g,n}(X,\beta)]^{\vir}}
\text{ev}_1^*\gamma_1\wedge\ldots\wedge\text{ev}_n^*\gamma_n.
$$
You will notice that there are many notations in this definition
that we have not defined yet. We will slowly compute the
Gromov--Witten invariants of $\CP^2$ (hence how many cubic
parameterized curves pass through $8$ points etc) and the
Gromov--Witten invariants of $X_{S^3}$. Along the way, we will
define the extra notations used in the above formula. For now when
you see ${[\wwbar{\calM}_{g,n}(X,\beta)]^{\vir}}$ you should just
think $M_{g,n}(X,\beta)$. We will see that the overline refers to a
compactification of this space later in this article. The
calligraphic font refers to the stack structure on the moduli space.
Intuitively the stack structure `adds' a group to each point of the
space in order to have a proper count of points taking symmetry into
account. The necessity of stacks is motivated in the first article
in \fullref{locstr},  and the definition of the moduli stack
$\wwbar{\calM}_{g,n}(X,\beta)$ is given in the second article of
this subsection. The general definition of a stack is given in
\fullref{app:a}, and the best example is contained in \fullref{mc}.
The square brackets with vir superscript indicate the virtual
fundamental class. This is motivated and defined in \fullref{vc}.
The short explanation is that when intersections are counted one
generally assumes that objects are in general position. However, one
can still get sensible answers when the objects are not in general
position provided one uses the correct virtual fundamental classes.

\begin{figure}[ht!]
\centering
\includegraphics[width=3truein]{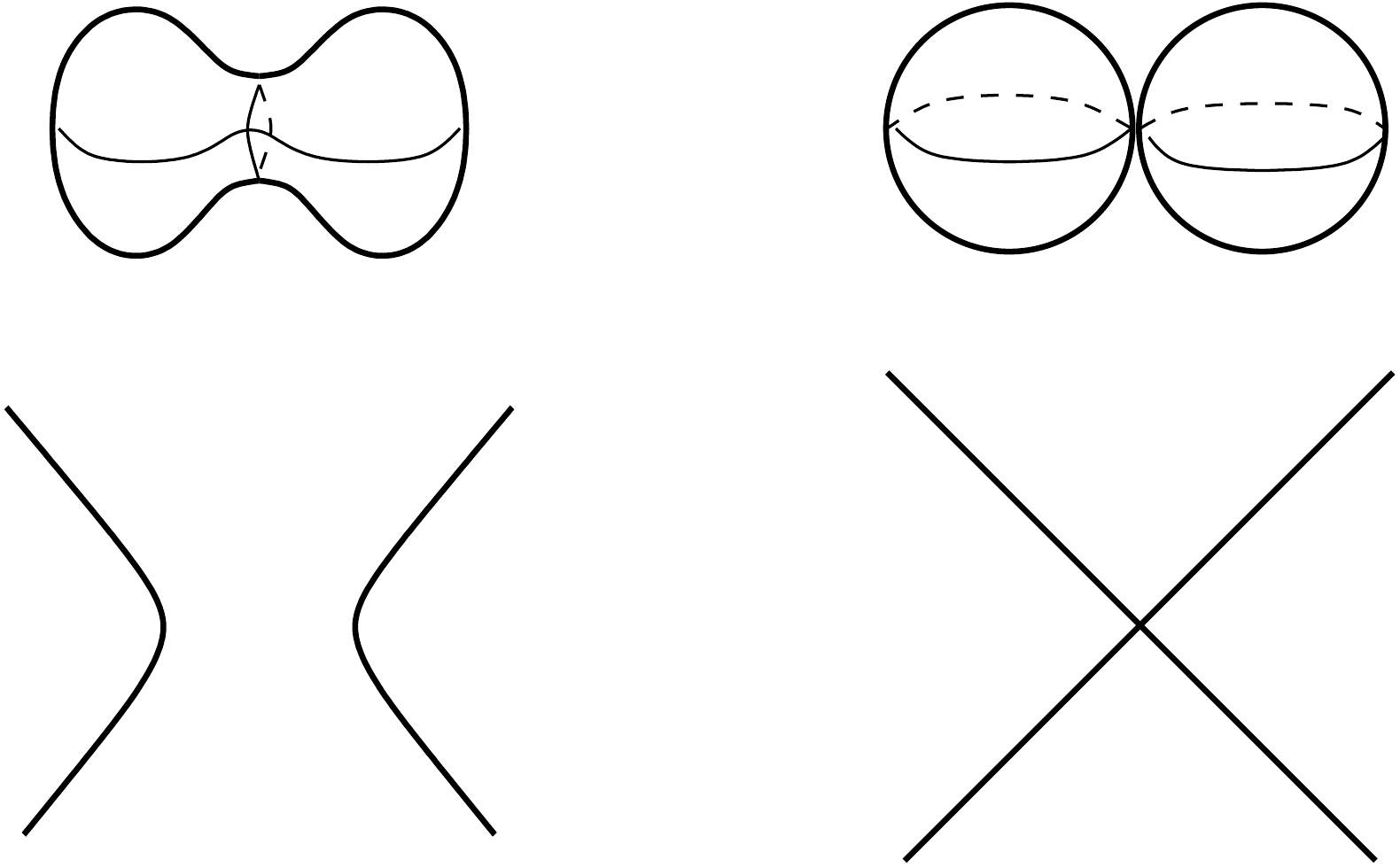}
\caption{Bubbling}\label{bubble}
\end{figure}

To make sense of the integral in the definition we need to have a
fundamental cycle to integrate over. This is easiest to establish
when the domain of integration is compact. As defined above the
coarse moduli spaces would not be compact because we insist that the
marked points be distinct. Even without considering marked points
these spaces will fail to be compact. To see the problem, consider
the family of degree two parameterized curves $u_n\co\CP^1\to\CP^2$
given by $u_n([s:t]):=[n^{-1}(s^2+t^2):2n^{-1}st:s^2-t^2]$. The
limit appears to be given by $u_\infty([s:t])=[0:0:s^2-t^2]$, but
this is not a well defined map (consider points where $s=\pm t$).
The image of the map $u_n$ in an affine chart is just the hyperbola,
$x^2-y^2=n^{-2}$. The geometric limit of this sequence is just a
pair of lines. This is a clue that motivates the correct definition
of a compactification of the moduli space. The correct limit is a
map from the one point union of two copies of the complex projective
line. The domains and (real part of the) images of $u_n$ and the
limit $u_\infty$ are displayed in \fullref{bubble}. To be precise
the limit is defined by
$$u_\infty\co\CP^1\times\{\pm 1\}\to\CP^2;\quad
  u_\infty([s:t],\epsilon):=[s:\epsilon s:t]; \quad \epsilon=\pm 1.$$
The splitting off of an extra component in the limit is called
bubbling.

\begin{figure}[ht!]
\includegraphics[width=4.7truein]{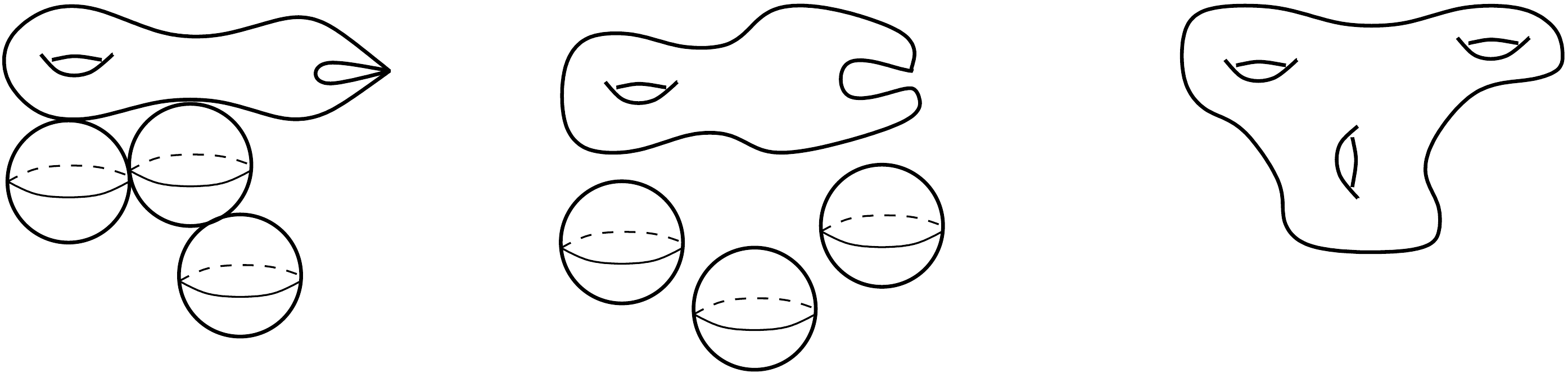}
\caption{Prestable curve, normalization and smoothing}\label{node}
\end{figure}

To describe domains of stable maps we take a disjoint union of
Riemann surfaces (called the normalization) and glue them together
along special points called nodes. The resulting `surface' needs to
be connected. These objects are called prestable curves or nodal
Riemann surfaces. The original surface without any identifications
is called the normalization. \index{normalization} The nodes
\index{nodes} are locally modeled on $\{(x,y)|xy=0\}$. There is also
a notion of a smoothing \index{smoothing} of a prestable curve
\index{prestable curve} obtained by replacing each node by
$\{(x,y)|xy=\epsilon\}$. An example of a marked prestable curve
together with its normalization and smoothing is shown in
\fullref{node}.
Precise definitions are given in the paper by Siebert on
Gromov--Witten invariants for general symplectic manifolds
\cite{Siebert}. In fact, this is a superb reference for the
definition of Gromov--Witten invariants in the symplectic category
for readers who are more familiar with differential topology than
algebraic geometry. While we are on the subject of references, we
should mention the book by Hori et al \cite{Hori}, the book by Cox
and Katz \cite{CK} and the little-known book edited by Aluffi
\cite{aluffi}.  These are the main references from which this
section on Gromov--Witten invariants was derived. For more
information see these sources and the references contained therein.

What follows is a more formal definition of a marked prestable curve.
\begin{defn}[Siebert \cite{Siebert}]
A marked prestable curve is a pair, $(\Sigma,p)$ where $\Sigma$ is a
reduced, compact, connected, one-dimensional, complex, projective
variety with no worse than ordinary double point singularities, and
$p$ is an $n$--tuple of pairwise distinct regular points.
\end{defn}
An ordinary double point is a point modeled on $\{(x,y)|xy=0\}$ (see Griffiths--Harris \cite{GH} and
Hartshorne \cite{hart}). Reduced means that there is no locally defined non-zero
holomorphic function with a power equal to zero. To see how such a function
could arise consider the variety defined by $x^2y=0$. This looks like a one point union of a pair of lines,
but it is not a reduced variety because $xy$ is a non-zero
holomorphic function with square zero. The one point union of a pair
of lines is represented by the reduced variety defined by $xy=0$.

A prestable map is a holomorphic map from a prestable curve to a
symplectic manifold or projective variety. The definition of a
morphism of marked pseudoholomorphic curves, and resulting notions
of equivalence, stability and automorphism group extend naturally to
prestable curves. This allows one to define a compactification of
the moduli space as follows.
\begin{defn}
The (compactified) coarse moduli space \index{c@(compactified)
coarse moduli space}\index{$M1$@$\wwbar{M}_{g,n}(X,\beta)$ the
(compactified) coarse moduli space} of genus $g$, $n$--marked curves
is
\begin{multline*}
\wwbar{M}_{g,n}(X,\beta) = \\
\bigl\{[\Sigma,p,u] ~|~ u_*[\Sigma]=\beta \text{ and } (\Sigma,p)
\text{ is a genus $g$, $n$--marked prestable curve} \\[-1ex]
\text{such that $[\Sigma,p,u]$ is stable}\bigr\}/\sim.
\end{multline*}
\end{defn}
It may seem weird that this space is compact when it is required that
the marked points be disjoint from each other and the nodes.  The limit
of a sequence where two or more marked points or nodes collide will be
described in the discussion of the boundary divisors and $\psi$ classes
in \fullref{GWco} on the cohomology of the moduli space. Basically,
the idea is that one or more new components bubble off at the point of
collision and the marked points move into these bubbles while staying
distinct. Before we talk about a compactification, we should introduce
a topology on the moduli space.  We will conclude this subsection with
a description of the topology used in the symplectic category.

In the symplectic case a topology on the space of stable maps can be
described via Gromov convergence \index{Gromov convergence}. The
main difficulty in defining it is that the limit map may have a
different domain (as we saw in the example depicted in \fullref{bubble}). To circumvent this difficulty define a resolution
$\kappa\co \tildeS\to\Sigma$ between two prestable curves to be a map
satisfying:

\begin{enumerate}
\item If $p\in\tildeS$ is a node, then $\kappa(p)$ is a node;

\item  If $q\in\Sigma$ is a node, then $\kappa^{-1}(q)$ is either a node or a circle disjoint from nodes;

\item If $V$ is any neighborhood of all nodes in $\Sigma$, then $\kappa\mid_{\kappa^{-1}(V)}$ is a diffeomorphism onto its image.
\end{enumerate}

One says that a sequence of stable maps $(\Sigma_k, u_k)$ Gromov
converges to a map $(\Sigma, u)$ if there is a sequence of
resolutions $\kappa_k\co \Sigma_k\to\Sigma$ such that for any
neighborhood $V$ of all nodes in $\Sigma$:

\begin{enumerate}
\item
$\bfig\morphism|b|[u_k\circ\kappa^{-1}_k`u; k\to\infty]\efig$
in $C^\infty(\Sigma\backslash V)$;

\item $\bfig\morphism|b|<600,0>[d\kappa_k{\circ}j_k{\circ}d\kappa^{-1}_k`j;
  k\to\infty]\efig$
in $C^\infty(\Sigma\backslash V)$, where $j_k$, $j$ are complex structures
on $\Sigma_k$, and $\Sigma$ respectively;

\item $\bfig\morphism|b|<800,0>[\text{Area}(u_k(\Sigma_k))`
  \text{Area}(u(\Sigma)); k\to\infty]\efig$.
\end{enumerate}

Equipped with the Gromov topology the moduli space is Hausdorff and
compact. One can check that the image homology class $\beta\in
H_2(X,\Z)$ and the (arithmetic) genus are preserved in this limit.
With this topology the coarse moduli space has the structure of a
generalization of a manifold called an orbifold. See Siebert
\cite{Siebert} and the references contained therein for more
details.
\begin{remark}\label{orbicheat}
For genus zero invariants of convex spaces such as $\CP^n$ the
moduli space is in fact a manifold. We will discuss the moduli space
as if it were a manifold until \fullref{locstr} where we explain
where the orbifold singularities arise.
\end{remark}

\subsection{The algebraic construction}

The construction of the coarse moduli space for projective varieties
works a bit differently. One can start working from the ground up by
studying some examples. The easiest example has $X$ equal to a point
(so $\beta=0$) and the genus equal to zero. This amounts to studying
configurations of $n$ distinct points in $\CP^1$ modulo equivalence
by degree one rational maps. Any degree one holomorphic map from
$\CP^1\to\CP^1$ takes the form $$\varphi([z:w])=[az+bw:cz+dw]$$
where $ad-bc\ne 0$. In affine coordinates on the domain and codomain
this takes the form $\varphi(z)=\frac{az+b}{cz+d}$. It is easy to
see that the following function for fixed $z_0$, $z_1$ and $z_2$ is
a linear fractional transformation taking $z_3=z_0$ to $0$,
$z_3=z_1$ to $1$ and $z_3=z_2$ to $\infty$.
\begin{defn}
The cross ratio is the function defined by $$
\gamma(z_0,z_1,z_2,z_3):=\frac{(z_1-z_2)(z_3-z_0)}{(z_0-z_1)(z_2-z_3)}.
$$
\end{defn}
It is easy to say what the cross ratio means -- it is nothing more
than the image of the fourth point under the unique linear
fractional transformation taking the first three points to $0$, $1$,
and $\infty$ respectively. We now follow the exposition of D
Salamon which utilizes cross ratios to realize $\wwbar
M_{0,n}(\text{pt},0)$ as projective algebraic varieties
(see Salamon \cite{salamon}). For each tuple $(i, j, k, \ell)$ of distinct
positive integers less than or equal to $n$ define a function,
$\gamma_{i,j,k,\ell}\co  M_{0,n}(\text{pt},0)\to \CP^1$ by
$\gamma_{i,j,k,\ell}([p]):=\gamma(p_i,p_j,p_k,p_\ell)$. One can
easily check that these functions satisfy the relations
\begin{eqnarray*}
&\gamma_{j,i,k,\ell}=\gamma_{i,j,\ell,k}=1-\gamma_{i,j,k,\ell},& \\
&\gamma_{i,j,k,\ell}\gamma_{i,k,j,\ell}-\gamma_{i,k,j,\ell}=
  \gamma_{i,j,k,\ell},& \\
&\gamma_{j,k,\ell,m}\gamma_{i,j,k,m}-\gamma_{j,k,\ell,m}\gamma_{i,j,k,\ell}
=\gamma_{i,j,k,m}-1.&
\end{eqnarray*}
As a particular case, we see that the cross ratio maps $M_{0,4}$
isomorphically to $\CP^1-\{0, 1, \infty\}$. It is natural to guess
that this extends to a bijection between $\wwbar{M}_{0,4}$ and
$\CP^1$, which in fact it does. The stable maps corresponding to $1$
and $\infty$ are displayed on the left of \fullref{WDVV} in the
next subsection. More generally, the cross ratios of the marked
points may be used to identify $\wwbar{M}_{0,n}$ with the
projective subvariety of $(\CP^1)^{\times N}$ specified by the
solutions to the above displayed equations in the $\gamma$. Here
$N=n(n-1)(n-2)(n-3)$ is just the number of possible distinct
$4$--tuples marked points. In the following aside we continue with a
very brief description of the algebraic construction of the coarse
moduli spaces of higher genus curves as projective varieties.

\begin{aside}{\small
We can now step things up a bit and consider moduli of higher genus
curves. Here our exposition follows that of D Mumford from
\cite{mum}. Let $\Sigma$ be a genus $g>1$ curve and $K_\Sigma$ be
the canonical bundle (top exterior power of the cotangent bundle).
Using the Riemann--Roch theorem one can compute that the dimension of
$H^0(\Sigma,K_\Sigma^3)$ is $5g-5$. The space
$H^0(\Sigma,K_\Sigma^3)$ is just the space of globally defined
holomorphic forms of the form $f(z)dz^{\otimes 3}$. Given a basis
for $H^0(\Sigma,K_\Sigma^3)$, say $\{\omega_k\}$, define a map to
$\CP^{5g-6}$ by $\phi(z)=[\omega_k(z)]$. Here we use any
trivialization of $K_\Sigma^3$ around $z$ to identify the
$\omega_k(z)$ with complex numbers. Changing the trivialization
clearly does not change the projective equivalence class. The
Weierstrass points of $\Sigma$ are defined to be those points in
$\Sigma$ for which the tangent plane to $\Sigma$ in $\CP^{5g-6}$
matches to order $5g-5$ or more. There are $g(5g-5)^2$ such points
counted with multiplicity, label them by $z_j$. Now take a large $N$
and consider the following set of functions from $5g-5$ element
subsets of $E=\{1,\ldots,g(5g-5)^g\}$ to the non-negative integers
$$
R=\bigl\{r:\{I\subset E||I|=5g-5\}\to \Z|r(J)\ge 0\ \text{for all $J$
and}\ \textstyle{\sum_{k\in I}} r(I)=N\bigr\}.
$$
Define an embedding $M_{g,0}\to \CP^{|R|-1}$ by
$$
[\Sigma]\mapsto \Biggl[\sum_{\sigma\in\text{perm}(E)}~\prod_{I\subset E,
|I|=5g-5} \left(\text{det}_{j\in
I}(\omega_k(z_j))\right)^{r(\sigma(I))}\Biggr].
$$
The determinants in the above expression lie in $K_\Sigma^{3N}$;
they can be interpreted as complex numbers by evaluation in any
trivialization.
\begin{exm}
Assuming that this construction works, jazz it up to define an
embedding of $M_{g,n}$ into a sufficiently large projective space.
This is difficult, but luckily it is not needed.
\end{exm}

There are other approaches to proving that $M_{g,n}$ and $\wwbar
M_{g,n}$ admit the structure of quasiprojective and projective
varieties respectively, but no way is easy. See Mumford \cite{GIT}
for the standard exposition.

The next step is to describe the structure of $\wwbar
M_{g,n}(\CP^r,d)$. The final step is to define $\wwbar
M_{g,n}(X,\beta)$ for general projective varieties $X$. These last
two steps are not so bad. Our exposition comes from the lectures by
Aluffi \cite{aluffi}. Given $[\Sigma,p]\in\wwbar{M}_{g,n +d(r+1)}$ such
that the divisors $(p_{n+kd+1}+\cdots +p_{n+kd+d})$ and $(p_{n+\ell
d+1}+\cdots +p_{n+\ell d+d})$ are linearly equivalent for $k,\ell =
0,\ldots, r$ and non-zero sections $s_k$ of the line bundle
associated to these equivalent divisors such that
$s_k(p_{n+kd+1})=s_k(p_{n+kd+d})=0$ one can associate a stable curve
$[u,\Sigma,q]\in\wwbar{M}_{g,n}(\CP^r,d)$ by $q_j=p_j$ for
$j=1,\ldots, n$ and $u(z)=[s_k(z)]$. It is not hard to see that two
pairs $([\Sigma,p],s)$ and $([\Sigma^\prime,p^\prime],s^\prime)$
produce the same stable map if and only if $s$ and $s^\prime$  agree
up to a constant factor and $[\Sigma,p]$ and
$[\Sigma^\prime,p^\prime]$ agree after a permutation in the marked
points fixing $p_j$ for $j=1,\ldots, n$ and each divisor
$(p_{n+kd+1}+\cdots p_{n+kd+d})$. The subset of $\wwbar{M}_{g,n
+d(r+1)}$ satisfying the equivalent divisor condition is a
subvariety, and the set of data that we have described here forms a
$(\C^\times)^{r+1}$ bundle over this subvariety. The quotient of
this bundle by the group generated by the change of scale and
permutations produces a quasiprojective variety that embeds into
$\wwbar{M}_{g,n}(\CP^r,d)$ as an open set. Of course we can embed
it in a different way by composing each map with a fixed holomorphic
isomorphism of $\CP^r$. The fact is that by choosing a finite number
of such isomorphisms we can completely cover ${\wwbar
M}_{g,n}(\CP^r,d)$. To see this, pick a basis, $\{t_k\}$, for
$H^0(\CP^r,\calO(-1))$. Then to any generic
$[u,\Sigma,q]\in\wwbar{M}_{g,n}(\CP^r,d)$ we associate
$([\Sigma,p],t_k\circ u)$, where $p$ is obtained by adjoining the
zeros of the $t_k\circ u$ to $q$.}
\end{aside}

\section{Cohomology of the moduli space}\label{GWco}
We see that our first definition of Gromov--Witten invariants is
just the evaluation of natural cohomology classes on the moduli
space. This leads one to ask if there are any other cohomology
classes on the moduli space. Indeed there are other interesting
classes. In this subsection we will define some of them and derive
several important recurrence relations between the Gromov--Witten
invariants. The definition of a new set of cohomology classes
appears in the first article and the recurrence relations are
described in subsequent ones. The paper by R Vakil \cite{vakil} is
also a good reference for this material.

Here we describe the new cohomology classes as Chern classes of natural
vector bundles over the moduli space. In general, as explained in \fullref{locstr} one needs a
generalization of vector bundles called orbibundles that allow finite quotient
singularities. However, in a number of examples singularities do not appear (\fullref{orbicheat}).
For now we assume that everything is smooth and introduce more general examples in
\fullref{locstr}.

\subsection{Gromov--Witten invariants and descendants}\label{321}

We recall the intersection \index{Chern class} \index{$C$@$c_1$
first Chern class} theory definition of the first Chern class of a
line bundle. There are many other possible definitions, see
Milnor--Stasheff \cite{MS}, Griffiths--Harris \cite{GH} and Bott--Tu
\cite{BT}. Given a line bundle $L\to X$ and two generic (transverse)
sections $\sigma_0, \sigma_1\co X\to L$, the first Chern class of
the line bundle may be defined to be the cohomology class Poincar\'e
dual to $\sigma_1^{-1}(\sigma_0(X))$.  As an example, consider the
tangent bundle to the $2$--sphere. A section of the tangent bundle
is nothing other than a vector field. We can (and generally will)
take $\sigma_0$ to be the zero section. We can take $\sigma_1$ to be
a vector field that flows up from the south pole to the north pole.
As a set $\sigma_1^{-1}(\sigma_0(S^2))$ consists of exactly two
points. One can see that the intersections are transverse, and
conclude that $c_1(TS^2)[S^2]$ must be $-2$, $0$ or $2$. It is in
fact $2$.

\begin{exm}
Write out consistent orientation conventions and verify that the
signed count of zeros implies that $c_1(TS^2)[S^2]=2$.
\end{exm}

\begin{exm}\label{Chernprop}
Prove the following properties of line bundles and their Chern classes.
\begin{enumerate}
\item $c_1(\underCC)=0$ (Here and elsewhere $\underV$ will denote
the trivial bundle with fiber $V$.)
\index{$V$@$\underV$ the trivial bundle with fiber $V$}
\item $c_1(L_1\otimes L_2)=c_1(L_1)+c_1(L_2)$.
\item $c_1(f^*L)=f^*c_1(L)$.
\item $L\otimes L^*\cong \underCC$.
\end{enumerate}
\end{exm}

There are $n$ distinguished line bundles defined over $\wwbar
M_{g,n}(X,\beta)$, denoted by $\calL_k$ \index{$L$@$\calL_k$ the
tautological bundles over moduli space} for $k=1,\ldots,n$.
Intuitively, these bundles are specified by an identification of the
fiber over each point as
$$
\calL_k|_{[u,\Sigma,p]}=T^*\Sigma_{p_k}.
$$
This allows one to define new cohomology classes on $\wwbar
M_{g,n}(X,\beta)$ and extend the definition of the Gromov--Witten
invariants.
\begin{defn}
The $\psi$--classes are defined by \index{$\psi_k:=c_1(\calL_k)$}
$\psi_k:=c_1(\calL_k)\in H^2(\wwbar{M}_{g,n}(X,\beta))$. The
descendant \index{descendant} Gromov--Witten invariants are defined
by
$$
\langle \tau_{a_1}(\gamma_1),\ldots, \tau_{a_n}(\gamma_n)
\rangle_{g,\beta}^X:=\int_{[\wwbar
{\calM}_{g,n}(X,\beta)]^{\vir}}
\psi_1^{a_1}\wedge\text{ev}_1^*\gamma_1\wedge\ldots\wedge\psi_n^{a_n}\wedge\text{ev}_n^*\gamma_n\,.
$$ \index{Gromov--Witten invariants}
Here \index{$\tau_{a}$ descendant insertion} $\gamma_1,\ldots,
\gamma_n\in H^*(X;\Q)$, $\beta\in H_2(X;\Z)$ and the integral of a
form over a space is defined to be zero if the degree of the form
does not match the dimension of the space.
\end{defn}
It is often possible to reduce the computation of Gromov--Witten
invariants on one space to a computation of descendant
invariants on a smaller space. In addition there are
recursion relations relating various descendant invariants.

\begin{figure}[ht!]
\centering
\labellist\small
\pinlabel {$1$} at 144 206
\pinlabel {$2$} at 147 148
\pinlabel {$3$} at 67 193
\pinlabel {$1$} at 360 184
\pinlabel {$2$} at 388 184
\endlabellist
\includegraphics[width=3truein]{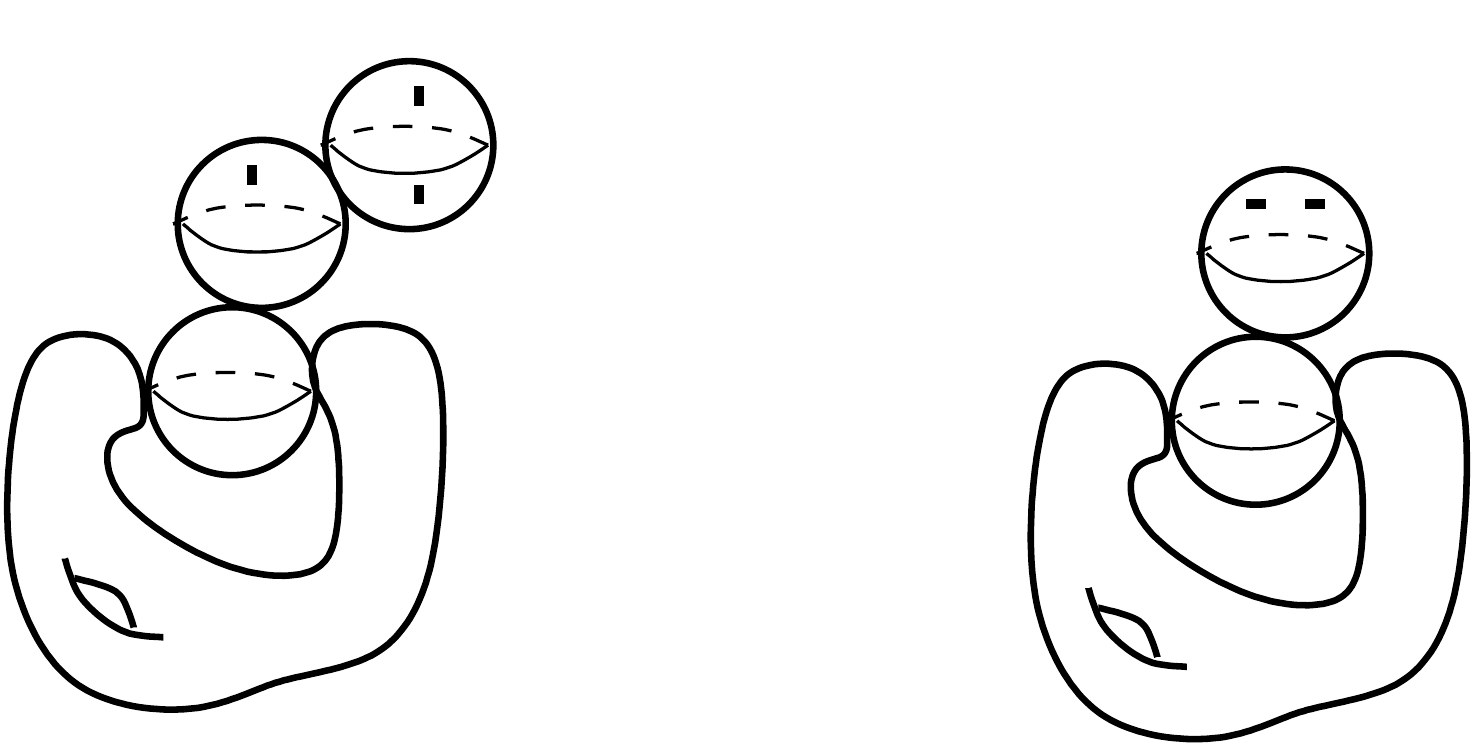}
\caption{Projection from $\wwbar{M}_{2,3}$ to $\wwbar{M}_{2,2}$}\label{proj}
\end{figure}

We should now make the definition of the bundle $\calL_k$ more
precise. To do so, we need to study the relationship between moduli
spaces with different numbers of marked points. There is a natural
projection from $\wwbar{M}_{g,n+1}(X,\beta)$ to $\wwbar
M_{g,n}(X,\beta)$ given by ignoring the final point. One subtle
point is that after deleting a marked point a stable map with marked
points may no longer be stable. This can be fixed by stabilization.
\index{$S$@st stabilization} \index{stabilization} The stabilization
of a prestable map, $\text{st}([u,\Sigma,p])$ is defined by
identifying to a point any component of the normalization of
$\Sigma$ on which an infinite subgroup of the automorphism group
acts effectively (that is, only the identity element fixes everything).
This gives,
$$
\pi\co  \wwbar{M}_{g,n+1}(X,\beta)\to \wwbar{M}_{g,n}(X,\beta); \
\pi([u,\Sigma,p_1,\ldots, p_{n+1}]):=\text{st}([u,\Sigma,p_1,\ldots,
p_n])\,.
$$
See \fullref{proj} to see the result of projection with nontrivial stabilization.

\begin{figure}[ht!]
\centering
\labellist\small
\pinlabel {$1$} at 22 6
\pinlabel {$2$} at 124 71
\pinlabel {$3$} at 58 69
\pinlabel {$1$} at 255 9
\pinlabel {$2$} at 346 113
\pinlabel {$3$} at 291 74
\pinlabel {$4$} at 382 114
\endlabellist
\includegraphics[width=4.3truein]{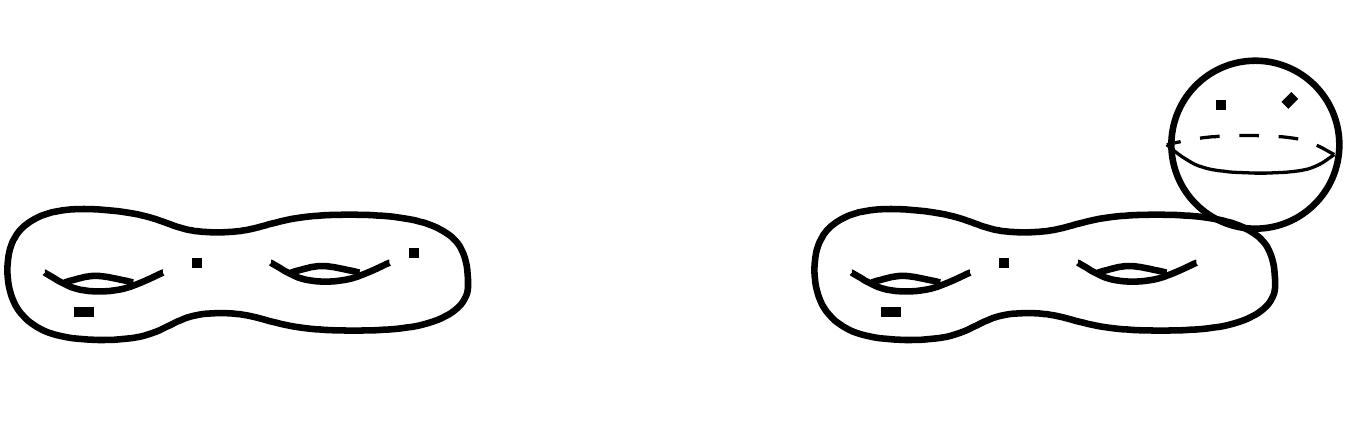}
\caption{Inclusion of $\wwbar{M}_{2,3}$ in $\wwbar{M}_{2,4}$
at $p_2$}\label{inc}
\end{figure}

There are also inclusion maps going in the other direction defined
as follows:
\begin{align*}
\rho_k&\co\wwbar{M}_{g,n}(X,\beta) \to \wwbar{M}_{g,n+1}(X,\beta); \\
\rho_k&([u,\Sigma, p_1,\ldots, p_n]):= \\
&\hspace{50pt}\bigl[\baru, \Sigma\cup_{p_k=[0:1]}\CP^1, p_1,\ldots, p_{k-1},
[1:1],p_{k+1}, \ldots, p_n, [1:0]\bigr],
\end{align*}
where $\baru$ is defined by $\baru|_{\Sigma}=u$ and
$\baru|_{\CP^1}=u(p_k)$. See \fullref{inc} to see the result of
inclusion at a marked point. \index{$\rho_k$ section of a family of
marked stable maps}

As an aside we can use these  last two figures to explain the limit
of a sequence in which one marked point approaches a node or
approaches a different marked point. If a third marked point is
added to one of the genus zero components of the curve on the right
side of \fullref{proj}, the limit as this point approaches the
node between the genus zero components will be the curve on the
left. If a fourth marked point is added to the curve on the left of
\fullref{inc}, the limit as this point approaches $p_2$ is the
curve on the right.

Setting $\calU=\wwbar{M}_{g,n+1}(X,\beta)$, the collection of maps,
$$(\pi\co \calU\to\wwbar{M}_{g,n}(X,\beta), \rho_k\co
\wwbar{M}_{g,n}(X,\beta)\to\calU, \text{ev}_{n+1}\co \calU\to X)
$$
is an example of a family of stable maps. The key property is that
for every $s_0\in\wwbar{M}_{g,n}(X,\beta)$, we have that
$[\text{ev}|_{\pi^{-1}(s_0)}, \pi^{-1}(s_0), \rho_1(s_0),\ldots,
\rho_n(s_0)]\in\wwbar{M}_{g,n}(X,\beta)$.
\begin{exm}
Given that $s_0$ has trivial automorphism group prove that $s_0$ is
isomorphic to the following stable map.
$$[\text{ev}|_{\pi^{-1}(s_0)}, \pi^{-1}(s_0), \rho_1(s_0),\ldots, \rho_n(s_0)]\,.$$ This
almost implies that $\calU$ is a `universal' family of stable maps.
What happens if $s_0$ has nontrivial automorphisms? See \fullref{calm}.
\end{exm}

We are now ready to give a precise definition of the bundles
$\calL_k$. This definition will only work as intended when stable maps have
no nontrivial automorphisms. To define $V$ in general one needs to use
the language of stacks. \index{$L$@$\calL_k$ the tautological
bundles over moduli space} For starters, define the subbundle $V$ of vertical vectors on
$\calU$ to be those vectors that are tangent to the fibers of the
projection $\pi$. Next, recall the definition of the pull-back of a vector bundle. Given a
vector bundle $\pi\co E\to Y$ and a map $f\co X\to Y$ the pull back is
defined by
$$
f^*E:=\{(x,e)\in X\times E | f(x)=\pi(e) \}.
$$
Now we have,

\begin{defn}
Let
$$
V:=\text{ker}(\pi_*\co T\calU\to T\wwbar{M}_{g,n}(X,\beta))\,.
$$
be the vertical subbundle. Then the tautological bundles over the moduli space are $\calL_k:=\rho_k^*V^*$,
where $V^*$ is the dual of $V$.
\end{defn}
We have seen that in the absence of automorphisms the inverse image
of a stable map $s_0\in\wwbar{M}_{g,n}(X,\beta)$ under the
projection from $\wwbar{M}_{g,n+1}(X,\beta)$ to $\wwbar{M}_{g,n}(X,\beta)$
is isomorphic to the original stable map $s_0$.
The image $\rho_k(s_0)$ corresponds to the marked point $p_k$ in the
isomorphic copy of $s_0$. It follows that the fiber of $\rho_k^*V$
over $s_0=[u,\Sigma,p]$ is isomorphic to $T_{p_k}\Sigma$. Thus the
fiber of $\calL_k$ over a stable map is isomorphic to the cotangent
space of the underlying curve at the associated marked point. The
analogous construction in the category of stacks will work when
there are automorphisms.

\subsection{Boundary divisors}\label{322}

We need a way to describe the cohomology of the moduli
spaces. By intersection theory we can identify a $k$--dimensional
cohomology class with a real codimension $k$ cycle. One sees that
the set of stable maps with one node has complex codimension $1$, so
may be used to define real dimension $2$ cohomology classes. The
classes defined in this way are called boundary divisors.
\index{$D$@$D(g,A,\alpha|h,B,\beta)$ boundary
divisor}\index{boundary divisors}
\begin{defn}
Given two disjoint sets $A$ and $B$ such that $A\cup B=\{1,\ldots,
n+m\}$ with order preserving bijections $j_a\co \{1,\ldots, n\}\to A$
and $j_b\co \{1,\ldots, m\}\to B$, non-negative integers $g$ and $h$
and second cohomology classes $\alpha$ and $\beta$ in $X$, we define
the boundary divisor $D(g,A,\alpha|h,B,\beta)$ to be the push-forward of
the fundamental cycle by the map,
\begin{align*}
\wwbar{M}_{g,n+1}(X,\alpha)\times\wwbar{M}_{h,m+1}(X,\beta) &\to
  \wwbar{M}_{g+h,n+m}(X,\alpha+\beta); \\
([u,\Sigma,p],[u^\prime,\Sigma^\prime,p^\prime]) &\mapsto
  [u\cup u^\prime,\Sigma\cup_{p_{n+1}=q_{m+1}}\Sigma^\prime, r],
\end{align*}
where $r_{j_a(k)}:=p_k$ and $r_{j_b(k)}:=q_k$.
\end{defn}

\begin{figure}[ht!]
\centering
\labellist\small
\pinlabel {$1$} at 193 130
\pinlabel {$2$} at 22 135
\pinlabel {$3$} at 223 145
\pinlabel {$4$} at 249 138
\pinlabel {$1$} at 610 135
\pinlabel {$2$} at 453 147
\pinlabel {$3$} at 637 148
\pinlabel {$4$} at 672 134
\pinlabel {$\alpha+\beta$} at 131 11
\pinlabel {$\alpha$} at 480 11
\pinlabel {$\beta$} at 647 11
\endlabellist
\includegraphics[width=4.3truein]{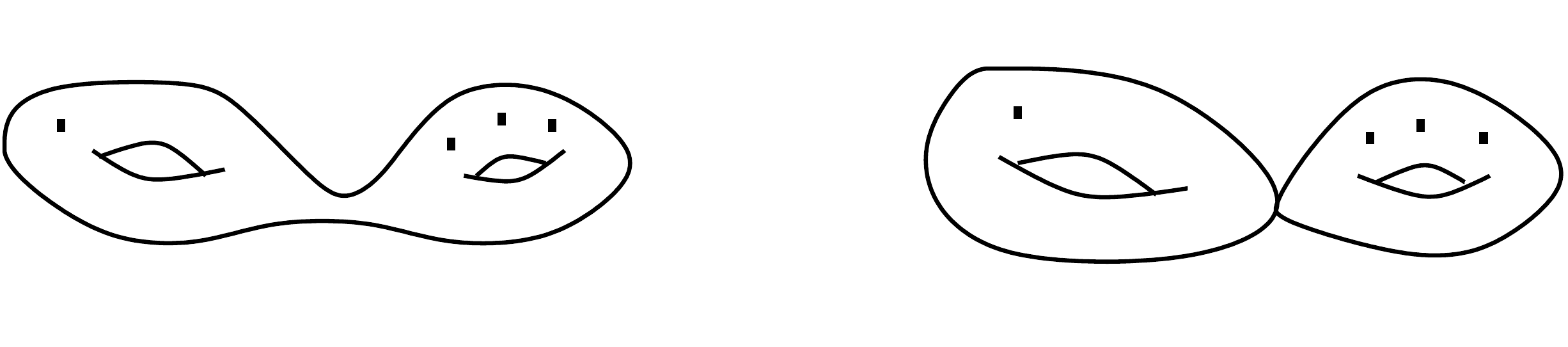}
\caption{Degeneration to the boundary divisor
$D(1,\{2\},\alpha|1,\{1,3,4\},\beta)$}\label{degeneration}
\end{figure}
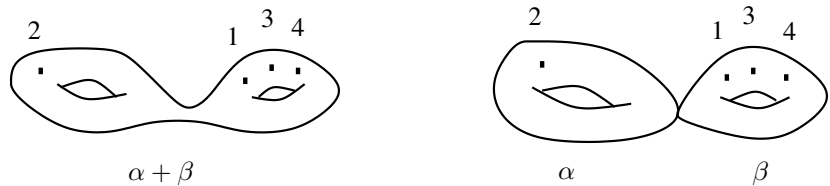

The boundary divisor $D(g,A,\alpha|h,B,\beta)$ corresponds to the
closure of the subset of moduli space of the set of stable maps with
one node joining two irreducible components having data $g,A,\alpha$
and $h,B,\beta$ respectively. See \fullref{degeneration}.
There is an additional type of boundary divisor corresponding to the
degeneration of a non-separating simple closed curve in the domain.
This divisor is typically denoted by $\Delta_0$.
\begin{exm}
Give a formal definition of $\Delta_0$.
\end{exm}

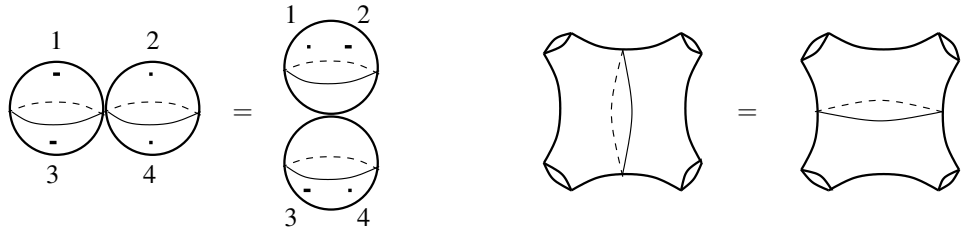
\begin{figure}[ht!]
\centering
\labellist\small
\pinlabel {$1$} at 41 151
\pinlabel {$2$} at 121 151
\pinlabel {$3$} at 38 40
\pinlabel {$4$} at 118 40
\pinlabel {$=$} at 195 90
\pinlabel {$1$} at 235 173
\pinlabel {$2$} at 295 173
\pinlabel {$3$} at 235 4
\pinlabel {$4$} at 295 4
\pinlabel {$=$} at 612 90
\endlabellist
\includegraphics[width=5truein]{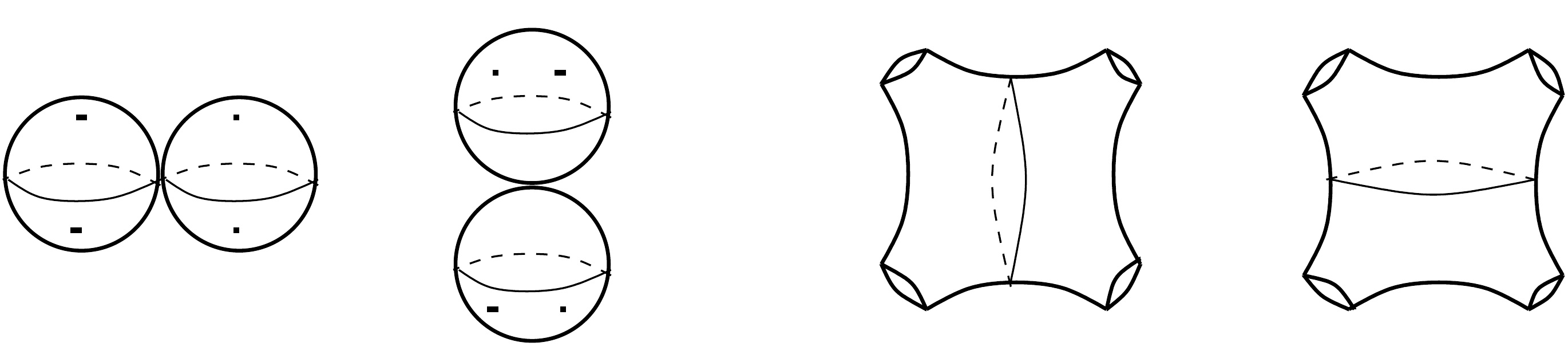}
\caption{The WDVV equation and the fusion rule}\label{WDVV}
\end{figure}

When factors in the definition of a boundary divisor are clear from
the context we will leave them out of the notation. For example
$D(\{1,3\}|\{2,4\})$ in $\wwbar{M}_{0,4}$ is short-hand notation
for $D(0,\{1,3\},0|0,\{2,4\},0)$, as there is no possible way to
have nontrivial genus or homology in this case. There is an
additional boundary divisor obtained by identifying two points on
one curve, or equivalently degenerating a non-separating simple
closed curve. This divisor is typically denoted by $\Delta_0$. The
map combining two curves into one nodal curve is very similar to
constructions gluing two surfaces with boundary along a boundary
component. This latter operation is common in the analysis of
topological quantum field theories. As an example consider the
divisors $D(\{1,2\}|\{3,4\})$ and $D(\{1,3\}|\{2,4\})$ in $\wwbar
M_{0,4}$. We have identified $\wwbar{M}_{0,4}$ with $\CP^1$ via
the cross ratio. Under this identification we are mapping the first
point to zero, the second to one and the third to infinity and then
looking at the coordinates of the last point. Thus
$D(\{1,2\}|\{3,4\})$ is identified with infinity and
$D(\{1,3\}|\{2,4\})$ is identified with one, so these divisors are
distinct. They are however linearly equivalent, which implies that
they represent the same cohomology class. The equality of these
cohomology classes is the starting point for the
Witten--Dijkgraaf--Verlinde--Verlinde equations (WDVV equations). The
WDVV equations are similar to a consequence of the fusion rule for
monoidal functors, see \fullref{WDVV}.

This similarity is one good reason to believe that there may be
some relationship between Gromov--Witten invariants and Chern--Simons
invariants.

The linear equivalence of these two simple divisors is more powerful
than it may first appear. The projection map from $\wwbar
M_{0,n}(X,\beta)$ to $\wwbar{M}_{0,4}$ defined via stabilization
whenever $n\ge 4$ translates this equivalence into many more moduli
spaces. To demonstrate the power of this equivalence, the next
exercise outlines a proof of Kontsevich's recursion for the genus
zero Gromov--Witten invariants of $\CP^2$. The answer to this
exercise is explained nicely in the book by Hori et al \cite{Hori}.
\begin{exm}\label{Nd}
Let $N_d$ be the Gromov--Witten invariant
$$\langle \text{pt},\ldots,\text{pt}\rangle_{0, d[\CP^1]}^{\CP^2}.$$
This is the number of degree $d$ parameterized curves that pass through
$3d-1$ generic points in the plane. Let $\pi\co \wwbar
M_{0,3d}(\CP^2,d[\CP^1])\to \wwbar{M}_{0,4}$ be the standard
projection. Finally, let $Y$ be a cycle in $\wwbar
M_{0,3d}(\CP^2,d[\CP^1])$ representing
$$\text{ev}_1^*[\CP^1]\text{ev}_2^*[\CP^1]\text{ev}_3^*[\text{pt}]\ldots\text{ev}_{3d}^*[\text{pt}].$$
We are denoting cohomology classes by their Poincar\'e duals. Of
course one should take a generic collection of two hyperplanes and
$3d-2$ points when writing a representative for $Y$.
\begin{enumerate}
\item Express $\pi^*D(\{1,2\}|\{3,4\})$ and $\pi^*D(\{1,3\}|\{2,4\})$ in terms of
boundary divisors on $\wwbar{M}_{0,3d}(\CP^2,d[\CP^1])$.
\item Show that $\#(Y\cap D(\{1,2\},0|\{3,\ldots,3d\},d[\CP^1]))=N_d$.
\item  Argue that $\#(Y\cap D(A,e|B,f))$ must be zero unless
$B$ has the `right' number of points. When $B$ does have the right
number of points construct a covering projection
\begin{multline*}
Y\cap D(A,e|B,f)\to
\bigl(\wwbar{M}_{0,3e-1}(\CP^2,e[\CP^1])\cap
  \text{ev}_1(q_1)\cap\ldots\text{ev}_{3e-1}(q_{3e-1})\bigr) \\
\times\bigl(\wwbar{M}_{0,3f-1}(\CP^2,f[\CP^1])\cap
  \text{ev}_1(r_1)\cap\ldots\text{ev}_{3f-1}(r_{3f-1})\bigr).
\end{multline*}
\item Express $\#(Y\cap\pi^*D(\{1,2\}|\{3,4\}))$ in terms of the numbers $N_e$. Hint: When
computing $\#(Y\cap D(A,e|B,f))$ recall that a degree $e$ curve
intersects a degree $f$ curve in $ef$ points when you are
enumerating the locations of the node and the two marked points that
lie on the lines.
\item Express $\#(Y\cap\pi^*D(\{1,3\}|\{2,4\}))$ in a similar way.
\item Use the linear equivalence to deduce a recurrence relation between the various $N_d$,
and test your recurrence by computing a few values. You should get $N_2=1$, $N_3=12$, and $N_4=620$.
\end{enumerate}
\end{exm}
If these Gromov--Witten invariants are put into a suitable
generating function the recurrence relation will translate into a
differential equation. One important example of this allows one to
conclude that the descendant Gromov--Witten invariants of a point
are encoded into the solutions of the KdV equation. This may be
considered as one of the first tests of Large $N$ Duality for the
case of zero-dimensional space. We however, will not emphasize this
aspect of Large $N$ Duality here. See the paper by Witten on 2D
gravity \cite{2D} for a discussion of this when it was still a
conjecture and the paper by Kontsevich for a proof \cite{KontAiry}.

\subsection{The string and dilaton equations}\label{323}

We now turn to a description of the recursion relations between the
$\psi$--classes. For clarity, we will denote the first $\psi$ class
on $\wwbar{M}_{g,n}(X,\beta)$ by $\psi_1^{(n)}$ when we are
comparing $\psi$--classes on $\wwbar{M}_{g,n}(X,\beta)$ and
$\wwbar{M}_{g,n+1}(X,\beta)$. We will use a similar notation for
the tautological bundles. We have by definition,
$$\psi_1^{(n+1)}-\pi^*\psi_1^{(n)}=c_1(\calL_1^{(n+1)}\otimes(\pi^*\calL_1^{(n)})^*)\,.
$$
To compute this Chern class, we just need to construct a section of the bundle. A section of
this bundle may be viewed as a vector bundle morphism,
$s\co \rho_1^*V^{(n+1)}\to \pi^*\rho_1^*V^{(n)}$. We can view tangent vectors to the
universal bundle as paths of stable maps. A natural morphism is given by
\begin{multline*}
s([u,\Sigma,p_1,\ldots,p_{n+1}],[u,\Sigma,p_1,\ldots,p_{n+1},q(t)]):=\\
([u,\Sigma,p_1,\ldots,p_{n+1}],\text{st}[u,\Sigma,p_1,\ldots,
  p_{n}],\text{\underline{st}}[u,\Sigma,p_1,\ldots,p_{n},q(t)]).
\end{multline*}
Here $q(t)$ is a smooth path in $\Sigma$ with $\lim_{t\to
0}q(t)=p_1$, and the above formula is correct when $t\ne 0$. The
underline on the second stabilization refers to the fact that it is
the result of stabilizing the $[u,\Sigma,p]$ and then placing the
$q(t)$ at the corresponding point of the stabilization.  The above
formula appears complicated because of the pull-backs in the
definitions of the two bundles. This bundle
map is identically zero over the boundary divisor
$D(0,\{1,n+1\},0|g,\{2,\ldots,n\},\beta)$ because the entire bubble
containing $p_1$ and $p_{n+1}$ collapses to a point when $p_{n+1}$
is forgotten. The resulting map is stabilized and thus the path
$q(t)$ becomes constant. In fact, this is the only way that this
section can be identically zero over some point. One concludes
$$
\psi_1^{(n+1)}-\pi^*\psi_1^{(n)}=
D(0,\{1,n+1\},0|g,\{2,\ldots,n\},\beta).
$$
The recursion between cohomology classes that we just derived
implies three important recursion relations between descendant
Gromov--Witten invariants. We will prove two of these recursions;
the string equation and the dilaton equation. This will require use
of the Thom isomorphism and the Euler class of a vector bundle.

Recall that relative cohomology classes may be represented as the
Poincar\'e duals of closed cycles. Given a real rank $r$ vector
bundle, $\pi\co E\to X$ over a closed base, there is a natural relative
cohomology class, $\Phi\in H^r(E,E-\sigma_0(X))$  (the Thom class)
given as the Poincar\'e dual of $\sigma_0(X)$.  The Thom isomorphism
is the map \index{$\Phi$ the Thom class} \index{Thom
class/isomorphism} \index{$E$@$e(E)$ the Euler class}\index{Euler
class}
\begin{align*}
H^k(X) &\to H^{r+k}(E,E-\sigma_0(X)) \\
\alpha &\mapsto \Phi\cup\pi^*\alpha.
\end{align*}
The Euler class of the vector bundle is the pull-back of the Thom
class by a section, $e(E):=\sigma_1^*\Phi\in H^r(X)$. Notice that
the Euler class of the real bundle underlying a complex line bundle
is the first Chern class of the line bundle. It is easier to put all
of this on a firm theoretical foundation by using cohomology;
however the conceptual picture is harder to follow. See for example
Milnor--Stasheff \cite{MS} Spanier \cite{span}.

Using this technology we will establish two important recursion
relations for Gromov--Witten invariants. The first is the string
equation, \index{string equation}
\[
\langle \tau_{a_1}(\gamma_1)\ldots
\tau_{a_n}(\gamma_n)\tau_0(X)\rangle_{g,\beta}^X = \sum_{k=1}^n
\langle \tau_{a_1}(\gamma_1)\ldots \tau_{a_k-1}(\gamma_k)\ldots
\tau_{a_n}(\gamma_n)\rangle_{g,\beta}^X\,,
\] and the second is the dilaton equation, \index{dilaton equation}
\[
\langle \tau_{a_1}(\gamma_1)\ldots
\tau_{a_n}(\gamma_n)\tau_1(X)\rangle_{g,\beta}^X =(2g-2+n) \langle
\tau_{a_1}(\gamma_1)\ldots
\tau_{a_n}(\gamma_n)\rangle_{g,\beta}^X\,.
\]
We can now prove the string equation (we denote cycles and their
Poincar\'e duals by the same symbols). To start the computation, we
use the definition of the descendant Gromov--Witten invariants and
the relation between the $\psi$--classes on the moduli of
$(n{+}1)$--marked curves and the moduli of $n$--marked curves. After
expanding the powers, we recognize that the first integral is
trivial because the moduli space of $n$--marked curves is lower
dimensional than the moduli space of $(n{+}1)$--marked curves. Any
cycle containing a factor of $D_iD_j$  for $i{\ne}j$ is empty
because the set of all stable maps with $p_i$ and  $p_{n+1}$
isolated in a single bubble is disjoint from the set of all stable
maps with just $p_j$ and $p_{n+1}$ in an isolated bubble, so the
last integral below is trivial. Using
$\text{ev}_k\circ\pi=\text{ev}_k$ on the remaining integral gives
\begin{align*}
\langle \tau_{a_1}&(\gamma_1)\ldots
\tau_{a_n}(\gamma_n)\tau_0(X)\rangle_{g,\beta}^X
  = \int_{[\wwbar{\calM}_{g,n+1}(X,\beta)]^{\vir}}
 \psi_1^{a_1}\ldots \psi_n^{a_n} \text{ev}_1^*\gamma_1
 \ldots\hbox{ev}_n^*\gamma_n\\
& =\int_{[\wwbar{\calM}_{g,n+1}(X,\beta)]^{\vir}}
  (\pi^*\psi_1+D_1)^{a_1}\ldots (\pi^*\psi_n+D_n)^{a_n}
  \hbox{ev}_1^*\gamma_1 \ldots\hbox{ev}_n^*\gamma_n\\
& =\int_{[\wwbar{\calM}_{g,n+1}(X,\beta)]^{\vir}} \pi^*(\psi_1^{a_1}\ldots
  \psi_n^{a_n} \hbox{ev}_1^*\gamma_1 \ldots\hbox{ev}_n^*\gamma_n)\\
& \ +\sum_{k=1}^n\sum_{p=1}^{a_k} \binom{a_k}{p}
  \int_{[\wwbar{\calM}_{g,n+1}(X,\beta)]^{\vir}}
  D_k^p\pi^*(\psi_1^{a_1}\ldots \psi_n^{a_n}\psi_k^{-p}
  \hbox{ev}_1^*\gamma_1 \ldots\hbox{ev}_n^*\gamma_n)\\
& \ +\sum \int_{[\wwbar{\calM}_{g,n+1}(X,\beta)]^{\vir}}
  \bigl(D_iD_j\ldots\bigr).
\end{align*}
We next change the integrand on the remaining integral to $N(D_k)$
because $D_k$ is supported in this tubular neighborhood. The
equation we use in the second step below,
$\int_E\Phi\beta=\int_X\sigma^*\beta$, is easy to understand from
the view of intersection theory. Any variation of a stable map in
the boundary divisor
$$D_k:=D(0,\{k,n+1\},0|g,\{1,\ldots,\hat k, \ldots, n\},\beta)$$
may be decomposed into variations that move
$p_k$ and $p_{n+1}$ out of a common bubble and variations that leave
these two points in a common bubble. The normal bundle of $D_k$
consists of those variations that move the points out of a common
bubble, that is, the vertical bundle restricted to $D_k$. It follows
that the Euler class of the normal bundle is given by
$$e(N(D_k)):=c_1(\calL_k^*)=-\psi_k.$$
Together with a combinatorial
identity from the binomial theorem, this allows us to complete the
derivation of the string equation.
\begin{align*}
& =\sum_{k=1}^n\sum_{p=1}^{a_k} \binom{a_k}{p} \int_{N(D_k)}
  D_k^p\pi^*(\psi_1^{a_1}\ldots
  \psi_n^{a_n}\psi_k^{-p} \text{ev}_1^*\gamma_1
  \ldots\text{ev}_n^*\gamma_n)\\[-0.5ex]
& =\sum_{k=1}^n\sum_{p=1}^{a_k} \binom{a_k}{p}
  \int_{[\wwbar{\calM}_{g,n}(X,\beta)]^{\text{\begin{small}vir\end{small}}}}
  \rho_k^*\bigl( D_k^{p-1}\pi^*(\psi_1^{a_1}\ldots \psi_n^{a_n}\psi_k^{-p}
  \text{ev}_1^*\gamma_1 \ldots\text{ev}_n^*\gamma_n)\bigr)\\[-0.5ex]
& =\sum_{k=1}^n\sum_{p=1}^{a_k} \binom{a_k}{p} \int_{[\wwbar
  {\calM}_{g,n}(X,\beta)]^{\text{\begin{small}vir\end{small}}}}
  e(N(D_k))^{p-1}\psi_1^{a_1}\ldots \psi_n^{a_n}\psi_k^{-p}
  \text{ev}_1^*\gamma_1 \ldots\text{ev}_n^*\gamma_n 
\end{align*}
\begin{align*}
& =\sum_{k=1}^n\sum_{p=1}^{a_k} \binom{a_k}{p} (-1)^{p-1}
  \int_{[\wwbar{\calM}_{g,n}(X,\beta)]^{\text{\begin{small}vir\end{small}}}}
 \psi_1^{a_1}\ldots \psi_n^{a_n}\psi_k^{-1} \text{ev}_1^*\gamma_1
 \ldots\text{ev}_n^*\gamma_n \\[-0.5ex]
& =\sum_{k=1}^n \int_{[\wwbar{\calM}_{g,n}(X,\beta)]^{\vir}}
 \psi_1^{a_1}\ldots \psi_n^{a_n}\psi_k^{-1} \text{ev}_1^*\gamma_1
 \ldots\text{ev}_n^*\gamma_n \\[-0.5ex]
& =\sum_{k=1}^n \langle \tau_{a_1}(\gamma_1)\ldots
  \tau_{a_k-1}(\gamma_k)\ldots \tau_{a_n}(\gamma_n)\rangle_{g,\beta}^X
\end{align*}

\begin{example}\label{psi3}
We will use the string equation to compute a couple of descendant
invariants. First consider the integral
$$
\int_{[\wwbar
{\calM}_{0,6}(\text{pt},0)]^{\vir}}\psi_1^2\psi_2\,.
$$
This integral is denoted by $\langle
\tau_2(\text{pt})\tau_1(\text{pt})(\tau_0(\text{pt}))^4\rangle_{0,0}^{\text{pt}}$.
Using the string equation we obtain
\begin{align*}
\langle \tau_2(\text{pt})\tau_1(\text{pt})
  (\tau_0(\text{pt}))^4\rangle_{0,0}^{\text{pt}}
&= \langle \tau_1(\text{pt})\tau_1(\text{pt})
  (\tau_0(\text{pt}))^3\rangle_{0,0}^{\text{pt}}
  +\langle \tau_2(\text{pt})(\tau_0(\text{pt}))^4\rangle_{0,0}^{\text{pt}}
  \\[-0.5ex]
&=3\langle \tau_1(\text{pt})(\tau_0(\text{pt}))^3\rangle_{0,0}^{\text{pt}}
=3\langle (\tau_0(\text{pt}))^3\rangle_{0,0}^{\text{pt}} = 3\,.
\end{align*}
More generally consider the integral
$$\int_{[\wwbar{\calM}_{0,n}(\text{pt},0)]^{\vir}}\psi_1^{a_1}
  \ldots\psi_k^{a_k},$$
where $a_1+\cdots+a_k=n-3$. Using the string equation this can be
reduced to a sum of terms with the power of one of the $\psi$
decremented. Repeating this, a total of $n-3$ subtractions without
combining like terms will lead to a sum of ones. Each term of this
sum will correspond to selecting $a_1$ of the subtractions, then
$a_2$ of the subtractions, etc. This shows that the value of the
integral is $\binom{n-3}{a_1 a_2 \ldots a_k}$. The symmetries in
this formula arising from rearranging the $a_j$ correspond to the
geometric operations on the moduli space obtained by rearranging the
marked points.
\end{example}

The derivation of the dilaton equation is similar.
\begin{multline*}
\langle \tau_{a_1}(\gamma_1)\ldots
  \tau_{a_n}(\gamma_n)\tau_1(X)\rangle_{g,\beta}^X
  = \int_{[\wwbar{\calM}_{g,n+1}(X,\beta)]^{\vir}} \psi_1^{a_1}\ldots
  \psi_n^{a_n}\psi_{n+1} \text{ev}_1^*\gamma_1
  \ldots\text{ev}_n^*\gamma_n\\[-1ex]
 =\int_{[\wwbar{\calM}_{g,n+1}(X,\beta)]^{\vir}}
  (\pi^*\psi_1+D_1)^{a_1}\ldots (\pi^*\psi_n+D_n)^{a_n}\psi_{n+1}
  \text{ev}_1^*\gamma_1 \ldots\text{ev}_n^*\gamma_n 
\end{multline*}
\begin{multline*}
 =\int_{[\wwbar{\calM}_{g,n+1}(X,\beta)]^{\vir}}
  \psi_{n+1}\pi^*(\psi_1^{a_1}\ldots \psi_n^{a_n}
  \text{ev}_1^*\gamma_1 \ldots\text{ev}_n^*\gamma_n)\\[-0.5ex]
 +\sum_{k=1}^n\sum_{p=1}^{a_k} \binom{a_k}{p}
  \int_{[\wwbar{\calM}_{g,n+1}(X,\beta)]^{\vir}}
  \psi_{n+1}D_k^p\pi^*(\psi_1^{a_1}\ldots \psi_n^{a_n}\psi_k^{-p}
  \text{ev}_1^*\gamma_1 \ldots\text{ev}_n^*\gamma_n)\\[-2ex]
 +\sum \int_{[\wwbar{\calM}_{g,n+1}(X,\beta)]^{\vir}}
  \bigl(D_iD_j\ldots\bigr)
\end{multline*}
Two new observations need to be used. The first is that
the evaluation of $\psi_{n+1}$ on the fiber of the projection is the
evaluation of the first Chern class of the cotangent bundle of a
surface (on the surface), which is $2g-2$. The second observation is
that the pull-back of $\psi_{n+1}$ under the natural map that
interchanges $p_k$ and $p_{n+1}$ is $\psi_k$. The restriction of
this map to the divisor $D_k$ is trivial, so these two bundles agree
over this divisor.
\begin{align*}
& =(2g-2)\int_{[\wwbar{\calM}_{g,n}(X,\beta)]^{\vir}}
  \pi^*(\psi_1^{a_1}\ldots \psi_n^{a_n}
  \text{ev}_1^*\gamma_1 \ldots\text{ev}_n^*\gamma_n)\\[-0.5ex]
& +\sum_{k=1}^n\sum_{p=1}^{a_k} \!\binom{a_k}{p}\!
  \int_{[\wwbar{\calM}_{g,n}(X,\beta)]^{\vir}} \rho_k^*\bigl(
  \psi_{n+1}D_k^{p-1}\pi^*(\psi_1^{a_1}\ldots \psi_n^{a_n}\psi_k^{-p}
  \text{ev}_1^*\gamma_1
  \ldots\text{ev}_n^*\gamma_n)\bigr) \\[-0.5ex]
& =(2g-2)\int_{[\wwbar{\calM}_{g,n}(X,\beta)]^{\vir}}
  \pi^*(\psi_1^{a_1}\ldots \psi_n^{a_n}
  \text{ev}_1^*\gamma_1 \ldots\text{ev}_n^*\gamma_n)\\[-0.5ex]
&+\sum_{k=1}^n\sum_{p=1}^{a_k} \binom{a_k}{p}
  \int_{[\wwbar{\calM}_{g,n}(X,\beta)]^{\vir}} \psi_k
  e(N(D_k))^{p-1}\psi_1^{a_1}\ldots \psi_n^{a_n}\psi_k^{-p}
  \text{ev}_1^*\gamma_1 \ldots\text{ev}_n^*\gamma_n \\[-0.5ex]
& =(2g-2)\int_{[\wwbar{\calM}_{g,n}(X,\beta)]^{\vir}}
  \pi^*(\psi_1^{a_1}\ldots \psi_n^{a_n}
  \text{ev}_1^*\gamma_1 \ldots\text{ev}_n^*\gamma_n)\\[-0.5ex]
&+\sum_{k=1}^n\sum_{p=1}^{a_k} \binom{a_k}{p} (-1)^{p-1}
  \int_{[\wwbar{\calM}_{g,n}(X,\beta)]^{\vir}}
 \psi_1^{a_1}\ldots \psi_n^{a_n} \text{ev}_1^*\gamma_1
 \ldots\text{ev}_n^*\gamma_n \\[-0.5ex]
&=(2g-2+n) \langle \tau_{a_1}(\gamma_1)\ldots
  \tau_{a_n}(\gamma_n)\rangle_{g,\beta}^X
\end{align*}
To conclude this subsection we should remark that there is a
different natural family of cohomology classes on these moduli
spaces that arise as the Chern classes of the Hodge bundle. This
bundle will appear after we discuss the deformation-obstruction
sequence in \fullref{defobsec}. In addition there is a third
important recurrence relation known as the divisor equation
(see Hori et al \cite{Hori}).

\section{Local structure of moduli spaces}\label{locstr}
We begin this section by pointing out what seems to be a paradox in
this theory. When the Gromov--Witten invariants of certain integral
classes are computed one can get non-integral answers. A simple
example can explain the origin of this paradox. This is how we begin
the following subsection.

\subsection{Orbifolds and ${\wwbar{M}}_{1,1}$}

\begin{figure}[ht!]
\centering
\labellist\small
\pinlabel {$\mathbb{Z}_2$} [l] at 431 78
\pinlabel {$\mathbb{Z}_3$} [t] at 374 38
\pinlabel {$\mathbb{Z}_4$} [b] at 374 160
\endlabellist
\includegraphics[width=3truein]{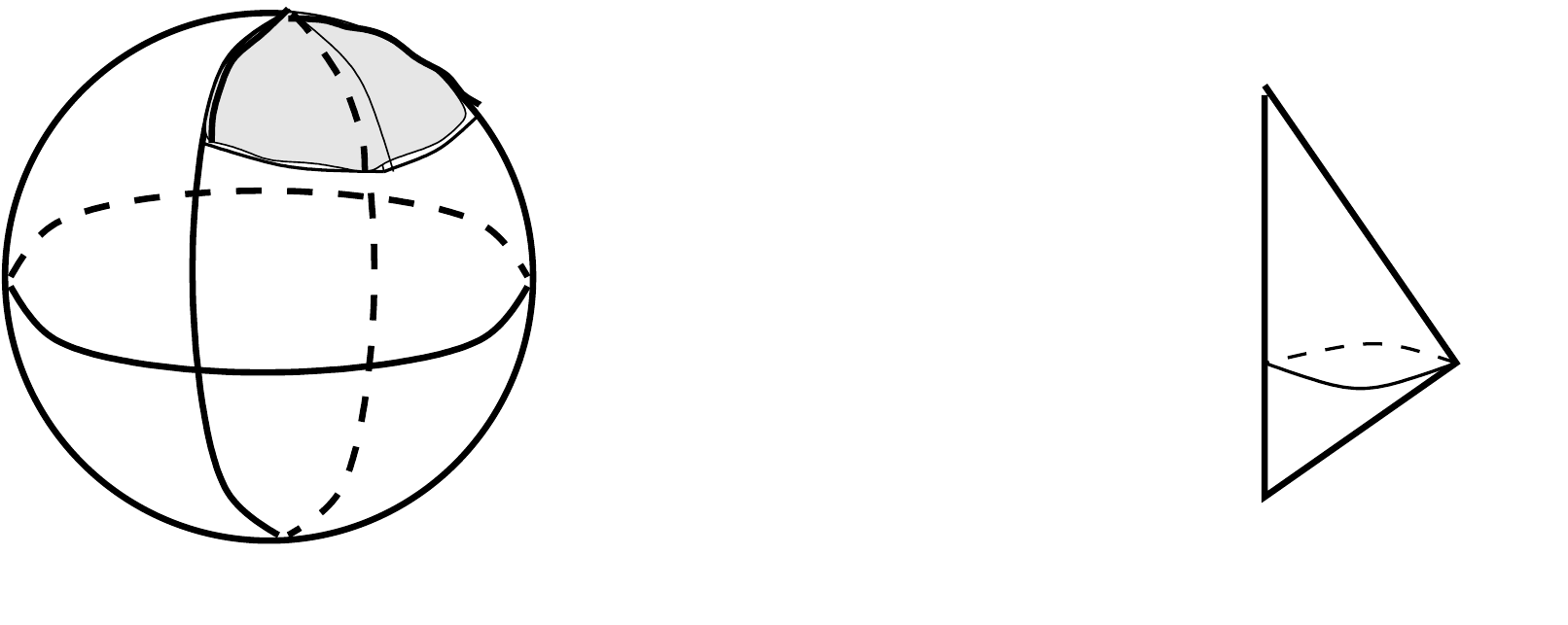}
\caption{The octahedral orbifold}\label{octspace}
\end{figure}

The group of orientation preserving isometries of an octahedron acts
on $S^2$ and on $TS^2$ in a natural way. The quotient of $S^2$ by
this group action is topologically a new copy of $S^2$ and the
degree of the quotient projection is $24$. Note that the projection
is not a covering projection, rather it is a branched cover,
branched over the points labeled with finite cyclic groups in
\fullref{octspace}. The cyclic group labels are just the stabilizer
groups of the points in the preimage of each branch point. By
keeping track of the local symmetry groups, one can develop a theory
of branched covers that is very close to the theory of covering
spaces. This gives rise to the notions of orbifolds, orbibundles and
their characteristic classes. Orbibundles are not locally trivial
over labeled points, but rather look (locally) like quotients of
such bundles by the label groups.

Thus an orbifold is locally a quotient of Euclidean space by a
\index{orbifold} finite group where one keeps track of the local
symmetry groups. Note that orbifolds are not manifolds with
singularities. Topologically, the points labeled with the finite
cyclic groups in \fullref{octspace} are locally homeomorphic to
$\C$ (although this need not be the case in general), but the
stabilizer groups make these points special.

The formal definition is a Deligne--Mumford stack over the category
DIFF. This is formalized in \fullref{app:a}, but we wanted to give
an intuitive description first. We work out further examples in
\fullref{mc} below.

Continuing with the quotient of $S^2$ by the octahedral group, note
that by the third part of exercise in \fullref{Chernprop}, we know that the
Chern class of a pull-back bundle is the pull-back of the Chern
class of the bundle. Since the first Chern class of $TS^2$ evaluates
to $2$ on the fundamental class of $S^2$ we would have to conclude
that the first Chern class of the quotient bundle of $TS^2$ by the
orientation preserving octahedral group would have to be
$\frac{1}{12}$. Of course there is a map $TS^2/\sim\to S^2/\sim$;
however, this map does not give $TS^2/\sim$ the structure of a
vector bundle. It gives it the structure of an orbibundle.
Chern classes of orbibundles are defined; however, they are not
always integral classes.

\begin{figure}[ht!]
\centering
\labellist\small
\pinlabel {$-2$} [t] at 145 24
\pinlabel {$-1$} [t] at 202 24
\pinlabel {$0$} [t] at 258 24
\pinlabel {$1$} [t] at 315 24
\pinlabel {$2$} [t] at 372 24
\endlabellist
\includegraphics[width=4truein]{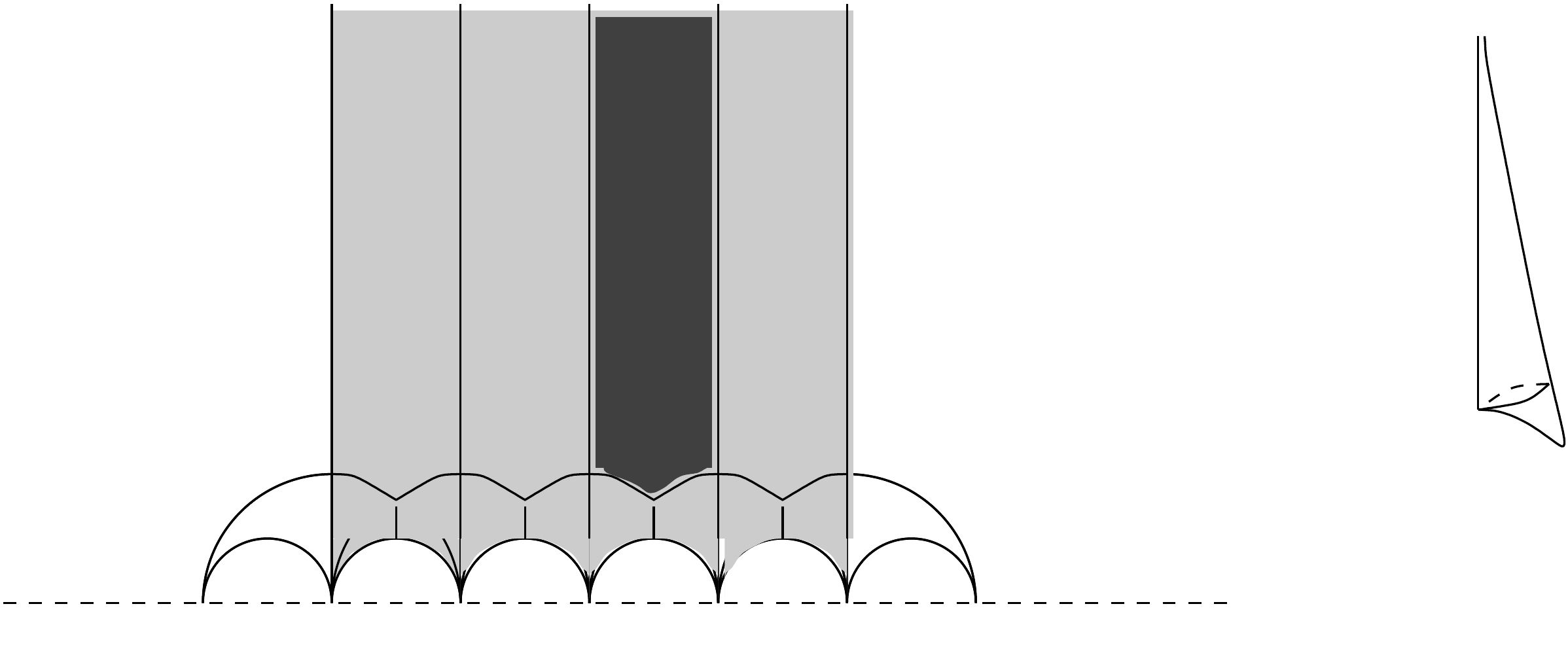}
\caption{The orbifold $M_{1,1}$}\label{m11}
\end{figure}
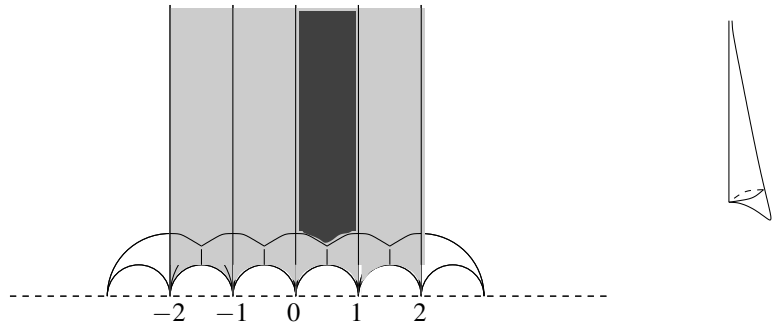

To see that this situation arises in Gromov--Witten theory
consider the structure of the moduli space $\wwbar{M}_{1,1}:=\wwbar
M_{1,1}(\text{point},0)$. It is a fact that every smooth genus $1$
Riemann surface is equivalent to the quotient of $\C$ by a lattice
(see Griffiths--Harris \cite{GH}). All of these admit a natural
group structure, so without loss of generality we may assume that the
marked point is the equivalence class of zero. Next, notice that any
biholomorphism between a pair of tori must be induced from a linear
endomorphism of $\C$. To see this, compose a given biholomorphism,
$\varphi\co T\to T^\prime$, with the projection $\C\to T$ and notice that
the resulting map lifts to a biholomorphism $\wbar \varphi\co \C\to\C$.
Now $(\wbar\varphi(z^{-1}))^{-1}$ has an isolated singularity at zero that
must be removable as this function is bounded near zero.  It follows
that the map $\wbar\varphi$ extends to a biholomorphism of $\CP^1$
taking infinity to infinity. Such must be a linear map when restricted
to $\C$. The map $z\mapsto -z$ is an automorphism of any torus, thus
every point in $\wwbar{M}_{1,1}$
has stabilizer containing $\Z_2$.
\begin{exm}
What is the stabilizer of the unique nodal curve in $\wwbar{M}_{1,1}$?
\end{exm}
We may assume that one generator of any lattice used to construct a
torus is $1$ and that the other generator is in the upper half plane
by applying a suitable linear transformation. Equivalently, the
second generator is chosen so that the generators of the lattice
form a positively oriented basis.
\begin{figure}[ht!]
\centering
\labellist\small
\pinlabel {$\omega$} [b] at 55 130
\pinlabel {$1$} [t] at 115 28
\endlabellist
\includegraphics[width=1.5truein]{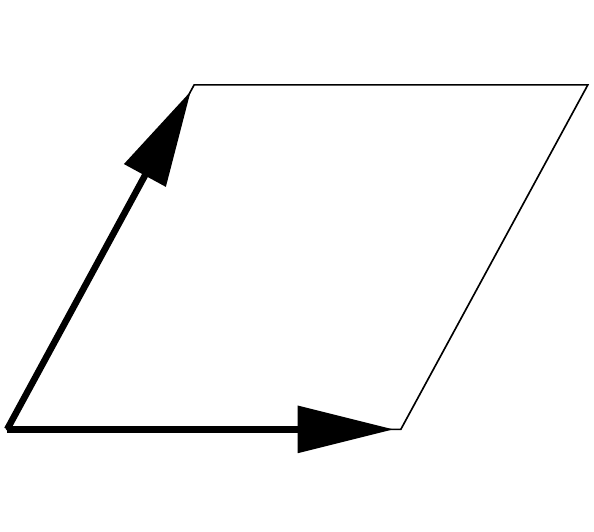}
\caption{Lattice generators}\label{funddom}
\end{figure}
It is standard to parameterize a once-marked torus by the location
of the second lattice point. The choice of second lattice point is
however not unique. Clearly, adding one to the second generator $z\mapsto z+1$ does
not change the torus. Similarly, interchanging the
two generators and changing the sign of one and then renormalizing
does not change the torus, $z\mapsto -z^{-1}$. These two operations
generate an action of $SL_2\Z$ on the upper half plane. The quotient
of the upper half plane by this action may be identified with
$M_{1,1}$; adding a single point at the cusp representing the nodal
marked torus gives $\wwbar{M}_{1,1}$. A fundamental domain for
this action and the quotient is displayed in \fullref{m11}.

We can now analyze the group of automorphisms of one of these marked
tori. Consider the torus labeled by $\omega:=e^{i\frac{\pi}{3}}$. We
have seen that any automorphism must be of the form $[z]\mapsto
[\alpha z]$. Since $1$ is equivalent to $0$ modulo the lattice it
must be mapped to a new lattice point, say $\alpha=a+b\omega$. This
implies that $\omega$ gets mapped to a lattice point,
$(a+b\omega)\omega=-b+(a+b)\omega$. This allows us to represent the
automorphism as a real linear transformation by a $2\times 2$
matrix. Since the map is invertible and takes the lattice to itself
the determinant must be $\pm 1$. Orientation preserving implies that
the determinant must be $1$, so
$$
1=\text{det}\begin{bmatrix}a&-b\\b&a+b\end{bmatrix}=a^2+ab+b^2=
\bigl(a+\tfrac12 b\bigr)^2+\bigl(\tfrac{\sqrt{3}}{2}b\bigr)^2.
$$
It follows that $(a,b)$ must be one of the six pairs, $\pm(1,0)$,
$\pm(0,1)$ or $\pm(1,-1)$. These correspond exactly to the rotations
generated by multiplication by powers of $\omega$. Notice that the
computation of the stabilizer of $\omega$ under the action of
$SL_2\Z$ on the upper half plane is exactly the same, thus the upper
half plane quotient accurately models the orbifold structure on
$M_{1,1}$. In particular notice that $SL_2\Z$ does not act
effectively on the upper half plane as every point is fixed by $-I$.
This corresponds to the fact that a generic elliptic curve has
$\Z_2$ automorphism group. Thus, $PSL_2\Z$ is definitely not the
right group to study from the point of view of algebraic geometry.
\begin{exm}
Compute the automorphism group of the marked torus labeled by $i$.
\end{exm}

We will now take a detour via exercises to prove that $SL_2\Z$ is
generated by two matrices and sketch a proof that
$\langle\tau_1(\text{pt})\rangle_{1,0}^{\text{pt}}=\frac{1}{24}$.
The matrices are
$$
S=\begin{pmatrix}0&-1\\1&0\end{pmatrix}, \quad
T=\begin{pmatrix}1&1\\0&1\end{pmatrix}.
$$
After the detour, we will give the correct definition of the moduli
stack. This encodes the orbifold structure of the moduli space.

Let $PSL_2\Z$ be the quotient of $SL_2\Z$ by $\pm 1$. This acts on
the upper half plane by linear fractional transformations. One can
define an orbifold fundamental group. This group should satisfy a
van Kampen theorem and the orbifold fundamental group of an
$X/\Gamma$ should be $\Gamma$ when $X$ is a simply connected
orbifold and the action is nice (see Ratcliffe \cite{thurston}).
\begin{exm}
Conclude that $\pi_1^{orb}(\H/PSL_2\Z)=PSL_2\Z$. Now use the
orbifold van Kampen theorem and \fullref{m11} to prove that this
is $\langle [S], [W] | [S]^2=1, [W]^3=1\rangle$, where $W=ST$ and
$[A]$ is the image of $A\in SL_2\Z$ under the natural projection.
\end{exm}
\begin{exm}
Check that the map from the group $\langle s, t|s^2=(st)^3,
s^4=1\rangle$ to $SL_2\Z$ taking $s$ to $S$ and $t$ to $T$ is a
well-defined group homomorphism. Show that the kernel of the
composition of this map with projection to $PSL_2\Z$ is $\Z_2$.
Conclude that $SL_2\Z\cong\langle s, t|s^2=(st)^3, s^4=1\rangle$.
\end{exm}
\begin{exm}
Give a definition of a line bundle in the orbifold category. Also
define the first Chern class of an orbifold bundle. You should see
that the Chern class may be computed by counting zeros with sign,
but one must divide each term by the order of the stabilizer group.
Also prove that this Chern class satisfies the usual multiplication
by degree for pull-backs.
\end{exm}
Now $\smash{\langle\tau_1(\text{pt})\rangle_{1,0}^{\text{pt}}}$ is
just the evaluation of the first Chern class of the tautological
bundle over $\wwbar{M}_{1,1}$. To do this we would like to have a
section of $\calL_1$. Such a section would associate a holomorphic
form to each marked torus. It is natural to try a section that would
associate $f(\tau)dz$ to $\C/\langle 1,\tau\rangle$. Transformations
by elements of $SL_2\Z$ fix the torus and so should fix the
holomorphic form. Multiplication by $c\tau+d$ is an isomorphism
$\C/\bigl\langle 1,\frac{a\tau+b}{c\tau+d} \bigr\rangle\to\C/\langle
c\tau+d,a\tau+b\rangle=\C/\langle 1,\tau\rangle$. It follows that we
need $f\bigl(\frac{a\tau+b}{c\tau+d}\bigr)=(c\tau+d)f(\tau)$ in order
for this section to be well defined.
\begin{defn}
The Dedekind eta function \index{Dedekind eta function} is the
analytic function on the upper half plane given by
$$
\eta(\tau):=e^{\frac{\pi i\tau}{12}}\prod_{n=1}^\infty (1-e^{2\pi in\tau})\,.
$$
\end{defn}
The following functional equation is proved in Apostol \cite{apostol}.
$$
\eta^2\left(\frac{a\tau+b}{c\tau+d}\right)=-i\epsilon(a,b,c,d)^2(c\tau+d)\eta^2(\tau)\,,
$$
where
$$
\epsilon(a,b,c,d):=\text{exp}\left\{\pi
i\left(\frac{a+d}{12c}+s(-d,c)\right)\right\}\,,
$$
and
$$
s(h,k):=\sum_{r=1}^{k-1}\frac{r}{k}\left(\frac{hr}{k}-\left[\frac{hr}{k}\right]-\frac12\right)\,.
$$
Thus $\tau\mapsto \eta^2(\tau)dz$ is almost a holomorphic section of
the tautological bundle $\calL_1$. In fact,
$(-i\epsilon(a,b,c,d)^2)^{12}=1$ so $\tau\mapsto \eta^2(\tau)dz$ is
a well-defined section on the pull-back of $\calL_1$ to a $12$--fold
branched cover of $\wwbar{M}_{1,1}$.
\begin{exm}
Label translates of the fundamental domain in \fullref{m11} by
the group elements in $PSL_2\Z$ used to translate them. Now compute
$-i\epsilon([A])^2$ for each group element and use this as a label.
Pick one translate with each of the $12$ different labels. The
collection of these translates is a fundamental domain for the
$12$--fold branched cover of $M_{1,1}$. Also label the adjoining
translates to figure out the identifications on the larger
fundamental domain.
\end{exm}
\begin{exm}\label{e12}
Show that tubular neighborhoods of the cusps in each of these $12$
labeled translates glue together to a once-punctured disk, and the
function $\eta^2(\tau)$ extends across this disk with a simple zero
at the center. (Use the group label to translate back to the small
fundamental domain without losing the $-i\epsilon([A])^2$ factor and
introduce $w=e^{2\pi i\tau}$ as a coordinate.)
\end{exm}
Adding the extra point described in \fullref{e12} gives the
$12$ fold cover together with a section of the pull-back of
$\calL_1$ to this cover. The negative of the identity matrix acts
trivially on every point of this cover. This means that the
stabilizer group of every point in the cover is $\Z_2$ so the
evaluation of the first Chern class of the pull-back bundle is
$\frac12$ of the number of zeros which is just $\frac12$. Thus since
Chern classes are multiplicative under covers
$\langle\tau_1(\text{pt})\rangle_{1,0}^{\text{pt}}=\frac{1}{24}$ as
claimed concluding our detour.

\subsection{Moduli stacks}\label{calm}
The proper structure to encode the orbifold idea in algebraic
geometry is a stack. One does not have to understand the definition
of a stack to get a feel for Gromov--Witten invariants so we do not
include the definition here. However most sources do not even define
a stack, so we have given the definition together with some
motivating examples in \fullref{app:a}. The main points to keep in
mind are that a stack adds extra structure that remembers the
stabilizer groups, and that the `universal' family is a stack.

A family of stable curves $\{\Sigma_{y}\}_{y\in Y}$ is encoded as a
map say $f$ from one space $X$ to a second space $Y$ such that for
every point $y\in Y$, $f^{-1}(y)$ can be thought of as a stable
curve. One thinks of this as the set of stable maps $f^{-1}(y)$
together with some topology linking everything together.

To define families of stable maps, we need a technical definition.
\begin{defn}
An $R$--module $M$ is {\it flat} if $0\to M\otimes_R A\to M\otimes_R
B \to M\otimes_R C\to 0$ is exact whenever $0\to A\to B\to C\to 0$
is exact. Given a morphism between analytic spaces, $f\co X\to Y$ the
ring of germs of analytic functions over a point $x_0\in X$ (denoted
by $\calO_X(x_0)$) is a module over the corresponding ring of germs
on $Y$, $\calO_Y(f(x_0))$. The map $f$ is called {\it flat} if
$\calO_X(x_0)$ is a flat $\calO_Y(f(x_0))$--module for every $x_0\in
X$.
\end{defn}\index{flat}\index{family}
Intuitively a map $f$ is flat if it has a nice fiber structure. See
Hartshorne \cite{hart} for a discussion and examples.
\begin{defn}\label{family}
A family of stable maps is a flat morphism $\pi\co \calV\to S$ together
with maps, $u\co \calV\to X$, and $\rho_k\co S\to\calV$ such that
$f\circ\rho_k=\text{id}_S$ and
$$
[u|_{\pi^{-1}(s_0)}, \pi^{-1}(s_0), \rho_1(s_0),\ldots,
\rho_n(s_0)]\in \wwbar{M}_{g,n}(X,\beta)
$$
for every $s_0\in S$. A morphism of families of stable maps is a
pair of maps, $G\co \calV\to\calV^\prime$ and $g\co S\to S^\prime$ that
intertwine the structure maps, $u$, $u^\prime$, $\rho_k$ and
$\rho_k^\prime$.
\end{defn}

A universal family of stable curves is one such that there is a
unique morphism of families from any given family into the universal
one. If one sticks to families over schemes or analytic spaces,
there is no universal family; see the discussion in Mumford
\cite{pic}, Harris--Morrison \cite{HM} or \fullref{mc}. The moduli
stack \index{moduli stack} is just a formal construction of a
universal family. The resulting construction generalizes the
category of schemes or analytic spaces.  Given a space $T$ one
constructs the contravariant functor $\underT$ that takes a scheme
$S$ to the set of all morphisms from $S$ to $T$. See Appendix
\ref{app:a} or \fullref{mc} for more information.

\begin{defn}
The moduli stack, $\wwbar\calM_{g,n}(X,\beta)$ is the functor
from the category of schemes (or analytic spaces) to the category of
sets taking a scheme $S$ to the set of equivalence classes of
families of genus $g$, $n$--marked stable curves representing $\beta$
in $H_2(X)$. (Recall the definition of a family of stable maps from
\fullref{family}).
\end{defn}
The definition of the moduli stack is very elegant and nicely avoids
questions about the topology of moduli space. However, one still has to
work to extend intersection theory to stacks; see Vistoli \cite{vistoli}.

\subsection{Deformation complexes}

Finally, we come to the main point of this section -- the local
structure of the moduli space. It is a general principle that
nondegenerate spaces or maps have the same local structure as their
linear approximations. The implicit function theorem is one version
of this principle. We are interested in a generalization to spaces
with group actions; see Bredon \cite{bredon}, Audin \cite{audin} and
Atiyah--Bott \cite{AB}. It will be convenient to consider right
$G$--spaces, ie assume that the group $G$ acts on the right in
contrast to the standard convention.

As a warm-up we shall study the following situation. Let $F\co X\to Y$ be
a smooth equivariant map of right $G$--spaces and let $y_0:=F(x_0)$
be a $G$ fixed point of $Y$. To get a local model of $F^{-1}(y_0)/G$ in a
neighborhood of $[x_0]$ we linearize the following sequence of maps,
$$
G\overset{L_{x_0}}{\lra} X\overset{F}{\lra} Y\,,
$$
where $L_{x_0}(g):=x_0g$. The linearization is
$$
T_1G\overset{TL_{x_0}}{\lra} T_{x_0}X\overset{TF}{\lra} T_{y_0}Y\,.
$$

\begin{exm}
Prove that the above sequence is a complex, that is, $TF\circ
TL_{x_0}=0$ given that $y_0$ is a fixed point. It is called the
deformation complex.
\end{exm}
\index{deformation complex} We will call the stabilizer group of
$x_0$ the automorphism group of $x_0$, $\text{Aut}(x_0):=\{g\in
G|x_0g=x_0\}$. The kernel of $TL_{x_0}$ is the zeroth cohomology of
the deformation complex and it is isomorphic to the Lie algebra of
the automorphism group, ${\mathfrak{aut}}(x_0)$. Provided $TF$ is
surjective, $F^{-1}(y_0)$ will be a manifold. The cokernel of $TF$
measures the failure of $TF$ to be surjective. This cokernel is the
second cohomology of the deformation complex and is called the
obstruction space, ${\mathfrak{ob}}(x_0)$. Assuming that the
obstruction space vanishes, $F^{-1}(y_0)$ is a manifold locally
homeomorphic to the kernel of $TF$. The quotient of this manifold by
the $G$ action can be locally identified with the first cohomology
of the deformation complex provided that the automorphism group is
trivial. This cohomology group is called the deformation space of
$x_0$ and is denoted by ${\mathfrak{def}}(x_0)$. The point is that
the exponential map applied to the orthogonal complement of the
image of $TL_{x_0}$ in the kernel of $TF$ is a local slice for the
action of $G$ on $F^{-1}(y_0)$. In other words, it intersects each
nearby $G$ orbit in exactly one point. The obstruction, deformation
and automorphism spaces glue together to form the obstruction,
deformation and automorphism bundles over the space $F^{-1}(y_0)/G$.
If the obstruction space is trivial at a point, but the automorphism
group is nontrivial, then there is a natural action of the
automorphism group on the deformation space. Furthermore,
$F^{-1}(y_0)/G$ is locally homeomorphic to the quotient of the
deformation space by the action of the automorphism group. If the
obstruction space is nontrivial, there is a map from the
deformation space to the obstruction space and the quotient
$F^{-1}(y_0)/G$ is locally homeomorphic to the quotient of the
inverse image of zero under this map by the automorphism group. This
special map is called a Kuranishi map. The Kuranishi map is
described in the following exercise for the case of trivial
automorphism group. \index{Kuranishi map}
\begin{exm}
Let $F\co H^1\oplus V\to H^2\oplus W$ be a smooth (non-linear) map
between linear spaces satisfying $F(0)=0$, $H^1=\text{ker}(T_0F)$,
and $H^2=\text{coker}(T_0F)$. Define a map $\Phi\co H^1\oplus H^2\oplus
V\to H^1\oplus H^2\oplus W$ given by
$$\Phi(x,y,z)=(x,y+F_1(x,z),F_2(x,z))\,.$$ Use the inverse function
theorem to prove that there is a locally defined inverse and smooth
maps $\psi\co H^1\to H^2$ and $\xi\co H^1\to V$ defined on open
neighborhoods of zero in $H^1$ such that $\Phi(x,\psi(x),\xi(x))=0$.
Conclude that $\psi^{-1}(0)$ is locally homeomorphic to $F^{-1}(0)$.
The map $\psi$ is the Kuranishi map.
\end{exm}
\begin{exm}
Apply these ideas to various orbit types in the quotient of $S^3$
(viewed as the unit sphere in $\C^2$) by the natural action of
$S^1\times S^1$.
\end{exm}

\subsection{Deformations of stable maps}

We will now apply these ideas to the moduli space of stable maps. A
stable map is specified by a complex structure on a surface, a
collection of marked points and a $J$--holomorphic map into a
symplectic manifold. Recall the relevant definitions from \fullref{coarse}. We will separate the deformations of a stable curve
into deformations of the marked points and complex structure and
deformations of the map with a fixed set of marked points and fixed
complex structure.

Begin by considering deformations of the map with the marked points
and complex structure fixed. By definition, a map $u\co \Sigma\to X$ is
$J$--holomorphic if $\dbar u=0$. The expression $\dbar u$ may be
applied to a tangent vector in $T_{x_0}\Sigma$ to produce a tangent
vector in $T_{u(x_0)}X$. This may be reinterpreted to say that
$\dbar u\in C^\infty(\wedge^1\Sigma\otimes u^*TX)$, ie $\dbar u$ is a
$1$--form with values in the pullback of the tangent bundle. The expression
$\dbar u$ extends to a map from the complexified tangent space of
$\Sigma$ to the complexified tangent space of $X$. It is completely
determined by an induced map from $\wedge^{0,1}\Sigma$ to $u^*TX$.
Thus we usually view $\dbar u\in C^\infty(\wedge^{0,1}\Sigma\otimes
u^*TX)$.
\begin{exm}
Let $V$ be a real vector space with almost complex structure $J$,
and let $V^\C:=V\otimes_\R\C$ be the complexification. The
(anti)holomorphic subspace is the $(-)i$--eigenspace of
$J^\C:=J\otimes 1$ acting on $V^\C$. These are denoted by
($V^{0,1}$) or $V^{1,0}$.
\begin{enumerate}
\item Show that the projection $P^{1,0}\co V^\C\to V^{1,0}$ given by
$$P(X)=\tfrac12(X-J\otimes iX)$$ restricts to an isomorphism of complex vector spaces
from $V$ to $V^{1,0}$ when $V$ is viewed as a complex vector space via $J$ and is
viewed as a subspace of $V^\C$ via $V\otimes 1$. (For this reason we often identify
the holomorphic subspace of $V^\C$ with $V$ by $P^{1,0}$.)
\item Let $\dbar u$ act on the complexified tangent space of $\Sigma$ in the natural way.
Show that it takes $\wedge^{1,0}\Sigma:=T\Sigma^{1,0}$ to $u^*TX^{0,1}$ (it takes
holomorphic vectors to antiholomorphic vectors.)
\item Conclude that $\dbar u$ is uniquely determined by the restriction of $P^{1,0}\dbar u$
to $\wedge^{0,1}\Sigma$.
\end{enumerate}
\end{exm}
Now consider a one-parameter family of maps passing through $u$, say
$v_t$. The derivative of this family (at $t=0$) associates a tangent
vector in $T_{u(x_0)}X$ to a point $x_0\in\Sigma$, so we may view
$\dot v_0\in C^\infty(u^*TX)$. The linearization of $\dbar$ maps
$\dot v_0$ to the derivative of $\dbar v_t$ at $t=0$. We will abuse
notation and denote this by $\dbar\dot v_0$. This may be viewed as
an element of $C^\infty(\wedge^{0,1}\Sigma\otimes u^*TX)$ as
explained in the previous exercise. An expression for the
linearization of $\dbar$ may be found in Salamon \cite{salamon},
McDuff--Salamon \cite{MdS} and Audin--Lafontaine \cite{AL}. We
represent it by
\begin{equation}\label{defu}
\dbar\co C^\infty(u^*TX)\to C^\infty(\wedge^{0,1}\Sigma\otimes u^*TX).
\end{equation}
The kernel of this map is the deformation space of maps with fixed
complex structure and marked points. It is denoted by
\index{$D$@$\mathfrak{def}(u)$ map deformation space}
\index{deformations} $\mathfrak{def}(u)$. The space of stable maps
with underlying marked curve $[\Sigma,p]$ is locally diffeomorphic
to ${\mathfrak{def}}(u)$ at $[u,\Sigma,p]$ provided that the
obstruction space (cokernel of the above map) ${\mathfrak{ob}}(u)$
\index{obstruction space} \index{$O$@${\mathfrak{ob}}(u)$ map
obstruction space} vanishes. There are no automorphisms that act on
maps with fixed marked domain.

We now consider automorphisms and deformations of the underlying
marked curve. An automorphism is just a holomorphic map
$\varphi\co \Sigma\to\Sigma$ fixing the marked points, so
$\dbar\varphi=0$ and $\varphi(p_k)=p_k$. Given a one-parameter
family of maps $\varphi_t$, we can differentiate the defining
conditions of the automorphism group to obtain the conditions
specifying infinitesimal automorphisms: $\dbar\dot\varphi_0=0$ and
$\dot\varphi_0|_{p_k}=0$. The nice way to encode these holomorphic
vector fields that vanish at the marked points is to use the kernel
of the following operator,
\begin{equation}\label{defs}
\dbar\co C^\infty(T\Sigma\otimes [-p])\to C^\infty(\wedge^{0,1}\Sigma\otimes T\Sigma\otimes
[-p])\,.
\end{equation}
Here $-p=-p_1-\cdots-p_n$ is negative the divisor associated to the
marked points and $[-p]$ is the complex line bundle associated with
this divisor; see Griffiths--Harris \cite{GH}. One may pick a (unique up to scale)
non-zero meromorphic section of this line bundle with a simple pole
at each marked point, say $s^-$. The map
$\dot\varphi_0\mapsto\dot\varphi_0 \otimes s^-$ identifies the
infinitesimal automorphisms ${\mathfrak{aut}}([\Sigma,p])$
\index{automorphism space} \index{$A$@${\mathfrak{aut}}([\Sigma,``p])$
curve automorphism space} with the kernel of the above map.

Since a complex structure satisfies $J^2=-I$, the derivative of a
one-parameter family of complex structures satisfies $\dot
J_0J+J\dot J_0=0$. For a holomorphic vector $X$ we have
$$
J(J\dot J_0X)=-J(\dot J_0JX)=-iJ\dot J_0X\,.
$$
Thus $\frac12J\dot J_0$ takes holomorphic vectors to antiholomorphic
vectors and may be interpreted as an element of
$\wedge^{0,1}\Sigma\otimes T\Sigma$ similar to the way that $\dbar
\dot v_0$ may be interpreted as an element of
$\wedge^{0,1}\Sigma\otimes u^*TX$. A one-parameter family of
deformations $J_t$ is trivial if there is a one-parameter family of
reparametrizations $\varphi_t$ so that $d\varphi_t\circ J=J_t\circ
d\varphi_t$.
\begin{exm}
Show that that $d\varphi_t\circ J=J_t\circ d\varphi_t$ implies that
$\frac12 J\dot J_0=\dbar\dot\varphi_0$.
\end{exm}
One also has to consider deformations of the marked points.
Surprisingly, deformations of the marked points may be modeled by
deformations of the complex structure that are fixed at the marked
points. This is easiest to understand on the Riemann sphere $\CP^1$
because $\CP^1$ has exactly one complex structure up to
reparametrizations.

Consider a collection of four or more \index{deformation of a marked
curve} points on $\CP^1$, say $p=(p_1,\ldots,p_n)$. Since we know
that there is a unique complex structure on $\CP^1$ up to
reparametrizations on $\CP^1$, given any family of complex
structures $J_t$ we may find a family of reparametrizations
$\varphi_t$ so that $d\varphi_t\circ J=J_t\circ d\varphi_t$. The
expression $\varphi_t(p)$ describes the associated one-parameter
family of marked points. This generalizes to deformations of the
complex structure and marked points on any marked curve. To
summarize, each deformation of the equivalence class of a marked
curve $[\Sigma,p]$ is uniquely specified by an element ($\frac12
J\dot J_0$) of the cokernel of the map in equation \eqref{defs}. The
obstruction space vanishes for dimensional reasons.

\begin{figure}[ht!]
\centering
\includegraphics[width=3truein]{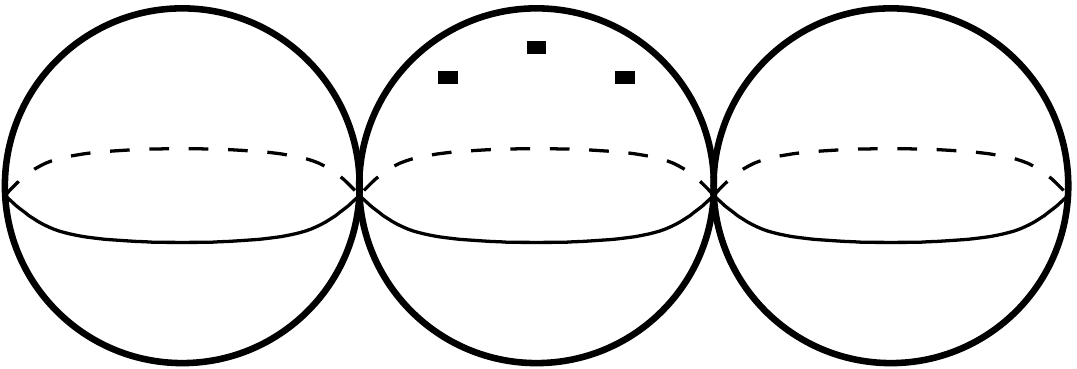}
\caption{A nodal curve with automorphisms and deformations}\label{defaut}
\end{figure}

One may think that it is impossible to have nontrivial
automorphisms and nontrivial deformations of the same marked
Riemann surface. This is true for smooth surfaces, but it is not
true for nodal curves. The marked curve in \fullref{defaut} has a
four-complex-dimensional family of automorphisms corresponding to
linear reparametrizations of each side bubble. In addition it has a
seven-complex-dimensional family of deformations, three for the
marked points in the center bubble, two more for the locations of
the nodes and two more for the resolutions of the side bubbles.

We can now assemble the deformation complexes of a map (equation
\eqref{defu}) and marked curve (equation \eqref{defs}) to obtain the
deformation complex of a stable map $[u,\Sigma,p]$. We first form
the double complex
\begin{equation}\label{double}
\bfig\barrsquare/->`<-`<-`->/<1500,500>[C^{\infty}(u^*TX)`
  C^{\infty}(\wedge^{0,1}\Sigma\otimes u^*TX)`
  C^{\infty}(T\Sigma\otimes{[-p]})`
  C^{\infty}(\wedge^{0,1}\Sigma\otimes T\Sigma\otimes{[-p]});
  \dbar`du(-\otimes s^+)`du(-\otimes s^+)`\dbar]\efig
\end{equation}
Here $s_+$ is a non-zero section of $[p]$ with a simple zero at each marked point.
A bit of thought allows one to conclude that the associated total complex,
$$
C^\infty(T\Sigma\otimes [-p]) \to
C^\infty(\wedge^{0,1}\Sigma\otimes T\Sigma\otimes [-p]) \oplus C^\infty(u^*TX) \to
C^\infty(\wedge^{0,1}\Sigma\otimes u^*TX) \,,
$$
is the deformation complex of the stable map $[u,\Sigma,p]$. The
kernel of the first map ($(\dbar, du(-\otimes s^+))$) encodes the
conditions that an infinitesimal automorphism corresponds to a
one-parameter family of holomorphic maps $\varphi_t$ fixing the
marked points and satisfying $u\circ\varphi_0=u$. The kernel of the
second map ($\dbar- du(-\otimes s^+)$) encodes the fact that each
map in a deformation of a stable map must be $J$--holomorphic with
respect to the corresponding structure ($J\circ du_t=du_t\circ
J_t$).

\subsection{Homological description of deformations}
Introducing a bit of homological algebra will allow us to formulate
the deformation complex for nodal curves, where the notion of
$C^\infty$ sections may be less clear, and will provide additional
computational tools for the study of the local structure of the
moduli space. The homological definition of the deformation complex
will also have the added benefit of defining bundles of
infinitesimal automorphisms, deformations and obstructions over the
moduli space as opposed to just defining vector spaces attached at
each point. The book by Weibel \cite{wei} is a good reference for the
homological algebra that we use. \index{injective
module}\index{injective resolution}
\begin{defn}
An $R$--module $I$ is called injective if the functor
$\Hom_R(-,I)$ is exact (that is, takes exact sequences to exact
sequences). An \emph{injective resolution} of an $R$--module $M$ is an exact
sequence,
$$
0\to M\to I_0\to I_1\to\cdots\,,
$$
where the $I_k$ are injective. Given a complex of $R$--modules,
$$A_*:= \cdots\to A_k\to A_{k+1}\to\cdots $$ and an $R$--module $B$, one may take an
injective resolution of $B$ say $0\to B\to I_*$, and form the double
complex $\Hom_R(A_i,I_j)$. The \emph{hyperext} group $\E
xt_R^k(A_*,B)$ \index{$E$@$\Ext_R^k(A_*,B)$ hyperext group} is by
definition the $k$th cohomology of the associated double complex.
\end{defn}
Notice that we have isomorphisms
$$
\begin{aligned}
C^\infty(u^*TX)\cong & \Hom_{\calO_\Sigma}(u^*\Omega_X, C^\infty(\wedge^0\Sigma)) \\
C^\infty(\wedge^{0,1}\Sigma \otimes u^*TX)\cong &
\Hom_{\calO_\Sigma}(u^*\Omega_X,C^\infty(\wedge^{0,1}\Sigma))
\\
C^\infty(T\Sigma\otimes [-p])\cong & \Hom_{\calO_\Sigma}(\Omega_\Sigma([p]),C^\infty(\wedge^0\Sigma)) \\
C^\infty(\wedge^{0,1}\Sigma\otimes T\Sigma\otimes [-p])\cong &
\Hom_{\calO_\Sigma}(\Omega_\Sigma([p]),C^\infty(\wedge^{0,1}\Sigma))\,.
\end{aligned}
$$
Here $\calO_\Sigma$ \index{$O$@$\calO_\Sigma$ sheaf of holomorphic
functions} is the sheaf of holomorphic functions on $\Sigma$, $\Omega_\Sigma(L)$ is the
sheaf of holomorphic differentials on $\Sigma$ taking values in a
line bundle $L$ \index{$O$@$\Omega_\Sigma(L)$ sheaf of holomorphic differentials
taking values in $L$} and \index{$O$@$\Omega_X$ sheaf of holomorphic
differentials} $\Omega_X$ is the sheaf of holomorphic differentials
on $X$. Any sheaf ${\mathcal E}$ may be realized as the sections of
an associated sheaf space $E\to X$. By definition, the pull-back of
a sheaf \index{pull-back sheaf} under a map $f\co Y\to X$ is
$f^*{\mathcal E}:=\Gamma(f^*E)$.
\begin{exm}
The third isomorphism listed above is defined by $$F(X\otimes
s)(\alpha\otimes t):=\alpha(X)st\,.$$ Find the other three
isomorphisms.
\end{exm}

We can now see how to represent the infinitesimal automorphisms,
deformations and obstructions as hyperext groups. Notice that a
partition of unity argument may be used to show that a sheaf of
$C^\infty$ sections is fine. This in turn implies that it is
$\Hom_{\calO_\Sigma}(du,-)$--acyclic thus the hyperext groups
can be computed with these sheaves; see Hartshorne \cite{hart}. Also recall
that the definition of an exact sequence of sheaves is defined by
the requirement that the sequence be exact at the level of germs.
This means that the usual de Rham complex is an exact sequence of
sheaves (even though the sequence of global sections fails to be
exact as measured by the de Rham cohomology). The point is that a
closed form becomes exact when it is restricted to a small enough
open set. It follows that the following is an acyclic resolution of
$\calO_\Sigma$:
$$
0\to \calO_\Sigma\to C^\infty(\wedge^0\Sigma)\to C^\infty(\wedge^{0,1}\Sigma)\to 0\,.
$$
We will use this resolution and the definition of the hyperext
groups to construct the deformation-obstruction complex.

\subsection{The deformation-obstruction sequence}\label{defobsec}
Using the four isomorphisms and the definition of the hyperext
groups we arrive at the following definitions which generalize our
proceeding discussion for smooth stable maps.
\begin{defn}
The spaces in the deformation-obstruction sequence of a stable map
are given by: \index{$D$@$\mathfrak{def}(u)$ map deformation
space}\index{$O$@${\mathfrak{ob}}(u)$ map obstruction
space}\index{$A$@${\mathfrak{aut}}([\Sigma,``p])$ curve automorphism
space} \index{$D$@${\mathfrak{def}}([\Sigma,p])$ curve deformation
space}\index{$A$@${\mathfrak{aut}}([u,\Sigma,p])$ automorphism
space}\index{$D$@${\mathfrak{def}}([u,\Sigma,p])$ deformation
space}\index{$O$@${\mathfrak{ob}}([u,\Sigma,p])$ obstruction space}
\begin{align*}
{\mathfrak{def}}(u):= & \Ext_{\calO_\Sigma}^0(u^*\Omega_X, \calO_\Sigma) \\
{\mathfrak{ob}}(u):= & \Ext_{\calO_\Sigma}^1(u^*\Omega_X, \calO_\Sigma) \\
{\mathfrak{aut}}([\Sigma,p]):= & \Ext_{\calO_\Sigma}^0(\Omega_\Sigma([p]),\calO_\Sigma) \\
{\mathfrak{def}}([\Sigma,p]):= & \Ext_{\calO_\Sigma}^1(\Omega_\Sigma([p]),\calO_\Sigma) \\
{\mathfrak{aut}}([u,\Sigma,p]):= & \Ext_{\calO_\Sigma}^0(u^*\Omega_X\to\Omega_\Sigma([p]),\calO_\Sigma) \\
{\mathfrak{def}}([u,\Sigma,p]):= & \Ext_{\calO_\Sigma}^1(u^*\Omega_X\to\Omega_\Sigma([p]),\calO_\Sigma) \\
{\mathfrak{ob}}([u,\Sigma,p]):= & \E
xt_{\calO_\Sigma}^2(u^*\Omega_X\to\Omega_\Sigma([p]),\calO_\Sigma)\,.
\end{align*}
\end{defn}

Homological algebra will lead to a long exact sequence relating
these groups. The following aside sketches the derivation of the
exact sequence.
\begin{aside}{\small
Since the hyperext groups are defined as the cohomology of the total
complex of a double complex and every double complex has an
associated spectral sequence, there is a spectral sequence related
to the hyperext groups. As the double complex only has two terms in
any direction, the spectral sequence will collapse at the $E_2$ term
and produce a long exact sequence. See Bott and Tu \cite{BT}, Brown
\cite{Brown} and Spanier \cite{span} for a
detailed discussion of spectral sequences. Briefly, to any double
complex $C^{p,q}$ one associates a total complex
$$TC^n:=\oplus_{k+\ell=n}C^{k,\ell}$$
filtered by
$$F^pTC^n:=\oplus_{k+\ell=n, k\ge p}C^{k,\ell}.$$
The zeroth page of the spectral sequence is defined by
$$E_0^{p,q}:=\frac{F^pTC^{p+q}}{F^{p+1}TC^{p+q}}.$$
Similarly the last page of the spectral sequence is defined by
$$E_\infty^{p,q}:=\frac{F^pH^{p+q}(TC^*)}{F^{p+1}H^{p+q}(TC^*)}.$$

\begin{exm}
Prove that for the spectral sequence of a double complex one has $E_0^{p,q}=C^{p,q}$.
\end{exm}
The differentials of the double complex may be used to define
differentials on the pages of the spectral sequence,
$d_n\co E_n^{p,q}\to E_n^{p+n,q+1-n}$. Here the pages may be
inductively defined by $E_{n+1}:=H(E_n,d_n)$. Applying this to the
double complex of equation \eqref{double} we see that the $E_1$--page
of the associated spectral sequence is given by taking the
cohomology in the vertical direction. Combined with the definition
of the deformation spaces this gives:
$$
E_1=\begin{array}{cc}
{\mathfrak{def}}([\Sigma,p])\quad & {\mathfrak{ob}}(u) \\
{\mathfrak{aut}}([\Sigma,p]) \quad & {\mathfrak{def}}(u).
\end{array}
$$
Taking cohomology in the horizontal direction gives:
$$
E_2=\begin{array}{cc} \text{ker}({\mathfrak{def}}([\Sigma,p])\to
{\mathfrak{ob}}(u)) \quad &
\text{coker}({\mathfrak{def}}([\Sigma,p])\to {\mathfrak{ob}}(u))  \\
\text{ker}({\mathfrak{aut}}([\Sigma,p])\to {\mathfrak{def}}(u))
\quad & \text{coker}({\mathfrak{aut}}([\Sigma,p])\to
{\mathfrak{def}}(u)) \,.
\end{array}
$$
Since all groups outside of this square are zero all other
differentials are zero and we conclude that $E_\infty=E_2$. Using
the definition of the $E_\infty$ page and the definition of the
deformation groups we conclude that
$$E_\infty^{0,1}=\frac{F^0H^{1}(TC^*)}{F^{1}H^{0 }(TC^*)},\qua
E_\infty^{1,0}= F^0H^{1}(TC^*) \qua\text{and}\qua
F^0H^{1}(TC^*)={\mathfrak{def}}([u,\Sigma,p]).$$
Combining this with
the above computation of the $E_\infty$ page gives the following
exact sequence:
$$
0\to \text{coker}({\mathfrak{aut}}([\Sigma,p])\to
{\mathfrak{def}}(u))\to \text{ker}({\mathfrak{def}}([\Sigma,p])\to
{\mathfrak{ob}}(u))\to 0\,.
$$

\begin{exm}
Continue in this way to prove that the following
deformation-obstruction sequence is exact.
\end{exm}}
\end{aside}
The deformation-obstruction sequence is
\begin{multline}\label{defob}
0\to {\mathfrak{aut}}([u,\Sigma,p]) \to
  {\mathfrak{aut}}([\Sigma,p]) \\
\to {\mathfrak{def}}(u) \to {\mathfrak{def}}([u,\Sigma,p]) \to
  {\mathfrak{def}}([\Sigma,p])\\
\to {\mathfrak{ob}}(u) \to {\mathfrak{ob}}([u,\Sigma,p]) \to 0\,.
\end{multline}
It is worth pointing out what the maps in the
deformation-obstruction sequence are. Any automorphism of
$[u,\Sigma,p]$ is an automorphism of $[\Sigma,p]$, so the first map
is just the inclusion. The second map is given by
$\dot\varphi_0\mapsto \frac{d}{dt}u\circ\varphi|_{t=0}$. The third
map is $\dot v\mapsto (\dot v, 0)$; the fourth is $(\dot v,B)\mapsto
B$; the fifth is $B\mapsto du(B\otimes s^+)$. The sixth map is just
the projection because both ${\mathfrak{ob}}(u)$ and
${\mathfrak{ob}}([u,\Sigma,p])$ are quotients of the same group, but
the latter is a quotient by a larger equivalence. A few more
comments will clarify these notions.

First consider the infinitesimal automorphisms and deformations of a
smooth marked curve. We have the following acyclic resolution of
$\calO_\Sigma(T\Sigma\otimes [-p])$:
$$
0\to\calO_\Sigma(T\Sigma\otimes [-p])\to C^\infty(T\Sigma\otimes [-p])\to
C^\infty(\wedge^{0,1}\Sigma\otimes T\Sigma\otimes [-p])\to 0\,.
$$
It follows that we may make the identifications
$$
\begin{aligned}
{\mathfrak{aut}}([\Sigma,p])\cong & H^0(\Sigma,\calO_\Sigma(T\Sigma\otimes [-p])) \\
{\mathfrak{def}}([\Sigma,p])\cong &
H^1(\Sigma,\calO_\Sigma(T\Sigma\otimes [-p]))\,.
\end{aligned}
$$
The simplest version of the Riemann--Roch theorem (see Forster
\cite{forster}) states that
$$
\text{dim}_\C H^0(\Sigma;L)-\text{dim}_\C
H^1(\Sigma;L)=c_1(L)[\Sigma]+1-g\,.
$$
This gives
\begin{align*}
\dim_\C &{\mathfrak{aut}}([\Sigma,p]) - \text{dim}_\C
  {\mathfrak{def}}([\Sigma,p]) \\
&= \dim_\C H^0(\Sigma;\calO_\Sigma(T\Sigma\otimes [-p])) -
  \dim_\C H^1(\Sigma;\calO_\Sigma(T\Sigma\otimes [-p])) \\
&= \deg(T\Sigma\otimes [-p])+1-g = 3-3g-n.
\end{align*}
This is in fact true for the deformations and infinitesimal
automorphisms of any marked curve. For the special case of
$\Sigma=\CP^1$ the bundle $T\CP^1\otimes [-p]$ is just the degree
$2-n$ bundle over the sphere. The \emph{Kodaira vanishing theorem}
(see Griffiths--Harris \cite{GH}) states that $H^q(X;\Omega^p(L))=0$ for $p+q>n$ when $L$
has positive degree. With $3$ or fewer marked points the Kodaira
vanishing theorem implies that
$H^1(\Sigma,\calO_\Sigma(T\Sigma\otimes [-p]))=0$ and we conclude
that $\text{dim}_\C\, {\mathfrak{aut}}([\Sigma,p])=3-n$ and
$\text{dim}_\C\, {\mathfrak{def}}([\Sigma,p])=0$. One may also
directly compute $H^0$ as homogenous degree $2-n$ polynomials in
this case. One sees that this matches perfectly with linear
fractional transformations fixing $3$ or fewer points.

With three or fewer marked points $H^1$ is trivial and $H^0$ can be
nontrivial. We will see that the situation with $4$ or more marked
points is just the opposite--$H^0$ is trivial and $H^1$ can be
nontrivial. To be more precise, recall that for smooth varieties
there is a non-degenerate pairing,
$$
H^{0,k}(X,\calO(E))\times H^{0,n-k}(X,\calO(E^*\otimes
\wedge^{n,0}X))\to \C\,,
$$
given by $(\alpha, \beta)\mapsto \int_X \alpha\wedge\beta$. Here the
$E$ and $E^*$ components pair to give a number and we have $k$
components $d\barz_i$ from the first term, $n-k$ from the second
term and $n$ components $dz_i$ from the coefficients in the second
term. This pairing gives the duality,
$$
H^{0,k}(X,\calA)\cong (H^{0,n-k}(X,\Omega_X\otimes\calA))^*\,,
$$
known as \emph{Kodaira--Serre duality} \cite{GH}. Here $\Omega_X$ is the
sheaf of top-dimensional holomorphic forms. \index{$O$@$\Omega_X$
top-dimensional holomorphic forms}(The bundle $\wedge^{n,0}X$ is
called the canonical bundle and is denoted by $K_X$. Applying this
duality to the deformations gives,
$$
{\mathfrak{def}}([\Sigma,p])\cong
H^0(\Sigma,\Omega_\Sigma\otimes\Omega_\Sigma([p]))^*\,.
$$
The elements of this last group take the form $f(z)\,dz\otimes dz$
in local coordinates and are called quadratic differentials.

\index{Kodaira--Serre duality} Kodaira--Serre duality generalizes to
more general projective varieties where there is a dualizing sheaf,
denoted by $\omega_X$, so that
$$
H^{0,k}(X,\calA)\cong (H^{0,n-k}(X,\omega_X\otimes\calA))^*\,. $$
See Hartshorne \cite{hart}. The dualizing sheaf leads to the last family of
cohomology classes that we will need. One can construct a bundle
over the moduli space, such that the fiber over any point
$[u,\Sigma,p]$ is just the space of sections of the dualizing sheaf
$H^0(\Sigma,\omega_\Sigma)$. This bundle is called the Hodge bundle
\index{$E$@$\E$ Hodge bundle}\index{dualizing sheaf} and it is
denoted by $\E$. The Chern classes of the Hodge bundle are called
\index{Hodge classes} Hodge classes. The formal definition of the
\index{$\omega_\Sigma$ dualizing sheaf} dualizing sheaf is as
follows (see Hori et al \cite{Hori}):
\begin{defn}
The \emph{dualizing sheaf} over a nodal curve $\Sigma$ denoted by $\omega_\Sigma$ is
the sheaf of meromorphic differentials that:
\begin{enumerate}
\item are holomorphic away from the nodes
\item have at worst a pole of order one at each node branch
\item have residues that sum to zero on each pair of node branches.
\end{enumerate}
\end{defn}
\begin{exm}
Give a formal definition of the Hodge bundle analogous to the
definition of the tautological bundles $\calL_k$.
\end{exm}

According to the previous application of the Riemann--Roch theorem,
for every extra marked point the difference between the dimensions
of the automorphisms and deformations decreases by one. Returning to
the example in \fullref{defaut}, we see that adding one point to
the left bubble would reduce the dimension of the automorphism group
by one because the restriction of the automorphism to the left
bubble would have to fix two points, not just one as before. However,
this would not change the space of deformations. One might think
that a possible position of the new marked point needs to be taken into account.
But any motion of the new point is a trivial deformation. Indeed, an
automorphism changing its position can can be applied to the space before
adding the point. A similar thing occurs with the addition of a second
marked point to the left bubble. When one adds a third marked point
to the left bubble, the dimension of the group of automorphisms will
not change from the dimension with two marked points because there
are no automorphisms acting nontrivially on the left bubble (such
an automorphism would have to fix at least three points). There
would be an extra deformation corresponding to changing the location
of the third point.

We now have a fairly good local description of the moduli stack.
Since the moduli space models stable curves there are no
infinitesimal automorphisms. Assuming that the obstruction space
vanishes, we see that the moduli space is locally isomorphic to the
quotient of ${\mathfrak{def}}([u,\Sigma,p])$ by a finite group. The
deformation-obstruction sequence can be combined with the
Riemann--Roch theorem to compute the dimension of
${\mathfrak{def}}([u,\Sigma,p])$. This is called the \index{virtual
dimension} virtual dimension of the moduli space.
\begin{exm}\label{virdim}
Recall that the alternating sum of the dimensions of an exact
sequence of spaces is zero. Use the Riemann--Roch theorem
$$\sum_{k=0}^n(-1)^k\text{dim}_\C\,H^k(X,E)=\text{Td}(TX)\text{ch}(E)[X]
$$
to compute $\text{dim}_\C \,{\mathfrak{def}}(u)  - \text{dim}_\C
\,{\mathfrak{ob}}(u)$. Combine this with our earlier computation of
$$\text{dim}_\C \,{\mathfrak{aut}}([\Sigma,p])
- \text{dim}_\C \,{\mathfrak{def}}([\Sigma,p])$$
 and the deformation-obstruction sequence to conclude
$$
\virdim_\C \,\wwbar{M}_{g,n}(X,\beta)=\int_\beta
c_1(TX)+(\text{dim}_\C \,{X}-3)(1-g)+n\,.
$$
\index{$V$@$\virdim_\C$ virtual dimension}
\end{exm}

\section{Localization}\label{locGW}

We have come a long way in our review of Gromov--Witten invariants.
We described all the technical elements in their definition with the exception
of the virtual fundamental class (addressed in a later subsection). We also performed a number of non-trivial computations via recursion or direct reasoning. In this section we will describe a new computational tool called
localization used to compute the Gromov--Witten invariants of the resolved conifold $X_{S^3}$. We start with a
general description of localization and an outline of the virtual
localization formula for Gromov--Witten invariants. Some readers may
prefer to skip down to the sample computations that we give for
$\CP^2$ after the general discussion.

\subsection{The Umkehrung}

Localization is a technique reducing a computation of an integral
over a higher-dimensional space to an integral over a
lower-dimensional space. Of course this is impossible in general,
but it is instructive to try.

Given an inclusion (or any map) $\iota_F\co F\hookrightarrow M$, we
have the well-known pull-back $\iota_F^*\co H^k(M)\to H^k(F)$, given on
the level of Poincar\'e duals by the inverse image. There is a less
well-known push-forward $\iota_{F!}\co H^k(F)\to H^{k+m-f}(M)$ defined
for any map between oriented manifolds. This map (called the
Umkehrung) is defined by the following diagram.
$$\bfig\barrsquare/->`->`<-`->/<1000,500>[H^k(F)`H^{k+m-f}(M)`
  H_{f-k}(F)`H_{f-k}(M);
  \iota_{F!}`PD`PD^{-1}`\iota_{F*}]\efig$$
Here $m$ is the dimension of $M$, $f$ is the dimension of $F$ and
$PD$ is Poincar\'e duality. We will apply these ideas to orbifolds.
For orbifolds everything goes through as in the manifold case
provided one uses rational coefficients.

There are nice descriptions of the Umkehrung \index{Umkehrung} for
fibrations and embeddings. As an example, if $\pi_M\co M\to\text{pt}$,
the map $\pi_{M!}\co H^m(M)\to H^0(\text{pt})$ is just given by
integration $f_!\alpha=\int_M\alpha$. When $\pi\co E\to M$ is a fiber
bundle, the Umkehrung is just integration over the fiber (see Bott and Tu
\cite{BT}).

In the case of an embedding $\iota_F\co F\hookrightarrow M$, we can
compute $\iota_{F!}1$ and get an interesting answer. The cycle dual
to $1$ in $F$ is just $F$ and so the cycle dual to $\iota_{F!}1$ is
just the image of $F$ in $M$ which is the zero section of the normal
bundle to $F$ in $M$. Recall that the image of the zero section of a
bundle is dual to the \emph{Thom class} \index{Thom class}.
The \emph{Euler class} \index{Euler class} of the bundle is the pull-back
under any section of the Thom class \cite{BT}. It follows that
$$
\iota_F^*\iota_{F!}1=e(N(F))\,.
$$

\begin{remark}
The Umkehrung for an \index{$F$@$F_{{\pling}}$ Umkehrung} arbitrary map
can be decomposed into one for an embedding and one for a
projection. Starting with a cycle in $F$ dual to the given
cohomology class one obtains a cycle in $F\times M$ by the natural
inclusion $F\hookrightarrow F\times M$ taking $x$ to
$(x,\iota_F(x))$. Given this cycle in $F\times M$ one obtains a
cycle in $M$ by composition with the projection. The class dual to
this final cycle is the value of the Umkehrung.
\end{remark}

If one could invert the Euler class, one might hope that the maps
$\frac{\iota_F^*}{e(N(F))}\co H^*(M)\to H^*(F)$ and
$\iota_{F!}\co H^*(F)\to H^*(M)$ would be inverses, thus reducing
integrals over $M$ to integrals over $F$. Of course this is too much
to expect in general. Unfortunately, cohomology classes of positive
degree on a finite-dimensional manifold are all nilpotent (since high
enough powers would be forms of degree larger than the dimension of
the manifold) and therefore not invertible.

\subsection{Equivariant cohomology}

It is much more reasonable to expect such a reduction to work if
everything is invariant under a group action because as we will see,
equivariant cohomology has more invertible elements. The insight of
Atiyah and Bott \cite{AB} was that this could be made to work when $F$ is the
fixed point locus of a torus action on $M$. The geometric
intuition behind the reduction of an integral over a larger set to
an integral over a smaller set for equivariant functions is that the
symmetry allows one to  sample their values at a smaller collection
of points.

We now need a brief review of \index{equivariant cohomology}
equivariant cohomology. If $G$ is any Lie group and $EG$ \index{$EG$
classifying bundle} \index{classifying bundle} is a contractible,
free, right, proper $G$--space and $M$ is any left $G$--space one can
form the \index{$EG\times_G M$ twisted product}twisted product
\index{twisted product} $EG\times_G M:=EG\times M/\sim$ where
$(eg^{-1},gm)\sim (e,m)$. The equivariant cohomology of $M$ is
defined to be
$$
H^*_G(M):=H^*(EG\times_G M)\,.
$$
\index{$H^*_G(M)$ equivariant cohomology} As an example take
$G=T^2$, then $EG=S^\infty\times S^\infty$ where $S^\infty$ is
viewed as the unit sphere is an infinite-dimensional complex space
and the action is given by
$(x_0,x_1)\cdot(\lambda_0,\lambda_1):=(x_0\lambda_0,x_1\lambda_1)$.
The equivariant cohomology of a point is computed as follows \cite{AB},
$$
H^*_{T^2}(\text{pt};\Q)=H^*(\CP^\infty\times\CP^\infty;\Q)\cong
\Q[\alpha_0,\alpha_1]\,.
$$
Notice that $\CP^\infty$ is infinite dimensional and $\alpha_k$ have
degree $2$. Given a representation $\mu\co G\to GL_n\C$ we get a vector
bundle $EG\times_G \C^n\to EG\times_G\text{pt}$, and associated
equivariant classes given by the Chern classes of this bundle
\index{$\mu_k\co T^n\to GL_1\C$ projection to $k$th factor}
$c_k(EG\times_G\C^n)$. If $\mu_k\co T^n\to GL_1\C$ is the projection to the
$k$th factor and $\mu_k^*$ is the dual representation then the first
Chern class of the line bundle associated to $\mu_k^*$ is the
standard generator $\alpha_k$ \index{$\alpha_k$ standard generator
of $H^*_{T^2}(\text{pt};\Q)$} of $H^*_{T^2}(\text{pt};\Q)$. This
corresponds to the fact that the first Chern class of the
tautological line bundle over projective space evaluates to $-1$ on
a standardly oriented generator of the second homology. In fact
homogeneous polynomials of degree $d$ may be regarded as sections of
the degree $d$ (as measured by the first Chern class) line bundle
over projective space and this is where all of the sign conventions
come from. \index{$T=T^{n+1}$ the $(n{+}1)$--torus}

A second example of equivariant cohomology is the equivariant
cohomology of the group with the natural left action. The result is
$$
H_G^*(G):=H^*(EG\times_GG)=H^*(EG)=\Q\,.
$$
Thus when the action was free the equivariant cohomology was trivial
and when the action was trivial the equivariant cohomology was
interesting. In some sense equivariant cohomology is generated by
the fixed point set of the action.

From here on forward we will work with torus actions only and $T$
will always denote a torus. In order to make the observation that
the equivariant cohomology is generated by the fixed point set more
precise and follow the reduction outline from the beginning of this
section we need to invert elements of the equivariant cohomology and
ultimately invert the Euler class of the normal bundle to the fixed
point locus. Notice that
$H^*_{T^{n+1}}(\text{pt})\cong\Q[\alpha_0,\ldots,\alpha_n]$ is an
integral domain. This is in stark contrast to the cohomology of
finite dimensional manifolds. We let
$$F^*_{T^{n+1}}\cong\Q(\alpha_0,\ldots,\alpha_n)
$$ be the associated fraction field. The obvious map $EG\times_GM\to EG\times_G\text{pt}$
induces a map $H^*_G(\text{pt})\to H^*_G(M)$ giving $H^*_G(M)$ the
structure of an $H^*_G(\text{pt})$--module. This leads us to the
first version of the theorem of Atiyah and Bott.
\begin{thm}
Let Fix be the fixed point locus of a torus $T$ action on $M$. The map
$ET\times_T\text{Fix}\to ET\times_TM$ induces an isomorphism
$$
\iota^*_{\text{Fix}}\co H^*_T(M)\otimes_{H^*_T(\text{pt})}F^*_T \to
H^*_T(\text{Fix})\otimes_{H^*_T(\text{pt})}F^*_T\,.
$$
\end{thm}
Atiyah and Bott actually state a more refined theorem that specifies
which elements need to be inverted in order to obtain an
isomorphism. The proof of this theorem is to use the formula $
\iota_F^*\iota_{F!}=e(N(F)) $ to see that the maps
$\iota_{\text{Fix}!}$ and $Q:=\sum_F\frac{\iota^*_F}{e(N(F))}$ are
inverses of each other. Here $F$ represents a component of Fix and
the sum is taken over all such components.

In order to apply these ideas to the integration of ordinary
cohomology classes on $M$, notice that the map $j\co M\to
ET\times_TM$ taking $x$ to the equivalence class of $(e_0,x)$
induces a map $H^*_T(M)\to H^*(M)$ via pull-back. A class in the
image of this map is called an equivariant class. One standard way
to construct equivariant classes is to start with an equivariant
vector bundle $E\to M$ and take characteristic classes. An
equivariant bundle is just a vector bundle together with a torus
action compatible with the bundle structure. To any such bundle one
can associate the pull-back of the induced bundle $ET\times_TE\to
ET\times_TM$. Any characteristic class of $ET\times_TE$ is then an
equivariant class. Given any equivariant class $\phi\in H^m(M)$ with
equivariant lift $\hat\phi\in H^m_T(M)$ we have
\begin{multline}\label{loc}
\int_M\phi= \pi_{M!}\phi = \pi_{M!}j^*\hat\phi =
\iota^*_{\text{pt}} \pi_{M!}\hat\phi = \iota^*_{\text{pt}}
\pi_{M!}\iota_{\text{Fix}!}\sum_F\frac{\iota^*_F\hat\phi}{e(N(F))} \\
= \iota^*_{\text{pt}}
\pi_{\text{Fix}!}\sum_F\frac{\iota^*_F\hat\phi}{e(N(F))} =
\sum_F\int_F\frac{\iota^*_F\hat\phi}{e(N(F))}.
\end{multline}
This is the standard way of thinking about the Atiyah--Bott localization
formula.

\subsection{Equivariant cohomology of $\CP^n$}

To use the localization formula \eqref{loc} to integrate a
cohomology class $\phi$ one must find an equivariant lift
$\what\phi$ that maps to $\phi$. The standard way to do this is
to express $\phi$ as a product of characteristic classes of vector
bundles and then extend the group action over the vector bundles.
Such an extension is called a \index{linearization of an
action/bundle} \emph{linearization of the action}. For example, we can
represent the hyperplane bundle over $\CP^n$ by
$(\C^{n+1}-\{0\})\times\C/\sim$ with
$(z_0,\ldots,z_n,\xi)\sim(wz_0,\ldots,wz_n,w\xi)$. Define a family
of linearizations of this bundle by
$\lambda\cdot(z_0,\ldots,z_n,\xi):=(\lambda_0z_0,\ldots,\lambda_nz_n,\lambda^{q_0}_0\ldots\lambda_n^{q_n}\xi)$.
Denoting the hyperplane bundle with this linearization by
$L_{q_0,\ldots, q_n}$, we have
$$
c_1(L_{q_0,\ldots, q_n})=h-q_0\alpha_0-\cdots-q_n\alpha_n\,.
$$
This serves to define equivariant cohomology classes of $\CP^n$. The
class $h$ \index{$H$@$h$ generator of $H_T^*(\CP^n)$ as an
$H_T^*$--module} is defined to be the first Chern class of the bundle
$L_{0,\ldots,0}$, and $\alpha_k$ is defined as a difference of Chern
classes. Let us outline a proof that
$$
H_T^*(\CP^n)=\Q[h,\alpha_0,\ldots,\alpha_n]/\left(\prod_{k=0}^n(h-\alpha_k)\right)\,.
$$

First consider a finite-dimensional model $X_{p,n,N}$ for $ET\times_T\CP^p$.
If $T=T^{n+1}$ acts on $\CP^p$ in the usual way then
\begin{multline*}
X_{p,n,N}:=\\
\bigl\{([x_0],\ldots,[x_n];[y_0:\ldots:y_p])\in(\CP^N)^{n+1}
  \times\CP^{(p+1)(N+1)-1}~\big|~\wedge^2
  \bigl[\begin{smallmatrix}x_k\\y_k\end{smallmatrix}\bigr]=0\bigr\}.
\end{multline*}
\index{$X_{p,n,N}$ a finite-dimensional model of $ET\times_T\CP^p$}
A specific case helps explain what this means. Take $n=2$, $p=2$,
and $N=4$. We will write $[x_0]=[x_0:x_1:x_2:x_3:x_4]$,
$[x_1]=[y_0:\ldots:y_4]$, $[x_2]=[z_0:\ldots:z_4]$ and
$$[y_0:\ldots:y_2]=[p_0:p_1:p_2:p_3:p_4:q_0:q_1:q_2:q_3:q_4:r_0:r_1:r_2:r_3:r_4]$$
The condition on the second wedge states that all $2\times 2$
determinants pairing $x$'s and $p$'s etc are zero, so
$$
\left|\begin{matrix} x_0 & x_3\\p_0 & p_3
\end{matrix}\right|=0,\qquad
\left|\begin{matrix} z_2 & z_3\\r_2 & r_3
\end{matrix}\right|=0,\ldots\,.
$$
To see where this comes from note that there is a natural action of
$S^1$ on $S^{2N+1}$ with quotient $\CP^N$. In the $N\to\infty$ limit
we see that $S^\infty$ is contractible, so $ES^1=S^\infty$ and
$BS^1=\CP^\infty$. Thus $(\CP^N)^{n+1}$ is a finite-dimensional
model for $BT$. The condition on the vanishing of the second
exterior powers implies that $y_0$ is proportional to $x_0$, etc.
This means that the inverse image of a point under the natural
projection $X_{p,n,N}\to (\CP^N)^{n+1}$ is a copy of $\CP^p$.
The class $h$ is Poincar\'e dual to
$$
H:=\{([x_0],\ldots,[x_n];[y_0:\ldots:y_p])\in X_{p,n,N}|\langle
x_p,y_p\rangle=0\}\,,
$$
and the class $\alpha_k$ is Poincar\'e dual to
$$
A_k:=\{([x_0],\ldots,[x_n];[y_0:\ldots:y_p])\in
X_{p,n,N}|(x_k)_N=0\}\,.
$$
These formulas can serve as alternative definitions of these classes.
The following exercise will prove that
$$
H_T^*(X_{p,n,N})=\Q[h,\alpha_0,\ldots,\alpha_n]/\left(\prod_{k=0}^n(h-\alpha_k,
\alpha_0^{N+1},\ldots,\alpha_n^{N+1})\right)\,.
$$

\begin{exm}
Prove that the cohomology rings of $X_{p,n,N}$ and $\CP^n$ take the
stated form. Note for example that $X_{2,2,4}$ is a
$28$--real-dimensional space with one $28$--cell, no $27$--cells and
$\smash{X_{2,2,4}^{(26)}}=X_{1,2,4}\cup\smash{\bigcup_{k=0}^4}A_k$.
In fact, by properly taking complements of unions of intersections
of the $A_k$ and $H$ cycles, $X_{p,n,N}$ may be decomposed into a
union of cells. It follows that the cohomology group of $X_{p,n,N}$
is as stated. The ring structure follows from the combinatorics of
the intersections of the  $A_k$ and $H$ cycles.
\end{exm}

\begin{exm}
Let $q_k$ be the point in $\CP^m$ with all coordinates other than
$z_k$ equal to zero. These points are fixed by the standard
$T$--action, so there is an equivariant class $\phi_k$ Poincar\'e
dual to $q_k$. By considering appropriate intersections of the $A_k$
and $H$ prove that \index{$\phi_k$ Poincar\'e dual to $q_k$}
$$
\phi_k=\prod_{j\ne k}(h-\alpha_j)\,.
$$
\end{exm}

Let $\iota_k\co \text{pt}\to\CP^m$ be the map with image $q_k$. We
compute
$$
\iota_k^*h=\int_{\CP^m}PD(q_k)h=\int_{\CP^m}PD(q_k)(h-\alpha_k)+
\int_{\CP^m}PD(q_k) \alpha_k=\alpha_k\,.
$$
The above computation may look a bit weird. It appears that we are
restricting a $2$--form to a point and getting a non-zero answer. The
thing to remember here is that we are working equivariantly so we
replace $q_k$ by $ET\times_Tq_k$ and $\CP^m$ by $ET\times_T\CP^m$,
and we have tensored with $\Q(\alpha_0,\ldots,\alpha_m)$.

\section{Localization computations of Gromov--Witten invariants}
As with many parts of this theory Kontsevich was the first to apply
localization to Gromov--Witten invariants \cite{konloc}. We are
following the exposition from Hori et al \cite{Hori} and Cox and
Katz \cite{CK}. To see how this applies to Gromov--Witten
invariants, notice that a $T$ action on $X$ will induce a $T$ action
on $\wwbar\calM_{g,n}(X,\beta)$ via composition
$\lambda\cdot[u,\Sigma,p]:=[\lambda u,\Sigma,p]$. In order to better
demonstrate how localization may be used to compute Gromov--Witten
invariants we will recompute some of the Gromov--Witten invariants
of $\CP^2$. This is not a very efficient way to compute these
invariants but it will allow us to explain the important points.

We should point out that thus far our discussion of localization
only applies to well-behaved spaces with torus actions. In the
Gromov--Witten setting this will only occur when the automorphism
group of every stable map in the relevant moduli space is trivial
and the obstruction space over each stable curve is zero as well,
that is, $\text{Aut}([u,\Sigma,p])=1$ and
$\mathfrak{ob}([u,\sigma,p])=0$. When the obstruction spaces are
still trivial but the automorphism groups are not, one must
generalize the localization formula to the orbifold setting. When
the obstruction spaces do not vanish, one must define virtual
fundamental cycles and prove that localization works for these
virtual cycles. In this section we will ignore the more technical points and
discuss localization in the context of Gromov--Witten invariants as
though everything were smooth and automorphism-free. Take heart, we will give an
introduction to virtual fundamental cycles later in this section and
(virtual) localization really does work in this context as was shown by Graber
and Pandharipande \cite{virloc}.

We have a standard $T^3$ action on $\CP^2$ given by
$$(\lambda_0,\lambda_1,\lambda_2)\cdot[z_0:z_1:z_2]:=
  [\lambda_0z_0:\lambda_1z_1:\lambda_2z_2].$$
The first step is to find the fixed point set of the $T$ action on
$\wwbar\calM_{g,n}(\CP^2,d)$. The subtle point is that if one
looks for maps as opposed to equivalence classes of maps fixed by
$T$ there will be none (unless $d=0$ in which case a constant map
with value a fixed point in $\CP^2$ would be fixed). Thus one must
remember to look for fixed equivalence classes. This means that
given a stable map $[u,\Sigma,p]$ and a $\lambda\in T$ one must find
an automorphism $\varphi$ of the underlying marked surface so that
$\lambda u=u\circ\varphi$.

\begin{example}\label{eg1}
An example of a genus zero degree three stable map fixed by the
natural $T$ action is $[u,\CP^1]$ where
$u([x_0:x_1])=[x_0^3:x_1^3:0]$. This is a stable map because the
only automorphisms are given by $\varphi([x_0:x_1])=[\xi x_0:x_1]$
where $\xi^3=1$. It is fixed by every $\lambda\in T$ because the map
$\varphi([x_0:x_1])=[\xi_0 x_0:\xi_1x_1]$ induces an equivalence
between $[\lambda\cdot u,\CP^1]$ and $[u,\CP^1]$ where
$\xi^3_0=\lambda_0$ and $\xi_1^3=\lambda_1$.
\end{example}
\begin{remark}
A stable map fixed by $T$ has an infinite number of left symmetries.
One should not confuse this with the requirement of at most a finite
number of right symmetries from the definition of a stable map. The
first is the $T$ action by post-composition. The second is
pre-composition by an automorphism of the marked surface.
\end{remark}

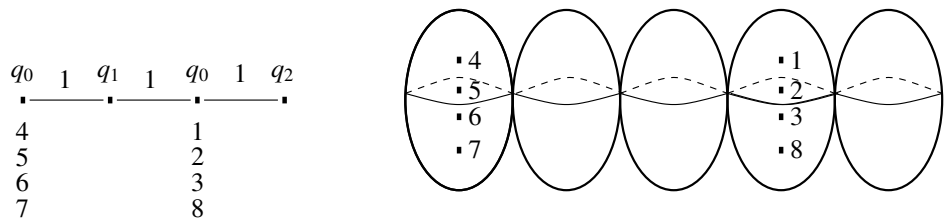
\begin{figure}[ht!]
\centering
\labellist\small
\pinlabel {$q_0$} [b] at 12 100
\pinlabel {$1$} [b] at 48 95
\pinlabel {$q_1$} [b] at 85 100
\pinlabel {$1$} [b] at 120 95
\pinlabel {$q_0$} [b] at 160 100
\pinlabel {$1$} [b] at 196 100
\pinlabel {$q_2$} [b] at 232 100
\pinlabel {$4$} [t] at 12 80
\pinlabel {$5$} [t] at 12 58
\pinlabel {$6$} [t] at 12 36
\pinlabel {$7$} [t] at 12 14
\pinlabel {$1$} [t] at 160 80
\pinlabel {$2$} [t] at 160 58
\pinlabel {$3$} [t] at 160 36
\pinlabel {$8$} [t] at 160 14
\pinlabel {$4$} [l] at 382 126
\pinlabel {$5$} [l] at 382 100
\pinlabel {$6$} [l] at 382 78
\pinlabel {$7$} [l] at 382 50
\pinlabel {$1$} [l] at 654 126
\pinlabel {$2$} [l] at 654 100
\pinlabel {$3$} [l] at 654 78
\pinlabel {$8$} [l] at 654 50
\endlabellist
\includegraphics[width=4.9truein]{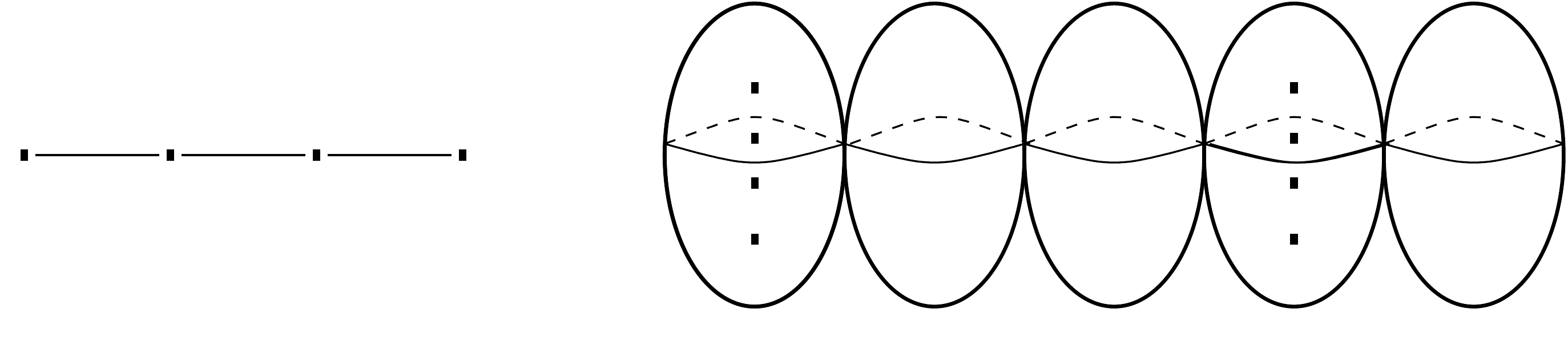} \caption{A labeled
graph and element}\label{fixstab}
\end{figure}

In general, a \index{$\Gamma$ component of Fix} component of the
fixed point set of the $T$ action on $\wwbar\calM_{g,n}(\CP^2,d)$
can be described by a labeled graph, denoted $\Gamma$. The labels on
the graphs and the correspondence between labeled graphs and
components of the fixed point set will be described in
\fullref{locform}. See \fullref{fixstab} for a labeled graph and
element of the corresponding fixed point component.

Applying the localization formula \eqref{loc} to the Gromov--Witten
invariants of $\CP^m$ gives
\begin{equation}\label{cpnlocfor}\bigl\langle h^{\ell_1}\ldots
h^{\ell_n}\bigr
  \rangle^{\CP^m}_{g,d[\CP^1]}=\sum_\Gamma
  \frac{1}{|\A_\Gamma|}\int_\Gamma\frac{\prod_{j=1}^n
\what{h}|_{k(u(p_j))}^{\ell_j}}{e(N_\Gamma^{\vir})}.\end{equation}
The notations used in this formula are explained further in the next
article. Here $\what h$ is an equivariant lift of \index{$H$@$\what
h$ equivariant lift of $h$} the class $h$. The standard lift is also
denoted by $h$, in which case $\what h|_{k(u_j)}=\alpha_{k(u_j)}$.
More generally, lifts of $h$ take the form $h+\sum a_j\alpha_j$ and
$(h+\sum a_j\alpha_j)|_{k(u_j)}=\alpha_{k(u_j)}+\sum a_j\alpha_j$.
We have to divide each term by a group factor to take into
\index{$N_\Gamma^{\vir}$ normal bundle to $\Gamma$} account the fact
that the moduli space is a stack (think orbifold). The group
$\A_\Gamma$ is the automorphism group of a generic element of the
fixed point set $\Gamma$. The automorphisms of the labeled graph
$\text{Aut}(\Gamma)$ \index{$A$@$\text{Aut}(\Gamma)$ automorphism
group of labeled graph} act on the edges of the graph and therefore
on the group $\prod_{e}\Z_{d(e)}$. The group $\A_\Gamma$
\index{$A$@$\A_\Gamma$ automorphism group of generic element of
component $\Gamma$} is defined as the following semidirect product:
$$
\A_\Gamma:=\text{Aut}(\Gamma)\ltimes
\left(\prod_{e}\Z_{d(e)}\right)\,.
$$
Formula \eqref{cpnlocfor} gives the Gromov--Witten invariants of
$\CP^m$ for any genus. For genus zero one just has to include the
automorphism groups; for higher genus one also has to apply virtual
localization. Of course, there is still work to do to get numbers
from this formula. In particular we still have to compute the Euler
class $e(N_\Gamma^{\vir})$.

The next step in the localization computation is the computation of
the Euler class of the normal bundle to each component of the fixed
point set. We have seen that the tangent space to the moduli stack
at a stable curve is ${\mathfrak{def}}([u,\Sigma,p])$.  Assuming
that $[u,\Sigma,p]$ is a fixed point of the $T$ action there will be
an induced action on the tangent space. This linear space decomposes
into a collection of irreducible $T$ representations. The tangent
space to the fixed point set at $[u,\Sigma,p]$ is the sum of the
trivial $T$ representations because the fixed point set is well, fixed.
The normal space is therefore the sum of the nontrivial
representations and is denoted by \index{$V^{\text{mov}``}$ maximal
nontrivial subrepresentation}
${\mathfrak{def}}([u,\Sigma,p])^{\text{mov}}$. In fact, we will see
that the $T$ action extends to all of the spaces in the
deformation-obstruction complex, so the representation
${\mathfrak{def}}([u,\Sigma,p])^{\text{mov}}$ may be deduced from
similar parts from the other terms in the deformation-obstruction
complex.

Recall that given an exact sequence of vector bundles
$$
0\to E_0\to E_1\to \cdots\to E_{2n+1}\to 0\,,
$$
one has the relation
$\prod_{k=0}^ne(E_{2k})=\prod_{k=0}^ne(E_{2k+1})$. Applying this to
the moving part of the deformation-obstruction sequence
\eqref{defob}, we obtain
$$
e({\mathfrak{def}}([u,\Sigma,p])^{\text{mov}}) =\frac{
e(\mathfrak{def}(u)^{\text{mov}})
e({\mathfrak{def}}([\Sigma,p])^{\text{mov}})
e({\mathfrak{ob}}([u,\Sigma,p])^{\text{mov}})}{
e(\mathfrak{aut}([\Sigma,p])^{\text{mov}})
e({\mathfrak{ob}}(u)^{\text{mov}})}\,.
$$
This last formula divides the computation of the Euler class of the
normal bundle required for localization calculations into more
tractable parts. Instead of analyzing all deformations of a stable
map we are able to analyze deformations of the map that fix the
marked curve, deformations of the marked curve, etc separately.
This program is carried out in detail in Sections \ref{seau},
\ref{tact}, \ref{sedef} and \ref{mape} where we compute
$e(\mathfrak{aut}([\Sigma,p])^{\text{mov}})$ first, describe the
torus actions second, compute
$e({\mathfrak{def}}([\Sigma,p])^{\text{mov}})$ third, and the ratio
$e({\mathfrak{def}}(u)^{\text{mov}})/e({\mathfrak{ob}}(u)^{\text{mov}})$
last. The term $e({\mathfrak{ob}}([u,\Sigma,p])^{\text{mov}})$ is
addressed in \fullref{fmc}. Before embarking on this program
we collect all of the resulting formulas and provide two examples.

\subsection{Representation of fixed point components by graphs}\label{locform}

We now describe the correspondence between components of the fixed
point set and labeled graphs. This article introduces the notation
used throughout this and the following sections. Recall that
the domain of a stable map is a prestable genus $g$ curve with
marked points. A stable map fixed by the standard torus action maps
each node of the curve to one of the points
\index{$Q$@$q_i:=[0:\ldots:1:\ldots:0]$}
$q_i:=[0:\ldots:1:\ldots:0]$, which are the points fixed  by the
standard action on $\CP^m$. Some irreducible components of the
prestable curve are mapped entirely to a single $q_i$.  We call such
components contracted components or ghost bubbles. Note that there
are no non-constant maps of higher genus curves into $\CP^m$ fixed
by the standard torus action. Therefore if a component is not
contracted it has to be a copy of $\CP^1$ that is mapped onto a
projective line containing exactly two of the $q_i$ fixed points.
Notice that any two stable maps in the same component of the fixed
point set have the same non-contracted $\CP^1$ configuration.

Thus we can describe a component of the fixed point set by a graph
with edges corresponding to non-contracted components and vertices
corresponding to the components of the preimages of the $q_i$. Each
edge $e$ is labeled with a positive integer $d(e)$ indicating the
degree of the map and each vertex is labeled with one of the fixed
points $q_i$. In addition, marked points are listed in columns under
the vertices (see \fullref{fixstab}).  Finally, contracted
components unlike non-contracted ones, can have higher genus and in
principle this has to be indicated at the vertices as well. We adopt
the convention that the absence of such a label corresponds to genus
$0$. In particular, in this subsection we are only concerned with
the Gromov--Witten invariants of rational curves and no genus
labeling is necessary. \fullref{glab} summarizes the notation used
to label these graphs.

\begin{remark}
One can visualize a flag by drawing an arrow on the edge; the source
of the arrow along with the edge is the flag. Notice that this use
of the term flag agrees with the usual definition in terms of
increasing sequences of subspaces (see Harris \cite{AG1} and
Griffiths--Harris \cite{GH}) because a vertex in
the graph corresponds to a point in $\CP^m$ which is a
$1$--dimensional subspace of $\C^{m+1}$ and an edge in the graph
corresponds to a line in $\CP^m$ which is a $2$--dimensional subspace
of $\C^{m+1}$.
\end{remark}
\begin{example}\label{graphnot}
For the graph in \fullref{fixstab} let $v$ denote the leftmost
vertex. Then $\text{val}(v)=1$, $d(v)=1$, $n(v)=4$, $k(v)=0$,
$k(v^\prime)=1$, $g(v)$=0, and $k(u(p_7))=0$.
\end{example}

\index{$V$@val($v$) valence of $v$} \index{$N$@$n(v)$  number of
marked points} \index{$K$@$k(v)$ index of fixed
point}\index{$G$@$g(v)$ genus of contracted
component}\index{$D$@$d(e)$ degree non-contracted component}

\begin{table}[ht!]
\begin{tabular}{|l|l|}
\hline
Notation & Description  \\
\hline $\val(v)$& The valence of the vertex $v$, that is, the number of
edges incident to \\ & it. If $\val(v)=1$ we let $v^\prime$ denote
the vertex on the other side of \\ & this single edge and if
$\text{val}(v)=2$
let $v_1$ and $v_2$ denote the other \\ & two vertices on the two edges. \\
$n(v)$ & The number of marked points listed under the vertex. Some
 \\ &  authors draw graphs with `legs' to indicate the marked points. \\
$k(v)$ & The index of the fixed point in $\CP^m$ corresponding to
the \\ & vertex,   for example if $v$ is labeled by $q_i$ then $k(v)=i$. \\
$\alpha_k$ & The first Chern class of the line bundle associated to
the dual of
\\&  the representation given by projection to the $k$th
factor of the \\ & torus.\\
$g(v)$ & The genus of the contracted component associated to the
vertex. \\& We set $g(v)=0$ if there is no such component at the
vertex.
\\
$d(e)$ & The degree of the map of the non-contracted component \\ &
corresponding to edge $e$. If $v$ has valence one $d(v)$ will be the
\\ & degree of the unique edge meeting $v$. If $v$ has valence two
$d(e_1)$ \\ & and $d(e_2)$ will
denote the degrees of the edges containing $v_1$ and \\ & $v_2$ respectively.\\
$F$ & A flag in the graph of a fixed point component, that is, a pair of
a \\ & vertex and an incident edge. We will use the flag to denote
the \\& corresponding vertex or point on the prestable curve or
image \\ & point
in $\CP^m$ without comment. \\
$k(u(p_j))$ & The label image of the $j$th marked point.
\\ \hline
\end{tabular}
\vskip.1in
\caption{Graph labels}
\label{glab}
\end{table}

We are now ready to present the formulas that are used in
localization computations of Gromov--Witten invariants of $\CP^m$.
\begin{remark}
We note that the index $j$ in \eqref{eT}--\eqref{eH1} runs over all
possible values, not just the $q_j$ depicted on the graph. For
example, if we are in $\CP^3$ one should take into account terms
with $j=3$ even when $q_3$ is not on the graph. See Remarks
\ref{can1} and \ref{can2} in the proofs.
\end{remark}

\subsection{Formulas used in localization}
\begin{equation}\label{eaut} e({\mathfrak{aut}}
([\Sigma,p])^{\text{mov}})=\prod_{\stackrel{
\text{\begin{tiny}val\end{tiny}}(v)=1}{\stackrel{n(v)=0}{g(v)=0}}}
\frac{\alpha_{k(v)}-\alpha_{k(v^\prime)}}{d(v)}\,.
\end{equation}
\begin{multline}\label{edefs}
e({\mathfrak{def}}([\Sigma,p])^{\text{mov}})=
\prod_{\stackrel{\val(v)=2}{\stackrel{n(v)=0}{g(v)=0}}}
\Bigl(\frac{\alpha_{k(v)}-\alpha_{k(v_1)}}{d(e_1)}+
\frac{\alpha_{k(v)}-\alpha_{k(v_2)}}{d(e_2)}\Bigr) \\[-4ex]
\prod_{\val(F)+n(F)+2g(v)>2}\Bigl(\frac{\alpha_{k(F)}-
\alpha_{k(v^\prime)}}{d(F)}-\psi_F\Bigr).
\end{multline}
\begin{equation}\label{edefobu}
\frac{e\bigl({\mathfrak{def}}(u)^{\text{mov}}\bigr)}
  {e\bigl({\mathfrak{ob}}(u)^{\text{mov}}\bigr)} =
\frac{e\bigl(H^0\bigl(\wwhat\Sigma,\calO_{\wwhat\Sigma}
  \bigl(\nu^*u^*T\CP^2\bigr)\bigr)^{\text{mov}}\bigr)}
{e\bigl(H^1\bigl(\wwhat\Sigma,\calO_{\wwhat\Sigma}
  \bigl(\nu^*u^*T\CP^2\bigr)\bigr)^{\text{mov}}\bigr)
\prod_ce\bigl(T_{u(c)}\CP^2\bigr)}.
\end{equation}
\begin{equation}\label{eT}
\prod_ce(T_{u(c)}\CP^2))=\prod_F\prod_{j\ne k(F)}(\alpha_{k(F)}-\alpha_j).
\end{equation}
\begin{multline}\label{eH0}
e\bigl(H^0\bigl(\wwhat\Sigma,
  \calO\bigl(u^*T\CP^2\bigr)\bigr)^{\text{mov}}\bigr)= \\[-1ex]
\prod_{v}\prod_{j\ne k(v)}(\alpha_{k(v)} -\alpha_j)
\prod_e
\frac{(-1)^{d(e)}(d(e)!)^2}{d(e)^{2d(e)}}\bigl(\alpha_{k(v)}
  -\alpha_{k(v^\prime)}\bigr)^{2d(e)} \\[-2ex]
\prod_{\stackrel{a+b=d}{j\ne k(v), k(v^\prime)}}
  \hspace{-1em}\Bigl(\frac{a}{d(e)}\alpha_ {k(v)}
  +\frac{b}{d(e)}\alpha_{ k(v^\prime)} -\alpha_j\Bigr).
\end{multline}
\begin{multline}\label{eH1}
e(H^1(\wwhat\Sigma,\calO(\nu^*u^*T\CP^2))^{\text{mov}})= \\[-3ex]
\prod_{\text{\begin{tiny}val\end{tiny}}(v)+n(v)+2g(v)>2}~\prod_{j\ne
k(v)}~\sum_{i=0}^{g(v)}c_i(\E^\vee)(\alpha_{k(v)}-\alpha_j)^{g-i}.
\end{multline}
\begin{multline}\label{eN}
e(N_\Gamma)=
e({\mathfrak{def}}([u,\Sigma,p])^{\text{mov}}) = \\
\frac{
e(\mathfrak{def}(u)^{\text{mov}})
e({\mathfrak{def}}([\Sigma,p])^{\text{mov}})
e({\mathfrak{ob}}([u,\Sigma,p])^{\text{mov}})}{
e({\mathfrak{aut}}([\Sigma,p])^{\text{mov}})
e({\mathfrak{ob}}(u)^{\text{mov}})}\,.
\end{multline}
\begin{equation}\label{GWloc}
\langle h^{\ell_1}\ldots
h^{\ell_n}\rangle^{\CP^m}_{g,d[\CP^1]}=\sum_\Gamma
\frac{1}{|\A_\Gamma|}\int_\Gamma\frac{\prod_{j=1}^n \what
h\big|_{k(u(p_j))}^{\ell_j}}{e(N_\Gamma)}\,.
\end{equation}

\begin{remark}
In the localization formula \eqref{GWloc} the sum is taken over all
labeled graphs with $n$ marked points having the correct genus and
degree. The integral is taken over the moduli space $\calM_{\Gamma}$
of all stable maps having the given graph. To evaluate these
integrals one has to expand $\frac{1}{e(N_\Gamma)}$ into an infinite
series in the $\psi$ classes and integrate the (finite number of) terms
of top degree. Since the values of the integral are in the
equivariant cohomology ring of $\CP^m$ the Chern classes $\alpha_k$
play the role of constants and can simply be factored out of the
integral. Thus, actual integration is only required when one of the
descendant classes $\psi_F$ from \eqref{edefs} is non-zero in which
case we can use the string and dilaton equations of \fullref{323}
to evaluate the integral. Sample computations in the next article
will clarify the details and should convince the reader that it is
possible to extract useful information from these cumbersome
formulas.
\end{remark}

\subsection{Small degree invariants of rational curves in
$\CP^2$}\label{N3}

In this article we demonstrate the localization formulas by
computing the genus zero invariants in degree one, two and three. We
start with
$$N_1:=\langle h^2h^2\rangle^{\CP^2}_{0,1[\CP^1]},$$
where $h$ is the Poincar\'e dual to the hyperplane class. We can
lift $h$ to an equivariant class that we denote by the same letter.
Intuitively, we are looking for the number of lines through two
generic points in $\CP^2$ so you can guess that the answer is one,
but it is instructive to see how localization produces this answer.

We are working with fixed stable maps in
$\wwbar\calM_{0,2}(\CP^2,[\CP^1])$. Since the overall degree is one,
our graph can only have one edge and it will be labeled with a $1$.
The two vertices are labeled with two of the three points
$q_0=[1:0:0]$, $q_1=[0:1:0]$, $q_2=[0:0:1]$. We also have to
distribute the $2$ marked points between the two vertices and there
are two essentially different ways to do this: to place the marked
points at the same vertex or at different vertices. This leads to
two different labeled graph types. There are a total of $12$ graphs
to consider, $6$ of type one and $6$ of type two. Due to the obvious
symmetry, contributions from all graphs within a type are similar.

For the first type we consider the graph labeled by $q_0$ and $q_1$
with the two marked points placed at $q_0$. The only automorphisms
are trivial so $|\text{Aut}(\Sigma_\Gamma)|=1$ and since both marked
points are mapped to $q_0$ we have
$\iota^*\text{ev}_1^*h^2=\alpha_0^2$ and
$\iota^*\text{ev}_2^*h^2=\alpha_0^2$.

Following \eqref{eaut}--\eqref{eN} we get
\begin{align*}
e({\mathfrak{aut}}([\Sigma,p])^{\text{mov}})
&=\frac{\alpha_1-\alpha_0}{1}=\alpha_1-\alpha_0\\
e({\mathfrak{def}}([\Sigma,p])^{\text{mov}})
&=1\cdot1\cdot\Bigl(\frac{\alpha_0-\alpha_1}{1}-\psi\Bigr)=\alpha_0-\alpha_1
\end{align*}
Recall that the $\psi$--classes are the first Chern classes of the
tautological bundles. In this case the components of the fixed point
set are points so all $\psi$--classes vanish.
\begin{align*}
&\prod_ce(T_{u(c)}\CP^2))
=(\alpha_0{-}\alpha_1)(\alpha_0{-}\alpha_2){\cdot}
  (\alpha_1{-}\alpha_0)(\alpha_1{-}\alpha_2)
e(H^0(\wwhat\Sigma,\calO(u^*T\CP^2))^{\text{mov}}) \\[-1ex]
&\quad=(\alpha_0{-}\alpha_1)(\alpha_0{-}\alpha_2)\cdot
  (\alpha_1{-}\alpha_0)(\alpha_1{-}\alpha_2)
\cdot
\frac{(-1)^1(1!)^2}{1^{2\cdot1}}(\alpha_0{-}\alpha_1)^{2\cdot1}
(\alpha_0{-}\alpha_2)(\alpha_1{-}\alpha_2)
\\
&e(H^1(\wwhat\Sigma,\calO(\nu^*u^*T\CP^2))^{\text{mov}})
=1,\qquad\text{since all vertices have genus $0$.}
\end{align*}
In fact we could have seen that
$e(H^1(\wwhat\Sigma,\calO(u^*T\CP^2))^{\text{mov}})= 1$ directly
since the fiber of the obstruction bundle ${\mathfrak{ob}}(u)$ is
$H^1(\Sigma,\calO_\Sigma (u^*T\CP^2))=0$. Using the
deformation-obstruction sequence \eqref{defob} this implies that
$e({\mathfrak{ob}}(u)^{\text{mov}})=0$ as well. Continuing,
\begin{align*}
\frac{e({\mathfrak{def}}(u)^{\text{mov}})}{e({\mathfrak{ob}}(u)^{\text{mov}})}
&=(\alpha_0-\alpha_1)^2(\alpha_0-\alpha_2)(\alpha_1-\alpha_2)\\
e(N_\Gamma)
&=\frac{\alpha_0-\alpha_1}{\alpha_1-\alpha_0}\cdot-(\alpha_0-\alpha_1)^2(\alpha_0-\alpha_2)(\alpha_1-\alpha_2)\\
&=(\alpha_0-\alpha_1)^2(\alpha_0-\alpha_2)(\alpha_1-\alpha_2)\,.
\end{align*}
The computation of the Euler class of the normal bundle for the
second graph type is analogous and we leave it as an exercise.
\begin{exm}
Repeat the above computation for the one-edge graph labeled with
$q_0$ and $p_1$ labeling one vertex and $q_1$ and $p_2$ labeling the
other. You should get
$$e(N_\Gamma)=-(\alpha_0-\alpha_1)^2(\alpha_0-\alpha_2)(\alpha_1-\alpha_2).$$
\end{exm}
The Euler classes of the normal bundles for the other graphs can be
obtained from the first two examples by symmetry considerations.
There are no $\psi$--classes to integrate so the integral in
\eqref{GWloc} can be dropped. It is convenient to first sum up all
the contributions from graphs labeled with $q_0, q_1$. There are
four such graphs and the remaining two can be obtained from the ones
we already computed by switching $\alpha_0$ and $\alpha_1$. By
\eqref{GWloc} the contribution to $N_1$ from these four graphs is
\begin{multline*}
\frac{1}{(\alpha_0-\alpha_1)^2(\alpha_0-\alpha_2)(\alpha_1-\alpha_2)}\cdot(\alpha_0^2\cdot\alpha_0^2+\alpha_1^2\cdot\alpha_1^2)\\
-\frac{1}{(\alpha_0-\alpha_1)^2(\alpha_0-\alpha_2)(\alpha_1-\alpha_2)}\cdot(\alpha_0^2\cdot\alpha_1^2+\alpha_1^2\cdot\alpha_0^2)\\
=\frac{(\alpha_0^2-\alpha_1^2)^2}{(\alpha_0-\alpha_1)^2(\alpha_0-\alpha_2)(\alpha_1-\alpha_2)}
\end{multline*}
There are eight more graphs to account for, four labeled with
$q_0,q_2$ and four labeled with $q_1,q_2$. The joint contributions
of each four are obtained from the last expression by applying the
obvious substitutions.
\begin{exm}
Add up the three fractions and get $N_1=1$ as expected.
\end{exm}

We can repeat the previous computation with different linearizations
(lifts to $H^*_T(\CP^2)$ of the class $h$). For example, since $h$
is the first Chern class of the $\calO(1)$ line bundle we can
construct lifts as the equivariant Chern classes of the same line
bundle with various group actions. The lift with the same name comes
from the first Chern class of the standard extension of the $T$
action on $\CP^2$ to this bundle. We may take the tensor product
with a trivial line bundle with the $\mu_1$--action to change the
action. This does not change the Chern class but it does change the
equivariant lift to $\what h=h-\alpha_1$. If we use
$(h-\alpha_1)(h-\alpha_2)$ as the linearization for the cohomology
class corresponding to the first marked point and
$(h-\alpha_0)(h-\alpha_2)$ for the second marked point then only one
of the $12$ graphs will contribute to the sum. Namely, the graph
with one vertex labeled with $q_0$ and $p_1$ and the other vertex
labeled with $q_1$ and $p_2$.
\begin{exm}
Repeat the computation of
$$\langle h^2h^2\rangle^{\CP^2}_{0,1[\CP^1]}$$
using the
$(h-\alpha_1)(h-\alpha_2)$ and $(h-\alpha_0)(h-\alpha_2)$
linearizations.
\end{exm}

Generally speaking components of the fixed point set of a torus
action on a moduli space of stable curves can be expressed as a
finite quotient of a product of moduli spaces of stable maps into
the fixed point set of the target manifold. The best way to
understand this is to work out some less trivial examples.

To demonstrate a less trivial localization computation we will
conclude the subsection with a computation of
$$N_2:= \langle h^2h^2h^2h^2h^2\rangle^{\CP^2}_{0,2[\CP^1]}$$
and some of the terms from the computation of
$$N_3:= \langle h^2h^2h^2h^2h^2h^2h^2h^2\rangle^{\CP^2}_{0,3[\CP^1]}.$$
It is convenient to use the linearization $(h-\alpha_1)(h-\alpha_2)$ on
the class associated to each marked point. This forces all of the marked
points to be mapped to $q_0$ if a graph is going to contribute to the sum.

\begin{figure}[ht!]
\centering
\labellist\small
\pinlabel {$I(k_0k_1)$} [l] at -20 220
\pinlabel {$k_0$} [t] at 114 215
\pinlabel {$k_1$} [t] at 168 215
\pinlabel {$\II(k_0k_1k_2)$} [l] at -20 90
\pinlabel {$k_0$} [t] at 88 80
\pinlabel {$k_1$} [t] at 142 80
\pinlabel {$k_2$} [t] at 195 80
\pinlabel {$\III(k_0k_1k_2k_3)$} [l] at 260 220
\pinlabel {$k_0$} [t] at 398 220
\pinlabel {$k_1$} [t] at 450 220
\pinlabel {$k_2$} [t] at 505 220
\pinlabel {$k_3$} [t] at 560 220
\pinlabel {$\IV(k_0k_1k_2k_3)$} [l] at 260 90
\pinlabel {$k_0$} [tl] at 477 87
\pinlabel {$k_1$} [t] at 477 34
\pinlabel {$k_2$} [l] at 508 123
\pinlabel {$k_3$} [r] at 445 123
\endlabellist
\includegraphics[width=4.5truein]{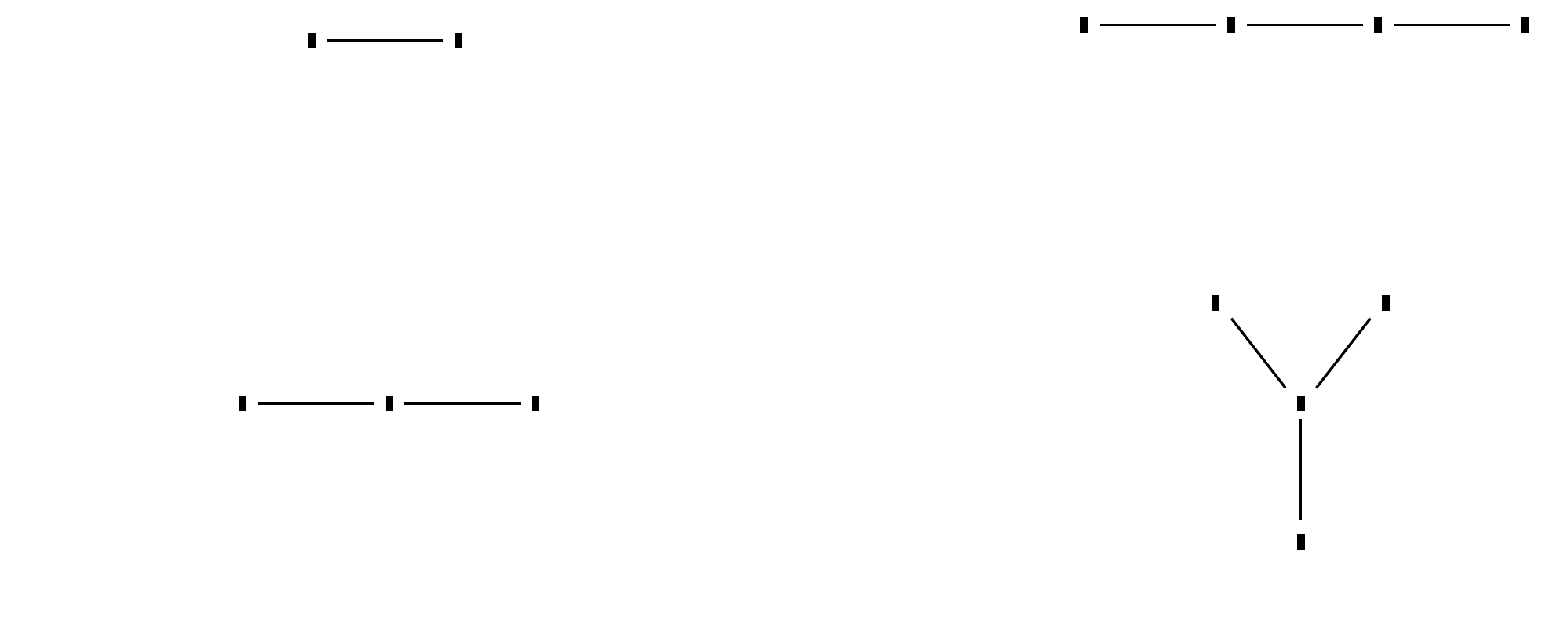} \caption{Graph
types}\label{gtype}
\end{figure}
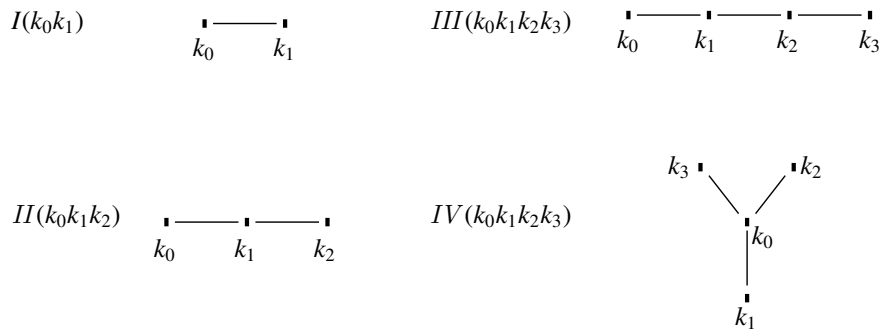

The graphs that contribute to the $N_2$ computation come in one of two
types. We label these two types by $I(k_0k_1)$ and $\II(k_0k_1k_2)$. The
corresponding graphs are displayed in \fullref{gtype} along with the
graphs that contribute to the $N_3$ computation.

Each graph type contains all of the graphs obtained by an admissible
labeling. An admissible labeling consists of assigning each $k_j$ a
value of $0$, $1$ or $2$ such that one of the $k_j$ is zero and no
adjacent two are equal. In addition an admissible labeling includes
an assignment of the five marked points to vertices labeled with a
zero. The vertex labeled with $k_j$ is really labeled by $q_{k_j}$;
the $k$--labels just avoid the double index notation and this is
consistent with our use of $k$ in the formulas.

\begin{exm}
Compute the contributions to $N_2$ of several graphs from
\fullref{gtype}. You can check your answer with the results listed
in \fullref{app:b}. Notice that the sum of all the contributions is
equal to one as we outlined in \fullref{Nd}.
\end{exm}

This last exercise may be difficult, but we feel that the interested
reader will learn more by doing it than just reading the answer. We
won't feel too guilty since the answer is in an appendix. We will
provide more explanation for some of the contributions to $N_3$.

The graphs that can contribute to $N_3$ come in one of four
different types. These types are displayed in \fullref{gtype}.
Once again each graph type contains all of the graphs obtained by an
admissible labeling as described in the $N_2$ case. This time an
admissible labeling includes an assignment of the eight marked
points to vertices labeled with a zero.

As an example of this notation the graph and curve displayed in
\fullref{fixstab} are of type $\III$.
\begin{example}
To be specific the graph from the figure is $\III(0102)$ with the
marked points $p_4$, $p_5$, $p_6$ and $p_7$ on the $k_0$ vertex. In
this case the moduli stack of the fixed point component is given by
$${\wwbar\calM}_\Gamma={\wwbar\calM}_{0,5}(\text{pt},0) \times
{\wwbar\calM}_{0,6}(\text{pt},0),$$
and the inclusion map
$\iota_\Gamma\co {\wwbar\calM}_\Gamma\to{\wwbar\calM}_{0,8}(\CP^2,3[\CP^1])$
is given by
\begin{multline*}
\iota_{\Gamma}\bigl([\Sigma,p_1,p_2,p_3,p_4,p_5],
  [\Sigma^\prime,p_1^\prime,p_2^\prime,p_3^\prime,
  p_4^\prime,p_5^\prime,p_6^\prime]\bigr) = \\
\Bigl[u,\left(\Sigma\cup\Sigma^\prime\right)
  \bigcup_{\stackrel{p_1=
  ([0:1],1)}{\stackrel{p_1^\prime=
  ([1:0],2)}{\scriptscriptstyle p_2^\prime=([0:1],3)}}}\bigl(\CP^1\times\{1,2,3\}\bigr)\big/\sim,
  p_3^\prime, p_4^\prime,p_5^\prime,p_2,p_3,p_4,p_5,p_6^\prime\Bigr],
\end{multline*}
where $([1:0],1)\sim([0:1],2)$, $([1:0],2)\sim([0:1],3)$ and the map
$u$ is given by
\begin{align*}
u([x_0:x_1],1)&=[x_1:x_0:0], & u([x_0:x_1],2)&=[x_0:x_1:0], \\
u([x_0:x_1],3)&=[x_1:0:x_0], & u|_{(\Sigma\cup\Sigma^\prime)}&=[1:0:0].
\end{align*}
\end{example}

Applying formulas \eqref{eaut}--\eqref{eN} to the labeled graph from
\fullref{fixstab} gives
\begin{multline*}
e(N_\Gamma)^{-1}=
-\tfrac12(\alpha_0{-}\alpha_1)^{-4}(\alpha_0{-}\alpha_2)^{-2}(\alpha_1{-}\alpha_2)^{-2} \\
\int_{[{\wwbar\calM}_{0,5}(\text{pt},0)]^{\vir}}(\alpha_0{-}\alpha_1{-}\psi_1)^{-1}
{\cdot}\int_{[{\wwbar\calM}_{0,6}(\text{pt},0)]^{\vir}}(\alpha_0{-}\alpha_1
{-}\psi_1)^{-1}(\alpha_0{-}\alpha_2{-}\psi_2)^{-1}
\end{multline*}
Now consider the second integral that appears in the above
expression. We compute
\begin{align*}
\int_{[\wwbar\calM_{0,n}(\text{pt},0)]^{\vir}}&(a_1-\psi_1)^{-1}(a_2-\psi_2)^{-1}\\
&=\int_{[\wwbar\calM_{0,n}(\text{pt},0)]^{\vir}}a_1^{-1}a_2^{-1}\left(\sum_{p=0}^\infty
a_1^p\psi_1^p\right) \left(\sum_{q=0}^\infty
a_2^q\psi_1^q\right)  \\
&=a_1^{-1}a_2^{-1}\sum_{p+q=n-3}a_1^{-p}a_2^{-q}\int_{[\wwbar\calM_{0,n}(\text{pt},0)]^{\vir}}
\psi_1^p\psi_2^q \\
&= a_1^{-1}a_2^{-1}\sum_{p+q=n-3}\binom{n-3}{p\ q}a_1^{-p}a_2^{-q}=a_1^{-1}a_2^{-1}(a_1^{-1}+a_2^{-1})^{n-3}.
\end{align*}
The last line in this computation used the result from \fullref{psi3}.

The graph $\III(0102)$ has no nontrivial automorphisms and each edge
has degree one, so this component of the fixed point set has only
trivial automorphisms. Combining the two previous computations with
the $(h-\alpha_1)(h-\alpha_2)$ linearization allows us to conclude
that the contribution of the graph $\III(0,1,0,2)$ with the marked
points $p_4$, $p_5$, $p_6$ and $p_7$ on the $k_0$ vertex is
$$
-\tfrac12(\alpha_0-\alpha_1)^{-4}(\alpha_0-\alpha_2)^{5}(\alpha_1-\alpha_2)^{-2}
\bigl((\alpha_0-\alpha_1)^{-1}+(\alpha_0-\alpha_2)^{-1}\bigr)^3.$$
Of course there are $\binom{8}{4}$ ways to choose four marked points
for the first vertex. This means that the contribution of such
graphs is
$$
-\frac12\binom{8}{4}
(\alpha_0-\alpha_1)^{-3}(\alpha_0-\alpha_2)^{2}(\alpha_1-\alpha_2)^{-2}
(2\alpha_0-\alpha_1-\alpha_2)^{3}.$$

\begin{exm}
Show that the contribution of the $\III(0102)$ graphs with $k$ marked
points on the first vertex is
$$-\frac12\binom{8}{k}
(\alpha_0-\alpha_1)^{-3}(\alpha_0-\alpha_2)^{k-2}(\alpha_1-\alpha_2)^{-2}
(2\alpha_0-\alpha_1-\alpha_2)^{7-k}.$$
Conclude that the total contribution to $N_3$ from all $\III(0102)$
graphs is
$$-\tfrac12
(\alpha_0-\alpha_1)^{-3}(\alpha_0-\alpha_2)^{-2}(\alpha_1-\alpha_2)^{-2}
(2\alpha_0-\alpha_1-\alpha_2)^{-1}(3\alpha_0-\alpha_1-2\alpha_2)^8.$$
\end{exm}

\begin{exm}
Explain why there is no graph corresponding to $\II(002)$. Explain
why there is no graph contribution corresponding to $I(21)$.
\end{exm}

\begin{exm}
Compute the contributions to $N_3$ of several graphs from \fullref{gtype}. Using some mathematical software it would be possible
to push this computation all the way and get the answer $N_3=12$.
\end{exm}

\section{Derivation of the Euler class formulas}

In this subsection we derive the general formulas that were used to
compute the various Euler classes used in localization computations.
We derive the formulas for $\CP^2$; however careful inspection shows
that these formulas are valid for $\CP^m$ as well.

\subsection{The Euler class of moving infinitesimal
automorphisms}\label{seau}

Before computing $e(\mathfrak{aut}([\Sigma,p])^{\text{mov}})$ we
need to understand the bundle of automorphisms. First consider the
automorphisms of a genus zero irreducible component with no marked
points mapped into $\CP^2$ by the composition of a degree $d$ map
from $\CP^1$ to $\CP^1$ and an inclusion of $\CP^1$ into $\CP^2$.
This is exactly the situation described in the following example. We
assume that this component is attached to the rest of the stable map
by a node at $[1:0]$.

\begin{example}
The $T$--action on ${\mathfrak{def}}(u)$ is a good example  to study
the induced $T$--action on a space in the deformation complex.
Consider the $T$--fixed stable map given by
$u([x_0:x_1])=[x_0^d:x_1^d:0]$. The action of $T$ on this is just by
pointwise multiplication of the coordinates. However, as we saw in
\fullref{eg1} it is not immediately clear from this description
why this map is fixed. We would like to include the
reparametrization demonstrating that the class of this map is built
into the definition of the action. The reparametrization requires a
$d$th root of elements of $T$, but this can not be done in a
consistent way. This may be corrected by considering the
representation where $T$ acts by multiplication of $d$th powers of
the elements of $T$. The first Chern classes of the line bundles
associated to the irreducible factors of this representation are $d$
times the first Chern classes of the original representation. It
follows that we can get the Chern classes of the line bundles in the
original representation by dividing by $d$. The induced action on a
$1$--parameter family of maps (that is, a deformation) is then given by
\begin{multline*}
(\lambda\cdot u_t)([x_0:x_1]):= \\
\bigl[ \lambda_0^du^0_t(\lambda_0^{-1}x_0,\lambda_1^{-1}x_1):
\lambda_1^du^1_t(\lambda_0^{-1}x_0,\lambda_1^{-1}x_1):
\lambda_2^du^2_t(\lambda_0^{-1}x_0,\lambda_1^{-1}x_1)\bigr].
\end{multline*}
Expressed in this way it is clear why this action fixes the map $u$.
\end{example}

One can now determine the  $T$--action on
${\mathfrak{aut}}([\Sigma,p])$, and apply it to a specific bubble.
Recall that the map ${\mathfrak{aut}}([\Sigma,p])\to
{\mathfrak{def}}(u)$ is given by $\varphi_t\mapsto u\circ\varphi_t$.
Assuming that the domain of a stable map has a bubble with a node at
$[1:0]$ and this bubble is mapped into $\CP^2$ by a degree $d$ map
(as described in \fullref{eg2}) any automorphism must restrict to
this bubble to a map of the form
$\varphi([x_0:x_1])=[ax_0+bx_1:x_1]$. In order for the map
${\mathfrak{aut}}([\Sigma,p])\to {\mathfrak{def}}(u)$ to be
equivariant with respect to the $d$th power action, we must have
$$
(\lambda\cdot
\varphi_t)([x_0:x_1]):=[a_tx_0+\lambda_0\lambda_1^{-1}b_tx_1:x_1]\,.
$$
A nice way to compute the first Chern class of the corresponding
subbundle is to use the Borel construction. Recall that one can
combine a principal $T$--bundle with a $T$ representation to get a
vector bundle, as explained by the following exercise.
\begin{exm}
Let $L_{n_0,\ldots,n_1}$ be the line bundle associated to the
representation $\lambda\cdot
z=\lambda_0^{n_0}\ldots\lambda_k^{n_k}z$ and show that
$$c_1(L_{n_0,\ldots,n_1})=n_0c_1(L_{1,0,\ldots,0})+\cdots+n_kc_1(L_{0,\ldots,0,1})\,.$$
\end{exm}
We conclude that the first Chern class of the corresponding
subbundle is $(\alpha_{[1:0:0]}-\alpha_{[0:1:0]})/d$,
$\alpha_k$ is the Chern class of the bundle associated to the
divisor $\phi_k$ described in the computation of the equivariant
cohomology of $\CP^n$.

A nice way to keep track of extra factors like $1/d$ involved in
computing the Chern classes is to introduce the notion of a virtual
representation. Virtual representations are just $\Q$--linear
combinations of ordinary representations. The virtual representation
on the automorphisms arising from the usual $T$--action is defined to
be $\frac{1}{d}$ times the $d$th power representation.

We can split the automorphisms of a marked curve into automorphisms
of each irreducible component of the marked curve and compute the
Euler class of the bundle of infinitesimal automorphisms from the
automorphisms of the components. In fact we can split each factor in
\eqref{eN} into a sum of line bundles. The irreducible
representations of $T$ are all one-dimensional. These
representations induce line bundles over the components of the fixed
point set of the moduli stack. The Euler class of a complex vector
bundle is just the top Chern class of the bundle (see Bott and Tu
\cite{BT}) and the Chern class of a sum is given by
$$c_{\text{top}}(E\oplus F) = c_{\text{top}}(E) c_{\text{top}}(F).$$
To get to $e(\mathfrak{aut}([\Sigma,p])^{\text{mov}})$ consider
other irreducible components of the fixed stable map labeled by
graphs as described around \fullref{fixstab}. Each irreducible
component of the stable map contributes a summand to the bundle
$\mathfrak{aut}([\Sigma,p])^{\text{mov}}$ and hence a factor to the
Euler class. The factor of $e({\mathfrak
aut}([\Sigma,p])^{\text{mov}})$ corresponding to an edge in the
graph representing a component of the fixed point set that has two
nodes or a node and a marked point is trivial. The induced action on
the factors of ${\mathfrak{aut}}([\Sigma,p])$ corresponding to any
contracted component (that is, one for which the stable map is constant)
is trivial. The above discussion implies that
$$e({\mathfrak{aut}}
([\Sigma,p])^{\text{mov}})=\prod_{\stackrel{
\val(v)=1}{\stackrel{n(v)=0}{\scriptscriptstyle g(v)=0}}}
\frac{\alpha_{k(v)}-\alpha_{k(v^\prime)}}{d(e)}.$$

\subsection{The $T$--action on the deformation
complex}\label{tact}

Before addressing the deformations of the underlying curve, we
explain the $T$--action on all of the terms in the deformation
complex in greater detail. We can use the homological algebra
introduced to derive the deformation-obstruction complex to give a
uniform treatment of the induced actions on the terms in the
complex.

Given an action $T\times X\to X$ and a stable map fixed by this
action, $u\co \Sigma\to X$ one can construct the $d$th power action
and the corresponding virtual $T$--actions on $X$ and $\Sigma$ making
$u$ equivariant as we did in \fullref{eg2}.

\begin{exm}\label{eg2}
Let $u\co \CP^1\cup_{[0:1]=[1:0]}\CP^1\to\CP^2$ be given by
$$u([x_0:x_1],1)=[x_0^{d_1}:x_1^{d_1}:0]\quad\text{and}\quad
u([x_0:x_1],2)=[0:x_0^{d_2}: x_1^{d_2}];$$
the original action on
$\CP^2$ is given by $\lambda\cdot
[z_0:z_1:z_2]=[\lambda_0z_0:\lambda_1z_1;\lambda_2z_2]$; then the
map $u$ is equivariant with respect to the action on $\CP^1$ given
by
$$\lambda\cdot
([x_0:x_1],1)=([\lambda_0^{d_2}x_0:\lambda_1^{d_2}x_1],1)
\quad\text{and}\quad
\lambda\cdot
([x_0:x_1],2)=([\lambda_0^{d_1}x_0:\lambda_1^{d_1}x_1],2),$$
the $d_1d_2$ power of the original action on $\CP^2$.
\end{exm}

Recall that the push-forward \index{push-forward} of a sheaf $\calA$
under a map $f\co X\to Y$ is given by \index{$F$@$f_*\calA$
push-forward} $f_*\calA(\calU):=\calA(f^{-1}(\calU))$. Let
$L_\lambda$ represent left multiplication by $\lambda$. We will
consider the push-forward by $L_\lambda$ of various sheaves. There
are natural transformations given by pull-back on the various
sheaves,
\begin{align*}
L^*_\lambda\co \Omega_X^1\longrightarrow & L_{\lambda *}\Omega_X^1 \\
L^*_\lambda\co u^*\Omega_X^1\longrightarrow & u^*L_{\lambda *}\Omega_X^1 \\
L^*_\lambda\co \calO_\Sigma\longrightarrow & L_{\lambda *}\calO_\Sigma \\
L^*_\lambda\co \Omega_\Sigma^1([p])\longrightarrow &
  L_{\lambda *}\Omega_\Sigma^1([p]) \\
\end{align*}
The map in the first line is just the pull-back of holomorphic
$1$--forms on $X$ and the rest are analogous. The push-forward just
formalizes the fact that the pull-back of a form over $U$ is a form
over $L_\lambda^{-1}(U)$. Any of these pull-backs will induce a
group action on the space of sections over any invariant set. In
particular, they induce actions on the space of global sections. Our
convention will be to use left group actions everywhere. Recall that
a right action may be turned into a left action by taking the
inverse of the group element. These actions will induce actions on
the ext-groups. Recall that
$${\mathfrak{aut}}([\Sigma, p])=\E
xt^0_{\calO_{\Sigma}}(\Omega_\Sigma(p),\calO_\Sigma)\cong
\Hom_{\calO_\Sigma}(\Omega_\Sigma(p),\calO_\Sigma).$$
Given $X\in
\Hom_{\calO_\Sigma}(\Omega_\Sigma(p),\calO_\Sigma)$, the left
action is given by
$$(\lambda\cdot X)(\theta):=L^*_{\lambda^*}(X(L^*_\lambda\theta)).$$
Working in the $([x_0:1],1)$--chart we obtain
$$(\lambda\cdot((\dot ax_0+\dot b)\partial_x))(dx_0)=L_{\lambda^*}^*((\dot
ax_0+\dot b)\partial_x(L^*_\lambda dx_0))=\dot
ax_0+\lambda_0^{d_2}\lambda_1^{-d_2}\dot b.$$
Note that $\lambda\cdot([x_0:1],1)=([\lambda_0^{d_2}\lambda_1^{-d_2}],1)$
This agrees with our earlier computation of the action.
\begin{exm}
Check that the action on ${\mathfrak{def}}([\Sigma,p])$ derived via
homological algebra agrees with the action described earlier.
\end{exm}

\subsection{The Euler class of moving deformations of the
curve}\label{sedef}

Following Cox and Katz \cite{CK}, we can use local models to analyze the
$T$--action on ${\mathfrak{def}}([\Sigma,p])^{\text{mov}}$. Most
marked points are on contracted components. The deformations
corresponding to moving these points are tangent to the fixed point
set of the action. The one exception is when there is exactly one
marked point labeling a vertex of valence one. This corresponds to a
marked point on one of the branch points of a degree $d$ cover of a
standard line in $\CP^2$. The space of deformations corresponding to
moving such a point is trivial because a genus zero curve with two
marked points has no deformations. All of the normal deformations
arise from resolutions of nodes. The local model here is
$\Sigma=\{(x,y)\in \C^2|xy=0\}$. Let \index{$I$@$\calI_\Sigma$ ideal
sheaf} $\calI_\Sigma=(xy)$ be the sheaf of algebraic functions on
$\C^2$ containing a factor of $xy$ (this is called the ideal sheaf)
and notice that we have the following exact sequence,
$$0\to\calI_\Sigma/\calI_\Sigma^2\to\Omega^1_{\C^2}|_\Sigma\to
  \Omega^1_\Sigma\to 0.$$
The associated long exact sequence of ext-groups reads
\begin{multline*}
\to\Ext^0_{\calO_\Sigma}(\Omega^1_{\C^2}|_\Sigma,\calO_\Sigma)
\to\Ext^0_{\calO_\Sigma}(\calI_\Sigma/\calI_\Sigma^2,\calO_\Sigma) \\
\to \Ext^1_{\calO_\Sigma}(\Omega^1_\Sigma,\calO_\Sigma) \to\E
xt^1_{\calO_\Sigma}(\Omega^1_{\C^2}|_\Sigma,\calO_\Sigma)\to\cdots.
\end{multline*}
Now, $\Omega^1_{\C^2}|_\Sigma$ is a free $\calO_\Sigma$--module, so
$\Ext^1_{\calO_\Sigma}(\Omega^1_{\C^2}|_\Sigma,\calO_\Sigma)=0$.
It follows that
\begin{align*}
\Ext^1_{\calO_\Sigma}(\Omega^1_\Sigma,\calO_\Sigma)&\cong\text{coker}(\Ext^0_{\calO_\Sigma}(\Omega^1_{\C^2}|_\Sigma,\calO_\Sigma) \to\Ext^0_{\calO_\Sigma}(\calI_\sigma/\calI_\Sigma^2,\calO_\Sigma)) \\
&\cong T_0(\C\times 0)\oplus T_0(0\times\C)\,.
\end{align*}
These local results can be combined to give
${\mathfrak{def}}([\Sigma,p])^{\text{mov}}$. We first introduce or
recall some notation.

Let $F$ refer to a flag in the graph of the stable map (that is, a
pair consisting of a vertex and incident edge). We will use the flag
to denote the vertex or the point on $\Sigma$ or the point on
$\CP^2$ corresponding to the vertex without comment. The
normalization of the surface is denoted by $\wwhat\Sigma$
\index{$\wwhat\Sigma$ normalization}\index{normalization} and the
node branches will be denoted by $b_1$ and $b_2$. Recall that
$\calL_k$ denotes the line bundle over the moduli space whose fiber
is the cotangent space of the corresponding curve.

Combining these local results gives,
$${\mathfrak{def}}([\Sigma,p])^{\text{mov}}=
\Biggl(
\bigoplus_{\stackrel{\val(v)=2}
  {\stackrel{n(v)=0}{\scriptscriptstyle g(v)=0}}}
  \bigl(T_{b_1}\wwhat\Sigma\oplus T_{b_2}\wwhat\Sigma\bigr)
\Biggr)
\Biggl(
\bigoplus_{\val(F)+n(F)+2g(v)>2}\bigl(T_{F}\wwhat\Sigma\oplus
\calL_F^*\bigr)
\Biggr).$$
This formula translates directly into a formula for the Euler class,
\begin{multline*}
e({\mathfrak{def}}([\Sigma,p])^{\text{mov}})=
\Biggl(\prod_{\stackrel{\val(v)=2}{\stackrel{n(v)=0}
  {\scriptscriptstyle g(v)=0}}}
  \Bigl(\frac{\alpha_{k(v)}-\alpha_{k(v_1)}}{d(e_1)}+
  \frac{\alpha_{k(v)}-\alpha_{k(v_2)}}{d(e_2)}\Bigr) \Biggr)\\[-2ex]
\Biggl(\prod_{\val(F)+n(F)+2g(v)>2}
  \Bigl(\frac{\alpha_{k(F)}-\alpha_{k(v^\prime)}}{d(e)}-\psi_F \Bigr)
\Biggr).
\end{multline*}

\subsection{Euler class associated to the
map}\label{mape}

We now turn to the computation of the Euler classes
of ${\mathfrak{def}}(u)^{\text{mov}}$ and $\mathfrak{ob}(u)^{\text{mov}}$. The following exact sequence serves to define
holomorphic functions on a nodal curve in terms of holomorphic
functions on the normalization $\nu\co \wwhat\Sigma\to\Sigma$:
$$0\to\calO_\Sigma\to\nu_*\calO_{\wwhat\Sigma}\to\oplus_c\calO_c\to0.$$
Here we use $c$ to denote the nodes (crossings) of $\Sigma$. We have
a related exact sequence for holomorphic sections of the pull-back
bundle $u^*T\CP^2$:
$$0\to\calO_\Sigma(u^*T\CP^2)\to\nu_*\calO_{\wwhat\Sigma}(\nu^*u^*T\CP^2)
  \to\oplus_cT_{u(c)}\CP^2\to0.$$
The associated cohomology long exact sequence reads,
\begin{multline*}
0\to H^0(\Sigma, \calO_\Sigma(u^*T\CP^2))\to
  H^0(\wwhat\Sigma,\calO_{\wwhat\Sigma}(\nu^*u^*T\CP^2)) \\
\to\oplus_cH^0(\Sigma,T_{u(c)}\CP^2)\to H^1(\Sigma,
  \calO_\Sigma(u^*T\CP^2))\\
\to H^1(\wwhat\Sigma,\calO_{\wwhat\Sigma}(\nu^*u^*T\CP^2))
  \to\oplus_cH^1(\Sigma,T_{u(c)}\CP^2)=0.
\end{multline*}
Notice that the fibers of ${\mathfrak{def}}(u)$ and
$\mathfrak{ob}(u)$  over $[u,\Sigma,p]$ are $H^0(\Sigma,
\calO_\Sigma(u^*T\CP^2))$ and $H^1(\Sigma, \calO_\Sigma(u^*T\CP^2))$
respectively. There are analogous sequences associated to any point
$[u^\prime,\Sigma^\prime,p^\prime]$ in the moduli space and the
spaces in these sequences glue together to define $T$--equivariant
vector bundles over the moduli space. This is covered in more detail
in \fullref{vc} below. The sequence of linear spaces generalizes
to a $T$--equivariant sequence of vector bundles. This implies that
$$\frac{e({\mathfrak{def}}(u)^{\text{mov}})}{e(\mathfrak{ob}(u)^{\text{mov}})}
=
\frac{e(H^0(\wwhat\Sigma,\calO_{\wwhat\Sigma}(\nu^*u^*T\CP^2))^{\text{mov}})}
{e(H^1(\wwhat\Sigma,\calO_{\wwhat\Sigma}(\nu^*u^*T\CP^2))^{\text{mov}})
\prod_ce(T_{u(c)}\CP^2)}.$$
Recall that we are computing equivariant Euler classes here. It
turns out that the possibility of a nontrivial bundle with
nontrivial action does not arise in our computations. Thus, there
are two different cases that we need to analyze. The first case is
when the bundle is possibly nontrivial and the action is trivial.
In this case the equivariant Euler class is just the usual Euler
class. The second case is when the bundle is trivial and the action
is nontrivial.  In this case, the action on a vector space induces
a bundle over $ET\times_T V$ and the equivariant Euler class is just
the product of the Chern classes associated to the irreducible
representations. Recall that $\alpha_k$ is the first Chern class of
the bundle associated to the dual representation to the projection
to the $k$th factor $\mu_k^*$.

The standard action of $T$ on $\CP^2$ may be expressed as
$\lambda\cdot(z_1,z_2)=(\lambda_1\lambda_0^*z_1,\lambda_2\lambda_0^*z_2)$
in the $(z_1,z_2)$--chart. It follows that
$$T_{[1:0:0]}\CP^2=\mu_0\otimes(\mu_1^*\oplus\mu_2^*) \quad\text{and}\quad
e(T_{[1:0:0]}\CP^2)=(\alpha_0-\alpha_1)(\alpha_0-\alpha_2),$$
so
$$\prod_ce(T_{u(c)}\CP^2))=\prod_F\prod_{j\ne k(F)}(\alpha_{k(F)}-\alpha_j).$$

\begin{remark}\label{can1}
The number of nodes at a vertex with $n(v)=g(v)=0$ is one less than
the valence of this vertex, so when we take this product over all
flags we should remove one flag from each such vertex. However, it
is easier to include all flags here and include canceling terms in
the expression for
$e(H^0(\wwhat\Sigma,\calO(u^*T\CP^2))^{\text{mov}})$ (see \fullref{can2}). We abuse notation by not introducing new notation for
these modifications to $\prod_ce(T_{u(c)}\CP^2))$ and
$e(H^0(\wwhat\Sigma,\calO(u^*T\CP^2))^{\text{mov}})$.
\end{remark}

The space
$H^0(\wwhat\Sigma,\calO_{\wwhat\Sigma}(\nu^*u^*T\CP^2))^{\text{mov}}$
splits into a sum of terms corresponding to the components of
$\wwhat\Sigma$. Clearly the bundle $u^*T\CP^2$ is trivial over
contracted components. Thus the contribution to the Euler class from
the contracted components is
$$\prod_{\val(v)+n(v)>2}\prod_{j\ne k(v)}(\alpha_{k(v)}-\alpha_j).$$
To compute the contribution to the Euler class of the non-contracted
components we will use the Euler sequence described in the next
exercise.
\begin{exm}
The subbundle of $T(\C^3-\{0\})$ generated by
$z_k\partial_k|{(z_0,z_1,z_2)}$ is invariant under the
multiplicative action of $\C^\times$. It therefore induces a bundle,
say $L$, over the quotient $\CP^2$. Check that the following is an
exact sequence of vector bundles (it is called the Euler sequence):
$$0\to L\to {\underCC}^3\to T\CP^2\to 0.$$
\end{exm}
The Euler sequence leads to the following sequence of sheaves over $\CP^2$:
$$0\to \calO_{\CP^2}\to\calO(1)\otimes{\underCC}^3\to T\CP^2\to 0$$
As a representative model of the non-contracted components, consider
the map $u\co \!\CP^1\!\to\CP^2$ given by $u([x_0:x_1])=[x_0^d:x_1^d:0]$.
Pulling back to $\CP^1$ via $u$ and taking cohomology gives
\begin{multline*}
0\to H^0(\CP^1,\calO_{\CP^1})\to H^0(\CP^1,\calO(d))\otimes{\underCC}^3\\
\to H^0(\CP^1,\calO(u^*T\CP^2))\to H^1(\CP^1,\calO_{\CP^1})\to\cdots.
\end{multline*}
Recall that $H^0(\CP^1,\calO(d))$ can be identified with the space
of degree $d$ homogeneous polynomials in $x_0$ and $x_1$. Each of
the terms in the last sequence is $T$--equivariant. Recall that the
$T$--action on $\C^3$ is given by $\lambda\cdot z=(\lambda_0
z_0,\lambda_1 z_1,\lambda_2 z_2)$ and the maps
$[z_0:z_1]\mapsto[z_0^{d-j}:z_1^j]$ generate all degree $d$
meromorphic functions on $\CP^1$. As a $T$--representation, we have
$$H^0(\CP^1,\calO(d))\otimes{\underCC}^3=
  \Biggl(\bigoplus_{j=0}^d\mu_0^{\otimes \upnfrac{d-j}{d}}
  \otimes\mu_1^{\otimes\unfrac{j}{d}}\Biggr)
  \otimes(\mu_0\oplus\mu_1\oplus\mu_2).$$
Thus taking equivariant Euler classes gives
\begin{multline*}
e(H^0(\CP^1,\calO(u^*T\CP^2)))=
  \prod_{k=0}^2\prod_{j=0}^d
  \Bigl(\frac{d-j}{d}\alpha_0+\frac{j}{d}\alpha_1-\alpha_k\Bigr)\\
= \prod_{j=1}^d\Bigl(\frac{j}{d}\alpha_1-\frac{j}{d}\alpha_0\Bigr)
  \prod_{j=0}^{d-1}\Bigl(\frac{d-j}{d}\alpha_0+\frac{j-d}{d}\alpha_1\Bigr)
  \prod_{j=0}^d\Bigl(\frac{d-j}{d}\alpha_0+\frac{j}{d}\alpha_1-\alpha_2\Bigr)\\
=\frac{(-1)^d(d!)^2}{d^{2d}}(\alpha_0-\alpha_1)^{2d}
  \prod_{j=0}^d\Bigl(\frac{d-j}{d}\alpha_0+\frac{j}{d}\alpha_1-\alpha_2\Bigr).
\end{multline*}
In general we have
\begin{multline*}
e(H^0(\wwhat\Sigma,\calO(u^*T\CP^2))^{\text{mov}})=
  \prod_{v}\prod_{j\ne k(v)}(\alpha_{k(v)} -\alpha_j)\\
\prod_e\Biggl(\frac{(-1)^{d(e)}(d(e)!)^2} {d(e)^{2d(e)}}(\alpha_{k(v)}
  {-} \alpha_{k(v^\prime)})^{2d(e)}
  \hspace{-1em}
  \prod_{\stackrel{a+b=d}{\scriptscriptstyle j\ne k(v), k(v^\prime)}}
  \hspace{-1em}
  \Bigl(\frac{a}{d(e)}\alpha_{k(v)} {+} \frac{b}{d(e)}\alpha_{k(v^\prime)}
  {-} \alpha_j\Bigr)\Biggr).
\end{multline*}

\begin{remark}\label{can2}
The first product in this expression should be over all $v$ such
that $\val(v)+n(v)+2g(v)>2$ because these are the vertices with
contracted components. However, taking the product over all vertices
exactly cancels the extra terms introduced in $\prod_ce(T_{u(c)}\CP^2))$
in \fullref{can1}.
\end{remark}

We have one remaining term to compute to finish our computation
\eqref{edefobu} of the equivariant Euler characteristic of the
normal bundle to the fixed point set, namely
$$e(H^1(\wwhat\Sigma,\calO(\nu^*u^*T\CP^2))^{\text{mov}}).$$
Notice that the Kodaira vanishing theorem implies that the cohomology
corresponding to non-contracted components vanishes, giving
\index{$V^\vee$ dual of $V$}
$$H^1(\wwhat\Sigma,\calO(\nu^*u^*T\CP^2))\cong\oplus_v
H^1(\wwhat\Sigma_v,\calO(\nu^*u^*T\CP^2)).$$
Now
\begin{align*}
H^1(\wwhat\Sigma_v,\calO(\nu^*u^*T\CP^2))&\cong
H^1(\wwhat\Sigma_v,\calO_{\Sigma_v}))\otimes T_{u(v)}\CP^2 \\
\text{and}\qquad H^1(\wwhat\Sigma_v,\calO_{\Sigma_v}))
  &\cong(H^0(\wwhat\Sigma_v,\calO_{\Sigma_v})))^\vee =\E^\vee.
\end{align*}
(Recall that $\E$ is the Hodge bundle.) Putting this together gives
\begin{align*}
H^1(\wwhat\Sigma_v,\calO(\nu^*u^*T\CP^2))&\cong \E^\vee\otimes
  T_{u(v)}\CP^2 \\
&\cong\E^\vee\otimes(\mu_1\otimes\mu_0^*\oplus\mu_2\otimes\mu_0^*) \\
&\cong\E^\vee\otimes\mu_1\otimes\mu_0^*\oplus\E^\vee\otimes\mu_2
  \otimes\mu_0^*.
\end{align*}
Here we are assuming that $u(v)=[1:0:0]$.
\begin{exm}\label{etensorl}
The splitting principle states that any formula for characteristic
classes that is valid for sums of line bundles is valid for
arbitrary bundles. Use the splitting principle to prove that
$$e(E\otimes L)=\sum_{i=0}^r c_{i}(E)c_1(L)^{r-i},$$
when $E$ is a rank $r$ vector bundle and $L$ is a line bundle.
\end{exm}
It follows that
\begin{multline*}
e(H^1(\wwhat\Sigma,\calO(\nu^*u^*T\CP^2))^{\text{mov}})= \\[-2ex]
\prod_{\val(v)+n(v)+2g(v)>2}~
  \prod_{j\ne k(v)}~
  \sum_{i=0}^{g(v)}c_i(\E^\vee)(\alpha_{k(v)}-\alpha_j)^{g-i}.
\end{multline*}

\begin{remark}
We derived the formulas for the factors of the Euler class of the
normal bundle to the components of the fixed point set for $\CP^2$.
However careful inspection shows that all of these formulas are
valid for $\CP^n$ without modification.
\end{remark}

\section{The virtual fundamental class}\label{vc}

There is one new ingredient that arises in this computation. Up to
now we have considered intersections in general  position. It is
still possible to do intersection theory when intersections are not
generic. We now describe this non-generic intersection theory.

Let $\pi\co E\to X$ be a vector bundle and let $\sigma_0$ be the zero
section. Let $F$ be a subbundle of $E$ and $\sigma\co X\to F$ be a
generic section of $F$. We may also consider $\sigma$ as a section
of $E$, but it will not be transverse to the zero section. Let
$Z:=\sigma^{-1}(\sigma_0)$ and consider the following exact sequence
of bundles,
$$0 \xrightarrow{\qquad} TZ \xrightarrow{\qquad} TX|_Z
  \xrightarrow{\qua d\sigma\qua} E|_Z \xrightarrow{\qquad} F^{\perp}|_Z
  \xrightarrow{\qquad} 0$$
Since $E$ is a vector bundle, we have a natural map $E\to TE$. We
also have a map $d\sigma_0\co TX\to TE$. One can check that
$TE|_{\sigma_0(X)}=d\sigma_0(TX)\oplus E|_{\sigma_0(X)}$. The map
labeled by $d\sigma$ in the above sequence is the projection of the
push-forward to $E|_{\sigma_0(X)}$. The bundle $F^\perp|_Z$ is
called the obstruction bundle. It is usually denoted by
$\mathfrak{ob}(Z)$. Notice that this agrees with our earlier
description of the obstruction bundle. If the section $\sigma$ were
generic then the fundamental class of $Z$ would be the Euler class
of $E$. We define the \index{$Z$@$[Z]^{\vir}$ the virtual
fundamental class of $Z$} \index{virtual fundamental class} virtual
fundamental class of $Z$ (denoted by $[Z]^{\vir}$) to be the
the Poincar\'e dual of the Euler class of $E$. We have
\begin{multline*}
~[Z]^{\vir} = \text{PD}(e(E)) = \text{PD} (e(F)\cup e(F^\perp))\\
=\text{PD}([Z]\cup e(F^\perp)) = [Z]\cap e(F^\perp) =[Z]\cap
e({\mathfrak{ob}}(Z)).
\end{multline*}
We now consider an example of the virtual fundamental class.
\begin{example}
Consider the self-intersection of a line in $\CP^2$. The usual way
to compute this is to perturb one copy of the line and then take the
intersection, but this is not necessary. Let
$L_1=\CP^3-\{[0:0:0:1]\}$ with projection $L_1\to\CP^2$ be the Chern
class $1$ line bundle over $\CP^2$. Setting $\sigma_1\co \CP^2\to L_1$
to be the section $\sigma_1([x:y:z])=[x:y:z:x]$, we see that
$\sigma_1^{-1}(\sigma_0(\CP^2))$ is just a standard line.  The
self-intersection of this line is just the zeros of the section
$\sigma=\sigma_1\oplus\sigma_1$ of $L_1\oplus L_1$. In this case it
is easy to see that $\sigma$ takes values in the diagonal $L_1$
subbundle and is transverse to the zero section of this subbundle.
This is exactly the situation described above, so we have
\begin{multline*}
\#(\CP^1\cap \CP^1)=\int_{[\sigma^{-1}(0)]^{\vir}}1 =
\int_{[\CP^1]\cap e({\mathfrak{ob}}(\CP^1))}1 \\
= \int_{\CP^1}e({\mathfrak{ob}}(\CP^1)) =
\int_{\CP^1}\Omega_{\CP^1}=1. \index{$X_{S^3}$ the resolved
conifold}\index{local $P^1$}
\end{multline*}
\end{example}

This same behavior happens in computations on many moduli spaces.
Our main example is the large $N$ dual of the $3$--sphere.
\begin{defn}
The {\it local} $P^1$ or {\it small resolution of the conifold} is
denoted by $X_{S^3}$ or $\calO(-1)\oplus\calO(-1)$. It is the total
space of the $\calO(-1)\oplus\calO(-1)$ complex vector bundle over
$\CP^1$. We have,
$$X_{S^3}:=
\biggl\{([z_0:z_1],w_0,w_1,w_2,w_3)\in \CP^1\times\C^4~\bigg|~\det
\begin{pmatrix}w_0&w_1\\z_0&z_1\end{pmatrix}=\det
\begin{pmatrix}w_2&w_3\\z_0&z_1\end{pmatrix}=0\biggr\}.$$
\end{defn}
The restriction of the symplectic form on $\CP^1\times\C^4$ is a
symplectic form on $X_{S^3}$. Let $q_1\co X_{S^3}\to \CP^1$ and
$q_2\co X_{S^3}\to\C^4$ be the projection maps. Clearly these maps
are holomorphic. Given a holomorphic map $u\co \Sigma\to X_{S^3}$ we
see that $q_2\circ u$ is holomorphic, therefore constant (see
\fullref{loui}.) Let $(w_1,w_2,w_3,w_4)$ be this constant value. For
positive degree $d$ the map $q_1\circ u$ must be surjective. Taking
$x_0\in\Sigma$ with $q_1\circ u(x_0)=[0:1]$ we see that $w_0=w_2=0$.
Taking $x_\infty\in\Sigma$ with $q_1\circ u(x_\infty)=[1:0]$ we see
that $w_1=w_3=0$.
\begin{exm}\label{loui}
Recall the definition of the $\wbar\partial$ operator from
\fullref{coarse}. Show that $\partial u$ and
$\wbar\partial u$ are perpendicular, and that
$$(|\partial u|^2-|\wbar\partial
u|^2)\,d\text{vol}_\Sigma=-g(du,J\circ du\circ
j))\,d\text{vol}_\Sigma=2u^*\omega_X.$$
Conclude that
$$\frac12\int_\Sigma|du|^2\,d\text{vol}_\Sigma=\int_\Sigma|\wbar\partial
  u|^2\,d\text{vol}_\Sigma+\int_\Sigma u^*\omega_X.$$
When $u$ is holomorphic, $\Sigma$ is closed and $u^*\omega_X$ is
exact this implies that $u$ is constant.
\end{exm}

Let $\sigma_0\co \CP^1\to X_{S^3}$ be the zero section. Our previous
discussion implies that the induced map
$$\sigma_0\co \wwbar{\mathcal M}_{g,0}(\CP^1,d[\CP^1])\to
\wwbar{\mathcal M}_{g,0}(X_{S^3},d[\CP^1])$$
is an isomorphism of stacks. According to \fullref{virdim},
$$\virdim(\wwbar{\mathcal M}_{g,0}(X_{S^3},d[\CP^1]))=0,$$
and
$$
\virdim(\wwbar{\mathcal
M}_{g,0}(\CP^1,d[\CP^1]))=2g-2+2d\,.
$$
We conclude that if the moduli space $\wwbar{\mathcal
M}_{g,0}(\CP^1,d[\CP^1])$ is unobstructed, then we are in exactly
the situation described by the excess intersection formula. In
\fullref{locstr}, we claimed that homological algebra would
allow one to extend the spaces from the deformation-obstruction
complex to bundles over the moduli space. This is similar to the
situation encountered earlier of glueing the spaces
$T_{p_k}^*\Sigma$ into the bundle ${\mathcal L}_k$ over the moduli
space. We use the universal curve over the moduli space here as we
did before. This time we need a construction from homological
algebra called higher direct image functors.

Given a left exact functor $F\co {\mathcal C}\implies{\mathcal D}$
and an injective resolution $A\to I^0\to I^1\to\cdots$ one defines
the right derived functors \index{right derived functors}
\index{$R^*(F)(A)$ right derived functors} of $F$ applied to $A$
(denoted $R^*(F)(A)$) to be the cohomology of the complex $F(I^0)\to
F(I^1)\to\cdots$. Given a map $f\co X\to Y$ and a sheaf ${\mathcal
A}$ over $X$ one defines the direct image sheaf
\index{$F$@$f_*{\mathcal A}$ push-forward} \index{direct
image}\index{higher direct image} $f_*{\mathcal A}$ over $Y$ by
$f_*{\mathcal A}(V):={\mathcal A}(f^{-1}(V))$. The \emph{higher direct
image functors} are just the right derived functors of the direct
image functor. \emph{Hyper-derived functors} are a generalization of derived
functors applicable to complexes of sheaves. One takes injective
resolutions of each sheaf in the complex and applies the functor to
obtain a double complex. The hyper-derived functor is just the
cohomology of the corresponding total complex. \emph{Hyper-higher direct
image functors} \index{hyper-higher direct image functors} ${\bf
R}^*(f_*)$ of complexes \index{$R$@$\mathbf{R}^*(f_*)$} may be defined via
the total space of the resulting complexes in the same way one defines
the hyper-ext functors.

We will apply the hyper-higher direct image functor to sheaves
arising from vector bundles. To any vector bundle we can associate
the sheaf of sections. The sheaf of sections of a holomorphic vector
bundle over $X$ is a locally free, finite rank sheaf of ${\mathcal
O}_X$--modules. Conversely, given a locally free, finite rank sheaf
${\mathcal E}$ of ${\mathcal O}_X$--modules one defines an
associated holomorphic vector bundle. Namely, there is an open cover
$\{U_\alpha\}$ of $X$ and isomorphisms $\varphi_\alpha\co {\mathcal
O}_X(U_\alpha)^n\to{\mathcal E}(U_\alpha)$. Define
$\psi_{\alpha\beta}:=\varphi_\alpha^{-1}\circ\varphi_\beta$ on the
overlaps and $E:=\perp\!\!\perp U_\alpha\times \C^n/\sim$, where
$(\alpha,x,z)\sim (\beta,x,\psi_{\beta\alpha}(x)z)$.
\begin{exm}\label{risbndl}
Show that $\Gamma(E)$ is naturally isomorphic to ${\mathcal E}$.
Show that
$${\mathcal O}_X^n\otimes_{{\mathcal O}_X}
  \left({\mathcal O}_{x_0}/{\mathfrak m}_{x_0}\right)\cong \C^n.$$
\end{exm}
We conclude from this exercise that the fiber of $E$ over a point
$x_0$ may be identified with ${\mathcal E}\otimes_{{\mathcal
O}_X}\left({\mathcal O}_{x_0}/{\mathfrak m}_{x_0}\right)$.

Now consider the universal curve over the moduli space. We have the
vertical bundle
$$\mathcal{V} \xrightarrow{\qquad} \mathcal{U}_X
  \xrightarrow{\qua\pi_X\qua} \wwbar{\mathcal{M}}_{g,n}(X,\beta)$$
with sections $\rho_k\co \wwbar{\mathcal M}_{g,n}(X,\beta)\to
{\mathcal U}_X$ and the diagram
$$\bfig
  \barrsquare<500,500>[\text{ev}_X^*TX`TX`\mathcal{U}_X`X;
  ```\text{ev}_X]
  \morphism<0,-500>[\phantom{\mathcal{U}_X}`
    \wwbar{\mathcal{M}}_{g,n}(X,\beta);\pi_X]
  \efig$$
Let ${\mathcal L}_{-\rho}$ be the line bundle associated to the
divisor
$$
-\rho_1(\wwbar{\mathcal
M}_{g,n}(X,\beta))-\cdots-\rho_n(\wwbar{\mathcal
M}_{g,n}(X,\beta))\,.
$$
Let $s$ be a section of this bundle with simple zeros along the
divisor. Define ${\mathcal V}_{-\rho}$ to be the bundle ${\mathcal
V}\otimes {\mathcal L}_{-\rho}$ and a bundle map
$\text{ev}_*(-\otimes s)\co {\mathcal V}_{-\rho}\to\text{ev}_X^*TX$.
\begin{defn}
The bundles in the deformation-obstruction sequence of a stable map
are naturally associated to the sheaves given by:
\begin{align*}
{\mathfrak{def}}(u) &:= {\bf R}^0(\pi_*)(\text{ev}^*TX) \\
{\mathfrak{ob}}(u) &:= {\bf R}^1(\pi_*)(\text{ev}^*TX) \\
{\mathfrak{aut}}([\Sigma,p]) &:= {\bf R}^0(\pi_*)({\mathcal V}_{-\rho}) \\
{\mathfrak{def}}([\Sigma,p]) &:= {\bf R}^1(\pi_*)({\mathcal V}_{-\rho}) \\
{\mathfrak{aut}}([u,\Sigma,p]) &:= {\bf R}^0(\pi_*)({\mathcal V}_{-\rho}\to\text{ev}^*TX) \\
{\mathfrak{def}}([u,\Sigma,p]) &:= {\bf R}^1(\pi_*)({\mathcal V}_{-\rho}\to\text{ev}^*TX) \\
{\mathfrak{ob}}([u,\Sigma,p]) &:= {\bf R}^2(\pi_*)({\mathcal
V}_{-\rho}\to\text{ev}^*TX) \,.
\end{align*}
\end{defn}

Before explaining why these bundles over the moduli stack have the
required fibers we should explain exactly what the various maps
actually are in this setting. In the above description we just
treated stacks as spaces. The easiest way to understand a stack at
this point is to use the definition that the moduli stack is the
contravariant functor from the category of schemes to sets that
associates the set of all equivalence classes of families of stable
maps over a scheme to a scheme.
\begin{quote} In short, a stack is just the collection of all families of stable
maps over a scheme.
\end{quote}

Several observations will clarify a correct way to think about
these constructions. The first observation is that the coarse moduli
space is the set that the moduli stack associates to a point, that
is, $\wwbar{\mathcal
M}_{g,n}(X,\beta)(\text{pt})=\wwbar{M}_{g,n}(X,\beta)$. The second
observation is that the universal curve over \index{universal curve}
\index{$U$@${\mathcal U}_X$ the universal curve} $\wwbar{\mathcal
M}_{g,n}(X,\beta)$ is just $\mathcal{U}_X=\wwbar{\mathcal
M}_{g,n+1}(X,\beta)$ even when the stable maps have nontrivial
automorphisms. The map
$$
\pi_X\co {\mathcal U}_X \to \wwbar{\mathcal M}_{g,n}(X,\beta)
$$
is just the natural transformation of functors that takes a family
of $(n{+}1)$--pointed curves  to the family of $n$--pointed curves
obtained by ignoring the last point (section) and stabilizing the
fibers. The third observation is that the stack associated to a
space just consists of all maps from schemes into the space, and the
map $\text{ev}_X\co \mathcal{U}_X\to X$ is just the natural
transformation that takes a family $[w\co \mathcal{W}\to
X,\mathcal{W}\to S, \rho]$ to the map $w\circ \rho_{n+1}\to X$. The
example in \fullref{mc} will continue with the idea that the
moduli stack is the collection of families of stable maps. The
following exercise is good practice translating constructions into
families.
\begin{exm}
Give a definition of the vertical bundle
$$\mathcal{V} \xrightarrow{\qquad} \mathcal{U}_X
  \xrightarrow{\qua\pi_X\qua} \wwbar{\mathcal{M}}_{g,n}(X,\beta)$$
by making $\mathcal{V}$ a stack that associates appropriate families
of vector bundles over families of stable maps.
\end{exm}

We now describe why the fibers of the bundles in the
deformation-obstruction sequence are the expected spaces. To apply
\fullref{risbndl} to the sheaves of the deformation-obstruction
complex, we need a theorem of Grauert (see Hartshorne \cite[page
33]{hart}). Recall the definition of flat morphism from \fullref{calm}.
\begin{thm}
If $f\co X\to Y$ is a flat morphism,  ${\mathcal F}$ is a coherent
sheaf over $Y$ and $\text{dim}_{{\mathcal O}_{y}/{\mathfrak
m}_{y}}H^k(X_y,{\mathcal F}_y)$ is constant, then
$R^k(f_*)({\mathcal F})$ is locally free of finite rank and
$$
R^k(f_*)({\mathcal F})\otimes_{{\mathcal O}_Y}\left({\mathcal
O}_{y}/{\mathfrak m}_{y}\right)\cong H^k(X_y,{\mathcal F}_y)\,,
$$
where $X_y$ is the fiber over $y$ and ${\mathcal F}_y={\mathcal
F}|_{X_y}$.
\end{thm}
\begin{quote}
This theorem gives us a good way to think about higher direct image
functors -- under nice conditions the higher direct image functors
associate to a family of spaces over $S$ a vector bundle over $S$
with fiber equal to the cohomology of the fiber in the original
family. \index{higher direct image functor}
\end{quote}

Notice that the condition that $H^k(X_y,{\mathcal F}_y)$ have
constant dimension is not satisfied in cases of interest to us. For
example, a degree three map from a surface of genus three to
$X_{S^3}$ has no nontrivial infinitesimal automorphisms of the
underlying surface. However there is a nodal surface in the same
moduli space consisting of the one point union of a surface of genus
two and a surface of genus zero mapped by degree two on the genus
two part and degree one on the genus zero part. This nodal surface
has a two-complex-dimensional space of infinitesimal automorphisms.
The next exercise compares the virtual fibers of the bundles from
the deformation-obstruction bundle complex with the spaces from the
deformation-obstruction complex.
\begin{exm}
Given an injective resolution ${\mathcal O}_\Sigma\to I^0\to
I^1\to\cdots$, show that $u^*TX\to\Hom_{{\mathcal
O}_\Sigma}(u^*\Omega_X,I^0)\to\cdots$ is an injective resolution of
$u^*TX$. Conclude that
$$
R^1(\pi_*)(\text{ev}^*TX)\otimes_{{\mathcal O}/{\mathfrak
m}}\left({\mathcal O}/{\mathfrak m}\right)\cong H^1(\Sigma,u^*TX)=\E
xt^1(u^*\Omega_X,{\mathcal O}_\Sigma)\,.
$$
Repeat this computation with the other bundles.
\end{exm}

Assume for the moment that the fiber
$$
\Ext^2(u^*\Omega_X,{\mathcal O}_\Sigma)\cong {\mathbb
H}^2(\Sigma,T\Sigma\to u^*TX)
$$
has constant dimension as $[u,\Sigma]$ varies in $\wwbar{\mathcal
M}_{g,0}(X_{S^3},d[\CP^1])$. (Here ${\mathbb H}^k$ represents
hypercohomology.) Then ${\mathfrak{ob}}([u,\Sigma])={\bf
R}^2(\pi_*)({\mathcal V}\to \text{ev}^*TX)$ satisfies the
assumptions of Grauert's theorem and therefore satisfies the
assumptions required for our definition of the virtual fundamental
class. Using the map $\sigma_0\co \wwbar{\mathcal
M}_{g,0}(\CP^1,d[\CP^1])\to\wwbar{\mathcal
M}_{g,0}(X_{S^3},d[\CP^1])$, we write
\begin{align*}
N_{g,d}:=\langle 1\rangle_{g,d[\CP^1]}^X &=
\int_{[\wwbar{\mathcal
M}_{g,0}(X_{S^3},d[\CP^1])]^{\text{\begin{tiny}vir\end{tiny}}}} 1 \\
&=\int_{[\wwbar{\mathcal M}_{g,0}(X_{S^3},d[\CP^1])]}e({\bf
R}^2(\pi_*)({\mathcal V}\to \text{ev}^*TX)) \\
&=\int_{[\wwbar{\mathcal
M}_{g,0}(\CP^1,d[\CP^1])]}\sigma_0^*e({\bf R}^2(\pi_*)({\mathcal
V}\to \text{ev}^*TX))\\
&=\int_{[\wwbar{\mathcal
M}_{g,0}(\CP^,d[\CP^1])]}e(\sigma_0^*{\bf R}^2(\pi_*)({\mathcal
V}\to \text{ev}^*TX))\,.
\end{align*}
We continue with the change of base theorem for higher direct image
functors. See Hartshorne \cite[page 255]{hart}.
\begin{thm}
Given a commutative diagram,
$$\bfig\barrsquare[W`X`Z`Y;v`g`f`u]\efig$$
with $u$ flat and a complex of coherent sheaves ${\mathcal A}$ one
has
$$
u^*{\bf R}^*(f_*)({\mathcal A})={\bf R}^*(g_*)(v^*{\mathcal A})\,.
$$
\end{thm}
We conclude that
$$
\sigma_0^*{\bf R}^2(\pi^X_*)({\mathcal V}\to \text{ev}_X^*TX)={\bf
R}^2(\pi^{\CP^1}_*)(\what\sigma_0^*{\mathcal V}\to
\what\sigma_0^*\text{ev}_X^*TX)\,.
$$
Here $\what\sigma_0\co {\mathcal U}_{\CP^1}\to{\mathcal U}_X$ is the
natural map. Since
$\sigma_0\circ\text{ev}_{\CP^1}=\text{ev}_X\circ\what\sigma_0$,
we know
$${\bf R}^2(\pi^{\CP^1}_*)(\what\sigma_0^*{\mathcal V}\to
\what\sigma_0^*\text{ev}_X^*TX)={\bf
R}^2(\pi^{\CP^1}_*)(\what\sigma_0^*{\mathcal V}\to
\text{ev}_{\CP^1}^*\sigma_0^*TX).$$
To go further, consider the exact sequence of bundles,
$$
0\to \text{ev}_{\CP^1}^*T\CP^1\to \text{ev}_{\CP^1}^*\sigma_0^*TX\to
\text{ev}_{\CP^1}^*({\mathcal O}(-1)\oplus{\mathcal O}(-1))\to 0\,.
$$
By the horseshoe lemma from homological algebra there are injective
resolutions $\text{ev}_{\CP^1}^*T\CP^1\to I^*$,
$\text{ev}_{\CP^1}^*\sigma_0^*TX\to J^*$ and
$\text{ev}_{\CP^1}^*({\mathcal O}(-1)\oplus{\mathcal O}(-1))$ that
form a short exact sequence of complexes
$$0\longrightarrow I^*\longrightarrow J^*\longrightarrow K^*\longrightarrow 0.$$

\begin{exm}
Let $0\to A\to B\to C\to 0$ be a short exact sequence of modules,
and let $A\to I^*$ and $C\to K^*$ be injective resolutions. Show
that there are maps $\varepsilon\co B\to I^0\oplus K^0$ and
$d^k\co I^k\oplus K^k \to I^{k+1}\oplus K^{k+1}$ such that $B\to
I^*\oplus K^*$ is an injective resolution that fits into a short
exact sequence of complexes  (see Weibel \cite{wei}).
\end{exm}
Let ${\mathcal V}\to L^*$ be an injective resolution and form the
following short exact sequence of complexes,
$$\bfig\iiixiii<700,500>{'0077}[L^0`L^0`0`
  L^1\oplus I^0`L^1\oplus J^0`K^0`
  I^1`J^1`K^1;```````````]\efig$$
Applying the direct image functor $\pi^{\CP^1}_*$ to each term and
writing out the associated long exact sequence of cohomology groups
produces a long exact sequence containing
\begin{multline*}
\to {\bf R}^0(\pi_*)\bigl(\text{ev}_{\CP^1}^*({\mathcal
O}(-1)\oplus{\mathcal O}(-1))\bigr)\to {\bf
R}^2(\pi_*)\bigl(\text{ev}_{\CP^1}^*{\mathcal
V}\to\text{ev}_{\CP^1}^*T\CP^1\bigr) \\
\to {\bf R}^2(\pi_*)\bigl(\text{ev}_{\CP^1}^*{\mathcal
V}\to\text{ev}_{\CP^1}^*\sigma_0^*TX\bigr)\to {\bf
R}^1(\pi_*)\bigl(\text{ev}_{\CP^1}^*({\mathcal O}(-1)\oplus{\mathcal
O}(-1))\bigr)\to 0.
\end{multline*}
Now look at the deformation-obstruction sequence for
$\wwbar{\mathcal M}_{g,0}(\CP^1,d[\CP^1])$,
\begin{multline*}
0\longrightarrow  {\mathfrak{aut}}([v,\Sigma]) \longrightarrow
  {\mathfrak{aut}}([\Sigma]) \longrightarrow {\mathfrak{def}}(v)
  \longrightarrow {\mathfrak{def}}([v,\Sigma]) \longrightarrow
  {\mathfrak{def}}([\Sigma]) \\
\longrightarrow {\mathfrak{ob}}(v) \longrightarrow
  {\mathfrak{ob}}([v,\Sigma]) \longrightarrow 0.
\end{multline*}
The fiber of ${\mathfrak{ob}}(v)$ is $H^1(\Sigma,v^*T\CP^1)$.  For
every stable map in the genus zero case this cohomology group is
zero (this property is called convexity). For $\CP^1$ convexity
follows from the Kodaira vanishing theorem (see Griffiths and Harris
\cite{GH}). This implies that
$$0={\mathfrak{ob}}([v,\Sigma])={\bf
R}^2(\pi_*)\bigl(\text{ev}_{\CP^1}^*{\mathcal
V}\to\text{ev}_{\CP^1}^*T\CP^1\bigr),$$
so that
$${\bf R}^2(\pi_*)\bigl(\text{ev}_{\CP^1}^*{\mathcal
V}\to\text{ev}_{\CP^1}^*\sigma_0^*TX\bigr)\cong {\bf
R}^1(\pi_*)\bigl(\text{ev}_{\CP^1}^*({\mathcal O}(-1)\oplus{\mathcal
O}(-1))\bigr),$$
by the previous exact sequence. Now,
$$\text{dim}_\C H^1(\Sigma,v^*({\mathcal O}(-1)\oplus{\mathcal
O}(-1)))=2d+2g-2$$
by the Riemann--Roch theorem together with the Kodaira vanishing
theorem. Since this is constant, we conclude that the excess
intersection formula holds and
$$
N_{g,d}=\int_{[\wwbar{\mathcal
M}_{g,0}(\CP^1,d[\CP^1])]}c_{2d+2g-2}\bigl({\bf
R}^1(\pi_*)\bigl(\text{ev}_{\CP^1}^*({\mathcal O}(-1)\oplus{\mathcal
O}(-1))\bigr)\bigr).
$$
The argument we used to get here was
valid for the $g=0$ case. In fact, one can show that this formula is
valid for all genera. When the dimension of the virtual fibers of
${\mathfrak {ob}}([u,\Sigma,p])$ is not constant, one can still
define a virtual fundamental class. This was done independently and almost simultaneously
in the spring of 1996 by Fukaya--Ono, Hofer--Salamon, Li--Tian, Ruan and Siebert,
see Salamon's lectures \cite{salamon},5 and references therein.
All of the authors originally worked in a symplectic setting.
Our description is closest to Li--Tian \cite{LT} who later extended their
results to an algebraic setting \cite{litian}. See also Behrend \cite{beh2}, Cox and Katz
\cite{CK} and Liu \cite{liu}.

Virtual fundamental classes are important in the computation of
degree zero Gromov--Witten invariants as well. The next exercises
address this situation.
\begin{exm}
Show that the obvious map $\sigma_0\co \wwbar{\mathcal
M}_{g,n}\times X\to \wwbar{\mathcal M}_{g,n}(X,0)$ is an
isomorphism of stacks.
\end{exm}

\begin{exm}
Let $\pi_M\co \wwbar{\mathcal M}_{g,0}\times X \to
\wwbar{\mathcal M}_{g,n}$ and $\pi_X\co \wwbar{\mathcal
M}_{g,0}\times X \to X$ be the projection maps and show that
$$
\int_{[\wwbar{\mathcal M}_{g,n}(X,0)]^{\vir}}\gamma
=\int_{[\wwbar{\mathcal M}_{g,n}]^{\vir}\times
X}c_{\text{top}}(\pi_M^*\E^\vee\otimes\pi_X^*TX)\cup\sigma_0^*\gamma\,.
$$
\end{exm}

\begin{exm}\label{xdeg0}
Use the splitting principle to show that $\text{dim}_\C X=3$ and
$c_1(TX)=0$ imply
$$
N_{g,0}(X):=\int_{[\wwbar{\mathcal M}_{g,0}(X,0)]^{\vir}} 1
=(-1)^g\chi(X) \int_{[\wwbar{\mathcal
M}_{g,0}]^{\vir}}c_{g-1}(\E)^3/2\,.
$$
\end{exm}

\section{The multiple cover formula in degree two}\label{mc}

In this subsection we describe how to compute the Gromov--Witten
invariants of the manifold $X_{S^3}$. We begin with a direct
computation of
$$N_{0,2}(X):=\langle 1\rangle^X_{0,2[\CP^1]}.$$
A localization computation of this same number is presented in
Cox and Katz \cite{CK}. We begin by analyzing the corresponding coarse
moduli space.

Recall from the previous subsection that any stable map into
$X_{S^3}$ factors through $\CP^1$, so the moduli stack of stable
maps to $X_{S^3}$ is isomorphic to the moduli stack of stable maps
to $\CP^1$. Given a stable map, one can construct a graph with
vertices corresponding to the maximal contracted components, edges
corresponding to the non-contracted components, and labels
corresponding to the marked points, images of the contracted
components, genera of the contracted components, and degrees of the
non-contracted components. Since the curves in our case have genus
zero, the corresponding graphs must be trees. Computing the Euler
characteristic of the resulting graph gives $1=\sum_v
(1-\frac12\text{valence}(v))$. It follows that we must either have a
vertex of valence zero or two vertices of valence one. In order to
be stable the map must have positive degree on each vertex with
valence less than three. It follows that the only stable curves in
the genus zero degree two moduli space have domain $\CP^1$ or two
copies of $\CP^1$ joined by a single node.

Now consider a degree two holomorphic map $u\co \CP^1\to\CP^1$. Locally
such a map has a power series representation, so it must be a
branched cover. This implies that
$$
2=\chi(\CP^1)=2\chi(\CP^1)-\sum_{y\in
S(u)}(2-|u^{-1}(y)|)=4-|S(u)|\,.
$$
It follows that such a map must have exactly two critical points and
two critical values. We will use the critical values to parametrize
these maps. Notice that pre-composition with a linear fractional
transformation cannot change the locations of the critical values,
so two maps with different critical values are different. Given any
two distinct points $p$, $q$ in $\CP^1$ we can take a linear
fractional transformation taking $[0:1]$ to $p$ and $[1:0]$ to $q$.
Composing the map $[z:w]\mapsto [z^2:w^2]$ with this linear
fractional transformation gives a degree two map with the desired
critical values. Notice that this is not an equivalence of stable
maps because such equivalences must be by pre-composition. We will
map a nodal curve to the image of the node.
\begin{exm}
Show that two stable maps with the same critical values are
equivalent, and show that the automorphism group of any stable curve
in this space is $\Z_2$. Conclude that the coarse moduli space is
isomorphic to the symmetric product of two copies of $\CP^1$. (This
is a reasonable example to use to understand the Gromov topology, or
the algebraic structure of the coarse moduli space.)
\end{exm}

The second symmetric power of $\CP^1$ may be identified with the
space of degree two polynomials up to scale. One associates the
roots of the polynomial to the polynomial. This gives the
isomorphism $\text{Sym}^2\CP^1\to \CP^2$ taking
$([z_0:z_1],[w_0:w_1])$ to $[z_1w_1:-z_0w_1-z_1w_0:z_0w_0]$. The
same map is an explicit isomorphism
$$
\Phi\co {\wwbar{M}}_{0,0}(\CP^1,2[\CP^1])\cong \CP^2\,,
$$
when the singular set of the map is $([z_0:z_1],[w_0:w_1])$. The
boundary divisor is just $D=\{[a:b:c]\in\CP^2|b^2-4ac=0\}$.

At the level of stacks the universal curve is just the moduli space
${\wwbar {\mathcal M}}_{0,1}(\CP^1,2[\CP^1])$ with
evaluation as the map to $\CP^1$ and projection as the map to
${\wwbar {\mathcal M}}_{0,0}(\CP^1,2[\CP^1])$. It is interesting
to compute the corresponding coarse moduli space. We claim that
${\wwbar{M}}_{0,1}(\CP^1,2[\CP^1])$ is isomorphic to
$\CP^2\times\CP^1$ with isomorphism taking $[u,\Sigma,p]$ to
$(\Phi(u),u(p))$. Given a point in $\CP^2\times\CP^1$ one can take a
stable map corresponding to the $\CP^2$--component and then pick a
point in the inverse image of the $\CP^1$--component as the marked
point. Such a point exists because the map has degree two. (If the
inverse image is a node, we add a ghost bubble containing the marked
point at the node.) If the $\CP^1$--component is one of the critical
values there is a unique choice for the marked point. Otherwise
there are two possibilities. However one finds that the resulting
marked stable curves are equivalent. For example, the point
$[0:1:0],[1:1]$ gives a stable map with critical points $[0:1]$ and
$[1:0]$ (the roots of $bz_0z_1=0$.) This map is just $[z:w]\mapsto
[z^2:w^2]$. The inverse image of $[1:1]$ is just $[\pm 1:1]$. The
reparametrization $[z:w]\mapsto [-z:w]$ takes one to the other. Most
points in ${\wwbar{M}}_{0,1}(\CP^1,2[\CP^1])$ have trivial
automorphism group. This is because the underlying stable map has
automorphism group $\Z_2$ and only the trivial automorphism will fix
the marked point unless the marked point is at a critical value. It
follows that the points along the divisor
$$
D_1=\{[a:b:c],[z:w]|az^2+bzw+cw^2=0\}
$$
have  automorphism group $\Z_2$.

For the universal curve ${\mathcal U}\to {\mathcal M}$ one should
have that the inverse image of a point $s_0\in{\mathcal M}$ is
isomorphic to $s_0$. This fails with the coarse moduli spaces. The
inverse image of a point in $\CP^2$ is a copy of $\CP^1$ and the
restriction of the evaluation map to this copy is just a degree one
map. This should be a degree two map. The reason for this failure is
that each map in this moduli space has automorphism group $\Z_2$. It
appears that one should take a double cover of $\CP^2\times\CP^1$
branched along the divisor $D_1$. The fiber of such a cover over a
point in $\CP^2$ would be a two-fold branched cover of $\CP^1$
branched over two points. This is also a copy of $\CP^1$, but the
induced evaluation map would have degree two as it should. The only
problem with this is that no such cover exists in the category of
schemes - such a map would restrict to a two-fold connected cover
over a simply-connected space. The map
$$
{\wwbar {\mathcal M}}_{0,1}(\CP^1,2[\CP^1])\to {\wwbar
{\mathcal M}}_{0,0}(\CP^1,2[\CP^1])
$$
has exactly this structure in the category of stacks.

To go further we have to use the power of stacks. Recall that one
description of a Deligne--Mumford stack is a contravariant functor
from SCHEME to SET. Any scheme produces such a functor. For $\CP^2$
one gets the functor $\underCP^2$ taking a scheme $S$ to
$\text{Mor}(S,\CP^2)$. To analyze the local symmetry groups we need
to use the fibered category structure. The objects in the associated
fibered category are ordered pairs consisting of an element of the
set associated to a scheme and the scheme. The points in the stack
are just those objects corresponding to the one point scheme. For
example,
$$
\llbracket 0:1:0\rrbracket :=([0:1:0]:\text{pt}\to\CP^2,\text{pt})\in\text{Ob}({\mathbf
D}^{\underCP^2})\,,
$$
has $\Hom(\llbracket 0:1:0\rrbracket ,\llbracket 0:1:0\rrbracket )=\{\text{id}\}$, so the local
automorphism group of a point in the stack associated to $\CP^2$ is
trivial as expected.

Let ${\mathcal M}:={\wwbar {\mathcal M}}_{0,0}(\CP^1,2[\CP^1])$.
The argument showing that the corresponding coarse moduli space is
isomorphic to $\CP^2$ shows that the points of this stack correspond
to the points of $\CP^2$. Consider the local automorphisms of a
point of this stack. Let
$$
z^2:=(u([z:w])=[z^2:w^2]\co \CP^1\to\CP^1,\pi\co \CP^1\to\text{pt},\text{pt})\in\text{Ob}({\mathbf
D}^{\underMM})\,.
$$
We have $\Hom(z^2,z^2)=\{\text{id},n\}$, where
$n([z:w])=[-z:w]$. A similar thing is true for every point in
$\mathcal M$, so every point in this stack has automorphism group
$\Z_2$.
\begin{exm}
Do a similar analysis in ${\wwbar {\mathcal
M}}_{0,1}(\CP^1,2[\CP^1])$ to show that the points in $D_1$ have
automorphism group $\Z_2$ and all others have trivial automorphism
group.
\end{exm}

We can interpret the fundamental cycle of the stack $\mathcal M$ to
be $\frac12[{\wwbar {M}}_{0,0}(\CP^1,2[\CP^1])]$
or just $\frac12[\CP^2]\in H_4(\CP^2;{\mathbb Q})$. The factor of
$1/2$ here is due to the $\Z_2$ automorphism group. To compute the
virtual fundamental class, we need to compute
$$c_{2d+2g-2}({\bf R}^1(\pi_*)(\text{ev}_{\CP^1}^*({\mathcal
O}(-1)\oplus{\mathcal O}(-1)))).$$
In this case $2d+2g-2=2$, so the Whitney sum formula gives
$$
c_{2}({\bf R}^1(\pi_*)(\text{ev}_{\CP^1}^*({\mathcal
O}(-1)\oplus{\mathcal O}(-1))))=c_{1}({\bf
R}^1(\pi_*)(\text{ev}_{\CP^1}^*({\mathcal O}(-1))))^2\,.
$$
Attempting to construct a two-fold cover of $\CP^2\times \CP^1$
branched along $D_1$ provides good motivation for the computation of
the above first Chern class. In general to construct a $p$--fold
cyclic cover branched along a divisor $D$, one constructs the line
bundle associated to $D$, say $L_D$, and takes a section vanishing
along $D$, say $\sigma_D$. If $L^{1/p}$ is a $p$th root of this
bundle in the sense that $(L^{1/p})^{\otimes p}\cong L_D$, the
desired cover will be
$$
\{\xi\in L^{1/p}|\xi^2=\sigma_D(\pi(\xi))^p\}\,.
$$
Applying this idea to the divisor $D_1$ in $\CP^2\times\CP^1$, we
cover $\CP^2\times\CP^1$ by charts $V_a$, $V_b$, and $V_c$
corresponding to $a=1$, $b=1$ and $c=1$. Define a bundle
$L_c^{1/2}=\C^3\times (\C^2-\{0\})/\sim$ where
$(a,b,\gamma,z_0,z_1)\sim(a, b,\lambda \gamma,\lambda z_0,\lambda
z_1)$, and similar bundles over $V_a$ and $V_b$. We have a section
of the tensor square of this bundle taking $(a,b,[z_0:z_1])$ to
$[a,b,az_0^2+bz_0z_1+z_1^2,z_0,z_1]$. The problem is that there is
no reasonable way to glue these pieces into a global bundle.
Continue anyway and define
$$
Q_c:=\{[a,b,\gamma,z_0,z_1]\in
L_c^{1/2}|\gamma^2=az_0^2+bz_0z_1+z_1^2\}\,
$$
with natural projection map $\pi\co Q_c\to\C^2$ and $v\co Q_c\to\CP^1$
given by
$$v([a,b,\gamma,z_0,z_1])=[z_0:z_1].$$

\begin{exm}
Show that $\pi$ is a flat morphism.
\end{exm}
We can see that $Q_c$ is a flat family of stable, genus zero, degree
two maps to $\CP^1$. We just need to check that the restriction of
$v$ to the inverse image of any point in $\C^2$ is such a stable
map. For example, we have an isomorphism $\CP^1\to\pi^{-1}(1,0)$
given by $[s:t]\mapsto [1,0,s^2-t^2,2st,s^2+t^2]$. The composition
of this map with the restriction of $v$ is a degree two map with
critical values at $[\pm i:1]$ as expected.
\begin{exm}
Identify the restriction of $v$ to $\pi^{-1}(0,0)$.
\end{exm}

Stacks should be considered as generalizations of schemes that
include orbifold information. Just as a manifold is defined via a
maximal atlas while specific manifolds are usually described by a
finite atlas, a specific stack can be described by a finite cover
while the general definition adds a condition analogous to
maximality. One can cover the stack ${\wwbar {\mathcal
M}}_{0,0}(\CP^1,2[\CP^1])$ by the family $Q_c$ together with two
other analogous families denoted by $Q_a$ and $Q_b$. \index{stack}

\begin{exm}
Construct analogous families $Q_a$ and $Q_b$ and use these three
families to conclude that ${\wwbar {\mathcal
M}}_{0,0}(\CP^1,2[\CP^1])$ is locally representable.
\end{exm}

We have been treating stacks as contravariant functors from the
category of schemes to sets, so for example the stack associated to
a scheme $T$ is the functor that takes a scheme $S$ to the set of
morphisms from $S$ to $T$. As explained in \fullref{app:a}, stacks
can also be viewed as fibered categories. The objects of the fibered
category ${\mathbf D}^{\underT}$ are just morphisms $u\co R\to T$,
similarly the objects of the fibered category version ${\mathbf
D}^{\mathcal{M}}$ of the stack ${\wwbar {\mathcal
M}}_{0,0}(\CP^1,2[\CP^1])$ are just families of stable maps $[v\co
V\to\CP^1, \pi\co V\to R]$.

Notice that any flat family $Q$ over a scheme $T$ (really,
$[q\co Q\to\CP^1, \pi\co Q\to T]$) living in ${\wwbar {\mathcal
M}}_{0,0}(\CP^1,2[\CP^1])(T)$ defines a map of stacks (covariant
functor) taking ${\mathbf D}^{\underT}$ to ${\mathbf
D}^{\mathcal{M}}$. This map takes a morphism $u\co R\to T$ in ${\mathbf
D}^{\underT}$ to the fiber product family
$[q\circ\text{pr}_1,Q\times_TR\to R]$ where $Q\times_TR:=\{(x,y)\in
Q\times R|\pi(x)=u(y)\}$ in${\mathbf D}^{\mathcal{M}}$.
\begin{exm}
Define the action of the functor on a morphism from $u\co R\to T$ to
$v\co S\to T$, that is, $w\co R\to S$ that satisfies $u=v\circ w$.
\end{exm}

\begin{quote}
An alternate way to think of stacks is as a categorical construction
of terminal objects. The moduli space should be a terminal object in
the category of families of stable maps. The problem is that
generally no such terminal object exists. When we stackify we
replace a scheme by a contravariant functor and a morphism by a
natural transformation and arrive at a category in which the
collection of all families itself turns into a terminal object.
\end{quote}

We now turn to the computation of the Chern classes in the formula
for the Gromov--Witten invariants from the last subsection. Recall
that Chern classes satisfy $c(f^*E)=f^*c(E)$, so we could compute a
Chern class of a bundle by computing the Chern class of the
pull-back under a finite branched cover and dividing by the degree
of the cover. Since $H^2(\CP^2;\Z)\cong\Z$,  to compute the first
Chern class of any bundle over $\CP^2$ it suffices to compute the
first Chern class of the restriction of the bundle to $\CP^1$. Even
though there is no bundle corresponding to the square root of the
line bundle associated to $D_1$, there is an object corresponding to
the result of taking the pull-back of the restriction to $\CP^1$ of
such a bundle under a two-fold cover. We apply this to the $\CP^1$
at $a=0$.

Define
$$\wwhat{Q}:=
  \bigl\{[\alpha,b,c,z_0,z_1]\in\C\times(\C^2-\{0\})^2~\big|~
  \alpha^2=b^2z_0z_1+c^2z_0z_1\bigr\}/\sim,$$
with $[\alpha,b,c,z_0,z_1]\sim[\lambda\mu\alpha,\lambda b,\lambda
c,\mu z_0,\mu z_1]$. This is the total space of a family with
projection $\pi\co \wwhat Q\to\CP^1$ taking $[\alpha,b,c,z_0,z_1]$ to
$[b:c]$ and evaluation $v\co \wwhat Q\to\CP^1$ given by
$v([\alpha,b,c,z_0,z_1]):=[z_0:z_1]$. Notice that the variables $b$
and $c$ here do not  correspond to the same variables as used
earlier. It might have been clearer to use $b_1$, $c_1$ here with
the relations $b=b_1^2$ and $c=c_1^2$ defining the two-fold branched
cover of the $\CP^1$ cycle. This would just complicate the notation
a bit further. The space $\wwhat Q$ is also a flat family of genus
zero, degree two stable maps. We wish to compute $c_{1}({\bf
R}^1(\pi_*)(v^*({\mathcal O}(-1))))$. Recall that ${\bf R}^1(\pi^*)$
is just a vector bundle over $\CP^1$ with fibers isomorphic to
$H^1(\pi^{-1}(-),v|^*({\mathcal O}(-1))))$. It helps to remember how
to compute sheaf cohomology at this point.
\begin{exm}
Define $L_n:=\C^2\times (\C-\{0\})/\sim$, with $(z_0,z_1,\zeta)\sim
(\lambda z_0,\lambda z_1,\lambda^n\zeta)$. Taken with the natural
projection to $\CP^1$, this is a line bundle. Let ${\mathcal O}(n)$
be the associated sheaf of sections. Using the standard $z_1\ne 0$,
$z_0\ne 0$ cover of $\CP^1$, compute the \v{C}ech cohomology groups
$\check{H}^*(\CP^1;{\mathcal O}(n))$ for various positive and
negative values of $n$.
\end{exm}

To compute ${\bf R}^1$ we use the  following description of it from
Hartshorne's book \cite{hart}.
\begin{thm}
If $\mathcal A$ is a sheaf over $X$ and $f\co X\to Y$, then
$R^k(f_*)({\mathcal A})$ is the sheaf associated to the presheaf
taking $V$ to $H^k(f^{-1}(V);{\mathcal A}|_{f^{-1}(V)})$.
\end{thm}
One can check that the $v$ pull-back of the $L_{-1}$ bundle over
$\CP^1$ to $\wwhat Q$ is the bundle defined by
$$v^*L_{-1}:=
  \bigl\{[\zeta,\alpha,b,c,z_0,z_1]\in\C^2\times(\C^2-\{0\})^2 ~|~
  \alpha^2=b^2z_0z_1+c^2z_0z_1\bigr\}/\sim,$$
with
$[\zeta,\alpha,b,c,z_0,z_1]\sim[\mu^{-1}\zeta,\lambda\mu\alpha,\lambda
b,\lambda c,\mu z_0,\mu z_1]$. Now work over the $c\ne 0$ chart
$U_c$ of $\CP^1$. We can take an open cover of the $\pi$--inverse
image of this chart consisting of $z_1\ne 0$ and $z_0\ne 0$.
Computing $R^1(\pi_*)(\mathcal{O}(v^*L_{-1}))(U_c)$ amounts to
computing the \v{C}ech cohomology of the inverse image of $U_c$. A
$0$--\v{C}ech cochain consists of an algebraic section of $v^*L_{-1}$
over the $z_1\ne 0$ chart and a section over the $z_0\ne 0$ chart. A
section over the first takes the form
$[f(\alpha,b,z_0),\alpha,b,1,z_0,1]$ where $f(\alpha,b,z_0)$ is
polynomial and we are using the obvious coordinates obtained by
setting $c$ and $z_1$ to one. A algebraic section over the $z_0\ne
0$ chart takes the form $[g(\alpha,b,z_0),\alpha,b,1,1,z_1]$ with
$g$ polynomial. Using the relation from the definition of $\wwhat
Q$ we can eliminate all second and higher powers of $\alpha$ and
write $f(\alpha,b,z_0)=f_0(b,z_0)+\alpha f_1(b,z_0)$. The polynomial
$g(\alpha,b,z_1)$ can be expressed similarly. A \v{C}ech $1$--cochain
is just a algebraic section on the overlap. The \v{C}ech coboundary
is just the difference of the restrictions of the sections from the
two large charts. In order to compute this difference we must write
each section in the same coordinates. Using the $z_1\ne 0$ chart we
can write the \v{C}ech coboundary as
$$\delta(f,g):=\bigl[(f_0(b,z_0)-z_0^{-1}g_0(b,z_0^{-1}))+
  \alpha(f_1(b,z_0)-z_0^{-2}g_1(b,z_0^{-1})),\alpha,b,1,z_0,1\bigr].$$
Notice that on the overlap $z_0\ne 0$ and $\alpha^2=bz_0+1$, so any
algebraic function on this overlap may be written in the form
$z_0^{-N}(F_0(b,z_0)+F_1(b,z_0)\alpha)$, with polynomial $F_0$ and
$F_1$.  Combining this with the expression for the coboundary
implies that the cokernel of $\delta$  is the $\C[b]$--module
generated by $z_0^{-1}\alpha$ and this is the definition of ${\bf
R}^1(\pi_*)(v^*({\mathcal O}(-1)))(\{c\ne 0\})$. This confirms our
theoretical arguments from the previous subsection that ${\bf
R}^1(\pi_*)(\text{ev}^*{\mathcal O}(-1))$ is a locally free,
finite-rank sheaf of ${\mathcal O}_S$--modules of the correct
dimension. On the $b\ne 0$ chart we will use $c$, $z_0$ ($z_1$) and
$\alpha$ with $b=1$ as our coordinates.
\begin{exm}
Show that ${\bf R}^1(\pi_*)(v^*({\mathcal O}(-1)))(\{b\ne 0\})$ is
the $\C[c]$--module generated by $z_0^{-1}\alpha$ in the given
coordinates.
\end{exm}
We must be careful when making the identifications between the $b\ne
0$ chart and the $c\ne 0$ chart. It might help to use a different
variable, say $\beta$ in place of $\alpha$ when describing the $b\ne
0$ chart. The correct way is to scale the $b=1$ answer to a $c=1$
answer using the equivalence from the definition with
$\lambda=c^{-1}$. It follows that the section of ${\bf
R}^1(\pi_*)(v^*({\mathcal O}(-1)))$ over the $c\ne 0$ chart given by
$z_0^{-1}\alpha$ extends to a meromorphic section over all of
$\CP^1$ given by $c^{-1}z_0^{-1}\alpha$ in the $b\ne 0$ chart. This
meromorphic section has exactly one simple pole, and no other poles
or zeros. It follows that $c_1({\bf R}^1(\pi_*)(v^*({\mathcal
O}(-1))))([\CP^1])=-1$. We conclude that
$$c_1({\bf R}^1(\pi_*)(\text{ev}^*({\mathcal
  O}(-1))))([\CP^1])=-\tfrac12.$$
Putting everything together gives
\begin{align*}
N_{0,2}&=\int_{[\wwbar{\mathcal
M}_{0,0}(\CP^1,2[\CP^1])]}c_{1}({\bf
R}^1(\pi_*)(\text{ev}_{\CP^1}^*{\mathcal O}(-1)))^2\\
&=c_{1}\bigl({\bf R}^1(\pi_*)\bigl(\text{ev}_{\CP^1}^*{\mathcal
O}(-1)\bigr)\bigr)^2\bigl(\tfrac12[\CP^2]\bigr)\\
&=\tfrac12\bigl(c_{1}\bigl({\bf R}^1(\pi_*)\bigl(\text{ev}_{\CP^1}^*{\mathcal
O}(-1)\bigr)\bigr)([\CP^1])\bigr)^2=\tfrac18\,.
\end{align*}
While this direct computation clarifies all of the ingredients in
the definition of the Gromov--Witten invariants, it is not very
practical for computing the general case. To compute the answer in
the general case, we return to localization.

\section{The full multiple cover formula via localization}\label{fmc}

In this section we apply localization to compute the Gromov--Witten
invariants of the resolved conifold $X_{S^3}$. We should really say
that we are using virtual localization. This is a generalization of
the localization formula that we have already explained in two
different directions. First one must apply localization in the stack
setting, and second one must apply it with virtual fundamental
classes. The correct generalization is given by Graber and Pandharipande
\cite{virloc}.

The standard torus action on $\CP^1$ extends to a torus action on
all of the spaces involved in the multiple cover formula for the
Gromov--Witten invariants of $\mathcal{O}(-1)\oplus\mathcal{O}(-1)$.
In fact this action extends to the bundle $\mathcal{O}(-1)$ in
several ways. Just as line bundles are classified by the first Chern
class, equivariant line bundles are classified by the equivariant
first Chern class. An action on a vector bundle compatible with an
action on the base is called a linearization. Using the
classification of equivariant bundles it is standard to label
linearizations by the associated equivariant first Chern class. Thus
we label the linearizations of $\mathcal{O}(-1)$ by
$n\alpha_0+m\alpha_1-h$.
\begin{exm}\label{fmce1}
Construct group actions on $\mathcal{O}(-1)$ corresponding to the
equivariant class $n\alpha_0+m\alpha_1-h$.
\end{exm}

The numbers that we need to compute are
$$N_{g,d}=\int_{[\wwbar{\mathcal
M}_{g,0}(\CP^1,d[\CP^1])]}c_{2d+2g-2}\bigl({\bf
R}^1(\pi_*)\bigl(\text{ev}_{\CP^1}^*({\mathcal O}(-1)\oplus{\mathcal
O}(-1))\bigr)\bigr).$$
Based on the computation of $N_3$ via localization in \fullref{N3} one might worry that this localization computation is going
to be very complicated. In fact, the computation is fairly
straightforward provided that one chooses the proper linearizations.
Faber and Pandharipande discovered that if one chooses $\alpha_0-h$
as the linearization on the first factor of ${\mathcal O}(-1)$ and
$\alpha_1-h$ on the second factor then only components of the fixed
point set corresponding to one graph contribute to $N_{g,d}$
\cite{FP}.

To see why only one graph can contribute let $[u,\Sigma,p]$ be a
stable map fixed by the group action. We start with the
normalization sequence
$$0\to\calO_\Sigma\to\nu_*\calO_{\wwhat\Sigma}\to\oplus_c\calO_c\to0.$$
We already used this sequence in computing the Euler class of the
normal bundle associated to deformations of the map in \fullref{mape}. As before we use $c$ to denote the nodes of $\Sigma$. We
have a related exact sequence for holomorphic sections of the
pull-back bundle $u^*\mathcal{O}(-1)$:
$$0\to u^*\calO(-1)\to\nu_*\nu^*u^*\calO(-1)\to\oplus_cL_{-1}|_{u(c)}\to 0.$$
The associated long exact sequence on cohomology reads,
$$\cdots\to H^0(\wwhat\Sigma,\nu^*u^*\calO(-1))
\to\oplus_nL_{-1}|_{u(n)}\to H^1(\Sigma, u^*\calO(-1)) \to\cdots.$$
There are bundles over the moduli space with fibers isomorphic to
the spaces in this exact sequence. We will use the fibers as names
for these bundles so for example, we will denote the bundle ${\bf
R}^1(\pi_*)\bigl(\text{ev}_{\CP^1}^*({\mathcal O}(-1))\bigr)$ by $H^1(\Sigma,
u^*\calO(-1))$.

The normalization is a union of smooth components
$\wwhat\Sigma=\coprod \Sigma_k$ and the cohomology
$H^0(\wwhat\Sigma,\nu^*u^*\calO(-1))$ is isomorphic to $\oplus
H^0(\Sigma_k,\nu^*u^*\calO(-1))$. The contribution
$H^0(\Sigma_k,\nu^*u^*\calO(-1))$ is trivial unless the component
$\Sigma_k$ is contracted to a point under the stable map in which
case it is isomorphic to $\C$. Since the original prestable curve is
connected we see that there is at least one node attached to each
contracted component and there are extra nodes if there is any
vertex in the graph associated to the stable curve with valence
grater than one.

It follows that the cokernel of the map
$$H^0(\wwhat\Sigma,\nu^*u^*\calO(-1)) \to\oplus_cL_{-1}|_{u(c)}$$
is only trivial if there are no vertices of valence greater than
one.

In the case the cokernel is nontrivial the long exact sequence of
bundles implies that $e\bigl({\bf
R}^1(\pi_*)\bigl(\smash{\text{ev}_{\CP^1}^*}({\mathcal O}(-1))\bigr)\bigr)$ contains a
factor of $e(L_{-1}|_{u(c)})$. In order for the stable map to be
equivariant $u(c)$ must either equal $q_0=[1:0]$ or $q_1$. Without
loss of generality it is equal to $q_0$. Now consider the factor
$e\bigl({\bf R}^1(\pi_*)\bigl(\smash{\text{ev}_{\CP^1}^*}
({\mathcal O}(-1))\bigr)\bigr)$
corresponding to the $\alpha_0-h$ linearization. The corresponding
linearization of $L_{-1}|_{u(c)}$ is trivial, implying that
$e(L_{-1}|_{u(c)})$ is trivial, so such a stable map cannot
contribute to the Gromov--Witten invariant.

It follows that the only components of the fixed point set that
contribute to the Gromov--Witten invariant correspond to graphs
consisting of one straight edge. Consider such a graph with genus
$g_k$ on the $q_k$ vertex and compute the Euler class in the
integrand and the Euler class of the normal bundle to the fixed
point set.

Start with the Euler class in the integrand. We have seen that the
cokernel of the map into the fibers over the nodes is trivial in
this case. On the other hand, the higher cohomology groups of any
skyscraper sheaf like $L_{-1}|_{u(c)}$ vanish. If follows
from the long exact cohomology sequence that
\begin{multline}\label{er1}
H^1(\Sigma,u^*\calO(-1))\cong
H^1(\wwhat\Sigma,\nu^*u^*\calO(-1))\\
\cong H^1(\Sigma_{g_0};L_{-1}|_{q_0})\oplus
H^1(\Sigma_{g_1};L_{-1}|_{q_1})\oplus H^1(\CP^1;u^*\calO(-1))\,.
\end{multline}
We now need to compute the equivariant Euler classes of these
bundles. For this we need to describe the group actions on the
relevant spaces. The degree $d$ line bundle over $\CP^1$ may be
described as equivalence classes $[z_0,z_1,\xi]$ where
$(az_0,az_1,a^d\xi)\sim (z_0,z_1,\xi)$ for $a\in \C-\{0\}$. One can
check that the action $$ \lambda\cdot
[z_0,z_1,\xi]=[\lambda_0z_0,\lambda_1z_1,\lambda_0^{n_0}\lambda_1^{n_1}\xi]
$$ is the linearization labeled by $dh-n_0\alpha_0-n_1\alpha_1$.
This is the answer to \fullref{fmce1}.
\begin{exm}
Show that the action on the fiber over $q_0:=[1:0]$ is
$(d-n_0)\alpha_0-n_1\alpha_1$. Use this to  show that the natural
linearization on $T\CP^1$ is $2h-\alpha_0-\alpha_1$. Also show that
the linearization on a tensor product of line bundles is the sum of
the corresponding linearizations.
\end{exm}

\begin{exm}
Recall that the holomorphic sections of the degree $d$ line bundle
over $\CP^1$ are just degree $d$ polynomials. Use this to show that
the linearization $dh-n_0\alpha_0-n_1\alpha_1$ on $\calO(d)$
 turns  $H^0(\CP^1;\calO(d))$ into the representation
 $$
\oplus_{k=0}^2[(k-n_0)\alpha_0+(d-k-n_1)\alpha_1]\,.
 $$
\end{exm}
Finally, recall that to get an action on the domain of the degree
$d$ map $[z_0:z_1]\mapsto [z_0^d:z_1^d]$ one must
pass to the $d$th power action on the codomain. This is kept track
of by dividing by $d$ at the end.

Kodaira--Serre duality implies that
$$
H^1(\CP^1;u^*\calO(-1))\cong
H^0(\CP^1;\calO(d)\otimes\calO(-2))^\vee\,,
$$
where the $\calO(d)$ inherits a linearization from $\calO(-1)$ and
$\calO(-2)$ has the natural linearization arising as the cotangent
bundle (dualizing sheaf) on $\CP^1$. Using the $\alpha_0-h$
linearization on $\calO(-1)$ the previous two exercises together
with the above remarks give the Euler class of the last summand of
equation \eqref{er1} as
\begin{multline}\label{er1.1}
e(H^1(\CP^1;u^*\calO(-1)))=
  \prod_{k=0}^{d-2}[(d-k-1)\alpha_0/d+(k-d+1)\alpha_1/d]\\[-2ex]
= (d-1)!d^{1-d}(\alpha_0-\alpha_1)^{d-1}.
\end{multline}
Similarly, with the $\alpha_1-h$ linearization we obtain
$$
e(H^1(\CP^1;u^*\calO(-1)))=
(-1)^{d-1}(d-1)!d^{1-d}(\alpha_0-\alpha_1)^{d-1}\,.
$$
Now Kodaira--Serre duality implies that
$$
H^1(\Sigma_{g_0};L_{-1}|_{q_0})\cong
H^0(\Sigma_{g_0};\omega_{\Sigma_{g_0}})^\vee\otimes
L_{-1}|_{q_0}\cong \E^\vee\otimes L_{-1}|_{q_0}\,.
$$
Recall that $\E$ is the Hodge bundle, which is by definition the
bundle over the moduli space with fiber over a point isomorphic to
the first cohomology of the curve representing the point with
coefficients in the dualizing sheaf $\omega_\Sigma$. By
\fullref{etensorl} we conclude that
\begin{equation}\label{er123}
e(H^1(\Sigma_{g_0};L_{-1}|_{q_0}))=\sum_{i=0}^{g_0}c_i(\E^\vee)t^{g_0-i}\,,
\end{equation}
where $t$ is the linearization (equivariant first Chern class) of
$L_{-1}|_{q_0}$. When the linearization on $L_{-1}$ is $\alpha_0-h$
we have $t=0$; when the linearization is $\alpha_1-h$ we have
$t=\alpha_1-\alpha_0$.

Now turn to the equivariant Euler class of the normal bundle to the
fixed point set. Consider a component of the fixed point set with
$0<g_0<g$. Formula \eqref{eaut} implies that
$e({\mathfrak{aut}}([\Sigma,p])^{\text{mov}})=1$, while formula
\eqref{edefs} gives
$$
e({\mathfrak{def}}([\Sigma,p])^{\text{mov}})=\left((\alpha_0-\alpha_1)/d-\psi_0\right)
\left((\alpha_1-\alpha_0)/d-\psi_1\right)\,.
$$
Similarly formula \eqref{eT} gives
$$
\prod_ce(T_{u(c)}\CP^2))=-(\alpha_0-\alpha_1)^2\,,
$$
formula \eqref{eH0} gives
$$
e(H^0(\wwhat\Sigma,\calO(u^*T\CP^2))^{\text{mov}})=
-(-1)^d(d!)^2d^{-2d} (\alpha_0 -\alpha_1)^{2d+2}\,,
$$
and formula \eqref{eH1} gives
\begin{multline*}
e(H^1(\wwhat\Sigma;\calO(\nu^*u^*T\CP^2))^{\text{mov}})=\\
\biggl(\sum_{i=0}^{g_0}c_i(\E^\vee)(\alpha_{0}-\alpha_1)^{g_0-i}\biggr)
\biggl(\sum_{i=0}^{g_1}c_i(\E^\vee)(\alpha_{1}-\alpha_0)^{g_1-i}\biggr).
\end{multline*}
Localization also works on stacks with virtual fundamental cycles
(see Graber and Pandharipande \cite{virloc}). The factor
$e({\mathfrak{ob}}([u,\Sigma,p])^{\text{mov}})$ does not need to be
computed because it is accounted for in the virtual fundamental
cycle. Automorphisms of generic elements of the fixed point
components do have to be taken into account. The virtual
localization formula reads (compare to \eqref{loc} and \eqref{cpnlocfor}):
$$
\int_{M^{\vir}}\phi=
\sum_F\frac{1}{|\A_F|}\int_{F^{\vir}}\frac{\iota^*_F\hat\phi}{e(N(F)^{\vir})}\,.
$$
Combining all of the above formulas together with \eqref{eN} and
\eqref{edefobu} shows that the contribution to the Gromov--Witten
invariant coming from the $g_0,g_1$ fixed point component is
\begin{small}
\begin{multline*}
d^{-1}(-1)^d(d!)^{-2}d^{2d} (\alpha_0{-}\alpha_1)^{-2d}(d{-}1)!
d^{-1}(\alpha_0{-}\alpha_1)^{d{-}1}(-1)^{d{-}1}(d{-}1)!d^{-1}
  (\alpha_0{-}\alpha_1)^{d{-}1}\\[-1ex]
\qua\biggl(\int_{[\wwbar{\calM}_{g_0,1}]^{\vir}}
  \sum_{i=0}^{g_0}c_i(\E^\vee)(\alpha_0{-}\alpha_1)^{g_0{-}1}c_{g_0}(\E^\vee)
\sum_{i=0}^{g_0}c_i(\E^\vee)(\alpha_1{-}\alpha_0)^{g_0{-}i}
  ((\alpha_0{-}\alpha_1)/d{-}\psi_0)^{-1}\biggr)\\[-1ex]
\qua\biggl(\int_{[\wwbar{\calM}_{g_1,1}]^{\vir}}
  \sum_{i=0}^{g_1}c_i(\E^\vee)(\alpha_0{-}\alpha_1)^{g_1{-}1}c_{g_1}(\E^\vee)
\sum_{i=0}^{g_1}c_i(\E^\vee)(\alpha_1{-}\alpha_0)^{g_1{-}i}
  ((\alpha_1{-}\alpha_0)/d{-}\psi_1)^{-1}\biggr).
\end{multline*}
\end{small}
Using the fact that $c_i(E^\vee)=(-1)^ic_i(E)$ in general together
with the relation $c(\E)c(\E^\vee)=1$ proved by Mumford for the
Hodge bundle \cite{mum5} this contribution can be simplified to
$$d^{2g-3}\int_{[\wwbar{\calM}_{g_0,1}]^{\vir}}c_{g_0}(\E)\psi^{2g_0-2}
\int_{[\wwbar{\calM}_{g_1,1}]^{\vir}}c_{g_1}(\E)\psi^{2g_1-2}.$$

\begin{exm}
Prove that the same formula is valid when either $g_0$ or $g_1$ is
zero provided
$\int_{[\wwbar{\calM}_{g_0,1}]^{\vir}}c_{g_0}(\E)\psi^{2g_0-2}$
is interpreted to be one when $g_0=0$.
\end{exm}

Expressions such as
$$b_{g_0}:=\int_{[\wwbar{\calM}_{g_0,1}]^{\vir}}c_{g_0}(\E)\psi^{2g_0-2}$$
are called Hodge integrals. The classes $c_k(\E)$ are called Hodge
classes and are often denoted by $\lambda_k$.

If a different pair of linearizations was chosen for the two $\calO(-1)$
factors, one would obtain a very different looking expression for
this Gromov--Witten invariant. This generates relations between the
various Hodge integrals that can be used together with similar relations
coming from an integral with a $\calO(1)$ factor to compute the Hodge
integrals. This is the approach taken by Faber and Pandharipande
\cite{FP}. The result obtained there is
$$\sum_{g=0}^\infty b_gs^{2g}=\left(\frac{s/2}{\sin(s/2)}\right).$$
Relations in the cohomology of moduli space allow one to express the
cubic Hodge integral from \fullref{xdeg0} in terms of the above
Hodge integrals. The answer  obtained in \cite{FP} for $g\ge 2$ is
$$\int_{[\wwbar{\mathcal M}_{g,0}]^{\vir}}c_{g-1}(\E)^3=
\frac{(1-2g)B_{2g}B_{2g-2}}{(2g-2)(2g)!}.$$
This can be combined with \fullref{xdeg0} to obtain the
following formula for the degree zero invariants of Calabi--Yau
$3$--folds:
\begin{equation}\label{Xdeg0}
N_{g,0}(X)=\frac{(-1)^{g-1}(2g-1)B_{2g}B_{2g-2}\,\chi(X)}{2(2g-2)(2g)!}\,.
\end{equation}

\section{The Gromov--Witten free energy}\label{FGW}

In this final section on the Gromov--Witten invariants, we cover a
way to package the Gromov--Witten invariants in case all of the
unmarked moduli spaces have zero virtual dimension. For a Calabi--Yau
$3$--fold, this assumption is true. This implies that the most
interesting invariants are of the form
$$
N_{g,\beta}(X):=\langle\ \rangle^X_{g,\beta}=\int_{[\wwbar
{\calM}_{g,n}(X,\beta)]^{\vir}}1\,.
$$
These are combined in the following definition.
\begin{defn}\label{fGW}
The Gromov--Witten free energy \index{Gromov--Witten free energy}
\index{$F$@$F^{GW}_X$ Gromov--Witten free energy} of a Calabi--Yau
3-fold $X$ is the following formal function depending on a complex
parameter $y$ and a cohomology class $t\in H^2(X;\C)$.
$$
F^{GW}_X(t,y):=\sum_{g=0}^\infty\sum_\beta N_{g,\beta}(X)e^{-\langle
t,\beta\rangle}y^{2g-2}\,.
$$
The restricted Gromov--Witten free energy \index{restricted
Gromov--Witten free energy} \index{$F$@$\wwhat F^{GW}_X$ restricted
Gromov--Witten free energy} is the sum taken over all non-zero
homology classes $\beta$.
\end{defn}
It might be interesting to formulate a free energy for more general
symplectic manifolds. The Gromov--Witten invariants of $X_{S^3}$ as
originally computed in the previous article are given by
$$
N_{g,d}:=N_{g,d[\CP^1]}=d^{2g-3}\sum_{g_0+g_1=g}b_{g_0}b_{g_1},
$$
where $b_g$ are the Hodge integrals computed by Faber and
Pandharipande presented in the form of the generating function
$$
\sum_{g=0}^\infty b_gs^{2g}=\left(\frac{s/2}{\sin(s/2)}\right)\,.
$$

\begin{exm}
Combine the last three displayed formulas to prove that the
restricted Gromov--Witten free energy of $X_{S^3}$ is
$$
\wwhat F^{GW}{X_{S^3}}=\sum_{d=1}^\infty \frac{1}{d}\left(2\sin
\frac{dy}{2}\right)^{-2}e^{-td}\,.
$$
\end{exm}
It is the full Gromov--Witten free energy that appears in the
gauge-string duality. The full Gromov--Witten free energy is simply
the sum of the restricted free energy with the degree zero invariant
computed earlier. After we give an overview of path integral
techniques we provide a heuristic argument for combining the
Gromov--Witten invariants into this generating function. This is
contained at the end of  \fullref{pCS}.

We can think of $t$ as a complex parameter if we identify
$H^2(X_{S^3};\C)$ with $\C$ via evaluation on $[\CP^1]$. The
Gopakumar--Vafa integrality conjecture states that the Gromov--Witten
invariants are uniquely specified by a set of integer invariants
called Gopakumar--Vafa invariants \index{Gopakumar--Vafa invariants}
(or BPS states) \cite{GV}. These BPS states \index{BPS
states}\index{$N$@$n^g_\beta$ BPS states} are denoted by $n^g_\beta$
and are supposed to be a count of embedded genus $g$,
$J$--holomorphic curves in the homology class $\beta$. The
Gopakumar--Vafa integrality conjecture \index{Gopakumar--Vafa
integrality conjecture} takes the exact form
$$
F^{GW}(X)=\sum_{g=0}^\infty\sum_\beta\sum_{d=1}^\infty n^g_\beta \frac{1}{d}\left(2\sin \frac{dy}{2}\right)^{2g-2}e^{-d\langle t,\beta\rangle}\,.
$$
Clearly $X_{S^3}$ satisfies the Gopakumar--Vafa conjecture with
$\smash{n^0_{[\CP^1]}}=1$ and the rest of the BPS states equal to zero.

\newpage
\part{Witten--Chern--Simons theory}
\setobjecttype{Part}
\label{CS}

There is a vast amount of literature on Witten's (quantum)
Chern--Simons theory, conformal field theory, quantum groups and so
on. We will try to outline the basic notions and definitions that
are needed to address Large $N$ Duality. The first subsection below
reviews background material about surgery and $3$--manifolds. The
easiest definition of Witten's Chern--Simons invariants is based on
skein theory (\fullref{app:c}). It is however, difficult to to
compute the resulting invariants or to show that they are well-defined.
For this reason we will work with the definition based on quantum
group. \fullref{pCS} describes the physical motivation for these
invariants. Path integral motivation leads to perturbative
invariants. Quantum field theory motivation leads to the
Reshetikhin--Turaev (exact) invariants. Definitions of the
Chern--Simons partition function and free energy based on the latter
are given in Sections \ref{modinv}, \ref{tilt} and \ref{csfree}
after several sections on background and motivation.

\section{Framed links and $3$--manifolds}

Chern--Simons theory provides topological invariants of
$3$--manifolds and framed links in $3$--manifolds.  There is a close
connection between framed links and $3$--manifolds, namely any
closed oriented $3$--manifold may be obtained by surgery on a framed
link in the $3$--sphere.

\begin{defn} A framed link \index{framed link} is an embedding of a finite disjoint union of
copies of $S^1\times D^2$ into a $3$--manifold. Two framed links are
considered equivalent if they are related by an ambient isotopy.
\end{defn}
Any ambient isotopy can be decomposed into elementary isotopies
called Reidemeister moves. Two framed links are isotopic if and
\index{Reidemeister moves} only if they are related by a finite
sequence of Reidemeister moves.
\begin{figure}[ht!]
\centering
\includegraphics[width=4truein]{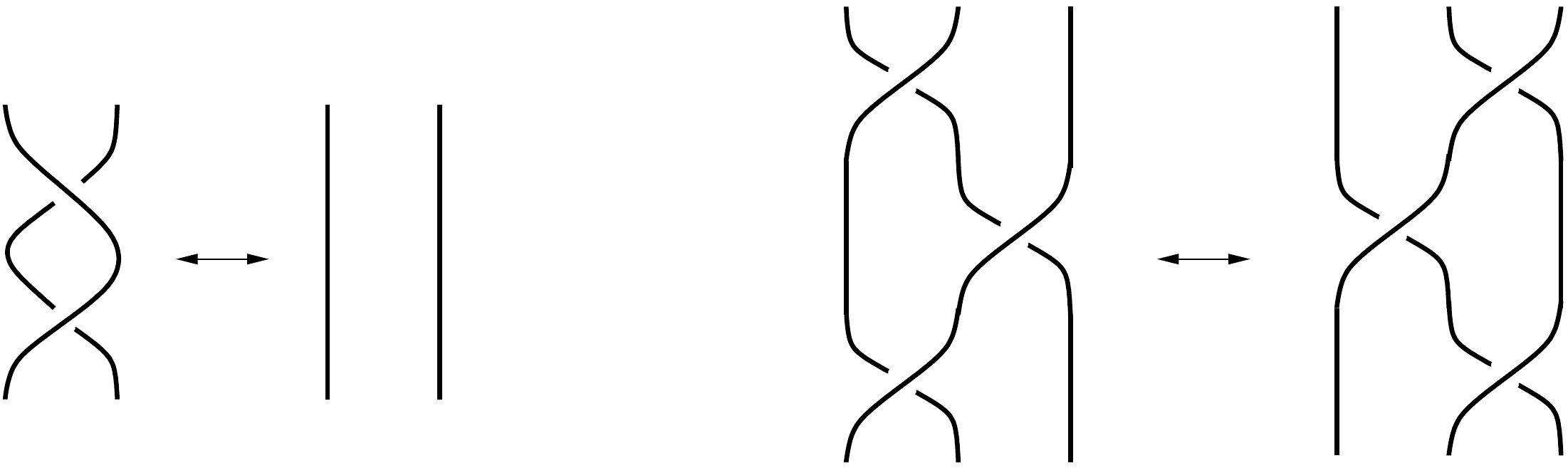}
\caption{Reidemeister moves II and III}\label{reidemeister}
\end{figure}
Framed links in $\R^3$ can be cut into elementary pieces called
tangles (see \fullref{lhopf}). It is often easier to analyze
elementary tangles.
\begin{defn}
A framed (or ribbon) tangle \index{tangle} is an embedding of a
finite disjoint union of copies of $S^1\times D^2$ and $[0,1]\times
D^2$ into $[0,1]\times \R^2$ taking $\{0,1\}\times D^2$ into
$\{0,1\}\times \R^2$. The embedding of $\{0,1\}\times D^2$ into
$\{0,1\}\times \R^2$ must be standard, depending on the number of
components landing on each boundary component so that ribbon tangles
may be stacked. Tangles are also considered up to isotopy.
\end{defn}

\begin{figure}[ht!]
\centering
\includegraphics[width=4truein]{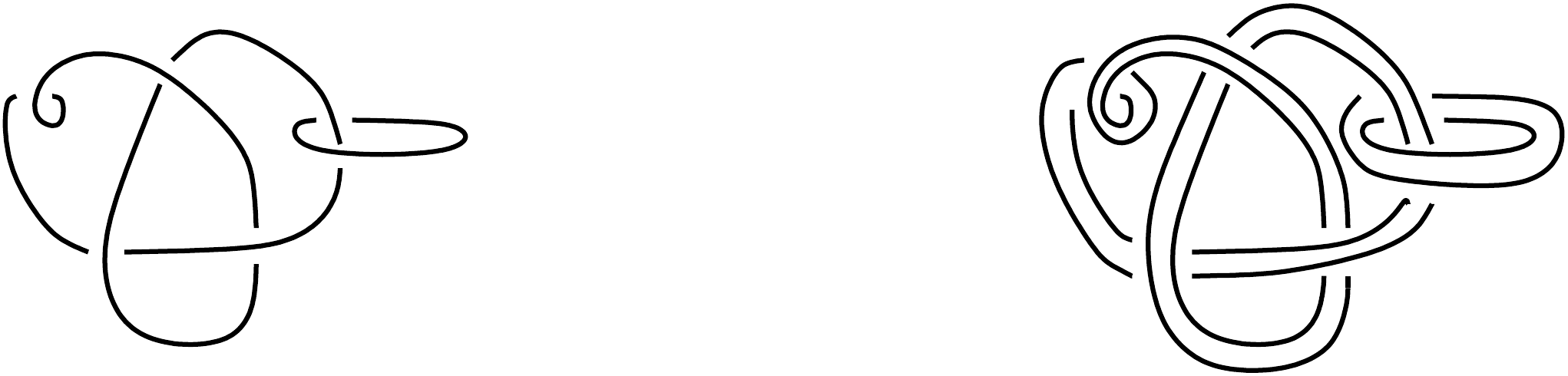}
\caption{Framed link projection and ribbon link}\label{framed-link}
\end{figure}

In order to move everything into closed $3$--manifolds it is standard
to work with the one point compactification of ${\mathbb R}^3$. This
is homeomorphic to the $3$--sphere. By general position we may assume
that any framed link in $S^3$ misses the point at infinity. All
framed links in $S^3$ can therefore be brought back to $\R^3$ and
represented by a projection of the cores $S^1\times \{0\}$ to a
plane keeping track of over-crossings and under-crossings.

To recover a framed link from the projection one first pushes the
over-crossings slightly above the plane to obtain an embedding
$\gamma\co \coprod S^1\hookrightarrow {\mathbb R}^3$. This is
then extended to an embedding $\what\gamma\co \coprod
S^1\times D^2\hookrightarrow {\mathbb R}^3$ by
$\what\gamma(t,x,y):= \gamma(t)+{\bf k}x+\dot\gamma(t)\times {\bf
k}y$. This convention is called the \index{blackboard framing}
blackboard framing. The same technique may be used to represent
ribbon tangles. \fullref{framed-link} displays the projection of
a framed link and the image of $S^1\times (\{0\}\times [0,1])$. This
second picture justifies the name `ribbon'.

Framed links in $S^3$ can be used to construct more complicated
$3$--manifolds by a process called surgery. \index{surgery}
\begin{defn}
Surgery on a framed link $\what\gamma$ refers to removing the
image of the $S^1\times B^2$'s and attaching the same number of
$S^1\times D^2$'s in such a way as to glue  $\{1\}\times S^1$ to
$\what\gamma(S^1\times\{(0,1)\})$.
\end{defn}
If you have never seen surgery before, the book by Rolfsen is a good
reference \cite{Rolfsen}. Other good references that are relevant to
this exposition are Prasolov--Sossinsky \cite{prasolov} and Kassel
\cite{kassel}.

\begin{figure}[ht!]
\centering
\includegraphics[width=3truein]{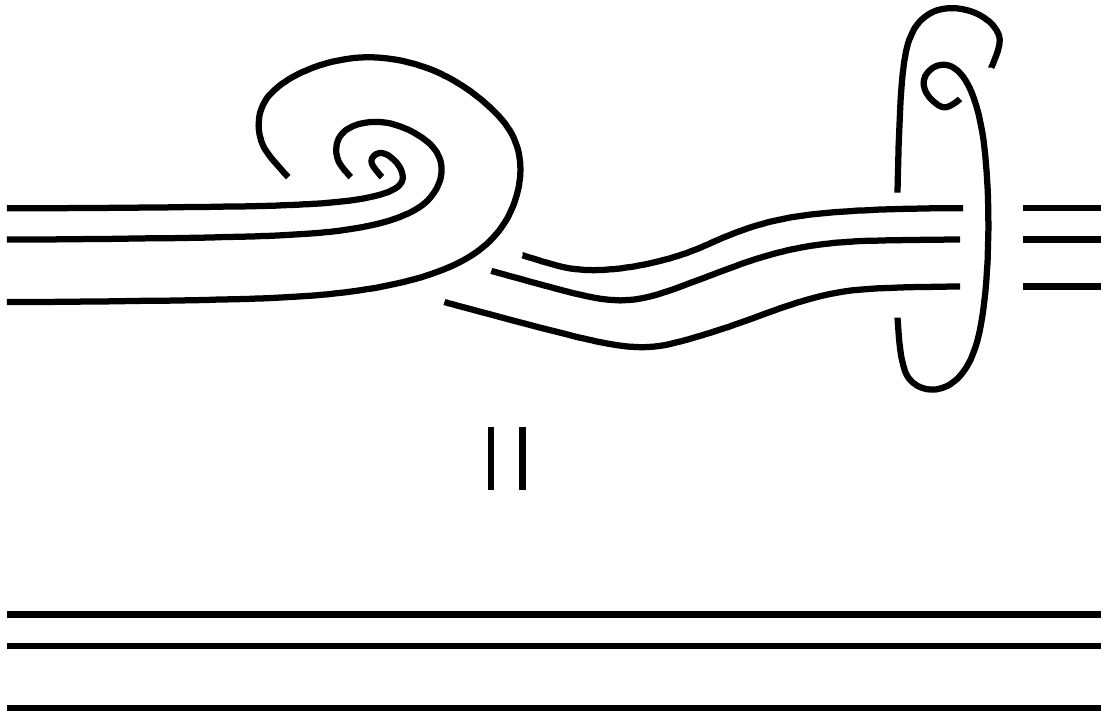}
\caption{Left Kirby move}\label{kirby-move}
\end{figure}

By a theorem of Lickorish and Wallace any closed oriented
$3$--manifold can be obtained by surgery on a framed link in $S^3$
\cite{Rolfsen}. It is not difficult to see that surgery on isotopic
framed links produces homeomorphic $3$--manifolds. What is less
obvious but still not difficult is that surgery on two framed links
related by an additional move called the Kirby move \index{Kirby
move} (\fullref{kirby-move}) still produces homeomorphic
$3$--manifolds. In fact, Kirby proved that two framed links
represent the same $3$--manifold if and only if they are related by
a sequence of Reidemeister and Kirby moves \cite{kirbysmovs,Rolfsen}.
Actually, what we are calling the Kirby move was introduced by Fenn and Rourke \cite{Rourkefenn}.
Kirby himself used an equivalent pair of moves: blow up/down and handle
slide \cite{kirbysmovs}. Blow up/down adds or removes an unlinked, $(\pm
1)$--framed
circle and handle slide tubes one component of a link to a parallel copy
of a second component. The names come from descriptions of
$3$--manifolds as boundaries of $4$--manifolds and the corresponding
moves on $4$--manifolds.
\begin{exm}
Show that surgery on the framed link with projection a simple circle
in the plane gives $S^2\times S^1$.
\end{exm}

A framed link in an arbitrary closed oriented $3$--manifold can thus
be represented by a pair of disjoint framed links in $S^3$. One of
them represents a surgery description of the $3$--manifold and the
other represents the framed link in it. It follows that to define an
invariant of $3$--manifolds or framed links in $3$--manifolds it suffices to define it on
a pair of framed links. Of course one then has to prove that this link invariant is
also invariant under the Kirby moves. There
are many invariants of links and framed links, some are described in
\fullref{app:c}. Most of them are not invariant under the Kirby moves
and therefore do not produce invariants of $3$--manifolds.
But Chern--Simons theory does generate $3$--manifold invariants along
these lines.

\section{Physical and heuristic descriptions}\label{pnhdesc}

The definition of the Chern--Simons invariants of $3$--manifolds was
motivated by quantum field theory (QFT) (see Deligne et al \cite{deligne}
and Zee \cite{Zee}). We can afford to be sketchy because this will
only serve as motivation for the mathematically rigorous definition of
the Reshetikhin--Turaev invariants that we use later. The mathematical
foundations of quantum field theory are not completely developed, but
the existing machinery and conjectural structure of QFT has produced and
motivated many remarkable mathematical theorems. The invariants that we
are discussing in this section are the expectation values of observables
of the QFT associated to the Chern--Simons action. These values do not
depend on any additional geometric structure and are therefore topological
invariants of the underlying $3$--manifold space-time. Such QFT's are
called topological quantum field theories TQFT's.

First recall the Lagrangian approach to quantum mechanics. In this
approach the partition function is equal to the integral of the
exponential of the action that is, $Z=\int e^{i\hbar^{-1}S}$. For
the quantum mechanics of a classical particle under the influence of
a conservative force with potential $V$ the action is $S=\int
\frac12 m \dot x(t)^2 - V(x(t))\,dt$. The classical equations of
motion are just the stationary curves of this functional. It is
exactly these stationary points that contribute to highest order in
the stationary phase approximation of $Z$. This generalizes in an
obvious way for extra degrees of freedom. The same framework
generalizes to continuum mechanics and field theory.

The fields in Chern--Simons theory are connections. A
\index{connection} connection can be viewed as a Lie algebra valued
$1$--form on the $3$--manifold. Think of the sum $A=A_i\theta^i$
where the matrices $A_i$ live in the Lie algebra and the $\theta^i$
are $1$--forms. If $R$ is a representation of a Lie algebra we
define a trace function on the algebra by
$\text{Tr}_R(A)=\text{Tr}(R(A))$. The Chern--Simons action of a
$U(N)$ connection is given by \index{Chern--Simons action}
\index{$CS$@CS Chern--Simons action}
$$
CS(A)=\frac{1}{4\pi}\int_M \text{Tr}_{\tableau{1}}\bigl(A\wedge
dA+\tfrac23 A\wedge A\wedge A\bigr)\,
$$
where $\tableau{1}$ is the defining representation of $U(N)$. Witten
suggested the simple idea that the average of a function of this
Chern--Simons action taken over all connections should be a
topological invariant. Indeed, nothing in the definition of $CS$
depends on any geometric data \cite{W1} except the connection which
is integrated out. This average is called the Chern--Simons
partition function. It is formally written as
$$
Z_k(M)=\int_{{\mathcal A}} e^{\frac{i}{2x}CS(A)}{\mathcal D}A\,,
$$
where ${\mathcal D}A$ is an as of yet undefined measure on the space
of connections and $x$ is the so-called string coupling constant.
\index{string coupling constant}\index{$X$@$x$ string coupling
constant} In the path integral expression one takes
$x=\frac{2\pi}{k}$ where $k$ is a positive integer called the level.
It turns out that after performing formal perturbative expansion in
$x$ one needs to `renormalize' it to $x=\frac{2\pi}{k+N}$, where $N$
is the rank of $U(N)$ to get the 'correct' answer. The explanation
for this shift  is not fully understood mathematically and
underscores subtleties of infinite-dimensional integration. This
shift comes from the interpretation of the signature of an operator
as an eta invariant, see Witten \cite{W1} and Atiyah \cite{At}.
Because of this some authors call the level $k$ and others will call
it $k+N$, so one must be careful when comparing different results in
the literature. We call the level $k$. \index{level}\index{$K$@$k$
level} \index{rank} \index{$N$ rank}

This formal representation of the partition function is an example
of a Feynman path integral (see Etingof \cite{ent}). This invariant
or various normalizations of it is more often denoted by $\tau(M)$
in the mathematical literature. The problem is that this `average'
is not well-defined because the space of connections is
infinite-dimensional and a translation-invariant measure does not
exist. However, there are ways to formally define invariants that
have most of the properties expected of this 'average'. It turns out
however that they do depend on an additional geometric structure on
a 3--manifold known as $2$--framing (trivialization of $T(TM)$ up to
homotopy) \cite{At,W1}. This phenomenon is called gravitational
anomaly \index{gravitational anomaly} by physicists and is sometimes
explained by the 'measure' ${\mathcal D}A$ not being purely
topological \cite{Zee}. Gravitation in physics is represented by a
background metric and a metric in its turn determines many additional
structures including a $2$--framing. Realistic quantum field theories
such as quantum chromodynamics do depend on metric or in physical
terms, are coupled to gravity. In (almost) topological quantum
Chern--Simons theory framing dependence can be seen as a lingering
ghost of this metric dependence.

This being said, $2$--framing is a very weak structure, so weak in fact
that every $3$--manifold admits a canonical one. Using it one can normalize the partition function so as to cancel out the framing dependence altogether. This is exactly the Reshetikhin--Turaev normalization of invariants that we adopt in Definitions \ref{taudef}, \ref{CSunnorm}. It differs from the physical normalization used by Witten \cite{W1} and this has implications for Large $N$ Duality. For instance, Ooguri and Vafa \cite{OV1} find a much better agreement between the gauge and string partition functions than we do. Unfortunately, there is no consistent definition of the 'physical normalization'. In examples it is usually derived ad hoc by comparing the exact answers to perturbative expansions, see eg Rozansky \cite{Roz}.

One can also `average' holonomies around colored framed links. A
colored framed link is a framed link with a group representation
associated to each component. The Chern--Simons invariants of
colored framed links are the expectation values of observables
constructed from the holonomy of \index{holonomy} connections. In
physical language these observables are called Wilson loop
operators. \index{Wilson loop operator} The holonomy is given by
$\text{Hol}_A(\gamma):={\rm P}\,\exp\, \oint_{\gamma} A$. More
explicitly this means that one solves the system of ODE's given by
$\frac{d}{dt}X(t)+A(\dot\gamma(t))=0$ with initial data $X(0)=I$.
Given this the holonomy is given by $\text{Hol}_A(\gamma):=X(1)$.
Here we are assuming that $\gamma(0)=\gamma(1)$.

A Wilson loop operator for one component is the trace of the
holonomy to a given connection along that component in a given
representation. That is $W^K_R(A):=\text{Tr}_R (\text{Hol}_A(K))$.
The link invariant associated to these Wilson loop operators is just
the vacuum expectation value (vev) \index{vacuum expectation value
(vev)} or correlation function: \index{correlation function}
$$
W_{R_1,\ldots, R_c}(L):=\frac{1}{Z_k(M)}\int_{{\mathcal A}}
e^{\frac{i}{2x}CS(A)}\prod W^{L_i}_{R_i}(A){\mathcal D}A.
$$
The invariants defined mathematically based on this motivation (up
to various different normalizations) are called colored Jones
polynomials for SU$(2)$ and the colored THOMFLYP polynomials for
SU$(N)$. They are sometimes also denoted by $J(L, R_1,\ldots,R_c)$
for $M=S^3$ or by $\tau(M,L)$ in general.

\begin{remark}\label{orhol}
This $W_{R_1,\ldots, R_c}(L)$ is an invariant of oriented framed
links, as changing the orientation inverts the holonomy.
\end{remark}

After a physical construction of invariants,
Witten went further and outlined ideas that led to one way of making
these invariants mathematically rigorous. He argued using skein relations that the
expectation values of Wilson loop operators are given by the Jones polynomial
of the corresponding links. Moreover, he gave an explicit prescription for
computing Chern--Simons partition functions of $3$--manifolds based on their
link surgery presentation, see Witten \cite{W1} and Axelrod--Della Pietra--Witten \cite{ADW}.
This was an amazing insight, but Witten's surgery prescription is a long way
from a mathematically rigorous definition of the invariants.
Reshetikhin and Turaev were first to devise a rigorous definition
based on quantum groups \cite{RT1,RT2}. This is the definition that
we will ultimately use.

There are two philosophically different ways to interpret the
expressions for these invariants: the perturbative approach and the
TQFT approach.  Each of these approaches can be formalized in
different ways. The perturbative approach is outlined in \fullref{pCS}
and the exact approach is outlined in \fullref{TQFT}. Witten's original
idea is also explained in the book of Atiyah \cite{At}.

\section{Perturbative Chern--Simons theory}\label{pCS}

Here we outline the perturbative approach to Chern--Simons theory
because it motivates many of the definitions and conjectures that appear
later. Numerous other authors have written expositions on perturbative
expansions (see Sawon \cite{sawon} in this volume, Bar-Natan \cite{bn0}
and Polyak \cite{polyak}). We include an overview here because it helps
motivate Large $N$ Duality. The perturbative approach is a generalization
of two ideas for finite-dimensional integrals: the stationary phase
approximation for oscillatory integrals, and a graphical calculus due
to Feynman for evaluating Gaussian integrals (see Etingof \cite{ent}).

Recall the stationary phase expansion, \index{stationary phase
expansion}
\begin{multline*}
\int_M  e^{i\lambda H}\,d\text{vol}_M=\\[-1ex]
\sum_{dH|_p=0} (2\pi \lambda^{-1})^{n/2}e^{\pi
i\text{sgn}(D^2H_p)}|\text{det}\,
D^2H_p|^{-1/2}e^{i\lambda H(p)}+O(\lambda^{-n/2-1}).
\end{multline*}
A theorem of Duistermaat and Heckman asserts that this is exact
(with no $O(\lambda^{-n/2-1})$ term) when $M$ is a symplectic
manifold and $H$ is an invariant Hamiltonian with only
non-degenerate critical points \cite{DH}. One nice proof of this is
based on the localization formula discussed in the section on
Gromov--Witten invariants. In fact, the localization formula was
originally discovered in an attempt to better understand why the
stationary phase approximation was exact (see Atiyah and Bott \cite{AB}).

Similarly, there is an infinite-dimensional symplectic structure on
the space of connections and the Chern--Simons action is invariant
under the action of the gauge group. So one might expect that the
stationary phase approximation is exact in this setting. The critical points
of the Chern--Simons action are flat connections
(see Baez and Muniain \cite{BzM}) and one can define a perturbative expansion about these flat
connections by analogy to the stationary phase approximation (see Bar-Natan\cite{bn0,bn2}).
For the unknot in $S^3$ the agreement between perturbative and exact invariants has been verified, see Bar-Natan--Garoufalidis--Rozansky--Thurston \cite{wheels}.
There are still interesting open questions related
to the appropriate interpretation of the full expansion on
nontrivial manifolds since such manifolds admit nontrivial flat
connections as critical points for the perturbative expansion.

To better understand the structure of the perturbative expansion we
consider a finite-dimensional Gaussian integral analog.
\begin{exm}\label{gtrick}
Recall that $f(a)=\int_{-\infty}^\infty e^{-ax^2}\,dx$ may be
evaluated by squaring it and then converting to polar coordinates.
By taking successive derivatives of $f(a)$ evaluate
$\int_{-\infty}^\infty x^{2n}e^{-ax^2}\,dx$.
\end{exm}
The expressions in the previous exercise become more complicated as
$n$ grows and are even more complicated for integrals over
higher-dimensional Euclidean spaces. Feynman added some slick
book-keeping machinery to produce an efficient method for computing
higher-dimensional Gaussian integrals. These integrals are analogous
to the path integrals that arise in Quantum Field Theory.

For finite-dimensional integrals the method is as follows. Let $Q$
be a symmetric bilinear form on $\R^n$ and let $V$ be a trilinear
form on $\R^n$. There is an obvious analogy between the following
Gaussian integral,
$$
Z  = \int_{\R^n} e^{-\hbar^{-1}(Q(x,x)/2+V(x,x,x)/6)}\, d^nx\,,
$$
and the path integral formally defining the Chern--Simons partition
function,
$$
Z_k(M)=\int_{{\mathcal A}} e^{\frac{i}{8\pi x}\left(\int_M
\text{Tr}_{\tableau{1}}(A\wedge
dA)+\int_M\text{Tr}_{\tableau{1}}(\frac23 A\wedge A\wedge A)\right)
}{\mathcal D}A\,.
$$
Substituting $x=\sqrt{\hbar}y$ and expanding the second exponential
gives
$$
\begin{aligned}
  & = \hbar^{\frac{n}{2}}\int_{\R^n} e^{-Q(y,y)/2}\cdot
   e^{-\sqrt{\hbar}V(y,y,y)/6}\, d^ny\\
  & = \hbar^{\frac{n}{2}} \sum^\infty_{m=0}\int_{\R^n} e^{-Q(y,y)/2}
   \frac{1}{6^{2m}(2m)!}   (-\sqrt{\hbar}V(y,y,y))^{2m}\, d^ny,
\end{aligned}
$$
the odd-order terms are missing from the above expression since their integrals
evaluate to $0$. A typical term in this sum may be evaluated by the
trick described in \fullref{gtrick} by diagonalizing the
quadratic form $Q$ or generalizing the trick to higher-dimensional
Gaussian integrals (see Sawon \cite{sawon} and Etingof \cite{ent}).
The result has the form
 \begin{equation}\label{e2.3.1}
 \hbar^{ \frac{n}{2}} \frac{(2\pi)^{\frac{n}{2}}}{(\det Q)^{\frac{1}{2}}}\quad
  \frac{\hbar^m}{6^{2m}(2m)!} \sum_\sigma W_\sigma,
\end{equation} where $\sigma$ represents a partition of the set
$\{1,\ldots,6m\}$ into two-element subsets encoded as a permutation
on $1,\ldots,6m$ and the $W_\sigma$ have the form
\begin{equation}\label{e2.3.2}
\sum_{j_1=1}^n{\cdots}\sum_{j_{6m}=1}^n\prod_{k=0}^{2m-1}V(e_{j_{3k+1}},
e_{j_{3k+2}},
e_{j_{3k+3}})\prod_{k=0}^{3m-1}Q^{-1}(e_{j_{\sigma(2k+1)}},
e_{j_{\sigma(2k+2)}}).
\end{equation}
with $ e_1,\ldots, e_n$ being the
standard basis in $\R^n$. For example, when $m=1$
\begin{equation}\label{e2.3.3}
\sum_{i_1=1}^n\!\cdots\!\sum_{i_{6}=1}^n V(e_{i_1},e_{i_2},e_{i_3})
V(e_{i_4},e_{i_5},e_{i_6}) Q^{-1}(e_{i_1},e_{i_2})
Q^{-1}(e_{i_3},e_{i_6}) Q^{-1}(e_{i_4},e_{i_5})
\end{equation}
is a typical term.
Terms with $Q^{-1}( e_{i_1}, e_{i_2})$ and $Q^{-1}( e_{i_2},
e_{i_1})$ are not distinguished, so in total we get
$\frac{1}{(3m)!}\binom{6m}{2,\ldots,2}=\frac{(6m)!}{(3m)! 2^{3m}}$
summands. This is cumbersome. Fortunately, there is a way due to
Feynman to represent such terms diagrammatically.

\begin{figure}[ht!]
\centering
\labellist
\pinlabel {$i_1$} [bl] at 61 64
\pinlabel {$i_2$} [tl] at 61 19
\pinlabel {$i_3$} [t] at 85 40
\pinlabel {$i_4$} [br] at 167 64
\pinlabel {$i_5$} [tr] at 167 19
\pinlabel {$i_6$} [t] at 140 40
\endlabellist
\includegraphics[width=3truein]{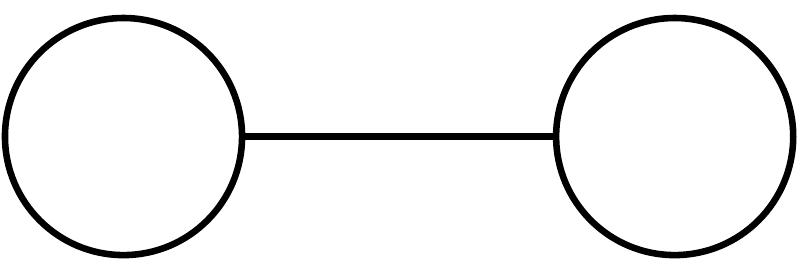}
\caption{Glasses}\label{sfig2.1}
\end{figure}

We construct a trivalent graph $\Gamma$ for each term $W_\sigma$
with a vertex for each $V$ in \eqref{e2.3.2} and an edge for each
$Q^{-1}$. The edge labeled $Q^{-1}(e_a,e_b)$ will connect to the
vertex or vertices containing $a$ or $b$. The graph corresponding to
example \eqref{e2.3.3} is then represented by the graph (Feynman
diagram) in \fullref{sfig2.1}. It is \index{Feynman diagram} easy
to see that different summands in \eqref{e2.3.1}  give the same
contributions as long as they have isomorphic graphs. So instead of
summing over all partitions $\sigma$ as in \eqref{e2.3.1} we may sum
over graphs as long as we factor in the number of different
partitions associated with each graph properly. Counting the number of
partitions associated to a given trivalent graph is an elementary
combinatorial problem with the answer (see Etingof \cite{ent})
\begin{equation}\label{e2.3.4}
\#\text{ partitions associated to }
  \Gamma=\frac{6^{2m}(2m)!}{|\text{Aut}(\Gamma)|},
\end{equation}
where $\text{Aut}(\Gamma)$ is the automorphism group of the
graph. Here we view a graph as a one-dimensional complex and an
automorphism must restrict to a linear map on each edge.

\begin{figure}[ht!]
\centering
\includegraphics[width=1truein]{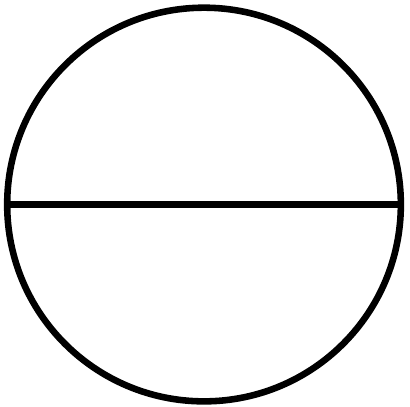}
\caption{Theta}\label{sfig2.2}
\end{figure}

\begin{example}
There are $15$ partitions of $1,\ldots,6$ into two-element sets, so
the sum in \eqref{e2.3.1} for $m=1$ would have $15$ terms. However,
there are exactly two trivalent graphs with $2m=2$ vertices (see
Figures \ref{sfig2.1} and \ref{sfig2.2}) so the corresponding sum over
graphs would only have two terms.  The automorphism groups of the
graphs in Figures \ref{sfig2.1} and \ref{sfig2.2} have orders $8$ and
$12$ respectively.
\end{example}

Notice that
\begin{equation*} \chi(\Gamma):=\text{\ \# vertices$(\Gamma)-$\#
edges$(\Gamma)$ }=2m-3m=-m.\end{equation*}
Thus the partition function
\index{$Z$ partition function} $Z$ can be rewritten as
\begin{equation}\label{e2.3.5} Z=\frac{(2\pi\hbar)^{{n}/{2}}}{(\det
Q)^{{1}/{2}}} \sum_{\Gamma\in\Xi}
\frac{\hbar^{-\chi(\Gamma)}}{|\text{Aut}(\Gamma)|} W(\Gamma),
\end{equation}
where the sum is taken over the set of all trivalent graphs
(including disconnected ones) $\Xi$ and $W(\Gamma)$ is the
contribution or \emph{Feynman amplitude} $W_\sigma$ \index{Feynman amplitude} of any
partition with graph $\Gamma$.

The sum
$$
F=\sum_{\Gamma\in\Xi^\prime}
\frac{\hbar^{-\chi(\Gamma)}}{|\text{Aut}(\Gamma)|} W(\Gamma),
$$
taken over the set of all connected trivalent graphs $\Xi^\prime$ is
called the {\it free energy}. \index{free energy}

\begin{exm}\label{freex}
Show that if $\Gamma_1\Gamma_2$ represents the disjoint union of $\Gamma_1$ and $\Gamma_2$ then one has
\begin{align*}
W({\Gamma_1\Gamma_2}) &= W({\Gamma_1}) W({\Gamma_2}),\\
\chi(\Gamma_1\Gamma_2) &= \chi(\Gamma_1)+\chi(\Gamma_2)\\
|\Aut(\Gamma^{n_1}_1\ldots\Gamma^{n_e}_e)|
    &=\prod^e_{i=1} |\Aut(\Gamma_i)|^{n_i} (n_i)!.
\end{align*}
Use exponentiation and series expansion to conclude that
$F=\ln(Z/Z_0)$ where $Z_0=(2\pi\hbar)^{{n}/{2}}/(\det Q)^{{1}/{2}}$.
\end{exm}
We conclude that it will be helpful to consider the natural
logarithm of the exact Chern--Simons invariants.

Of course the notions of partition function and free energy were
introduced in the field of statistical thermodynamics before they
ever appeared in Quantum Field Theory (see Schr\"odinger
\cite{schrod}). In statistical mechanics, the probability that a
state at energy $E$ in a given system will be occupied is given by
$e^{-E/(kT)}$. The partition function in this context is defined to
be the following integral over phase space:
$Z=\int_{T^*Q}e^{-\frac{H}{kT}}\,d\text{vol}$, where $H$ is the
Hamiltonian. In this context the free energy is defined by $F=-kT\ln
Z$.

\begin{figure}[ht!]
\centering
\labellist\Large
\pinlabel {$0$} at 540 100
\endlabellist
\includegraphics[width=3truein]{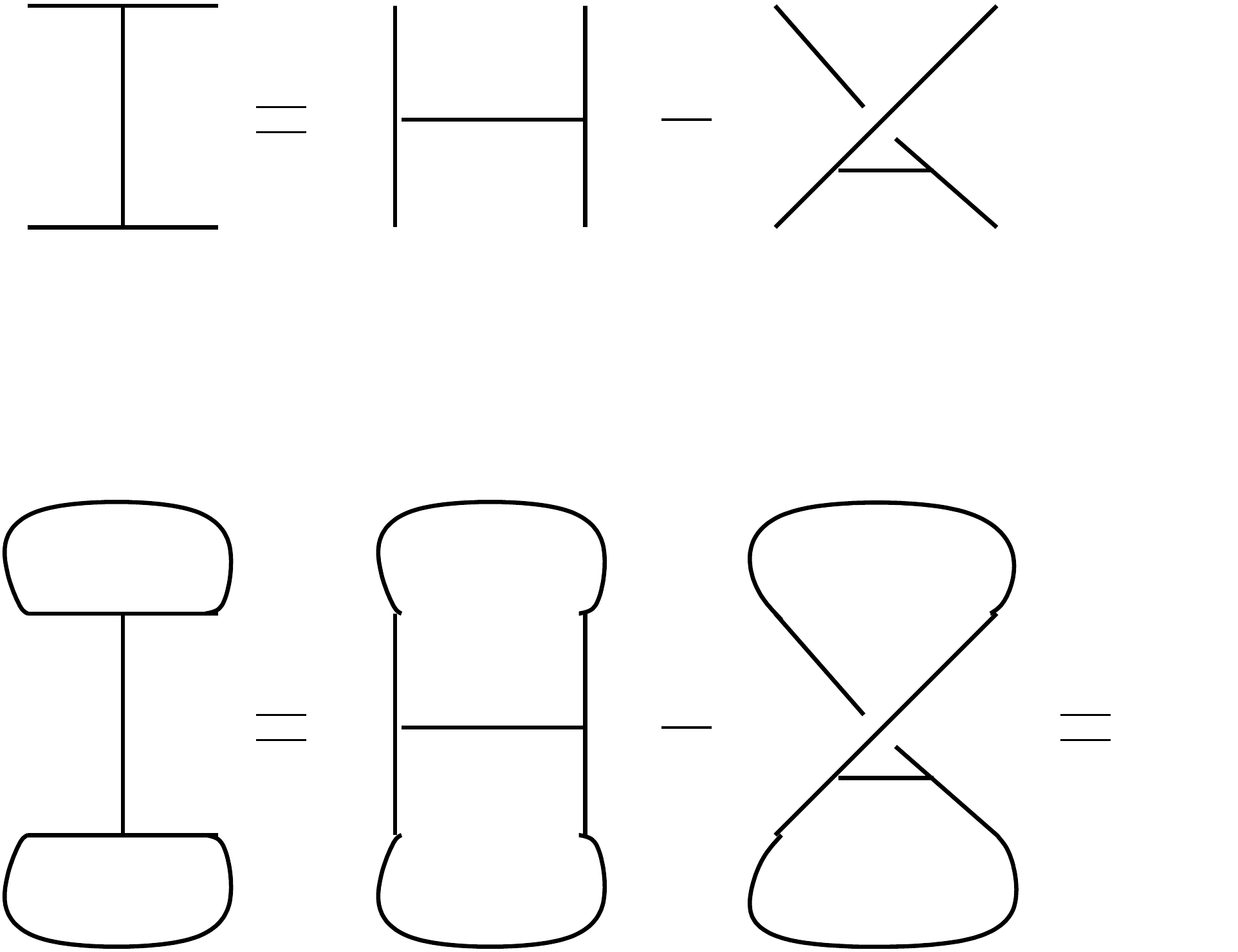} \caption{The IHX
relation}\label{IHX}
\end{figure}

The finite-dimensional analogy may be taken further and used to
motivate definitions of similar invariants of $3$--manifolds and
framed links in $3$--manifolds from Chern--Simons theory. The common
feature of all of these invariants is that they are expressed as
sums over graphs analogous to equation \eqref{e2.3.5}. Thus it makes
sense to introduce the free algebra generated by all trivalent
graphs with multiplication being disjoint union as in
\fullref{freex}. Two different combinations of graphs can have the
same contribution so we introduce relations in the algebra to
identify such combinations.
\begin{example}
One such relation is the IHX relation displayed in \fullref{IHX}.
The meaning of this relation is that the contribution of any graph
containing the piece on the left can be replaced by the contribution
of the difference of the two graphs containing the pieces on the
right. As \fullref{IHX} shows this relation implies that the
contribution of the `glasses' graph is trivial.
\end{example}

It is natural to consider the quotient algebra by these relations to
reduce the expression for the invariants even further. In fact this
sum can be taken over a basis for the quotient algebra. The
resulting expression for the free energy takes the form
\begin{equation}\label{freeform} F(M)=\sum_{\Gamma\in\text{\em B}}
{\hbar^{-\chi(\Gamma)}} {\mathcal W}^G(\Gamma)\,\text{FA}(M,\Gamma),
\end{equation}
where ${\mathcal W}^G$ is a homomorphism from the algebra of graphs
to the complex numbers called a weight system, {\em B} is a basis
for the algebra of graphs and $\text{FA}(M,\Gamma)$ is the Feynman
amplitude associated to the graph and the $3$--manifold.

\begin{figure}[ht!]
\centering
\labellist
\pinlabel {$0$} at 38 20
\pinlabel {$1$} at 455 20
\endlabellist
\includegraphics[width=4truein]{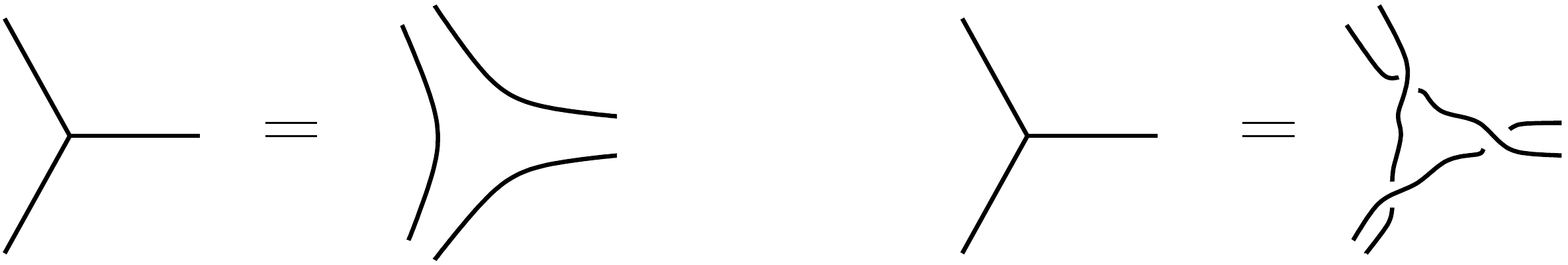} \caption{Fat
Graphs}\label{sfig2.8}
\end{figure}

We are particularly interested in the weight system for $U(N)$ (see
Bar-Natan \cite{bn2}). This weight system applied to a graph can be
computed as a sum over labelings. Given a graph $\Gamma$ label each
of the vertices with $0$ or $1$, then fatten the graph according to the
rules in \fullref{sfig2.8}.
The graph then turns into what is called a `fat graph' which
topologically represents a Riemann surface with boundary.  Let $g$
be the genus of the surface and $h$ the number of boundary
components. Also, let $\l$ denote a labeling of the graph,
$\Gamma_{(\l)}$ the labeled graph,  ${\Lambda}_{(\l)}$ its fattened
version, and $|\l|$ the sum of all labels in $\l$; then
\begin{equation}\label{uwt} \mathcal{W}^{U(N)}(\Gamma)=\sum_\l  (-1)^{|\l|}
N^{h({\Lambda}_{(\l)} )}. \end{equation}

\begin{figure}[ht!]
\centering
\labellist
\pinlabel {$0$} [r] at 30 118
\pinlabel {$0$} [l] at 187 118
\pinlabel {$0$} [r] at 690 118
\pinlabel {$1$} [l] at 845 118
\endlabellist
\includegraphics[width=4.5truein]{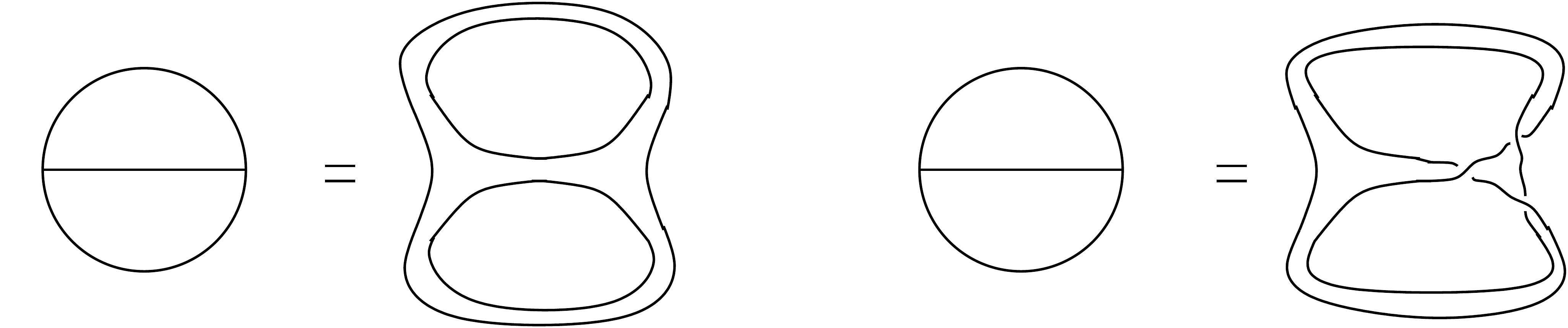} \caption{Fat
theta}\label{sfig2.9}
\end{figure}
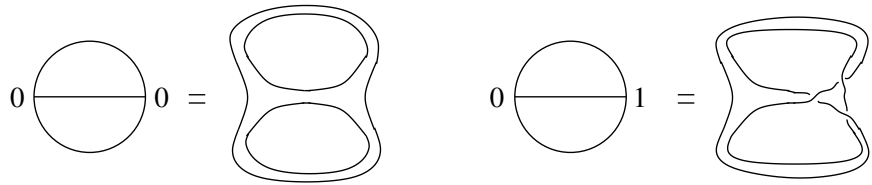

\begin{example}
To compute the U$(N)$ weight of the `theta' graph consider the fat
graphs corresponding to the labels $(0,0)$ and $(0,1)$ displayed in
\fullref{sfig2.9}. We see that the $(0,0)$ labeled graph has
genus zero and three boundary components, and the $(0,1)$ labeled
graph has genus one and one boundary component. The $(1,1)$ and
$(1,0)$ labeled graphs can be constructed similarly. Thus, for the
`theta' graph $\mathcal{W}^{U(N)}(\Gamma)=N^3-N-N+N^3=2N^3-2N$.
\end{example}

Since the Euler characteristic is a homotopy invariant,
$2g-2+h=-\chi({\Lambda}_{\l})=-\chi(\Gamma)$ for any labeling.
Combining equations \eqref{freeform} and \eqref{uwt} one sees that
the $U(N)$ Chern--Simons free energy has the form
$$
F^{\text{CS}}_M=\sum_{g=0}^\infty\sum_{h=1}^\infty
x^{2g-2+h}N^hF_{g,h}(M)\,.
$$
Here $x$ is playing the role of $\hbar$ in our finite-dimensional
Gaussian integral. The free energy takes the form of the free energy
of an open string theory. One expects that there is such a theory
of `instantons at infinity' (degenerate curves) but there is no mathematical
definition of this theory.

The fat graphs lead to closed surfaces as explained/conjectured by
't Hooft  and these in turn lead to the $J$--holomorphic curves that
are counted on the Gromov--Witten side. See Ooguri--Vafa \cite{OV1}
for this idea applied to Chern--Simons theory. 't Hooft suggested to
`sum over all holes' in this sum to obtain a `closed string'
expansion \cite{tH}.  This means introducing a new parameter $t=xN$
and combining all summands with like powers of $h$. Denoting
$$
F_g(M)=\sum_{h=1}^\infty t^hF_{g,h}(M)\,,
$$
we obtain
$$
F^{\text{CS}}_M=\sum_{g=0}^\infty x^{2g-2} F_g(M)\,.
$$
\begin{quote} This expression for the Chern--Simons free energy has the
structure of the free energy of a closed string theory and is one
reason to believe that there may be some relationship between string
theory and Chern--Simons invariants.
\end{quote}

We explain the structure of the free energy of a closed string
theory in more detail at the end of this section. On the string
theory side, the `instantons at infinity' live in the cotangent
bundle to the $3$--manifold and are open strings. The cotangent
bundle then undergoes a geometric transition where the boundaries of
the open strings are collapsed to points giving closed strings on
the  manifold on the other side of the transition. The manifold
across the transition from the cotangent bundle to $S^3$ is the
resolved conifold. Thus, one expects that the Chern--Simons free
energy is the same as the free energy of a closed string theory on
the resolved conifold. There is a mathematically defined closed
string theory on the resolved conifold, namely the Gromov--Witten
theory. Identifying it as the correct dual theory completes the
physical derivation of the duality and was the major contribution of
Gopakumar and Vafa \cite{GV}, see also Ooguri--Vafa \cite{OV1}.

\begin{quote}
To summarize, the first step is to describe the Chern--Simons
invariants via fat graphs that should be considered as open strings.
The second step is to sum over the `holes' via a geometric
transition to obtain a closed string theory.
\end{quote}

If the expression for $F^{\text{CS}}_M$ is rewritten using $N$ as a
parameter in place of $x=tN^{-1}$ one obtains a $(1/N)$--expansion
and the $g=0$ terms will dominate for large $N$.  When the
gauge-string duality holds for the leading terms (genus zero
contributions) it is said to hold in the large $N$ limit. This was
the case originally studied by 't Hooft. Witten realized that the
duality would be exact without considering the large $N$ limit for
topological (metric independent) field theories, see Witten
\cite{Wcss}.

The physical ideas here have been encoded into various
mathematically defined perturbative Chern--Simons invariants
\cite{bn2}. Kontsevich defined a universal Vassiliev invariant for
links taking the form of a rational linear combination of trivalent
graphs in an algebra generated by trivalent graphs with a few simple
relations, \cite{KontU1,KontU2}. See Bar-Natan \cite{bn1} for a good
time reading about these invariants and see \cite{bn2} by the same
author for a more typical overview.  Schematically this invariant
takes the form
\begin{multline*}
Z(L)= \\
\sum_{m=0}^\infty\sum_{(z_1,z_1^\prime)
  \ldots,(z_m,z_m^\prime)}E({(z_1,z_1^\prime),
\ldots,(z_m,z_m^\prime)};L)
  \Gamma({(z_1,z_1^\prime),\ldots,(z_m,z_m^\prime)};L),
\end{multline*}
where ${(z_1,z_1^\prime),\ldots,(z_m,z_m^\prime)}$ are $m$ pairs of
points on the link, and
$$\Gamma({(z_1,z_1^\prime),\ldots,(z_m,z_m^\prime)};L)$$
is the chord diagram
representing the locations of these points on the circles in $L$.
Furthermore, $E({(z_1,z_1^\prime),\ldots,(z_m,z_m^\prime)};L)$ is
some expression computed as an integral or via intersection theory
from the link $L$ and pairs of the points. To get numerical
invariants one just applies an algebra homomorphism from the algebra
of chord diagrams to the complex numbers. Such homomorphisms are
called weight systems \cite{bn2}.

There is a similar set of $3$--manifold invariants. A universal
$3$--manifold invariant of this type (now called the LMO invariant)
was introduced by Le, Murakami and Ohtsuki \cite{LMO}, and expressed
like the Feynman expansion of a Gaussian integral in \cite{bar} by Bar-Natan, Garoufalidis, Rozansky and
Thurston. As in the case of links, the universal $3$--manifold
invariant is a weighted sum of graphs.  This invariant is an element
in an algebra obtained as a quotient of the algebra freely generated
by all trivalent graphs.

\index{free energy}
\index{Gromov--Witten free energy}
It is interesting to see how string theory considerations suggest
that the Gromov--Witten free energy will take the form
$F^{\text{GW}}=\sum F_g y^{2g-2}$. The action for a simple model of
a vibrating string is
$$
S=\frac12\iint u_{tt}-u_{ss}\,ds\,dt\,.
$$
Notice the close similarity between this action and the Dirichlet
functional
$$
D=\frac12\iint u_{tt}+u_{ss}\,ds\,dt\,.
$$
As we saw in \fullref{loui}, the minima of the Dirichlet
functional are exactly the $J$--holomorphic maps. One common feature
of all string theories is the existence of an internal degree of
freedom $s$ in addition to the time $t$ appearing in the action.
This means that one must consider collections of surfaces in string
theory where one considered paths in ordinary field theory.

To discretize a path one just subdivides the interval. To develop a
discrete model of surfaces it is natural to triangulate the
surfaces. Notice that the dual graph of a triangulation is a
$3$--valent graph similar to the Feynman diagrams encountered in path
integrals. In fact neighborhoods of these dual graphs are the `fat
graphs' that we just discussed above. One can turn the process that
we used in this subsection backwards and write out a partition
function that would have these `fat graphs' in its perturbative
expansion. The most obvious candidate is the matrix integral
$$
Z=\int e^{-N \text{Tr}(\frac12 M^2 + w M^3)} \,dM\,,
$$
where the integral is taken over the space of all $N\times N$
Hermitian matrices.

When one performs a perturbative expansion on this partition
function the important things to notice about each term in the
expansion are
\begin{enumerate}
\item Each vertex contributes a factor of $wN$.
\item Each edge corresponds to a propagator and contributes a factor
of $N^{-1}$.
\item Each face (of the dual complex to the triangulation)
contributes a factor of $N$ (for the sum of the indices).
\end{enumerate}
It follows that each term in the expansion has order
$w^VN^{V-E+F}=w^VN^{2-2g}$. Thus the free energy is also a sum of
terms of order $w^VN^{V-E+F}=w^VN^{2-2g}$ because it is just the sum
of the contributions from the connected graphs. It follows that the
free energy can be written as
$$
F=\sum_g F_g N^{2-2g}\,.
$$
Let us summarize the main conclusions of this section. It is natural
to package Chern--Simons perturbative invariants into a formal
partition function, then take the natural logarithm, introduce a new
variable $t=xN$ and expand into a power series in $x$. Furthermore,
this function called the free energy will take the form $F=\sum
F_gx^{2g-2}$ identical to the form of the free energy in any closed
string theory.

In the next section we review physical motivations behind
Topological Quantum Field Theory (TQFT) and describe in detail a
simplified version of it that can be nicely packaged into the
language of ribbon and modular categories invented by N Reshetikhin
and V Turaev. For a particular choice of categories coming from the
theory of quantum groups these constructions produce the celebrated
Reshetikhin--Turaev or quantum invariants of framed links and
$3$--manifolds that can be seen as a mathematical formalization of
Witten's Chern--Simons path integrals.

\section{Modular categories and topological invariants}
\label{TQFT}

Chern--Simons theory is a special type of quantum field theory called
a topological quantum field theory (TQFT). Just as classical
mechanics may be described in the Lagrangian or Hamiltonian
frameworks any quantum field theory may be described in these two
frameworks. The Lagrangian framework leads to path integrals and the
Hamiltonian approach leads to canonical quantization (see Simms and
Woodhouse \cite{wod}). We have already discussed the Lagrangian approach and the resulting
perturbative Chern--Simons invariants. We will use the Hamiltonian
framework for formal definitions and computations that we do from
here on.

\subsection{The Hamiltonian approach to TQFT}
In the Hamiltonian approach, one begins with a symplectic manifold
called the phase space. This is the collection of all positions and
momenta. The mechanical system is specified by specifying an energy
or Hamiltonian function denoted by $H$ on this space. The
evolution of the mechanical system is given by Hamilton's equations
(see Arnol'd \cite{arnold}). A quantization of a classical mechanical
system is a map from the space of smooth functions on the phase space to
the Hermitian operators on a Hilbert space.  The relation between the
Hamiltonian and Lagrangian is given by
$$
\langle \phi_0|e^{itH}|\phi_1\rangle \,=\, \int_{\phi(0)=\phi_0,
\phi(1)=\phi_1} e^{iL(\phi)}{\mathcal D}\phi,
$$
Here $|\phi_k\rangle$ are elements of the Hilbert space,
$\langle\phi_k|$ are the functionals obtained from the elements via
the pairing on the Hilbert space (bras and kets), $H$ is the
Hamiltonian, $L$ is the Lagrangian and the right hand side is a path
integral (that is, the formal integral over the space of paths).

One can reinterpret this last expression by calling $|\phi_1\rangle$
a state and considering the Hilbert space as the space of states.
The operator $e^{itH}$ just represents evolution through $t$ units
of time. Then one sees that the last displayed expression expresses
the evolution of the state $|\phi_1\rangle$ through time as a path
integral.

\begin{figure}[ht!]
\centering
\labellist\small
\pinlabel {Glue} [b] at 370 10
\endlabellist
\includegraphics[width=4truein]{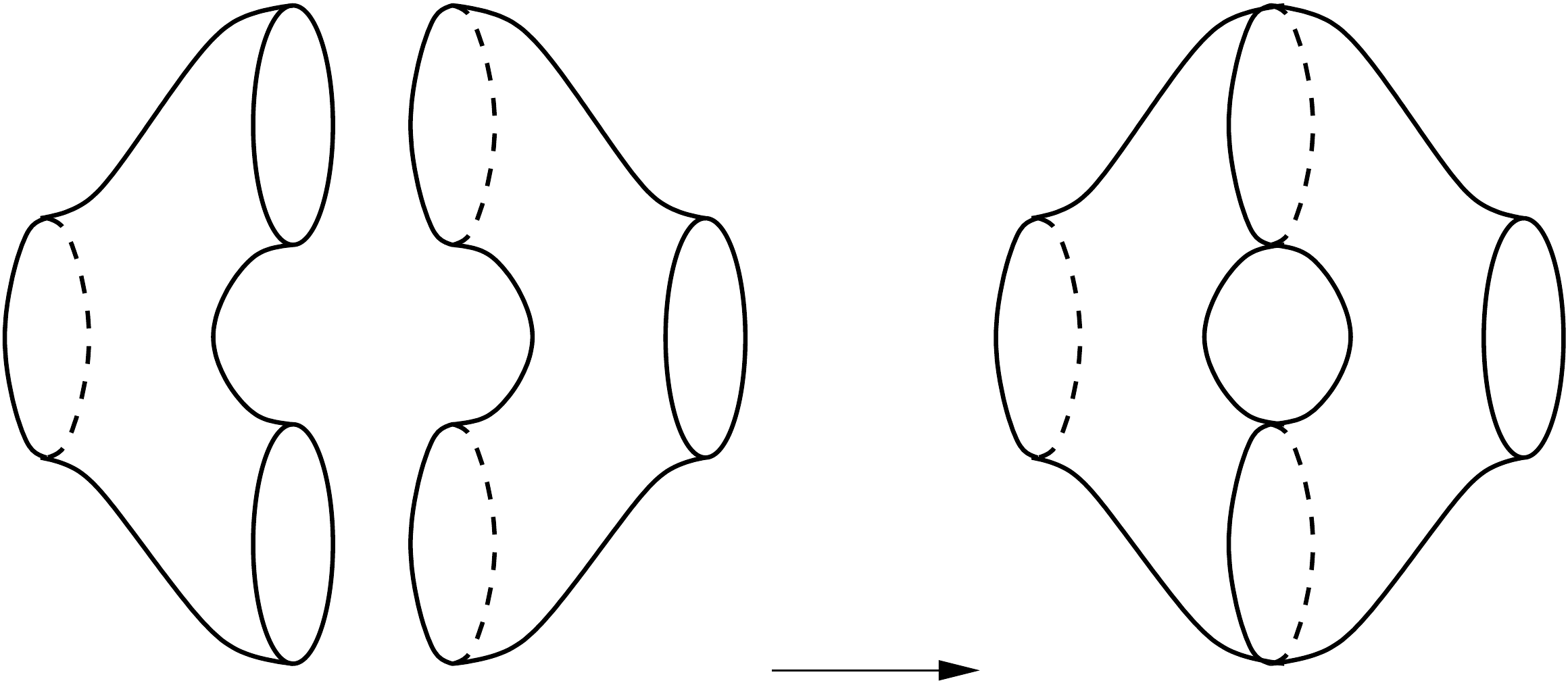} \caption{Gluing
for topological quantum field theories}\label{tqftfig}
\end{figure}

Geometrically this suggests considering manifolds with two boundary
components and associating a Hilbert space to each boundary
component. The Hilbert space associated to one boundary component
corresponds to the initial states and the Hilbert space associated
to the other component corresponds to the states after evolving
through time. The path integral should give an operator from the
Hilbert space associated to one boundary component to the Hilbert
space associated to the other boundary component. Gluing two
manifolds along a common boundary as in \fullref{tqftfig}
corresponds to taking the composition of the corresponding
operators. The Hilbert space associated to an empty boundary should
just be $\C$; thus operators corresponding to closed manifolds can
be interpreted as complex numbers.

Formalizing these ideas leads one to the notion of a topological
\index{topological quantum field theory} quantum field theory.
Slices of $3$--manifolds have $2$--dimensional boundary components.
Here slice means cobordism (that is, a $3$--manifold with boundary
components $\Sigma_\pm$). One views a product cobordism $[0,1]\times
\Sigma$ as a space-time with two spatial dimensions. Hence such
theories are called $(2{+}1)$--dimensional TQFT's. A
$(2{+}1)$--dimensional TQFT associates a Hilbert space ${\mathcal
H}_\Sigma$ to any closed surface and a bounded linear map,
${\mathcal H}_{\Sigma_-}\to{\mathcal H}_{\Sigma_+}$ to every
cobordism. This formalism was suggested by Segal \cite{Segal0,sea2,sea}
and axiomatized by Atiyah \cite{At}.

\begin{aside}{\small
Recall that the TQFT approach has origins in the Hamiltonian
framework. This is a brief outline of the Hamiltonian description of
Chern--Simons theory. In Chern--Simons theory the Lagrangian is just
the Chern--Simons functional and the corresponding Hamiltonian is
zero. The phase space used for the Hamiltonian description of
Chern--Simons theory is the space of flat connections modulo gauge
equivalence over a surface  $\Sigma$ of genus $g$. Equivalently this
is the space of all connections modded out by the complexified gauge
group (the geometric invariant theory picture) or the symplectic
quotient of the space of connections by the gauge group (the
symplectic reduction picture). Denote this moduli space of flat
connections by ${\mathcal M}_\Sigma$.

In order to construct the associated Hilbert space for Chern--Simons
theory, we need to introduce a line bundle over the moduli space
${\mathcal M}_\Sigma$. Quillen's determinant line bundle is the
complex line bundle over ${\mathcal M}_\Sigma$ with fiber
$$
{\mathcal
L}_A=\bigwedge^{\text{top}}\left(\text{ker}(d_A)\right)^*\otimes\bigwedge^{\text{top}}\left(\text{coker}(d_A)\right),
$$
where $d_A\co \Omega^{0,0}(\Sigma,E)\to\Omega^{0,1}(\Sigma,E)$ is
the covariant derivative associated to the flat connection $A$ on
the bundle $E$. The associated Hilbert space is then ${\mathcal
H}_\Sigma:=H^0({\mathcal M}_\Sigma, {\mathcal L}^{\otimes k+N})$.
One should note that sometimes the quantity $k+N$ is called the
level and is sometimes denoted by $k$. For general Lie groups the
two different notions of level are related by the so-called dual
Coxeter number. In this paper we will always use $k$ to be the level
as used in the definition of the string coupling constant. We
discuss this in greater detail in \fullref{app:e}. More information
on the above description of Witten--Chern--Simons theory may be
found in Axelrod--Della Pietra--Witten \cite{ADW}, Hu \cite{hu}, Di
Francesco--Mathieu--S\'en\'echal \cite{DMS} and Kohno \cite{kohno}.}
\end{aside}
While the motivation for this approach is fairly straightforward,
formally constructing invariants in this manner is very complicated.
There is an alternative approach that is more difficult to motivate
but slightly easier technically.

\subsection{Link invariants in a $U(1)$ theory}\label{u1link}
Instead of cutting $3$--manifolds into cobordisms we use the fact
that every $3$--manifold can be expressed as surgery on a framed link
to reduce our considerations to framed links. Just as every
$3$--manifold can be cut into cobordisms every framed link can be cut
into elementary framed tangles. This is easier to draw and
conceptualize.

To pass from a TQFT describing $3$--manifold invariants to a
corresponding theory for framed tangles the notion of a
$(2{+}1)$--dimensional TQFT was enhanced by Reshetikhin, Turaev and
many others \cite{RT1}. One includes Wilson loops (framed links)
into the manifolds. Formalizing the entire picture with Wilson loops
in general is fairly complicated because one has to consider framed
tangles in general $3$--manifolds with surfaces of arbitrary genus as
boundaries.

\begin{figure}[ht!]
\centering
\labellist\small
\pinlabel {bottom} [b] at 300 0
\pinlabel {top} [b] at 1020 0
\pinlabel {\rotatebox{-90}{$p$}} at 333 63
\pinlabel {\rotatebox{-90}{$-p$}} at 333 223
\pinlabel {\rotatebox{-90}{$-q$}} at 385 95
\pinlabel {\rotatebox{-90}{$q$}} at 385 200
\pinlabel {\rotatebox{-90}{$q$}} at 435 170
\pinlabel {\rotatebox{-90}{$p$}} at 435 235
\pinlabel {\rotatebox{-90}{$p$}} at 520 170
\pinlabel {\rotatebox{-90}{$q$}} at 520 235
\pinlabel {\rotatebox{-90}{$q$}} at 620 195
\pinlabel {\rotatebox{-90}{$p$}} at 805 240
\pinlabel {\rotatebox{-90}{$-q$}} at 695 120
\pinlabel {\rotatebox{-90}{$q$}} at 695 190
\pinlabel {\rotatebox{-90}{$q$}} at 760 120
\pinlabel {\rotatebox{-90}{$-q$}} at 760 200
\pinlabel {\rotatebox{-90}{$-p$}} at 910 120
\pinlabel {\rotatebox{-90}{$p$}} at 910 190
\pinlabel {\rotatebox{-90}{$p$}} at 988 120
\pinlabel {\rotatebox{-90}{$-p$}} at 988 205
\endlabellist
\includegraphics[width=5truein]{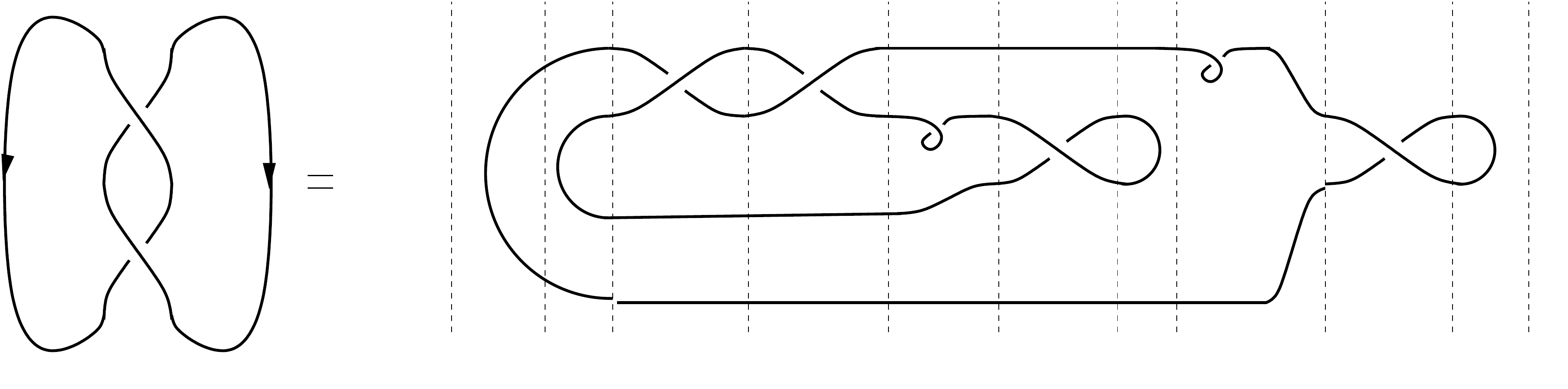} \caption{The left
Hopf link}\label{lhopf}
\end{figure}
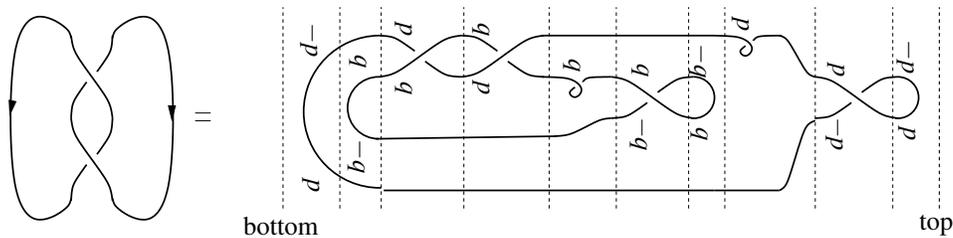

It is easier to formalize this picture for the special case of
framed tangles in $S^2\times [0,1]$. These are usually viewed as
framed tangles in $\R^2\times [0,1]$ as in \fullref{lhopf}. This
proves to be sufficient because any framed link can be obtained by
stacking such tangles and any $3$--manifold can be obtained by
surgery on a framed link.

Before describing the general case in the next section we are going
to introduce the cutting and pasting idea in this section in the
setting of U$(1)$ Chern--Simons theory. In the U$(1)$ theory the
cubic term in the Chern--Simons action vanishes because purely
imaginary $1$--forms anti-commute. Thus one expects considerable
simplification in the U$(1)$ case.

According to the discussion on topological quantum field theories
with framed tangles, one should be able to cut a framed link into
elementary pieces and associate invariants to each piece. This is
indeed the case. Any framed link can be cut into cups, caps, right
crossings, left crossings and twists (see \fullref{lhopf} for an
example and \fullref{mtcsec} for formal definitions).

The simple U$(1)$ theory at level $2m+1$ produces an invariant of
colored, oriented framed links. Appropriately drawn link diagrams
can be oriented by putting an upward pointing arrow on the left
branch of each cup. Right and left crossings are then defined so
that the crossings in the simple diagram for the left Hopf link from
\fullref{lhopf} are left crossings and the opposite crossings are
right crossings. Left and right \index{left crossing} \index{right
crossing} refer to the direction of the strand that goes under the
crossing. Coloring is captured by labeling every cup with a number
$p$ from $\Z_{2m+1}$ subject to a compatibility condition.
The evaluation rule is very simple. In the obvious analogy with
Feynman rules, every $p$--labeled cup on a framed link creates a pair
$\pm p$ with $p$ on the left and $-p$ on the right. The rest of the
link receives labels by extension once we choose a `color' $p$ for
each cup (see \fullref{lhopf}). The compatibility condition is
that every cap annihilates such a $\pm p$ pair.
 Once the link is
colored, the $\text{U}(1)$ theory  at level $2m+1$  associates the
following numbers to elementary pieces:
\begin{center}
\begin{tabular}{rcllp{3em}rcll}
$|$ & $=$ & $1$ & id & &
$\times_{p,q}$ & $=$ & $e^{2\pi i pq/(2m+1)}$ & right crossing \\
$\cup$ & $=$ & $1$ & cup & &
$\times^{-1}_{p,q}$ & $=$ & $e^{-2\pi i pq/(2m+1)}$ & left crossing \\
$\cap$ & $=$ & $1$ & cap & &
$\theta_p$ & $=$ & $e^{2\pi i p^2/(2m+1)}$ & twist
\end{tabular}
\end{center}
The colors $p,q$ can be any natural numbers from $0$ to $2m$. The final colored invariant
$J_{p_1,\ldots,p_l}(L)$ of the framed link is just the product of
all numbers assigned to the pieces. In analogy to Feynman rules, the
quantum invariant $F(L)$ is just the sum of the colored Jones
polynomials over all possible ways to color the link.
\begin{example}\label{lhopfeg}
Consider the left Hopf link $L$ from \fullref{lhopf} with the
larger component colored by $p$ and the smaller one colored by $q$.
Multiplying from the bottom up we compute
$$
J_{p,q}(L)=\times^{-1}_{p,q}\cdot\times^{-1}_{q,p}\cdot\theta_q\cdot\times_{q,-q}\cdot\theta_p\cdot\times_{p,-p}
=e^{-4\pi i pq/(2m+1)}\,.
$$
Therefore, $F(L)=\sum^{2m}_{p,q=0}e^{-4\pi i pq/(2m+1)}=2m+1$.
\end{example}
The `color and multiply' rule is not very sophisticated. In
particular, it implies that any two links with the same elementary
pieces have the same invariant no matter how those pieces are
assembled. Of course, this makes it easier to prove independence of
the presentation of a link by a regular projection but it also
yields rather weak invariants. We want evaluation rules that are
still independent of the projection but `see' much more structure of
the link.  The right balance for link invariants is struck in the
notion of ribbon categories. There colors are replaced by simple
objects (think irreducible representations), numbers by morphisms
(think linear maps) and multiplication by composition and tensor
product. Accordingly, the evaluation rules become more involved.

To get invariants of $3$--manifolds one takes weighted sums of
colored invariants such as the one computed in the previous
example. The hard part is to make sure that the final answer does
not depend on how the manifold is expressed as a framed link and how
the framed link is decomposed into elementary tangles. The required
axioms are formalized in the notion of a strict modular category. In
the rest of this subsection we introduce ribbon categories first,
then discuss modular categories and evaluation rules, and finally
the quantum invariants that arise from these categories.

\subsection{Ribbon categories}\label{catrib}

As explained in the previous subsection the invariants that we are
considering are naturally defined for framed tangles in $\R^2\times
[0,1]$ and satisfy formal gluing axioms. We will follow the version
of these invariants due to Reshetikhin and Turaev \cite{RT2} as
explained by Bakalov and Kirillov \cite{BK} and Turaev \cite{T}.

It is important to keep the basic idea in mind. Every $3$--manifold
may be expressed as surgery on a framed link.
\begin{quote}
Associate a simple invariant to each elementary piece of a framed
link diagram and define the final invariant to be an algebraic
combination of all of the elementary pieces. Of course, one must
know that different ways to assemble the same manifold give the same
final invariant.
\end{quote}
Axiomatizing assembly rules with outcomes independent of a link
presentation leads to the notion of modular categories.

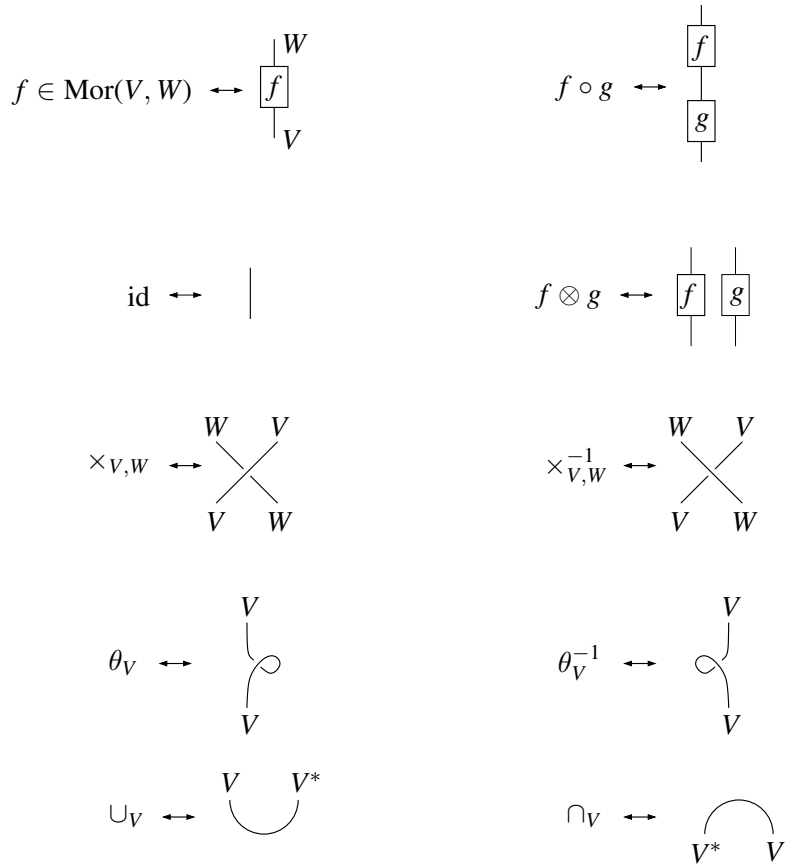
\begin{figure}[ht!]
\centering
\labellist
\pinlabel {$f\in\Mor(V,W)$} [r] at 135 630
\pinlabel {$f$} at 198 633
\pinlabel {$W$} [l] at 196 670
\pinlabel {$V$} [l] at 196 590
\pinlabel {$f\circ g$} [r] at 485 633
\pinlabel {$f$} at 552 667
\pinlabel {$g$} at 552 605
\pinlabel {$\id$} [r] at 100 460
\pinlabel {$f\otimes g$} [r] at 475 460
\pinlabel {$f$} at 542 460
\pinlabel {$g$} at 581 460
\pinlabel {$\times_{V,W}$} [r] at 100 320
\pinlabel {$W$} [b] at 148 337
\pinlabel {$V$} [b] at 201 337
\pinlabel {$V$} [t] at 148 287
\pinlabel {$W$} [t] at 201 287
\pinlabel {$\times_{V,W}^{-1}$} [r] at 480 318
\pinlabel {$W$} [b] at 533 337
\pinlabel {$V$} [b] at 587 337
\pinlabel {$V$} [t] at 533 287
\pinlabel {$W$} [t] at 587 287
\pinlabel {$\theta_V$} [r] at 90 155
\pinlabel {$V$} [b] at 176 187
\pinlabel {$V$} [t] at 176 118
\pinlabel {$\theta_V^{-1}$} [r] at 475 155
\pinlabel {$V$} [b] at 576 187
\pinlabel {$V$} [t] at 576 118
\pinlabel {$\cup_V$} [r] at 95 27
\pinlabel {$V$} [b] at 160 40
\pinlabel {$V^*$} [b] at 223 40
\pinlabel {$\cap_V$} [r] at 475 27
\pinlabel {$V^*$} [t] at 555 12
\pinlabel {$V$} [t] at 612 12
\endlabellist
\includegraphics[width=4truein]{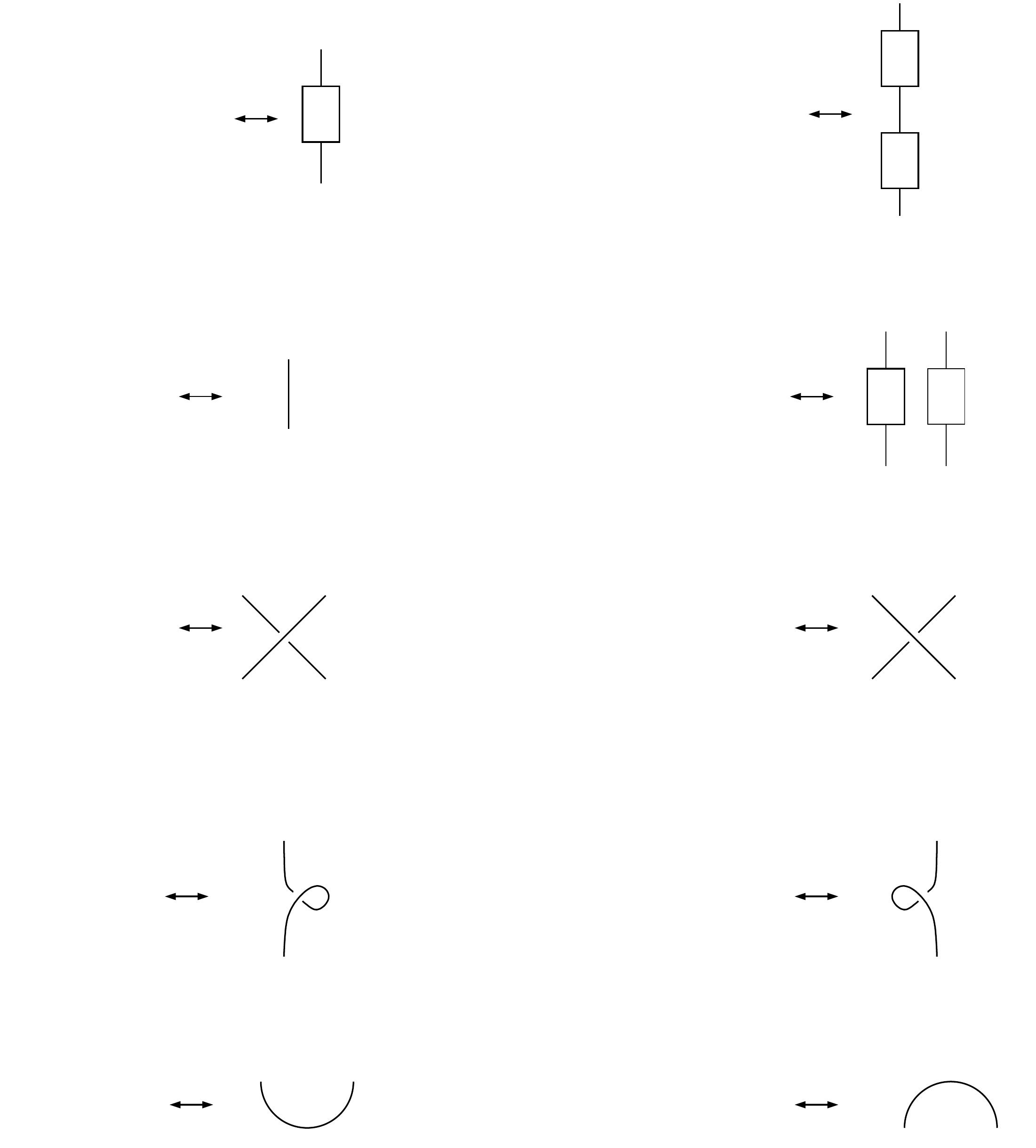} \caption{Strict
ribbon category}\label{MTC}
\end{figure}

The basic outline is as follows. Any framed link may be constructed
out of elementary building blocks that look like $\times_{U,V}$,
$\theta_V$, $\cap_V$ and $\cup_V$ from \fullref{MTC}. The trick is
to associate an invariant to each of these and then write out all of
the axioms that correspond to changing the height function isotopy
of the link, or Kirby moves. In particular, one colors and orients
the link. To each generic horizontal line one associates a tensor
product of objects and their duals according to the sign of the
intersection of the link with the horizontal line. The elementary
building blocks between the horizontal lines induce morphisms
between the associated objects. We use the standard convention that
combining elementary building blocks side-by-side corresponds to
taking the tensor product of the associated morphisms.

There are many types of categories related to modular
categories. It is helpful to consider a related category with an
easier definition first. The categories related to framed link
invariants (as opposed to $3$--manifold invariants) are ribbon categories.
\begin{defn}
A strict ribbon category \index{strict ribbon category} is a
category with a unit object ${\bbone}$, tensor product functor
$\otimes$, families of isomorphisms $\times$, $\theta$ called the
\index{$\times$ braiding} braiding \index{braiding} and the
\index{$\theta$ twist} twist \index{twist} respectively, and a
duality triple \index{duality triple} $(*,\cup,\cap)$
\index{$(*,\cup,\cap)$ duality triple} satisfying axioms 1 through
12 under \eqref{axmtc} below.
\end{defn}

\begin{remark}
We are slightly changing notation from Turaev \cite{T} and other
references. The correspondence is $\times_{U,V}=c_{U,V}$,
$\theta_V=\theta_V$, $\cup_V=b_V$ and $\cap_V=d_V$.
In fact, our notation is often used to represent objects in a colored ribbon category.
There is a natural functor taking the colored ribbon morphisms to morphisms in
a strict modular category. We feel that no confusion will arise by
using the same notation for both, and the ribbon notation is more
descriptive.
\end{remark}

\begin{remark}
One can obtain a new strict ribbon category by replacing the
braiding and twist by their inverses.
\end{remark}

The first example of a strict ribbon category is the category of
representations of a group.
\begin{example}\label{repg}
Let $\text{REP}_G$ be the category of representations of a Lie
group. The objects are representations $\rho\co G\to \text{Aut}(V)$.
The unit object is the trivial representation. Morphisms are
equivariant linear maps $f(gv)=gf(v)$. The dual representation $V^*$
has the dual space to $V$ as the representation space and the action
is given by $(g\varphi)(v):=\varphi(g^{-1}v)$. The pairing $\cup$ is
the standard duality pairing for vector spaces and $\cap$ the
copairing given by $1\mapsto \sum e_k\otimes e^k$ where $\{e_k\}$
and $\{e^k\}$ are dual bases. The image of $1$ under the copairing
is sometimes called the Casimir element. The tensor product is the
standard one in the category of finite-dimensional vector spaces
with action given by tensor product of actions as described in
\fullref{app:d}. The braiding is given by $\times_{V,W}(v\otimes
w)=w\otimes v$, and the twist is given by $\theta_V=\text{id}_V$. As
is, this is not a strict ribbon category because equality signs in
the axioms are only canonical isomorphisms. For example $V$ is not
equal to $V\otimes \C$. In truth, we should be talking about
equivalence classes of representations rather than representations
themselves. This requires some tweaking in the notion of morphisms
and ribbon operations that can be done in a standard way by the
Mac\,Lane coherence theorem \cite{mac}.
\end{example}
The category of representations cannot be used to construct
nontrivial link invariants following the procedure given in the next
article because the braiding is its own inverse (this category
cannot detect the difference between right and left crossings.) The
nontrivial invariants described in the previous subsection can be
seen to arise from a strict ribbon category.
\begin{example}\label{u1categ}
Construct a category $\calU(1)_{2m+1}$ \index{$U$@ U$(1)$
ribbon/modular category} with objects elements of $\Z_{2m+1}$, unit
object equal to $0$, only the zero morphism between unequal objects
and the endomorphisms of any object equal to $\C$. The tensor
product is given by
$$(f\co p\to q)\otimes(g\co r\to s):=fg\co p+r\to q+s,$$
and the braiding given by
$$\times_{p,q}=e^{2\pi ipq/(2m+1)}\co p+q\to q+p.$$
The twist is given by
$$\theta_{p}=e^{2\pi ip^2/(2m+1)}\co p\to p.$$
The duality pairing is $\cap_p\co  p+(-p)\to 0$ and the copairing
is $\cup_p=1\co 0\to p+(-p)$. One sees that this weird category is
strict and satisfies all of the axioms for a strict ribbon category.
\end{example}
\begin{exm}
Prove that in a strict ribbon
category $\times_{V,{\bbone}}=\text{id}_V=\times_{{\bbone},V}$ and
the following Yang--Baxter equation holds \index{Yang--Baxter equation}
\begin{multline*}
(\times_{V,W}\otimes\text{id}_U)\circ(\text{id}_V\otimes\times_{U,W})\circ
(\times_{U,V}\otimes\text{id}_W)\\
=(\text{id}_W\otimes\times_{U,V})\circ(\times_{U,W}\otimes\text{id}_V)\circ(\text{id}_U\otimes\times_{V,W})
\,.
\end{multline*}
\end{exm}

There is a way to represent the axioms and other formulas with
ribbon operations graphically making them much more intuitive.
\index{strict ribbon category}
The correspondence between elementary operations and graphs is
depicted on \fullref{MTC}. To draw the picture corresponding to a
formula start from the right (as in reading Arabic) and draw the
pieces corresponding to each expression between two compositions on
the same horizontal level while moving from the bottom up and
putting the obtained pieces on top of each other. Conversely, given
a labeled graph draw horizontal lines not intersecting any
crossings, twists, maxima or minima (see \fullref{lhopf}). Next
write the expression for each elementary piece on the same level and
tensor them. Finally assemble the expressions right to left by
compositions moving from the bottom to the top of the graph. The
graphs cannot be labeled arbitrarily; the labels on successive
horizontal levels must match. If a ribbon expression produces a
graph with mismatching labels it is nonsensical: in categorical
language you would be trying to compose morphisms with targets and
sources that do not match.

The main advantage of using graphs is that one can immediately see
if two expressions in a ribbon category are equal.
\begin{quote} If two (correctly composed) expression graphs represent isotopic framed
tangles the expressions are in fact equal and the equality can be
established using the axioms; see Bakalov and Kirillov \cite{BK}.
\end{quote}
Framing is important here.  Even though we draw pictures with
strands one should actually think of them as very thin ribbons so
that the twists (depicted as curls on the strands) cannot be undone.
Later we will slightly enhance the notation to allow arrows on the
strands but for now this will suffice.

\begin{figure}[ht!]
\centering
\labellist
\pinlabel {$\times$} [b] at 25 0
\pinlabel {$\theta$} [b] at 175 0
\pinlabel {$\cap$} [b] at 345 0
\pinlabel {$\cup$} [b] at 500 0
\endlabellist
\includegraphics[width=4truein]{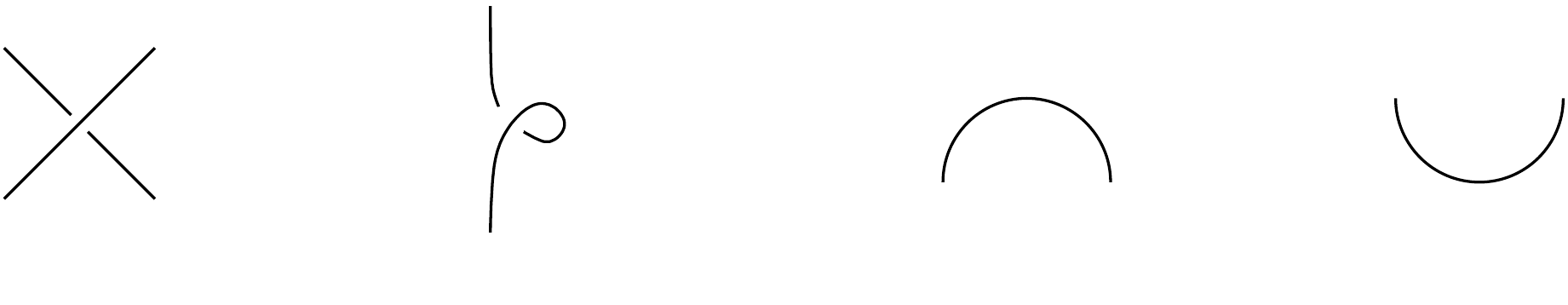} \caption{Category
of ribbon tangles}\label{ribtan}
\end{figure}
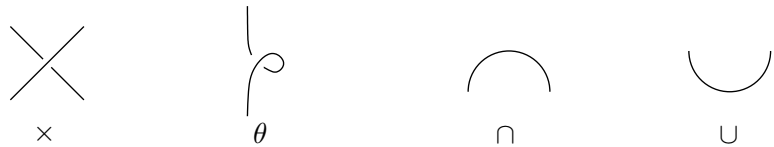

The motivating example of a ribbon category is the category of
ribbon tangles.
\begin{example}
The objects of the category of ribbon tangles \index{category of
ribbon tangles} are just non-negative integers, with zero
representing the unit object. The morphisms between $n$ and $m$ are
just the isotopy classes of framed tangles from $n$ points to $m$
points. Our graphs of expressions are just plane projections (with
indication of under- and over-crossings and twists) of framed
tangles with labels being natural numbers. The tensor product is the
sum and the dual of $n$  is $n$ itself. The braiding, twist and
duality for single strands are displayed in \fullref{ribtan}; in
general one just has to put $n$ and $m$ strands parallel to the one
or two depicted. Invariants produced by this category are complete
but useless: the invariant of an isotopy class is the isotopy class
itself. Fortunately, there are more interesting examples.
\end{example}

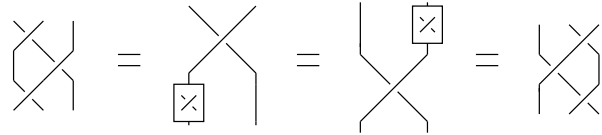
\begin{figure}[ht!]
\centering
\labellist
\pinlabel {Axiom 6} [b] at 125 200
\pinlabel {$U\otimes V$} [t] at 33 290
\pinlabel {$W$} [t] at 99 290
\pinlabel {$U$} [t] at 162 290
\pinlabel {$V$} [t] at 185 290
\pinlabel {$W$} [t] at 210 290
\pinlabel {Axiom 7} [b] at 468 200
\pinlabel {$f$} at 372 275
\pinlabel {$g$} at 425 275
\pinlabel {$g$} at 517 350
\pinlabel {$f$} at 570 350
\endlabellist
\includegraphics[width=4truein]{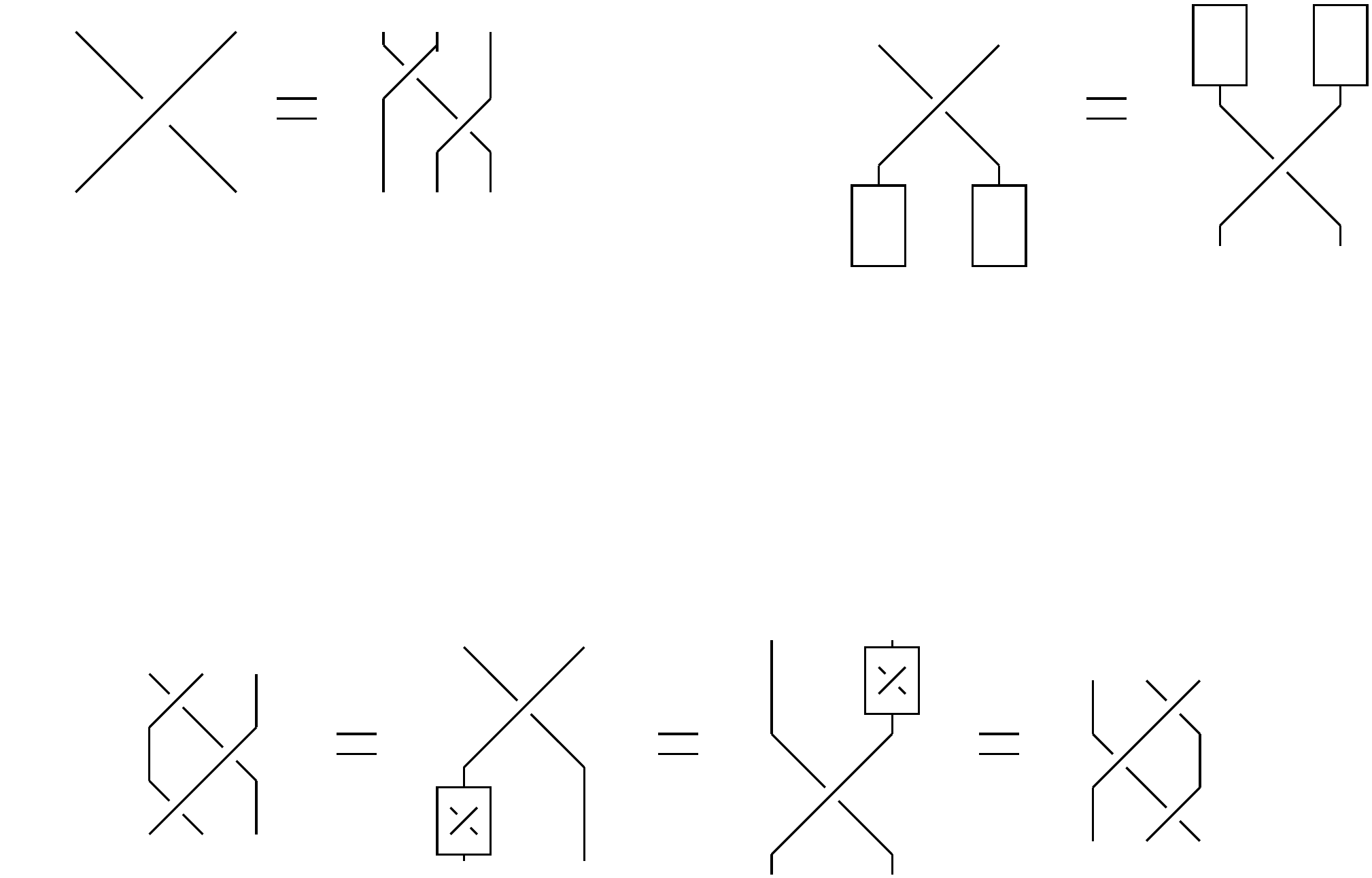} \caption{Axioms 6 and
7}\label{axioms}
\end{figure}

The pictures for axioms 6 and 7 are displayed in \fullref{axioms}. These axioms correspond exactly to elementary
isotopies. For example, a combination of axioms 6 and 7 implies the
third Reidemeister move as shown in \fullref{axioms}.

\begin{exm}
Draw diagrams representing each of the braiding, twist and duality
axioms. Remember that it is standard to compose maps from the bottom up in a
diagram.
\end{exm}

\subsection{Modular categories}\label{mtcsec}

While looking at axioms 1--12 you may have noticed that the
remaining ones 13--17 have a different flavor. Rather than
describing algebraic properties of operations they describe the
global structure of a category. The idea is that the first twelve
axioms take care of the Reidemeister moves and produce link
invariants. Surgeries on isotopic links produce diffeomorphic
$3$--manifolds but diffeomorphic $3$--manifolds are also produced by
non-isotopic links related by Kirby moves. The invariants we are
looking for must take the same values on such pairs of links. The
Kirby move has a more complex structure than the Reidemeister moves
and one requires the modular axioms to account for it.
\begin{quote}
The point is that any ribbon category generates many framed link
invariants. There is a framed link invariant for each object in the
category or more generally one can define invariants of framed links
colored by objects in the category. The hope is that an appropriate
linear combination of the resulting framed link invariants will be
invariant under the Kirby move and thus define a $3$--manifold
invariant.
\end{quote}

Before stating the axioms of a strict modular category we need to
define a few terms used in the axioms. {\sf Axiom 13} simply says
that we can add morphisms with the same sources and targets and this
addition behaves the same way as linear operators on vector spaces.
Categories satisfying axiom 13 are called preabelian. Denote
$\text{End}(V):=\text{Mor}(V,V)$ and notice that
$(\text{End}({\bbone}),+,\circ)$ is a ring. Moreover, it is a
commutative ring.

\begin{exm}
Use $f\circ g=(f\otimes \text{id}_{\bbone})\circ(\text{id}_{\bbone}\otimes
g)$ to prove that End$({\bbone})$ is commutative based on axioms 1--13.
Hint: recall that $\otimes$ is functorial.
\end{exm}

\begin{exm}
Prove that $\theta_{{\bbone}}^2=\theta_{{\bbone}}$ and
$$
\theta_{V\otimes
W}=\times_{W,V}\circ\times_{V,W}\circ(\theta_V\otimes\theta_W)\,.
$$
\end{exm}

\begin{figure}[ht!]
\centering \labellist \pinlabel {\rotatebox{-90}{$f$}} at 65 77
\pinlabel {bottom} [b] at 30 10 \pinlabel {top} [b] at 180 10
\endlabellist
\includegraphics[width=3truein]{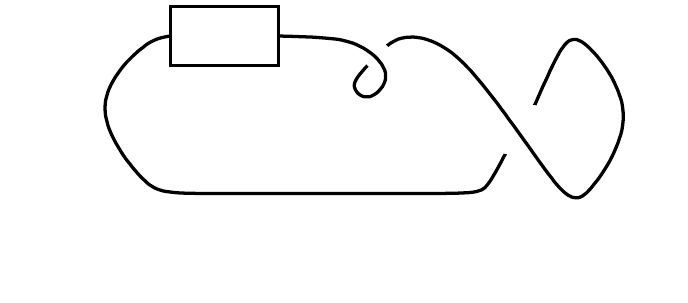} \caption{The
quantum trace}\label{qtrace}
\end{figure}

In the light of the above we will sometimes omit the composition
sign between morphisms when composing them. The analogy with
vector spaces goes further: one can define `traces' of endomorphisms
and `dimensions' of objects.
\begin{defn}\label{Qtrace}
The quantum trace of $f\in\text{End}(V)$ is \index{quantum
trace}\index{$T$@$\text{Tr}_q(f)$ quantum trace}
$$
\text{Tr}_q(f):=\cap_V\circ\times_{V,V^*}\circ(\theta_V\otimes\text{id}_V)
\circ(f\otimes\text{id}_V)\cup_V\in\text{End}({\bbone})\,.
$$
The corresponding graph is shown in \fullref{qtrace}. The quantum
dimension of an object is \index{quantum dimension}
\index{$D$@$\text{dim}_q(V)$ quantum dimension}
$\text{dim}_q(V):=\text{Tr}_q(\text{id}_V)$.
\end{defn}
\begin{exm}\label{qtrprop}
Prove that $\text{Tr}_q(fg)=\text{Tr}_q(gf)$ for all $f$ and $g$
where the compositions make sense, $\text{Tr}_q(f)=f$ for all
$f\in\text{End}({\bbone})$, and $\text{Tr}_q(f\otimes
g)=\text{Tr}_q(f)\text{Tr}_q(g)$ for all morphisms $f$ and $g$.
Conclude that for any pair of objects $\text{dim}_q(V\otimes
W)=\text{dim}_q(V)\text{dim}_q(W)$.
\end{exm}
\begin{example}\label{tracedim}
In the category $\text{REP}_G$ of \fullref{repg} we have
\def\Sum{\sum\nolimits}
\begin{align*}
\text{Tr}_q(f)&=\cap_V\circ\times_{V,V^*}\circ(\theta_V\otimes\text{id}_V)
\circ(f\otimes\text{id}_V)\Bigl(\Sum_k e^k\otimes e_k\Bigr)\\
&=\cap_V\circ\times_{V,V^*}\circ(\theta_V\otimes\text{id}_V)\Bigl(\Sum_k
f(e^k)\otimes e_k\Bigr)\\
&=\cap_V\Bigl(\Sum_k e_k\otimes f(e^k)\Bigr) =\Sum_k \langle e_k,
f(e^k)\rangle=\text{Tr}(f)
\end{align*}
so the quantum trace reduces to the ordinary one and
$\text{dim}_q(V)=\text{dim}(V)$.
\end{example}
\begin{exm}
In the category $\calU(1)_{2m+1}$ of \fullref{u1categ} show that
$\text{Tr}_q(f)=f$.
\end{exm}
In preabelian categories one can define an analog of irreducible representations.
\begin{defn}
To say that an object $V$ is {\it simple} means that the map
\index{simple}
 End$({\bbone})\to\text{Mor}(V,V)$ given by $f\mapsto
f\otimes\text{id}_V$ is an isomorphism. To say that such a category
\index{dominated} is {\it dominated by} $\{V_\lambda\}_{\lambda\in
I}$ means that for every $V\in\text{Ob}({\mathcal V})$ there are
morphisms $f_r\co V_{\lambda_r}\to V$ and $g_r\co V\to V_{\lambda_r}$ such
that $\text{id}_V=\sum f_r\circ g_r$.
\end{defn}
\begin{remark}
Domination is the preabelian analog of semi-simplicity in Abelian
categories. The point is that while one can take sums of objects in
Abelian categories one can only take sums of morphisms in preabelian
categories.
\end{remark}
The Schur lemma implies that in the category of representations of a
Lie group $\text{REP}_G$ the simple objects are the irreducible representations.
However, any nontrivial Lie group has infinitely many simple
objects. The category of finite-dimensional representations of a
semi-simple Lie group is dominated by the simple objects
because any representation decomposes into irreducible ones.
In the category from \fullref{u1categ} every object is simple
and there are only finitely many objects. Clearly, this category is
dominated by its simple objects.

The fact that $\text{End}({\bbone})$ is a commutative ring allows one to
define the matrix with the following entries (to be used in axiom
17):
\begin{equation}\label{smat}
\wtilde{s}_{\lambda\mu}:=\text{Tr}_q(\times_{\mu,\lambda}\circ\times_{\lambda,\mu})
\end{equation}
From \fullref{qtrprop} we see that $\wtilde{s}$ is symmetric.
\begin{exm}\label{smatex}\index{$S$@$\tilde s$ quantum $s$ matrix}
Show that the $\tilde s$--matrix defined in \eqref{smat} is
\begin{multline*}
\tilde
s_{\lambda\mu}=\cap_{\lambda\otimes\mu}\circ\times_{\lambda\otimes\mu,\mu^*\otimes\lambda^*}
\circ
(\theta_{\lambda\otimes\mu}\otimes\text{id}_{\mu^*\otimes\lambda^*})
\circ(\times_{\mu\otimes\lambda}\otimes
\text{id}_{\mu^*\otimes\lambda^*})\\
\circ (\times_{\lambda\otimes\mu}\otimes
\text{id}_{\mu^*\otimes\lambda^*}) \circ\cup_{\lambda\otimes\mu} \in
\text{End}({{\bbone}}).
\end{multline*}
 A graphical representation of the $\tilde s$--matrix is given by
the right Hopf link as displayed in \fullref{mtcchar}.
\end{exm}
We are now ready to define a strict modular category. We refer our
readers to Bakalov and Kirillov \cite{BK} and to Turaev \cite{T} for
further description, properties and related definitions. The reference
\cite{BK} uses modular tensor categories as opposed to strict modular
categories. The difference is that the underlying category in a modular
tensor category is Abelian rather than just preabelian, that is, direct
sums of objects are defined. Even though our main example is in fact
a modular tensor category, at this point it is expedient to ignore the
additional structure.

\begin{defn}
A strict modular category \index{strict modular category} is a
category with a tensor product functor, natural isomorphisms
$\times$, $\theta$, a duality $(*,\cup,\cap)$ together with an
indexed collection of special objects $\{V_\lambda\}_{\lambda\in I}$
satisfying (all 17 of) the axioms under \fullref{axmtc} below.
\end{defn}
In particular any strict modular category is a strict ribbon
category. The ingredients in a strict modular category are listed in
\fullref{mtci}. This table also includes the ingredients of the
$\calU(1)_{2m+1}$ category from \fullref{u1categ} as an
illustration.
\begin{table}
\begin{small}\label{mtci}
\begin{tabular}{|l|l|} \hline
Strict modular category ingredients & Example  \\ \hline
A category $\mathcal V$ & Ob$({\calU(1)_{2m+1}})=\Z_{2m+1}$\\
 & Mor$(p,q)= \C$ if $p=q$ and $0$ otherwise \\
A tensor product $\otimes\co {\mathcal V}\times{\mathcal V}\Rightarrow
{\mathcal V}$ & $(f\co p\to q)\otimes (g\co r\to s)=fg\co p+r\to q+s$ \\
A unit ${{\bbone}}\in \text{Ob}({\mathcal V})$ & $0\in\Z_{2m+1}$ \\
A braiding $\times_{U,V}\co U\otimes V\to V\otimes U$
& $\times_{p,q}=e^{2\pi ipq/(2m+1)}\co p+q\to q+p$ \\
A twist $\theta_V\co V\to V$ & $\theta_{p}=e^{2\pi ip^2/(2m+1)}\co p\to p$ \\
A duality pairing $\cap_V\co V^*\otimes V\to {{\bbone}}$ & $\cap_p\co  p+(-p)\to 0$ \\
A copairing $\cup_V\co {{\bbone}}\to V\otimes V^*$ & $\cup_p=1\co 0\to p+(-p)$ \\
A finite collection of simple objects &
$V_p=p$, $I=\Z_{2m+1}$ \\
$\{V_\lambda\}_{\lambda\in I}$ & \\ \hline
\end{tabular}
\end{small}
\vskip.1in
\caption{Strict modular category ingredients}
\end{table}

\subsection{Axioms defining a strict modular category}
\label{axmtc}
\quad

\textbf{Tensor axioms}
\begin{description}
\item[Axiom 1] ${V}\otimes{\bbone}={\bbone}\otimes V=V$.
\item[Axiom 2] $U\otimes (V\otimes W)=(U\otimes V)\otimes W$.
\item[Axiom 3] $f\otimes\text{id}_{\bbone}=\text{id}_{\bbone}\otimes f=f$.
\item[Axiom 4] $f\otimes(g\otimes h)=(f\otimes g)\otimes h$.
\end{description}
\textbf{Braiding axioms}
\begin{description}
\item[Axiom 5] $\times_{U, V\otimes
W}=(\text{id}_V\otimes\times_{U,W})\circ(\times_{U,V}\otimes\text{id}_W)$.
\item[Axiom 6] $\times_{U\otimes V,W}=(\times_{U,W}\otimes\text{id}_V)\circ
(\text{id}_U\otimes\times_{V,W})$.
\item[Axiom 7] $(g\otimes
f)\circ\times_{U,W}=\times_{V,Z}\circ(f\otimes g)$ for any morphisms $f\co
U\to V$, $g\co W\to Z$.
\end{description}
\textbf{Twist axioms}
\begin{description}
\item[Axiom 8] $\theta_{V\otimes
W}=\times_{W,V}\circ\times_{V,W}\circ(\theta_V\otimes\theta_W)$.
\item[Axiom 9] $f\circ\theta_U=\theta_V\circ f$ for any morphism $f\co U\to V$.
\end{description}
\textbf{Duality axioms}
\begin{description}
\item[Axiom 10]
$(\text{id}_V\otimes\cap_V)\circ(\cup_V\otimes\text{id}_V)=\text{id}_{V}$.
\item[Axiom 11]
$(\cap_V\otimes\text{id}_{V^*})\circ(\text{id}_{V^*}\otimes\cup_V)=\text{id}_{V^*}$.
\item[Axiom 12]
$(\theta_V\otimes\text{id}_{V^*})\circ\cup_V=(\text{id}_V\otimes\theta_{V^*})\circ\cup_V$.
\end{description}
\textbf{Modular axioms}
\begin{description}
\item[Axiom 13] Mor$(V,W)$ are abelian groups and
$\circ\co \text{Mor}(V,W)\times\text{Mor}(U,V)\to\text{Mor}(U,W)$ is
bilinear.
\item[Axiom 14] ${\mathcal V}$ is dominated by a finite collection
$\{V_\lambda\}_{\lambda\in I}$.
\item[Axiom 15] There is $0\in I$ such that $V_0={\bbone}$.
\item[Axiom 16] For every $\lambda\in I$ there is a $\lambda^*\in I$
such that $V_{\lambda^*}\cong V_\lambda^*$.
\item[Axiom 17] The matrix $\tilde s_{\lambda,\mu}=\text{Tr}_q(\times_{\mu,\lambda}\circ\times_{\lambda,\mu})$
is non-singular.
\end{description}
The tensor axioms are not difficult to understand if one keeps the
example of vector spaces in mind in which case ${\bbone}$ is just
the underlying field. Graphical representations help one understand
the braiding, twist and duality axioms. As we emphacized in
\fullref{catrib} they simply catalogue elementary transformations of
tangle diagrams. In fact, one can forget about them and work with
diagrams directly. Axioms 7 and 9 just restate the naturality of
$\times$ and $\theta$ but we included them for the sake of
diagrammatic interpretation. For instance, axiom 9 means that one
can slide the twist through any morphism. The next example is
intended to help illustrate the modular axioms.
\begin{example}
The category $\text{REP}_G$ for an infinite group $G$ fails to be strict
modular for the trivial reason of having infinitely many simple
objects (irreducible representations). But even when $|G|<\infty$
this category is not strict modular. It comes very close though: the
only thing that goes wrong is the non-degeneracy axiom 17. Indeed,
for any pair of objects (representation spaces) we have by \fullref{tracedim}:
$$
\text{Tr}_q(\times_{W,V}\circ\times_{V,W})=\text{Tr}_q(\text{id}_{V\otimes
W}) =\text{dim}_q(V\otimes W) =\text{dim}(V)\text{dim}(W)
$$
and the structure matrix
$\wtilde{s}_{\lambda\mu}=\text{dim}(V_\lambda)\text{dim}(V_\mu)$
always has rank $1$. We get non-degeneracy for $|G|=1$ that is, the
trivial group, but this is a rather trivial example.

On the other hand, the category $\calU(1)_{2m+1}$ from \fullref{u1categ}
is both strict modular and nontrivial. In fact, every
object is simple so it is definitely dominated by simple objects,
and there are only $2m+1<\infty$ of them. By \fullref{tracedim}
$$\wtilde{s}_{pq}=\text{Tr}_q(e^{2\pi ipq/(2m+1)}\circ e^{2\pi
ipq/(2m+1)})=e^{4\pi ipq/(2m+1)}=z^{pq}$$
with $z:=e^{4\pi i/(2m+1)}$.
\end{example}
\begin{exm}
Verify that the category $\calU(1)_{2m+1}$ from  \fullref{u1categ} and the \fullref{mtci} satisfies the definition of a
a strict modular category. Hint: notice that $\text{det}\ z^{pq}$ is
a Vandermonde determinant.
\end{exm}
We extend the honor of being named a number to the elements of
$\text{End}({\bbone})$ which is, after all, a commutative ring. In all
examples of interest to us $\text{End}({\bbone})=\C$ anyway.
\begin{defn}\label{charnum}
The \emph{characteristic numbers} \index{characteristic numbers}
\index{$d_\lambda:=\text{Tr}_q(\text{id}_\lambda)$}
\index{$p_\lambda^\pm:=\text{Tr}_q(\theta^{\pm 1}_\lambda)$} of a
strict modular category are given by
$d_\lambda:=\text{Tr}_q(\text{id}_\lambda)$,
$p_\lambda^\pm:=\text{Tr}_q(\theta^{\pm 1}_\lambda)$, where
$\lambda$ indexes simple objects. The numbers
$p^\pm:=\sum_{\lambda\in I}
p^\pm_\lambda=\sum_{\lambda\in I} \theta^{\pm 1}_\lambda d_\lambda$
are called the twists and $\calD:=\bigl(\sum d_\lambda^2\bigr)^{\unfrac12}$
\index{$D$@$\calD:=\bigl(\sum d_\lambda^2\bigr)^{\unfrac12}$}
\index{quantum diameter} the quantum diameter (also rank or dimension; see
Brugui\`eres \cite{Brug} and M\"uger \cite{Mug}) of a category.
\end{defn}
\begin{remark}\label{qdiam}
Note that the characteristic `numbers' are defined in any ribbon
category not just modular categories. The numbers $p^\pm,\calD$ are
instrumental in making sure that framed link invariants defined by a
modular category are invariant under Kirby moves and hence define
$3$--manifold invariants. The square root we need to take to define
$\calD$ may not exist in the ring End$({\bbone})$ and, when it does
may not be unique. It is possible to extend an arbitrary strict
modular category so that this square root does exist and there is no
important difference between choosing any of the two roots (see
Turaev \cite{T}). Numerical values of invariants do however depend
on a choice of quantum diameter. In examples of interest to us
$\End({\bbone})=\C$ and $\calD^2$ is a positive real number so we
agree to always choose the positive square root. For instance, in
$\calU(1)_{2m+1}$ there are $2m+1$ objects of dimension one each so
$\calD=(2m+1)^{\unfrac12}$. For $\text{REP}_G$ with $G$ finite the sum of
the squares of dimensions of irreducible representations is $|G|$ by
the Burnside theorem (see Fulton and Harris \cite{fulton-harris}) so
$\mathcal{D}=|G|^{\unfrac12}$.
\end{remark}

\subsection{Coloring, double duals and the arrow convention}\label{color}

\begin{figure}[ht!]
\centering
\labellist\small
\pinlabel {$\cup_V^*$} [r] at 40 192
\pinlabel {$V^*$} [b] at 112 206
\pinlabel {$V$} [b] at 178 206
\pinlabel {$V^*$} [b] at 250 242
\pinlabel {$V$} [b] at 288 242
\pinlabel {$V$} [r] at 245 160
\pinlabel {$V^*$} [l] at 290 160
\pinlabel {$:=$} at 215 190
\pinlabel {$\cap_V^*$} [r] at 420 192
\pinlabel {$V$} [t] at 480 170
\pinlabel {$V^*$} [t] at 546 170
\pinlabel {$V$} [t] at 615 140
\pinlabel {$V^*$} [t] at 660 140
\pinlabel {$V^*$} [r] at 615 225
\pinlabel {$V$} [l] at 660 225
\pinlabel {$:=$} at 580 190
\pinlabel {$V$} [r] at 35 25
\pinlabel {$V^*$} [l] at 77 25
\pinlabel {$:=$} at 55 25
\pinlabel {$V$} [r] at 205 35
\pinlabel {$V^*$} [b] at 338 37
\pinlabel {$V$} [b] at 393 35
\pinlabel {$:=$} at 300 25
\pinlabel {$V$} [r] at 485 35 \pinlabel {$V$} [b] at 600 35
\pinlabel {$V^*$} [b] at 658 35 \pinlabel {$:=$} at 570 25
\endlabellist
\includegraphics[width=4truein]{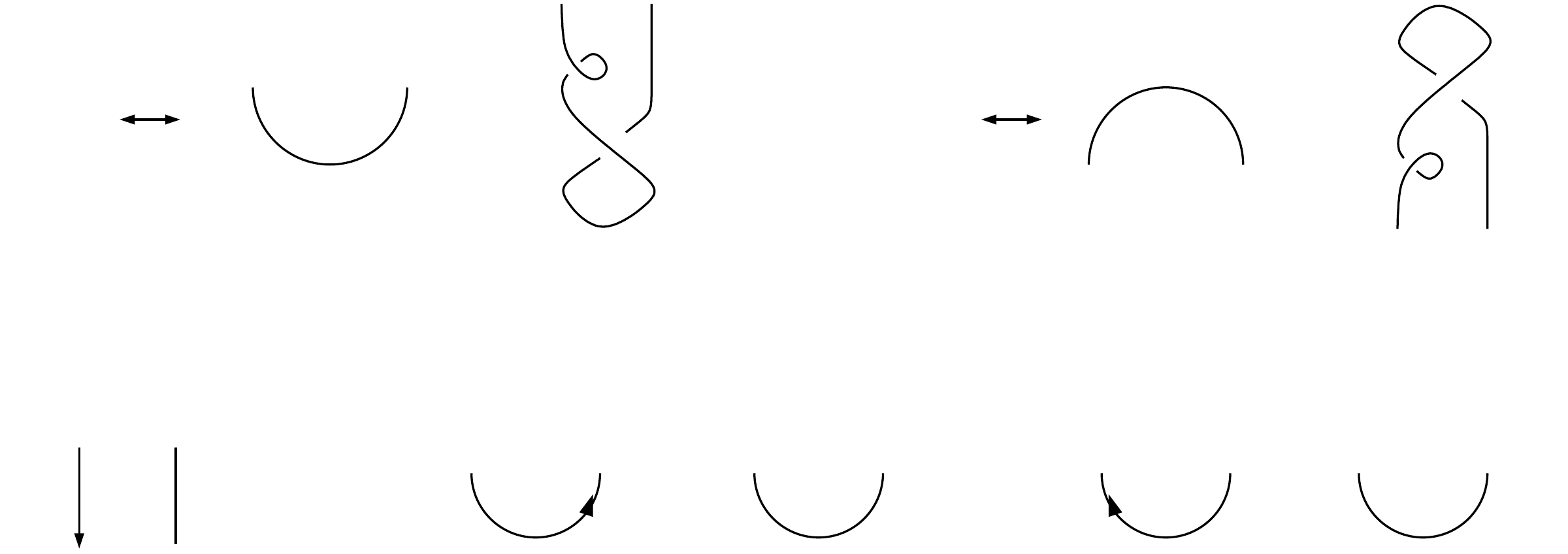} \caption{Dual
cups and caps and the arrow convention}\label{arrules}
\end{figure}
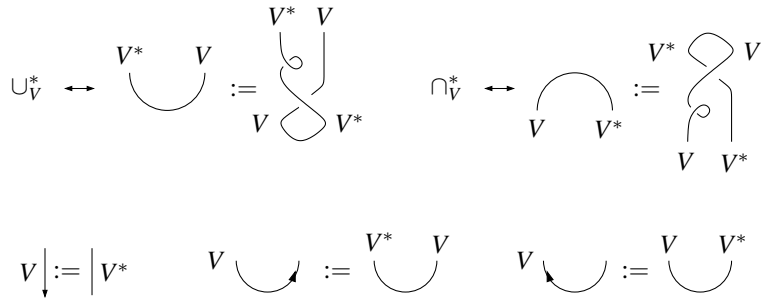

Before introducing invariants of links and $3$--manifolds we augment
the graphic notation introduced in \fullref{ribtan}. You may
notice that it is impossible to label a simple circle coherently
with the rules we have so far  because $\cup_V$ and $\cap_V$ put
together have mismatching labels $V,V^*$. We can isotope the circle
as in \fullref{qtrace} (without $f$) so that it can be labeled
consistently. Isotoping links into shapes that allow labeling every
time leads to cumbersome pictures as in \fullref{lhopf}. Instead
we can simply add a dual pairing and copairing (cups and caps) to
the notation. Namely, define $\cup^*_V\co {{\bbone}}\to V^*\otimes V$ and
$\cap^*_V\co V\otimes V^*\to{{\bbone}}$ by (see \fullref{arrules}):
$$
\begin{aligned}
\cup^*_V &=(\theta^{-1}_{V^*}\otimes\text{id}_V)\circ\times^{-1}_{V,V^*}\circ\cup_V\\
\cap^*_V &=\cap_V\circ\times_{V,V^*}\circ(\theta_V\otimes\text{id}_{V^*}).
\end{aligned}
$$
Using the dual pairing and copairing  \index{dual
pairing}\index{dual copairing} \index{$\cap^*_V$ dual
pairing}\index{$\cup^*_V$ dual copairing} we can re-express the
quantum trace of $f\co V\to V$ as
$\text{Tr}_q(f)=\cap^*_V\circ(f\otimes\text{id}_{V^*})\circ\cup_V$
and simplify all the definitions that use it.

\begin{exm}\label{dcup}
Show that $\cup^*_V$ and $\cap^*_V$ in the category
$\calU(1)_{2m+1}$ from \fullref{u1categ} are multiplications by
$1$.
\end{exm}
\begin{figure}[ht!]
\centering
\labellist\small
\pinlabel {$\theta = $} [r] at 45 75
\pinlabel {$\delta$} at 85 90
\endlabellist
\includegraphics[width=1.4truein]{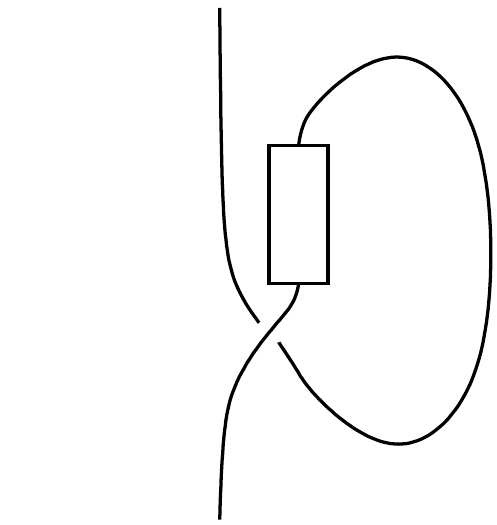}
\caption{Twist from the double dual isomorphism}\label{twistfromdual}
\end{figure}
Using $\cup^*_V$ and $\cap^*_V$ we can label any diagram
consistently with $V$ and $V^*$ only. The final simplification is to
get rid of $V^*$ as well by placing arrows on the strands and
agreeing that an up arrow  on a strand in a $V$--labeled component
corresponds to the $V$ label and a down arrow on the same component
corresponds to the $V^*$ label. With the arrow convention all we
have to do in order to label an entire link or tangle is to orient
every component and label it with an object in a single place, that is,
{\it color} it. To evaluate a colored graph one simply has to
transform it into the ribbon expression according to the rules on
Figures \ref{arrules} and \ref{ribtan} and tensor and compose all the
morphisms.
\begin{exm}\label{colmult}
Show that for the category $\calU(1)_{2m+1}$ from \fullref{u1categ} this reduces to the `color and multiply' rule of
\fullref{u1link}.
\end{exm}
\begin{figure}[ht!]
\centering
\labellist\small
\pinlabel {$\delta$} at 207 60
\pinlabel {$\delta$} at 379 37
\pinlabel {$\delta$} at 558 77
\endlabellist
\includegraphics[width=4truein]{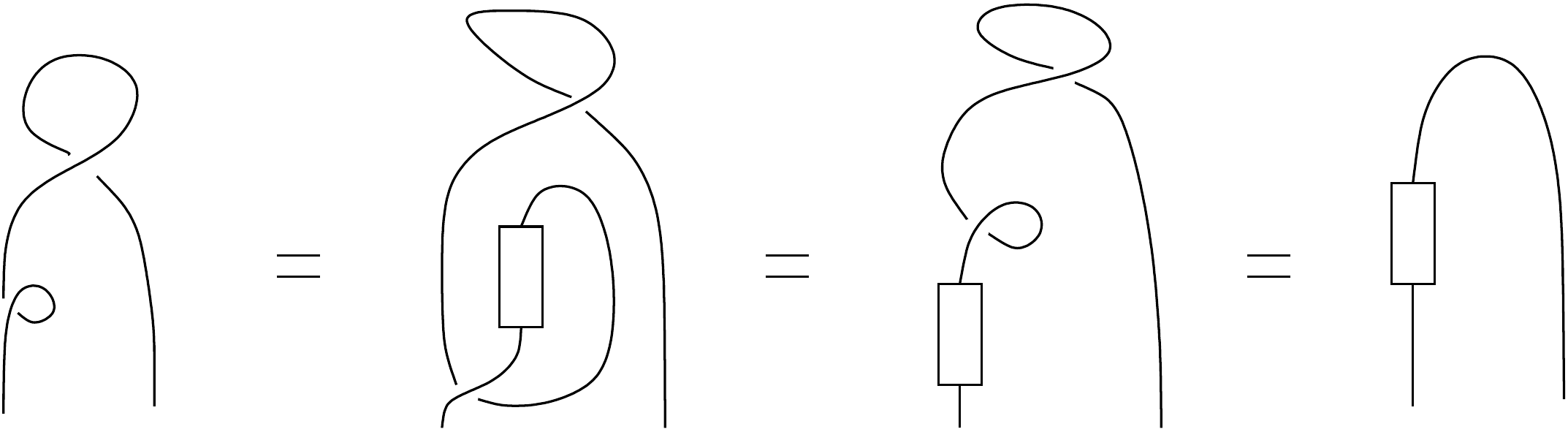} \caption{Dual pairing via the double dual isomorphism}\label{capddual}
\end{figure}
One may wonder why double and higher duals in the category have not
been discussed. What we are doing of course is implicitly
identifying $V^{**}$ with $V$. In any strict modular category it is
possible to do this explicitly.
\begin{exm}
Show that there exists a natural isomorphism
$\smash{V\xrightarrow{\delta_V}V^{**}}$ in any strict modular category such
that the twist is given by \fullref{twistfromdual} or
$$
\theta_V=(\text{id}_{V^*}\otimes\cap_V)\circ(\text{id}_V\otimes\delta_V\otimes\text{id}_{V^*})
\circ(\times_{V,V}\otimes\text{id}_{V^*})\circ(\text{id}_V\otimes\cup_V).
$$
The right way to do this exercise is to draw a colored ribbon tangle
representing a map from $V$ to $V^{**}$ built from elementary
pieces. The answer is \index{double dual isomorphism}
\index{$\delta$@$\delta_V$ double dual isomorphism}
\begin{multline*}
\delta_V=(\cap_V\otimes
\text{id}_{V^{**}})\circ(\times_{V,V^*}\otimes\text{id}_{V^{**}})\circ(\theta_V\otimes\text{id}_{V^{**}})
\circ(\text{id}_{V}\otimes\times_{V^{**},V^*})\\[-0.5ex]
\circ(\text{id}_{V}\otimes\theta^{-1}_{V^{**}}\otimes\text{id}_{V^*})
\circ(\text{id}_{V}\otimes\times_{V^*,V^{**}})\circ(\text{id}_{V}\otimes\cup_{V^*})\circ\theta_V.
\end{multline*}
\end{exm}
We call $\delta$ from the exercise the double dual isomorphism.
Using this isomorphism we can streamline the definitions of the dual
pairing and copairing (see \fullref{capddual})
$$
\begin{aligned}
\cup^*_V &=(\text{id}_{V^*}\otimes\delta_V^{-1})\circ\cup_{V^*}\\
\cap^*_V &=\cap_{V^*}\circ(\delta_V\otimes\text{id}_{V^*}).
\end{aligned}
$$
The best part about knowing $\delta_V$ explicitly is that by using
it one can evaluate a number of graphs without evaluating any
braidings in the process. We we will take full advantage of
this fact when discussing the modular categories coming from quantum
groups because as in many other nontrivial modular categories, it is the braiding that is
the hardest to compute. In particular, we get the following braiding-free
formula for the quantum trace (\fullref{qtrace}):
$$
\Tr_q(f):=\cap_{V^*}\circ(\delta_Vf\otimes\id_{V^*})\circ\cup_V.
$$
\begin{exm}\label{exddual}
Show that the double dual in $\calU(1)_{2m+1}$ from \fullref{u1categ} is the identity map $\C\to\C$ and the double dual in
$\text{REP}_G$ is the standard isomorphism between a vector space and its
double dual $\delta_V(v)(\varphi):=\varphi(v)$.
\end{exm}
\begin{figure}[ht!]
\centering
\labellist\small
\pinlabel {$d_{\lambda}$} [r] at 35 60
\pinlabel {$\lambda$} [r] at 90 35
\pinlabel {$p^+$} [r] at 195 60
\pinlabel {$p^-$} [r] at 354 60
\pinlabel {$\wtilde{s}_{\lambda\,\mu}$} [r] at 525 60
\pinlabel {$\lambda$} [r] at 564 100
\pinlabel {$\mu$} [l] at 587 87
\endlabellist
\includegraphics[width=4truein]{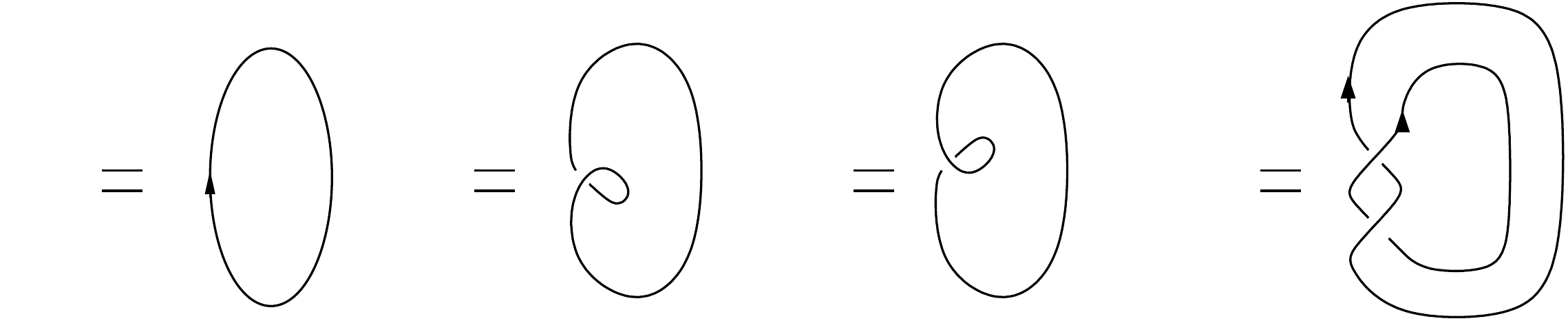}
\caption{Characteristic numbers of a  strict modular category}\label{mtcchar}
\end{figure}
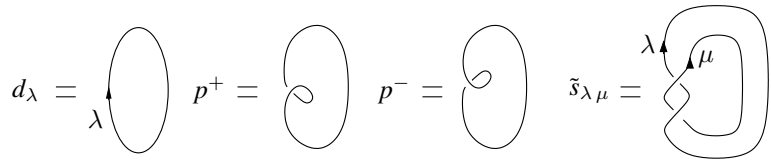

Our conventions thus far only allow evaluation of graphs that are
completely colored (labeled). But usual framed links do not have
representation labels attached to their components. One can make
sense of unlabeled graphs as well. The trick is to come up with
evaluation rules such that the resulting expressions do not change
at least under the Reidemeister moves (for link invariants) and even
better, do not change under the Kirby moves (for $3$--manifold
invariants). There are many different ways to color a particular
graph and the hope is that an appropriate linear combination of the
resulting evaluations will be invariant under the appropriate moves.
It turns out that there is an essentially unique way to evaluate
unlabeled graphs or unlabeled closed components in a graph that
guarantees invariance under the Kirby moves (see Turaev and Wenzl
\cite{TW}). We therefore introduce the following important convention.

\begin{quote}
In any diagram we sum over all ways of labeling unlabeled closed
components by simple objects $\lambda$ with multiplicity $d_\lambda$.
\end{quote}

Examples using our extended graphic notation to visualize the
definitions of the $\wtilde{s}$--matrix and the characteristic
numbers of a modular category are shown on \fullref{mtcchar}. The
next article explains in detail how to construct framed link
invariants from strict ribbon categories and how to construct
$3$--manifold invariants from strict modular categories.

\subsection{Invariants from modular categories}\label{modinv}

To see the correspondence between the Chern--Simons invariants and strict
modular categories one needs to describe how to build $3$--manifold
invariants from a strict modular category. We first define the
colored Jones polynomial of an oriented framed link $L$ with $c(L)$
components. \index{$c(L)$ number of components}Color the components
with objects $V_1,\ldots, V_{c(L)}$ from the strict ribbon
category. Notice that a diagram with only closed components
represents the element of $\text{End}({{\bbone}})$ obtained by taking
the composition of all of the basic morphisms corresponding to the
elementary framed tangles occurring between horizontal lines as in
\fullref{lhopf}. The axioms of a strict modular category ensure
that this element is invariant under all elementary isotopies of the
link.
\begin{defn}\label{Jones}
The colored Jones polynomial \index{colored Jones polynomial}
\index{$J_{V_1,\ldots, V_{c(L)}}(L)$ colored Jones polynomial} of
an oriented framed link is the element $J_{V_1,\ldots,
V_{c(L)}}(L)$ corresponding to the morphism represented by the
labeled link diagram.
\end{defn}
\begin{remark}
This invariant is not a polynomial. The name comes from the fact
that the invariant associated to the category of tilting modules is
closely related to the classical Jones polynomial \cite{Jn}
(which by the way, is not a polynomial either).
\end{remark}
\begin{defn}
Given a framed link (not oriented) $L$ define the invariant
\index{$F$ invariant}
$$
F(L)=\sum_{\lambda_1,\ldots,\lambda_{c(L)}\in
I}J_{\lambda_1,\ldots,\lambda_{c(L)}}(L)d_{\lambda_1},\ldots,d_{\lambda_{c(L)}}\,.
$$
Just as any isotopy of an ordinary link can be expressed as a
composition of elementary isotopies (Reidemeister moves I, II and
III) any isotopy of ribbon tangles with height function can be
decomposed into elementary isotopies. According to Turaev \cite{T7}, Freyd
and Yetter \cite{FY}, and Reshetikhin and Turaev \cite{RT1},
the elementary isotopies for ribbon tangles with height function are
Reidemeister moves II and III (see \fullref{reidemeister})
together with the moves displayed in \fullref{htreid}.
\begin{figure}[ht!]
\centering
\includegraphics[width=4truein]{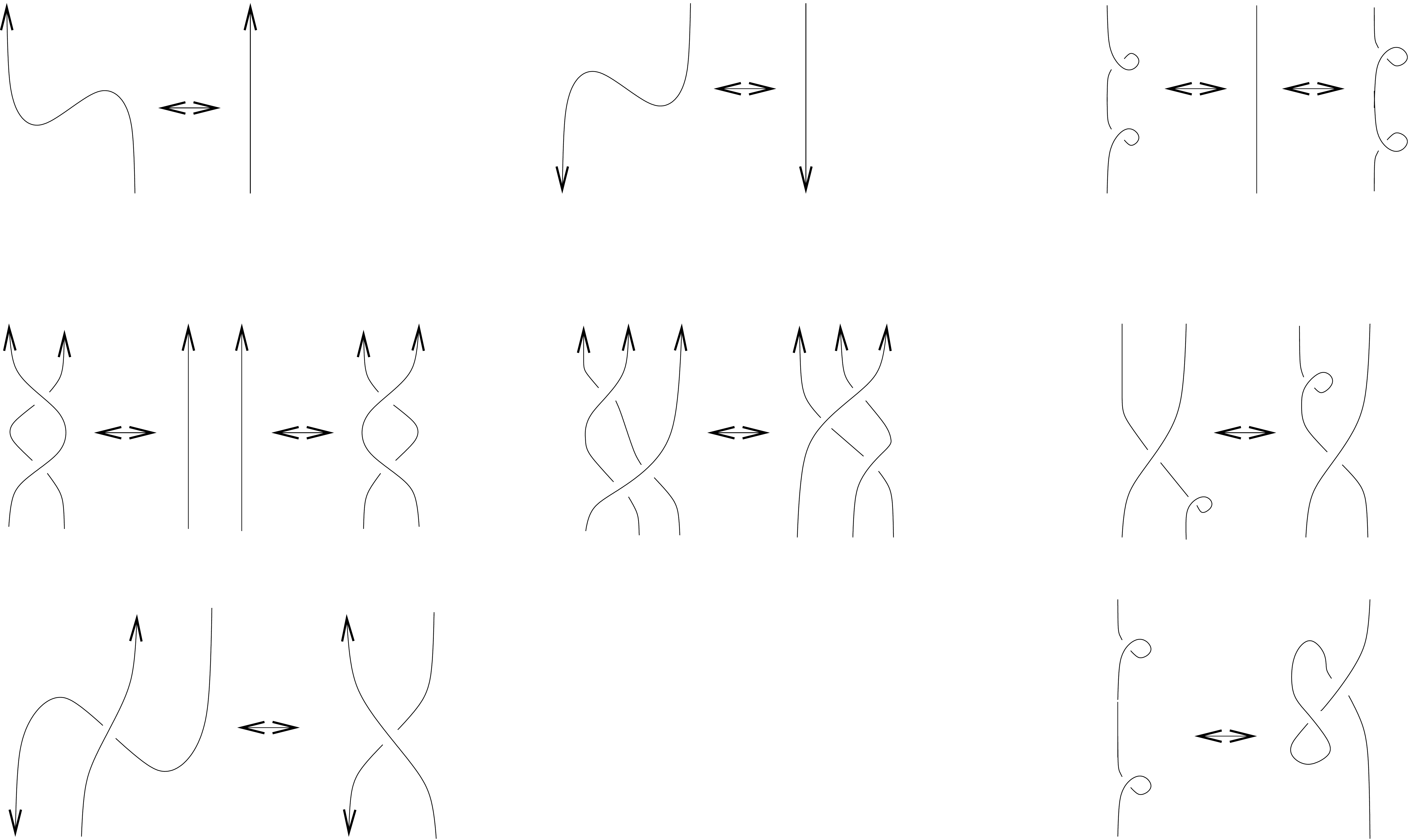}
\caption{Elementary ribbon isotopies}\label{htreid}
\end{figure}
\begin{exm}
Show that the invariant $F(L)$ is well defined when $I$ is any
collection of objects in a ribbon category, ie that it is
invariant under elementary isotopies.
\end{exm}

More generally given a framed link $L_M\cup L$ where $L$ is colored
and oriented define
$$
F(L_M,L)=\sum_{\lambda_1,\ldots,\lambda_{c(L_M)}\in
I}J_{\lambda_1,\ldots,\lambda_{c(L_M)},V_1,\ldots,V_{c(L)}}(L_M\cup
L)\,.
$$
\end{defn}

\begin{remark}
This is well defined independently of orientations on $L_M$ because
labeling a component with $\lambda$ is equivalent to changing the
orientation on the component and labeling it with $\lambda^*$. By
axiom 16 the sum is symmetric with respect to taking duals. This is
also consistent with the sum over all colorings of unlabeled
components convention.
\end{remark}

Let $M$ be the manifold obtained by surgery on a framed link $L_M$.
Even though $F(L_M)$ is invariant under the Reidemeister moves as
is, it is not invariant under the Kirby moves. Luckily, invariants
of non-isotopic links defining the same manifold are related in a
very simple manner and we can multiply $F(L_M)$ by normalizing
factors that cancel out this dependence.

Given a framed link, $L_M$, one defines a linking matrix
\index{$N$@$n_{ij}(L)$ linking matrix} \index{linking matrix}
$n_{ij}(L)$ with off-diagonal entries equal to the linking numbers
of the corresponding components and with diagonal entries equal to
the self-linking or writhe of the corresponding component
(see Prasolov and Sossinsky \cite{prasolov} and Rolfsen \cite{Rolfsen}).
\begin{example}
The linking matrix of the framed left Hopf link from \fullref{lhopf} is given on the left below and the linking matrix of the
framed link from \fullref{framed-link} is given on the right.
$$
\begin{pmatrix}0&-1\\-1&0\end{pmatrix}\qquad\qquad
\begin{pmatrix}2&\pm1\\\pm1&0\end{pmatrix}
$$
\end{example}

Let $\sigma(L_M)$ \index{$\sigma(L)$ signature of $L$} denote the
signature, that is, the number of positive eigenvalues minus the number
of negative eigenvalues of the linking matrix. With these notations
one can show that the quantity $\tau(M)$ defined below is a topological
invariant (see Bakalov and Kirillov \cite{BK} and Turaev \cite{T}). This
means that it does not depend on the link used to represent $M$ (so any
sequence of positive or negative Kirby moves or their inverses may be
applied to $L_M$ without changing the answer).
\begin{defn}\label{taudef}
The Reshetikhin--Turaev invariant \index{Reshetikhin--Turaev
invariant}\index{$\tau()$ Reshetikhin--Turaev invariant} of a
$3$--manifold is given by
$$
\tau(M):=(p^-)^{\sigma(L_M)}{\calD}^{-\sigma(L_M)-c(L_M)-1}F(L_M)\,.
$$
\end{defn}

It is also possible to define invariants of oriented colored framed
links in $3$--manifolds. A colored framed link in a $3$--manifold can
be represented by a link in $S^3$ with some of the components
labeled with objects from a strict modular category. Let $L_M$
denote the sublink consisting of unlabeled components. We assume
that $L_M$ is a surgery presentation of $M$. Each of the labeled
components represent a component of the framed link $L$ in $M$. The
labels correspond to the representations labeling Wilson loops in
the heuristic description.

The more general invariant is defined by
$$
\tau(M,L):=(p^-)^{\sigma(L_M)}{\mathcal
D}^{-\sigma(L_M)-c-1}F(L_M,L)\,.
$$
We define,
\begin{equation}\label{zwdef}
Z(M):=\tau(M,\emptyset), \qquad\text{and}\qquad
W_{R_1,\ldots,R_c}(L):=\tau(S^3,L)/\tau(S^3).
\end{equation}\index{$Z(M):=\tau(M,\emptyset)$}\index{$W_{R_1,\ldots,R_c}(L):=\tau(S^3,L)/\tau(S^3)$}
to be the mathematical interpretations of the physically motivated
invariants discussed earlier.

It may appear that we have a reasonably simple definition of
partition functions in \eqref{taudef}; however it depends on the
choice of a strict modular category and we do not have any other
than $\calU(1)_{2m+1}$ yet. The invariants corresponding to
$\calU(1)_{2m+1}$ are heuristically the same as the invariants
'defined' via the path integral by integrating holonomies over a
space of U$(1)$ connections. Since U$(1)$ is Abelian, the cubic
terms in the Chern--Simons invariant vanish so the resulting theory
is what physicists call a Gaussian theory. This is the case where
the path integral would be easiest to formalize mathematically but
it leads to fairly weak invariants. In the next subsection we define
quantum groups because their representations lead to more
complicated strict modular categories that in turn lead to more
interesting quantum invariants.

We now compute the invariant of $S^3$ in three different ways. Using
the empty link we see that $\tau(S^3)=\calD^{-1}$. The following examples
compute the same invariant from different link presentations.
Note that there is some algebra involved in establishing that the obtained
values are the same.
\begin{example}
Let $L$ be the link (twisted unknot) that defines $p^-$ in \fullref{mtcchar}. We compute $F(L)$ and $\tau(M_L)$ in the category
$\calU(1)_{2m+1}$. Coloring the single component by $p$ and
orienting the component clockwise we have
$$
J_p(L)=\cap_{p^*}\circ(\theta^{-1}_p\otimes\text{id}_{p^*})\circ\cup_p=\theta^{-1}_p=e^{-2\pi
ip^2/(2m+1)},
$$
where $p$ takes values $0,1,\ldots,2m$. Since $d_p=1$ for all $p$,
\def\Sum{\sum\nolimits}
$$
F(L)=\Sum^{2m}_{p=0}d_pJ_p(L)=\Sum^{2m}_{p=0}e^{-2\pi ip^2/(2m+1)}.
$$
Since $L$ only has one component the linking `matrix' is just the
self-linking number which is $-1$ because of the twist. The
signature is also $-1$. The quantum diameter of $\calU(1)_{2m+1}$ is
$\mathcal{D}=(2m+1)^{\unfrac12}$ (see \fullref{qdiam}). By the very
choice of $L$ we have $p^-=F(L)$. Thus
$$
\tau(M_L)=(p^-)^{-1}{\calD}^{1-1-1}F(L)={\calD}^{-1}=(2m+1)^{-\unfrac12}.
$$
\end{example}

\begin{example}
A slightly more difficult example is the left Hopf link from \fullref{lhopf}. In fact, we already computed
$F(L)=\sum^{2m}_{p,q=0}e^{-4\pi ipq/(2m+1)}=2m+1$ back in \fullref{lhopfeg}. Self-linking numbers in this case are both $0$ (no
twists) and the linking numbers between the components are $-1$
(orientation!). Hence the linking matrix is
$$
\text{lk}(L)=\begin{pmatrix}0&-1\\-1&0\end{pmatrix}
$$
and $\sigma(L)=1-1=0$. Assembling the results we get
$$\tau(M_L)=(p^-)^0\mathcal{D}^{-3}F(L)=(2m+1)^{-\unfrac32}
\sum^{2m}_{p,q=0}e^{-4\pi ipq/(2m+1)}=(2m+1)^{-\unfrac12}.$$
\end{example}
\begin{exm}
Use the Reidemeister and Kirby moves to show that all three links
employed above represent the same $3$--manifold, namely the $3$--sphere.
\end{exm}
To fully describe how the notion of a strict modular category
relates to the Witten's Chern--Simons invariants one should describe
the associated topological quantum field theory. The general outline
is as follows. Diagrams with only closed components may be cut into
parts with marked surface boundaries. A Hilbert space will then be
associated to each such marked surface and a bounded linear map will
be associated to a cobordism. The Hilbert space associated to a
genus $g$ surface marked with colors (simple elements of the
category) $V_1,\ldots, V_n$ is ${\mathcal
H}=\Hom({{\bbone}},(\bigoplus_{j\in I}\lambda_j\otimes
\lambda_j^*)^{\otimes g}\otimes V_1\otimes\cdots\otimes V_n)$. There
are many technical details that must be addressed in order to do
this properly. A very careful explanation is given by Turaev
\cite{T}.
\begin{exm}
Compute the $\mathcal{U}(1)_{2m+1}$ invariant of the $3$--manifold
represented by the framed link from \fullref{framed-link}.
\end{exm}
\begin{exm}
Describe how to compute the $W$ invariant of any framed link for
this strict modular category. (There is a fairly simple formula in
terms of the linking matrix and colors.)
\end{exm}

\section{Quantum groups and their representations}\label{RTinv}

The ingredients in a strict modular category look like the
representations of some algebraic object. This is indeed one of the
best methods to construct strict modular categories. In this
subsection we consider the representations of deformations of Lie
algebras called quantum groups. Deformations are constructed to
arrive in a non-commutative setting. Without such a deformation the
morphism that would be attached to a right crossing would be the
same as the one attached to a left crossing. The invariants
constructed from the resulting representations are the original
Reshetikhin--Turaev invariants.

\subsection{Quantum groups at roots of unity}

The axioms of  a strict modular category are very complicated, so
without some additional motivation it would be difficult to
construct an interesting modular category. Luckily the study of
symmetry in quantum mechanics led to very similar structures and it
was possible to construct interesting modular categories from
quantum groups.

The category that produces the usual quantum or Reshetikhin--Turaev
invariants is a category of representations of certain quantum
groups. The latter are algebraic objects that replace classical
groups in description of symmetries on quantum (or noncommutative)
spaces. The term `quantum group' is used rather loosely and is
usually reserved for deformations of algebras associated to
classical groups. In particular, they are not groups; instead they
generalize classical group algebras and enveloping algebras rather
than groups.

We begin this section with a brief description of the ideas that led
to the discovery of quantum groups. As one would guess from the name
the original motivation for the
definition of a quantum group comes from quantum mechanics. However
there is also a strong analogy between quantum groups and algebraic
groups and much of the theory of quantum groups was first developed
for algebraic groups. To quantize a mechanical system via canonical
quantization one tries to find an embedding of a deformation of the
algebra of observables into the linear operators on a Hilbert space
(see Simms and Woodhouse \cite{wod}). The failure of operators
corresponding to various observables to commute can be interpreted
as the uncertainty principle. Often the Hilbert space can be taken
to be a space of functions on the configuration space.

\begin{example}
For a free particle moving in $1$--dimensional space the observables
are functions of position $x$ and momentum $p$. A possible Hilbert
space to associate to this system is the space of $L^2$ functions in
the variable $x$. One then associates multiplication by $x$, denoted
by $L_x$, to the observable $x$ and the operator
$L_p=-i\hbar\frac{\partial}{\partial x}$ to the observable $p$.
Notice that one has $[L_x,L_p]=i\hbar$. Thus $L_x$ and $L_p$ do not
commute. However in the classical limit $\hbar\to 0$ one recovers
the classical algebra of observables. In most everyday situations
this is a reasonable approximation because $\hbar=1.055\, \times\,
10^{-34}$ joule-sec.
\end{example}

For the case of maximal symmetry the configuration space is a group.
When the underlying configuration space is a group the space of
functions on the group can be given the structure of a Hopf algebra.
Thus Hopf algebras appear naturally in quantum mechanics as algebras
of functions on groups.

Recall that a group can be described as a tuple $(G,\mu\co G\times G\to
G, e\co {\bbone}\to G, n\co G\to G)$ such that the following diagrams commute.
$$\bfig
  \barrsquare<600,500>[G{\times}G{\times}G`G{\times}G`G{\times}G`G;
    \text{id},\mu`\mu,\text{id}`\mu`\mu]
  \efig\qquad\bfig
  \barrsquare<600,500>[G{\times}\bbone`G{\times}G`G`G;
    \text{id},e`\cong`\mu`\text{id}]
  \efig\qquad\bfig
  \barrsquare/->`<-`->`->/<600,500>[G{\times}G`G{\times}G`G`G;
   \text{id},n`\delta`\mu`e\circ p]
  \efig$$
  Here $\delta\co G\to G\times G$ is the diagonal map.
If $A$ is an algebra of functions on a finite group $G$, the
multiplication $\mu$ identity $e$ and inverse $n$ on $G$ will induce
a comultiplication $\Delta\co A\to A\otimes A$, counit
$\varepsilon\co A\to\C$ and antipode $\gamma\co A\to A$ on $A$ such
that the following diagrams commute.
$$\bfig
  \barrsquare/<-`<-`<-`<-/<700,500>[
    A{\otimes}A{\otimes}A`A{\otimes}A`A{\otimes}A`A;
    \text{id}{\otimes}\Delta`\Delta{\otimes}\text{id}`\Delta`\Delta]
  \efig\quad\bfig
  \barrsquare/<-`<-`<-`<-/<700,500>[
    A{\otimes}\C`A{\otimes}A`A`A;
    \text{id}{\otimes}\varepsilon`\cong`\Delta`\text{id}]
  \efig\quad\bfig
  \barrsquare/<-`->`<-`<-/<700,500>[
    A{\otimes}A`A{\otimes}A`A`A;
    \text{id}{\otimes}\gamma`\mu`\Delta`\iota{\circ}\varepsilon]
  \efig
$$
The comultiplication \index{$\Delta$ comultiplication}
\index{comultiplication} $\Delta$ and counit \index{$\varepsilon$
counit} $\varepsilon$ \index{counit} are algebra homomorphisms and
the antipode \index{antipode}\index{$\gamma$ antipode} is an
anti-homomorphism ($\gamma(ab)=\gamma(b)\gamma(a)$.) This is
essentially the \index{Hopf algebra}  definition of a Hopf algebra.
\begin{defn}\label{Hopfalg}
A Hopf algebra $\calA$ over a field $\F$ is an associative algebra
with additional algebra homomorphisms (counit, coproduct)
$\ve\co \calA\to\F$, $\Delta\co \calA\otimes\calA$ and an algebra
antihomomorphism (antipode) $\gamma\co \calA\to\calA$ satisfying axioms
that dualize the usual axioms for the unit, product and inverse in a
group (see the examples below, Chari and Pressley \cite{CP} or Majid
\cite{Maj} for the complete list).
\end{defn}

For quantization one needs to deform the algebra of functions into a
non-commutative algebra. Now consider the case of the group
$SL_2\C$. The natural action of this group on $\C^2$ will help
motivate the correct deformation. The algebra of functions on $\C^2$
is just the ring of polynomials $\C[x,y]$ and the algebra of
functions on $SL_2\C$ is $\C[a,b,c,d]/(ad-bc-1)$. Recall that a map
between spaces induces a map on the associated function algebras
going in the opposite direction. Now the natural action of $SL_2\C$
on $\C^2$ by matrix multiplication induces the corresponding map of
function algebras
$$
\C[x,y]\to(\C[a,b,c,d]/(ad-bc-1))\otimes\C[x,y].
$$
To arrive at a non-commutative deformation of the function algebra
of $SL_2\C$ (that is, $\C[a,b,c,d]/(ad{-}bc{-}1)$), start with the
non-commutative complex plane. The function algebra of the
non-commutative plane is $\C\{x,y\}\llbracket h\rrbracket
/(xy-e^{-h}yx)$. Here $\C\{x,y\}$ is the free algebra on two
generators over $\C$ and $R\llbracket h\rrbracket$ refers to formal
power series in the variable $h$ with coefficients in $R$. In order
to keep track of all of the important algebraic structure one must
consider the Hopf algebra structure on the resulting algebra.

There is essentially a unique Hopf algebra that respects the map of
function algebras induced by the action of $SL_2\C$ on $\C^2$ with
the algebra of $\C^2$ replaced by the non-commutative version and
the determinant replaced by $ad-e^{-h}bc$. This algebra is denoted
by $SL^q_2\C$. Such algebras or their duals are called quantum
groups (sic!). The quantum group $U_q({\mathfrak sl}_2{\mathbb C})$
with $q=e^h$ is `dual' to the deformed function algebra $SL^q_2\C$.
See Chari and Pressley \cite[Chapter 7]{CP} for more details.

Because the axioms for operations and co-operations in a Hopf
algebra are symmetric one can define a dual Hopf algebra $\calA^*$
(in a couple of ways) switching them, that is, the product in the
dual comes from the coproduct in the original, etc. We will be
mostly interested in the duals of the deformed group algebras such
as $SL^q_2\C$. For simply connected Lie groups these duals or
distribution algebras can be described as deformations of the
universal enveloping algebras $U(\g)$ of the corresponding Lie
algebras $\g$. Our main example comes from deforming $U(\Sl_N\C)$.

Many Lie algebras are matrix algebras with bracket given by
$[X,Y]=XY-YX$. Any Lie algebra ${\mathfrak g}$ can be embedded in a
unital associative algebra so that the Lie bracket takes this form.
Recall that the universal enveloping algebra \index{$U(\g)$
universal enveloping algebra}\index{universal enveloping algebra}
$U(\g)$ is the quotient of the tensor algebra
$\oplus_{n\geq0}\g^{\otimes n}$ by the ideal generated by $x\otimes
y-y\otimes x-[x,y]$ with $x,y\in\g$. In a form more germane to
quantum generalizations $U(\g)$ can be described by generators and
relations. Namely, if $x_1,\ldots,x_n$ generate $\g$ with relations
given in terms of brackets then $U(\g)$ has a presentation with the
same set of generators with relations obtained by replacing all
brackets with the corresponding commutators. For example $[x,y]$ is replaced by
$xy-yx$. We will usually write brackets even when working in $U(\g)$
interpreting them as commutators.

Before discussing the deformations recall the presentation of the
enveloping algebra.
\begin{exm}
Compute $[e,f]$, $[L,e]$ and $[L,f]$ for the following matrices.
$$e=\left(\begin{matrix}0&1\\0&0\end{matrix}\right), \qquad
f=\left(\begin{matrix}0&0\\1&0\end{matrix}\right), \qquad
L=\left(\begin{matrix}1&0\\0&-1\end{matrix}\right).$$
\end{exm}
\begin{exm}\label{quadraticserre}
Express $[X,[X,Y]]$ as a product of $X$ and $Y$ factors. Let
$E_{i,j}$ be the square matrix with a $1$ in the $i,j$ entry and
zeros elsewhere. Compute $[E_{i,i+1},E_{j,j+1}]$,
$[E_{i,i+1},E_{j+1,j}]$ and $[E_{i,i+1},E_{j,j}-E_{j+1,j+1}]$.
\end{exm}

\begin{example}\label{UslN}
Recall from Fulton and Harris \cite{fulton-harris}, Humphreys
\cite{Hum} or \fullref{app:d} that $\Sl_N\C$ is generated by
$e_i,f_i,\alpha_i^\vee$ and $i=1,\ldots,N-1$. Here $\alpha_i^\vee$
represent the simple coroots of $\Sl_N\C$ and $e_i,f_i$ are the
corresponding positive and negative root vectors. In other words
$e_i=E_{i,i+1}$, $f_i=E_{i+1,i}$ and
$\alpha_i^\vee=E_{i,i}-E_{i+1,i+1}$. The relations in $\Sl_N\C$ and
$U(\Sl_N\C)$ are:
\begin{equation}\label{slNrel}
\begin{aligned}
{[\alpha^\vee_i,\alpha^\vee_j]}&=0 &
  [e_i,f_j]&=\delta_{ij}\alpha^\vee_j,\\
[\alpha^\vee_i,e_j]&=a_{ij}e_j &
  [e_i,e_j]&=[f_i,f_j]=0, & \quad |i-j|&\ne 1\\
[\alpha^\vee_i,f_j]&=-a_{ij}f_j &\quad
  [e_i,[e_i,e_j]] &=[f_i,[f_i,f_j]] =0, & j&=i\pm1,
\end{aligned}
\end{equation}
where
$$a_{ij}=\text{Tr}(\alpha^\vee_i\alpha^{\vee\dagger}_j)=\begin{cases} 0 &|i-j|>1\\ -1 & |i-j|=1\\2 & i=j\end{cases}$$
is the Cartan matrix \index{Cartan matrix}\index{$a_{ij}$ Cartan
matrix} of $\Sl_N\C$. The last pair of relations are known as the
Serre relations. The Serre relations look more familiar if written
in the associative form
\begin{equation}\label{Serre}
e_i^2e_j-2e_ie_je_i+e_je_i^2=f_i^2f_j-2f_if_jf_i+f_jf_i^2=0,\ j=i\pm1.
\end{equation}
$U(\Sl_N\C)$ also has a  Hopf algebra structure given by
$$
\ve(x)=0,\qquad\gamma(x)=-x,\qquad\Delta(x)=1\otimes x+x\otimes 1.
$$
For $\Sl_2\C$ there is only one generator of each type
$e,f,\alpha^\vee$, the Cartan matrix is $1\times 1$, that is, the number
$2$ and the Serre relations trivialize so \eqref{slNrel} reduces to
\begin{equation*}
[\alpha^\vee,e]=2e,\qquad[\alpha^\vee,f]=-2f\qquad[e,f]=\alpha^\vee.
\end{equation*}
\end{example}

Following the motivation above we now explain the quantum deformation $U_q(\Sl_2\C)$.
According to V Drinfeld \cite{Drin} the first two relations in \eqref{slNrel}
should stay intact whereas the last one should become
$$[e,f]=\frac{\sinh(\hbar\alpha^\vee)}{\sinh(\hbar)},$$
where $\hbar$ is a formal parameter (the Planck constant). To make
sense of this we would have to consider formal power series in
$\hbar$ with coefficients in $U(\Sl_2\C)$ which is not very
convenient. Fortunately, there is a bypass due to M Jimbo who
suggested setting $q=e^{\hbar}$ and introduced two new generators
$q^{\pm\alpha^\vee}$ so that the last relation becomes
$$[e,f]=\frac{q^{\alpha^\vee}-q^{-\alpha^\vee}}{q-q^{-1}}.$$
Now we only have to extend the field of coefficients from $\C$ to $\C(q)$
(rational functions in $q$). Of course we now need to eliminate
$\alpha^\vee$ form the first two relations above which leads to
$$q^{\alpha^\vee}eq^{-\alpha^\vee}=q^2e,\qquad
  q^{\alpha^\vee}fq^{-\alpha^\vee}=q^{-2}f.$$
Note that $q^{\pm\alpha^\vee}$ are indeed new generators and by no
means a variable $q$ `taken to the power' $\pm\alpha^\vee$ (some
authors denote them $K^{\pm1}$ to prevent confusion; see Chari and
Pressley \cite{CP} and Majid \cite{Maj}). It is interesting that the
representations of $U(\Sl_2\C)$ and $U(\Sl_2\C)[[\hbar]]$ are in
bijective correspondence \cite{CP}.  Moreover, any deformation of
$U(\Sl_2\C)$ and more generally $U(\g)$ as an associative algebra is
trivial, that is it produces an isomorphic algebra \cite{CP}. It is
in the Hopf algebra structure that the difference between $U$ and
$U_q$ becomes essential. For $U_q(\Sl_2\C)$ the deformed
co-operations are given on the generators by
\begin{align*}
\ve(q^{\pm\alpha^\vee})&=1, &
  \gamma(q^{\pm\alpha^\vee})&=q^{\mp\alpha^\vee}, &
  \Delta(q^{\pm\alpha^\vee})&=q^{\pm\alpha^\vee}\otimes q^{\pm\alpha^\vee},\\
\ve(e)&=0, &
  \gamma(e)&=-eq^{-\alpha^\vee}, &
  \Delta(e)&=e\otimes q^{\alpha^\vee}+1\otimes e,\\
\ve(f)&=0, &
  \gamma(f)&=-q^{\alpha^\vee}f, &
  \Delta(f)&=f\otimes 1+q^{-\alpha^\vee}\otimes f.
\end{align*}
As Hopf algebras $U(\Sl_2\C)$ and $U_q(\Sl_2\C)$ are not isomorphic.

One would expect that $U_q(\mathfrak{sl}_2{\mathbb C})$ is a
deformation of $U(\mathfrak{sl}_2{\mathbb C})$ but this is not
quite true. What is true (see Kassel \cite{kassel}) is
$$U(\mathfrak{sl}_2{\mathbb C})=U_{q=1}(\mathfrak{sl}_2{\mathbb
  C})/(q^{\alpha^\vee}-1).$$
We now give the general definition for $U_q(\Sl_N\C)$ following
Chari and Pressley \cite{CP}.
\begin{defn}\label{UqslN}
The (rational form of the) Drinfeld--Jimbo quantum group
$U_q(\Sl_N\C)$ \index{Drinfeld--Jimbo quantum group}
\index{$U_q(\Sl_N\C)$
 quantum group} is the Hopf algebra generated as an
associative algebra over $\C(q)$ by the generators
\index{$q^{\pm{\alpha^\vee}_i},e_i,f_i$ generators of
$U_q(\Sl_N\C)$} $q^{\pm{\alpha^\vee}_i},e_i,f_i$ with
$i=1,\ldots,N-1$ satisfying the relations
\begin{equation}\label{Uqalg}
\begin{array}{lll}
q^{{\alpha^\vee}_i}q^{-{\alpha^\vee}_i}{=}q^{-{\alpha^\vee}_i}q^{{\alpha^\vee}_i}{=}1,
\qua &
  [e_i,f_j]{=}\delta_{ij}\frac{q^{{\alpha^\vee}_i}-q^{-{\alpha^\vee}_i}}{q-q^{-1}}, & \\
q^{{\alpha^\vee}_i}e_jq^{-{\alpha^\vee}_i}{=}q^{a_{ij}}e_j, &
  [e_i,e_j]{=}[f_i,f_j]{=}0, & |i-j|\ne 1,\\
q^{{\alpha^\vee}_i}f_jq^{-{\alpha^\vee}_i}{=}q^{-a_{ij}}f_j, &
  e_i^2e_j{-}(q{+}q^{-1})e_ie_je_i{+}e_je_i^2{=}0, & j{=}i\pm1,\\
& f_i^2f_j{-}(q{+}q^{-1})f_if_jf_i{+}f_jf_i^2{=}0, & j{=}i\pm1.
\end{array}
\end{equation}
Here $a_{ij}=\text{Tr}(\alpha^\vee_i\alpha^{\vee\dagger}_j)$ is the
Cartan matrix  of $\Sl_N\C$. On the generators the counit, the
antipode and the coproduct are given by
\begin{equation}\label{Uqcoalg}
\begin{aligned}
\ve(q^{\pm{\alpha^\vee}_i})&=1, &
  \gamma(q^{\pm{\alpha^\vee}_i})&=q^{\mp{\alpha^\vee}_i}, &
  \Delta(q^{\pm{\alpha^\vee}_i})&=q^{\pm{\alpha^\vee}_i}\otimes q^{\pm{\alpha^\vee}_i},\\
\ve(e_i)&=0, &
  \gamma(e_i)&=-e_iq^{-{\alpha^\vee}_i}, &
  \Delta(e_i)&=e_i\otimes q^{{\alpha^\vee}_i}+1\otimes e_i,\\
\ve(f_i)&=0, &
  \gamma(f_i)&=-q^{{\alpha^\vee}_i}f_i, &
  \Delta(f_i)&=f_i\otimes 1+q^{-{\alpha^\vee}_i}\otimes f_i.\\
\end{aligned}
\end{equation}
\end{defn}

\index{$\Delta$ comultiplication}
\index{comultiplication}\index{$\varepsilon$ counit} \index{counit}
\index{antipode}\index{$\gamma$ antipode}
\begin{remark}
Notice that for $\mathfrak{sl}_N\C$ the coroots and the roots can be
identified. The reason for using the notation that we used here is
that it generalizes automatically to a quantum group constructed
from any semisimple Lie algebra.
\end{remark}

The comultiplication is used to construct the tensor product in the
strict modular category constructed from representations of
quantum groups, the counit leads to the unit, the antipode leads to
the duality and the square of the antipode leads to the twist.
Constructing the braiding is a bit tricky. We address these issues
later in this section.

One other tricky point is making sure that there are only a finite
number of simple objects. This is accomplished by specializing to a
root of unity; however, this is not as easy as one might guess. It
is tempting to consider $U_q(\Sl_N\C)$ over $\C$ rather than $\C(q)$
by specializing $q$ to a particular complex number $z$. When $z$ is
not a root of unity one obtains a Hopf algebra with generators and
relations given by \eqref{Uqalg} and \eqref{Uqcoalg}, and $q$ replaced by
$z$ (not in $q^{\pm\alpha_i}$ since those are names of generators).
As associative algebras $U_z(\Sl_N\C)$ and  $U(\Sl_N\C)$ are
isomorphic and thus have identical representation theories. As nice
as this might be it means an infinite number of simple objects
(irreducible representations) and no hope of modularity. To
understand why the case $z=\eps,\ \eps^l=1$ with $l$ an integer is
exceptional we need to introduce the notions of $q$--integers,
$q$--factorials and $q$--binomials.
\begin{defn}\label{qnum}
For any $n\in\Z$ define the $q$--integers  \index{quantum integers}
(quantum integers) as
$$[n]_q:=\frac{q^{n}-q^{-n}}{q-q^{-1}}\in\C(q)$$
\index{$[n]_q:=(q^{n}-q^{-n})/(q-q^{-1})\in\C(q)$} and for $n\geq0$
define the $q$--factorials \index{quantum factorial} as
\index{$[n]_q\pling:=[n]_q[n-1]_q\ldots[2]_q[1]_q$}
$$[n]_q!:=[n]_q[n-1]_q\ldots[2]_q[1]_q.$$
For any pair of integers \index{$\left[n\atop
m\right]_q:=\frac{[n]_q\pling}{[n-m]_q\pling[m]_q\pling}$}\index{quantum
binomial coefficient} $0\leq m\leq n$ define the $q$--binomial
coefficient as
\begin{equation}\label{qbin}
\left[n\atop m\right]_q:=\frac{[n]_q!}{[n-m]_q![m]_q!}.
\end{equation}
\end{defn}
\begin{exm}
Verify that $[n]_q=q^m[n-m]_q+q^{-(n-m)}[m]_q$ and use it to
derive the Pascal recurrence relation for $q$--binomials:
$$
\left[n\atop m\right]_q=q^m\left[{n-1}\atop m\right]_q+q^{-(n-m)}\left[{n-1}\atop{m-1}\right]\,.
$$
\end{exm}
Obviously, in the limit $q\to1$ the $q$--numbers turn into the
ordinary integers, factorials and binomials. Many formulas in
representation theory involve multiplication by $n!$.  Under quantum
deformation these become multiplications by $[n]_q!$. However, if
$\eps^{l}=1$ and we specialize to $q=\eps$ then $[l^\prime]_\eps!=0$
for $\eps^{l^\prime}=\pm 1$ and some irreducible representations
will become reducible. For $\eps=\pm 1$ the defining relations in
\eqref{Uqalg} do not even make sense as written due to division by
$\eps-\eps^{-1}=0$. But even for $l>2$ simply setting $q=e^{2\pi
i/l}$ will not produce the `right' version of the quantum group. The
right version was introduced by G Lusztig \cite{Lus}.

Lusztig notes that
$[n]_q=\sum_{k=0}^{n-1}q^{n-1-2k}\in\Z[q,q^{-1}]$, that is, it is a Laurent
polynomial in $q$, and therefore so is $[n]_q!$. Less obviously, the
$q$--binomial $\left[n\atop m\right]_q$ is also a Laurent polynomial.
\begin{exm}
Prove the last claim.\\
Hint: Use the Pascal recurrence and proceed by induction on $n$ and
$m$.
\end{exm}
The Laurent polynomials play the same role in $\C(q)$ as the
integers play in $\C$. More importantly, for any
$\pi\in\Z[q,q^{-1}]$ the value $\pi(z)$ is well-defined for any
$z\ne0$, in particular, $q$ can be specialized even to roots of
unity in Laurent polynomials. This suggests defining an integral
form  of $U_q(\g)$ before specializing to roots of unity. For the
classical enveloping algebras such a form is known as the Kostant
$\Z$--form (see Humphreys \cite{Hum}) and it uses the divided powers
$\frac{e^n_i}{n!},\frac{f^n_i}{n!}$ as generators. This motivates
the following definition (see Lusztig \cite{Lus}):
\begin{defn}\label{divpw}
The divided powers in $U_q(\Sl_N\C)$ are \index{divided powers}
\index{$e_i^{(n)}:=\frac{e^n_i}{[n]_q\pling}$}
\index{$f_i^{(n)}:=\frac{f^n_i}{[n]_q\pling}$} defined by
$$e_i^{(n)}:=\frac{e^n_i}{[n]_q!},\qquad
f_i^{(n)}:=\frac{f^n_i}{[n]_q!}.$$
\end{defn}
The trick now is to rewrite the relations \eqref{Uqalg} in terms of
the divided powers and make sure that their coefficients are Laurent
polynomials. This is indeed the case and we give some of the
relations below (the full list occupies an entire page in Chari and
Pressley \cite{CP}).
\begin{equation}\label{Uqdiv}
\begin{gathered}
\begin{aligned}
q^{\alpha_i}e^{(n)}_jq^{-\alpha_i}&=q^{na_{ij}}e^{(n)}_j, &
  e^{(m)}_ie^{(n-m)}_i&=\left[n\atop m\right]_qe^{(n)}_i,\\
q^{\alpha_i}f^{(n)}_jq^{-\alpha_i}&=q^{-na_{ij}}f^{(n)}_j, &
  f^{(m)}_if^{(n-m)}_i&=\left[n\atop m\right]_qf^{(n)}_i, \\
e^{(n)}_if^{(m)}_j&=f^{(m)}_je^{(n)}_i,\qquad i\ne j,
\end{aligned}\\
e^{(m)}_ie^{(n-m)}_i=\Sum_{j=1}^{\min[n,n-m]}f^{(n-m-j)}_i\left[q^{\alpha_j};\
2j-n\atop j\right]_qe^{(m-j)}_i\,.
\end{gathered}
\end{equation}
In the last formula we used a new notation
\index{$\left[q^{\alpha_i};\ c\atop j\right]_q$ integral form
elements}
\begin{equation}\label{cengen}
\left[q^{\alpha_i};\ c\atop j\right]_q:=\Prod_{k=1}^j\frac{q^{c+1-k}q^{\alpha_i}-q^{-(c+1-k)}q^{-\alpha_i}}{q^k-q^{-k}}
\end{equation}
with $c\in\Z$ and $j\in\Z_{\geq0}$. These new elements are not
generators, they can be expressed via $e_i^{(n)}$, $f_i^{(n)}$ with
$\Z[q,q^{-1}]$ coefficients (see Lusztig \cite{Lus}). We mention them because they
will play an important role in the representation theory later. Once
we have the new relations we can forget about the origin of the
divided powers and treat them as formal symbols that satisfy the
relations \eqref{Uqdiv}. Indeed, we have to do this since \fullref{divpw} makes no sense for $q=\eps$ a root of unity.
\begin{defn}[Quantum groups at roots of unity]\label{Uqres}
The \index{restricted integral form} restricted integral form
\index{$U_{\Z[q,q^{-1}]}^{\text{res}}(\Sl_N\C)$ restricted integral
form} $\smash{U_{\Z[q,q^{-1}]}^{\text{res}}}(\Sl_N\C)$ of the
quantum group $U_q(\Sl_N\C)$ is the Hopf algebra generated over
$\Z[q,q^{-1}]$ by $q^{\pm\alpha^\vee_i},e_i^{(n)},f_i^{(n)}$.  The
corresponding quantum group at a root of unity
$U_\eps^{\text{res}}(\Sl_N\C)$ is obtained by specializing $q$ to
$\eps$ and changing the coefficients to $\C$. Formally,
$$
U_\eps^{\text{res}}(\Sl_N\C):=U_{\Z[q,q^{-1}]}^{\text{res}}(\Sl_N\C)\otimes_{\Z[q,q^{-1}]}\C\,,
$$
where $\Z[q,q^{-1}]$ acts on $\C$ in the obvious way with
$\pi\mapsto\pi(\eps)$. \index{$U_\eps^{\text{res}}(\Sl_N\C)$ quantum group
at $\eps$}
\end{defn}
\begin{remark}
The algebra $\smash{U_{\Z[q,q^{-1}]}^{\text{res}}(\Sl_N\C)}$ has a
presentation with generators $q^{\pm\alpha_i}$, $e_i^{(n)}$,
$f_i^{(n)}$ and relations including \eqref{Uqalg} and co-operations
induced from \eqref{Uqcoalg}.
\end{remark}

At this point the whole digression on divided powers may seem
superfluous: why not simply take the subalgebra of $U_q(\Sl_N\C)$
generated over $\Z[q,q^{-1}]$ by $q^{\pm\alpha^\vee_i},e_i,f_i$ and
then specialize to $\eps$? Intuitively, the difference is due to the
following. In $U_q(\Sl_N\C)$ we have the equality
$e_i^n=\smash{[n]_q!e_i^{(n)}}$. Since this equality only contains a
Laurent polynomial it continues to hold in
$\smash{U_{\Z[q,q^{-1}]}^{\text{res}}(\Sl_N\C)}$ and therefore in
${U_\eps^{\text{res}}(\Sl_N\C)}$. However, if $\eps^l=1$ then
$[n]_\eps!=0$ for $n\geq2l$ and the higher powers of $e_i,f_i$
vanish while the divided powers survive! Otherwise there is no way
to define divided powers and the quantum group simply looses part of
the structure. On a bright side, for $0\leq n<\unfrac{l}{2}$ we have
$[n]_\eps!\ne0$ and the quantity
$e_i^{(n)}=\unfrac{e^n_i}{[n]_\eps!}$ is well-defined even in
$\smash{U_\eps^{\text{res}}}(\Sl_N\C)$. This simplifies computations
with the divided powers in this range and allows one to use simpler
relations \eqref{Uqalg} instead of \eqref{Uqdiv}.

\subsection{Representations of $U_\eps^{\text{res}}(\Sl_N\C)$ and tilting modules}\label{UeslN}

The category of all finite-dimensional representations of
$U_\eps^{\text{res}}(\Sl_N\C)$ is still too large to be modular. What we
need is the category $\Tilt$ of the `reduced tilting modules'. This
is a suitable subquotient of the category of representations of
$U_\eps^{\text{res}}(\Sl_N\C)$. This means that we will consider only some
of the representations (the tilting modules) and construct $\Tilt$
by quotients of these by `negligible' parts (reduced). In this
subsection we define the subcategory of tilting modules. Here the word
module simply means representation space of $U_\eps^{\text{res}}(\Sl_N\C)$.
We define the category of reduced tilting modules in the next
subsection and defer the modular structure until even later.

Much of the early work on these representations was done by algebraists
(see Lusztig \cite{Lus} and Andersen \cite{And}) interested in the
representation theory of algebraic groups in positive characteristic. They
introduced the terminology which became standard. We adopt it here
despite the fact that it is not customary in the representation theory
of Lie groups and algebras that are the closest classical analogs. Since
prime roots are the most interesting for algebraic groups the theory
was originally developed for odd roots of unity. This makes the algebra
somewhat easier (see \cite{And} and Chari and Pressley \cite{CP}). On
the other hand, for topological applications one has to consider the
even roots of unity. In particular, the correspondence with the $SU_N$
Chern--Simons theory at level $k$ requires the order of the root to
be $l=2(k+N)$. More seriously, the category $\Tilt$ for odd roots of
unity is not modular, the non-degeneracy Axiom 17 fails(see Sawin \cite{Saw}).
This problem can actually be fixed by `modularization' but
this was realized much later and involves additional technicalities
(see Brugui\`eres \cite{Brug}).

To keep track of the difference between the even and odd cases we
define
$$l^\prime:=\begin{cases} l& l \text{ even}\\
  l/2 & l \text{ odd}\end{cases}.$$
Said in a different way $l^\prime$ is the smallest positive integer such
that $\eps^{l^\prime}=\pm 1$ where $\eps$ is a primitive $l$th root of
unity. \index{$l^\prime$ $l$ or $l/2$}

For a while in the 1990s there existed a well-developed
representation theory for odd roots of unity that did not lead to a
non-degenerate $\wtilde{s}$--matrix. There was also a non-degeneracy
proof for even roots (see Turaev and Wenzl \cite{TW}) but no
corresponding representation theory. Thus in papers and monographs
written in the 1990s authors either did not treat $3$--manifold
invariants at all \cite{CP} or implicitly assumed that the
representation theory transfers from the odd case to the even (see
the work of Bakalov, Kirillov, Reshetikhin and Turaev
\cite{BK,KRT,T}). There is still no single source where all the
required algebraic facts are stated and proved in the correct
generality.  We will state results in a form that works for both
cases but a reader interested in connections to Chern--Simons theory
may safely assume everywhere that $l=2(k+N)$.

Representations of the classical enveloping algebra $U(\g)$ are of
course the `same' as those of $\g$ meaning that every representation
of the latter extends to one of the former. The representation
theory of $U_q(\g)$ after specializing to $q=\eps$ includes some
subtleties. To stimulate intuition we begin by recalling how a
classical irreducible representation of $\Sl_2\C$ with the highest
weight $\lambda\in\Lambda_w^+$ is constructed. Looking at Appendix
\ref{app:d} first for notation and basic results from classical
representation theory may help.
\begin{example}\label{sl2}
Let $\g=\Sl_2\C$ and $\lambda=\omega=\frac12\alpha$, where
$\alpha=E_{11}-E_{22}$ is the (unique) simple root of $\Sl_2\C$ and
$\omega$ is the corresponding fundamental weight (here and below we
identify the Cartan subalgebra $\h$ with its dual $\h^*$ via the
Killing form). Recall from \fullref{UslN} that the generators
and relations for $\Sl_2\C$ are $e,f,\alpha$:
\begin{equation}\label{sl2rel}
[\alpha,e]=2e,\qquad[\alpha,f]=-2f,\qquad[e,f]=\alpha.
\end{equation}
Let $u_0$ be the highest weight vector. By definition
\begin{align*}
eu_0 &=0\\
\alpha u_0 &=(\omega,\alpha)u_0=\tfrac12(\alpha,\alpha)u_0=u_0.
\end{align*}
Now set $u_1:=fu_0$ and compute using \eqref{sl2rel}
\begin{align*}
eu_1 &=efu_0=([e,f]+fe)u_0=\alpha u_0+0=u_0\\
\alpha u_1 &=\alpha fu_0=([\alpha,f]+f\alpha)u_0=-2fu_0+fu_0=-u_1.
\end{align*}
An analogous computation shows that setting $u_2:=fu_1$ leads to
$eu_2=0$ and $\alpha u_2=-3u_2$. It follows that we could set
$fu_1=0$ and obtain a well-defined representation. Furthermore the
representation must be simple because the orbit of any non-zero vector under
$\SL_2\C$ generates the whole space. Thus, the representation space
$V_\omega$ is spanned by $u_0,u_1$. In the $u_0,u_1$ basis we have
$$e=\left(\begin{matrix} 0 & 1\\ 0 & 0 \end{matrix}\right),
\qquad\alpha=\left(\begin{matrix} 1 & 0\\ 0 & -1
\end{matrix}\right), \qquad f=\left(\begin{matrix} 0 & 0\\ 1 & 0
\end{matrix}\right),$$ which one easily recognizes as the defining
representation of $\Sl_2\C$.
\end{example}
In general, given generators $e_i,f_i,\alpha_i$ of $\Sl_N\C$ and a
dominant weight $\lambda\in\Lambda_w^+$ with highest weight vector
$u_0$ one has $e_iu_0=0$, $\alpha_i u_0=(\lambda,\alpha_i)u_0$. This
is exactly the classical Verma module $\what V_\lambda$ where $u_0$
is the equivalence class of $1$.  We keep generating new vectors
$u_{i_1,\ldots,i_m}:=f_{i_m}\ldots f_{i_1}u_0$ in lexicographic
order and keep computing the action of the generators on the new
vectors using the commutation relations. When we find the largest
proper set of vectors generated by these and the $\Sl_N\C$ action it
will be a maximal proper ideal of the Verma module and the quotient
will be the maximal abelian quotient. Said differently, we set the
action of $f_i$ to $0$ on the last vectors produced and obtain an
irreducible representation $V_\lambda(\Sl_N\C)$. More formally, we
take the quotient of the space spanned by all (infinitely many)
$u_0,u_{i_1,\ldots,i_m}$ by the maximal invariant subspace under the
action of $e_i,f_i,\alpha_i$. In \fullref{sl2} the maximal proper
invariant subspace is spanned by $u_2,u_3,\ldots$ and this is why
taking the quotient reduced to setting $fu_1=0$.
\begin{exm}\label{funs3bas}
Show that the basis of $V_{\omega_1}(\Sl_3\C)$ formed as above is
$u_0,u_1,u_{12}$ by computing the action of $e_i,f_i,\alpha_i$ as
described above.
\end{exm}
Now we wish to apply the same approach to $U_q(\Sl_2\C)$.
\begin{example}\label{Uqsl2}
The generators now are $q^{\pm{\alpha^\vee}},e,f$ and the relations
\eqref{sl2rel} get replaced by
\begin{equation}\label{Uqsl2rel}
q^{{\alpha^\vee}}eq^{-{\alpha^\vee}}=q^2e,\qquad
q^{{\alpha^\vee}}fq^{-{\alpha^\vee}}=q^{-2}f,
\qquad[e,f]=\frac{q^{{\alpha^\vee}}-q^{-{\alpha^\vee}}}{q-q^{-1}}.
\end{equation}
As in \fullref{sl2} we have
\begin{align*}
q^{\alpha^\vee} u_0&=q^{(\omega,{\alpha^\vee})}u_0=qu_0\\
eu_0&=0,\\
fu_0&=:u_1\\
q^{\alpha^\vee} u_1&=q^{\alpha^\vee} fq^{-{\alpha^\vee}}q^{\alpha^\vee} u_0=q^{-2}f\cdot qu_0=qu_0,\\
eu_1&=efu_0=([e,f]+fe)u_0=\frac{q^{\alpha^\vee}-q^{-{\alpha^\vee}}}{q-q^{-1}}u_0+0=u_0.
\end{align*}
An analogous computation for $u_2:=fu_1$ shows that $q^{\alpha^\vee}
u_2=q^{-3}u_2$ and $eu_2=0$. Hence we should set $fu_1=0$ then
$u_0,u_1$ form a basis of $\calV_{\omega}^q(\Sl_2\C)$.

It is customary for $\Sl_2\C$ to use the basis of the divided
powers $v_i:=\frac1{i!}u_i=\frac1{i!}f^iu_0$ as a canonical one.
One can show along the above lines (see Kassel \cite{kassel}) that in this
basis the action of $\Sl_2\C$ and $U(\Sl_2\C)$ in the representation
$V_{(m-1)\omega}(\Sl_2\C)$ is given by:
\begin{align*}
{\alpha^\vee} v_i &=(m-2i)v_i,\\
ev_i &=(m-i+1)\,v_{i-1},\quad i=0,1,\ldots,m-1\\
fv_i &=(i+1)\,v_{i+1}.
\end{align*}
This generalizes straightforwardly to the quantum case, where
ordinary numbers are replaced by $q$--numbers from \fullref{qnum}. The corresponding representation
$\calV_{(m-1)\omega}^q(\Sl_2\C)$ of $U_q(\Sl_2\C)$ is given by
Chari and Pressley \cite{CP}
\begin{align}\label{RepUqsl2}
q^{\pm{\alpha^\vee}}v_i &=q^{\pm(m-2i)}v_i,\notag\\
ev_i &=[m-i+1]_q\,v_{i-1},\quad i=0,1,\ldots,m-1\\
fv_i &=[i+1]_q\,v_{i+1}.\notag
\end{align}
\end{example}
Representations of interest to us are constructed from the so-called
Weyl modules that are the $\eps$ analogs of the $\calV_\lambda^q$
representations from \fullref{Uqsl2}. Their construction is
largely parallel to the construction of the corresponding classical
representations but with important caveats.
\begin{defn}
Let $\Lambda_r$ be the root lattice of $\SL_N\C$ i.e. the lattice
generated by all $\alpha_i=E_{ii}^*-E_{i+1i+1}^*$ and let
$\phi\co\Lambda_r\to\Z_2$ be a homomorphism. A
$U_q(\SL_N\C)$--weight is a pair $(\lambda,\phi)$ with $\lambda$
an ordinary $\SL_N\C$ weight and $\phi$ a homomorphism. A weight
vector in a $U_q(\SL_N\C)$ representation is a vector
$v_{(\lambda,\phi)}$ such that
$q^{\alpha_i^\vee}v_{(\lambda,\phi)}=\phi(\alpha_i)q^{\lambda(\alpha_i^\vee)}v_{(\lambda,\phi)}$.
A type I representation is one with $\phi$ being the trivial homomorphism.
\end{defn}
\begin{remark}
Only the type I representations of $U_q(\SL_N\C)$ have classical
analogs. From here forward we only consider type I representations.
\end{remark}
\begin{defn}[Weyl modules]\label{Wmod} Let $\lambda$ be a dominant
weight and let $\mathcal{I}_\lambda$ be the left ideal of
$U_q(\SL_N\C)$ generated by $e_i$ and
$q^{\pm\alpha_i^\vee}-q^{\pm\lambda(\alpha_i^\vee)}$. The Verma
module is the quotient \index{$V$@$\wwhat{\calV}_\lambda^q(\Sl_N\C)$
Verma module} \index{Verma module}
\[\wwhat{\calV}_\lambda^q(\Sl_N\C):=U_q(\SL_N\C)/\mathcal{I}_\lambda\,.\]
Denote by $\calV_\lambda^q(\Sl_N\C)$ the quotient of
$\wwhat{\calV}_\lambda^q(\Sl_N\C)$ by the maximal invariant subspace
with the induced action. The restricted integral form of this
representation $\calV_\lambda^{q,res}(\Sl_N\C)$
\index{$V$@$\calV_\lambda^{q,res}(\Sl_N\C)$ Verma quotient} is the
$\smash{U_{\Z[q,q^{-1}]}^{\text{res}}}(\Sl_N\C)$ submodule of
$\wwhat{\calV}_\lambda^q(\Sl_N\C)$ generated by $1$.  The Weyl
module \index{$\calW_\lambda^\eps(\Sl_N\C)$ Weyl module}\index{Weyl
module} $\calW_\lambda^\eps(\Sl_N\C)$is the vector space over $\C$
generated from $\calV_\lambda^{q,res}(\Sl_N\C)$ by changing
coefficients from $\Z[q,q^{-1}]$ to $\C$ ($\pi\mapsto\pi(\eps)$)
with the action of $U_\eps^{\text{res}}(\Sl_N\C)$ obtained by
specializing $q$ to $\eps$. Formally,
$$
\calW_\lambda^\eps(\Sl_N\C):=\calV_\lambda^{q,res}(\Sl_N\C)\otimes_{\Z[q,q^{-1}]}\C.
$$
\end{defn}
This definition seems a bit convoluted but it is in essence parallel
to the definition of $U_\eps^{\text{res}}(\Sl_N\C)$ itself. We are
trying to avoid the trivialization of the powers of $e_i,f_i$ by
making sure that the divided powers are present as `independent'
quantities. The space $\wwhat{\calV}_\lambda^q(\Sl_N\C)$ is just the
vector space over $\C(q)$ generated by $u_0,u_{i_1,\ldots,i_m}$ for
$i_k=1,\ldots,N-1$ with the action of $U_q(\Sl_N\C)$ determined by
$$
e_iu_0=0,\quad q^{\pm\alpha_i}u_0=q^{\pm(\lambda,\alpha)}u_0, \quad
f_iu_{i_1,\ldots,i_m}=u_{i_1,\ldots,i_m,i}\,.
$$
To get to the Weyl modules we restrict to an integral form then
specialize coeficients to $\C$.
\begin{example}\label{Uesl2}
Recall from \fullref{Uqsl2} that
$\calV_{(m-1)\omega}^q(\Sl_2\C)$ is spanned by $v_0,\ldots,$
$v_{m-1}$. Iterating the action \eqref{RepUqsl2} we obtain
$$
e^nv_i=[m-i+1]_q\ldots[m-i+n]_q\,v_{i-n}=\frac{[m-i+n]_q!}{[m-i]_q!}\,v_{i-n}.
$$
and analogously for $f^n$. Using \fullref{divpw} of the divided powers and
\eqref{qbin} of $q$--binomials we get
\begin{equation}
\label{RepUesl2}
\begin{aligned}
e^{(n)}v_i &=\left[m-i+n\atop n\right]_q\,v_{i-n}, \\
f^{(n)}v_i &=\left[i+n\atop n\right]_q\,v_{i+n}.
\end{aligned}
\end{equation}
Since $q$--binomials are in $\Z[q,q^{-1}]$ we see that
$v_0,\ldots,v_{m-1}$ also form a basis of
$\calV_{(m-1)\omega}^{\text{q,res}}(\Sl_2\C)$ with the action given by
\eqref{RepUesl2}. Thus by \fullref{Wmod} the Weyl module
$\calW_{(m-1)\omega}^\eps(\Sl_2\C)$ is spanned by the same vectors
and the action of $U_\eps^{\text{res}}(\Sl_2\C)$ on them is given by
replacing $q$ with $\eps$ in \eqref{RepUesl2}. It can be shown
(see Chari and Pressley \cite{CP}) that
$\calW_{(m-1)\omega}^\eps(\Sl_2\C)$ 
has an invariant subspace $\calW_{\inv}$ unless $m\leq l^\prime$ or
$m\equiv l^\prime-1\mod l^\prime$; see \fullref{ex1727}. This means that
for other values of $m$ not only is $\calW_{(m-1)\omega}^\eps(\Sl_2\C)$
not irreducible but it is not even a direct sum of irreducibles. Thus,
complete reducibility of the classical representations that still holds
for $U_q(\Sl_N\C)$ is lost for $U_\eps^{\text{res}}(\Sl_N\C)$.
\end{example}
This example demonstrates an important method of doing computations
in $U_\eps^{\text{res}}(\Sl_N\C)$ that avoids using its complicated
relations directly (and this is the reason we did not give a
complete list of them in \eqref{Uqdiv}). This idea will be used
again and again in the sequel.

\begin{quote}
To obtain an equality in $U_\eps^{\text{res}}(\Sl_N\C)$ perform all
the computations in $U_q(\Sl_N\C)$ using \eqref{Uqalg} and rewrite
the end result in terms of the divided powers so that it only
contains Laurent polynomials as coefficients. Then specializing $q$
to $\eps$ gives an equality in $U_\eps^{\text{res}}(\Sl_N\C)$.
\end{quote}

Before dealing with the loss of complete reducibility we have to
address a more basic problem with the definition of a weight space.
Recall from \fullref{app:d} that classically a vector $v$ has a
weight $\lambda\in\Lambda_w$ if $\alpha_i v=(\lambda,\alpha_i)v$ for
all simple roots $\alpha_i$. For $q$ an indeterminate this
generalizes straightforwardly to the quantum case by setting
$q^{\pm\alpha_i}v=q^{\pm(\lambda,\alpha_i)}v$ instead. This also
works when we specialize $q$ to a generic complex number $\eps$. If
however $\eps$ is a root of unity and $\beta\in\Lambda_r$ we have
$\eps^{(\lambda+l\beta,\alpha_i)}=\eps^{(\lambda,\alpha_i)}$ and
$\lambda$ would only be defined modulo $l\Lambda_r$.

The underlying reason for the weight ambiguity is that for roots of
unity the maximal Abelian subalgebra of $U^{\text{res}}_\eps$ is no longer
generated by $q^{\pm\alpha_i}$ \cite{CP}. The additional generators
are the ones we already met in \eqref{cengen}
$\smash{\sbinom{q^{\alpha_i};0}{l^\prime}_{\eps}}$.
Note that substituting $\eps$ into
`definition' \eqref{cengen} leads to a meaningless expression. To
make sense of these elements in $U^{\text{res}}_\eps$ one has to reexpress
them in terms of the divided powers. Below we incorporate these new
elements into the definition of a weight space so that the weight is
now well-defined.
\begin{defn}\label{wspace}
Let $\calV$ be a representation space of $U^{\text{res}}_\eps(\Sl_N\C)$. We
say that $v\in\calV$ is a weight vector \index{weight vector} with
weight $\lambda\in\Lambda_w$ if
$$q^{\pm\alpha_i}v=\eps^{\pm(\lambda,\alpha_i)}v, \qquad
  \sbinom{q^{\alpha_i};0}{l^\prime}_{\eps} v
  =\sbinom{(\lambda,\alpha_i)}{l^\prime}_{\eps} v,$$
where on the right we have the $q$--binomial coefficient
\eqref{qbin} rewritten as a Laurent polynomial and specialized to
$\eps$. We denote the subspace of vectors in $\calV$ with weight
$\lambda$ by $\calV^\lambda$ and call it the $\lambda$--weight space.
\index{$V$@$\calV^\lambda$ weight space}  \index{weight space} If
$\calV=\oplus_{\lambda\in\Lambda_w}\calV^\lambda$ the righthand side
is called the weight space decomposition of $\calV$.
\end{defn}
\begin{remark}
Notice that we are using $V_\lambda$ to represent the irreducible
representation with heighest weight $\lambda$ and $V^\lambda$ to
denote the $\lambda$--weight space of a representation $V$. Generally
we will use calagraphic fonts to denote representations of quantum
groups and roman fonts to denote representations of classical
algebras.
\end{remark}

By \eqref{Uqdiv} we have
\begin{multline*}
q^{\alpha_j}(e_i^{(n)}v)=(q^{\alpha_j}e_i^{(n)}q^{-\alpha_j})q^{\alpha_j}v
=q^{na_{ji}}q^{(\lambda,\alpha_i)}v\big|_{q=\eps} \\
=\eps^{n(\alpha_i,\alpha_i)+(\lambda,\alpha_j)}
=\eps^{(\lambda+n\alpha_i,\alpha_j)}v.
\end{multline*}
As we explained this in itself does not mean that $e_i^{(n)}v$ has
the weight $\lambda+n\alpha_i$ but one can show using the full list
of relations in Lusztig \cite{Lus} or Chari and Pressley \cite{CP} that indeed
$$
e_i^{(n)}(V^\lambda)\subseteq V^{\lambda+n\alpha_i},\quad
f_i^{(n)}(V^\lambda)\subseteq V^{\lambda-n\alpha_i}.
$$
With this notation we have a very important result that follows from
the definitions by a deformation argument \cite{CP}:
\begin{quote}
Weight space decompositions of
$V_\lambda,\calV_\lambda^q,\calW_\lambda^\eps$ are the same, that
is, their weights and the dimensions of their weight spaces are
equal.
\end{quote}
Since weight space decompositions of classical representations are
well-known this observation comes handy when performing computations
with representations of quantum groups.

Circumventing the lack of complete reducibility is not as simple.
Ultimately, we will have to restrict to the class of admissible
representations for which complete reducibility still holds. This
will lead to the desired finite number of irreducibles. To proceed
in this direction we need to introduce the notions of the dual of a
representation and of the tensor product of representations. Recall
that given vector spaces $\calU,\calV$ with a linear action of a Lie
group $G$ the actions on $\calV^*$ and $\calU\otimes\calV$ are given
by
$$
gf(v):=f(g^{-1}v),\quad g(u\otimes v):=gu\otimes gv.
$$
Quantum groups $U_q(\Sl_N\C),U_\eps^{\text{res}}(\Sl_N\C)$ are neither Lie
groups nor Lie algebras but Hopf algebras so the inverse and the
tensor action are replaced by the antipode and the coproduct
respectively.
\begin{defn}\label{dualtensor}
Given representations $\calA\to\End(\calU)$, $\calA\to\End(\calV)$
of a Hopf algebra $\calA$ (see \fullref{Hopfalg}) the dual
representation on $\calV^*$ \index{$V$@$\calV^*$ dual
representation} is given by $af(v):=f(\gamma(a)v)$ and the tensor
representation on $\calU\otimes\calV$ \index{$U$@$\calU\otimes\calV$
product representation} is given by $a(u\otimes
v):=\Delta(a)(u\otimes v)$. The unit representation
${\bbone}:=\calA\to\End(\C)$ is given by $az:=\ve(a)z$.
\end{defn}
This definition works because $\gamma,\Delta$ are an
antihomomorphism and a homomorphism respectively. For Lie group
algebras and universal enveloping algebras with their usual Hopf
structure we recover the standard definitions of the unit (trivial
representation) dual and the tensor product. This provides the unit,
dual and tensor product in the modular category that we are
defining.

Recall from \fullref{app:d} that for classical Lie algebras all
irreducible representations are indexed by dominant weights. Since
obviously the dual to an irreducible is an irreducible we get a
duality involution on the set of dominant weights $\Lambda_w^+$. For
$\Sl_N\C$ it can be easily described explicitly (see Fulton and
Harris \cite{fulton-harris} or Humphreys \cite{Hum}). Let
$w_0=\bigl(\begin{smallmatrix} 1 & 2 & \ldots & N-1 & N\\ N & N-1 &
\ldots & 2 & 1\end{smallmatrix}\bigr)$ be the order-reversing
permutation then $V_\lambda^*\simeq V_{-w_0(\lambda)}$ and this is
an involution because $w_0^2=\id$.
\begin{exm}
Check that this indeed works for a couple of representations of
$\Sl_3\C$.
\end{exm}
More importantly for us, this carries over to the quantum groups, in
particular,
\begin{equation}\label{irdual}
\calV_\lambda^{q*}(\Sl_N\C)\simeq\calV_{-w_0(\lambda)}^q(\Sl_N\C).
\end{equation}
The next definition introduces the type of representations that we
will use to build our modular category. These `tilting' modules were
originally introduced and studied in the context of algebraic
groups. In fact many of the proofs refer to facts that are true by
analogy with results from algebraic groups, and we do not know of a
reference that addresses the representation theory of quantum groups
that we need without assuming familiarity with algebraic groups. For
example Chari and Pressley refer to algebraic groups as the
`classical' case \cite{CP}. See H\,H Andersen \cite{And} for more
history and some important results related to these modules.
\begin{defn}[Tilting modules]\label{Tmod}
A representation $V$ of $U^{\text{res}}_\eps(\Sl_N\C)$ is said to
have a Weyl filtration if there exists an increasing sequence of
invariant subspaces $0=V_0\subset V_1\subset\cdots\subset V_r=V$
with $V_i/V_{i-1}\simeq\calW_{\lambda_i}^\eps(\Sl_N\C)$ for some
dominant weights $\lambda_i\in\Lambda_w^+$. A representation $V$ is
called a tilting module \index{tilting module} if both it and its
dual $V^*$ have Weyl filtrations. The category of all tilting
modules is denoted by $\mathcal{T}\bthin\mathit{ilt}_\eps(\Sl_N\C)$
\index{$T$@$\mathcal{T}\bthin\mathit{ilt}_\eps(\Sl_N\C)$ tilting
modules}
\end{defn}
Note that if we had complete reducibility, the existence of the Weyl
filtration is equivalent to $V$ being a direct sum of
$\calW_{\lambda_i}^\eps$--s. However, such direct decomposition does
not hold in general, that is, the filtration 'tilts'. In particular,
the tilting modules still are not completely reducible and we need a
condition weaker than irreducibility to describe the `elementary'
tilting modules.
\begin{defn}\label{Indec}
A representation $V$ is called indecomposable \index{indecomposable
representation}if it does not split into a direct sum of two proper
invariant subspaces $V\neq V_1\oplus V_2$.
\end{defn}
It follows from \fullref{Uesl2} that in the $\Sl_2\C$ case every
Weyl module is indecomposable. Also every Weyl module obviously has
a trivial Weyl filtration. However, not every Weyl module is
tilting. The problem is that its dual does not necessarily have a
Weyl filtration. This is closely related to the fact that some Weyl
modules are reducible. The following exercise provides a good
example to think about when studying these issues.
\begin{exm}
\label{ex1727}
Take $N=2$ and $k=3$ giving $l=10$ and $\eps=e^{\pi i/5}$. Show that
the subspace of $\calW_{7\omega_1}^{\eps}(\Sl_2\C)$ generated by
$f^{(3)}u_0$ and $f^{(4)}u_0$ is an invariant subspace
($[5]_{\eps}=0$) so that $\calW_{7\omega_1}^{\eps}(\Sl_2\C)$ is
reducible. Use this invariant subspace to explicitly construct a
non-split extension
\[
0\to
\calW_{7\omega_1}^{\eps}(\Sl_2\C)\to\calQ_{7\omega_1}^\eps\to\calW_{\omega_1}^{\eps}(\Sl_2\C)\to
0\,.
\]
The module $\calQ_{7\omega_1}^\eps$ is the unique
indecomposable tilting module of weight $7\omega_1$.
\end{exm}

Nonetheless, we have the following major theorem originally due to
H\,H Andersen \cite{And} for $l$ odd. The case of even $l$ should
follow the general arguments of \cite{And} but it is not stated
explicitly there.
\begin{thm}\label{decomp}
Direct sums and summands, duals and tensor products of tilting
modules are again tilting modules. For every $\lambda\in\Lambda_w^+$
there is a unique indecomposable tilting module
$\calQ_{\lambda}^\eps$ with the highest weight $\lambda$ and
one-dimensional highest weight space. Moreover, every tilting module
$V$ admits a decomposition
$$
V=\oplus_{\lambda\in\Lambda_w^+}\left(\calQ_{\lambda}^\eps\right)^{\oplus m_{\lambda}(V)}
$$
with multiplicities $m_{\lambda}(V)$ canonically determined by $V$
and only finitely many of them non-zero. The dual is given by
$\calQ_{\lambda}^{\eps*}\cong\calQ_{-w_0(\lambda)}^\eps$.
\end{thm}
Thus, tilting modules form a subcategory
$\mathcal{T}\!\mathit{ilt}_\eps(\Sl_N\C)$ of the category of
representations of $U_\eps^{\text{res}}(\Sl_N\C)$ with morphisms
being equivariant maps $f(av)=af(v)$, which is closed under duality
and tensor products and is dominated by indecomposable objects
$\calQ_{\lambda}^\eps$. In a way, the tilting modules resemble the
classical representations much more than general representations of
$U_\eps^{\text{res}}(\Sl_N\C)$ do. However, there are two
difficulties that prevent $\mathcal{T}\!\mathit{ilt}_\eps(\Sl_N\C)$
from being modular as is: $\calQ_{\lambda}^\eps$ are not exactly
simple objects being indecomposable but not irreducible and there
are still infinitely many of them. Based on ideas from physics or
more likely ideas from algebraic groups Andersen was able to resolve
both problems by discarding tilting modules of quantum dimension
$0$. This is explained in the next subsection.

\subsection{Quantum dimensions and the Weyl alcove}

We do not have a ribbon structure on
$\mathcal{T}\!\mathit{ilt}_\eps(\Sl_N\C)$ yet, so defining quantum
traces and dimensions as it was done in \fullref{mtcsec} is not
possible. Recall however that quantum traces were reinterpreted in
\fullref{color} in terms of double duals. Namely, given a double
dual isomorphism $\smash{V\xrightarrow{\delta_V}V^{**}}$ one can set
$\Tr_q(f):=\cap_{V^*}\circ(\delta_V\circ
f\otimes\id_{V^*})\circ\cup_V$ and $\dim_q(V):=\Tr_q(\id_{V})$,
where $\cap_{V^*},\cup_V$ are the standard pairing and copairing for
vector spaces. Any candidate for $\delta_V$ must of course be
equivariant under the $U_\eps^{\text{res}}$ action. The standard
identification of $V$ and $V^{**}$ in the category of
finite-dimensional vector spaces $v\mapsto[v]$ with
$[v](\varphi):=\varphi(v)$ is not equivariant. Indeed, we have
$$
(a[v])(\varphi)=[v](\gamma(a)\varphi)=(\gamma(a)\varphi)(v)=\varphi(\gamma^2(a)v)=(\gamma^2(a)[v])(\varphi)
$$
and one easily sees from \eqref{Uqcoalg} that $\gamma^2\ne\id$. The
idea is to `fix' $[\cdot]$ to make it equivariant.
\begin{prop}\label{qddual} Let $\rho:=\frac12
\sum_{\alpha\in\Delta^+}\alpha$ be the Weyl weight (see Appendix
\ref{app:d}). Then for $a$ in $U_q(\Sl_N\C)$ or in
$U_\eps^{\text{res}}(\Sl_N\C)$ we have
\begin{equation}\label{Qgdual}
\gamma^2(a)=q^{2\rho}aq^{-2\rho}
\end{equation}
and the \index{double dual isomorphism}\index{$\delta$@$\delta_V$
double dual isomorphism} double dual map $\delta_V\co V\to V^{**}$
given by $v\mapsto [q^{2\rho}v]$ is an equivariant isomorphism.
Consequently,
\begin{equation}\label{Qgtr}\index{quantum trace}
\Tr_q(f)=\Tr(q^{2\rho}f),
\end{equation}
where $\Tr$ is the usual trace.
\end{prop}
\begin{proof}
Since both sides of \eqref{Qgdual} are algebra homomorphisms it
suffices to check the equality on the generators
$q^{\pm{\alpha^\vee}_i},e_i,f_i$. We have
$\gamma^2(q^{\pm{\alpha^\vee}_i})=q^{\mp{\alpha^\vee}_i}$ by
\eqref{Uqcoalg} and the equality is obvious. For $e_k$ we have from
\eqref{Uqalg} and \eqref{Uqcoalg}:
$$
\gamma^2(e_k)=\gamma(-e_kq^{-{\alpha^\vee}_k})=-\gamma(q^{-{\alpha^\vee}_k})\gamma(e_k)=q^{{\alpha^\vee}_k}e_kq^{-{\alpha^\vee}_k}=q^2e_k.
$$
On the other hand
$q^{2\rho}e_kq^{-2\rho}=q^{(\alpha_k,2\rho)}e_k=q^2e_k=\gamma^2(e_k)$
as claimed. (Recall from \fullref{app:d} that $\rho=\sum \omega_i$.)
In $U_q(\Sl_N\C)$ the equality for the divided powers $e_k^{(n)}$
follows by dividing both sides by $[n]_q!$. The case of $f_k^{(n)}$
is analogous. Since \eqref{Qgdual} does not involve any rational
functions of $q$ it remains valid in
$U_{\Z[q,q^{-1}]}^{\text{res}}(\Sl_N\C)$ and therefore specializes
to $U_\eps^{\text{res}}(\Sl_N\C)$. The equivariance is now
straightforward
$$
a\delta_V(v)=a[q^{2\rho}v]=[\gamma^2(a)q^{2\rho}v]=[q^{2\rho}av]=\delta_V(av).
$$
Finally, for the quantum trace we get
\begin{multline*}
\Tr_q(f):=\cap_{V^*}\circ(\delta_V\circ f\otimes\id_{V^*})
\biggl(\sum_k v_k\otimes v^k\biggr)
=\cap_{V^*}\biggl(\sum_k[q^{2\rho}f(v_k)]\otimes v^k\biggr)\\
=\sum_kv^k\biggl(q^{2\rho}f(v_k)\biggr)=\Tr(q^{2\rho}f).
\end{multline*}
This completes the proof.
\end{proof}

\begin{remark}\label{someinterest}
Something interesting is happening here. We are starting with a
representation of a deformation of the universal enveloping algebra
of ${\mathfrak{sl}}_N$. In order to define the appropriate
deformation  we introduced elements such as $q^\beta$ in place of
the coroots. Furthermore, because we are exponentiating ($q=e^{2\pi
i/l}$ after specialization giving $q^\beta=e^{2\pi i\beta/l}$) these
are actually elements of the group SU$(N)$ so the characteristic
numbers that we compute here can be expressed in terms of characters
of the group SU$(N)$. The Weyl character formula can then be used to
compute the resulting characters.
\end{remark}

Recall from the previous article that the Weyl modules
$\calW_\lambda^{\eps}(\Sl_N\C)$ have the same weight space
decompositions as the classical representation spaces $V_\lambda$.
The quantum dimension is the quantum trace of the identity, and it
is given by the Weyl character formula (see \fullref{app:d}).
\begin{cor} The quantum dimensions of the space $\calW_\lambda^{\eps}$
is given by
\begin{equation}\label{Wdim}\index{quantum
dimension}\index{$D$@$\text{dim}_q(V)$ quantum dimension}
\dim_q\!\calW_\lambda^{\eps}=\prod_{\alpha\in\Delta^+}
\frac{\eps^{(\lambda+\rho,\alpha)}-\eps^{-(\lambda+\rho,\alpha)}}{\eps^{(\rho,\alpha)}-\eps^{-(\rho,\alpha)}}.
\end{equation}
\end{cor}
\begin{exm}
Since $\eps^{-1}=\wbar\eps$ it is obvious that
$\dim_q\!\calW_\lambda^{\eps}$ is real. Prove that in fact
$\dim_q\!\calW_\lambda^{\eps}\geq0$ if $\eps$ is a primitive root of
unity.
\end{exm}
It is obvious from \eqref{Wdim} that
$\dim_q\!\calW_\lambda^{\eps}=0$ if and only if
$(\lambda+\rho,\alpha)$ is divisible by $l^\prime$ for some
$\alpha\in\Delta^+$. Therefore, for $\lambda$ in the range
$0<(\lambda+\rho,\alpha)<l^\prime$ for all positive roots $\alpha$
the quantum dimension of $\calW_\lambda^{\eps}$ is non-zero.
\begin{example}\label{Wlsl2}
For the Weyl modules $\calW_{(m-1)\omega}^\eps(\Sl_2\C)$ from
\fullref{Uesl2} the Weyl weight $\rho:={\alpha}/2$, where
$\alpha$ is the only positive (and also simple) root of $\Sl_2\C$.
Therefore
$$
\bigl((m-1)\omega+\tfrac{\alpha}2,\alpha\bigr)=m-1+\tfrac12(\alpha,\alpha)=m
$$
and $\dim_q\!\calW_{(m-1)\omega}^\eps(\Sl_2\C)=0$ if and only if
$l^\prime$ divides $m$. The condition
$0<(\lambda+\rho,\alpha)<l^\prime$ is in turn equivalent to
$0<m<l^\prime$. Note that by what we mentioned in \fullref{Uesl2} when $m=l^\prime$ the Weyl module is still irreducible
even though its quantum dimension is already $0$.
\end{example}
It turns out that the Weyl modules with the highest weights in the
range $0<(\lambda+\rho,\alpha)<l^\prime$ are the most important ones
for our purposes.
\begin{defn}[Weyl alcove]\label{Walcove}
The (open) Weyl alcove \index{open Weyl alcove}\index{$C^l$ open
Weyl alcove} is a subset of the Cartan subalgebra $\h$ of $\Sl_N\C$
defined by
$$
C^l:=\{x\in\h\mid\,0<(x+\rho,\alpha)<l^\prime,\ \text{for all}\
\alpha\in\Delta^+\},
$$
where $l^\prime=l$ for $l$ odd and $=l/2$ for $l$ even. We also
denote \index{$\Lambda_w^l:=C^l\cap\Lambda_w$ Weyl alcove}
\index{Weyl alcove} the set of dominant weights in the Weyl alcove
by $\Lambda_w^l:=C^l\cap\Lambda_w$ and call this the Weyl alcove
when there is no confusion. The closure \index{$C$@$\wwbar{C}^l$
closed Weyl alcove} $\wwbar{C}^l$ \index{closed Weyl alcove}is
called the closed Weyl alcove. We will often abuse notation by
saying that a Weyl module is or is not in the Weyl alcove according
to the location of its highest weight.
\end{defn}
We already saw that the modules in the alcove have non-zero quantum
dimensions. In fact, much more is true (see Andersen \cite{And},
Andersen and Paradowski \cite{AP}, Chari and Pressley \cite{CP}, and Sawin
\cite{Saw}):
\begin{thm}\label{alcove}
The Weyl alcove $\Lambda_w^l$ contains precisely all the (highest
weights of) Weyl modules $\calW_\lambda^\eps$ that are both
irreducible and have non-zero quantum dimensions. An indecomposable
module \index{$\calQ_{\lambda}^\eps$ indecomposable
module}\index{indecomposable module} $\calQ_{\lambda}^\eps$ (see
\fullref{decomp}) has non-zero quantum dimension if and only if
$\lambda\in\Lambda_w^l$ in which case
$\calQ_{\lambda}^\eps=\calW_\lambda^\eps$. If
$\lambda\in\Lambda_w^l$ then
$\lambda^*:=-w_0(\lambda)\in\Lambda_w^l$.
\end{thm}
Recall that $\calQ_\lambda^{\eps*}\simeq\calQ_{\lambda^*}^\eps$ so
the Weyl alcove is closed under duality. In practice for $\SL_N\C$,
it is more convenient to write the defining condition of the alcove
as follows
\begin{equation}\label{pralcove}
\begin{cases}
0<(x+\rho,\alpha_i), & \text{for all simple roots $\alpha_i$}\\
(x+\rho,\theta)<l^\prime, & \text{for the highest positive root $\theta$.}
\end{cases}
\end{equation}
We proceed with some sample computations that make \eqref{pralcove}
more explicit.
\begin{example}\label{emalcove}
Recall from \fullref{app:d}  that
$L_i:=E_{ii}-\frac1NI\,i=1,\ldots,N-1$ form a basis in $\h=\h^*$ of
$\Sl_N\C$. Given a weight it is convenient to write
$\lambda:=\sum_{i=1}^N\lambda_iL_i$ by setting $\lambda_N:=0$. The
highest root is $\theta=E_{11}-E_{NN}=\sum_{i=1}^{N-1}\alpha_i$ with
$\alpha_i=E_{ii}-E_{i+1\,i+1}=L_i-L_{i+1}$. One easily checks that
$(L_j,\alpha_i)=\delta_{ij}-\delta_{i+1\,j}$ and since
$\rho=\sum_{i=1}^{N-1}\omega_i$ we have
\begin{align*}
(\lambda+\rho,\theta)&=\sum_{i=1}^{N-1}\sum_{j=1}^{N}\lambda_j(L_j,\alpha_i)
+\sum_{i=1}^{N-1}\sum_{j=1}^{N-1}(\omega_j,\alpha_i)\\
&=\sum_{i=1}^{N-1}\sum_{j=1}^{N}\lambda_j(\delta_{ij}-\delta_{i+1\,j})
+\sum_{i=1}^{N-1}\sum_{j=1}^{N-1}\delta_{ij} \\
&=\sum_{i=1}^{N-1}(\lambda_i-\lambda_{i+1}+1) \\
&=\lambda_1+N-1.
\end{align*}
Hence the second condition in \eqref{pralcove} is equivalent to
$\lambda_1\leq l^\prime-N$. By analogous computation the first
condition reduces to $\lambda_i\geq\lambda_{i+1}$ for all $i$. This
means that $\lambda_i$ is a partition and the weights $\lambda$ from
the alcove are always dominant: $\Lambda_w^l\subset\Lambda_w^+$.
Thus \eqref{pralcove} reduces to $\lambda$ being dominant and
$0\leq\lambda_i\leq l^\prime-N$ for all $i$. Note that in the case
of Chern--Simons theory where $l=2(k+N)$ this leads to an
$N$--independent condition $0\leq\lambda_i\leq k$. Another direct
computation shows that
$\lambda^*=\sum_{i=1}^N(\lambda_1-\lambda_{N-i+1})L_i$.
\end{example}
\begin{exm}
Show that in the case of $\Sl_2\C$ we have $\lambda^*=\lambda$, that is, all Weyl modules are self-dual.
\end{exm}
\begin{exm}\label{wLi}
The symmetric group $W=\mathfrak{G}_N$ acts on $E_{ii}$ by $\sigma
E_{ii}=E_{\sigma(i)\sigma(i)}$. Compute the induced action on $L_i$:
$$
\sigma L_i=\begin{cases}
L_{\sigma(i)},&\sigma(i)<N\\
-\sum_{j=1}^{N-1}L_j,&\sigma(i)=N
\end{cases}
$$
and use it to derive a formula for $\lambda^*=-w_0(\lambda)$.
\end{exm}
\begin{exm}\label{nidual}
Recall from \fullref{app:d} that one can also express any weight as
a sum of fundamental weights $\lambda=\sum_{i=1}^{N-1}n_i\omega_i$
and $n_i\geq0$ for the dominant weights. Show that \eqref{pralcove}
is equivalent to $n_i\geq0$ and $\sum_{i=1}^{N-1}n_i\leq
l^\prime-N$. Also in terms of $n_i$s we have
$\lambda^*=\sum_{i=1}^{N-1}n_{N-i}\omega_i$
\end{exm}
\begin{exm}
In the Cartan subalgebra $\h$ of $\Sl_N\C$ neither $\alpha_i$ nor
$\omega_i$ form an orthonormal basis. For $\Sl_3\C$ orthonormalize
$\alpha_1,\alpha_2$ into $u_1,u_2$ and show that
$$\alpha_1=\sqrt{2}u_1,\qquad
  \alpha_2=\sqrt{2}\bigl(-\tfrac12u_1+\tfrac{\sqrt{3}}{2}u_2\bigr),\qquad
  \theta=\alpha_1+\alpha_2=
  \sqrt{2}\bigl(\tfrac12u_1+\tfrac{\sqrt{3}}{2}u_2\bigr).$$
Given $x=x_1u_1+x_2u_2$ describe the condition $x\in C^l$ in terms
of $x_1,x_2$ and draw a picture of $C^l$ (cf. \fullref{translatedI} that shows $\rho+C^l$).
\end{exm}
The properties of the Weyl alcove indicate how to proceed with the
construction of a modular category. First of all, we finally have a
finite set $\Lambda_w^l$ of simple objects and they are indeed
simple because alcove Weyl modules are irreducible. To make them
dominate our category it suffices to consider only tilting modules
that decompose into direct sums of alcove Weyl modules. Since dual
alcove modules are still in the alcove the duals to such tilting
modules will again decompose into direct sums of the alcove Weyl
modules. This idyllic picture is spoiled by the behavior of the
tensor product: the tensor product of two alcove modules may have
non-alcove modules in its decomposition. This means that if we want
a category with `tensor products' the usual tensor product has to be
redefined. The idea is to discard the submodule of the tensor
product that has quantum dimension zero. Since
$\calQ_{\lambda}^\eps$ has positive dimension if and only if it is
an alcove Weyl module keeping only the latter should do the trick.
\begin{defn}\label{redtensor} Let
$V=\oplus_{\lambda\in\Lambda_w^+}\left(\calQ_{\lambda}^\eps\right)^{\oplus
m_{\lambda}(V)}$ be a decomposition of a tilting module $V$ into the
indecomposables as given by \fullref{decomp}. Then its reduction
is defined by keeping only summands with highest weight in the Weyl
alcove, that is,
\index{$V$@$\wwbar{V}:=\oplus_{\lambda\in\Lambda_w^l}\left(\calQ_{\lambda}^\eps\right)^{\oplus
m_{\lambda}(V)}$ reduction of $V$}\index{reduction}
$$\wwbar{V}:=\smash{\oplus_{\lambda\in\Lambda_w^l}}
\bigl(\calQ_{\lambda}^\eps\bigr)^{\oplus m_{\lambda}(V)}.$$
A tilting module is said to be reduced if
\index{$T$@$\Tilt(\Sl_N\C)$ reduced tilting modules}
\index{negligible} $\wwbar{V}=V$ and negligible if
$\wwbar{V}=0$. The reduced tensor product of \index{reduced
tensor product} two reduced modules $U,V$ is defined
\index{$U\rotimes V:=\overline{U\otimes V}$ reduced tensor product}
by
$$U\rotimes V:=\overline{U\otimes V}.$$
\end{defn}
\begin{remark}
Actually, we only defined $\wwbar{V}$ up to isomorphism. For a
strict category one needs a canonical construction of it which goes
as follows (see Chari and Pressley \cite{CP}). Let $\widetilde{V}$ be the maximal reduced
submodule of $V$ and $V^\prime$ its maximal negligible submodule
then
$\wwbar{V}=\widetilde{V}/(V^\prime\cap\widetilde{V})$
is the canonical representative of the reduction. In physics
literature reduced tensor product is often called fusion tensor
product and rules for computing it are called fusion rules.
\end{remark}
The following theorem also due to Andersen \cite{And} (see also \cite{CP})
indicates that the new tensor product behaves `properly'.
\begin{thm}\label{prodok}
If $V$ is any tilting module and $U$ is negligible then so are
$U\otimes V$ and $V\otimes U$. Thus, if $U,V,W$ are reduced we have
canonical isomorphisms $U\rotimes V\simeq V\rotimes U$ and
$(U\rotimes V)\rotimes W\simeq U\rotimes (V\rotimes W)$
\end{thm}
\begin{exm}\label{neglmorph}
A morphism $U\xrightarrow{f}V$ of tilting modules is called neglible
\index{neglible morphism} if it factors through a neglible tilting
module. An equivalence class of equivariant morphisms is called a
non-negligible morphism. \index{non-negligible morphism} Prove that
$f$ is neglible if and only if for any morphism $V\xrightarrow{g}U$
we have $\Tr_q(fg)=\Tr_q(gf)=0$.
\end{exm}
\begin{remark}
Note that we do not lose any information about link invariants by
discarding tilting modules of quantum dimension zero. It follows
from \fullref{Jones} that the colored invariant
$J_{V_1,\ldots,V_k}(L)$ is $0$ if any one of $V_i$ is negligible.
Indeed, we can always elongate the strand colored by $V_i$ so that
it cups under and caps over the rest of the link. The entire link
then reduces to a morphism $\smash{V_i\xrightarrow{f_L}V_i}$ with
$J_{V_1,\ldots,V_k}(L)=\Tr_q(f_L)$. But $f_L$ is negligible along
with $V_i$ and so its quantum trace is zero. This serves as
topological justification for discarding the negligible modules.
\end{remark}

\begin{defn}\label{nonneg}
We denote the category of reduced tilting modules with equivariant
linear maps by $\Tilt(\Sl_N\C)$.
\end{defn}
\begin{remark}
There is an alternative construction of this category used eg in
Bakalov--Kirillov \cite{BK}. One keeps the ordinary tensor product but
takes morphisms to be equivalence classes of equivariant linear maps
modulo negligible ones (so called non-negligible morphisms). The category
so defined turns out to be isomorphic to $\Tilt(\Sl_N\C)$, see Chari--Pressley \cite{CP}.
\end{remark}
The category $\Tilt(\Sl_N\C)$ together with the duality and reduced
tensor products is exactly what we were looking for.  However, two
key ingredients are still missing: the braiding and the twist. The
latter can be restored from the former using dual cups and caps as
in \fullref{color}. But to get a braiding we need to introduce a
new structure on $U_\eps^{\text{res}}(\Sl_N\C)$ known as quasitriangular
structure or the $R$--matrix. After defining the braiding and twist
we will see that for even roots of unity $\Tilt(\Sl_N\C)$ will turn
out to be modular.

\subsection{$R$--matrices and braiding}\label{Rbraid}
We are making good progress defining a concrete nontrivial modular
tensor category. By considering appropriate representations of
$U_\eps^{\text{res}}(\Sl_N\C)$ as representations of an associative
algebra, we were able to construct the objects and morphisms of an
abelian category $\Tilt$ dominated by a finite collection of simple
objects. The Hopf algebra structure of $U_\eps^{\text{res}}(\Sl_N\C)$, in
particular the counit, the antipode and the coproduct  were
responsible for providing additional structures in the category of
representations, namely the unit object, the duality and the tensor
product. In this article we define the appropriate braiding in our
category.

The braiding requires something beyond the Hopf algebra structure.
As motivation, recall that the (trivial) braiding on $\text{REP}_G$
from \fullref{repg} is just the flip $\times_{U,V}(u\otimes
v)=v\otimes u$. Of course, given any Hopf algebra $\calA$ we can
define a braiding $\sigma$ on $\calA\otimes\calA$ by $a\otimes
b\mapsto b\otimes a$ and this will induce a trivial `braiding' on
any pair of representations. The catch is to find a necessarily
nontrivial braiding compatible with the rest of the structure -- in
particular leading to the same definitions of quantum traces and
dimensions as in the previous article. Drinfeld showed that for a
finite-dimensional $\calA$ the following notion works
\cite{Drin} (see also Chari--Pressley \cite{CP} and Majid \cite{Maj}).
\begin{defn}
A Hopf algebra $\calA$ is called quasitriangular if there \index{$R$
$R$--matrix} \index{quasitriangular} exists an invertible element
$R$ (called the universal $R$--matrix) living in a certain
completion $\calA\what\otimes\calA$ of the tensor square
$\calA\otimes\calA$ that satisfies the following axioms
\begin{equation}\label{qtrian}
\begin{aligned}
\sigma\circ\Delta&=R\Delta R^{-1},\\
(\Delta\otimes\id)(R)&=(\sigma\otimes\id)(1\otimes R)\cdot(1\otimes R)\\
(\id\otimes\Delta)(R)&=(\sigma\otimes\id)(1\otimes R)\cdot(R\otimes
1).
\end{aligned}
\end{equation}
Here $\sigma(x\otimes y):=y\otimes x$. If $\calA$ is
finite-dimensional one can take
$\calA\what\otimes\calA=\calA\otimes\calA$.
\end{defn}
When a universal $R$--matrix exists it is essentially unique. Given a
pair of representations $\rho_U,\rho_V\co \calA\to\End(U),\End(V)$ the
universal $R$--matrix induces specializations
$$R_{U,V}\co U\otimes V\to U\otimes V;\quad R_{U,V}(x\otimes y):=(\rho_U\otimes\rho_V)(x\otimes y)\,.$$
and the braiding will be composition with the flip $\times_{U,V}:=\sigma\circ R_{U,V}$.

\begin{remark}
We have not defined the completed tensor product
$\calA\what\otimes\calA$ because we will not use it. We will express
the universal $R$--matrix as an infinite formal sum and it will not
matter what space it lives in. The point is that when it is applied
to any product of reduced tilting modules it will reduce to a finite
sum. Infinite-dimensional Hopf algebras such as $U_q(\Sl_N\C)$,
$U_\eps^{\text{res}}(\Sl_N\C)$ rarely contain a universal
$R$--matrix in $\calA\otimes\calA$.  For our examples one really has
to go out of one's way to construct a big enough extension
$\calA\what\otimes\calA$ to contain $R$. In other words we are
treating the universal $R$--matrix as the Cheshire cat whose smile
lingers on after the cat is gone. After all, tilting modules are
finite-dimensional and all we need are the induced $R$--matrices on
finite-dimensional representations (the smile).
\end{remark}

\begin{aside}{\small
There is a standard way to construct a matrix that satisfies the
conditions from the definition of a quasitriangular Hopf algebra.
The opposite algebra of a Hopf algebra $A$ is the algebra
$A_{\text{op}}$ with multiplication given by
$a\cdot_{\text{op}}b=b\cdot a$. One first constructs an object
living in $A_{\text{op}}\otimes A^*$ satisfying the correct
relations and then takes the image of this object in the correct
completed tensor square. If $A$ is any Hopf algebra let $\{a_i\}$ be
a basis, $\{a^i\}$ be the dual basis and set $R=\sum a_k\otimes
a^k$. We will see that this satisfies the required property in
$A_{\text{op}}\otimes A^*$. Write the comultiplication in $A$ in
index notation as $\Delta(a_k)=\sum \gamma^{ij}_ka_i\otimes a_j$.
This determines the product in $A^*$ by
$$
(a^ia^j)(a_k)=(a^i\otimes a^j)(\Delta(a_k))=\gamma^{ij}_k,
$$
so $a^ia^j=\gamma^{ij}_ka^k$. The product in $A_{\text{op}}$ is
given by $a\cdot_{\text{op}}b=ba$ and the comultiplication is the
same $\Delta$. Now compute
\begin{align*}
(\Delta\otimes\text{id})(R) &= \sum \gamma^{ij}_ka_i\otimes
a_j\otimes a^k \\
&= \sum \delta_{\ell m}^{ij} a_\ell\otimes
a_m\otimes \gamma^{ij}_ka^k \\
&= \sum \delta_{\ell m}^{ij} a_\ell\otimes a_m\otimes a^ia^j \\
&= \sum  a_i\otimes a_j\otimes a^ia^j =(\sigma\otimes\id)(1\otimes
R)\cdot(1\otimes R).
\end{align*}
The proof that the $R$--matrix that we are about to define satisfies
the same equation is that there is a homomorphism taking the $R$--matrix
constructed here to the complicated expression given below. This idea is
the heart of the quantum double construction.  Working out the details
takes from page 128 to 273 of Chari and Pressley \cite{CP}.
Drinfeld won a Fields medal in part for the work formalizing the idea of the
quantum double. The quantized enveloping algebras $U_q(\Sl_N\C)$,
$U_\eps^{\text{res}}(\Sl_N\C)$ are not quantum doubles themselves but they
are quotients of such and the $R$--matrices can be computed by
`projection' \cite{CP}, see also Rosso \cite{Ros}. Roughly speaking,
if $U_q^+$($U_q^-$) denote the subalgebras generated by
$q^{\pm\alpha_i},e_i$($q^{\pm\alpha_i},f_i$) then $U_q^-\simeq
U_q^{+*}$ and $U_q$ is a quotient of $D(U_q^+)\simeq U_q^+\otimes
U_q^-$. Note that $U_q^+$($U_q^-$) are quantum deformations of
algebras of the upper (lower) triangular matrices hence the name
quasitriangular.}
\end{aside}

For $U_q(\Sl_2\C)$ a formal computation following the previous aside
yields \cite{CP,Maj}:
\begin{equation}\label{Rsl2}
R=q^{\alpha\otimes\omega}\sum_{n=0}^\infty q^{n(n-1)/2}
\frac{(q-q^{-1})^{n}}{[n]_q!}e^{n}\otimes f^{n},
\end{equation}
where the second term acts on $x\otimes y$ in the obvious way and
$$
q^{\alpha\otimes\omega}(x\otimes y):=q^{(\alpha,\lambda)(\omega,\mu)}x\otimes y,
$$
when $x,y$ are weight vectors with weights $\lambda,\mu$
respectively. As before $\alpha$ is the simple root of $\Sl_2\C$ and
$\omega={\alpha}/{2}$ is the corresponding fundamental weight. This
formula is to be understood as follows: on any finite-dimensional
representation $e^{n},f^{n}$ act nilpotently and only finite number
of terms in the sum \eqref{Rsl2} are non-zero. Thus the induced
matrices $R_{U,V}$ are well-defined for any pair of
finite-dimensional representations even though the `universal $R$'
itself is just a formal expression.
\begin{example}\label{R3sl2} Recall from \fullref{Uesl2} that the action of $U_q(\Sl_2\C)$ on the representation
$\calV_{(m-1)\omega}^q(\Sl_2\C)$ spanned by $v_1,\ldots,v_{m-1}$ is
given by
\begin{equation}\label{Acsl2}
q^{\pm\alpha}v_i =q^{\pm(m{-}2i)}v_i,\quad
e^nv_i =\frac{[m{-}i{+}n]_q!}{[m{-}i]_q!}v_{i-n},\quad
f^nv_i =\frac{[m{-}i{+}n]_q!}{[i]_q!}v_{i+n}.
\end{equation}
Obviously $e^n$ acts as $0$ on $v_i$ for $n>i$ and $f^n$ acts as $0$
on $v_i$ for $n>m-1-i$. Therefore the sum in \eqref{Rsl2} applied to
$v_i\otimes v_j$ truncates at $\min\{i,m-1-j\}$. By \eqref{Acsl2}
$v_i$ is a weight vector with the weight $(m-2i)\omega$ and
therefore
\begin{multline*}
q^{\alpha\otimes\omega}(v_{i-n}\otimes v_{j+n})=
q^{(\alpha,(m-2i+2n)\omega)(\omega,(m-2j-2n)\omega)}v_{i-n}\otimes v_{j+n}\\
=q^{\frac12(m-2i+2n)(m-2j-2n)}v_{i-n}\otimes v_{j+n}.
\end{multline*}
Substituting this into \eqref{Rsl2} we  explicitly get
\begin{multline}\label{Rijsl2}
R(v_{i}\otimes v_{j})= \\
\sum_{n=0}^{\min\{i,m-1-j\}}
  \hspace{-1em}
  q^{\frac12((m{-}2i{+}2n)(m{-}2j{-}2n){+}n(n{-}1))}
  \frac{(q{-}q^{-1})^{n}}{[n]_q!}
  \frac{[m{-}i{+}n]_q![j{+}n]_q!}{[m{-}i]_q![j]_q!}v_{i-n}{\otimes}v_{j+n}.
\end{multline}
To make this more transparent we will compute this expression for $m=3$ and some pairs $i,j$.
\begin{align*}
R(v_{0}\otimes v_{0}) &=q^{\frac92}v_{0}\otimes v_{0}\\
R(v_{1}\otimes v_{1}) &=q^{\frac12(3-2)(3-2)}v_{1}\otimes v_{1}
+q^{\frac12(3-2+2)(3-2-2)}(q-q^{-1})[3]_q[2]_qv_{0}\otimes v_{2}\\
&=q^{\frac12}v_{1}\otimes v_{1}
+q^{-\frac32}(q-q^{-1})\frac{(q^3-q^{-3})(q^2-q^{-2})}{(q-q^{-1})^2}v_{0}\otimes v_{2}\\
&=q^{\frac12}v_{1}\otimes v_{1}+q^{-\frac32}(q^3-q^{-3})(q+q^{-1})v_{0}\otimes v_{2}\\
R(v_{0}\otimes v_{2}) &=q^{\frac123(3-4)}v_{0}\otimes v_{2}=q^{-\frac32}v_{0}\otimes v_{2}.
\end{align*}
\end{example}
\begin{exm}
Compute $R(v_{i}\otimes v_{j})$ for $m=2$, that is, for
$V=\calV_{\omega}^q(\Sl_2\C)$. Show that in the lexicographic basis
$v_{0}\otimes v_{0},v_{0}\otimes v_{1},v_{1}\otimes
v_{0},v_{1}\otimes v_{1}$ of $V\otimes V$ the $R_{V,V}$ is given by
the matrix
$$\left(\begin{matrix} q^2 & 0 & 0 & 0\\ 0 & 1 & q^2{-}q^{-2} & 0\\
0 & 0 & 1 & 0\\ 0 & 0 & 0 & q^2\end{matrix}\right).$$
\end{exm}
\begin{remark}\label{fracpw}
A reader may notice that expressions in \fullref{R3sl2} involve
fractional powers of $q$ that do not belong to $\C(q)$ that is
technically our field of coefficients. In fact this is a direct
consequence of having the term $q^{\alpha\otimes\omega}$ in
\eqref{Rsl2} that produces $(\omega,\mu)$th power of $q$ when acting
on $u\otimes v$ with $v$ having weight $\mu$. Since
$\omega=\frac{\alpha}{2}$ this number is potentially a half-integer.
In general, the analogous term produces $\frac1{\det(a_{ij})}$
powers of $q$, where $a_{ij}$ is the Cartan matrix. For $\Sl_N\C$ an
elementary computation shows that $\det(a_{ij})=N$ and for the
$R$--matrix formula to make sense we have to extend the coefficients
to $\C(q^{1/N})$. Formally, we now have to work in
$U_{q^{1/N}}(\Sl_N\C):=U_{q}(\Sl_N\C)\otimes_{\C(q)}\C(q^{1/N})$
instead of $U_{q}(\Sl_N\C)$. If we apply the same process to
$U_{q^{1/N}}(\Sl_N\C)$ that was used in \fullref{Uqres} to obtain
$U_\eps^{\text{res}}(\Sl_N\C)$ from $U_{q}(\Sl_N\C)$ we will get the
same algebra assuming that $q^{1/N}$ is specialized to the principal
value of $\smash{\eps^{\frac1N}:=e^{\frac{\pi i}{N(k+N)}}}$. This
allows us to continue doing all computations with the indeterminate
$q$ as explained in \fullref{UeslN}, rewrite the results in terms of
the divided powers and Laurent polynomials and only then specialize
to $\eps$.
\end{remark}
The formula \eqref{Rsl2} can not be specialized to
$U_\eps^{\text{res}}(\Sl_2\C)$ directly because of the $[n]_q!$ in the
denominator. This is easy to fix following our general ideology by
noticing that
$$
\frac{e^{n}\otimes f^{n}}{[n]_q!}=[n]_q!\,e^{(n)}\otimes f^{(n)}.
$$
Thus for $U_\eps^{\text{res}}(\Sl_2\C)$ formula \eqref{Rsl2} is replaced by
\begin{equation}\label{Resl2}
R=q^{\alpha\otimes\omega}\sum_{n=0}^\infty \eps^{n(n-1)/2}
(\eps-\eps^{-1})^{n}[n]_\eps!\,e^{(n)}\otimes f^{(n)}.
\end{equation}
The inconvenience of computing the action of the divided powers for
this formula can be bypassed as follows. Recall that we are
interested in the Weyl modules $\smash{\calW_{(m-1)\omega}^\eps(\Sl_2\C)}$
where $m$ lies in the range $1\leq m<(k{+}N)$ (the Weyl alcove). In this
range $\smash{\frac{e^{n}\otimes f^{n}}{[n]_\eps!}}$ still makes sense and
\eqref{Rsl2} with $q=\eps$ can be used equivalently.

The appearance of fractional powers is not the only nuisance we have
to deal with when moving on to $U_{q^{1/N}}(\Sl_N\C)$. To write the
$R$--matrix for $N>2$ we need analogs of the root vectors $e_i,f_i$
for non-simple positive roots $\alpha$. Unlike the classical case
there is no canonical way to introduce such. A correct
generalization comes from the following classical observation
\cite{CP}. Let $W(=\mathfrak{G}_N)$ be the Weyl group of $\Sl_N\C$
(=the symmetric group) generated by reflections $s_i=s_{\alpha_i}$
in the hyperplanes orthogonal to the simple roots $\alpha_i$. In the
basis $E_{kk}$ these act as the transpositions $s_k=(k\ k+1)$. It is
clear from this description that the Weyl group acts transitively on
the root vectors $e_{ij}:=E_{ii}-E_{jj}$.

More explicitly, each element in the Weyl group admits a
presentation $\sigma=s_{i_1}\ldots s_{i_\nu}$ as a `word' in
generators. A word representing an element is called reduced if it
has the shortest possible length $\nu$  (such a word may not be
unique). Let $w_0\in W$ be the order-reversing permutation that
we already met in connection with duality and
$w_0=s_{i_1}\ldots s_{i_\nu}$ be a reduced word for it. Then
each positive root occurs exactly once in the following sequence
\begin{equation}\label{proots}
\beta_1:=\alpha_{i_1},\,\beta_2:=s_{i_1}(\alpha_{i_2}),\ldots,\beta_\nu:=s_{i_1}\ldots
s_{i_{\nu-1}}(\alpha_{i_\nu}).
\end{equation}
This gives a natural enumeration of the set of all positive roots.
The standard choice of a reduced word for $\Sl_N\C$ is
\begin{equation}\label{rword}
w_0=s_1(s_2s_1)(s_3s_2s_1)\ldots(s_{N-1}\ldots s_2s_1)
\end{equation}
and it gives the anti-lexicographic (read from right to left)
\index{anti-lexicographic root order} enumeration for the roots,
namely
\begin{equation}\label{Eroots}\index{$\beta_k$ ordered positive
roots}
\begin{aligned}
\beta_1&=E_{11}-E_{22}, & \beta_2&=E_{11}-E_{33}, &
\beta_3&=E_{22}-E_{33}, \\
\beta_4&=E_{11}-E_{44},\ldots, & \beta_6&=E_{33}-E_{44},\ldots, \\
\beta_{\nu-N+2}&=E_{11}-E_{NN},\ldots, & \beta_{\nu}&=E_{N-1\,N-1}-E_{NN}.
\end{aligned}
\end{equation}

\begin{exm}
Verify the last claim.
\end{exm}

We now consider how to define standard root vectors in the quantum
setting. For this we need analogs of reflections $s_i$ for quantum
groups. They are given by the Lusztig automorphisms \cite{CP,Lus}.
\begin{defn}\index{Lusztig automorphisms} \index{$T_\sigma$ Lusztig
automorphism} Define algebra automorphisms
$T_i\co U_{q^{1/N}}(\Sl_N\C)\to U_{q^{1/N}}(\Sl_N\C)$ by the following
action on the generators
\begin{equation}\label{Lusaut}
\begin{aligned}
i&=j & |i-j|&>2 & j&=i\pm1\\ \hline
T_iq^{\pm\alpha_i}&=q^{\mp\alpha_i} &
  \quad T_iq^{\pm\alpha_j}&=q^{\pm\alpha_j}, &
  \quad T_iq^{\alpha_j}&=q^{\alpha_j}q^{\alpha_i},
  \hbox{\vrule height12pt depth0pt width0pt}\\[-1ex]
T_ie_i&=-f_iq^{\alpha_i} &
  T_ie_j&=e_j, &
  T_ie_j&=q^{-1}e_je_i-e_ie_j,\\[-0.5ex]
T_if_i&=-q^{\alpha_i}e_i &
  T_if_j&=f_j, &
  T_if_j&=qf_if_j-f_jf_i.
\end{aligned}
\end{equation}
For $\sigma\in W=\mathfrak{G}_N$ represented by the  reduced word
$\sigma=s_{i_1}\ldots s_{i_\nu}$ set
$T_\sigma:=T_{i_1}\ldots T_{i_\nu}$. The operators $T_i=T_{s_{i}}$,
$T_\sigma$ are called the Lusztig automorphisms.
\end{defn}
Note that $T_i$ and $T_\sigma$ are only algebra, not Hopf algebra
automorphisms (they do not preserve the coproduct). For general
permutations $T_\sigma$ are well-defined due to the following result
of G Lusztig \cite{Lus}.
\begin{thm}\label{LusT}
The action of $T_\sigma$ depends only on $\sigma$ and not on the
choice of a reduced expression for it. All $T_\sigma$ are invertible
and $U_{q^{1/N}}^{\text{res}}(\Sl_N\C)$ is invariant under all of them.
Furthermore, the following relations hold
\begin{equation}\label{braid}
\begin{aligned}
T_iT_j &=T_jT_i,\qquad|i-j|>1\\
T_iT_{i+1}T_i &=T_{i+1}T_iT_{i+1}.
\end{aligned}
\end{equation}
\end{thm}
Relations \eqref{braid} are nothing other than the defining
relations for the braid group $\calB_N$ on $N$ strands. In this
group multiplication is concatenation of braids and inverse is the
mirror image. The $t_i$ generator of the braid group  corresponds to
a simple crossing between the $i$th and $(i{+}1)$st strands
\cite{CP,Maj}. \fullref{LusT} implies in particular that
$t_i\mapsto T_i$ defines a representation of the braid group in
$U_{q^{1/N}}(\Sl_N\C)$. Also note that the forgetful map
$\calB_N\to\mathfrak{G}_N$ that only keeps track of the permutation
on the strands is an infinite cover of the Weyl group of $\Sl_N\C$.

We now introduce the quantized versions of the root vectors
corresponding to non-simple roots.
\begin{defn}\label{Rootgen}\index{$E$@$e^{(n)}_{\alpha}$ root vector power}
Let $\alpha\in\Delta^+$ be a positive root and $\alpha=\beta_k$
from \eqref{proots}. Then set
\begin{equation}\label{Rgen}\index{$F$@$f^{(n)}_{\alpha}$ root vector power}
e^{(n)}_{\alpha}:=T_{i_1}\ldots T_{i_{k-1}}(e^{(n)}_{i_k}),\quad
f^{(n)}_{\alpha}:=T_{i_1}\ldots T_{i_{k-1}}(f^{(n)}_{i_k}),
\end{equation}
where $i_1,\ldots,i_k$ are as in the reduced expression
\eqref{rword} for $w_0$. Naturally, we denote
$e_{\alpha}:=e^{(1)}_{\alpha}$, $f_{\alpha}:=f^{(1)}_{\alpha}$.
\end{defn}
\begin{exm}
Show that if $\alpha=\alpha_i$ is simple then
$e^{(n)}_{\alpha}=e^{(n)}_i$, and $f^{(n)}_{\alpha}=f^{(n)}_i$ as
expected.
\end{exm}
With this notation we can now write down formulas for the
$R$--matrices for $U_{q^{1/N}}(\Sl_N\C)$ and
$U_{q^{1/N}}^{\text{res}}(\Sl_N\C)$ that generalize
\eqref{Rsl2} and \eqref{Resl2}:
\begin{equation}
\label{RslN}\index{$R$ $R$--matrix}
\begin{aligned}
R&=q^{\sum_{i=1}^{N-1}\alpha_i\otimes\omega_i}
  \hspace{-1em}\sum_{n_1,\ldots,n_\nu=0}^\infty
  ~\prod_{\alpha\in\Delta^+}^{\gets} q^{\frac{n_k(n_k-1)}2}
  \frac{(q-q^{-1})^{n_k}}
  {[n_k]_q!}\,e^{n_k}_{\alpha}\otimes f^{n_k}_{\alpha},\\
&=q^{\sum_{i=1}^{N-1}\alpha_i\otimes\omega_i}
  \hspace{-1em}\sum_{n_1,\ldots,n_\nu=0}^\infty
  ~\prod_{\alpha\in\Delta^+}^{\gets}
  q^{\frac{n_k(n_k-1)}2}
  (q-q^{-1})^{n_k}[n_k]_q!\,e^{(n_k)}_{\alpha}\otimes f^{(n_k)}_{\alpha}
\end{aligned}
\end{equation}
The product in \eqref{RslN} is not commutative and should be
computed `in reverse' to \eqref{Eroots}, that is,
$\beta_{\nu},\ldots,\beta_1$ so that $\beta_1$ term is applied to
the tensor product first (hence the $\gets$). The $R$--matrix for
$U_\eps^{\text{res}}(\Sl_N\C)$ is obtained from the second line in
\eqref{RslN} by replacing $q$ with $\eps$ \cite{CP}.
\begin{example}\label{Rsl3} In $\Sl_3\C$ the simple roots
are $\alpha_1=E^*_{11}-E^*_{22}$ and $\alpha_2=E^*_{22}-E^*_{33}$ and the
only other positive root is $\alpha_1+\alpha_2=E^*_{11}-E^*_{33}$.
In the ordering of \eqref{Eroots} we have
$\beta_1=\alpha_1,\beta_2=\alpha_1+\alpha_2,\beta_3=\alpha_2$. Hence
$e^{(n)}_{\beta_1}=e^{(n)}_1,e^{(n)}_{\beta_3}=e^{(n)}_2$ and
$e^{(n)}_{\beta_2}=T_1(e^{(n)}_2)$. For $n=1$ we get from
\eqref{Lusaut}
\begin{equation*}
T_1(e_2)=q^{-1}e_2e_1-e_1e_2=-e_1e_2+q^{-1}e_2e_1
\end{equation*}
and because $T_1$ is an automorphism
\begin{equation}\label{T1e2}
T_1(e_2^2)=T_1(e_2)^2=(e_1e_2e_1e_2+q^{-2}e_2e_1e_2e_1)-q^{-1}(e_1e_2^2e_1+e_2e_1^2e_2).
\end{equation}
We now want to rewrite this in terms of the divided powers and
Laurent polynomials in accordance with the general strategy from
\fullref{fracpw}. Taking into account that $q+q^{-1}=[2]_q$ we
have by the quantum Serre relations from \eqref{Uqalg}
$e_2e_1e_2=\frac1{[2]_q}(e_2^2e_1+e_1e_2^2)$ and therefore
\begin{align*}
e_1(e_2e_1e_2)+q^{-2}(e_2e_1e_2)e_1
&=\frac{e_1}{[2]_q}(e_2^2e_1+e_1e_2^2)+q^{-2}(e_2^2e_1+e_1e_2^2)
  \frac{e_1}{[2]_q}\\
&=\frac{e_1^2e_2^2}{[2]_q}+q^{-2}
  \frac{e_2^2e_1^2}{[2]_q}+(1+q^{-2})
  \frac{e_1e_2^2e_1}{[2]_q}\\
&=\frac{e_1^2e_2^2}{[2]_q}+q^{-2}
  \frac{e_2^2e_1^2}{[2]_q}+q^{-1}e_1e_2^2e_1.
\end{align*}
Substituting into \eqref{T1e2} yields
$$T_1(e_2^2)=\frac{e_1^2e_2^2}{[2]_q}+q^{-2}
  \frac{e_2^2e_1^2}{[2]_q}-q^{-1}e_2e_1^2e_2,$$
and dividing both sides by $[2]_q$
\begin{equation}\label{T1e22}
e^{(2)}_{\beta_2}=T_1(e_2^{(2)})=
  e_1^{(2)}e_2^{(2)}-q^{-1}e_2e_1^{(2)}e_2+q^{-2}e_2^{(2)}e_1^{(2)}.
\end{equation}
Analogously,
\begin{equation}\label{T1f22}
f^{(2)}_{\beta_2}=T_1(f_2^{(2)})=
  q^{2}e_1^{(2)}e_2^{(2)}-qf_2f_1^{(2)}f_2+f_2^{(2)}f_1^{(2)}.
\end{equation}
Formulas for $U_\eps^{\text{res}}(\Sl_N\C)$ can now be obtained by
setting $q=\eps$.
\end{example}

\begin{exm}
Prove by induction that
$$e^{(n)}_{\beta_2}=\sum_{k=0}^nq^{-k}e_2^{(k)}e_1^{(n)}e_2^{(n-k)}$$
and derive an analogous formula for $f^{(n)}_{\beta_2}$, see
\cite{Lus}.
\end{exm}

\begin{example}\label{funsl3r}
We will now compute the action of the $R$--matrix of the fundamental
representation in the $\mathfrak{sl}_3\C$ case. Recall from
\fullref{funs3bas} that the fundamental representation
$\calV_{\omega_1}^q(\Sl_3\C)$ has a basis of weight vectors
$v_1:=u_0$, $v_2:=f_1u_0$, $v_3:=f_2v_2$ having weights $\omega_1$,
$\omega_1-\alpha_1$ and $\omega_1-\alpha_1-\alpha_2$ respectively.
This specifies the action of the $q^{\alpha_i}$. The remainder of
the action is specified by $f_iv_j=\delta_{ij}v_{j+1}$ and
$e_iv_j=\delta_{i+1\, j}v_{j-1}$. Continuing with the computation
from the previous example we have
\begin{align*}
e_{\beta_1}&=e_1, & e_{\beta_2}&=-e_1e_2+q^{-1}e_2e_1, & e_{\beta_3}&=e_2,\\
f_{\beta_1}&=f_1, & f_{\beta_2}&=-f_2f_1+qf_1f_2, & f_{\beta_3}&=f_2.
\end{align*}
Since $v_1$ is annihilated by $e_1$ and $e_2$, only the $q^{\sum
\alpha_i\otimes\omega_i}$ factor of the $R$ matrix acts on vectors
of the form $v_1\otimes x$ (see equation \eqref{RslN}). Since the
weight of the fundamental representation is
$\omega_1=\tableau{1}=\frac23E_{11}^*-\frac13E_{22}^*-\frac13E_{33}^*$,
it follows that
\begin{align*}
R(v_1\otimes v_1)
  &=q^{\langle\alpha_1,\omega_1\rangle\langle\omega_1,\omega_1\rangle}
  v_1\otimes v_1=q^{2/3}v_1\otimes v_1, \\
R(v_1\otimes v_2)
  &=q^{\langle\alpha_1,\omega_1\rangle\langle\omega_1,\omega_1-\alpha_1\rangle}
  v_1\otimes v_2=q^{-1/3}v_1\otimes v_2, \\
R(v_1\otimes v_3)
  &=q^{\langle\alpha_1,\omega_1\rangle\langle\omega_1,\omega_1-\alpha_1\rangle}
  v_1\otimes v_3=q^{-1/3}v_1\otimes v_3.
\end{align*}
Now
$$e_{\beta_1}v_2=v_1, \qquad e_{\beta_2}v_2=0, \qquad e_{\beta_3}v_2=0,$$
so the only terms that contribute to the $R$--matrix evaluated on
vectors of the form $v_2\otimes x$ will correspond to the root
$\beta_1$. We compute
$$f_{\beta_1}v_1=v_2, \qquad f_{\beta_1}v_2=0, \qquad e_{\beta_1}v_3=0.$$
It follows that
\begin{align*}
R(v_2\otimes v_1)
  &=q^{\sum\alpha_i\otimes\omega_i}(v_2\otimes v_1+(q-q^{-1})v_1\otimes v_2)\\
  &=q^{\langle\alpha_1,\omega_1-\alpha_1\rangle\langle\omega_1,
  \omega_1\rangle+\langle\alpha_2,\omega_1-\alpha_1\rangle\langle\omega_2,
  \omega_1\rangle}v_2\otimes v_1 \\
&\qquad+(q-q^{-1})q^{\langle\alpha_1,\omega_1\rangle\langle\omega_1,
  \omega_1-\alpha_1\rangle+\langle\alpha_2,\omega_1\rangle\langle\omega_2,
  \omega_1-\alpha_1\rangle}v_1\otimes v_2, \\
&=q^{-1/3}v_2\otimes v_1+q^{-1/3}(q-q^{-1})v_1\otimes v_2,\\
R(v_2\otimes v_2)
  &=q^{\langle\alpha_1,\omega_1-\alpha_1\rangle\langle\omega_1,
  \omega_1-\alpha_1\rangle+\langle\alpha_2,
  \omega_1-\alpha_1\rangle\langle\omega_2,\omega_1-\alpha_1\rangle}
  v_2\otimes v_2\\
&=q^{2/3}v_2\otimes v_2, \\
R(v_2\otimes v_3)
&=q^{\langle\alpha_1,\omega_1-\alpha_1\rangle\langle\omega_1,
  \omega_1-\alpha_1-\alpha_2\rangle+\langle\alpha_2,
  \omega_1-\alpha_1\rangle\langle\omega_2,\omega_1-\alpha_1-\alpha_2\rangle}
  v_2\otimes v_3\\&=q^{-1/3}v_2\otimes v_3.
\end{align*}
Now
\begin{align*}
e_{\beta_1}v_3&=0, & e_{\beta_2}v_3&=-v_1, & e_{\beta_3}v_3&=v_2, \\
f_{\beta_2}v_1&=-v_3, & f_{\beta_2}v_2&=0, & f_{\beta_2}v_3&=0, \\
f_{\beta_3}v_1&=0, & f_{\beta_3}v_2&=v_3, & f_{\beta_3}v_3&=0.
\end{align*}
It follows that
\begin{align*}
R(v_3{\otimes}v_1)
  &=q^{\sum\alpha_i{\otimes}\omega_i}(v_3{\otimes}v_1{+}(q{-}q^{-1})v_1{\otimes}v_3)\\
&=q^{\langle\alpha_1,\omega_1{-}\alpha_1{-}\alpha_2\rangle\langle\omega_1,
  \omega_1\rangle{+}\langle\alpha_2,\omega_1{-}\alpha_1
  {-}\alpha_2\rangle\langle\omega_2,\omega_1\rangle}v_3{\otimes}v_1 \\
  &\quad+(q{-}q^{-1})q^{\langle\alpha_1,\omega_1\rangle\langle\omega_1,
  \omega_1{-}\alpha_1{-}\alpha_2 \rangle}v_1{\otimes}v_3, \\
&=q^{-1/3}v_3{\otimes}v_1{+}q^{-1/3}(q{-}q^{-1})v_1{\otimes}v_3,\\
R(v_3{\otimes}v_2)
  &=q^{\sum\alpha_i{\otimes}\omega_i}(v_3{\otimes}v_2{+}(q{-}q^{-1})v_2{\otimes}v_3)\\
  &=q^{\langle\alpha_1,\omega_1{-}\alpha_1{-}\alpha_2\rangle\langle\omega_1,
  \omega_1{-}\alpha_1\rangle{+}\langle\alpha_2,\omega_1{-}\alpha_1{-}\alpha_2\rangle
  \langle\omega_2,\omega_1{-}\alpha_1\rangle} v_3{\otimes}v_2 \\
&\quad+(q{-}q^{-1})q^{\langle\alpha_1,\omega_1{-}\alpha_1\rangle\langle\omega_1,
  \omega_1{-}\alpha_1{-}\alpha_2 \rangle{+}\langle\alpha_2,
  \omega_1{-}\alpha_1\rangle\langle\omega_2,
  \omega_1{-}\alpha_1{-}\alpha_2\rangle}v_2{\otimes}v_3, \\
&=q^{-1/3}v_3{\otimes}v_2{+}q^{-1/3}(q{-}q^{-1})v_2{\otimes}v_3,\\
R(v_3{\otimes}v_3)
  &=q^{\langle\alpha_1,\omega_1{-}\alpha_1{-}\alpha_2\rangle\langle\omega_1,
  \omega_1{-}\alpha_1{-}\alpha_2\rangle{+}\langle\alpha_2,
  \omega_1{-}\alpha_1{-}\alpha_2\rangle\langle\omega_2,
  \omega_1{-}\alpha_1{-}\alpha_2\rangle}v_3{\otimes}v_3\\
&=q^{2/3}v_3{\otimes}v_3.
\end{align*}
Ordering the basis for $\calV_{\omega_1}^q(\Sl_3\C)\otimes
\calV_{\omega_1}^q(\Sl_3\C)$ lexicographically we may write the
matrix for $R$ as
$$R=q^{-1/3}\left[\begin{matrix} q&&&&&&&&\\
  &1&&(q{-}q^{-1})&&&&&\\ &&1&&&&(q{-}q^{-1})&&\\ &&&1&&&&&\\ &&&&q&&&&\\
  &&&&&1&&(q{-}q^{-1})&\\ &&&&&&1&&\\ &&&&&&&1&\\ \\ &&&&&&&&q
  \end{matrix}\right].$$
\end{example}
This last example can be generalized to the fundamental
representation $\calV_{\omega_1}^q(\Sl_N\C)$ for all $N$.
\begin{exm}\label{NfunR}
In general the fundamental representation has a basis with highest
weight vector $v_1$ and $v_{j+1}:=f_jv_j$. Show that in this basis
we have, $\omega_1(e_j)=E_{j\, j+1}$, $\omega_1(f_j)=E_{j+1\, j}$,
$\omega_1(q^{\alpha_i^\vee})=(q-1)E_{ii}+(q^{-1}-1)E_{i+1\,i+1}+I$,
and
$$R=q^{-1/N}\biggl(q\sum E_{ii}\otimes E_{ii}+\sum_{i\ne
j}E_{ii}\otimes E_{jj}+(q-q^{-1})\sum_{i<j}E_{ij}\otimes
E_{ji}\biggr).$$
\end{exm}
The answer to this exercise can be found on page 277 of \cite{CP}.
It is important because it provides the link between the quantum
invariants as we are defining them, and the THOMFLYP polynomial and
skein theory as sketched in \fullref{app:c}.

It follows from \fullref{126} that the representation
$$\calV^q_{\tableau{4 2}}(\Sl_4\C)$$
for $U_q({\mathfrak{sl}}_4{\mathbb C})$ has $126$ nontrivial weight
spaces, we know that the dimension of
$$\calV^q_{\tableau{4 2}}(\Sl_4\C)$$
is  $126$. It follows that the $R$--matrix is a
$15876\times 15876$ matrix. As the reader can see, computing the
$R$--matrices explicitly is no mean feat even for simple
representations. The next exercise is one of the few remaining cases
that can be worked out reasonably by hand.
\begin{exm}
Compute a basis for the representation
$$\calV_{\tableau{1 1}}^q(\Sl_3\C).$$
Compute the action of $q^{\alpha_i}$, $e_i^{(n)}$,
$f_i^{(n)}$, and the $R$--matrix.
\end{exm}

Recall that when $U_q(\Sl_N\C)$ acts on a vector space $V$, we can
define an action on the dual $V^*$ by
$$(ag)(v):=g(\gamma(a)v).$$
To see how this works, consider the action of $U_q(\Sl_2\C)$ on the
dual of the representation
$$\calV_{\tableau{1}}^q(\Sl_2\C).$$

\begin{example}\label{dualeg}
Recall that the representation
$$\calV_{\tableau{1 1}}^q(\Sl_2\C)$$
has a basis $v_1$, $v_2$. We denote the dual basis by $v^1$, $v^2$.
Since $(ev^1)(v_1){=}v^1(\gamma(e)v_1)=0$ and
$(ev^1)(v_2)=v^1(\gamma(e)v_2)=-q$, we have $ev^1=-qv^2$. Similarly,
$ev^2=0$, $fv^1=0$ and $fv^2=-q^{-1}v^1$. We can now see how the
$R$--matrix acts on terms that include functionals such as
$\omega_1\otimes\omega_1^*$ in the ${\mathfrak{sl}}_2$ case. We have
\begin{align*}
R(v_1\otimes v^1)&=q^{-1/2}v_1\otimes v^1 \\
R(v_1\otimes v^2)&=q^{1/2}v_1\otimes v^2 \\
R(v_2\otimes v^1)&=q^{1/2}v_2\otimes v^1 \\
R(v_2\otimes v^2)&=q^{-1/2}v_2\otimes v^2
-q^{-3/2}(q-q^{-1})v_1\otimes v^1.
\end{align*}
\end{example}

Recall that we can define link and framed link invariants from a
ribbon category. The $R$--matrix is the main tool to build a ribbon
category from the representations of a quantum group. In fact the
category of type I representations of $U_{q^{1/N}}(\Sl_N\C)$ becomes
a ribbon category if we use the braiding given by
$\times_{u,V}=\sigma\circ R_{U,V}$. We denote the (framed) link
invariants resulting from this category by $W^{{\mathfrak
sl}_N}_\Lambda$ \index{$W^{{\mathfrak{sl}}_N}_\Lambda$ quantum
framed link invariant} where $\Lambda$ is a collection of
representations (one for each component of the link). We next
compute the invariant
$W^{{\mathfrak{sl}}_2}_{\tableau{1},\tableau{1}}(\text{left Hopf
link})$ from the definition. First compute the composition of the
top three morphisms from \fullref{lhopf}. We have
\begin{align*}
v_1\otimes v^1 &\mapsto q^{3/2}v_1\otimes v^1\mapsto qv^1\otimes
v_1\mapsto q \\
v_2\otimes v^2 &\mapsto q^{3/2}v_2\otimes v^2\mapsto
q^{3/2}[q^{-1/2}v^2\otimes v_2-q^{-3/2}(q-q^{-1})v^1\otimes
v_1]\mapsto q^{-1}.
\end{align*}
We also have $v_1\otimes v^2\mapsto 0$ and $v_2\otimes v^1\mapsto
0$. Notice that this recreates the double dual isomorphism.  To
continue the computation of the invariant for the left Hopf link it
is helpful to write out the matrix for the braiding in the
lexicographic basis  for $\calV^q_{\tableau{1}}(\Sl_2\C)\otimes
\calV^q_{\tableau{1}}(\Sl_2\C)$ from \fullref{NfunR}. This is our
first sample computation of a braiding.
\begin{example}\label{breg}
We have
$$\times_{\tableau{1},\tableau{1}}{=}q^{-1/2}\!\!\left[\begin{matrix}
q&&&\\&0&1&\\&1&(q{-}q^{-1})& \\&&&q\end{matrix}\right], \quad
\times^{-2}_{\tableau{1},\tableau{1}}{=}q\!\!\left[\begin{matrix}
q^{-2}&&&\\&(q^{-1}{-}q)^2{+}1&(q^{-1}{-}q)&\\&(q^{-1}{-}q)&1&
\\&&&q^{-2}\end{matrix}\right].$$
\end{example}
We can combine our computations of the twist, double dual pairing
and braiding to compute the invariant of the Hopf link. We use
$v_i^k$ to denote $v_i\otimes v^k$ and $v_{ij}^{k\ell}$ to denote
$v_i\otimes v_j\otimes v^k\otimes v^\ell$. The computation is
contained in the following example.
\begin{example}\label{Wlhopfeg}
Following the element $1$ up through the morphisms as depicted in
\fullref{lhopf} gives
\begin{align*}
1 &\mapsto v_1^1+v_2^2  \\
&\mapsto  v_{11}^{11}+v_{12}^{21}+v_{21}^{12}+v_{22}^{22} \\
&\mapsto  q^{-1}v_{11}^{11}+(q(q-q^{-1})^2+q)v_{12}^{21}
  +(1-q^2)v_{21}^{21}
+(1-q^2)v_{12}^{12}+qv_{21}^{12} +q^{-1}v_{22}^{22} \\
&\mapsto v_1^1+((q-q^{-1})^2+1)v_1^1+q^2v_2^2+q^{-2}v_2^2 \\
&\mapsto q+q(q-q^{-1})^2+q+q+q^{-3} = q^3+q+q^{-1}+q^{-3}.
\end{align*}
We conclude that
$$W^{{\mathfrak{sl}}_2}_{\tableau{1},\tableau{1}}(\text{left Hopf
link})=q^3+q+q^{-1}+q^{-3}.$$
\end{example}

We can use this type of computation to derive an interesting
recurrence relation for the framed invariants. Notice that the twist
in the fundamental representation is given by
\begin{align*}
\times_{\tableau{1},\tableau{1}}(v_{ii})&=q^{-1/N}qv_{ii}\\
\times_{\tableau{1},\tableau{1}}(v_{ij})&=q^{-1/N}v_{ji}, \quad i<j\\
\times_{\tableau{1},\tableau{1}}(v_{k\ell})&=q^{-1/N}(v_{\ell
k}+(q-q^{-1})v_{k\ell}), \quad k>\ell.
\end{align*}
It follows that
\begin{align*}
\times^{-1}_{\tableau{1},\tableau{1}}(v_{ii})&=q^{1/N}q^{-1}v_{ii}\\
\times^{-1}_{\tableau{1},\tableau{1}}(v_{ji})&=q^{1/N}v_{ij}, \quad i<j\\
\times^{-1}_{\tableau{1},\tableau{1}}(v_{\ell
k})&=q^{1/N}(v_{k\ell}-(q-q^{-1})v_{\ell k}), \quad k>\ell.
\end{align*}

\begin{exm}\label{preskein}
Use the previous computations to show that
$$
q^{1/N}\times_{\tableau{1},\tableau{1}}-q^{-1/N}\times_{\tableau{1},\tableau{1}}^{-1}=
(q-q^{-1})\text{id}_{\tableau{1},\tableau{1}}\,.
$$
\end{exm}
Notice that this implies that the framed link invariant
$W^{{\mathfrak{sl}}_N}$ (all representations taken to by
fundamental) satisfies a skein relation. Also notice that the
invariant defined via quantum groups is provably invariant under
isotopies. A framed link invariant is very close to being a link
invariant; it is invariant under Reidemeister moves two and three
from \fullref{reidemeister}. It follows from \fullref{fanom} below
that
$$W^{{\mathfrak{sl}}_N}(\text{closure}(\theta_|\circ f))=
q^{N-1/N}W^{{\mathfrak{sl}}_N}(\text{closure}(f))\,.$$ Here closure
is just the quantum trace in the category of framed tangles; see
\fullref{qtrace}. This motivates the following definition of the
THOMFLYP polynomial.
\begin{defn}\label{thomflypdef}
The THOMFLYP polynomial is defined by \index{$P(L)$ numerical
THOMFLYP polynomial}\index{THOMFLYP polynomial}
$$P(L):=q^{(1/N-N)\sum_{i=1}^{c(L)}\sum_{j=1}^{c(L)}n_{ij}(L)}W^{{\mathfrak{sl}}_N}(L)\,.$$
Here $c(L)$ is the number of components of $L$ and $n_{ij}(L)$ is
the linking matrix of $L$.
\end{defn}
\begin{remark}
To define $3$--manifold invariants we have to go to the category or
reduced tilting modules $\Tilt(\Sl_N\C)$. Of course this is also a
ribbon category that generates link invariants. The resulting
invariants are just the evaluation of $P(L)$ at $q^{1/N}=e^{\pi
i/(N(k+N))}$ In this way $P(L)$ is a function of $N$ and a primitive
even root of unity $\eps=e^{\pi i/(k+N)}$. One can show \cite{CP}
that for each $L$ there is a unique rational function in variables
$\varlambda^{1/2}$, $\mathfrak{q}^{1/2}$ denoted by $\mathcal{P}(L)$
\index{$P$@$\mathcal{P}(L)$ THOMFLYP polynomial} such that
$P(L)=\mathcal{P}(L)(\varlambda=\eps^{2N},\mathfrak{q}=\eps^2)$.
\end{remark}
\begin{exm}
From \fullref{unknot} we know that
$$\mathcal{P}(\text{unknot})=\frac{\varlambda^{1/2}-\varlambda^{-1/2}}{\mathfrak{q}^{1/2}-\mathfrak{q}^{-1/2}}$$
Prove that the THOMFLYP polynomial satisfies the following skein
relation
$$\varlambda^{1/2}\mathcal{P}(\text{closure}(\times_{|,|}\circ f))-
\varlambda^{-1/2}\mathcal{P}(\text{closure}(\times_{|,|}^{-1}\circ
f))=
(\mathfrak{q}^{1/2}-\mathfrak{q}^{-1/2})\mathcal{P}(\text{closure}(f)).$$
\end{exm}
This (skein) recurrence relation together with the value of the
unknot specifies the link invariant uniquely. In fact one can turn
the entire process around and use the skein relation and
normalization as the definition of the link invariants. It is not
obvious that such a definition is well-formed. One must show that
different ways of applying the skein relation lead to the same
answers and one must show the resulting quantity is invariant under
the Reidemeister moves.  In addition to the approach to  the
definition that we used via quantum groups there are several
different approaches to do this. The first are described in the
original paper \cite{HOMFLY}.

Much more is true -- the full Reshetikhin--Turaev invariants of
links in general $3$--manifolds can be recovered from the skein
theory. This is described in \fullref{app:c}.

\begin{exm}
Compute $\mathcal{P}(\text{Right Hopf link})$ from the definition
and via the skein relation and reconcile the two answers.
\end{exm}

\section{Reshetikhin--Turaev invariants from quantum groups}

We are now in a position to combine all of the facts we derived
about representations of quantum groups into a definition
of a nontrivial strict modular category. This is the category of reduced tilting
modules
at even roots of unity based on the Lie algebra $\Sl_N\C$. Historically,
Reshetikhin and Turaev constructed their invariants directly \cite{RT1,RT2}
and only later
Turaev \cite{T} streamlined their construction into a notion of modular category.
In this section we will finally define these invariants. First we review the main
features of
$\Tilt(\Sl_N\C)$, then compute its
characteristic numbers, the $\wtilde{s}$--matrix and finally the
Reshetikhin--Turaev invariant of $S^3$, aka the Chern--Simons partition function.

\subsection{Modular category of reduced tilting modules at even roots of
unity}\label{tilt}

In this article we summarize the modular category structure on $\Tilt(\Sl_N\C)$ when $\eps$
is a primitive {\it even} root of unity.  This is the modular
category that produces the Reshetikhin--Turaev (quantum) invariants
of $3$--manifolds along the lines of \fullref{modinv}. After
discussing the main ingredients we compute some of the
characteristic numbers of this category (see \fullref{charnum}). It
is assumed throughout that $\smash{\eps=e^{\frac{2\pi i}l}}$, where
$l$ is {\it even}. To draw a connection with the Chern--Simons
theory set $k:=l/2-N$ for the level of the theory. Then $l=2(k+N)$,
$l^\prime=k+N$ and $\smash{\eps=e^{\frac{\pi i}{k+N}}}$. One should
keep in mind that there is no mathematical definition of Witten's
path integral, hence no rigorous way to compare Witten's invariants
to the Reshetikhin--Turaev ones and the above connection is merely
conjectural.

\begin{description}
\item[{Objects}] The objects are reduced tilting modules $V$, that is,
representations of  $U_\eps^{\text{res}}(\Sl_N\C)$ (\fullref{Uqres})
that are finite direct sums of Weyl modules (\fullref{Wmod}) with
the highest weights in the Weyl alcove $\Lambda_w^l$ (see
\fullref{Walcove} or \textbf{Index set} below):
$V=\oplus_{\lambda\in\Lambda_w^l}\smash{\bigl(\calW_{\lambda}^\eps\bigr)^{\oplus
m_{\lambda}(V)}}$.
\item[{Morphisms}] The morphisms are equivariant linear maps $f(av)=a(f(v))$.
\item[{Unit object}] The unit object is $\mathds{1}:=\C$ with the trivial
representation structure $az:=\varepsilon(a)z$, where $\varepsilon$
is the counit \eqref{Uqcoalg} of $U_\eps^{\text{res}}(\Sl_N\C)$.
\item[{Dual objects}] The dual objects are the vector space duals $V^*$ of tilting
modules with the action given by the antipode \eqref{Uqcoalg}
$(a\varphi)(v):=\varphi(\gamma(a)z)$.
\item[{Pairing and copairing}] The pairing and copairing are the usual ones from
the category of finite-dimensional vector spaces
$$\cap_V(\varphi\otimes v):=\varphi(v),\qquad
  \cup_V(1):=\sum_k\,v_k\otimes v^k,$$
where $v_k,v^k$ are dual bases in $V,V^*$ respectively.
\item[{Tensor product}] The tensor product is the reduced tensor product $\rotimes$
from \fullref{dualtensor}, that is, the maximal invariant
subspace $U\rotimes V$ of the ordinary tensor product $U\otimes V$
that is a direct sum of Weyl modules with highest weights in the
Weyl alcove. The action restricts from the action on the usual tensor product
$a(u\otimes v):=\Delta(a)(u\otimes v)$ on $U\otimes V$.
\item[{Braiding}] The braiding is $\times_{U,V}:=\sigma\circ R_{U,V}$, where $\sigma(u\otimes v):=v\otimes u$
and $R_{U,V}$ is the restriction of the `universal $R$--matrix' given
by \eqref{RslN} to the ordinary tensor product of $U,V$.
\item[{Twist}] The twist is obtained from the braiding, duality and double duality as in \fullref{twistfromdual},
explicitly
$$
\theta_V
:=(\text{id}_V\otimes\cap^*_V)\circ(\times_{V,V}\otimes\text{id}_{V^*})\circ(\text{id}_V\otimes\cup_V)\,,
$$
where $\cap_V^*(v\otimes\varphi):=\varphi(q^{2\rho}v)$.
\item[{Index set}] The index set of simple objects $I$ is the set of dominant
weights (see \fullref{app:d}) in the open Weyl alcove
$$
I:=\Lambda_w^l=\{\lambda\in\Lambda_w^+\mid\,(\lambda+\rho,\alpha)<l/2,\
\text{for all}\ \alpha\in\Delta^+\}.
$$
The unit object is indexed by $\lambda=0$.
\item[{Index involution}] The index involution is given by $\lambda^*:=-w_0(\lambda)$
(see \fullref{emalcove}), where $w_0$ is the order-reversing
permutation and the action of the symmetric group $W=\mathfrak{G}_N$
on weights is described in \fullref{app:d}. In particular, $0^*=0$.
The dual weight is the weight of the dual simple object:
$\calW_{\lambda}^{\eps}(\Sl_N\C)^*=\calW_{\lambda^*}^\eps(\Sl_N\C)$.
\item[{Simple objects}] The simple objects are the Weyl modules $\calW_{\lambda}^\eps(\Sl_N\C)$
with highest weights in the Weyl alcove $\lambda\in\Lambda_w^l$
(alcove Weyl modules).
\end{description}

For convenience we list the salient points of the definition in
\fullref{cint}.\hfill\newline
\begin{table}[ht!]
\label{cint}
\begin{tabular}{|l|l|}
\hline
Category \hbox{\vrule height12pt depth0pt width0pt}
  &  $\Tilt(\Sl_N\C)$, $\eps:=e^{{2\pi i}/l}$, $l$ even\\
Objects & Reduced tilting modules
$V=\oplus_{\lambda\in\Lambda_w^l}\left(\calW_{\lambda}^\eps\right)^{\oplus m_{\lambda}(V)}$\\
Morphisms & Equivariant linear maps $f(au)=af(u)$\\
Unit object & $\C$ with $az:=\varepsilon(a)z$\\
Dual object & Dual vector space with $(a\varphi)(v):=\varphi(\gamma(a)z)$\\
Tensor product & Reduced product with $a(u\otimes v):=\Delta(a)(u\otimes v)$\\
Pairing,copairing & $\cap_V(\varphi\otimes v):=\varphi(v),\cup_V(1):=\sum_k\,v_k\otimes v^k$\\
Braiding & $\times_{U,V}(u\otimes v):=\sigma(R(u\otimes v))$\\
Twist & $\theta_V:=(\text{id}_V\otimes\cap^*_V)\circ(\times_{V,V}\otimes\text{id}_{V^*})\circ(\text{id}_V\otimes\cup_V)$\\
Index set & Weyl alcove $I=\Lambda_w^l$\\
Index involution & $\lambda^*:=-w_0(\lambda)$\\
Simple objects &  Alcove Weyl modules $\{\calW_{\lambda}^\eps\}_{\lambda\in\Lambda_w^l}$\\
 \hline
\end{tabular}
\vskip.1in
\caption{Modular category of reduced tilting modules at even roots of unity}
\end{table}
\index{$T$@$\Tilt(\Sl_N\C)$ reduced tilting modules} \index{$\times$
braiding}\index{$\theta$ twist} \index{$(*,\cup,\cap)$ duality
triple} \index{${\bbone}$ unit}

As we described in \fullref{color} one can define a dual pairing and copairing in any modular category.
In $\Tilt(\Sl_N\C)$ we have, $\cap_V^*(v\otimes\varphi):=\varphi(q^{2\rho}v)$ and
$\cup_V^*(1):=\sum_k\,v^k\otimes q^{-2\rho}v_k$, where $\rho$ is the
Weyl weight of $\Sl_N\C$ (see \fullref{app:d}).

This is not yet a strict category because equality signs in the axioms
are only canonical isomorphisms. For example, $V$ is not equal to
$V\otimes \C$. In truth, we should be talking about isomorphism
classes of tilting modules rather than tilting modules themselves, or
selecting canonical representatives of those classes.
One faces the same nuance in the category of representations of classical Lie groups.
Formally, some tweaking in the notions of morphisms and ribbon
operations is required. It can be done in a standard way by the Mac\,Lane
coherence theorem \cite{mac} but in practice one can safely ignore the difference.

\begin{figure}[ht!]
\centering
\labellist\small
\pinlabel {$\delta$} at 23 410
\pinlabel {$\delta$} at 221 410
\pinlabel {$\delta$} at 244 410
\pinlabel {$\delta$} at 480 410
\pinlabel {$\theta$} at 480 344
\pinlabel {$\theta$} at 32 64
\pinlabel {$\theta$} at 61 64
\pinlabel {$\theta$} at 171 120
\pinlabel {$\theta$} at 171 47
\pinlabel {$\delta$} at 341 141
\pinlabel {$\theta$} at 329 24
\pinlabel {$\delta$} at 514 142
\pinlabel {$\theta$} at 503 35
\endlabellist
\includegraphics[width=4truein]{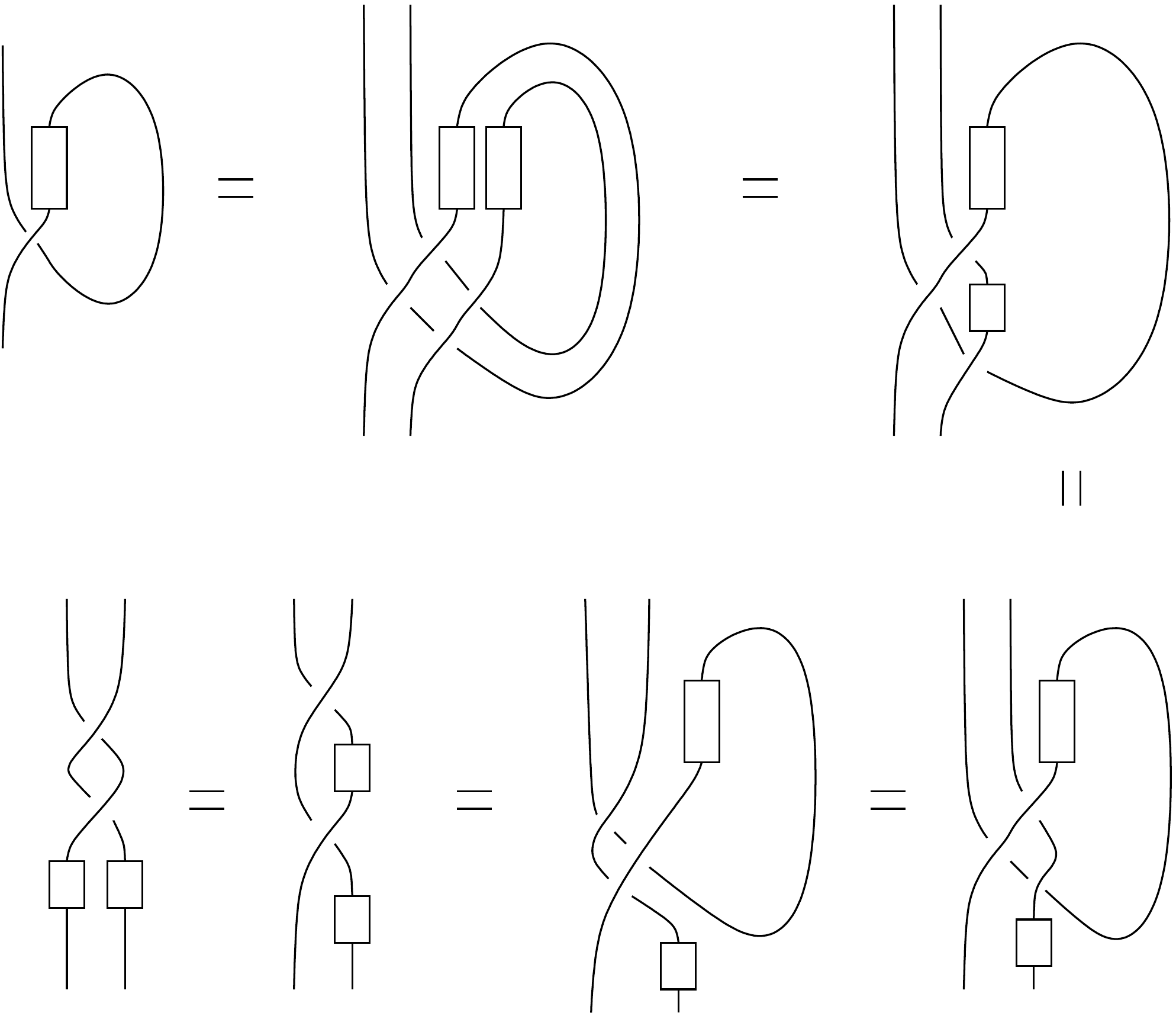} \caption{Proof of
axiom 8}\label{ax8}
\end{figure}
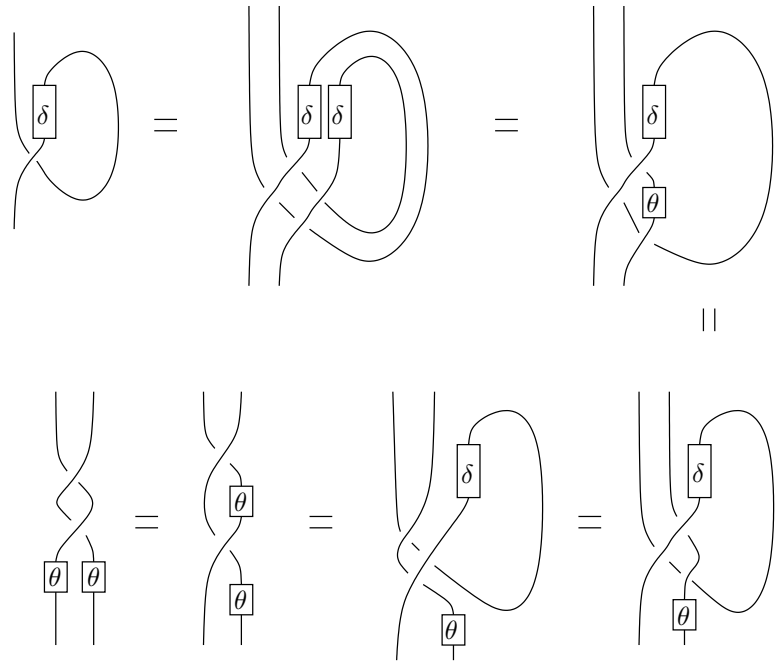

To see that the ingredients do indeed form a strict modular
category, we will sketch proofs of some of the axioms \eqref{axmtc}. The tensor
axioms (Axioms 1--4) essentially follow from \fullref{prodok}.
The braiding axioms (Axioms 5--7) follow from the defining
properties \eqref{qtrian} of the $R$--matrix. There is a subtlety
here since \eqref{RslN} is only a formal expression but
\eqref{qtrian} holds if $R$ is restricted to a pair of
finite-dimensional representations \cite{CP}.
\begin{exm}
If $R=\sum a\otimes b$ write $R_{13}=\sum a\otimes 1\otimes b$ and
$R_{23}=\sum 1\otimes a\otimes b$. Compute the left hand side and
the right hand side of Axiom 6. Conclude that Axiom 6
follows from
$$
(\Delta\otimes\text{id})(R)=R_{13}R_{23}\,.
$$
\end{exm}
To verify Axiom 7 write $R=\sum a\otimes b$ and compute
\begin{multline*}
(g\otimes f)\circ \times_{V,W}(s\otimes t)
  =(g\otimes f)\Bigl(\sum (bt)\otimes (as)\Bigr)\\
  =\sum (bg(t))\otimes(af(s))=
  \times_{V^\prime,W^\prime}\circ(f\otimes g)(s\otimes t).
\end{multline*}
The twist axioms and the last duality axiom (Axioms 8, 9 and 12) can
be verified by manipulating tangle diagrams. One uses the graphical
definition of the twist from \fullref{twistfromdual}. For example
Axiom 8 is verified in \fullref{ax8}. The remaining duality axioms
(Axioms 10 and 11) are trivial computations. Axioms 13--16 follow from Andersen's
\fullref{alcove} and the construction of $\Tilt(\Sl_N\C)$. See Andersen \cite{And},
Andersen and Paradowski \cite{AP}, Chari and Pressley \cite{CP} and Sawin \cite{Saw} for more details.
The only remaining axiom is 17, the non-degeneracy of the
$\wtilde{s}$--matrix. It will be proved as a byproduct of computing the Chern--Simons
partition function of $S^3$ in \fullref{cscompute}.

As we already mentioned the link between the Reshetikhin--Turaev and
Witten--Chern--Simons invariants is only conjectural. In the absence
of a formal definition for the latter we can simply identify them
with the former for the category $\Tilt(\Sl_N\C)$.
\begin{defn}\label{slntaudef}
The Reshetikhin--Turaev invariant corresponding to $U_q(\Sl_N\C)$ at
level $k$ is defined to be the invariant $\tau$ from definition
\eqref{taudef} arising from the category $\Tilt(\Sl_N\C)$. It is
denoted by $\tau^{\Sl_N\C}(M)$. The Chern--Simons partition function
is a different name for the same thing. It is denoted by
$Z^{CS}(M):=\tau^{\Sl_N\C}(M)$. \index{Chern--Simons partition
function}\index{$Z^{CS}(M):=\tau^{\Sl_N\C}(M)$ Chern--Simons
partition function}\index{$\tau^{\Sl_N\C}(M)$ $\Sl_N\C$
Reshetikhin--Turaev invariant}\index{$\Sl_N\C$ Reshetikhin--Turaev
invariant}
\end{defn}

As a warm-up to the computation of $Z^{CS}(S^3)$ we compute the
characteristic numbers (\fullref{charnum}) of $\Tilt(\Sl_N\C)$ here.
We already computed the number $d_\lambda$ in \eqref{Wdim}, but we
would like to simplify it. Recall how this went. The morphism
associated to an unknot labeled with $\lambda$ sends $1$ to
$d_\lambda$. Let $\{v_i\}$ and $\{v^i\}$ be bases for the
representation $\calV^q_\lambda$ and its dual respectively. The
morphism is a composition of two morphisms and we see that
$$1\mapsto \sum v_i\otimes v^i\mapsto \sum
  v^i(\lambda(q^{2\rho})v_i)=\text{Tr}(\lambda(q^{2\rho}))
  =\chi_\lambda(q^{2\rho}).$$
It follows that $d_\lambda=\chi_\lambda(q^{2\rho})$. The Weyl
character formula states that the characters are just Schur
polynomials, which in turn are ratios of determinants. In this case
the determinants are Vandermonde determinants. The appearance of
group characters was explained in \fullref{someinterest}.

Since the simple objects are the Weyl modules and the Weyl modules
are specializations of $\calV^q_\lambda$ we have,
\begin{equation*}
d_{\lambda}=\dim_q\!\calW_\lambda^{\eps}=\prod_{\alpha\in\Delta^+}
\frac{\eps^{(\lambda+\rho,\alpha)}-\eps^{-(\lambda+\rho,\alpha)}}{\eps^{(\rho,\alpha)}-\eps^{-(\rho,\alpha)}}.
\end{equation*}
Substituting $\eps=e^{\frac{\pi i}{k+N}}$ yields
$$\eps^{(\mu,\alpha)}-\eps^{-(\mu,\alpha)}
  =2i\,\frac{e^{\frac{\pi i}{k+N}(\mu,\alpha)}
  -e^{-\frac{\pi i}{k+N}(\mu,\alpha)}} {2i}
  =2i\,\sin\left(\frac{\pi(\mu,\alpha)}{k+N}\right)$$
for any weight $\mu$. For $\Sl_N\C$ the positive roots are
$\alpha_{ij}=E_{ii}-E_{jj}$ with $i<j$; therefore, we explicitly get
\begin{equation}\label{delambda}\index{$d_\lambda:=\text{Tr}_q(\text{id}_\lambda)$}
d_{\lambda}=\prod_{i<j}\frac{\sin\left(\frac{\pi(\lambda+\rho,\alpha_{ij})}{k+N}\right)}{\sin\left(
\frac{\pi(\rho,\alpha_{ij})}{k+N}\right)}
\end{equation}

\begin{exm}\label{raij}
Note that $\alpha_{ij}=\sum_{k=i}^{j-1}\alpha_{k}$ and
$\rho=\sum_{i=1}^{N-1}\omega_i$. Use the fact that
$\alpha_{i},\omega_j$ are biorthogonal to derive that
$(\rho,\alpha_{ij})=j-i$.
\end{exm}
\begin{exm}\label{unknot}
Show that $d_{\omega_1}=(\eps^N-\eps^{-N})/(\eps-\eps^{-1})$. This
is the invariant of the unknot in the fundamental representation.
\end{exm}

In principal to compute other invariants it appears that we have to
compute some braiding morphisms. Notice that when the dimension of a
representation is $n$, the $R$--matrix will be a $n^2$ by $n^2$
matrix. Clearly it is not practical to write out many $R$--matrices.
Fortunately, sometimes one can get away with knowledge of its action
only on elements of a special form, namely $u_\lambda\otimes v\in
U\otimes V$ with $u_\lambda$ being the highest weight vector in $U$
with the weight $\lambda\in\Lambda_w^+$. This observation is due to
V Turaev and H Wenzl \cite{TW} and will come handy for the
computations of $p^+_\lambda$ and $\wtilde s_{\lambda\mu}$.
\begin{prop}\label{Rhw} Let $u_\lambda$ ($v_\lambda$) be the highest weight vector of $U$ ($V$) and let $v$ ($u$) be arbitrary. Then
\begin{align}
R_{U,V}(u_\lambda\otimes v) &=u_\lambda\otimes q^\lambda v,\label{Rul}\\
R_{U,V}(u\otimes v_\lambda) &=q^\lambda u\otimes
v_\lambda+z_{<\lambda}\label{Rvl},
\end{align}
where $z_{<\lambda}$ is an element of $U\otimes V_{<\lambda}$ and
$V_{<\lambda}$ is the sum of weight subspaces of $V$ with weights
$<\lambda$.
\end{prop}
\begin{proof}
By definition of the highest weight we have
$e^{(n)}_{\alpha_i}u_\lambda=0$ for all $i$. We also get
$\smash{e^{(n)}_{\alpha}u_\lambda=0}$ for all positive roots from
\eqref{Rgen}. Hence the only term that survives in the sum of
products in \eqref{RslN} is the one that corresponds to $\nu=0$,
$n_1,\ldots,n_\nu=0$ and the whole sum reduces to $1\otimes 1$.
Therefore, we only have to compute
$\smash{q^{\sum_{i=1}^{N-1}\alpha_i\otimes\omega_i}}(u_\lambda\otimes
v)$. Since both sides of \eqref{Rul} are linear in $v$ it suffices
to prove it for the case when $v=v_\mu$ is a weight vector with
weight $\mu$. We have
\begin{align*}
q^{\sum_{i=1}^{N-1}\alpha_i\otimes\omega_i}(u_\lambda\otimes v_\mu)
&=q^{\sum_{i=1}^{N-1}(\lambda,\alpha_i)(\mu,\omega_i)}u_\lambda\otimes
v_\mu\\
&=q^{\sum_{i,j=1}^{N-1}(\lambda,\alpha_i)(\mu,\omega_j)(\omega_i,\alpha_j)}u_\lambda\otimes v_\mu\\
&=q^{(\sum_{i=1}^{N-1}(\lambda,\alpha_i)\omega_i,\sum_{i=1}^{N-1}(\lambda,\omega_j)\alpha_j)}u_\lambda\otimes
v_\mu =q^{(\lambda,\mu)}u_\lambda\otimes v_\mu
\end{align*}
since $\omega_i,\alpha_j$ are biorthogonal (see 
\fullref{app:d}). But by the definition of the weight
$q^{(\lambda,\mu)}v_\mu=q^\lambda v_\mu$ and we are done with the
first formula. The proof of the second formula proceeds analogously
but $v_\lambda$ is acted upon by $\smash{f^{(n)}_{\alpha}}$ and
$\smash{f^{(n)}_{\alpha}v_\lambda}$ is no longer zero. Instead, it
has a weight lower than $\lambda$ except in the case when $n=0$.
Thus, up to lower weight vectors in $V$ the sum of products in
\eqref{RslN} again reduces to $1\otimes 1$ and
$\smash{q^{\sum_{i=1}^{N-1}\alpha_i\otimes\omega_i}}$ reduces as
above to $q^\lambda$ acting on $u$.
\end{proof}

\begin{exm}
Derive formulas analogous to \eqref{Rul} and \eqref{Rvl} for the
case when $u_\lambda$, $v_\lambda$ are the lowest weight vectors.
\end{exm}

Now we compute $p_{\lambda}^\pm:=\Tr_q(\theta_{\lambda}^{\pm1})$.
Since the twist commutes with all morphisms (by Axiom 9) and
$\calW_\lambda^{\eps}$ is irreducible for $\lambda$ in the Weyl
alcove (\fullref{alcove}) by the Schur lemma
$\theta_{\lambda}^{\pm1}$ must be a scalar operator and it suffices
to compute it on a single vector in $\calW_\lambda^{\eps}$. In view
of the simple form of the $R$--matrix on pairs containing the
highest weight vector (\fullref{Rhw}) we choose the highest weight
vector $v_{\lambda}\in\calW_\lambda^{\eps}$:
\begin{align*}
\theta_\lambda v_{\lambda} :=&
  (\text{id}_\lambda\otimes\cap^*_\lambda)\circ
  (\times_{\lambda,\lambda}\otimes\text{id}_{\lambda^*})
  \biggl(\sum_k\,v_{\lambda}\otimes v_k\otimes v^k\biggr)\\
=&(\text{id}_\lambda\otimes\cap^*_\lambda)
  \biggl(\sum_k\times_{\lambda,\lambda}
  (v_{\lambda}\otimes v_k)\otimes v^k\biggr)\\
=&(\text{id}_\lambda\otimes\cap^*_\lambda)
  \biggl(\sum_k(q^{\lambda}v_k\otimes v_\lambda)\otimes v^k\biggr),
  \quad\text{by \eqref{Rul}}\\
=&\sum_k\, v^k(q^{2\rho}v_{\lambda})q^{\lambda}\,v_k
  =\sum_k\,\eps^{(2\rho,\lambda)}v^k(v_{\lambda})\,q^{\lambda}v_k\\
=&\eps^{(2\rho,\lambda)}q^{\lambda}\biggl(\sum_k\,v^k(v_{\lambda})v_k\biggr)
=\eps^{(2\rho,\lambda)}q^{\lambda}v_\lambda
=\eps^{(2\rho,\lambda)}\eps^{(\lambda,\lambda)}v_\lambda
  =\eps^{(2\rho+\lambda,\lambda)}v_\lambda.
\end{align*}
The case of $\theta_{\lambda}^{-1}$ is analogous so
\begin{equation}\label{tmat}
\theta_{\lambda}^{\pm1}=\eps^{\pm(2\rho+\lambda,\lambda)}\id_\lambda.
\end{equation}
Taking quantum traces on both sides yields.
\begin{equation}\label{pelambda}
p_{\lambda}^\pm=\eps^{\pm(2\rho+\lambda,\lambda)}\Tr_q(\id_\lambda)=e^{\pm\frac{\pi i}{k+N}(2\rho+\lambda,\lambda)}d_{\lambda}.
\end{equation}

\begin{exm}\label{fanom}
We have seen that $\theta_{\omega_1}$ is just multiplication by a
scalar. This scalar is called the framing anomaly. \index{framing
anomaly} Show that the framing anomaly is $\eps^{N-1/N}$.
\end{exm}

\subsection{$\wtilde{s}$--Matrix and Chern--Simons partition function
for $S^3$}\label{cscompute}

\begin{figure}[ht!]
\centering
\labellist\small
\pinlabel {$\lambda$} [r] at 21 120
\pinlabel {$\mu$} [r] at 21 69
\endlabellist
\includegraphics[width=.7truein]{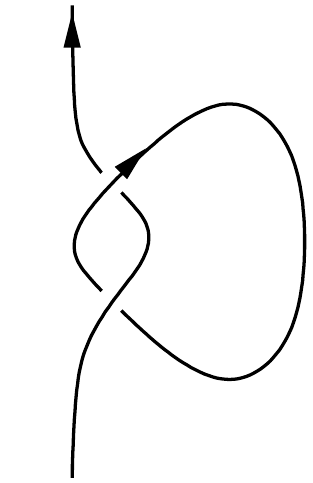}
\caption{The meridian morphism $\Gamma_{\lambda\mu}$}\label{meridian}
\end{figure}

In this subsection we compute the entries $\wtilde{s}_{\lambda\mu}$
of the structure matrix from Axiom 17 and we compute the partition
function of $S^3$. Following Turaev and Wenzl \cite{TW}  to compute
$\wtilde{s}_{\lambda\mu}$ we introduce the meridian morphism
$\Gamma_{\lambda\mu}\co
\calW_\lambda^{\eps}\to\calW_\lambda^{\eps}$, see
\fullref{meridian}.
One can see by inspection that
$\wtilde{s}_{\lambda\mu}=\Tr_q(\Gamma_{\lambda\mu})$. Again,
$\Gamma_{\lambda\mu}$ commutes with any morphism since we can slide
the latter up and down the $\lambda$ strand, and
$\Gamma_{\lambda\mu}$ acts as a scalar since it is an endomorphism
of an irreducible representation. It is therefore sufficient to
compute it on a highest weight vector $v_\lambda$. Let $\{u_i\}$ and
$\{u^i\}$ be dual bases for $V_\mu$ and $V_\mu^*$. Reading the
ribbon expression off \fullref{meridian} we have
\begin{align*}
v_\lambda &\mapsto \sum_i v_\lambda\otimes u_i\otimes u^i
\mapsto\sum_i
(q^\lambda u_i)\otimes v_\lambda\otimes u^i \\
&\mapsto \sum_i v_\lambda(q^{2\lambda}u_i+\text{lower
weight})\otimes u^i\mapsto \sum_i v_\lambda
u^i(q^{2(\rho+\lambda)}u_i)\\
&=\bigl(\text{Tr}\,\mu(q^{2(\rho+\lambda)})\bigr)v_\lambda=\chi_\mu(q^{2(\rho+\lambda)})v_\lambda.
\end{align*}
Here we abuse notation denoting by $\mu$ the irreducible representation of $\Sl_N\C$ indexed by the dominant weight $\mu$ so that $\text{Tr}\,\mu(\cdot)$ is the ordinary trace in this representation. Taking the quantum trace of this $\Gamma_{\lambda\mu}$ morphism implies that
\begin{equation}\label{selm}
\wtilde{s}_{\lambda\mu}=\Tr_q(\lambda(\Gamma_{\lambda\mu}))=\chi_\lambda(q^{2\rho})\chi_\mu(q^{2(\rho+\lambda)})=\chi_\mu(q^{2(\rho+\lambda)})d_\lambda.
\end{equation}
In particular, $\wtilde{s}_{00}=1$.

Since $S^3$ can be obtained from itself by a surgery on the empty
link $\emptyset$, we compute
$$
Z^{CS}(S^3):=\tau^{\Sl_N\C}(S^3):=(p^-)^{\sigma(\emptyset)}{\calD}^{-\sigma(\emptyset)-c(\emptyset)-1}F(\emptyset)={\calD}^{-1}\,,
$$
where recall $\calD^2:=\sum_{\lambda\in I}d_\lambda^2$ and
$p_-=\sum_{\lambda\in I} \theta^{-1}_\lambda d_\lambda$. Note that
there are many other links that produce $S^3$ and using expressions for
them produces universal algebraic relations among the characteristic
numbers of any modular category. Since we already computed $d_\lambda$
in \eqref{delambda} it appears that we just have to attach a sum and
a square root to the formula. The problem is that we are ultimately
interested in the Chern--Simons free energy and this is the logarithm
of the partition function. Hence a multiplicative, not an additive,
expression is desirable. This can be treated of course as a problem in
special functions theory, but we prefer a more conceptual approach based
on an understanding of the geometry of the Weyl alcove. The trick is to
represent the Weyl alcove as a fundamental domain of a group action;
see Kirillov \cite{Kir}, Kac \cite{Kac}, Samelson \cite{Samelson} and
Sawin \cite{Saw}. In addition this will allow us to finally verify the
non-degeneracy of the $\wtilde{s}$--matrix.

\begin{defn}
Given an integer $l\in\Z$ define the affine Weyl group $\AW$ of
$\Sl_N\C$ as the group of isometries of its Cartan subalgebra $\h$
generated by reflections about the hyperplanes
$(x,\alpha_i)=kl^\prime$ for every integer $k\in\Z$ and every simple
root $\alpha_i$ (as before $l^\prime$ is $l$ for $l$ odd and $l/2$
for $l$ even). We also define the translated action of $\AW$ as
$\wtilde{w}\cdot x:=\wtilde{w}(x+\rho)-\rho$, where $\rho$ is
the Weyl weight.
\end{defn}
\begin{remark}
Notice that we use both actions of elements of the Weyl group. In
the Weyl character formula the usual action given by reflections
perpendicular to simple roots is used. In the proof of
\fullref{nondegsmat} we use the translated action. We will always
use the notation $\sigma(\lambda)$ for the usual action and
$\sigma\cdot(\lambda)$ for the translated action.
\end{remark}
Setting $k=0$ in the definition we get the reflections $s_i$ that
generate the (ordinary) Weyl group $W\subset\AW$. Another type of
elements of $\AW$ is obtained by performing the reflection about
$(x,\alpha_i)=0$ followed by the reflection about
$$(x,\alpha_i)=l^\prime=l^\prime\frac{(\alpha_i,\alpha_i)}2=
  \Bigl(\frac{l^\prime\alpha_i}2,\alpha_i\Bigl).$$
This is easily seen to be the translation by $l^\prime\alpha_i$ and
we identified the subgroup of translations
$l^\prime\Lambda_r\subset\AW$. Every other reflection can be
performed by translating the corresponding hyperplane to the origin,
reflecting and then translating back. Thus we arrive at the
presentation
\begin{equation}\label{AWfact}
\AW=W\ltimes l^\prime\Lambda_r,
\end{equation}
where in the semi-direct product $W$ acts on $l^\prime\Lambda_r$ by
conjugation. In particular, $l^\prime\Lambda_r$ is a normal subgroup of $\AW$.
\begin{figure}[ht!]
\centering
\labellist\small
\pinlabel {$\mathcal{H}$} [b] at 70 282
\pinlabel {$(123)$} at 180 290
\pinlabel {$(13)$} at 130 210
\pinlabel {$(132)$} at 180 120
\pinlabel {$(12)$} at 280 290
\pinlabel {$(23)$} at 280 120
\pinlabel {$\rho{+}1$} at 320 210
\endlabellist
\includegraphics[width=4truein]{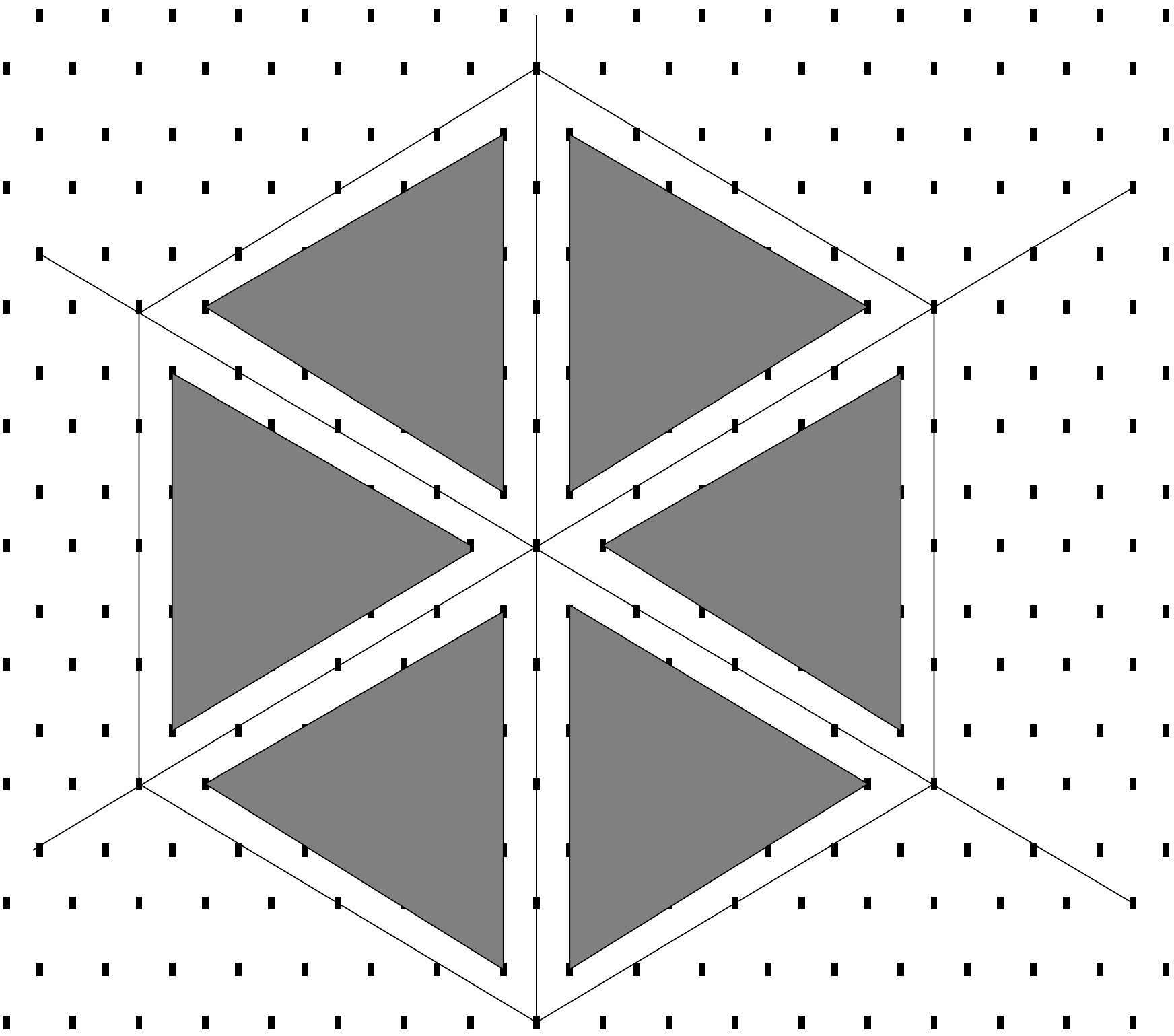} \caption{Affine
Weyl group and $\wwbar{C}^l$}\label{translatedI}
\end{figure}
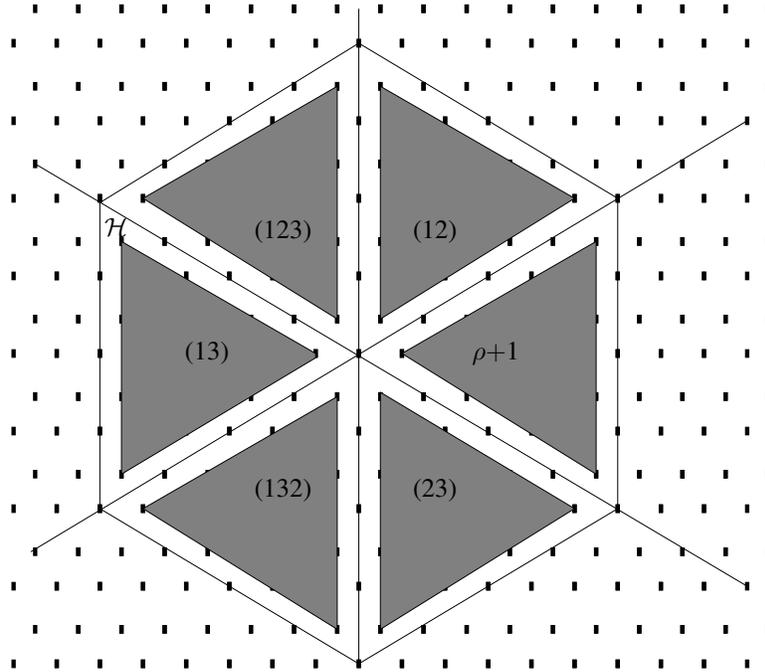

Recall that a closed set $D$ is a fundamental domain
of a continuous group action if every orbit intersects it
and no orbit intersects its interior more than once.
The subgroup of translations comprises integer multiples of
$l^\prime\alpha$ for $\alpha\in\Delta^+$ and therefore
one of its fundamental domains is the polygon
$$
\mathcal{H}:=\{x\in\h\mid\,(x,\pm\alpha)\leq l^\prime,\ \text{for
all}\ \alpha\in\Delta^+\}.
$$
Of course, any translate of $\mathcal{H}$ is also a fundamental
domain of $l^\prime\Lambda_r$. Note that by definition the
translated action of $l^\prime\Lambda_r$ is the same as the usual
one while $\mathcal{H}-\rho$ is invariant under the translated
action of $W$. The fundamental domain for the ordinary action of $W$
is well known to be the Weyl chamber (see Humphreys \cite{Hum} or Fulton
and Harris \cite{fulton-harris})
$$
\Lambda^+:=\{x\in\h\mid\,(x,\alpha)\geq 0,\ \text{for all}\
\alpha\in\Delta^+\}
$$
and hence $\Lambda^+-\rho$ serves as a fundamental domain for the
translated action.
\begin{exm}
Let $G=A\ltimes B$ with $B$ being a normal subgroup. Suppose $D_A$,
$D_B$ are fundamental domains of $A$ and $B$ respectively, and $D_B$
is invariant under $A$. Prove that $D:=D_A\cap D_B$ is a fundamental
domain of $G$.
\end{exm}
Intersecting $\mathcal{H}-\rho$ with $\Lambda^+-\rho$ (see \fullref{translatedI})
we establish
\begin{cor}\label{fundom}
The closed Weyl alcove
$$
\wwbar{C}^l=\mathcal{H}\cap\Lambda^+-\rho=\{x\in\h\mid\,0\leq(x+\rho,\alpha)\leq
l^\prime,\ \text{for all}\ \alpha\in\Delta^+\},
$$
is a fundamental domain for the translated action of the affine Weyl group $\AW$ on $\h$.
\end{cor}
After these geometric preliminaries let us move on to a
simplification of the quantum diameter $\calD^2:=\sum_{\lambda\in
I}d_\lambda^2$. First note that
$$
\wtilde{s}_{0\lambda}=\wtilde{s}_{\lambda0}=\Tr_q(\times_{0,\lambda}\circ\times_{\lambda,0})=\Tr_q(\id_\lambda)=d_\lambda.
$$
This suggests considering the square of the $\wtilde{s}$--matrix,
indeed
\begin{equation}\label{s00}
(\wtilde{s}^2)_{00}=\sum_{\nu\in\Lambda_w^l}\wtilde{s}_{0\nu}\wtilde{s}_{\nu0}
=\sum_{\nu\in\Lambda_w^l}d_\nu^2=\calD^2.
\end{equation}
Non-degeneracy of the $\wtilde{s}$--matrix was originally
established by Turaev and Wenzl in a roundabout way based on results
of Kac and Petersen from the theory of affine Lie algebras
\cite{TW}. The direct proof based on an explicit computation of
$\wtilde{s}^2$ that we present here is due to A Kirillov
\cite{BK,Kir}. As a bonus, it gives a nice formula for the quantum
diameter.  This proof requires the orthogonality relation for
characters of finite groups sketched in the next exercise. Another,
purely algebraic proof is given by M M\"uger \cite{Mug}.
\begin{exm}
Define $\langle \phi, \psi\rangle:=|G|^{-1}\sum_{g\in
G}\phi(g)\psi(g)^*$ for functions defined on an arbitrary finite
group. Prove that
$\langle\chi_\lambda,\chi_\mu\rangle=\delta_{\lambda\mu}$. This is
similar to \fullref{orch} from \fullref{app:d}.
\end{exm}

\begin{thm}\label{nondegsmat}
Let $l$ be an even integer. Then $(\wtilde{s}^2)_{\lambda\mu}=\calD^2\delta_{\lambda\mu^*}$ and
\begin{equation}\label{Ds00}
\calD^2=(-1)^{w_0}\frac{|\Lambda_w/l^\prime\Lambda_r|}{\delta_0(q^{2\rho})^2},
\end{equation}
where $w_0$ is the order-reversing permutation and $\delta_0$ is the
Weyl denominator \eqref{Wdenom}. In particular, $\calD\ne0$ and the
$\wtilde{s}$--matrix is non-degenerate.
\end{thm}
\begin{proof}
We first transform expression \eqref{selm} for
$\wtilde{s}_{\lambda\nu}$ using the Weyl character formula
\eqref{Wform}
\begin{equation}\label{prest}
\begin{aligned}
\wtilde{s}_{\lambda\nu}&=\chi_\nu(q^{2(\rho+\lambda)})d_\lambda
=\chi_\nu(q^{2(\rho+\lambda)})\chi_\lambda(q^{2\rho})\\
&=\frac{\delta_\nu(q^{2(\rho+\lambda)})}{\delta_0(q^{2(\rho+\lambda)})}
\frac{\delta_\lambda(q^{2\rho})}{\delta_0(q^{2\rho})}
=\frac{\delta_\nu(q^{2(\rho+\lambda)})\,\delta_\lambda(q^{2\rho})}
{\delta_\lambda(q^{2\rho})\,\delta_0(q^{2\rho})},\quad\text{by \eqref{dqlrho}}\\
&=\frac1{\delta_0(q^{2\rho})}\sum_{\sigma\in
W=\mathfrak{G}_N}(-1)^\sigma\,\eps^{2(\sigma(\lambda+\rho),\,\nu+\rho)}\,.
\end{aligned}
\end{equation}
Notice that while $\wtilde{s}_{\lambda\nu}$ is originally only
defined for weights $\lambda$, $\nu$ in the Weyl alcove
$\Lambda_w^l$, the right side of equation \eqref{prest} makes sense
for all weights. We extend the definition of
$\wtilde{s}_{\lambda\nu}$ to all weights by equation \eqref{prest}.
We can now replace the summation over $\Lambda_w^l$ in \eqref{s00}
by summation over a more convenient set, the quotient group
$\Lambda_w/l^\prime\Lambda_r$. This is done in two steps. First, we
know from \fullref{alcove} that the Weyl modules with the highest
weights in $\wwbar{C}^l\backslash C^l$ have quantum dimensions $0$
and \fullref{neglmorph} implies
$\wtilde{s}_{\lambda\nu}=\Tr_q(\times_{\lambda,\nu}\circ\times_{\nu,\lambda})=0$
for $\nu\in \wwbar{C}^l\backslash C^l$ so
\begin{equation*}
(\wtilde{s}^2)_{\lambda\mu}=\sum_{\nu\in\Lambda_w^l}\wtilde{s}_{\lambda\nu}\wtilde{s}_{\nu\mu}
=\sum_{\nu\in\wwbar{\Lambda}_w^l}\wtilde{s}_{\lambda\nu}\wtilde{s}_{\nu\mu},
\end{equation*}
where $\wwbar{\Lambda}_w^l:=\Lambda_w\cap\wwbar{C}^l$. For the
second step notice that the value of $\wtilde{s}_{\lambda\nu}$ is
only changed by the sign of the orientation of the transformation
when $\lambda$ or $\nu$ are acted on by elements of the affine Weyl
group via the translated action. Indeed applying a translation
$l^\prime\alpha_i$ from the affine Weyl group to $\nu$ does not
change the value of $\wtilde{s}_{\lambda\nu}$ since
\[\eps^{2(\sigma(\lambda+\rho),\,\nu+l^\prime\alpha_i+\rho)}
=\eps^{2(\sigma(\lambda+\rho),\,\nu+\rho)}
\eps^{2(\sigma(\lambda+\rho),l^\prime\alpha_i)}=\eps^{2(\sigma(\lambda+\rho),\,\nu+\rho)}\,,
\]
and $\alpha_i$ pairs with any weight to give an integer and
$\eps^{2l^\prime}=1$. Applying the translated action of an element
of the Weyl group to $\lambda$ gives
\begin{equation}\label{awsign}
\begin{aligned}
\wtilde{s}_{\tau\cdot\lambda\nu}
&=\frac1{\delta_0(q^{2\rho})}\sum_{\sigma\in
W=\mathfrak{G}_N}(-1)^\sigma\,\eps^{2(\sigma(\tau\cdot\lambda+\rho),\,\nu+\rho)}\\
&=\frac1{\delta_0(q^{2\rho})}\sum_{\sigma\in
W=\mathfrak{G}_N}(-1)^\sigma\,\eps^{2(\sigma(\tau(\lambda+\rho)),\,\nu+\rho)}\\
&=\frac{(-1)^{\tau}}{\delta_0(q^{2\rho})}\sum_{\sigma\in
W=\mathfrak{G}_N}(-1)^\sigma\,\eps^{2(\sigma(\lambda+\rho),\nu+\rho)}\,.
\end{aligned}
\end{equation}
By \fullref{fundom} and \eqref{awsign} summation over
$\wwbar{\Lambda}_w^l$ can be replaced by summation over equivalence
classes in $\Lambda_w/\AW$ and furthermore in view of \eqref{AWfact}
\begin{equation}\label{newsum}
\sum_{\nu\in\wwbar{\Lambda}_w^l}\wtilde{s}_{\lambda\nu}\wtilde{s}_{\nu\mu}
=\sum_{\nu\in\Lambda_w/\AW}\wtilde{s}_{\lambda\nu}\wtilde{s}_{\nu\mu}
=\frac1{|W|}\sum_{\nu\in\Lambda_w/l^\prime\Lambda_r}\wtilde{s}_{\lambda\nu}\wtilde{s}_{\nu\mu}.
\end{equation}
This gives
\begin{equation}\label{dubsum}
(\wtilde{s}^2)_{\lambda\mu}=\frac1{|W|\,\delta_0(q^{2\rho})^2}\sum_{\nu\in\Lambda_w/l^\prime\Lambda_r}
\sum_{\sigma_1,\sigma_2\in
W}(-1)^{\sigma_1\sigma_2}\,\eps^{2(\sigma_1(\lambda+\rho)+\sigma_2(\mu+\rho),\,\nu+\rho)}.
\end{equation}
To simplify the last double sum we change the order of summation and
look closer at the maps $\eps^{2(\kappa,\,\cdot\,)}$, where
$\kappa\in\Lambda_w$. Since $\eps^{2l^\prime}=\eps^l=1$ they are
well-defined as maps $\Lambda_w/l^\prime\Lambda_r\to\C^*$ and in
fact are characters of the Abelian group
$\Lambda_w/l^\prime\Lambda_r$. Since $\eps$ is a primitive $l$th
root of unity $\eps^{2(\kappa,\,\cdot\,)}$ defines the trivial
character if and only if $(\kappa,\nu)$ is divisible by $l/2$ for
any $\nu\in\Lambda_w$. Since $l$ is even and $l^\prime=l/2$ this is
equivalent to $\kappa\in l^\prime\Lambda_r$. The orthogonality
relation for characters now yields
\begin{multline*}
(\eps^{2(\kappa,\,\cdot\,)},1)=\sum_{\nu\in\Lambda_w/l^\prime\Lambda_r}\eps^{2(\kappa,\,\nu)}
=\sum_{\nu\in\Lambda_w/l^\prime\Lambda_r}\eps^{2(\kappa,\,\nu+\rho)}
=\begin{cases} 0,\qquad\quad\ \ \kappa\notin l^\prime\Lambda_r\\
|\Lambda_w/l^\prime\Lambda_r|,\,\kappa\in l^\prime\Lambda_r.\end{cases}
\end{multline*}
This implies that the terms in the double sum \eqref{dubsum} are $0$
unless $\sigma_1(\lambda+\rho)+\sigma_2(\mu+\rho)\in
l^\prime\Lambda_r$. But by \fullref{alcove}:
$\mu+\rho=-w_0(\mu^*+\rho^*)=-w_0(\mu^*+\rho)$ since $\rho^*=\rho$.
Therefore, the following are equivalent.
\begin{align*}
&\sigma_1(\lambda+\rho)+\sigma_2(\mu+\rho)\in l^\prime\Lambda_r \\
&\sigma_1(\lambda+\rho)-\sigma_2w_0(\mu^*+\rho)\in
l^\prime\Lambda_r \\
&\lambda+\rho\in
\sigma_1^{-1}\sigma_2w_0(\mu^*+\rho)+l^\prime\Lambda_r=:\sigma(\mu^*+\rho)+l^\prime\Lambda_r \\
&\lambda\in
\sigma(\mu^*+\rho)-\rho+l^\prime\Lambda_r=\sigma\cdot\mu^*+l^\prime\Lambda_r \\
&\lambda=\wtilde{w}\cdot\mu^*,\,\text{for $\wtilde{w}\in\AW$.}
\end{align*}
But $\lambda,\mu^*\in\Lambda_w^l\subset C^l$ are within a
fundamental domain of the translated action and the last condition
can only be satisfied when $\lambda=\mu^*$ and $\wtilde{w}=1$ or
equivalently $\sigma_1=\sigma_2w_0$. In particular,
$(\wtilde{s}^2)_{\lambda\mu}=0$ for $\lambda\ne\mu^*$. When
$\lambda=\mu^*$ expression \eqref{dubsum} reduces to
\begin{align*}
(\wtilde{s}^2)_{\lambda\mu}&=\frac1{|W|\,\delta_0(q^{2\rho})^2}
\sum_{w_1\in W}(-1)^{2w_1+w_0}\,|\Lambda_w/l^\prime\Lambda_r|\\
&=\frac1{|W|\,\delta_0(q^{2\rho})^2}(-1)^{w_0}|W||\Lambda_w/l^\prime\Lambda_r|
=(-1)^{w_0}\frac{|\Lambda_w/l^\prime\Lambda_r|}{\delta_0(q^{2\rho})^2}.
\end{align*}
The formula for the quantum diameter now follows directly from
\eqref{s00}.
\end{proof}

Using the formula for the Weyl denominator, we obtain the following
formula for the quantum diameter that generalizes to any semisimple
Lie algebra.
$$\index{$D$@$\calD:=(\sum d_\lambda^2)^{\frac12}$}\index{quantum
diameter} \calD =i^{w_0}|\Lambda_w/l^\prime\Lambda_r|^{1/2}
\prod_{\alpha\in\Delta^+}(\eps^{(\rho,\alpha)}-\eps^{-(\rho,\alpha)})^{-1}\,.
$$
We can make this formula  more explicit for $\Sl_N\C$. Recall that
the level is $k:=l^\prime-N$. By \eqref{dqrho} and
\fullref{raij}
\begin{multline*}
\delta_0(q^{2\rho})=\prod_{\alpha\in\Delta^+}(\eps^{(\rho,\alpha)}-\eps^{-(\rho,\alpha)})
=\prod_{i<j}\,2i\frac{e^{\frac{\pi i(\rho,\alpha_{ij})}{k+N}}-e^{-\frac{\pi i(\rho,\alpha_{ij})}{k+N}}}{2i}\\
=i^{\frac{N(N-1)}2}\prod_{i<j}2\,\sin\left(\frac{\pi(j-i)}{k+N}\right)
=i^{\frac{N(N-1)}2}\prod_{j=1}^{N-1}2\,\sin\left(\frac{\pi j}{k+N}\right)^{N-j}.
\end{multline*}

\begin{exm}
Prove that $|\Lambda_w/(k+N)\Lambda_r|=N(k+N)^{N-1}$, and show that
$(-1)^{w_0}=(-1)^{N(N-1)/2}$.
\end{exm}
We arrive at the following.
\begin{cor} The Chern--Simons partition function for the three-sphere is
\begin{equation}\label{D-1}
Z^{CS}(S^3)=\calD^{-1}=N^{-1/2}(k+N)^{(1-N)/2}
\prod_{j=1}^{N-1}\left(2\sin\left(\frac{\pi j}{k+N}\right)\right)^{N-j}.
\end{equation}
\end{cor}
\begin{remark}
When $l$ is odd the level $k=l/2-N$ is only a half-integer and
Witten's heuristic argument for the invariance of $Z^{CS}$ breaks
down. Thus, there is no `physical' reason to expect that
$U_\eps^{\text{res}}(\Sl_N\C)$ produces topological invariants for $\eps$ a
primitive odd root of unity. The above proof of non-degeneracy of
the $\wtilde{s}$--matrix also fails for $l$ odd. Indeed, it is
sufficient that $\kappa\in \frac{l}2\Lambda_r$ as opposed to
$\kappa\in l\Lambda_r=l^\prime\Lambda_r$ for
$\eps^{2(\kappa,\,\cdot\,)}$ to define the trivial character and we
can not reduce the double sum completely. As a matter of fact, the
$\wtilde{s}$--matrix {\it can be degenerate} in this case. This is
due to appearance of nontrivial transparent \cite{Brug} (also
called degenerate \cite{Saw} or central \cite{Mug}) simple objects
in $\Tilt(\Sl_N\C)$. These are the Weyl modules $\calW_{\tau}^\eps$
such that $R_{\lambda,\tau}=R_{\tau,\lambda}^{-1}$ for all alcove
weights $\lambda$. Obviously, the trivial object is transparent in
any ribbon category. The name comes from the fact that
$\times_{\lambda,\tau}=\times_{\tau,\lambda}$ and any strand can be
pushed through $\tau$--colored ones without changing the invariant.
As a result, transparent objects create a row
$\wtilde{s}_{\tau\lambda}$ in the $\wtilde{s}$--matrix that is
proportional to the row of the trivial object
$\wtilde{s}_{0\lambda}$ making the matrix degenerate. When $l$ is
odd and $N$ is even the fundamental weight $\omega_{\frac N2}$
indexes a nontrivial transparent object in $\Tilt(\Sl_N\C)$. In
particular, for $\Sl_2\C$ reduced tilting modules at odd roots do
not form a modular category (see Sawin \cite{Saw}). Moreover, under mild
technical assumptions one can prove that a category satisfying
Axioms 1--16 of \eqref{axmtc} (such categories are called premodular)
with $\calD\ne0$ is modular if and only if it does not have
nontrivial transparent simple objects (see Brugui\'eres \cite{Brug} and
M\"uger \cite{Mug}).
Surprisingly, $\Tilt(\Sl_N\C)$ can be 'modularized' by a quotient
construction and made to produce nontrivial invariants even when
$l$ is odd except when $N=2\mod4$ \cite{Saw}. At present, a
geometric or physical explanation for this phenomenon is lacking.
\end{remark}

We conclude with some remarks about computing the Reshetikhin--Turaev
invariants for more general $3$--manifolds. There are two different
ways that may be used to organize many such computations. The first
method is via the shadow invariants of V Turaev \cite{T}. The
second method is to expand the theory into a full topological
quantum field theory. Most computations in the physical literature
use the topological quantum field theory viewpoint.

Recall the twist matrix
$\wtilde{t}_{\lambda\mu}:=\eps^{(2\rho+\lambda,\lambda)}\delta_{\lambda\mu}$
from \eqref{tmat}. For many $3$--manifolds such as circle bundles
over Riemann surfaces one can avoid using the $R$--matrices directly
by utilizing the fact that  renormalized $\wtilde{s}$,
$\wtilde{t}$ matrices form a projective representation of the
modular group $\text{SL}_2\Z$ that appears in many surgeries.
Namely, set $s:=\calD^{-1}\wtilde{s}$ and
$t:=\zeta^{-1}\wtilde{t}$ with $\zeta:=(p^+/p^-)^\frac16$ then
$s,t$ can be taken as the images of the standard generators $S,T$ of
$\text{SL}_2\Z$
$$
S=\left(\begin{array}{cc}0&-1\\1&0\end{array}\right), \qquad
T=\left(\begin{array}{cc}1&1\\0&1\end{array}\right)\,.
$$

\begin{exm}
Show that $s,t$ satisfy the standard relations of $\text{SL}_2\Z$:
$(st)^3=s^2,\,s^2t=ts^2,\,s^4=1$, see Bakalov and Kirillov \cite{BK,Kir}.
\end{exm}
In terms of the $s$--matrix we have
\begin{align*}
d_{\lambda}&=\calD s_{0\lambda}=\frac{s_{0\lambda}}{s_{00}}\\
Z^{CS}(S^3)&=\calD^{-1}=s_{00}
\end{align*}
These identities are frequently used in practical computations.

>We combine equations \eqref{prest} and \eqref{Ds00} in this subsection into the following corollary.
\begin{cor} In any
semisimple Lie algebra we have the following explicit formula   for
the $s$--matrix and $\wtilde{t}$--matrix.
\index{$T$@$s:=\calD^{-1}\wtilde{s}$ $s$--matrix}
\index{$T$@$\wtilde{t}$--matrix}
\begin{align*}
s_{\lambda\nu}=i^{-|\Delta_+|}|\Lambda_w/l^\prime\Lambda_r^\vee|^{-1/2}&\sum_{\sigma\in
W} (-1)^\sigma \eps^{2\langle\!\langle
\sigma(\lambda+\rho),\nu+\rho\rangle\!\rangle}\,.\\
\wtilde{t}_{\lambda\mu}:=\eps^{\langle\!\langle
2\rho+\lambda,\lambda\rangle\!\rangle}\delta_{\lambda\mu}\,.
\end{align*}
Here $\Lambda_r^\vee$ is the coroot lattice and $\langle\!\langle
\cdot,\cdot\rangle\!\rangle$ is the Killing form normalized so that
$\langle\!\langle \alpha,\alpha\rangle\!\rangle=2$ for short roots.

Specializing to $\Sl_N\C$ gives
$$
s_{\lambda\nu}=N^{-1/2}(k+N)^{(1-N)/2}(i)^{-N(N-1)/2}\sum_{\sigma\in{\mathfrak
G}_N}(-1)^\sigma
\eps^{2\langle\rho+\nu,\sigma(\rho+\lambda)\rangle}\,.
$$
\end{cor}

\begin{remark}
The sign in the exponent $(i)^{-N(N-1)/2}$ differs from that of
\cite{BK} since we are using the right-handed Hopf link to define
the $s$--matrix (not the left-handed one.) The sign used here and the
use of the right-handed Hopf link are correct.
\end{remark}


In the next section we compute the free energy of the $3$--sphere and
perform the final comparison between the Gromov--Witten free energy
and Chern--Simons free energies.

\newpage
\part{Comparisons and recent developments}
\setobjecttype{Part}
\label{part3}

\section{Comparison of the free energies}\label{comparison}

In this section we reach the final goal of this paper by comparing
the Gromov--Witten and Chern--Simons free energies for $X_{S^3}$ and
$S^3$. We start by recalling some necessary facts and formulas from
the theory of special functions that are used in transforming the
expressions derived earlier in the paper.  We apply these formulas
to reduce both energies to a similar form and perform the
comparison. Even though the energies do not match exactly the
difference does not contain any positive powers of
$x:={2\pi}/{(k+N)}$ (called the string coupling constant and denoted
$g_s$ in the physical literature) which means that we have an exact
match for the counting invariants.

\subsection{Bernoulli numbers and special functions}\label{specfs}

After generating functions for the invariants are computed the claim
of the Gopakumar--Vafa Large $N$ Duality reduces to a problem in the
theory of special functions. Aside from Bernoulli numbers we need
the Eisenstein functions, the Riemann zeta-function, polylogarithms
and the Barnes function. Rather than simply referring the reader to
various sources for necessary formulas we briefly review the
definitions and relationship between them in this subsection.

Since free energies are essentially the natural logarithms of
partition functions it helps to present all terms in partition
functions as products. In particular, the Chern--Simons partition
function contains sines and we start by deriving the Euler product
formula for the sine function. There are different ways to do this
but we choose the one using Eisenstein functions since we need them
later, see Weil \cite{Weil}.
\begin{defn}\label{Eisn}
The Eisenstein functions are defined by $$ E_k(z):=\text{{\rm
v.p.}}\!\!\!\sum^\infty_{n=-\infty} (z+n)^{-k}, $$ where we use the
Eisenstein convention:
$$
\text{{\rm
v.p.}}\!\!\!\sum^\infty_{n=-\infty}f(n):=\lim_{N\to\infty}\sum^N_{n=-N}f(n)
$$
\end{defn}
The main properties of the Eisenstein functions are collected in the
following exercises:
\begin{exm}\label{44}
Show that $E_1(z)-\pi\cot(\pi z)$ is a bounded entire function that
evaluates to zero at $z=\frac12$ and use Liouville's theorem to
conclude that it is identically zero.
\end{exm}
\begin{exm}
From the definition, it follows that $E_k^\prime(z)=-kE_{k+1}(z)$,
so all of these functions may be found by differentiating the
cotangent function. Using this show that $E_2(z)=\pi^2\csc^2(\pi z)$
and $E_3(z)=E_1(z)E_2(z)$.
\end{exm}
\begin{exm}
Set $s(z)=z\prod_{n=1}^\infty \bigl(1-\frac{z^2}{n^2}\bigr)$ and prove that
$s^\prime(z)/s(z)=E_1(z)$.
\end{exm}
\begin{exm}\label{siniss}
Express $\frac{d}{dz}(s(z)^2E_2(z))$ using just $s(z)$, $E_1(z)$ and $E_2(z)$. Conclude that it is zero.
\end{exm}
It is not hard to see that $s(0)=0$ and $s^\prime(0)=1$. It follows that $$\lim_{z\to 0}\frac{\pi s(z)}{\sin \pi z}=1\,.$$
However we know from \fullref{siniss} that $\left(\frac{\pi s(z)}{\sin \pi z}\right)^2$ is constant. Thus
$$
\sin(\pi z)=\pi z\prod_{n=1}^\infty \Bigl(1-\frac{z^2}{n^2}\Bigr).
$$
To transform multiple sums that appear in the Chern--Simons free
energy we need closed formulas for power sums of natural numbers.
These formulas involve the Bernoulli numbers.
\begin{defn}\label{Bern} The Bernoulli numbers \index{$B_k$ Bernoulli numbers}
$B_k$ are defined by \index{Bernoulli numbers} their generating function:
$$  \frac{z}{e^z-1}=\sum_{k=0}^\infty B_k\frac{z^k}{k!}\,. $$
\end{defn}
Again we collect the main properties we need in an exercise, see
Andrews--Askey--Roy \cite{AAR}.
\begin{exm}\label{coth}
Show that
$\frac{z}{2}\coth\bigl(\frac{z}{2}\bigr)=\frac{z}{2}+\frac{z}{e^z-1}$ and
conclude that $B_1=-\frac12$ and the rest of the odd Bernoulli
numbers are zero. Also compute $B_0$, $B_2$, $B_4$ and $B_6$.
\end{exm}
\begin{exm}
Expand $\smash{\frac{e^{Nz}-1}{z}\frac{z}{e^z-1}}$ as a sum of exponentials,
and then expand the exponentials into power series to obtain
$N+\smash{\sum_{n=1}^\infty\sum_{j=1}^{N-1}j^p\frac{z^n}{n!}}$.
\end{exm}
\begin{exm}
Expand each factor of $\frac{e^{Nz}-1}{z}\frac{z}{e^z-1}$ in a power
series and multiply the two resulting power series to obtain
$$
\sum_{n=0}^\infty\sum_{j=0}^n\begin{pmatrix}n\\j\end{pmatrix}B_{j}\frac{N^{n-k+1}}{n-k+1}\frac{z^n}{n!}\,.
$$
\end{exm}
\begin{exm}\label{sumpower}
Compare the expansions from the two previous exercises to obtain a
formula for the sum of powers of the first $N-1$ natural numbers.
Use the fact that all of the odd Bernoulli numbers other than $B_1$
vanish together with the definition of the binomial coefficients to
conclude
\begin{align*}
\sum_{j=1}^{N-1}j^{2p+1}&=-\frac{N^{2p+1}}{2}+\sum_{k=0}^p\begin{pmatrix}2p+1\\2k\end{pmatrix} \frac{B_{2k}N^{2p-2k+2}}{2p-2k+2}\,, \\
\sum_{j=1}^{N-1}j^{2p}&=- \frac{N^{2p}}{2}+\sum_{k=0}^p \begin{pmatrix}2p+1\\2k\end{pmatrix}\frac{B_{2k}N^{2p-2k+1}}{2p+1}\,.
\end{align*}
\end{exm}
Bernoulli numbers are closely related to values of the celebrated
Riemann zeta function at even integers.
\begin{defn}\label{zeta}\index{$\zeta(z)$ the Riemann zeta function}
The Riemann zeta function is defined by (Andrews--Askey--Roy \cite{AAR}):
$$
\zeta(z):=\frac{1}{\Gamma(z)}\int_0^\infty\frac{u^{z-1}}{e^u-1}\,du=\sum_{n=1}^\infty
n^{-z}\,,
$$ where $\Gamma(z)$ is the usual gamma function of Euler
$$
\Gamma(z):=\int_0^\infty e^{-t}t^{z-1}\,dt\,.
$$
\end{defn}
\begin{exm}\label{zetber}
Use integration by parts to derive the usual relation between the
gamma function and the factorial. Expand $(e^u-1)^{-1}$ in powers of
$e^{-u}$ and substitute into the definition of the zeta function to
obtain the formula $\zeta(z)=\sum_{n=1}^\infty n^{-z}$. Use the
generating function definition of the Bernoulli numbers in the
definition of the zeta function to obtain for $n\ge 1$,
\begin{equation}
\label{zeta2n}
\zeta(2n)=(-1)^{n+1}\frac{(2\pi)^{2n}B_{2n}}{2(2n)!}\,.
\end{equation}
\end{exm}
The gamma function can be meromorphically extended to the entire
complex plane via $\Gamma(z+1)=z\Gamma(z)$ and the zeta function can be
meromorphically extended to the entire complex plane via the functional
equation discovered for the zeta function by Riemann himself (see Hardy
and Wright \cite{HW}):
$$\zeta(1-k)=2(2\pi)^{-k}\cos(\pi k/2)\Gamma(k)\zeta(k)\,.$$
The Chern--Simons partition function contains a factor that can be
identified with the volume of SU$(N)$. This volume can be expressed
in terms of the Barnes function.

\begin{defn}
The Barnes function \index{$G_2(\cdot)$ the Barnes function} is
defined by
$$
G_2(z+1)=(2\pi)^{z/2}e^{-\left(z(z+1)+\gamma
z^2\right)/2}\prod_{k=1}^\infty
\left((1+z/k)^ke^{(z^2-z)/2k}\right)\,,
$$
where $\gamma$ is the Euler constant given by $\gamma:=\lim_{n\to
0}\left(\sum_{k=1}^nk^{-1}-\ln n\right)$.
\end{defn}

\begin{remark}
For integers greater than $1$ one can show that
$G_2(N):=\prod_{j=1}^{N-2}j!$.
\end{remark}

The following asymptotic expansion for $G_2$ in terms of Bernoulli
numbers is given by Adamchik \cite{adm}:
\begin{multline}\label{Barnes}
\ln(G_2(N+1))= \\
\tfrac12 N^2\ln N-\tfrac34 N^2-\tfrac{1}{12}\ln N
  -N\zeta^\prime(0)+\zeta^\prime(-1)+\sum_{g=2}^\infty
  \frac{B_{2g}}{2g(2g-2)}N^{2-2g}
\end{multline}

\begin{remark}
The formula in \cite{adm} has a negative sign in front of the sum.
Careful inspection shows that a sign was lost when the expression
following line (29) in this paper was substituted into equation (20)
from this paper.
\end{remark}

Finally, manipulations with the Gromov--Witten free energy require
the use of polylogarithms.
\begin{defn}\label{polylog}\index{$\text{Li}_p(z)$ polylogarithm
function}
The polylogarithm functions are defined by (see eg \cite{CLZ})
\begin{equation}
\text{Li}_p(z):=\sum_{n=1}^\infty n^{-p}z^n\,. \end{equation}
\end{defn}
One can see by inspection that $\text{Li}_1(z)=-\ln(1-z)$and $z\frac{d}{dz}\text{Li}_{p+1}(z)=\text{Li}_{p}(z)$. Therefore, polylogarithms indexed by positive integers have a logarithmic branch point at $z=1$. On the other hand, polylogarithms with negative integers are meromorphic in the complex plane and relate to the Eisenstein functions by a change of variables and renormalization.
\begin{exm}
Expand $\pi\cot(\pi z)$ as a power series in the exponential
$e^{-2\pi iz}$ (see \cite{AAR}) and differentiate the expression obtained by
combining this with \fullref{44} and the definition of the
polylogarithm to obtain for $q\ge 1$
$$
\text{Li}_{-q}(e^{-2\pi i z})=\frac{q!}{(2\pi
i)^{q+1}}\sum_{n=-\infty}^\infty (n+z)^{-q-1}
=\frac{q!}{(2\pi i)^{q+1}}E_{q+1}(z)\,.
$$
\end{exm}
We will need power series expansions for $\text{Li}_{3-2g}(e^{-t})$ at $0$. Note also that these functions are manifestly periodic with the period $2\pi i$. For
$g\ge 2$ combining the last exercise with the negative power
binomial theorem gives
\begin{equation}
\begin{aligned}\label{lipower}
\text{Li}_{3-2g}(e^{-t})&=(2g{-}3)!t^{2-2g}{+}
  \frac{(2g{-}3)!}{(2\pi i)^{2g-2}}
  \sum_{\stackrel{n\ne 0}{n=-\infty}}^\infty
  (1{+}(t/2\pi i n))^{2-2g}n^{2-2g}\\
&=(2g{-}3)!t^{2-2g}\\
&\quad+(2g{-}3)!\sum_{\stackrel{n\ne 0}{n=-\infty}}^\infty
  \sum_{h=0}^\infty\binom{2g{+}h{-}3}{h}(-t)^h(2\pi i)^{2-2g-h}n^{2-2g-h}\\
&= (2g{-}3)!t^{2-2g} \\
&\quad+(2g{-}3)!\sum_{\stackrel{h\ \text{even}}{h\ge 0}}2
\binom{2g{+}h{-}3}{h}(-t)^h(2\pi i)^{2-2g-h}\zeta({2g{+}h{-}2})\,.
\end{aligned}
\end{equation}
As mentioned above the functions $\text{Li}_1(e^{-t})$ and $\text{Li}_3(e^{-t})$ have a branch point at $t=0$ and can not be expanded into a Taylor series. However, this is easy to fix by adding logarithmic 'counterterms' that render them holomorphic in a neighborhood of $0$. For instance, $\text{Li}_1(e^{-t})=-\ln(1-e^{-t})$ behaves like $-\ln t$ near $0$ so the sum $\text{Li}_1(e^{-t})+\ln t$ is holomorphic. Of course, by writing $\ln t$ we are implicitly fixing a branch of the logarithm and thus the polylogarithm as well. The corresponding expansions are derived in the next lemma.
\begin{lemma}\label{Li13}
The polylogarithms admit the following series expansions:
\begin{align}
\text{Li}_1(e^{-t})+\ln t
&=t/2+\sum_{\stackrel{\scriptstyle m\ \text{even}}{m\ge 2}}
\frac{2}{m}(2\pi)^{-m}\zeta(m)(it)^{m}, \label{Li1}\\[-1ex]
\text{Li}_3(e^{-t})+\frac{t^2}{2}\ln t
&=\zeta(3)-\zeta(2)t+3 t^2/4+t^3/12 \label{Li3} \\[-1ex]
&\quad-\sum_{\stackrel{\scriptstyle m\ \text{even}}{m\ge 4}}
\frac{2}{m(m-1)(m-2)}(2\pi)^{2-m}\zeta(m-2)(it)^{m}.\notag
\end{align}
\end{lemma}

\begin{proof}
The function
$f(t)=\ln(t)+\text{Li}_1(e^{-t})=\ln\left(t(1-e^{-t})^{-1}\right)$
is clearly holomorphic at zero. Now compute
\[
\begin{aligned}
f^\prime(t)&=t^{-1}\left(1-\frac{t}{e^t-1}\right)\\
&=\frac12-\sum_{n=1}^\infty \frac{B_{2n}}{(2n)!}t^{2n-1}\,.
\end{aligned}
\]
Since $f(0)=0$ we can integrate and use formula \eqref{zeta2n} relating
the
Bernoulli numbers to the zeta function see that
\[
\begin{aligned}
f(t)&=\tfrac12 t-\sum_{n=1}^\infty \frac{B_{2n}}{(2n)!(2n)}t^{2n}\\
&=\tfrac12 t+\sum_{n=1}^\infty \frac{(2\pi)^{-2n}}{n}\zeta(2n)(it)^{2n}\,.
\end{aligned}
\]
For the second equality set
$g(t)=\text{Li}_3(e^{-t})+\frac12 t^2\ln t -\tfrac{3}{4}t^2$. One
easily computes $g(0)=\zeta(3)$,
$$g^\prime(t)=-\text{Li}_2(e^{-t})+t\ln t-t,$$
giving $g^\prime(0)=-\zeta(2)$ and
$$g^{\prime\prime}(t)=\text{Li}_1(e^{-t})+\ln t
  =\tfrac12 t+\sum_{\stackrel{\scriptstyle m \text{ even}}{m\ge 2}}
  \frac{2}{m}(2\pi)^{-m}\zeta(m)(it)^{m}.$$
Integrating twice as above with the computed constants of integration,
reindexing the sum with $h=m+2$ and using the definition of $g(t)$
gives
\begin{multline*}
\text{Li}_3(e^{-t})+\tfrac12 t^2\ln t =\zeta(3)-\zeta(2)t+\tfrac34 t^2+
\tfrac{1}{12}t^3\\
  -\sum_{\stackrel{h \text{ even}}{h\ge 4}}
  \frac{2}{h(h-1)(h-2)}(2\pi)^{2-h}\zeta(h-2)(it)^{h}.
\end{multline*}
This completes the proof.
\end{proof}

\subsection{The Chern--Simons free energy}\label{csfree}

Recall from equation \eqref{D-1} that the
$SU(N)$--Witten--Reshetikhin--Turaev invariant of $S^3$ is given by
$$
Z(S^3)=\tau^{{\mathfrak{sl}}_N\C}_k(S^3)=N^{-1/2}(k+N)^{(1-N)/2}\prod_{j=1}^{N-1}
\left[2\sin\left(\frac{\pi j}{k+N}\right)\right]^{N-j}\,.
$$
Recall that the Chern--Simons free energy is given by
$$F_M=\ln(Z^{\text{U$(N)$}}(M)/Z_0^{\text{U$(N)$}}(M))\,.$$
The formal mathematical definition of the $U(N)$ partition function
has not yet been agreed upon. In addition there is no mathematical
definition of $Z_0^{\text{U$(N)$}}(M)$. To check Large $N$ Duality
we define the unnormalized Chern--Simons free energy to be the
logarithm of the SU$(N)$--Witten--Reshetikhin--Turaev invariant.
\begin{defn}\label{CSunnorm}
The unnormalized Chern--Simons free energy \index{unnormalized
Chern--Simons free energy}\index{Chern--Simons free
energy}\index{$F^{\text{CS}}_M$ unnormalized Chern--Simons free
energy} is
$$F^{\text{CS}}_M=\ln(\tau^{{\mathfrak{sl}}_N\C}_k(M))\,.$$
\end{defn}
It follows immediately that
\begin{equation}\label{FS3}
\begin{aligned}
F^{\text{CS}}_{S^3}=& \frac{1-N}{2}\ln(k+N)-\tfrac12\ln
N+\sum_{j=1}^{N-1}
(N-j)\ln\Bigl[2\sin\Bigl(\frac{\pi j}{k+N}\Bigr)\Bigr]\\[-1ex]
=& \frac{1-N}{2}\ln(k+N)-\tfrac12\ln
N+\sum_{j=1}^{N-1}(N-j)\ln\Bigl(\frac{2\pi j}{k+N}\Bigr)\\[-1ex]
&\qquad+ \sum_{j=1}^{N-1}(N-j)\biggl(\sum_{n=1}^\infty \ln\Bigl(1-\frac{j^2}{n^2(k+N)^2}
 \Bigr)\biggr)\,.
\end{aligned}
\end{equation}
Here we have used the product expansion of the sine function derived
just below Exercises \ref{44}--\ref{siniss}. We will now
concentrate on the last sum in this expression. We call it the
\index{perturbative Chern--Simons free
energy}\index{$F^{\text{pert}}$ perturbative Chern--Simons free
energy} perturbative Chern--Simons free energy, and denote it by
$F^{\text{pert}}$. Using the coupling constant $x= {2\pi}/{(k+N)}$,
replace $k+N$ in $F^{\text{pert}}$, expand the logarithm in a series
($\ln(1-z)=-z-\frac12 z^2-\frac13 z^3-\cdots$) and use
$\zeta(s)=\sum_{n=1}^\infty n^{-s}$  to obtain
\begin{multline*}
F^{\text{pert}} {=} \sum_{j=1}^{N-1}(N{-}j)
  \biggl(\sum_{n=1}^\infty \ln\Bigl(1{-}
  \frac{x^2j^2}{4\pi^2n^2}\Bigr)\biggr)\\[-1ex]
= -\!\sum_{j=1}^{N-1}(N{-}j)\biggl(\sum_{n=1}^\infty
  \sum_{p=1}^\infty\frac{x^{2p}j^{2p}}{4^p\pi^{2p}pn^{2p}}\biggr)
{=} -\!\sum_{j=1}^{N-1}(N{-}j)\biggl(\sum_{p=1}^\infty
  \frac{x^{2p}j^{2p}}{4^p\pi^{2p}p}\zeta(2p)\biggr).
\end{multline*}
Recall from \fullref{pCS} on perturbative Chern--Simons theory
that we expect based on intuition from the path integral that the
free energy should have an interesting 't Hooft expression of the
form $F=\sum_{g}\sum_{h} x^{2g-2+h}N^hF_{g,h}$. This motivates our
next manipulations. Applying the two formulas derived in
\fullref{sumpower} gives
\begin{align*}
F^{\text{pert}}=&\sum_{p=1}^\infty\sum_{k=0}^p\begin{pmatrix}2p+1\\2k\end{pmatrix}
\frac{B_{2k}}{p(2p-2k+2)}(2\pi)^{-2p}\zeta(2p)N^{2p-2k+2}x^{2p}\\[-1ex]
\
&\qquad-\sum_{p=1}^\infty\frac{1}{2p}(2\pi)^{-2p}\zeta(2p)N^{2p+1}x^{2p}\\[-1ex]
\
&\qquad-\sum_{p=1}^\infty\sum_{k=0}^p\begin{pmatrix}2p+1\\2k\end{pmatrix}
\frac{B_{2k}}{p(2p+1)}(2\pi)^{-2p}\zeta(2p)N^{2p-2k+2}x^{2p}\\[-1ex]
\
&\qquad+\sum_{p=1}^\infty\frac{1}{2p}(2\pi)^{-2p}\zeta(2p)N^{2p+1}x^{2p}\\[-1ex]
=&\sum_{p=1}^\infty\sum_{k=0}^p
\begin{pmatrix}2p+1\\2k\end{pmatrix}\frac{(2k-1)B_{2k}}{p(2p+1)(2p-2k+2)}(2\pi)^{-2p}\zeta(2p)N^{2p-2k+2}x^{2p}\,.
\end{align*}
We now re-index the double sum setting $g=k$ and $h=2p-2k+2$. Notice that
\begin{multline*}
\begin{pmatrix}2p+1\\2k\end{pmatrix}\frac{2k-1}{p(2p+1)(2p-2k+2)}
=\begin{pmatrix}2p-1\\2p-2k+2\end{pmatrix}\frac{1}{k(2k-2)} \\[-1ex]
=\begin{pmatrix}2g+h-3\\h\end{pmatrix}\frac{1}{g(2g-2)},
\end{multline*}
for $g>1$. We conclude that
\begin{multline}\label{fpert}
F^{\text{pert}} {=} \sum_{g=2}^\infty
  \sum_{\stackrel{\scriptstyle h\ \text{even}}{h\ge 2}}
  \binom{2g{+}h{-}3}{h} \frac{B_{2g}}{g(2g-2)}(2\pi)^{2-2g-h}
  \zeta(2g{-}2{+}h)N^{h}x^{2g-2+h}  \\[-1ex]
+\!\sum_{\stackrel{\scriptstyle h\ \text{even}}{h\ge 2}}
  \frac{1}{6h}(2\pi)^{-h}\zeta(h)N^{h}x^{h}
  {-}\!\sum_{\stackrel{\scriptstyle h\ \text{even}}{h\ge 4}}
  \frac{2}{h(h{-}1)(h{-}2)}(2\pi)^{2-h}\zeta(h{-}2)N^{h}x^{h-2}.
\end{multline}
The second sum in the above expression is the $g=1$ term and the
third sum is the $g=0$ term. In the original paper of Gopakumar and
Vafa the last two sums were rewritten using polylogarithm identities
\eqref{Li3} and \eqref{Li1} (see Gopakumar and Vafa \cite[(3.7), (3.10)]{GV}).
We prefer to keep them as they are and simplify the polylogarithms
on the Gromov--Witten side instead.

The term
$$\sum_{j=1}^{N-1}(N-j) \ln\left(\frac{2\pi j}{k+N}\right)=
  \frac{N(N-1)}{2}\ln x+\sum_{j=1}^{N-1}(N-j)\ln j,$$
from the Chern--Simons free energy is called the Barnes term \cite{GV}. This is
because one has
$$\sum_{j=1}^{N-1}(N-j)\ln j =\ln\biggl(\prod_{j=1}^{N-1}
  (j!)\biggl)=\ln(G_2(N+1)).$$
We use the asymptotic expansion \eqref{Barnes} for the Barnes function
to analyze this term in the final comparison.

\subsection{The Gromov--Witten free energy}

Recall from \fullref{FGW} that the restricted Gromov--Witten
free energy is
\begin{multline}\label{GWdef}
\wwhat{F}^{GW}_{X_{S^3}}=\sum_{g=0}^\infty\sum_{d=1}^\infty
N_{g,d}y^{2g-2}e^{-td}\\[-3ex]
=\sum_{d=1}^\infty \frac{1}{d}\Bigl(2\sin
\frac{dy}{2}\Bigr)^{-2}e^{-td}
=\sum_{d=1}^\infty \frac{1}{4d}\csc^2 \Bigl(\frac{dy}{2}\Bigr)e^{-td},
\end{multline}
where $t$ is assumed to be in the right half-plane for the
series to converge.
\begin{remark}
The formula derived in \fullref{zetber} allows one to compute
the signs of the Bernoulli numbers. Using these signs one can check
that the following formulas that we give for the Gromov--Witten
invariants agree with the formulas from Faber and Pandharipande \cite{FP}.
\end{remark}
We can rewrite equation \eqref{GWdef} using the polylogarithms
\eqref{polylog} and \fullref{coth},
\begin{align}
\wwhat F^{GW}_{X_{S^3}}&=\sum_{d=1}^\infty\frac{1}{4d}\csc^2
(\frac{dy}{2})e^{-td}\notag\\
&=-\sum_{d=1}^\infty\frac{\partial}{\partial y}\left(\frac{1}{2d^2}\cot(\frac{dy}{2})\right)e^{-td} \notag\\
&=-\sum_{d=1}^\infty\frac{\partial}{\partial y}\left(\frac{i}{2d^2}\coth(\frac{idy}{2})\right)e^{-td} \notag\\
&=-\sum_{d=1}^\infty\frac{1}{d^3}\frac{\partial}{\partial y}\left(\frac{idy}{2}\coth(\frac{idy}{2})/y\right)e^{-td} \notag\\
&=-\sum_{d=1}^\infty\frac{1}{d^3}\frac{\partial}{\partial y}\left(\sum_{g=0}^\infty(-1)^g\frac{B_{2g}}{(2g)!}d^{2g}y^{2g-1}\right)e^{-td} \notag\\
&=\sum_{g=0}^\infty\sum_{d=1}^\infty\frac{1}{d^{3-2g}}(-1)^{g-1}\frac{(2g-1)B_{2g}}{(2g)!}y^{2g-2}e^{-td} \label{FPform}\\
&=\sum_{g=0}^\infty\left((-1)^{g-1}{(2g-1)B_{2g}}
\text{Li}_{3-2g}(e^{-t})/{(2g)!}\right)y^{2g-2}\label{liform}
\\
&=\sum_{g=0}^\infty \wwhat F_g^{X_{S^3}}y^{2g-2}\,.
\end{align}
Comparing \eqref{FPform} with the definition of the Gromov--Witten free energy \eqref{GWdef} we see that
$$
N_{g,d}(X_{S^3})={d^{2g-3}}(-1)^{g-1}(2g-1)\frac{B_{2g}}{(2g)!}\,.
$$
The last line in the above equation serves to define the genus $g$
contribution to the restricted Gromov--Witten energy. In general the
expansion of the free energy can be expressed as a sum of
polylogarithms. The result is described in the following exercise.
\begin{exm}
Expand the $p=0$ term in the following expression as we did above,
then express the sines in the remaining terms using exponentials and
simplify with the binomial formula
$$
\wwhat F^{GW}_{X}=\sum_{p=0}^\infty\sum_\beta\sum_{d=1}^\infty
n^p_\beta \frac{1}{d}\left(2\sin
\frac{dy}{2}\right)^{2p-2}e^{-d\langle t,\beta\rangle}\,.
$$
The answer you should get is
\begin{align*}
&\what{F}^{GW}(X)=
\sum_{g=0}^\infty\biggl(\sum_\beta
\Bigl(n^0_\beta(-1)^g\frac{(2g{-}1)B_{2g}}{(2g)!}\\
&+\!\sum_{p=1}^\infty\sum_{j=0}^{2p-2}n^p_\beta(-1)^{p+g}\binom{2p{-}2}{j}
\frac{(1{-}p{+}j)^{2g-2}2g(2g{-}1)}{(2g)!}\Bigr)\text{Li}_{3-2g}(e^{-\langle
t,\beta\rangle})\biggr)y^{2g-2}.
\end{align*}
\end{exm}
Combining the polylogarithm formula for the Gromov--Witten free
energy \eqref{liform} with the power series expansion of
$\text{Li}_{3-2g}$ for $g\ge 2$ \eqref{lipower} gives
\begin{multline}\label{Fg}
\wwhat{F}_g^{X_{S^3}}=\frac{B_{2g}}{2g(2g{-}2)}(it)^{2-2g}\\[-1ex]
+\frac{B_{2g}}{g(2g-2)}\sum_{\stackrel{\scriptstyle h\ \text{even}}{h\ge
0}} \binom{2g{+}h{-}3}{h}(2\pi )^{2-2g-h}\zeta({2g{+}h{-}2})(it)^h\,.
\end{multline}
We next need to consider the $g=1$ and $g=0$ terms. The
polylogarithm formula together with \fullref{Li13} gives
\begin{equation}\label{F1}
\wwhat F_1^{X_{S^3}}=t/24-\tfrac{1}{12}\ln t+\sum_{\stackrel{\scriptstyle h\
\text{even}}{h\ge 2}} \frac{1}{6h}(2\pi)^{-h}\zeta(h)(it)^{h}\,.
\end{equation}
and
\begin{multline}\label{F0}
\wwhat F_0^{X_{S^3}}=
  \zeta(3)-\zeta(2)t+3 t^2/4+t^3/12 -\frac{t^2}{2}\ln t\\[-1ex]
-\sum_{\stackrel{\scriptstyle h\ \text{even}}{h\ge 4}}
  \frac{2}{h(h-1)(h-2)}(2\pi)^{2-h}\zeta(h-2)(it)^{h}.
\end{multline}
Note that in the definition of the Gromov--Witten free energy
\eqref{GWdef} the sum over degrees begins with $d=0$ and we have not
included $d=0$. The degree $0$ are constant maps and there is a
question as to whether they should be included. The answer turns out
to be `yes'. Indeed, notice that the second term in \eqref{Fg} is
almost completely identical to the first term in \eqref{fpert}. The
only difference is that summation there starts at $h=2$ as opposed
to $h=0$. The extra $h=0$ term reads
$$\frac{B_{2g}}{g(2g-2)}(2\pi )^{2-2g}\zeta({2g-2}).$$
Using the result from \fullref{zetber} together with the fact
that $\chi(X_{S^3})=2$ we compute
$$
\frac{B_{2g}}{g(2g-2)}(2\pi
)^{2-2g}\zeta({2g-2})=\frac{(-1)^{g}(2g-1)B_{2g}B_{2g-2}\,\chi(X_{S^3})}{2(2g-2)(2g)!}\,.
$$
But for $g\geq2$ this is exactly the negative of the degree zero
invariants $N_{g,0}(X)$ computed in equation \eqref{Xdeg0}! Genus
$0$ and $1$ contributions line up as well. There are no contracted
genus zero or one stable curves fixed by the torus action so
$N_{0,0}=N_{1,0}=0$. Combining the formulas we get the degree zero
term
\begin{equation}\label{h0term}
N_{g,0} = -\frac{B_{2g}}{g(2g-2)}(2\pi )^{2-2g}\zeta({2g-2})
\end{equation}
that exactly cancels the extra $h=0$ term. In other words, a major
discrepancy between the Gromov--Witten and the Chern--Simons free
energies takes care of itself if we include the degree zero terms (as
we should have from the beginning).  Also note that the sum over genus
would diverge without the cancelation afforded by the degree zero term.
This is why the definition of the Gromov--Witten free energy includes
the constant terms.

As we explained ${F}_1=\wwhat F_1$, ${F}_0=\wwhat F_0$ and
\eqref{Fg} turns into
\begin{equation}\label{Fg^}
\begin{aligned}
{F}_g(X_{S^3})&=
  \frac{B_{2g}}{2g(2g-2)}(it)^{2-2g}\\
&+\frac{B_{2g}}{g(2g-2)}
  \sum_{\stackrel{\scriptstyle h \text{ even}}{h\ge 2}}
  \binom{2g+h-3}{h}(2\pi )^{2-2g-h}\zeta({2g+h-2})(it)^h.
\end{aligned}
\end{equation}
We are all set for the final comparison of the free energies on both sides of the duality.

\subsection{The final comparison}

As we warned in the introduction the match between the Gromov--Witten
and the Chern--Simons free energies will not be exact. The
discrepancy may be due to the fact that as physicists insist, we
should really consider the $\text{U}(N)$ not $\text{SU}(N)$
Chern--Simons theory which is expected to insert some additional
normalizing factors into the partition function (see Mari\~no
\cite{Mar,Mar2}).
Combining equations \eqref{FS3}, \eqref{fpert}) and the asymptotic
expansion for the Barnes term  \eqref{Barnes} gives the following
expression for the unnormalized Chern--Simons free energy,
\begin{align*}
F^{\text{CS}}_{S^3}(N,x)&=\tfrac12{N(N-1)}\ln
  x+\tfrac12(1-N)\ln(k+N)+\tfrac12{N^2}\ln N\\
&-\tfrac12\ln N-\tfrac34 N^2-\tfrac{1}{12}\ln
  N-\zeta^\prime(0)N+\zeta^\prime(-1)\\
&-\sum_{\stackrel{\scriptstyle h\ \text{even}}{h\ge 4}}
  \frac{2}{h(h-1)(h-2)}(2\pi)^{2-h}\zeta(h-2)N^{h}x^{h-2} \\
&+\sum_{\stackrel{\scriptstyle h\ \text{even}}{h\ge 2}}
  \frac{1}{6h}(2\pi)^{-h}\zeta(h)N^{h}x^{h}\\
&+\sum_{g=2}^\infty \sum_{\stackrel{\scriptstyle h\ \text{even}}{h\ge 2}}
  \binom{2g{+}h{-}3}{h}\frac{B_{2g}}{g(2g-2)}(2\pi)^{2-2g-h}
  \zeta(2g{-}2{+}h)N^{h}x^{2g-2+h}\\
&+\sum_{g=2}^\infty\frac{B_{2g}}{2g(2g-2)}N^{2-2g}.
\end{align*}
In the same way we combined terms to get the Chern--Simons free
energy, \fullref{fGW} and equations \eqref{F0}, \eqref{F1} and
\eqref{Fg^} give the following expression for the full Gromov--Witten
free energy of the resolved conifold,
\begin{align*}
F^{\text{GW}}_{X_{S^3}}(t,y)=&\tfrac{1}{24}t-\tfrac{1}{12}\ln t
  +\zeta(3)y^{-2}-\zeta(2)ty^{-2}+3 t^2y^{-2}/4\\
&+t^3y^{-2}/12 -\tfrac12 {t^2y^{-2}}\ln t\\
&-\!\sum_{\stackrel{\scriptstyle h\ \text{even}}{h\ge 4}}
  \frac{2}{h(h-1)(h-2)}(2\pi)^{2-h}\zeta(h-2)(it)^{h}y^{-2}\\
&+\!\sum_{\stackrel{\scriptstyle h\ \text{even}}{h\ge 2}}
  \frac{1}{6h}(2\pi)^{-h}\zeta(h)(it)^{h}\\
&+\!\sum_{g=2}^\infty\frac{B_{2g}}{g(2g-2)}\sum_{\stackrel{\scriptstyle h\
  \text{even}}{h\ge 2}} \!\binom{2g{+}h{-}3}{h}(2\pi
  )^{2-2g-h}\zeta({2g{+}h{-}2})(it)^hy^{2g-2} \\
&+\!\sum_{g=2}^\infty\frac{B_{2g}}{2g(2g{-}2)}(it)^{2-2g}y^{2g-2}.
\end{align*}
Note that some of the extra terms that appeared `on the Chern--Simons
side' in the original paper \cite{GV} show up `on the Gromov--Witten
side' with opposite signs in our presentation. This is because we
chose not to reexpress the genus $0$ and $1$ contributions in the
Chern--Simons free energy via the polylogarithm identities. By
inspection, under the substitution $it=\pm Nx$ and $y=x$ all the
infinite sums match exactly. In light of the complicated definitions
and expressions for the free energies this is a remarkable
coincidence. Notice that the sums represent exactly the perturbative
part of the Chern--Simons free energy and thus contain all of the
information about the perturbative invariants.

Analytically, we are comparing series expansions of two functions
near the origin $(t,y)=(0,0)$. It may seem odd that we should choose
the origin since $(t,y)=\bigl(\smash{\frac{2\pi iN}{k+N}},
\smash{\frac{2\pi}{k+N}}\bigr)$
converges to $(2\pi i,0)$ at large $N$. However, as one can see for
example
from \eqref{GWdef} the free energy is periodic in $t$ with period
$2\pi i$ so the coefficients are the same as the coefficients at the
origin. Another issue is that originally in \eqref{GWdef} we assumed
the real part of $t$ to be positive. The problem with analytically
continuing the free energies to a punctured neighborhood of the
origin is that the logarithmic terms in both expressions are
ambiguous. However, for us the free energies are just a bookkeeping
device for the invariants on both sides of the duality. Since
logarithmic and other mismatching terms outside the infinite sums
carry no apparent geometric information they do not pose a serious
problem. Finally, notice that the infinite sums are real-valued
which allows us to package the comparison into the following nice
form.
\begin{thm}
The full Gromov--Witten free energy and the unnormalized Chern--Simons
free energy are related by
$$\text{Re}\bigl(F^{\text{GW}}_{X_{S^3}}(iNx,x)-\,
F_{S^3}^{\text{CS}}(N,x)\bigr)=\tfrac{5}{12}\ln
x+\zeta(3)x^{-2}-\tfrac12\ln(2\pi)-\zeta^\prime(-1)\,.
$$
\end{thm}
\begin{proof}
Combining the expressions for the free energies gives,
\begin{multline*}
\text{Re}\bigl(F^{\text{GW}}_{X_{S^3}}(iNx,x)\bigr)-F_{S^3}^{\text{CS}}(N,x) =\\
\tfrac12 (N{-}1)\ln(k{+}N){-}\tfrac12{N(N{-}1)}\ln x{-}\tfrac12{N^2}\ln N
{+}3N^2/4{+}\tfrac{1}{12}\ln N{+}\zeta^\prime(0)N{-}\zeta^\prime(-1)\\
{-}\text{Re}\bigl(\tfrac{1}{12}\ln(ixN)\bigr) {+}\zeta(3)x^{-2}
{-}\tfrac34 N^2{+}\tfrac12{N^2}\text{Re}(\ln(ixN)).
\end{multline*}
Using the fact that $x=\frac{2\pi}{k+N}$ to write $\ln(k+N)$ and the
value $\zeta^\prime(0)=-\frac12\ln(2\pi)$ from Sondow \cite{sund} allows
one to combine like terms further to obtain the result.
\end{proof}
We conclude our presentation of the Gopakumar--Vafa duality with this
remarkable equality.
\begin{remark}
Large $N$ duality is said to be exact when the full free energies
are equal. It is said to hold to a leading order when the genus zero
terms agree. Given this the term $\zeta(3)x^{-2}$ in the comparison
theorem is slightly disturbing. This term prevents the duality from
holding at the level of the genus zero contributions.
In the physical theory this term is canceled by additional genus zero
corrections in degree zero.  Ooguri and Vafa \cite{OV1} obtained a perfect
agreement of the two sides using the \emph{physical normalization} of the
$S^3$ Chern--Simons free energy (which can only be computed on a case by
case basis comparing exact answers to perturbative expansions).  Thus, we
expect that one will obtain an exact correspondence after enough examples
have been computed to find a general form of the correct normalization.
\end{remark}

\section{New results (2003--2006)}\label{last}

In this section we describe some recent directions of research
inspired by the Large $N$ Duality and discuss some difficulties and
open problems encountered within them. Obviously this account is
biased by our background and interests, and we apologize in advance
for any inaccuracies and/or omissions. As in the history \fullref{02} the dates in the text refer to arxiv submissions while
references are given wherever possible to journal publications.

\subsection{Computations of the Gromov--Witten invariants}\label{last1}

Computational verification of the Gopakumar--Vafa Large $N$ Duality
depends largely on one's ability to compute the Gromov--Witten
invariants for as large a class of threefolds as possible. Toric
threefolds seem to be natural candidates to start with since
holomorphic torus actions are a part of their definition and the
full power of virtual localization can be applied. However, as the
sample computations in \fullref{N3} show, the complexity of
expressions obtained through virtual localization often grows very
rapidly with degree and genus and quickly becomes unmanageable.

Aganagic, Mari\~no and Vafa  introduced an interesting way to attack
this computation for local toric Fano surfaces \cite{AMV}. Iqbal
found a nice reformulation of these results \cite{Iqbal} and Zhou
gave a mathematical proof of the results. This Aganagic, Mari\~no
and Vafa paper led to a breakthrough by Aganagic, Klemm, Mari\~no
and Vafa in \cite{vertex}, where an effective algorithm was offered
that produces explicit combinatorial answers (without Hodge
integrals, etc) for all toric Calabi--Yau threefolds. The idea is
that any toric Calabi--Yau threefold (which is necessarily
non-compact) can be presented by a labeled planar trivalent graph
that can be cut into trivalent vertices `with legs' representing
$\C^3$ patches. Labels on the edges provide the gluing data that
specifies the threefold. The topological vertex is an explicit
function  $C_{\vec \mu,\vec n}(\lambda)$ of the edge labels at each
vertex, three partitions $(\mu^1,\mu^2,\mu^3)$ and three integers
$(n^1,n^2,n^3)$ associated to each vertex of the graph. The
generating function of the Gromov--Witten invariants of the
threefold can then be written as a `state sum' of these $C_{\vec
\mu,\vec n}(\lambda)$ taken over additional labelings. The authors
of \cite{vertex} provided an explicit combinatorial expression for
$C_{\vec \mu,\vec n}(\lambda)$ based on a derivation assuming large
$N$ duality.

This algorithm has been almost proved mathematically by
Li, Liu, Liu and Zhou \cite{Mvertex} based on gluing formulas for
relative Gromov--Witten
invariants. However, mathematical redefinition leads to a seemingly
different expression $\widetilde{C}_{\vec \mu,\vec n}(\lambda)$ for
the topological vertex. The equality
$$
C_{\vec \mu,\vec n}(\lambda)=\widetilde{C}_{\vec \mu,\vec
n}(\lambda)
$$
has been verified for the case when one of the partitions $\mu^i$ is
empty or when all partitions have length $\leq6$ but the general
case remains open.

Another interesting class of Calabi--Yau threefolds is given by local
curves. Those are the total spaces of rank 2 complex vector bundles
$N$ over a complex curve $\Sigma$ with $c_1(N)=2g(\Sigma)-2$, for example
the resolved conifold is a local $\CP^1$. Although not toric in
general these threefolds always admit `degenerate' holomorphic torus
actions that leave the entire zero section fixed (as opposed to
isolated points in the usual case). In \cite{locurv} J Bryan and R
Pandharipande used relative Gromov--Witten invariants and the TQFT
approach to construct a recursive algorithm that computes the
invariants of any local curve.

Next in complexity is the case of local surfaces, that is, total spaces
of canonical bundles $K_S$ to complex surfaces $S$. The work of
D-E Diaconescu, B Florea and others on the invariants of local
del Pezzo surfaces \cite{dP} culminated in the joint work with N
Saulina \cite{Lvertex} that extends the toric topological algorithm
to the case of local ruled surfaces with a finite number of
reducible fibers. As in the case of local curves the authors make
use of degenerate torus actions that fix finitely many curves and
augment the toric formalism by the corresponding correction terms.
The derivation of the combinatorial formulae uses physical arguments
as in \cite{vertex} and mathematical justification of the ruled
vertex is an open problem. Another open problem is to generalize
this algorithm to arbitrary Calabi--Yau threefolds with degenerate
torus actions.

From the point of view of Large $N$ Duality it also important to
understand the pre-duals of the above threefolds, that is, analogs of
$T^*S^3$ for the resolved conifold and identify the correct pre-dual
theories. In the known examples as originally considered by Aganagic and
Vafa \cite{AV} several $2$--cycles are collapsed and then replaced by
Lagrangian $3$--cycles via resolving a singular deformation. These
pre-duals are thus of a more general form than $T^*M$ \cite{AMV,dP}.
The corresponding theories combine elements of both the Chern--Simons
and the Gromov--Witten theories in agreement with Witten's original
idea that the Chern--Simons theory on $M$ is the correct
`Gromov--Witten theory' on $T^*M$ (see Witten \cite{Wcss} and Grassi and
Rossi \cite[Appendix 9]{GR}). It was the formalism from \cite{AMV} that led to the
discovery of the topological vertex.

\subsection{Intermediate theories}\label{last2}

In this section we discuss some theories that have recently emerged
and could provide a bridge between the two sides of the large $N$
duality. The most developed of these is the Donaldson--Thomas theory
that has already been used to prove some of the duality's structural
predictions.

Gopakumar and Vafa predicted in \cite{GV} that the properly
normalized partition function of the Gromov--Witten invariants
$Z_X(\lambda,v)$ on a Calabi--Yau threefold $X$ is a rational
function of the variable $q=-e^{i\lambda}$ and expands into a series
in $q$ with integral coefficients (BPS states or Gopakumar--Vafa
invariants). Intuitively the integers should correspond to counts of
embedded curves in $X$. Classically, embedded curves are described
by ideal sheaves on $X$, that is, torsion-free rank one sheaves  with
trivial determinants. S Donaldson and R Thomas introduced in
\cite{DT} a new class of invariants $\widetilde{N}_{\chi,\beta}$
that count the number of ideal sheaves with a given holomorphic
Euler characteristic $\chi$ and the associated curve class $\beta$.
The Donaldson--Thomas theory is `better' than the Gromov--Witten
theory in the sense that no orbifolds occur as moduli spaces and the
numbers $\widetilde{N}_{\chi,\beta}$ are integers. In \cite{MNOP1,MNOP2} D
Maulik, N Nekrasov, A Okounkov and R Pandharipande conjectured
that the Donaldson--Thomas partition function
$\widetilde{Z}_X(q,v)=\sum_{\chi,\beta}\widetilde{N}_{\chi,\beta}q^\chi
v^\beta$ turns into $Z_X(\lambda,v)$ after the change of variables
$q=-e^{i\lambda}$. This automatically implies the rationality and
integrality predictions for all Calabi--Yau threefolds. Moreover, the
authors offered the `equivariant vertex' algorithm analogous to the
topological vertex for computing the Donaldson--Thomas invariants.
Combined with \cite{Mvertex} and \cite{locurv} their result proves
the Gromov--Witten/Donaldson--Thomas duality for toric Calabi--Yau
threefolds and local curves respectively. Obviously it is desirable
to extend the duality to local surfaces and more general threefolds.

It is expected that the Donaldson--Thomas theory admits a
gauge-theoretic interpretation and if so it could serve as a link in
a chain mathematically connecting Gromov--Witten theory to
Chern--Simons theory. In particular there are some promising
connections discovered between the BPS states and Yang--Mills theory in
two dimensions \cite{bh2D} and in four dimensions \cite{bh4D}.

Another possible intermediary is the symplectic field theory (SFT)
of Y Eliashberg, A Givental and H Hofer \cite{SFT}. In general
SFT studies invariants of a contact manifold $\mathcal{C}$ by
considering moduli of pseudoholomorphic curves in its
symplectization $\mathcal{C}\times\R$. One can naturally associate a
contact manifold to any $3$--dimensional manifold $M$, namely the
cosphere bundle $S(T^*M)$ with symplectization
$S(T^*M)\times\R\simeq T^*M\backslash M$. Unlike $T^*M$ itself
$T^*M\backslash M$ does admit nontrivial pseudoholomorphic curves
that may serve as Witten's `instantons at infinity' \cite{Wcss}.
Another attractive trait of SFT is that it reconstructs some knot
invariants, for example the Alexander polynomial, from the Gromov--Witten
invariants (see L Ng's paper \cite{ng} in this volume). The main challenge in
applying SFT to Large $N$ Duality is the scarcity of effective
algorithms for computing the invariants.

We should also mention an older approach to proving large $N$
duality suggested by B Acharya \cite{Ach} and developed by
M Atiyah, J Maldacena and C Vafa \cite{Mflop} (see also Grassi and
Rossi \cite{GR}) by lifting both sides to the
M--theory on a $7$--dimensional manifold with $G_2$ holonomy.
However, so far the M--theory approach (M for mystery) has not been
very fruitful mathematically because the geometry of $G_2$ manifolds
is much less understood than that of Calabi--Yau threefolds.

\subsection{Construction of large $N$ duals}\label{last3}

There are very few known large $N$ dual pairs despite the large
number of known Calabi--Yau threefolds with computable Gromov--Witten
invariants. For a while after 1998 the original
$T^*S^3/\calO(-1)\oplus\calO(-1)$ example remained the only one. In
2001 F Cachazo, K Intriligator and C Vafa  constructed a family of
examples \cite{CIV} with the deformed conifold $T^*S^3$ replaced by
the following hypersurface in $\C^4$:
$$
W'(x)^2+f(x)+y^2+z^2+w^2=0.
$$
Here $W(x)$, $f(x)$ are polynomials of degrees $n$, $n-2$
respectively and the deformed conifold is recovered for
$W(x)=x^2/2$, $f=\text{const}$. However, they do not provide new
examples of the form $T^*M$, which is where it is easiest to compute
the Chern--Simons side. In fact, the only known examples of this form
come from the spherical quotients.

The first published version of Large $N$ Duality for the lens spaces
$L(p,1)=S^3/\Z_p$ appeared in Giveon, Kehagias and Partouche\cite{GKP}
(see also Mari\~no \cite{matrix} for other credits). The idea is to extend
the group $\Gamma=\Z_p$ action to $T^*S^3$, pull it through the transition
to $\calO(-1)\oplus\calO(-1)$ and then resolve the resulting quotient
singularity (see Halmagyi, Okuda and Yasnov \cite{HOY} for details). The dual quotient is a
bundle over $\CP^1$ fibered by surface singularities $\C^2/\Gamma$
and one can obtain a threefold resolution by resolving the surface
singularities in each fiber. For the resolved quotient to be
Calabi--Yau the resolution must be crepant (see Harris \cite{AG1}). This restricts
the list of groups $\Gamma$ to finite subgroups of SU($2$) (in
particular this is why the lens spaces $L(p,q)$ with $q\neq 1$ are
out). The geometry of the transition for such quotients and more
general fibrations of surface ADE--singularities that also include
the \cite{CIV} examples is studied in \cite{ADE}. One may be able to
lift the SU($2$) restriction by considering orbifold Gromov--Witten
invariants. The crepant resolution conjecture proved for the affine
ADE--singularities by F Perroni \cite{crep} says roughly that the
orbifold Gromov--Witten invariants of a variety are equal to the
ordinary Gromov--Witten invariants of its crepant resolution (when
such exists). Thus conjecturally Large $N$ Duality for the quotients
with $\Gamma\subset\text{SO}$(4) should relate the Chern--Simons
invariants of $S^3/\Gamma$ to the orbifold Gromov--Witten invariants
of $\calO(-1)\oplus\calO(-1)/\Gamma$.

Computationally only the $L(p,1)$ example has been considered so far
\cite{matrix,HOY}. Here the Chern--Simons partition function is known
(see for example Hansen and Takata \cite{HT}, Garoufalidis and Mari\~no \cite{GM}
and Turaev \cite{T}) and the Gromov--Witten free energy is
computable via the topological vertex since the dual is toric. The
duality statement is more complicated than in the $S^3$ case and one
has to split the Chern--Simons partition function into contributions
from different (gauge classes of) flat connections before comparing
to the Gromov--Witten side \cite{matrix}.  In addition to the
analysis in \cite{matrix}, there is an indirect check of the duality
based on mirror symmetry in \cite{HOY}. In general, the precise
meaning of the duality for more general manifolds remains an open
question.

This justifies interest in explicit computations of the LMO
invariant (see Le, Murakami and Ohtsuki \cite{LMO} and Bar-Natan,
Garoufalidis, Rozansky and Thurston \cite{bar}) that is believed to capture the
contribution of the trivial connection into the full Chern--Simons
partition function. M Mari\~no showed in \cite{Smatrix} using
physical considerations related to mirror symmetry that the LMO
invariant can be expressed as a perturbed Gaussian matrix integral
for Seifert fiber spaces. This reduces the computation to a solvable
matrix model (see Aganagic, Klemm, Mari\~no and Vafa \cite{matrix},
Fiorenza and Murri \cite{MPROP} and Mari\~no \cite{Mar2}). The paper \cite{Smatrix}
also presented some sample computations. Recently S Garoufalidis
and M Mari\~no \cite{GM} gave a mathematical derivation of the
matrix model for general rational homology spheres based on the
{\AA}rhus integral presentation of the LMO invariant \cite{bar}. An
interesting open question is just how much of the Gromov--Witten
theory on the dual can be recovered from the LMO invariant.

There is also a duality involving $SO(N)$ or $Sp(N)$ Chern--Simons
theories discussed in the physics literature \cite{3,4,5,6}.

Finally we mention the symplectic surgery approach of I Smith and
R Thomas to constructing large $N$ duals \cite{STY}. Their work
suggests that in many cases such duals ought to be `non-K\"ahler
Calabi--Yaus', that is, non-K\"ahler symplectic manifolds $X$ with
$c_1(X)=0$. If so, this explains why so few duals to cotangent
bundles are known despite the abundance of known K\"ahler Calabi--Yau
threefolds.

\newpage
\part{Appendices}
\setobjecttype{Part}

\appendix
\section{Stacks}
\setobjecttype{App}
\label{app:a}

In this section we follow the excellent exposition covering stacks
in Metzler \cite{Met}, and just add a couple of motivating examples. Stacks
were introduced to encode the structure of an orbifold in the
category of schemes, but may also be used to define orbifolds in the
smooth, topological and analytic categories. We will provide
examples in the smooth category. The main geometric objects that we
consider may be represented by a category with additional structure,
a contravariant functor or a covariant functor. We will conclude
this section with a list of properties that such a covariant functor
satisfies if and only if it comes from an object of the category in
a natural way. The generalized objects that are used in
Gromov--Witten theory are just covariant functors that satisfy  a
subset of these properties.

Recall that an orbifold is a space locally modeled on the quotient
of $\R^n$ by a finite group action together with additional data to
measure the stabilizer subgroups. As a first example, consider the
quotient of $\C$ by the natural left action of $\Z_3$ by
multiplication by cube roots of unity. We can encode this as a
category $X$ with objects Ob$(X):=\C$ and arrows
Ar$(X):=\Z_3\times\C$ where we consider
$(\omega,z)\in\text{Mor}(z,\omega z)$. The underlying space of the
associated orbifold is the quotient of the objects obtained by
identifying those objects connected by a morphism. The stabilizer
group of a point $z\in\text{Ob}(X)$ is just Mor$(z,z)$. This is a
special category because every morphism has an inverse. Such a
category is called a groupoid. In fact this has the structure of a
smooth groupoid. The structure maps in this example are given by
(source -- $s\co \text{Ar}(X)\to\text{Ob}(X)$, $s(\omega,z)=z$;
target -- $t\co \text{Ar}(X)\to\text{Ob}(X)$, $t(\omega,z)=\omega
z$; inverse -- $i\co \text{Ar}(X)\to\text{Ar}(X)$,
$i(\omega,z)=(\omega^{-1},\omega z)$; composition -- $m\co
\text{Mor}(\omega z,\theta\omega z)\times\text{Mor}(z,\omega z)
\to\text{Mor}(z,\theta\omega z) $, $m(f,g)=f\circ g$)
\begin{defn}
A smooth groupoid is a category $X$ with invertible morphisms such
that Ob$(X)$ and Ar$(X)$ are smooth manifolds and the various
structure maps are smooth.
\end{defn}
One good example to keep in mind is the smooth groupoid associated
to any atlas on a smooth manifold. Given an atlas,
$\calA=\{\varphi_\alpha\co U_\alpha\to V_\alpha\}$ define a smooth
groupoid with Ob$(X^\calA):=\coprod_\alpha V_\alpha$ and
Ar$(X^\calA):=\coprod_{\alpha, \beta} \varphi_\alpha(U_\alpha\cap
U_\beta)$ with the obvious structure maps and smooth structures.
\begin{exm}
Combine the example of the $\Z_3$ quotient of $\C$ with the smooth
groupoid associated to an atlas to define a smooth groupoid modeling
an orbifold with underlying space homeomorphic to $S^2$, one point
with stabilizer $\Z_3$, one point with stabilizer $\Z_2$ and the
rest of the points having trivial stabilizer.
\end{exm}

We now turn to the second way to encode a geometric object -- a
contravariant functor. Let ${\mathbf {DIFF}}$ be the category of
smooth manifolds and ${\mathbf {SET}}$ be the category of sets.
Given a smooth manifold $M$, we define a contravariant functor
$\underM\co {\mathbf {DIFF}}\implies {\mathbf {SET}}$ by
$\underM(N):=C^\infty(N,M)$ and $\underM(f\co N\to
P):=f^*\co C^\infty(P,M)\to C^\infty(N,M)$. It is possible to
reconstruct the original manifold (up to diffeomorphism) from the
associated functor. We therefore think of a contravariant functor
$\calM \co {\mathbf {DIFF}}\implies {\mathbf {SET}}$ as a generalized
manifold. (Recall that we described the moduli stack as a
contravariant functor from ${\mathbf {SCHEME}}$ to ${\mathbf
{SET}}$.) However, we will have to add some restrictions in order to
have a reasonable family of generalizations. When we add these
restrictions to give the formal definition of a stack we will use a
third description of geometric objects. This third description will
generalize the first two frameworks.

Given a contravariant functor $\calM \co {\mathbf {DIFF}}\implies
{\mathbf {SET}}$ one can define a new category ${\mathbf D}^\calM$
with
\begin{align*}
\Ob({\mathbf D}^\calM)&:=\textstyle{\coprod_{N}}\calM(N) \\
\text{and}\quad
\Mor((\alpha,N),(\beta,P))&:=\{a\in C^\infty(N,P)|\calM(a)(\beta)=\alpha\}.
\end{align*}
One then defines a
covariant functor  $F^\calM\co  {\mathbf D}^\calM\implies {\mathbf
{DIFF}}$ by $F^\calM(a\co  (\alpha,N)\to (\beta,P)):= a\co N\to P$. When
the contravariant functor is of the form ${\underM}$ the
associated covariant functor will satisfy a number of special
properties. The first property that it will satisfy is that it will
be a fibered category.
\begin{defn}
A fibered category over ${\mathbf C}$ is a covariant functor
$F\co {\mathbf D}\implies {\mathbf C}$ such that
\begin{enumerate}
\item For every $f\co C_0\to C_1\in\text{Ar}({\mathbf C})$ and
every $D_1$ such that $F(D_1)=C_1$ there is an arrow $g\co D_0\to D_1$
such that $F(g\co D_0\to D_1)=f\co C_0\to C_1$.
\item If $F(g_1)\circ f=F(g_0)$, then there is a unique $g\in\text{Ar}({\mathbf D})$ such that $F(g)=f$.
\end{enumerate}
\end{defn}

\begin{exm}
Check that $F^{\underM}\co  {\mathbf D}^{\underM}\implies
{\mathbf {DIFF}}$ is a fibered category.
\end{exm}

It is also possible to construct a fibered category associated to a
smooth groupoid. Let $X$ be a smooth groupoid and define a category
${}^X\!{\mathbf D}$ with
\begin{align*}
&\text{Ob}(^X{\mathbf D}):=\coprod_{N\in\text{Ob}({\mathbf {DIFF}})}
  C^\infty(N,\text{Ob}(X)), \\
&\text{Mor}(f\co N\to\text{Ob}(X),g\co P\to\text{Ob}(X)) := \\
&\qquad\bigl\{(\varphi\co N\to P,h\co N\to\text{Ar}(X))~|~s(h(p))=f(p)
  \text{ and } t(h(p))=g(\varphi(p)) \bigr\}.
\end{align*}
The associated covariant functor is given by,
$$
{}^X\!F(\varphi\co N\to P,h\co N\to\text{Ar}(X)):=\varphi\co N\to P\,.
$$

\begin{exm}
Check that ${}^X\!F\co  {}^X\!{\mathbf D}\implies {\mathbf {DIFF}}$ is a
fibered category. (Recall that the arrows in a groupoid have
inverses.)
\end{exm}

Let $F\co {\mathbf D}\implies {\mathbf C}$ be a fibered category and
$C$ be an object of ${\mathbf C}$. We can define a fibered category
over $C$, denoted $F_C$, by Ob$(F_C):=\{D\in\text{Ob}({\mathbf
D})|F(D)=C\}$ and Ar$(F_C):=\{\varphi\co D_0\to
D_1|F(\varphi)=\text{id}_C\}$. One can see that $F_C$ is a groupoid.

\begin{exm}
Prove that $F_C$ is a groupoid.
\end{exm}

One often wishes to put additional structure on a category. A good
motivating example is the category of open sets of a topological
space with inclusions as arrows. In this setting, one would like to
axiomatize the properties of open covers of the original space. This
leads to the notion of a Grothendieck topology and the notion of a
site. Once a category has a notion of coverings one can define an
analogue of a sheaf. This is one way to introduce stacks. We do not
need the definitions of Grothendieck topologies or sites, but we do
use the following definition that we quote from Metzler \cite{Met} to encode
the notion of a covering.

\begin{defn}\label{def:basis}
  A  {basis} for a Grothendieck topology on a category
  $\CC$ is a function $K$ which assigns to every object $C$ of $\CC$
  a collection $K(C)$ of families of arrows with target $C$,
  called  {covering families}, such that
  \begin{enumerate}
\item if $f\co  C' \to C$ is an isomorphism, then $\{f\}$ is a covering family;
\item (stability) if $\{f_{i}\co  C_{i} \to C \}$ is a covering family, then for
      any arrow $g\co  D \to C$, the pullbacks $C_{i} \times_C D$ exist and
      the family of pullbacks $\pi_{2}\co  C_{i} \times_C D \to D$
      is a covering family (of $D$);
\item (transitivity) if $\{f_{i}\co  C_{i} \to C \: |\: i \in I\}$
      is a covering family and for each $i \in I$, one has a
      covering family $\{g_{ij}\co  D_{ij} \to C_{i} \: | \: j \in I_{i}\}$,
      then the family of composites
      $\{f_{i} g_{ij}\co  D_{ij} \to C \: | \: i \in I, j \in I_{i}\}$
      is a covering family.
\end{enumerate}
\end{defn}

In addition to the usual notion of an open cover in the category of
open sets, we obtain an example of a basis for a Grothendieck
topology on ${\mathbf {DIFF}}$ by considering collections of open
embeddings whose images cover a given manifold. This example will be
used in our definition of prestack and stack.

We now have two different constructions of special fibered
categories. One starts with a manifold and passes to an associated
contravariant functor and then to the associated fibered category.
The second passes from a manifold to a smooth groupoid to the
associated fibered category. The fibered categories constructed in
either of these ways satisfy additional conditions summarized in the
definition of a prestack. We take the characterization of a prestack
given in \cite[Lemma 25]{Met} as the definition of a prestack.

\begin{defn}
A fibered category $F\co \CC\implies \DD$ is a prestack if and only if
there is a unique arrow $\psi \co x \to y$ filling in the dotted arrow
in every diagram of the following form with $F(\psi) = 1$. In this
diagram $C \in \text{Ob}(\CC)$, $\{C_{\alpha} \to C\}$ is a cover of
$C$ and $x,y$ are objects of $\DD$ that map to $C$ and
$x_\alpha\in\DD$ map to $C_\alpha$ under $F$. In addition we denote
fibered products as $C_{\alpha\beta}:=C_\alpha\times_C C_\beta$
($x_{\alpha\beta}$).
$$\xymatrix@=5pt{
     & & x_{\beta} \ar[dr] \ar[ddd] &
      & & & & & & C_{\beta} \ar[dr] \ar@{=}[ddd] & \\
  x_{\alpha \beta} \ar[urr] \ar[dr] \ar[ddd] & & & x \ar@{.>}[ddd]^{\psi}
    & & & & C_{\alpha \beta} \ar[urr] \ar[dr] \ar@{=}[ddd] & & & C \ar@{=}[ddd] \\
     & x_{\alpha} \ar[urr] \ar[ddd] & & & \ar@{~>}[rr]^{F}
       & & & & C_{\alpha} \ar[urr] \ar@{=}[ddd] & & \\
     & & y_{\beta} \ar[dr] &
      & & & & & & C_{\beta} \ar[dr] & \\
  y_{\alpha \beta} \ar[urr] \ar[dr] & & & y
   & & & & C_{\alpha \beta} \ar[urr] \ar[dr] & & & C \\
     & y_{\alpha} \ar[urr] & &
     & & & & & C_{\alpha} \ar[urr] & & }$$
\end{defn}

\begin{exm}
Check that $F^{\underM}$ and ${}^X\!F$ are prestacks.
\end{exm}

We finally come to the definition of a stack. The definition we give
is not the usual one. It is instead the characterization given in
Metzler \cite[Lemma 32]{Met}. We chose to use this as the definition
because it was the quickest way to a clean definition. For the more
typical description of a stack in terms of descent data see the full
exposition in \cite{Met} or Behrend \cite{Beh}.

\begin{defn}
 A prestack $F\co  \DD \implies \CC$ is a stack
 if and only if
  for every cover $\{C_{\alpha} \to C\}$ in $\CC$
  and $x_{\alpha}$ mapping to $C_\alpha$ satisfying the following commutative diagram for all index triples
$$
\xymatrix@=5pt{
    & x_{\alpha \beta} \ar[r] \ar[dr] & x_{\alpha}
      & & & & & C_{\alpha \beta} \ar[r] \ar[dr] & C_{\alpha} \\
    x_{\alpha \beta \gamma} \ar[ur] \ar[r] \ar[dr]
       & x_{\alpha \gamma} \ar[ur] \ar[dr] & x_{\beta}
      & \ar@{~>}[rr]^{F} & & & C_{\alpha \beta \gamma} \ar[ur] \ar[r] \ar[dr]
       & C_{\alpha \gamma} \ar[ur] \ar[dr] & C_{\beta} \\
    & x_{\beta \gamma} \ar[ur] \ar[r] & x_{\gamma}
      & & & & & C_{\beta \gamma} \ar[ur] \ar[r] & C_{\gamma}
  }
$$
there is an object $x$ in $\DD$ mapping to $C$
  and arrows $\{x_{\alpha} \to x\}$ filling in the commutative diagram
  $$
    \xymatrix@=5pt{
        & x_{\alpha} \ar[dr] & & & & & & C_{\alpha} \ar[dr] &  \\
      x_{\alpha \beta} \ar[ur] \ar[dr] & & x
       & \ar@{~>}[rr]^{F} & & & C_{\alpha \beta} \ar[ur] \ar[dr] & & C \\
        & x_{\beta} \ar[ur]  & & & & & & C_{\beta} \ar[ur] &
    }
  $$
\end{defn}

\begin{exm}
Prove that $F^{\underM}$ is a stack.
\end{exm}

One may be expecting an exercise to prove that ${}^X\!F$ is a stack.
However, the fibered category associated to a typical groupoid is
not a stack. A good example to consider is the groupoid that arises
from the standard two-chart atlas of $\CP^1$. In this case the
smooth groupoid has Ob$(X)=\C\disj \C$ corresponding to the two
charts and Ar$(X)=\C\disj\C\disj\C^\times\disj\C^\times$
corresponding to the overlaps. The objects of the domain category of
the associated fibered category are smooth maps from smooth
manifolds into Ob$(X)$. The problem arises when one considers a map
into $\CP^1$ with image that is not contained in one of the
coordinate charts, for example $x\co \C\to\CP^1$ given by
$x(z):=[z-1:z+1]$. By restricting to the charts and their overlaps
we obtain objects of the domain of the fibered category
$x_\pm\co \C-\{\mp 1\}\to\CP^1$ given by $x_\pm(z):=\smash{\frac{z\mp 1}{z\pm
1}}$. These objects map to the charts under the fibered category and
the charts form a cover of $\CP^1$. The pull-backs of the cover and
the $x_\pm$ satisfy the diagram in the hypothesis of the definition
of a stack but not the required extension property. It is possible
to stackify a prestack by declaring objects to be equivalence
classes of objects over elements of covers as in the hypothesis of
the definition. This is similar to the sheafification of a presheaf.
See Metzler \cite{Met} for details. When the prestack associated to the
groupoid associated to an atlas is stackified, the domain category
has objects that correspond exactly to smooth maps into the given
manifold. Thus the constructions of ${}^X\!F$ and $F^{\underM}$
agree when one starts with a smooth manifold and stackifies.
\begin{exm}
Prove that the pull-backs of the cover and the $x_\pm$ satisfy the
diagram in the hypothesis of the definition of a stack.
\end{exm}

We now need to consider maps of stacks. As motivation consider what
can be constructed between the stacks associated to a pair of
manifolds from a map between the manifolds, $f\co M\to N$. Recall that
the objects of the domain category of the stack associated to $M$
are just maps $\alpha\co P\to M$. Such a map can be taken to
$f\circ\alpha\co P\to N$, which is an object of the domain category
associated to $N$. This extends to arrows in the natural way to give
a covariant functor $A_f\co \DD^{\underM}\implies\DD^{\underN}$. This motivates the definition of a map between stacks given
below.
\begin{defn}
A map between stacks $F\co \DD\implies\CC$ and $G\co {\mathbf
E}\implies\CC$ say $A\co F\to G$ is just a covariant functor
$A\co \DD\implies{\mathbf E}$ such that $F=G\circ A$.
\end{defn}
There are some technical issues that arise when one wishes to define
isomorphism of stacks, which are best addressed with $2$--categories,
see \cite{Met}.

We now quote some definitions of some properties of stacks and maps
of stacks from \cite{Met}.
\begin{defn}\label{def:monoepi}
  Let $A\co  F' \to F$ be a map of fibered categories over $\CC$.
  We say $A$ is a  {monomorphism}
  if, for every object $C$ of $\CC$, the functor $A_{C}\co  F'_{C} \to F_{C}$
  on fibers is fully faithful.
\end{defn}

\begin{defn}
  We say a map of fibered categories
  $A\co F' \to F$
  is  {covering} (French `couvrant') if, for every object
  $C$ of $\CC$, and every object $D$ of $F_{C}$,
  there is a covering family $\{f_{i}\co C_{i} \to C \}$
  and for every $i$ an object $D_{i}$ of $F'_{C_{i}}$ such that
  $A(D_{i}) \cong D|_{C_{i}}$.

  If $F'$ and $F$ are stacks, then we refer to a covering map
  as an  {epimorphism}.
\end{defn}

For the next definitions we need the definition of the fibered
product or pull-back of a pair of maps of stacks. Given $A\co F\to H$
and $B\co G\to H$ one defines the pull-back stack $F\times_HG$ to be
the stack with domain category having objects
$\text{Ob}(\DD_F)\times_{\text{Ob}(\DD_H)}\text{Ob}(\DD_G)$, that
is, ordered pairs $(x,y)$ such that $A(x)=B(y)$. We can generalize
the construction of a contravariant functor ${\underM}$
associated to a smooth manifold to any object $C$ in any category
$\CC$. If the category is reasonably well behaved the associated
covariant functor $F^{\underC}$ will be a stack.

The following definition generalizes the notions of manifold and a
submersion between smooth manifolds.
\begin{defn}\label{def:repmap}
A stack $F$ is representable if and only if it is isomorphic to a
stack of the form $F^{\underC}$ for some object $C$ of the base
category $\CC$.
  Let $F, G$ be stacks over $\CC$. We say that a map $A\co F \to G$
  is  {representable} if for every object $C$ of $\CC$ and map $B\co
F^{\underC} \to G$,
  the pull-back stack $F^{\underC} \times_{G} F$ is representable.
\end{defn}
In fact if $f\co M\to N$ is a smooth map, the associated map of stacks
$A_f\co F^{\underM}\to F^{\underN}$ is representable if and
only if $f$ is a submersion, see Metzler \cite{Met}.

\begin{defn}\label{def:locallyrep}
  Let $F$ be a stack over $\CC$. We say $F$ is  {locally representable}
  if there is an object $C$ of $\CC$ and a representable epimorphism
  $A\co  F^{\wwbar C} \to F$.
\end{defn}

\begin{defn}\label{def:propertiesofrepmaps}
  Let $P$ be a property of maps in $\CC$ that is stable under
  pullback. We say that a representable map $A\co F \to G$
   {has property $P$} if, for every object $C$ of $\CC$ and
  map $B\co  F^{\underC} \to G$, the projection $B^{*}A\co  F^{\underC} \times_{G} F \to C$ has property $P$.
\end{defn}

We now come to the result characterizing representable stacks, see \cite{Met}.
\begin{thm}
A stack $F$ is equivalent to a stack of the form $F^{\underM}$ if and only if
\begin{description}
\item[A1] The stack $F$ is locally representable by a map $A\co F^{\underN}\to F$.
\item[A2] The map $\Delta\co F\to F\times F$ is proper.
\item[DM] The map $A$ is \'etale.
\item[R1] The stack has trivial automorphisms.
\item[R2] The map $\Delta\co F\to F\times F$ is a closed embedding.
\end{description}
\end{thm}

We now come to the definition of a Deligne--Mumford stack, and the conclusion of this appendix.
\begin{defn}
An Artin stack is a stack that satisfies {\bf A1} and {\bf A2}. A
Deligne--Mumford stack is an Artin stack that satisfies {\bf DM}
(usually assumed to be over the category ${\mathbf {SCHEME}}$). An
orbifold is a Deligne--Mumford stack over the category ${\mathbf
{DIFF}}$.
\end{defn}

\begin{exm}
Prove that the example of the orbifold with underlying space $S^2$
and two non trivial stabilizers $\Z_3$ and $\Z_2$ really is an
orbifold. Analyze the automorphisms of this stack.
\end{exm}

\section{Graph contributions to $N_2$}
\setobjecttype{App}
\label{app:b}

In this appendix we list the contributions of all of the fixed point
components to the localization computation of $N_2$. We label graphs
according to the conventions depicted in \fullref{gtype}. The
first contribution is:
$$
I(01)=-32 (\alpha_0 - \alpha_1)^{-2} (\alpha_0 - \alpha_2)^4
(\alpha_1 - \alpha_2)^{1} (\alpha_0 + \alpha_1 - 2\alpha_2)^{-1}\,.
$$
The contribution from $I(02)$ can be obtained by exchanging $1$ and
$2$ in the above expression. The next contribution is:
$$
\II(010)=8 (\alpha_0 - \alpha_1)^{-2} (\alpha_0 - \alpha_2)^3
(\alpha_1 - \alpha_2)^{-1}\,.
$$
The next is:
$$
\II(101)=8 (\alpha_0 - \alpha_1)^{2} (\alpha_0 - \alpha_2)^4
(\alpha_1 - \alpha_2)^{-2}\,.
$$
The contributions from $\II(020)$ and $\II(202)$ can be obtained by
exchanging $1$ and $2$ in the above expressions. The next
contribution is:
$$
\II(012)=(\alpha_0 - \alpha_1)^{-1} (\alpha_0 - \alpha_2)^3 (\alpha_1
- \alpha_2)^{1}(2 \alpha_1 - \alpha_0 - \alpha_2)^{-1}\,.
$$
The final contribution is:
$$
\II(102)=- (\alpha_0 - \alpha_1)^{1} (\alpha_0 -
\alpha_2)^{-1}(\alpha_1 - \alpha_2)^{-2}(2 \alpha_0 - \alpha_1 -
\alpha_2)^4\,.
$$
A tedious simplification gives
$$
N_2=I(01)+I(02)+\II(010)+\II(020)+\II(101)+\II(202)+\II(012)+\II(102)=1\,.
$$
This is a remarkable check of the localization formula.

\section{Quantum invariants from skein theory}
\setobjecttype{App}
\label{app:c}

When Jones introduced his polynomial invariant, he was motivated by
representations of the braid group arising from operator algebras.
Shortly thereafter the theory of quantum groups started taking off,
and Reshetikhin and Turaev started work on defining link and
$3$--manifold invariants based on quantum groups. The papers of
Witten provided additional inspiration and helped them devise their
invariants. Since link invariants only have to be invariant under
Reidemeister but not Kirby moves, a modular structure is not
required and one can get by using ribbon categories that are much
easier to construct. For instance, type I representations of
$U_q(\Sl_N\C)$ with $q$ a free variable rather than a fixed number
form a ribbon category. We denote the resulting colored link
invariants by $\smash{W^{{\mathfrak sl}_N}_\Lambda}$, where
$\Lambda$ is a collection of representations, one for each component
of the link. We should note that mathematicians usually label these
invariants by the Lie algebra as we do while physicists label the
invariants by the corresponding Lie group. There is an easier way
suggested by H Wenzl \cite{wenzl} to get the same invariant using
skein relations. Even though this definition is easy to understand
from first principles it is next to impossible to compute with. In
this appendix we first review Wenzl's elementary construction and
then indicate how the same invariant can be obtained from quantum
groups.

\begin{figure}[ht!]
\centering
\labellist\small
\pinlabel {$L_+$} [r] at 30 56
\pinlabel {$L_-$} [r] at 325 56
\pinlabel {$L_0$} [r] at 607 56
\endlabellist
\includegraphics[width=4truein]{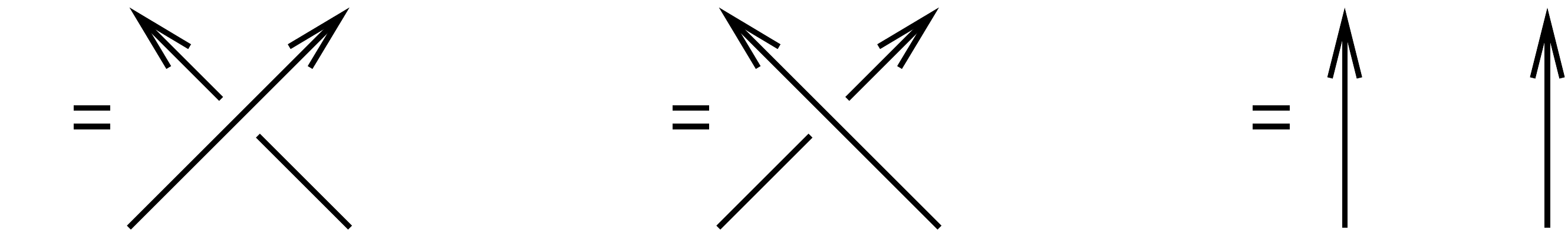} \caption{Terms in
the skein relation}\label{skein}
\end{figure}

The starting point for Wenzl's definition is the THOMFLYP
polynomial. (This generalization of the Jones polynomial was
simultaneously discovered by several different authors and described
in a joint paper \cite{HOMFLY}. The first acronym was HOMFLY and
this is still in common use. It has since been realized that two
other authors also deserve credit, so some people append PT to get
HOMFLYPT. We prefer to use the pseudonym.) This polynomial invariant
of links is defined by the recurrence relation
$$
\varlambda^{\frac12}\mathcal{P}(L_+)-\varlambda^{-\frac12}\mathcal{P}(L_-)=(\mathfrak{q}^{\frac12}-\mathfrak{q}^{-\frac12})\mathcal{P}(L_0),
$$
together with the normalization
$$\mathcal{P}(\text{unknot})=\frac{\varlambda^{\frac12}-
  \varlambda^{-\frac12}}{\mathfrak{q}^{\frac12}-\mathfrak{q}^{-\frac12}},$$
where $L_\pm$ and $L_0$ are three links that differ exactly in the
neighborhood of one crossing as in \fullref{skein}.
The THOMFLYP polynomial is not in fact a polynomial -- rather it is
a rational function of the formal variables $\mathfrak{q}^{\frac12}$
and $\varlambda^{\frac12}$. The Jones polynomial is recovered by
substituting $\varlambda=\mathfrak{q}^2$ into the THOMFLYP
polynomial.

Working backwards from \fullref{thomflypdef}, it is clear
that the framed link invariant should be defined by
\begin{equation}\label{WTHOMFLYP}
W^{{\mathfrak{sl}}_N}_{\tableau{1},\ldots,\tableau{1}}(L)=\eps^{(N-1/N)\sum_{i=1}^c\sum_{j=1}^c
n_{ij}(L)}\mathcal{P}(L)\,,
\end{equation}
where $n_{ij}(L)$ are the entries of the linking matrix. By
inserting a right-handed twist (the configuration labeled by
$\theta_V$ in \fullref{MTC}) into an arbitrary link diagram to
get $L_+$ and a left-handed twist to get $L_-$, one can compute that
$$\mathcal{P}(L0)=(\mathfrak{q}^{1/2}-\mathfrak{q}^{-1/2})
  (\varlambda^{1/2}-\varlambda^{-1/2})^{-1}\mathcal{P}(L)$$
where $L0$ is obtained from $L$ by adjoining a completely
unlinked and unknotted component. The physics literature often
uses a slightly different normalization. Namely,
$$
W^{\text{SU}(N)}_{\tableau{1},\ldots,\tableau{1}}(L)=\varlambda^{\sum_{i=1}^c\sum_{j=1}^c
n_{ij}(L)}\mathcal{P}(L)\,,
$$

\begin{exm}\label{exw1}
Assume zero framings (self-linking numbers). Compute
$$
W^{\text{SU}(N)}_{\tableau{1},\tableau{1}}(\text{left Hopf link})=
\left(\frac{\varlambda^{\frac12}-\varlambda^{-\frac12}}{\mathfrak{q}^{\frac12}-\mathfrak{q}^{-\frac12}}\right)^2
+\varlambda^{-1}-1\,.
$$
Compute $W^{\text{SU}(N)}_{\tableau{1},\tableau{1}}(\text{right Hopf
link}) $  and $W^{\text{SU}(N)}_{\tableau{1}}(\text{left (2,3) torus
knot})$. (The left Hopf link is depicted in \fullref{lhopf})
\end{exm}
The answer to the left $(2,3)$ torus knot may be found in
Mari\~no \cite{Marino-enu}.

There is a different invariant of framed links that may be
constructed out of The THOMFLYP polynomial. It is
$$
P^\alpha(L):= \alpha^{\sum_{i=1}^c\sum_{j=1}^c n_{ij}}P(L).
$$
To go further, we specialize by evaluating these polynomial
invariants at $\mathfrak{q}^{\frac12}{=}e^{{i\pi}/{(N+k)}}$,
$\varlambda^{\frac12}{=}e^{{iN\pi}/{(N+k)}}$.  Let $\gamma$ denote a
multi-index of length $c(L)$, that is, a $c(L)$--tuple of positive
integers. We will let $|\gamma|$ denote the sum of the components of
$\gamma$, and define $L^\gamma$ to be the framed link obtained from
$L$ by replacing the $p^{\text{th}}$ component of $L$ by $\gamma_p$
parallel copies. Define a number by
$$
\Gamma:=
\biggl(\frac{1}{N}\lim_{p\to\infty}\biggl(\frac{\varlambda^{\frac12}-\varlambda^{-\frac12}}{\mathfrak{q}^{\frac12}
-\mathfrak{q}^{-\frac12}}\biggr)^{-p}\mathcal{P}((\text{$1$--framed
unknot})^p)\biggr)^{-1}.
$$
Finally, Wenzl defines an invariant by
\begin{multline*}
\tau^{SU(N)}_k(M,\emptyset):= \\
N^{-1}\left(\frac{\Gamma}{|\Gamma|}\right)^{\sigma(L_M)-c(L)+1}
\lim_{p\to\infty}p^{-c(L)}\Gamma^{c(L)}\sum_{\text{max}(\gamma)\le
p} \Theta^{-|\gamma|}\mathcal{P}(L_M^\gamma),
\end{multline*}
where
$$
\Theta:=\left(\frac{\varlambda^{\frac12}-\varlambda^{-\frac12}}{\mathfrak{q}^{\frac12}-\mathfrak{q}^{-\frac12}}\right).
$$
The advantage of this definition is that all of the ingredients are
elementary. The disadvantages are that one has to work to prove that
it is well-defined, and it is not obvious how one can compute based
on this definition. Notice that the $1$--framed unknot squared is
just the right Hopf link, and then try to compute
$\mathcal{P}((\text{$1$--framed unknot})^3)$. This is fairly
difficult. Computing the invariant for higher cables without some
trick looks hopeless. Wenzl proves that this invariant is indeed
well-defined and is equal to the Reshetikhin--Turaev invariant
(see Wenzl \cite{wenzl}).

To prove that Wenzl's invariant is equal to the quantum group
invariant one notices that $\calQ_\lambda^q\otimes
\calQ_\mu^q=\sum_{\nu\in I}N^\nu_{\lambda\mu}\calQ_\nu^q$, where
$N^\nu_{\lambda\mu}$ are constants determined by the structure of a
modular category. Via repeated application of this formula one may
replace computations in arbitrary representations by computations
for various cables in the fundamental representation and then apply
the skein relation.

\section{Representation theory of Lie groups and Lie algebras}
\setobjecttype{App}
\label{app:d}

This appendix reviews some aspects of classical representation
theory that are used in this paper. The information we summarize
here can be found in the wonderful books by Fulton and Harris
\cite{fulton-harris} and Humphreys \cite{Hum}. We will denote the
$N\times N$ unitary matrices ($A^\dagger A=I$) by $U(N)$.
\index{$U$@U$(N)$ unitary group} \index{$S$@SU$(N)$ special unitary
group} The special unitary matrices are those with unit determinant,
this group will be denoted by SU$(N)$. The Lie algebras of these two
groups are the algebra of Hermitian matrices ($A^\dagger+A=0$),
denoted ${\mathfrak u}_N$, and trace-free Hermitian matrices,
denoted by ${\mathfrak su}_N$, respectively with the standard
bracket ($[A,B]=AB-BA$) as a product. The complexifications of these
two algebras are the collection of all complex matrices
(${\mathfrak{gl}}_N{\mathbb C}$) and the subalgebra of trace-free
matrices (${\mathfrak{sl}}_N{\mathbb C}$)
\index{${\mathfrak{sl}}_N{\mathbb C}$ trace-free matrices}
respectively. Let $E_{ij}$ \index{$E_{ij}$ matrix generator} denote
the matrix with a $1$ in the $j^{\text{th}}$ column of the
$i^{\text{th}}$ row and zeros elsewhere.

\subsubsection*{Young diagrams and irreducible representations}
\addcontentsline{toc}{subsection}{Young diagrams and irreducible
representations}

A representation \index{representation} of a group (algebra) is a
homomorphism into the automorphisms (endomorphisms) of a vector
space. All of the groups and algebras defined above admit an obvious
representation with vector space ${\mathbb C}^N$ called the
fundamental or defining representation.
\begin{remark}
Of course when we talk about differentiating a representation we are
talking about a representation of a Lie group. Representations of
groups and algebras are closely related but there is not an exact
correspondence. For example, the fundamental representation of
$GL_N\C$ can be differentiated to produce a representation of the
algebra ${\mathfrak{gl}}_N\C$. The same representation can be
conjugated to give a different representation of $GL_N\C$; however,
differentiating the conjugate representation will not produce a
corresponding algebra representation because it will not be complex
linear. In the other direction one can see that there is a
two-complex dimensional representation of the Lie algebra
${\mathfrak so}_3$ via the isomorphism ${\mathfrak su}_2$. However
there is no corresponding representation of the three-dimensional
rotation group SO$(3)$. The classical groups that we are considering
are all matrix groups, so each may be viewed as a subset of
$\C^{N^2}$ via the fundamental representation. The automorphisms of
any complex vector space can be described as a subset of complex
$m$--space in the same way. It therefore makes sense to consider
polynomial representations, that is, those that can be described by
polynomials. For simply-connected matrix groups there is a perfect
correspondence between finite-dimensional polynomial representations
\index{polynomial representation} of the group and finite
dimensional representations of the Lie algebra given by
differentiation.
\end{remark}

There are several obvious ways to construct new representations from
(sets of) old representations: the dual, direct sum, tensor product
etc. If $\rho_k\co G\to \text{Aut}(V_k)$ for $k=1, 2$ are group
representations, then the tensor product representation
$\rho_1\otimes\rho_2\co G\to \text{Aut}(V_1\otimes V_2)$ is given by
$(\rho_1\otimes\rho_2)(g)(x\otimes y)=\rho_1(g)(x)\otimes
\rho_2(g)(y)$. The tensor product of algebra representations
$\mu_k\co {\mathfrak g}\to \text{End}(V_k)$ for $k=1, 2$ is given by
$(\mu_1\otimes\mu_2)(A)(x\otimes y)=\mu_1(A)(x)\otimes y + x\otimes
\mu_2(A)(y)$. This is obtained by differentiating the tensor product
of group representations.

\begin{figure}[ht!]
\centering
\labellist\small
\pinlabel {$1$} at 307 43
\pinlabel {$1$} at 334 43
\pinlabel {$1$} at 361 43
\pinlabel {$1$} at 390 43
\pinlabel {$2$} at 307 16
\pinlabel {$2$} at 334 16
\endlabellist
\includegraphics[width=4truein]{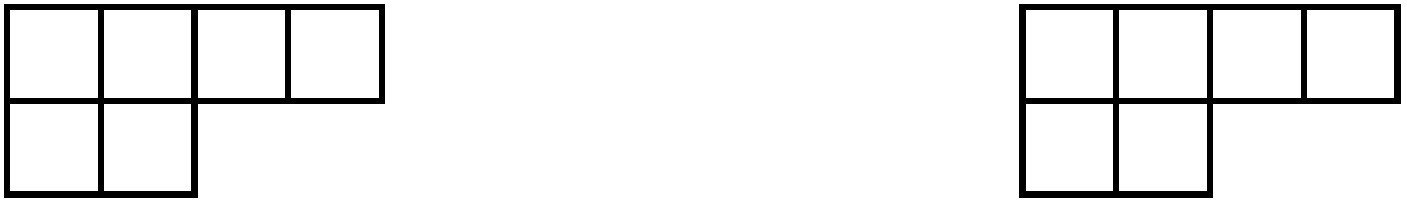}
\caption{Young diagram}\label{fig-young}
\end{figure}

Additional representations can be constructed by symmetric or
anti-symmetric tensors. All of the symmetries that we need can by
constructed using Young symmetrizers \index{$c_\lambda$ Young
symmetrizer}\index{Young symmetrizer} arising from \index{Young
diagram} \index{$\lambda$ Young diagram} Young diagrams. A typical
Young diagram is displayed on the left in \fullref{fig-young}. A
Young diagram represents a partition of a positive integer,
$\ell=\ell_1+\ell_2+\cdots+\ell_k$ with $\ell_j\ge \ell_{j+1}$. It
is standard to denote such a partition by $\lambda$ and we use the
same notation for the Young diagrams as well. The Young diagram in
\fullref{fig-young} corresponds to the partition $6=4+2$. A Young
tableau is a Young diagram filled in with natural numbers according
to some rules. Young diagrams describe vector spaces of
representations while Young tableau describe special basis vectors
in those spaces, namely weight vectors. The Young tableau
\index{Young tableau} describing what we will later define as the
highest weight associated to a Young diagram is obtained by filling
in a diagram in a specific manner (top row with ones, second row
with two's, etc) starting in the upper left corner and ending in
the lower right one (see the right side of \fullref{fig-young}).

The point is that a Young diagram encodes an endomorphism of the
$\ell$--fold tensor product of any vector space. Construct a specific
Young tableau by filling the Young diagram with the numbers $1$
through $\ell$ filling in rows starting in the upper left. If
$\lambda$ is a Young diagram, $a_\lambda$ will denote the
endomorphism of the tensor product that takes a tensor product of
vectors to the sum of the permutated tensor products by permutations
preserving the numbers in the rows of the Young tableau. Another
endomorphism, $b_\lambda$, is defined analogously but with the
alternating sum (according to the signs of permutations) over
permutations that preserve numbers in the columns of the Young
tableau. In the example from the figure,
$$
b_\lambda(e_{1,2,3,4,5,6})=e_{1,2,3,4,5,6}-e_{5,2,3,4,1,6}-e_{1,6,3,4,5,2}+e_{5,6,3,4,1,2},
$$
where $e_{1,2,3,4,5,6}$ denotes $e_1\otimes e_2\otimes e_3\otimes
e_4\otimes e_5\otimes e_6$ and $\{e_k\}$ is a basis for $V$. The map
$a_\lambda$ can be computed similarly; it is a sum of $48$ terms.

The Young symmetrizer is the composition of these two endomorphisms
$c_\lambda=a_\lambda\circ b_\lambda$ (a sum of $192$ terms in our
example). The image of a Young symmetrizer $c_\lambda$ applied to a
tensor product of $\ell$ copies of the fundamental representation is
a \index{$V_\lambda$ representation corresponding to $\lambda$}
denoted by $V_\lambda$. Computing a basis for $V_\lambda$ directly
from the definition is a little cumbersome. We will first consider
two special cases corresponding to a diagram with one row or one
column.

\begin{example}
For a diagram with just one row, the partition is just $\ell=\ell$.
The $a_\lambda$ endomorphism is just the sum over all permutations
and the $b_\lambda$ is the identity map. It follows that
$V_{\ell=\ell}$ is just the symmetric product of $\ell$ copies of
${\mathbb C}^N$ ( $\text{Sym}^\ell{\mathbb C}^N$). Similarly, for a
diagram with just one column the partition is just $\ell=1+1+\cdots
+1$. This time $a_\lambda$ is trivial and $b_\lambda$ is just the
alternating sum over all permutations so $V_{\ell=1+1+\cdots +1}$ is
just the $\ell{\text{--th}}$ exterior power ($\bigwedge^\ell{\mathbb
C}^N$).
\end{example}
\begin{exm}
Define the natural action of any of the classical matrix groups or
algebras on the set of homogeneous polynomials of degree $\ell$.
Show that this representation is isomorphic to
$\text{Sym}^\ell{\mathbb C}^N$.
\end{exm}
\begin{remark}
Whereas the $\ell{\text{--th}}$ symmetric power is nontrivial for
all $\ell$, the $\ell{\text{--th}}$ exterior power is trivial unless
$\ell\le N$. Furthermore, the group representation on
$\wedge^N{\mathbb C}^N$ is given by the determinant, so the
corresponding $SU(N)$ representation is trivial unless $\ell<N$.
Likewise the algebra representation on the top exterior power is
given by the trace so the $\mathfrak{su}_N$ and $\mathfrak{sl}_N\C$
representations are trivial unless $\ell<N$. The same comments hold
for any Young diagram: The $\text{GL}_N\C$, $\mathfrak{gl}_N\C$ and
U($N$) representations are trivial unless the diagram has less than
or equal to $N$ rows and the $\text{SL}_N\C$, $\mathfrak{sl}_N\C$
and SU($N$) representations are trivial unless the diagram has
strictly less than $N$ rows.
\end{remark}

When $\lambda$ has just one row, the set of vectors represented by
all possible ways of filling in the Young diagram with a
non-decreasing sequence of integers between $1$ and $N$ inclusive is
a basis for $V_{\lambda}$. When $\lambda$ has just one column, the
set of vectors represented by all possible ways of filling in the
Young diagram with an increasing sequence of integers between $1$
and $N$ inclusive is a basis for $V_{\lambda}$. In general,
$V_\lambda$ has a basis consisting of all ways of filling the Young
diagram with a  sequence of integers between $1$ and $N$ inclusive
so that the numbers do not decrease as one reads across rows and so
that the numbers strictly increase as one reads down columns. This
is exactly the rule specifying the Young tableau alluded to earlier.

The action of a matrix on the vector space $V_\lambda$ can be
described easily: the matrix acts on the tensor product of the
vectors with the given indices either using the group or algebra
action as appropriate, the answer is fully expanded and written as a
combination of terms in the standard non-decreasing/increasing
order.

As we go further, we will concentrate on finite-dimensional
representations. A finite-dimensional \index{indecomposable}
representation is called decomposable  if it can be expressed as a
nontrivial direct sum. It is indecomposable otherwise. A
representation is reducible if it contains a nontrivial
subrepresentation and irreducible otherwise. \index{irreducible
representation}  In all of the cases that we consider, every
finite-dimensional representation will be a direct sum of
irreducible representations. The irreducible representations
correspond to Young diagrams. It is not immediately clear that the
representations $V_\lambda$ are irreducible or that every
irreducible representation is of this form, but this is indeed the
case. See the book by Fulton and Harris for the complete story
\cite{fulton-harris}.

To be specific irreducible polynomial representations of
$\text{GL}_N\C$, $\mathfrak{gl}_N\C$ and U($N$) are indexed by Young
diagrams (partitions) with less than or equal to $N$ rows and
irreducible polynomial representations of $\text{SL}_N\C$,
$\mathfrak{sl}_N\C$ and SU($N$) are indexed by Young diagrams
(partitions) with strictly less than $N$ rows. It is not accidental
that irreducible representations for the different groups and
algebras here are indexed by the same sets. The fact is that every
irreducible representation of $\text{GL}_N\C$ restricts to an
irreducible representation of U$(N)$, which is its maximal compact
subgroup and all of the latter are obtained in this manner. The same
thing happens with $\text{SL}_N\C$,  SU($N$), ${\mathfrak{gl}}_N\C$
and ${\mathfrak{sl}}_N\C$.

We can encode our discussion up to this point in a definition.
\begin{defn}
Let $\lambda$ be a Young diagram with rows of length $\ell_1,
\ldots, \ell_N$ and columns of length $m_1, \ldots, m_{\ell_1}$. The
associated standard tableau is obtained by filling in the diagram
with the numbers one through $\ell:=\ell_1+\cdots+\ell_N$ across
rows starting in the upper left. The associated tableaux are
obtained by filling in the diagram with some numbers between one and
$N$ such that numbers are non-decreasing along rows and strictly
increasing down columns. Let $A_\lambda={\mathfrak
G}_{\ell_1}\times\cdots\times{\mathfrak G}_{\ell_N}$ embedded as the
subgroup of the permutation group $\mathfrak{G}_\ell$ preserving the
rows of the standard tableau. Let $B_\lambda={\mathfrak
G}_{m_1}\times\cdots\times{\mathfrak G}_{m_{\ell_1}}$ embedded as
the subgroup of the permutation group $\mathfrak{G}_\ell$ preserving
the columns of the standard tableau. The associated elements of the
group ring are $a_\lambda:=\sum_{\sigma\in A_\lambda}\sigma$,
$b_\lambda:=\sum_{\sigma\in B_\lambda}(-1)^\sigma\sigma$ and the
Young symmetrizer $c_\lambda=a_\lambda\circ b_\lambda$. The Specht
module $V_\lambda$ is the image of the tensor power
$\smash{V^{\otimes\ell}}$ under the natural action of the Young
symmetrizer.
\end{defn}

This is starting to get complicated. Every irreducible
representation is generated by what are called weight vectors. The
weight vectors in one of the $V_\lambda$ representations are just
the vectors represented by Young tableau.
\begin{example}\label{c1123.1}
For example in the ${\mathfrak{sl}}_4{\mathbb C}$ representation
given by $V_{\tableau{4 2}}$ \index{${\tableau{4 2}}$ Young diagram}
the vector obtained by filling the first row of the diagram with the
sequence $1, 1, 2, 3$ and the second row with $2, 4$ is a weight
vector. It can be written as $c_{\tableau{4 2}}(e_{1,1,2,3,2,4})$.
Let's see how several different elements of
${\mathfrak{sl}}_4{\mathbb C}$ act on the vector from our example.
The matrix $E_{21}$ maps $e_1$ to $e_2$ and the rest of the vectors
to $0$. So our weight vector will get mapped to
$c_{\tableau{4 2}}(e_{2,1,2,3,2,4})+c_{\tableau{4
2}}(e_{1,2,2,3,2,4})=c_{\tableau{4 2}}(e_{1,2,2,3,2,4})$.
The matrix
$E_{14}$ would replace the $4$ in the second row by a $1$, but this
is just the zero vector. For a more complicated example notice that
\begin{align*}
E_{21}c_{\tableau{4 2}}(e_{1,1,2,2,3,4}) &=c_{\tableau{4
2}}(e_{2,1,2,2,3,4})+c_{\tableau{4 2}}(e_{1,2,2,2,3,4})\\
&= 2c_{\tableau{4 2}}(e_{1,2,2,2,3,4})+c_{\tableau{4
2}}(e_{1,2,2,4,2,3})-c_{\tableau{4 2}}(e_{1,2,2,3,2,4})\,.
\end{align*}

\end{example}
\begin{exm}\label{126}
Verify that dim$(V_{\tableau{4 2}})=126$.
\end{exm}

We will concentrate on Lie algebra representations for a while. We
can see more about these representations once we know more about Lie
algebras. The Cartan subalgebra is a maximal abelian (every bracket
is zero) subalgebra. The Cartan subalgebra is usually denoted by
${\mathfrak h}$. \index{Cartan subalgebra}\index{$H$@${\mathfrak h}$
Cartan subalgebra}
\begin{exm}
Show that it is possible to simultaneously diagonalize a set of
Hermitian matrices exactly when the matrices in the set commute.
\end{exm}
\begin{example}
Using the previous exercise by changing the basis one can take the
subset of diagonal matrices as the Cartan subalgebra. The Cartan
subalgebra of ${\mathfrak{gl}}_N{\mathbb C}$ is generated by the
matrices $E_{ii}$. For ${\mathfrak{sl}}_N{\mathbb C}$ we can take
$E_{ii}-E_{i+1 i+1}$ as a basis for the Cartan subalgebra.
\end{example}
The dimension of the Cartan subalgebra is called the rank of the
group or algebra. It is denoted by $r$ in general and is $N-1$ for
${\mathfrak{sl}}_N{\mathbb C}$.

Just as eigenvalues are important invariants of matrices, weights
are important invariants of a representation.
\begin{defn}
The weights \index{weight}\index{weight vector} of an arbitrary
representation $\mu\co {\mathfrak g}\to \text{End}(V)$ are linear
functionals $\omega\co {\mathfrak h}\to {\mathbb C}$ that have an
associated non-zero weight vector, $v_\omega\in V$, satisfying
$\mu(H)v_\omega = \omega(H)v_\omega$ for every $H\in {\mathfrak h}$.
The set of all weight vectors (including zero) is called the weight
space and is denoted by $V^\omega$. The set of all weights of all
representations is called the weight lattice and is denoted
$\Lambda_w$. \index{weight lattice}\index{$\Lambda_w$ weight
lattice}
\end{defn}
Thus we see that weights are just eigenvalues of families of
operators, weight vectors are just the corresponding eigenvectors
and weight spaces are just the corresponding eigenspaces. The sum of
weights corresponding to the same representation may not be a weight
of the same representation; however, it is a weight in the tensor
product of the representation with itself. It helps to consider the
weights of all possible representations together. These form a
subgroup of the dual to the Cartan subalgebra $\mathfrak{h}^*$.

\begin{exm}
Prove that the set of weights is a subgroup of the dual ${\mathfrak
h}^*$. Hint: think about the tensor product and dual of
representations.
\end{exm}

\begin{exm}
Verify that Young tableaux always represent weight vectors.
\end{exm}

\begin{example}
Recall the weight vector $v_\lambda= c_{\tableau{4
2}}(e_{1,1,2,2,3,4})$ from \fullref{c1123.1}. The elements of
the Cartan subalgebra act as follows: $E_{11}-E_{22}$ multiplies our
weight vector by zero since there are two ones and two twos,
$E_{22}-E_{33}$ multiplies the vector by $1=2-1$, and
$E_{33}-E_{44}$ multiplies the vector by $0=1-1$. Thus $v_\lambda$
is indeed a weight vector.
\end{example}

It is helpful to write the weights in coordinates. Recall the
standard $E_{ii}$ basis of the Cartan subalgebra of
$\mathfrak{gl}_N\C$. Let $E^*_{ii}$ denote the dual basis and set
$$I^*:=\sum_{i=1}^NE^*_{ii},,
$$
\index{$I^*:=\sum_{i=1}^NE^*_{ii}$} then the Cartan subalgebra of
$\mathfrak{sl}_N\C$ contains exactly the matrices annihilated by
$I^*$. Therefore we can identify its dual
$\mathfrak{h}^*_{\mathfrak{sl}_N}$ with the subspace of
$\mathfrak{h}^*_{\mathfrak{gl}_N}$ that annihilates $I$. Thus
$\mathfrak{h}^*_{\mathfrak{sl}_N}$ is generated by
$E^*_{ii}-E^*_{i+1 i+1}$. It is not at all obvious, but the weight
lattices of $\mathfrak{gl}_N\C$, $\mathfrak{sl}_N\C$ are generated
over $\Z$ by $E^*_{ii}$ and
$$L_i:=E^*_{ii}-\frac1NI^*$$\index{$L_i:=E^*_{ii}-\frac1NI^*$}
respectively.

\begin{example}
Given this notation, the weight corresponding to the weight vector
of our $c_{\tableau{4 2}}(e_{1,1,2,2,3,4})$ example is
$E^*_{11}+E^*_{22}-\frac24I^*=L_1+L_2$.
\end{example}
\begin{exm}
Assuming that the irreducible representations are exactly
$V_\lambda$ \index{$V_\lambda$ irrep with highest weight $\lambda$}
show that the weight lattice of ${\mathfrak{sl}}_N{\mathbb C}$ is
generated by $L_i:=E^*_{ii}-\frac{1}{N}I^*$.
\end{exm}
Every Lie algebra has one special representation called the adjoint
representation.
\begin{defn} The adjoint representation \index{adjoint representation}\index{$A$@ad adjoint representation} is
$\text{ad}\co {\mathfrak g}\to \text{End}({\mathfrak g})$ given by
$\text{ad}(X)(Y)=[X,Y]$. The weights of the adjoint representation
are called roots. Roots are typically denoted by $\alpha$. The
corresponding weight vectors are called root vectors. \index{root
vector} The set \index{root} of roots is denoted by \index{$\Delta$
roots} $\Delta$ and the subgroup of ${\mathfrak h}^*$ generated by
the roots is called the root lattice and is denoted
\index{$\Lambda_r$ root lattice} $\Lambda_r$. \index{root lattice}
\end{defn}
The adjoint representation is not necessarily irreducible. If it is the algebra is called simple. For example,
$\mathfrak{sl}_N\C$ is simple and $\mathfrak{gl}_N\C$ is not.
\begin{exm}
Find the Young diagram corresponding to the adjoint representation
of ${\mathfrak sl}_N{\mathbb C}$. You can check your answer by
computing the dimension.
\end{exm}

It is standard to pick a vector in ${\mathfrak h}$ that is not
annihilated by any non-zero weight to measure the heights of
\index{height function}\index{$H$@Ht height function} weights and
roots. This vector can be chosen so that any weight evaluated on
this vector gives a real number. We can take the matrix
$$\text{Ht}:=\pi^{N-1}E_{11}+\pi^{N-2}E_{22}+\cdots E_{NN}-
\frac{\pi^N-1}{\pi-1}I\,,$$ where $I=\sum E_{ii}$ to define heights
for ${\mathfrak{sl}}_N{\mathbb C}$. \index{highest weight vector}
\begin{exm}
The highest weight vector in our sample representation is then the
one displayed in \fullref{fig-young}.
\end{exm}
\begin{defn}
The roots that evaluate to a positive number are called the positive
roots. The set of positive roots is denoted by \index{$\Delta^+$
positive roots}\index{positive root} $\Delta^+$. The positive roots
that cannot be written as a sum of positive roots are called simple
roots. \index{simple root}
\end{defn}
The negative of any root is also a root. Any Lie algebra may be
written as
$$
{\mathfrak g}={\mathfrak h}\oplus \bigoplus_{\alpha\in\Delta^+}
({\mathfrak g}_{-\alpha}\oplus {\mathfrak g}_\alpha)\,,
$$
where ${\mathfrak g}_{\pm\alpha}$ is the eigenspace of ${\mathfrak
h}$ \index{$G$@${\mathfrak g}_{\pm\alpha}$ root eigenspace} with
eigenvalue $\pm\alpha$.
\begin{defn}
A weight vector in a representation that is annihilated by
${\mathfrak g}_+:=\bigoplus_{\alpha\in\Delta^+}{\mathfrak
g}_{\alpha}$ is called a highest weight vector. The corresponding
weight is called a highest weight.
\end{defn}
The good news is that this is all very explicit for
${\mathfrak{sl}}_N{\mathbb C}$.
\begin{example}
The roots of ${\mathfrak{sl}}_N{\mathbb C}$ are just
\index{$\alpha_{ij}:=E_{ii}^*-E^*_{jj}$ roots}
$\alpha_{ij}:=E_{ii}^*-E^*_{jj}$. The positive roots are those with
$i<j$, and the simple roots are those with $j=i+1$. We use the
notation $\alpha_i=E^*_{ii}-E^*_{i+1 i+1}$ for simple roots. The
root vectors are $e_{ij}:=E_{ij}$ \index{$E$@$e_{ij}:=E_{ij}$} and
$f_{ij}:=E_{ji}$ \index{$F$@$f_{ij}:=E_{ji}$} The \index{root
vector} \index{$\alpha_i=E^*_{ii}-E^*_{i+1 i+1}$ simple roots}
subspace ${\mathfrak{sl}}_N{\mathbb C}_+$ is generated by $E_{ij}$
with $i<j$, thus the weight vector from our $c_{\tableau{4
2}}(e_{1,1,2,2,3,4})$ example is not a highest weight vector since
it is not annihilated by $E_{12}$.
\end{example}
\begin{remark}
The vector in a representation with highest weight as measured by
the height Ht is a highest weight vector. In general the highest
weights of a fixed representation are the weights with maximal
height in each of the irreducible factors.
\end{remark}
The notation ${\mathfrak
g}_-:=\bigoplus_{\alpha\in\Delta^+}{\mathfrak g}_{-\alpha}$ will
also be useful. If $v_\lambda$ is a highest weight vector in a
representation, the space $V_\lambda$ (notice the duplicate notation
from our discussion with Young diagrams) that is generated by
$({\mathfrak h}\oplus{\mathfrak g}_{-})v_\lambda$ will be an
irreducible subrepresentation and the original representation will
be the direct sum of all such. This is a very powerful idea that we
should elaborate. We begin with an example.
\begin{example}\label{wtspace}
For $\mathfrak{sl}_N\C$ the space $\mathfrak{sl}_N\C_-$ is spanned
by the matrices $E_{ij}$ with $i>j$. For the representation
described in \fullref{c1123.1} the highest weight vector
$v_{\tableau{4 2}}:=c_{\tableau{4 2}}(e_{1,1,1,1,2,2})$ is clearly a
weight vector. The fact that it is a highest weight vector allows
one to reconstruct the entire representation. In order to construct
a basis for the representation one should just start multiplying the
vector by the matrices $E_{ij}$ with $i>j$. For example, one should
add $E_{21}c_{\tableau{4 2}}(e_{1,1,1,1,2,2})=4c_{\tableau{4
2}}(e_{1,1,1,2,2,2})$ to the basis. Relations in the Lie algebra
determine the action of the rest of the Lie algebra on this vector.
The first thing that one can see is that it is a weight vector. This
is because the bracket of any element of the Cartan subalgebra with
$E_{21}$ is a multiple of $E_{21}$, for example
$[E_{11}-E_{22},E_{21}]=2E_{21}$. This implies that
$$
\begin{aligned}
(E_{11}-E_{22})E_{21}v_{\tableau{4
2}}&=[E_{11}-E_{22},E_{21}]v_{\tableau{4
2}}+E_{21}(E_{11}-E_{22})v_{\tableau{4
2}}\\
&= 2E{21}v_{\tableau{4 2}} +2E{21}v_{\tableau{4 2}}\,.
\end{aligned}
$$
Using the bracket $[E_{12},E_{21}]=E_{11}-E_{22}$ one can compute
$E_{12}E_{21}\smash{v_{\tableau{4 2}}}$. By continuing in this way adding
additional vectors such as $E_{32}E_{21}v_{\tableau{4 2}}$ and
$E_{21}^2v_{\tableau{4 2}}$ etc. The entire representation may be
reconstructed.
\end{example}

There are two important things to remember from this example. The
first is that any finite-dimensional representation is a direct sum
of weight spaces: $V=\oplus_\omega V^\omega$. The  $V^\omega$ are
not representations of the full Lie algebra, they are only
representations of the Cartan subalgebra. This weight space
decomposition completely determines the action of the Cartan
subalgebra according to $\mu(H)v=\omega(H)v$. This follows from the
natural generalization of $[E_{11}-E_{22},E_{21}]=2E_{21}$. In fact
when $E\in\mathfrak{g}_\alpha$, one has $E^nV^\omega\subseteq
V^{\omega+n\alpha}$.

The second thing to remember is how the entire representation could
be reconstructed from the highest weight vectors. To formalize this
define the universal enveloping algebra $U\mathfrak{g}$ to be the
associative unital algebra generated by $\mathfrak{g}$ subject to
the relations $[X,Y]=XY-YX$. For $\mathfrak{sl}_N\C$ this is just
the matrix algebra structure. The point is that in a general Lie
algebra one cannot multiply elements; one can only take brackets.
The universal enveloping algebra allows one to multiply elements.
Let $U_\pm$ be the subalgebras generated by
$\mathfrak{h}\oplus\mathfrak{g}_\pm$ respectively. Specifying a
weight vector $v_\lambda$ with weight $\lambda$ specifies a
$U_+$--module structure on $\C$ thought of as the vector space
generated by $v_\lambda$. The formal way to extend this to a
representation of the universal enveloping algebra (and thus to the
Lie algebra) is as the $U\mathfrak{g}$--module
$U\mathfrak{g}\otimes_{U_+}\C$.
\begin{remark}\label{weylwtrk}
We can summarize these two points by stating that every irreducible
representation is generated by words in $U_-$ multiplied by a
highest weight vector and the weight space decomposition can be used
to determine which words act trivially (see \fullref{sl2} and the discussion after it).
\end{remark}
Summarizing the above discussion we get the following description of
the irreducible representations in the cases of most interest to us.
\begin{quote}
Irreducible complex representations of $\mathfrak{gl}_N\C$ are
indexed by Young diagrams (partitions) with less than or equal to
$N$ rows and irreducible complex representations of
$\mathfrak{sl}_N\C$ are indexed by Young diagrams (partitions) with
strictly less than $N$ rows.
\end{quote}
In the main body of the paper we mention that general
finite-dimensional representations of quantum groups neither
decompose into direct sums of weight spaces nor contain highest
weight vectors. There is however a class of sub-representations
called tilting modules that do possess these desirable properties.

Following the ideas outlined above and in \fullref{wtspace}
leads to a classification of all finite-dimensional representations
of any Lie algebra, in particular $\mathfrak{gl}_N\C$, $\mathfrak{sl}_N\C$
(see Fulton and Harris \cite{fulton-harris}).
\begin{exm}
Prove that every finite-dimensional representation can be decomposed
as a sum of the $V_\lambda$ uniquely.
\end{exm}

\subsubsection*{Dominant weights and the Weyl character formula}
\addcontentsline{toc}{subsection}{Dominant weights and the Weyl
character formula}

As nice as the above description might be it is not convenient for
generalization to quantum groups. We now briefly review a different
(but closely related) approach due to E Cartan that uses dominant
weights instead of Young diagrams.

Any Lie algebra inherits a bilinear form according to
$$
\langle
X,Y\rangle=\text{constant}\cdot\text{Tr}(\text{ad}(X)\circ\text{ad}(Y))\,.
$$
This pairing is called the Killing form. \index{Killing form} A Lie
algebra is called semisimple exactly when this form is
nondegenerate.
\begin{exm}
Prove that every simple algebra, that is, one with irreducible adjoint
representation, is semisimple.
\end{exm}
There are two standard ways to normalize it. It induces a form on
${\mathfrak h}^*$ and one can require that
$\langle\alpha,\alpha\rangle=2$ for short roots or for long roots.
In the case of ${\mathfrak{sl}}_N{\mathbb C}$ all roots have the
same length, so there is a standard interpretation of the Killing
form. In fact it is just given by $\langle
A,B\rangle=\text{Tr}(AB^\dagger)$. Using the Killing form we may
identify the Cartan subalgebra with its dual. Under this
identification, the coroots \index{coroot}\index{$\alpha_i^\vee$
coroots} are defined by $\alpha^\vee=2\langle
\alpha,\alpha\rangle^{-1}\alpha$ (so in ${\mathfrak{sl}}_N{\mathbb
C}$ there is no difference between the coroots and the roots.) The
coroot lattice $\Lambda_r^\vee$ \index{$\Lambda_r^\vee$ coroot
lattice} is the lattice generated by the coroots. A coroot
$\alpha^\vee$ is simple if the corresponding root $\alpha$ is
simple. If $\alpha_i$ are the positive simple roots then the basis
biorthogonal to $\alpha^\vee_i$ is denoted by $\omega_i$ (that is,
$(\alpha_i^\vee,\omega_j)=\delta_{ij})$ and its elements are called
the \index{$\omega_i$ fundamental weights}\index{fundamental
weights} fundamental weights. The sum of fundamental weights $\rho:=
\sum_i\omega_i$ also plays an important role and is called the
\index{$\rho:= \sum_i\omega_i$ Weyl weight}\index{Weyl weight} Weyl
weight.
\begin{exm}\label{Ww}
Show that the fundamental weights for ${\mathfrak{sl}}_N{\mathbb C}$
are given by $\omega_k=\sum_{i=1}^kL_i=\sum_{i=1}^k
E_{ii}^*-\frac{k}{N}I^*$.
\index{$\omega_k=\sum_{i=1}^kL_i=\sum_{i=1}^k
E_{ii}^*-\frac{k}{N}I^*$} That is,
$\langle\alpha_i^\vee,\omega_j\rangle=\delta_{ij}$. Also show that
the Weyl weight can be written as $\rho:=\frac12
\sum_{\alpha\in\Delta^+} \alpha$, which is
$\sum_{i=1}^{N-1}(N-i)L_i$ for ${\mathfrak{sl}}_N{\mathbb C}$.
\end{exm}
\begin{exm}
Show that $L_k$ are (all) the weights of the defining representation of $\Sl_N\C$.
\end{exm}
We denote the highest root by $\theta$. The dual Coxeter number of a
Lie algebra is $h^\vee=\langle\rho,\theta\rangle+1$. For ${\mathfrak
sl}_N{\mathbb C}$ the highest root is $E_{11}^*-E_{NN}^*$
\index{$\theta=E_{11}^*-E_{NN}^*$} and the dual Coxeter number is
just $N$. In the main body of the paper we work only with $\Sl_N\C$.
Thus, we often assume that the Cartan subalgebra is identified with
its dual and drop $~^*$ and $~^\vee$ from notation.
\begin{exm}
Show that the Weyl weight is given by
$$\rho=\frac12 \sum_{k=1}^N
(N+1-2k)E_{kk}^*=\frac12\sum_{k=1}^N(N+1-2k)L_k$$ for
${\mathfrak{sl}}_N{\mathbb C}$.
\end{exm}
The Weyl group \index{Weyl group} is the group generated by
reflections in hyperplanes perpendicular to the roots. The
reflection corresponding to a simple root is
$s_i(\beta)=\beta-2\frac{\langle\beta,\alpha_i\rangle}
{\langle\alpha_i,\alpha_i\rangle}\alpha_i$. For
${\mathfrak{sl}}_N{\mathbb C}$ the Weyl group is the permutation
group \index{${\mathfrak G}_N$ permutation group} ${\mathfrak G}_N$.
The reflections $s_i$ act on the $E_{kk}^*$ via the permutation $(i
\ i+1)$.
\index{$s_i(\beta)=\beta-2\frac{\langle\beta,\alpha_i\rangle}
{\langle\alpha_i,\alpha_i\rangle}\alpha_i$ $s_i=(i \ i+1)$ on
$E_{kk}$}

It is time to establish a correspondence between the positive weights of
$\mathfrak{sl}_N\C$ and partitions (Young diagrams). Note that the
matrices $L_i:=E^*_{ii}-\frac1NI^*$, $1\leq i \leq N-1$ may also
serve as a basis of $\Lambda_w$ different from the basis of
fundamental weights $\omega_i$. These are related by,
$\omega_i=\sum_{j=1}^iL_j$. Given a highest weight $\lambda$ we have
$$\lambda=\sum_{i=1}^{N-1}n_i\omega_i=\sum_{i=1}^{N-1}n_i\sum_{j=1}^iL_j
  =\sum_{j=1}^{N-1}\biggl(\sum_{i=j}^{N-1}n_i\biggr)
  L_j=\sum_{j=1}^{N-1}\lambda_jL_j.$$
Since the $n_i$ are non-negative numbers we have
$\lambda_j\geq\lambda_{j+1}$, in other words the vector
$(\lambda_1,\ldots,\lambda_{N-1})$ is a partition. This is the
partition that corresponds to the highest weight $\lambda$ and we
abuse notation by using the same letter for the weight, partition
and Young diagram. It is possible to construct a Young diagram from
a sum of fundamental weights directly. Namely, if
$\lambda=\sum_{i=1}^{N-1}n_i\omega_i$ then the diagram has $n_{N-1}$
columns with $N-1$ boxes, $n_{N-2}$ columns with $N-2$ boxes, etc.
Since the $\lambda_i$ are the numbers of boxes in the rows this
provides a simple graphical method of converting sums of fundamental
weights into partitions and vice versa. Since we already know that
Young diagrams index irreducible complex representations of
$\Sl_N\C$ it makes sense to distinguish their counterparts among
weights.
\begin{defn}
A dominant weight \index{dominant weight} is a linear combination of
fundamental weights with non-negative integer coefficients. If
$\lambda=\sum_{i=1}^{N-1}n_i\omega_i$ is a dominant weight we denote
$\ell(\lambda):=\max\{i|\, n_i>0\}$ the length of $\lambda$ and
$|\lambda|:=\sum_{i=1}^{N-1}in_i$ the volume of $\lambda$. The set
of all dominant weights is denoted \index{$\Lambda_w^+$ dominant
weights} $\Lambda_w^+$.
\end{defn}
\begin{exm}
Show that the Young diagram corresponding to $\lambda$ has
$\ell(\lambda)$ rows and $|\lambda|$ boxes. For every non-negative
weight find a vector $H\in\mathfrak{h}$ such that
$|\lambda|=(\lambda,H)$. Hint: express $H$ as a sum of coroots and
use the biorthogonality relation. In case you are wondering: no,
$H\ne Ht$.
\end{exm}
The language of dominant weights allows one to characterize
irreducible representations for all semisimple Lie algebras
\begin{quote}
Any dominant weight is the highest weight of some irreducible
representation and any highest weight is dominant. Irreducible
complex representations of a complex semisimple Lie algebra, for example
$\Sl_N\C$ are indexed by  dominant weights.
\end{quote}

We now elaborate on this a bit.
\begin{defn}
If $\lambda$ is a dominant weight we let $\mathcal{I}_\lambda$ be
the left ideal of $U(\Sl_N\C)$ generated by $e_i$ and
$\alpha_i^\vee-\lambda(\alpha_i^\vee)$. The classical Verma module
associated to $\lambda$ is
\[
\what V_\lambda:=U(\Sl_N\C)/\mathcal{I}_\lambda\,.
\]
\end{defn}
The irreducible representation $V_\lambda$ that we constructed using
the Young diagram $\lambda$ is isomorphic to the maximal irreducible
quotient of $\what V_\lambda$. In particular the maximal irreducible
quotient is finite dimensional. Also every finite dimensional
representation decomposes into a direct sum of irreducible
representations. In addition the irreducible representations are
simple. One nice proof of all of these facts uses the Weyl unitary
trick. This goes as follows.

Any finite dimensional representation of $\Sl_N\C$ also gives a
representation of $\mathfrak{su}_N$ by restriction. This in turn
induces a representation of $SU(N)$ by exponentiation. It follows
from averaging over $SU(N)$ that any finite dimensional
representation of $\Sl_N\C$ admits an invariant inner product. This
quickly implies that any representation is a sum of simple
representations and that every representation has a weight. By
repeated action of the $f_i$ on any weight vector one can find a
highest weight vector. This is similar to the computation in example
\ref{wtspace}. The same computation shows that any two irreducible
representations with the same highest weight are isomorphic.
\begin{exm}
Write out details for the facts outlined in the previous paragraph.
Compare with example \fullref{sl2} and exercise \fullref{funs3bas}
from the main text.
\end{exm}

It is instructive to draw a picture of the weight lattice
$\Lambda_w$ and the subset of dominant weights inside of the dual to
the Cartan subalgebra. For $\mathfrak{sl}_N\C$ this dual can be
canonically identified with the subalgebra itself via the Killing
form. Under this identification the lattice in the background of
\fullref{translatedI} is the weight lattice for
$\mathfrak{sl}_3\C$.
\begin{exm}
Make a larger picture of the weight lattice for $\mathfrak{sl}_3\C$.
Given that the lattice points at $2/\sqrt{3}$ and
$(2/\sqrt{3})e^{\pi i/3}$ (viewing the picture in the complex plane)
correspond to the fundamental weights $\omega_1$ and $\omega_2$
respectively, plot $L_i$, the roots, simple roots, Weyl vector and
the set of dominant weights your picture.
\end{exm}

The final thing we need from the representation theory is the Weyl
character formula. Strangely enough, we need it for Lie group rather
than Lie algebra representations. This is because the form of quantum
groups (more precisely, quantized enveloping algebras) that we are
using is `partially integrated' and analogs of elements in the Cartan
subalgebra of $\mathfrak{sl}_N\C$ belong to $\text{SU}(N)$ rather
than $\mathfrak{sl}_N\C$. Of course, complex representation spaces of
$\text{SU}(N)$ and $\mathfrak{sl}_N\C$ can always be identified since
$\mathfrak{sl}_N\C$ is the complexified Lie algebra of $\text{SU}(N)$
(see Fulton and Harris \cite{fulton-harris}). The difference is only in
the operators whose traces are taken for characters.

We sketch a proof of the Weyl character formula here, but want to
mention that there is a different proof presented in the book by
Simon \cite{simon}.

Let $\lambda\co G\to\text{Aut}(\C^N)$ be a group representation then
its character is defined to be \index{character}
\index{$\chi_\lambda(g):=\text{tr}(\lambda(g))$ character}
$$
\chi_\lambda(g):=\text{tr}(\lambda(g)).
$$
A function $f\co G\to\C$ is called a class function if it is constant
on conjugacy classes, that is, $f(gxg^{-1})=f(x)$ for all $x$ and $g$.
Every character is a class function \index{class function} due to
the cyclic property of traces.
\begin{exm}
Verify the following basic facts about characters of unitary
representations: $\chi_{V\oplus W}=\chi_V+\chi_W$,
$\chi_{V^*}=\chi_V^*$ and $\chi_{V\otimes W}=\chi_V\chi_W$.
\end{exm}
We will need to use inner products and orthogonality relations for
characters. Let $G$ be a compact group and $dg$ be the Haar
(normalized bi-invariant) measure on it. For two continuous
complex-valued functions define the inner product
$$
\langle \psi,\varphi\rangle:=\int_{G}\psi(g)\wbar{\varphi}(g)\,dg\,,
$$
where $\wbar{ }$ stands for complex conjugation. We reserve $(\ ,\ )$ for
the inner product in $\text{U}(N)$. In the case of the symmetric
group $dg$ is just the normalized counting measure and
$$
\langle\psi,\varphi\rangle:=\frac1{N!}\sum_{\sigma\in{\mathfrak
G}_N}\psi(\sigma)\wbar{\varphi}(\sigma)\,.
$$
The symmetric group $\mathfrak{G}_N$ can be treated as a subgroup of
$\text{U}(N)$ by identifying permutations with permutation matrices.
Thus any function on $\text{U}(N)$ restricts to the symmetric group.

Given two representations $V$ and $W$ one can define a third as
Hom$(V,W)$ with action $(Af)(x):=A(f(A^*x))$. Let
Hom$(V,W)^{\text{\begin{small}U\end{small}}(N)}$ be the linear
subspace fixed by the action of U$(N)$.
\begin{exm}\label{orch}
Check that the transformation of
Hom$(V,W)^{\text{\begin{small}U\end{small}}(N)}$ defined by
$\Psi(f)=\sum_{\sigma\in{\mathfrak G}_N}\sigma f/N!$ is the identity
transformation. Notice that when $V$ and $W$ are irreducible
representations the dimension of
Hom$(V,W)^{\text{\begin{small}U\end{small}}(N)}$ is $1$ if they are
isomorphic and zero otherwise. By taking the trace of $\Psi$
conclude that
$\langle\chi_\lambda,\chi_\mu\rangle=\delta_{\lambda\mu}$.
\end{exm}
\begin{exm}
By constructing a continuous analog of $\Psi$ from the previous
problem prove that $(\chi_\lambda,\chi_\mu)=\delta_{\lambda\mu}$.
\end{exm}
\begin{remark}
It is a fact that the characters form a complete orthonormal basis
for the $L^2$ class functions on a compact Lie group. It follows
that the two norms are in fact the same.
\end{remark}
We will now describe the irreducible characters of $\text{U}(N)$ and
derive the Weyl character formula. Since every unitary matrix is
conjugate to a diagonal one it suffices to define any class function
just on the latter. Of course we have to make sure that if two
diagonal matrices are conjugate to each other our function takes the
same value on both.  Two diagonal matrices are conjugate in
$\text{U}(N)$ if and only if their diagonal entries differ by a
permutation. Hence any class function is symmetric in the
eigenvalues of matrices and conversely, any symmetric function can
be extended by conjugation to a class function on the entire group.

This means that we can describe characters as symmetric functions on
the eigenvalues of matrices. It is also useful to consider
alternating functions, that is, those that satisfy
$\omega(x_{\sigma(1)},\ldots,x_{\sigma(N)})=(-1)^\sigma\omega(x_1,\ldots,x_N)$
for any permutation $\sigma$. One such function of particular
interest to us is the Vandermonde determinant
$\delta_0(x):=\text{det}(x_j^{N-i})$ (the lower index is a
coordinate and the upper index is a power). The point is that ratios
(or products) of alternating functions are symmetric.
\begin{exm}\label{vandetform}
Show that $\text{det}(x_j^{N-i})=\prod_{i<j}(x_i-x_j)$.
\end{exm}

Now we are ready to introduce the symmetric functions that
correspond to the characters of $\text{U}(N)$. Set
\index{$\delta_\lambda(x):=\text{det}(x_j^{\lambda_i+N-i})$}
$\delta_\lambda(x):=\text{det}(x_j^{\lambda_i+N-i})$. Note that this
is an alternating polynomial and $\delta_0$ is the Vandermonde
determinant.
\index{$S$@$S_\lambda(x):={\delta_\lambda(x)}/{\delta_0(x)}$ Schur
polynomial}\index{Schur polynomial}
\begin{defn}
The ratios
$$
S_\lambda(x):=\frac{\delta_\lambda(x)}{\delta_0(x)}=\frac{\text{det}(x_j^{\lambda_i+N-i})}{\text{det}(x_j^{N-i})}
$$
are symmetric polynomials called the Schur polynomials.
\end{defn}
These were named in honor of I Schur who discovered their
connection to character theory in 1901.
\begin{exm}
Prove that these fractions are indeed polynomials.
\end{exm}
Just as the characters of irreducible representations are
orthonormal, the Schur polynomials are orthonormal:
$$
\langle S_\lambda,S_\mu\rangle=\delta_{\lambda\mu}\,.
$$
In fact we will see that the Schur polynomials evaluated on the
eigenvalues of a matrix are equal to the corresponding characters
applied to the matrix. Our proof will use the above orthogonality of
the Schur polynomials. A proof of this orthogonality independent of
the character formula can be found in Fulton and Harris \cite{fulton-harris}.

To prove that the Schur polynomials are the characters we will need
to use the complete symmetric polynomials defined by
$$
H_m(x):=\sum_{|\lambda|=m}x^\lambda.
$$
Alternatively, they can be described by a generating function:
\begin{exm}\label{Vandermonde}
Show that $\prod_{j=1}^N(1-tx_j)^{-1}=\sum_{m=0}^\infty H_m(x)t^m$.
\end{exm}
The following formula will prove very useful for manipulating
Vandermonde-type determinants.
\begin{lemma}\label{detf}
Let $a_i$ $i=1,\ldots,N$ be a decreasing sequence of integers and
set
$$ I_a=\{\,b\,|\,b_1\ge a_1>b_2\ge\cdots
>b_N\ge a_N\geq0\}\,.
$$ Then
$$
\text{det}(x_j^{a_i})\prod_{j=1}^N(1-x_j)^{-1}=\sum_{b\in
I_a}\text{det}(x_j^{b_i})\,.
$$
\end{lemma}

\begin{proof}
The proof uses induction on $N$. The case $N=1$ is standard. Using
cofactor expansion along the first row gives
\begin{align*}
\def\Sum{\sum\nolimits}
\def\Prod{\prod\nolimits}
\text{det}(x_j^{a_i})&\Prod_{j=1}^N(1{-}x_j)^{-1}\\
&=\Sum_{k=1}^N(-1)^{k+1}x_k^{a_1}(1{-}x_k)^{-1}\text{det}(x_j^{a_i})_{i\ne
1,j\ne
k}\Prod_{j\ne k}(1{-}x_j)^{-1}\\
&=\Sum_{k=1}^N(-1)^{k+1}\Sum_{b_1\ge a_1}\Sum_{b_2\ge a_2>\cdots}
\text{det}(x_j^{b_i})_{i\ne 1, j\ne k} \\
&=\Sum_{b_1\ge a_1, b_2\ge a_2>\cdots}\text{det}(x_j^{b_i})=
\Sum_{b\in I_{a}}\text{det}(x_j^{b_i}).
\end{align*}
The last equality follows by grouping terms with $b_1>b_2$ with
terms with $b_2>b_1$ when $b_2\ge a_1$.
\end{proof}

The next lemma employs the previous relation to prove the so-called
Pieri formula.
\begin{lemma}\label{perri}
$$
S_\lambda(x)H_m(x)=\sum_{\stackrel{\scriptstyle
|\nu|=m}{0\le\nu_i\le\lambda_{i-1}-\lambda_i}} S_{\lambda+\nu}(x)\,.
$$
\end{lemma}

\begin{proof}
We have
$$\begin{aligned}
\sum_{m=0}^\infty S_\lambda(x)H_m(x)t^m &=
\left(\text{det}(x_j^{\lambda_i+N-i})\prod_{j=1}^N(1-tx_j)^{-1}\right)/\text{det}(x_j^{N-i})\\
&=
t^{-|\lambda|}\left(\text{det}(tx_j^{\lambda_i+N-i})\prod_{j=1}^N(1-tx_j)^{-1}\right)/
\text{det}(tx_j^{N-i})\\
&=
t^{-|\lambda|}\left(\sum_\nu\text{det}(tx_j^{\lambda_i+\nu_i+N-i})\right)/
\text{det}(tx_j^{N-i})\\
&=
\left(\sum_\nu\text{det}(x_j^{\lambda_i+\nu_i+N-i})t^{|\nu|}\right)/
\text{det}(x_j^{N-i})\\
&=\sum_{0\le\nu_i\le\lambda_{i-1}-\lambda_i}S_{\lambda+\nu}(x)t^{|\nu|}\,.
\end{aligned}$$
Here $|\lambda|=\sum \lambda_i$, we used \fullref{detf} with
$a_i=\lambda_i+N-1$ and set $\nu_i=b_i-a_i$. The condition that
$b\in I_a$ is equivalent to $\nu$ being a partition with
$0\le\nu_i\le\lambda_{i-1}-\lambda_i$.
\end{proof}

\begin{remark} The Pieri formula has a very nice graphical
interpretation -- the product $S_\lambda H_m$ is the sum of Schur
polynomials with Young diagrams obtained from the Young diagram of
$\lambda$ in all ways of adding $m$ boxes to it with no two new
boxes in the same column. We use this interpretation later for
computations.
\end{remark}
We can now prove that Schur polynomials are the characters for the
representations having a Young diagram with just one row. This will
provide the base of induction for the general proof.
\begin{lemma}\label{smbox} Let the weight $m\omega_{1}$ also denote the corresponding Young
diagram with $m$ boxes in one row, the corresponding partition and
representation, then $S_{m\omega_{1}}=H_m=\chi_{m\omega_{1}}$.
\end{lemma}

\begin{proof}
The first equality comes from just applying the Pieri
formula with $\lambda=0$. To get the second equality note that a
basis of $V_{m\omega_{1}}$ is formed by the vectors
$c_{m\omega_{1}}(e_1^{\alpha_1}\otimes\cdots\otimes e_N^{\alpha_N})$
with $|\alpha|=m$. Here $c_{m\omega_{1}}$ is the Young symmetrizer
and $e_i$ form a basis in $V=\C^N$. Let $g\in\text{U}(N)$ act on
$\C^N$ with eigenvalues $x_j$ and eigenvectors $e_j$; then
$m\omega_{1}(g)c_{m\omega_{1}}(e_1^{\alpha_1}\otimes\cdots\otimes
e_N^{\alpha_N})=x^\alpha
c_{m\omega_{1}}(e_1^{\alpha_1}\otimes\cdots\otimes e_N^{\alpha_N})$.
It follows that
$$
\chi_{m\omega_{1}}(g)=\text{tr}({m\omega_{1}}(g))=\sum_{|\alpha|=m}
x^\alpha=H_m(x)\,,
$$
where $x^\alpha:=x_1^{\alpha_1}\ldots x_N^{\alpha_N}$.
\end{proof}

The Weyl character formula for U$(N)$ is the simple equality
$$\chi_\lambda=S_\lambda.$$
The general proof proceeds by induction on the number of boxes in
the Young diagram. \fullref{smbox} proves the equality for all
diagrams with one row. For example
$$\chi_{\tableau{4}}=S_{\tableau{4}}.$$
To illustrate the proof for more rows consider the case of four
boxes. Using the induction hypothesis and the Pieri formula
\eqref{perri} with $\lambda=3\tableau{1}$, $m=1$ and
$\lambda=2\tableau{1}$, $m=2$ respectively gives
\begin{align*}
\chi_{\tableau{3}}\chi_{\tableau{1}}&=S_{\tableau{4}}+S_{\tableau{3
1}}\\[-0.5ex]
\chi_{\tableau{2}}\chi_{\tableau{2}}&=S_{\tableau{4}}+S_{\tableau{3
1}}+S_{\tableau{2 2}}.
\end{align*}
Next, using the Pieri formula \eqref{perri} with
$\lambda=2\tableau{1}$ and $m=1$, multiplying the result by
$\chi_{\tableau{1}}=H_1$, and then using the Pieri formula again
gives the next formula. Repeating the process further gives the
formula after that.
\begin{align*}
\chi_{\tableau{2}}\chi_{\tableau{1}}\chi_{\tableau{1}}&=S_{\tableau{4}}
+2S_{\tableau{3 1}}+S_{\tableau{2 2}}+S_{\tableau{2 1 1}}\\[-0.5ex]
\chi_{\tableau{1}}\chi_{\tableau{1}}\chi_{\tableau{1}}\chi_{\tableau{1}}&=
S_{\tableau{4}}+3S_{\tableau{3 1}}+2S_{\tableau{2 2}}+3S_{\tableau{2
1 1}}+S_{\tableau{1 1 1 1}}\,.
\end{align*}
Since a product of characters is the character of the tensor
product, which in turn can be written as a sum of irreducible
representations, it follows that all our character products can be
written as sums of irreducible characters. For example we can write
$\chi_{\tableau{3}}\chi_{\tableau{1}}=\sum_\mu n_\mu\chi_\mu$ for
some numbers $n_\mu$. Since we need these numbers for different
character products we have to double index them. It is convenient to
use as the second index the diagram obtained by stacking the boxes
in the product one under another. For instance, the coefficients of
the above sum should be denoted $n_{\tableau{3 1}\,\mu}$.

Notice that the representation $V_{\tableau{3 1}}$ is a subspace of
$V_{\tableau{3}}\otimes V_{\tableau{1}}$ by definition and has to
appear at least once in the irreducible decomposition. It follows
that $n_{\tableau{3 1}\,\tableau{3 1}}\ge 1$. Writing the previous
five displayed equations in the matrix form will make the notation
and the argument even clearer.
$$
\left[\begin{array}{ccccc} 1&0&0&0&0\\1&1&0&0&0\\1&1&1&0&0\\
1&2&1&1&0\\1&3&2&3&1
\end{array}\right]\left[\begin{array}{c}S_{\tableau{4}}\\S_{\tableau{3
1}}\\
S_{\tableau{2 2}}\\ S_{\tableau{2 1 1}}\\ S_{\tableau{1 1 1
1}}\end{array}\right]= \left[\begin{array}{c}\sum n_{\tableau{4}\,\mu}\chi_\mu \\
\sum n_{\tableau{3 1}\,\mu}\chi_\mu \\\sum n_{\tableau{2 2}\,\mu}\chi_\mu \\
\sum n_{\tableau{2 1 1}\,\mu}\chi_\mu \\\sum n_{\tableau{1 1 1
1}\,\mu}\chi_\mu\end{array}\right].
$$
Inverting the matrix implies that we may write $S_\lambda=\sum
m_{\lambda\mu}\chi_\mu$ for some integers $m_{\lambda\mu}$. Taking
the norm gives
$$
1=\langle S_\lambda,S_\lambda\rangle=\sum m_{\lambda\mu}^2\,.
$$
Working from the top down we see that $m_{\lambda\lambda}\ge 1$ at
one level implies that $m_{\lambda\mu}=\delta_{\lambda\mu}$ and
implies that $m_{\lambda\lambda}\ge 1$ holds at the next level down.
It follows that $\chi_\lambda=S_\lambda$. There is a similar matrix for the Schur functions with $\lambda$
having any fixed arbitrary number of boxes.

So far we only considered irreducible characters of $\text{U}(N)$
but now we can get the $\text{SU}(N)$ case for free. Recall that
irreducible representations of $\text{SU}(N)$ are indexed by
partitions of length $<N$. But we can treat such a partition as
being of length $\leq N$ by setting $\lambda_N=0$. This produces an
irreducible representation of $\text{U}(N)$ that restricts to the
original one of $\text{SU}(N)$. Since the characters are just traces
the same formula gives them for $\text{SU}(N)$ as long as we treat
$\lambda$ as an $N$--vector with $\lambda_N=0$.
\begin{thm}[Weyl character formula for $\text{SU}(N)$]\label{uWeyl}
Let $\lambda$ denote any dominant weight (partition) of length $\leq
N$ and the corresponding irreducible representation of
$\text{SU}(N)$. Then
$$
\chi_{\lambda}=S_{\lambda}=\frac{\delta_\lambda}{\delta_0}
$$
\end{thm}
\begin{exm}
The coefficients $n^\nu_{\lambda\mu}$ defined by
$\chi_\lambda\chi_\mu=\sum^\nu_{\lambda\mu}n^\nu_{\lambda\mu}\chi_\nu$ are the
classical Racah coefficients. Compute $n^\nu_{3\omega_1\omega_1}$.
\end{exm}
There is a version of the Weyl character formula valid for any
semisimple Lie algebra. Even for $\text{SU}(N)$ the formula is more
useful to us in this general form. This also gives us a chance to
introduce notation that comes in handy when considering quantum
groups. Introduce a parameter $q=e^{iz}$, $z\in\C$, then the
expression $q^\beta:=\exp(iz\beta)$ is defined for all complex
numbers and complex-valued matrices, for example by the power
series. If $\beta\in\mathfrak{g}$ for some matrix Lie algebra then
$q^\beta$ is in the corresponding Lie group and on weight vectors of
a representation its action is given by $q^\beta
v_\gamma=q^{\langle\gamma,\beta\rangle}v_\gamma$. The Weyl
denominator \index{$\delta_\lambda(q^\beta)$ Weyl denominator}
\index{Weyl denominator} is defined to be
\begin{equation}\label{Wdenom}
\delta_\lambda(q^\beta):=\sum_{w\in W}(-1)^w q^{\langle w(\lambda+\rho),\beta\rangle}\,.
\end{equation}\index{Weyl character formula}

\begin{thm}[Weyl character formula]\label{Weyl}
The character $\chi_\lambda$ of an irreducible representation of a
semisimple Lie group is given by
\begin{equation}\label{Wform}
\chi_\lambda=\frac{\delta_\lambda}{\delta_0}.
\end{equation}
\end{thm}
\begin{proof}[Proof for $\text{SU}(N)$] It suffices to
consider diagonal matrices $\beta=\sum_{i=1}^{N}\beta_iE_{ii}$ with
real entries and the trace  $\sum_{i=1}^{N}\beta_i=0$. It is
straightforward to check that $\beta_i=\langle L_i,\beta\rangle$ and
$q^\beta=\text{diag}(x_1,\ldots,x_N)\in\text{SU}(N)$, where we put
$x_i:=q^{\beta_i}$. The Weyl denominator becomes
\begin{align*}
\delta_\lambda(q^\beta):=&\Sum_{\sigma\in\mathfrak{G}_N}(-1)^\sigma
q^{\langle \sigma(\lambda+\rho),\beta\rangle}\\
=&\Sum_{\sigma\in\mathfrak{G}_N}(-1)^\sigma
q^{\sum_{i=1}^N(\lambda_i+\rho_i)\langle \sigma(L_i),\beta\rangle} \\
=&\Sum_{\sigma\in\mathfrak{G}_N}(-1)^\sigma\Prod_{i=1}^N
  \left(q^{\langle(L_{\sigma(i)},\beta\rangle}\right)^{\lambda_i+\rho_i}\\
=&\Sum_{\sigma\in\mathfrak{G}_N}(-1)^\sigma
  x_{\sigma(1)}^{\lambda_1+\rho_1}\ldots
  x_{\sigma(N)}^{\lambda_N+\rho_N} \\
=&\text{det}(x_j^{\lambda_i+\rho_i}) \\
=&\text{det}(x_j^{\lambda_i+N-i}).
\end{align*}
The last equality is the result of \fullref{Ww}. Our claim now
follows directly from \fullref{uWeyl}.
\end{proof}

As an application let us compute $\chi_\lambda(q^{2\rho})$, where as
before $\rho$ is the Weyl weight. This quantity turns up as a
`quantum dimension' in \eqref{alcove}. First, for any $\beta$
we have
\begin{multline}\label{dqrho}
\delta_0(q^{2\beta})=\sum_{w\in\mathfrak{G}_N}(-1)^w q^{(w(\rho),2\beta)}
=\det\left(q^{(L_j,2\beta)(N-i)}\right)\\
=\prod_{i<j}\left(q^{(\beta,L_i-L_j)}-q^{-(\beta,L_i-L_j)}\right)
=\prod_{\alpha\in\Delta^+}(q^{(\beta,\alpha)}-q^{-(\beta,\alpha)}).
\end{multline}
The third equality is the Vandermonde determinant identity from \fullref{Vandermonde}.
Now specifically for $q^{2\rho}$:
\begin{multline}\label{dqlrho}
\delta_\lambda(q^{2\rho})=\sum_{w\in\mathfrak{G}_N}(-1)^w q^{(w(\lambda+\rho),2\rho)}
=\sum_{w\in\mathfrak{G}_N}(-1)^w q^{(2(\lambda+\rho),w(\rho))}\\
=\delta_0(q^{2(\lambda+\rho)})=\prod_{\alpha\in\Delta^+}(q^{(\lambda+\rho,\alpha)}-q^{-(\lambda+\rho,\alpha)}).
\end{multline}
Combining the last two formulas we prove the following.
\begin{cor} Let $\chi_\lambda$ be the character
of the irreducible representation of $\Sl_N\C$ with the highest
weight $\lambda$ and $\rho$ be the Weyl weight. Then
\begin{equation}\label{qchar}
\chi_\lambda(q^{2\rho})=\frac{\delta_\lambda(q^{2\rho})}{\delta_0(q^{2\rho})}
=\prod_{\alpha\in\Delta^+}\,\frac{q^{(\lambda+\rho,\alpha)}-q^{-(\lambda+\rho,\alpha)}}{q^{(\rho,\alpha)}-q^{-(\rho,\alpha)}}\,.
\end{equation}
\end{cor}

\section{Exact invariants from conformal field theory}\label{cft}
\setobjecttype{App}
\label{app:e}

There is a simple idea that leads to the mathematical definition of
conformal field theory. To understand this idea, consider the loop
group of a Lie group. The loop group is the space of all maps from
$S^1$ to the group. Now consider all formal power series with values
in the associated Lie algebra. An element of this `loop algebra'
produces an element of the loop group by a two step process. First
consider the formal variable to be an element of $S^1$, so the
formal power series is a map from $S^1$ to the Lie algebra. Second,
exponentiate the answer to obtain a map from $S^1$ to the group. In
this way problems in infinite dimensional geometry may be translated
into problems in algebra.

Let ${\mathbb C}((t)):=\{\sum_{k=-M}^\infty a_kt^k|a_k\in{\mathbb
C}\}$, and ${\mathbb C}\llbracket t\rrbracket :=\{\sum_{k=0}^\infty
a_kt^k|a_k\in{\mathbb C}\}$; then the affine Lie algebra associated
to ${\mathfrak g}$ is defined to be
$$
\what {\mathfrak g}:= ({\mathfrak g}\otimes {\mathbb C}((t)))\oplus
{\mathbb C}K,
$$
with bracket
$$
[A\otimes f+\alpha K,B\otimes g+\beta K]:= [A,B]\otimes fg+\langle
A,B\rangle\text{Res}_{t=0}(f^\prime g)K.
$$
The affine Lie algebra may be written as a direct sum of subalgebras,
$$
\what {\mathfrak g}= ({\mathfrak g}\otimes t{\mathbb C}\llbracket t\rrbracket )\oplus
{\mathfrak g}\oplus
{\mathbb C}K \oplus ({\mathfrak g}\otimes t^{-1}{\mathbb C}[t^{-1}]).
$$
Any representation of ${\mathfrak g}$, say $V$, is a $({\mathfrak
g}\otimes t{\mathbb C}\llbracket t\rrbracket )\oplus {\mathfrak g}\oplus {\mathbb
C}K$--module with $({\mathfrak g}\otimes t{\mathbb C}\llbracket t\rrbracket )$ acting
as zero, ${\mathfrak g}$ acting as usual, and ${\mathbb C}K$ acting
as $k\,\text{id}$. The number $k$ is called the level. We will
always assume that the level is a positive integer. If $R\subseteq
S$ and $M$ is an $R$--module, then $\text{Ind}_R^SM:=M\otimes_R S$ is
an $S$--module. The Weyl module at level $k$ with highest weight
$\lambda$ is defined to be $\text{Ind}_{({\mathfrak g}\otimes
t{\mathbb C}\llbracket t\rrbracket )\oplus {\mathfrak g}\oplus {\mathbb
C}K}^{\what
{\mathfrak g}}V_\lambda$ and is denoted by $V^k_\lambda$. Restrict
attention to ${\mathfrak sl}_N{\mathbb C}$ for simplicity. An
integrable module at level $k$ is an ${\mathfrak{sl}}_N{\mathbb
C}$--module, so that $K$ acts by multiplication by $k$ and
$(E_{ii}-E_{jj})\otimes t^n$ acts locally nilpotently. It turns out
that the category of level $k$ integrable modules ${\mathcal
O}_k^{int}$ is a strict modular category. To define the flip
$\times_{V,W}$, one uses a creative way to attach integrable modules
to a Riemann surface at a number of points to obtain what is called
the space of conformal blocks. This gives rise to the bundle of
conformal blocks over the moduli space of marked Riemann surfaces of
genus $g$. This bundle admits a projectively flat connection, and
one can solve the parallel transport (a system of differential
equations called the Knizhnik--Zamolodchikov equations) to see the
effect of interchanging points marked with $V$ and $W$ to define the
flip. This is not the obvious structure and the resulting invariants
are far from the obvious ones. See Bakalov and Kirillov \cite{BK}, Di
Francesco, Mathieu and S\'en\'echal \cite{DMS} and Kohno \cite{kohno}
for more information.

Before moving on it is worthwhile to describe the simple objects in
this category. The affine Weyl group at level $k$ is the semidirect
product $W^a_k:=W\ltimes k\Lambda_r^\vee$ acting on ${\mathfrak
h}^*$ by $(s,\beta^\vee)(\gamma):=s(\gamma)+\beta^\vee$, where $W$
is the Weyl group of the associated Lie algebra. The interior of a
fundamental domain for this action is given by
\begin{equation}\label{colors}
I=\{\lambda\in {\mathfrak h}^*|\langle\lambda+\rho,\alpha_i^\vee\rangle>0, \ \langle\lambda+\rho,\theta^\vee\rangle<k+h^\vee\}.
\end{equation}
The simple objects in ${\mathcal O}_k^{int}$ are in one to one correspondence with the weights in $I$. Given a weight $\lambda\in I$,
$$
L_\lambda^k:=V_\lambda^k:=V^k_\lambda/(U(\what {\mathfrak g})(v_\theta\otimes t^{-1})^{k-\langle \lambda,\theta^\vee\rangle+1}v_{\lambda,k}.
$$
Here $U(\cdot)$ denotes the universal enveloping algebra
and $v_\lambda,k$ is the highest weight vector of $V^k_\lambda$.

Even though many physicists approach the Witten--Chern--Simons invariants
through conformal field theory, we will approach these invariants via
quantum groups since it is the fastest way to supply the definitions. It
has been shown that the two approaches agree (see Tsuchiya, Ueno and
Yamada \cite{TUY} and Bakalov and Kirillov \cite{BK}).

\newpage

\input main_a.ind

\end{document}